\documentclass[11pt,twoside]{article}
\DeclareFontFamily{U}{matha}{\hyphenchar\font45}
\DeclareFontShape{U}{matha}{m}{n}{
      <5> <6> <7> <8> <9> <10> gen * matha
      <10.95> matha10 <12> <14.4> <17.28> <20.74> <24.88> matha12
      }{}
\DeclareSymbolFont{matha}{U}{matha}{m}{n}
\DeclareFontFamily{U}{mathx}{\hyphenchar\font45}
\DeclareFontShape{U}{mathx}{m}{n}{
      <5> <6> <7> <8> <9> <10>
      <10.95> <12> <14.4> <17.28> <20.74> <24.88>
      mathx10
      }{}
\DeclareSymbolFont{mathx}{U}{mathx}{m}{n}

\DeclareMathSymbol{\obot}{2}{matha}{"6B}
\DeclareMathSymbol{\bigobot}{1}{mathx}{"CB}
\usepackage{amsthm,epsfig,a4}

\usepackage{verse} %Gedichte
\usepackage[ansinew]{inputenc}
\usepackage{amsmath,amssymb}
\let\originallhook\lhook
\usepackage{MnSymbol}
\newcommand{\lhook}{\mathrel{\raise.018ex\hbox{$\originallhook$}}}
\usepackage{amsfonts}
\usepackage{latexsym}
\usepackage{esdiff}
\usepackage{oldgerm} %\textfrak{Fraktur Buchstaben} \textgoth{gotische Schrift} \textswab{Schwabacher Schrift}

\newfont{\suet}{suet14}
\newfont{\schwell}{schwell}
\DeclareTextFontCommand{\textsuet}{\suet}
\DeclareTextFontCommand{\textschwell}{\schwell}

\usepackage{lscape}
\usepackage{framed} %Boxen um Objekte
\usepackage{tikz}
\usetikzlibrary{arrows.meta} % für erweiterte Pfeiltypen
\usepackage{makeidx} %Index
\makeindex
\def\:{\colon\thinspace}
\def\q{\nolinebreak \hfill $\Box$} %q.e.d.%
\def\A{\\[3mm] \indent} %Beweisende Neuer Absatz%
\def\D{\displaystyle}
\def\N{\mathbb{N}}
\def\Z{\mathbb{Z}}
\def\Q{\mathbb{Q}}
\def\R{\mathbb{R}}
\def\C{\mathbb{C}}
\def\K{\mathbb{K}}
\def\L{\mathcal{L}}
\def\B{\mathbb{B}}
\DeclareMathSymbol{\Bbbk}{\mathord}{AMSb}{"7C}
\def\k{\Bbbk}
\def\U{\mathcal{U}^{\circ}\!}
\def\d{{\rm d}}
\def\e{{\rm e}}
\def\id{{\rm id}}
\def\inv{{\rm inv}}
\def\SOS{\overline{\mathfrak{S}}}
\def\SUS{\underline{\mathfrak{S}}}
\def\SZk{\mathfrak{Z}_k}
\def\OIf{\int\limits_0^{\underline{1}} f(t) \, \d t}
\def\UIf{\int\limits_{\overline{0}}^1 f(t) \, \d t}
\def\If{\int\limits_0^1 f(t) \, \d t}
\usepackage{german,t1enc}
\renewcommand{\i}{{\rm i}}
\theoremstyle{plain}
\newtheorem*{Lemma*}{Lemma}

\newtheorem{Lemma}{Lemma}[section]
\newtheorem*{Satz*}{Satz}
\newtheorem{Satz}[Lemma]{Satz}
\newtheorem{HS}[Lemma]{Hauptsatz}
\newtheorem*{Kor*}{Korollar}
\newtheorem{Kor}[Lemma]{Korollar}
\newtheorem*{Zusatz}{Zusatz}
\newtheorem{UA}{Aufgabe}
\theoremstyle{definition}
\newtheorem{Nr}[Lemma]{}
\newtheorem{Def}[Lemma]{Definition}
\newtheorem*{Def*}{Definiton}
\newtheorem*{Bsp*}{Beispiel}
\newtheorem{Bsp}[Lemma]{Beispiel}
\newtheorem{Bem}[Lemma]{Bemerkung}
\newtheorem*{Bem*}{Bemerkung}

\allowdisplaybreaks
\begin{document}
\pagenumbering{roman}
\begin{titlepage}
\vspace*{15mm}
\begin{center}
{\huge {\fontencoding{U} \fontfamily{yinit} \fontseries{m} \fontshape{n} \selectfont F}unktionalanalysis} \\
{\huge {\fontencoding{U} \fontfamily{yinit} \fontseries{m} \fontshape{n} \selectfont T}eil I}
\end{center}
\vspace*{35mm}
\begin{center}
Christoph Bock\\
\small{Department Mathematik\\ Universität Erlangen-Nürnberg\\ Cauerstraße 11\\ 91058 Erlangen\\ E-Mail: \verb+bock(at)mi.uni-erlangen.de+}
\end{center}
\vfill
\begin{center}
%\today
\end{center}
\end{titlepage}
\thispagestyle{empty}
\cleardoublepage
\newpage
\thispagestyle{empty}

$\,$
\\[9cm]
\begin{center}
{\bfseries Für meinen Vater}
\end{center}
\thispagestyle{empty}
\clearpage
\newpage
\thispagestyle{empty}
$\,$
\newpage
\section*{Abstract}
%
%
%\pagenumbering{roman}
%
%
Roughly speaking, \emph{functional analysis} is the study of vector spaces of arbitrary dimension over the field $\K \in \{\R,\C\}$, and the continuous linear mappings between such spaces.
Naturally, the notion of continuity requires a topology -- or more specifically, a norm -- on these vector spaces, bringing both analytic and algebraic tools into play.

The name functional analysis originates from the early attempts to extend calculus to functionals defined on function spaces.
Results from functional analysis offer powerful tools for solving problems in the theory of (partial) differential equations, in complex analysis, and for the formulation of quantum mechanics.
However, the aim of these pages is not to deal with such applications.

This first part is primarily concerned with the intrinsic properties of certain classes of spaces, namely almost metric spaces, normed vector spaces and algebras, spaces of continuous and and $p$-integrable functions (for $p \in {]}0, \infty{]}$), as well as reflexive, uniformly convex, and Hilbert spaces, rather than with the study of mappings between them.

If you notice any errors, please contact the author at the e-mail address \verb+bock(at)mi.uni-erlangen.de+.

These pages will be updated and extended from time to time, and the latest version can always be found on the author’s website via the following link \verb+https://www.mathematik.uni-erlangen.de/~bock/FAna_Teil_I.pdf+.
Please note that the numbering may still change compared to the current version.
\cleardoublepage
\newpage %\thispagestyle{empty}
\section*{Vorwort}
%
%
%\pagenumbering{roman}
%
%
\emph{Funktionalanalysis} bedeutet, grob gesagt, die Untersuchung beliebig-di\-men\-sio\-na\-ler Vektorräume und der stetigen linearen Abbildungen zwischen solchen, wobei der Begriff der Stetigkeit natürlich eine Topologie oder etwas spezieller eine Norm benötigt.
Der Name \emph{Funktionalanalysis} rührt daher, daß in den Anfängen der Theorie die Analysis auf Funktionale von Funktionenräumen ausgeweitet wurde.
Funktionalanalytische Resultate ergeben Möglichkeiten, Probleme der \emph{(Partiellen) Differentialgleichungen}, der \emph{Funktionentheorie} zu lösen und die \emph{Quantenmechanik} zu formulieren.
Ich verzichte hier allerdings größtenteils, auf die Anwendungen einzugehen.

Kapitel 1 behandelt fastmetrische Räume und topologische Eigenschaften jener.
Soweit es möglich ist, führe ich die Theorie auf dem Niveau topologischer Räume.
Bzgl.\ vieler Ergebnisse muß man sich allerdings auf den Fall fastmetrischer Räume beschränken.
Höhepunkt des Kapitels sind der Nachweis der Existenz einer bis auf Isometrie eindeutigen Vervollständigung (fast-)metrischer Räume, eine Charakterisierung kompakter Teilmengen fastmetrischer Räume sowie die Sätze von \textsc{Arzel\`{a}-Ascoli} und von \textsc{Baire}.
Das Kapitel kann auch als eine kleine Einführung in die Grundstrukturen der mengentheoretischen Topologie verstanden werden, wobei der Begriff des fastmetrischen Raumes im Mittelpunkt steht.

Im zweiten Kapitel werden normierte Vektorräume und Algebren betrachtet, welche zusammen mit ihren stetigen linearen Abbildungen die eigentlichen Objekte der Funktionalanalysis bilden.
Herausragende Ergebnisse sind die Existenz einer bis auf Isometrie eindeutigen Vervollständigung normierter Vektorräume bzw.\ Algebren zu sogenannten Banachräumen bzw.\ -algebren, eine Charakterisierung endlich-dimensionaler normierter Vektorräume, die Sätze von \textsc{Hahn-Banach}, zu deren Beweis das Lemma von \textsc{Zorn}, dessen Beweis in Anhang \ref{FAnaA} zu finden ist, benötigt wird, und die fundamentalen Prinzipien der gleichmäßigen Beschränktheit und der offenen Abbildung.

Die wichtigsten Räume in der Funktionalanalysis sind Räume stetiger bzw.\ integrierbarer Funktionen, denen die drei folgenden Kapitel gewidmet sind.
Die herausragenden Ergebnisse des dritten Kapitels sind, daß die im unendlichen verschwindenden stetigen Funktionen, die auf lokal-kom\-pak\-ten Hausdorff-Räumen definiert sind, in den stetigen Funktionen mit kompaktem Träger bzgl.\ der Supremumsnorm dicht liegen, wofür eine Version des Lemmas von \textsc{Urysohn}, deren Beweis in Anhang \ref{FAnaB} zu finden ist, benötigt wird, und der Satz von \textsc{Stone-Weierstraß}.
Kapitel 4 stellt zunächst eine leistungsfähige Integrationstheorie zur Verfügung, leistungsfähig in dem Sinne, daß das dargestellte Lebesguesche Integral, welches nach \textsc{H.\ L.\ Lebesgue} benannt ist, unter gewissen Grenzwertprozessen abgeschlossen ist.
Sodann werden im fünften Kapitel die Lebesgueschen Räume $L^p$ der $p$-integrierbaren Funktionen eingeführt und u.a.\ die Ungleichungen von \textsc{Hölder} und \textsc{Minkowski} sowie der Satz von \textsc{Riesz-Fischer} über die Vollständigkeit der Lebesgueschen Räume bewiesen.

Ein bedeutendes funktionalanalytisches Prinzip besteht darin, die Untersuchung eines normierten Vektorraumes mit dem Studium seines topologischen Dualraumes zu verbinden.
In diesem ist die Einheitsvollkugel bzgl.\ der durch die Operatornorm induzierten Topologie allerdings nur im endlich-dimensionalen Falle kompakt, und dies erfordert die Entwicklung einer speziellen Methode für die Funktionalanalysis.
In Kapitel \ref{FAna6} wird u.a.\ die sog.\ schwache-$*$-Topologie eingeführt und gezeigt, daß die Einheitsvollkugel des topologischen Dualraumes eines normierten Vektorraumes bzgl.\ dieser stets kompakt ist -- d.i.\ der Satz von \textsc{Banach-Alaoglu}, dessen Beweis den in Anhang \ref{FAnaC} präsentierten Satz von \textsc{Tychonoff} verwendet.
Mit diesem Rüstzeug kann ich im weiteren Verlauf des Kapitels reflexive Räume charakterisieren, wobei ein normierter Vektorraum genau dann reflexiv heißt, wenn seine kanonische Einbettung in den topologischen Bidualraum surjektiv ist.
An einer Stelle wird hierbei der Satz von \textsc{Eberlein-\v{S}mulian} wesentlich ausgenutzt und deswegen skizzenhaft bewiesen.
Eine der Charakterisierungen ermöglicht es außerdem, zu zeigen, daß die Approximationsaufgabe in reflexiven Räumen mindestens eine Lösung besitzt.

Das siebente Kapitel behandelt Räume, in deren Vervollständigung die Approximationsaufgabe eindeutig lösbar ist, nämlich gleichmäßig konvexe Räume.
Es werden die Sätze von \textsc{Clarkson} und \textsc{Milman}, die mit dem Satz von \textsc{Riesz-Fischer} die Reflexivität der Lebesgueschen Räume $L^p$ für $p \in {]}1, \infty{[}$ ergeben, sowie der Darstellungssatz von \textsc{Riesz} in der Version für Lebesguesche Räume, der im Falle $p \in {[}1, \infty{[}$ den topologischen Dualraum eines Lebesgueschen Raumes $L^p$ bis auf Isometrie als seinen konjugierten Lebesgueschen Raum $L^q$ angibt, welcher durch $\frac{1}{p} + \frac{1}{q} = 1$ charakterisiert ist, bewiesen.

Im achten Kapitel werden mit Hilberträumen spezielle gleichmäßig konvexe Räume studiert.
Jene sind per definitionem Banachräume, deren Norm durch ein Skalarprodukt induziert wird.
Letzteres ermöglicht die Definition des geometrischen Aspektes der Orthogonalität und einen Hilbertraum in einen beliebigen abgeschlossenen Untervektorraum und dessen orthogonales Komplement, welches ebenfalls ein abgeschlossener Untervektorraum ist, zu spalten.
Eine solche orthogonale Spaltung wird verwendet, um den Darstsellungssatz von \textsc{Riesz} in der Version für Hilberträume zu beweisen, welcher besagt, daß diese in kanonischer Weise selbstdual sind.
Des weiteren wird der Begriff einer Hilbertbasis eingeführt und gezeigt, daß jeder Hilbertraum bis auf Isometrie durch die Mächtigkeit einer solchen Basis eindeutig bestimmt ist.

I.d.R.\ werden dem Leser zum Abschluß eines Kapitels Übungsaufgaben gestellt, die bewußt gelassene Lücken in vorherigen Beweisführungen schließen oder die Theorie vervollständigen.
\A
Ich erhebe keinen Anspruch auf Originalität.
Bei der Erstellung der Kapitel dienten mir vor allem die Vorlesungen \cite{Reck}, die ich während meines Studiums bei \textsc{H.\ Reckziegel} gehört und von denen ich hier entscheidend profitiert habe, sowie das Buch \cite{Hirz}, das auf Grundlage einer Vorlesung von \textsc{F.\ Hirzebruch} geschrieben wurde, als Quelle.
Des weiteren fanden die Bücher von \textsc{K.\ Floret} \cite{Flo}, \textsc{R.\ V.\ Kadison} und \textsc{J.\ R.\ Ringrose} \cite{Kad}, \textsc{S.\ Lang} \cite{Lang} sowie \textsc{W.\ Rudin} \cite{Rud} Verwendung.
Bei topologischen Resultaten habe ich in \textsc{J.\ R.\ Munkres}' Buch \cite{Munkres} Anleihe genommen.
Die Lebesguesche Integrationstheorie, die auf \textsc{F.\ Riesz} sowie \textsc{B.\ v.\ Sz.\ Nagy} zurückgeht und von\linebreak \textsc{P.\ Dombrowski} erheblich verallgemeinert wurde, entstammt einer Mitschrift der Vorlesungen \cite{Dom} des letztgenannten sowie den Vorlesungen \cite{Henke} meines Diplomvaters \textsc{W.\ Henke}, bei dem ich die Theorie erlernt habe; in der Fachliteratur findet sich diese leider kaum.
Während der Ausarbeitung des vierten Kapitels bin ich auf das Buch \cite{Poeschel} von \textsc{J.\ Pöschel} aufmerksam geworden, und auch dieses hatte Einfluß auf vorliegende Seiten.
Ich habe auch weitere Literatur, die mir teilweise nicht mehr im einzelnen präsent ist, gelesen und bei meiner Darstellung und Beweisführung verwendet.
\A
Für den Austausch über die Materie sowie Literaturhinweise bin ich den Herren \textsc{J.\ Bock}, \textsc{M.\ Bohn}, \textsc{K.\ Conrad}, \textsc{P.\ Dombrowski} (\dag \, 2025), \textsc{W.\ Henke}, \textsc{M.\ A.\ Nieper-Wißkirchen}, \textsc{J.\ Pöschel}, \textsc{H.\ Reckziegel} und \textsc{P.\ Schwahn} dankbar. 
\A
Es sei angemerkt, daß ich mir sicher bin, daß mir Fehler unterlaufen sind.
Für Hinweise auf solche oder Kritik bin ich dankbar.
Kontakten können Sie mich per E-Mail an \verb+bock(at)mi.uni-erlangen.de+.
Bei der vorliegenden Ver\-sion zur \emph{Funktionalanlysis I} handelt es sich nur um ein Preprint, das nicht gegengelesen wurde und im Laufe der Zeit noch verändert und ergänzt wird -- insbesondere kann sich die Numerierung ändern.
Die jeweils aktuelle Version finden Sie unter \verb+https://www.mathematik.uni-erlangen.de/~bock/FAna_Teil_I.pdf+. 

\begin{flushleft}
Erlangen, im Herbst 2025 \hfill \raisebox{-1ex}{\schwell{Christoph Bock}}\\
\end{flushleft}
\cleardoublepage
\newpage %\thispagestyle{empty}
\tableofcontents
\cleardoublepage
\section{Fastmetrische Räume} \label{FAna1}
\pagenumbering{arabic}
\subsection*{Grundlagen} \addcontentsline{toc}{subsection}{Grundlagen}
\begin{Def}[(Fast-)Metrische Räume] \label{FA.1.1} 
Sei $X$ eine Menge.
\begin{itemize}
\item[(i)] Eine Abbildung $d \: X \times X \to [0, \infty]$ nennen wir eine \emph{Fastmetrik auf $X$}\index{Fastmetrik} genau dann, wenn gilt
\begin{eqnarray*}
({\rm D}1) & \forall_{x,y \in X} \, d(x,y) = 0 \Longleftrightarrow x = y, & \mbox{(Positiv-Definitheit),}\\
({\rm D}2) & \forall_{x,y \in X} \, d(x,y) = d(y,x),& \mbox{(Symmetrie),}\\
({\rm D}3) & \forall_{x,y,z \in X} \, d(x,z) \le d(x,y) + d(y,z), & \mbox{(Dreiecksungleichung),}\index{Ungleichung!Dreiecks-}
\end{eqnarray*}
wobei wir $t + \infty := \infty + t := \infty$ für alle $t \in [0, \infty]$ setzen.
Das Paar $(X,d)$ heißt dann ein \emph{fastmetrischer Raum}\index{Raum!fastmetrischer}.

\item[(ii)] Ist $d$ wie in (i) sogar eine Abbildung $X \times X \to {[}0, \infty{[}$, so heißt $d$ eine \emph{Metrik auf $X$}\index{Metrik (siehe auch Fastmetrik)} und das Paar $(X,d)$ ein \emph{metrischer Raum}\index{Raum!metrischer}.
\end{itemize}

\begin{Bem*} 
Oft werden wir einen (fast-)metrischen Raum mit einem einzigen Symbol -- z.B.\ $X$ -- bezeichnen.
Wir schreiben dann $\boxed{|X|}$ (oder einfach $X$) für die $X$ zugrundeliegende Menge und $\boxed{d_X}$ (oder einfach $d$) für die Metrik von $X$, also
$$ X = ( \underbrace{|X|}_{= X}, \underbrace{d_X}_{= d} ).$$
\end{Bem*}

\begin{Bsp*} $\,$
\begin{itemize}
\item[1.)] Eine Teilmenge eines (fast-)metrischen Raumes ist in kanonischer Weise ebenfalls ein (fast-)metrischer Raum. 
Wir versehen sie im folgenden stets mit dieser induzierten (Fast-)Metrik.
\item[2.)] Sind $\K \in \{\R,\C\}$, $n \in \N_+$ sowie $k \in \{1,2\}$, so werden durch
$$ \forall_{x = (x_1, \ldots,x_n), y = (y_1, \ldots, y_n) \in \K^n} \, \boxed{d_k(x,y)} := \left( \sum_{i=1}^n |x_i - y_i|^k \right)^{\frac{1}{k}} $$
und
$$ \forall_{x = (x_1, \ldots,x_n), y = (y_1, \ldots, y_n) \in \K^n} \, \boxed{d_{\infty}(x,y)} := \max \, \{ |x_i - y_i| ~ | ~ i \in \{1, \ldots, n\} \} $$
bekanntlich Metriken auf $\K^n$ definiert.

Im Falle $n=1$ gilt $\boxed{d_{\K}} := d_1 = d_2 = d_{\infty}$.
$d_{\K}$ heißt \emph{Standardmetrik für $\K$}.
\item[3.)] Seien $M$ eine nicht-leere Menge und $X$ ein fastmetrischer Raum.
Dann wird durch
$$ \forall_{f,g \in X^M} \, \boxed{d_{\infty}(f,g)} := \sup \{ d(f(p),g(p)) \, | \, p \in M \} $$
eine Fastmetrik, die sog.\ \emph{Supremumsfastmetrik\index{Fastmetrik!Supremums-} auf $X^M$}, definiert.

{[} (M1) und (M2) sind trivial. (M3) folgt daraus, daß für alle $f,g,h \in X^M$ und $p \in M$
$$ d(f(p),g(p)) \le  d(f(p),h(p)) + d(h(p),g(p)) \le d_{\infty}(f,h) + d_{\infty}(h,g) $$
gilt. {]}

Im Falle $M = \{1, \ldots, n\}$ und $X = \K$, wobei $n \in \N_+$ und $\K \in \{\R, \C\}$, identifiziert man $n$ mit $\{1, \ldots, n\}$, und die Definitionen von $d_{\infty}$ in 2.) und 3.) simmen überein.
\item[4.)] Sind $n \in \N_+$ und $(X_1,d_1), \ldots, (X_n,d_n)$ (fast-)metrische Räume, so wird durch
$$ \forall_{x=(x_1, \ldots, x_n), y=(y_1, \ldots, y_n) \in \bigtimes_{i=1}^n X_i} \, d(x,y) := \sum_{i=1}^n d_i(x_i,y_i) $$
offenbar eine Metrik, die sog.\ \emph{Produkt(fast-)metrik\index{Produkt(fast-)metrik}\index{Fastmetrik!Produkt-} auf $\bigtimes_{i=1}^n X_i$}, definiert. 
\end{itemize}
\end{Bsp*}
\end{Def}

\begin{Lemma} \label{FA.1.L} 
Sei $X$ ein metrischer Raum. Dann gilt für alle $x, \tilde{x}, y, \tilde{y} \in X$
$$ |d(x, \tilde{x}) - d(y, \tilde{y})| \le d(x,y) + d(\tilde{x}, \tilde{y}). $$
\end{Lemma}

\textit{Beweis.} Aus der Dreiecksungleichung folgt
$$ d(x, \tilde{x}) \le d(x,y) + d(y, \tilde{y}) +  d(\tilde{y}, \tilde{x}), $$
also
$$ d(x, \tilde{x}) - d(y, \tilde{y})  \le d(x,y) +  d(\tilde{x}, \tilde{y}) $$
und ebenso
$$ d(y, \tilde{y}) - d(x, \tilde{x})  \le d(x,y) +  d(\tilde{x}, \tilde{y}). $$
\q

\begin{Nr}[Fastmetrische Räume als topologische Räume] \label{FA.1.2} 
Sei $X$ ein fastmetrischer Raum.
\begin{itemize}
\item[(i)] Für jedes $x \in X$ und jedes $\varepsilon \in \R_+$ heißt die Menge
$$ \boxed{U_{\varepsilon}(x)} := \{ y \in X \, | \, d(x,y) < \varepsilon \} $$
die \emph{$\varepsilon$-Umgebung von $x$}.\index{Umgebung!$\varepsilon$-}

Es gilt:
\begin{itemize}
\item[(a)] $\forall_{x \in X} \forall_{\varepsilon \in \R_+} \, x \in U_{\varepsilon}(x)$,
\item[(b)] $\forall_{x \in X} \forall_{\varepsilon_1, \varepsilon_2 \in \R_+} \, U_{\varepsilon_1}(x) \cap U_{\varepsilon_2}(x) = U_{\min \{ \varepsilon_1, \varepsilon_2 \}}(x)$,
\item[(c)] $\forall_{x,y \in X} \forall_{\varepsilon \in \R_+} \, y \in U_{\varepsilon}(x) \Longrightarrow U_{\varepsilon - d(x,y)}(y) \subset U_{\varepsilon}(x)$.
\end{itemize}

{[} (a),(b) sind trivial, und (c) folgt aus der Dreiecksungleichung. {]}
\item[(ii)] Eine Teilmenge $U$ von $X$ heißt \emph{offen in $X$}\index{Menge!offene} genau dann, wenn gilt
$$ \forall_{x \in U} \exists_{\varepsilon \in \R_+} \, U_{\varepsilon}(x) \subset U.$$

\begin{Bsp*} 
Sind $x \in X$ und $\varepsilon \in \R_+$, so ist $U_{\varepsilon}(x)$ nach (i) (c) offen.
\end{Bsp*}
\pagebreak
\item[(iii)] Sei
$$\mathcal{T} := \left\{ U \subset X \, | \, U \mbox{ offen in } X \right\} \subset \mathfrak{P}(X).$$
Dann gilt offenbar
\begin{eqnarray*}
({\rm T}1) && \emptyset, X \in \mathcal{T}, \\
({\rm T}2) && \left( I \mbox{ beliebige Menge} \, \wedge \, \forall_{i \in I} \, U_i \in \mathcal{T} \right) \Longrightarrow \bigcup_{i \in I} U_i \in \mathcal{T}, \\ 
({\rm T}3) && U_1, U_2 \in \mathcal{T} \Longrightarrow U_1 \cap U_2 \in \mathcal{T}.
\end{eqnarray*}
\end{itemize}
\end{Nr}

\begin{Def}[Topologische Räume, offene Mengen, abgeschlossene Mengen, Umgebungen, Hausdorff-Räume] \label{FA.1.3} $\,$
\begin{itemize}
\item[(i)] Ein Paar $(X,\mathcal{T})$, bestehend aus einer Menge $X$ zusammen mit einer Teilmenge $\mathcal{T}$ von $\mathfrak{P}(X)$, das (T1) - (T3) erfüllt, heißt ein \emph{topologischer Raum}\index{Raum!topologischer}. 
$\mathcal{T}$ nennt man dann \emph{Topologie für $X$}\index{Topologie} und die Elemente von $\mathcal{T}$ die \emph{offenen Mengen\index{Menge!offene} des topologischen Raumes $(X,\mathcal{T})$}.

\begin{Bem*} 
Auch einen topologischen Raum bezeichnen wir oft mit einem einzigen Symbol, z.B.\ $X$.
Wir schreiben dann wieder $\boxed{|X|}$ (oder einfach $X$) für die $X$ zugrundeliegende Menge und $\boxed{{\rm Top}(X)}$ für die Topologie von $X$, also
$$ X = ( \underbrace{|X|}_{= X}, {\rm Top}(X) ).$$
\end{Bem*}

\begin{Bsp*} 
Ist $X$ eine Menge, so heißt $\{ \emptyset, X \}$ die \emph{triviale Topologie für $X$}\index{Topologie!triviale} und $\mathfrak{P}(X)$ die \emph{diskrete Topologie für $X$}\index{Topologie!diskrete}.
\end{Bsp*}

\begin{Bem*} $\,$
\begin{itemize}
\item[1.)] Ist $X$ ein (fast-)metrischer Raum, so wird durch \ref{FA.1.2} (iii) eine Topologie auf $X$ definiert, die sog.\ \emph{(fast-)metrische Topologie}\index{Topologie!(fast-)metrische}.
Jeder (fast-)metrische Raum ist also in kanonischer Weise ein topologischer Raum. 
Wir betrachten jeden (fast-)metrischen Raum als topologischen Raum mit dieser kanonischen Topologie.
\item[2.)] In \cite[Chapter 6]{Munkres} werden notwendige und hinreichende Bedingungen dafür angegeben, daß ein topologischer Raum \emph{metrisierbar}\index{Raum!topologischer!metrisierbarer} ist, d.h.\ eine Metrik besitzt so, daß die metrische Topologie mit der gegebenen übereinstimmt.
\item[3.)] Zu jeder Fastmetrik existiert eine Metrik derart, daß die kanonischen Topologien übereinstimmen, vgl.\ \ref{FA.1.18} unten.
\end{itemize}
\end{Bem*}
\item[(ii)]
Sei $X$ ein topologischer Raum.
\begin{itemize}
\item[(a)] Eine Teilmenge $A$ von $X$ heißt \emph{abgeschlossen in $X$}\index{Menge!abgeschlossene} genau dann, wenn $X \setminus A$ offen ist.

Aus (T1) - (T3) folgt leicht, daß $\emptyset, X$, die Vereinigung endlich vieler abgeschlossener Mengen und der Schnitt beliebig vieler abgeschlossener Mengen abgeschlossen ist.

\item[(b)] Ist $x \in X$, so heißt eine Teilmenge $U$ von $X$ \emph{Umgebung von $x$ in $X$}\index{Umgebung} genau dann, wenn $U$ offen in $X$ ist und $x \in U$ gilt.

Mit $\boxed{\U(x,X)}$ bezeichnen wir die Menge aller Umgebungen von $x$ in $X$.

\item[(c)] $X$ heißt \emph{Hausdorff-Raum}\index{Raum!Hausdorff-} oder \emph{hausdorffsch} genau dann, wenn zu $x,y \in X$ mit $x \neq y$ Umgebungen $U \in \U(x,X)$ und $V \in \U(y,X)$ existieren derart, daß $U \cap V = \emptyset$ gilt.
\end{itemize}
\end{itemize}
\end{Def}

\begin{Bem*} 
Für topologische Räume hat man sog.\ \emph{Trennungsaxiome}\index{Trennungs!-axiome} formuliert, die in einem gewissen Sinne Auskunft darüber geben, ,,wie viele`` offene Mengen es gibt, die Punkte oder allgemeiner Teilmengen voneinander ,,trennen``.
Das \emph{zweite Trennungsaxiom} ist genau die Hausdorff-Eigenschaft, die auch die \emph{Punktetrennungseigenschaft}\index{Punktetrennungseigenschaft}\index{Trennungs!-eigenschaft von Punkten} genannt wird.
\end{Bem*}

\begin{Satz} \label{FA.1.Punkt} 
Ist $X$ ein Hausdorff-Raum, so ist $\{x\}$ für jedes $x \in X$ abgeschlossen in $X$.
\end{Satz}

\textit{Beweis.} Sei $x \in X$. 
Dann existieren zu jedem $y \in X \setminus \{x\}$ Umgebungen $U_y \in \U(x,X)$ und $V_y \in \U(y,X)$ mit $U_y \cap V_y$, und es gilt $X \setminus \{x\} = \bigcup_{y \in X \setminus \{x\}} V_y$, welches nach (T2) eine offene Menge ist. \q

\begin{Bem*} 
Das \emph{erste Trennungsaxiom}\index{Trennungs!-axiome} für einen topologischen Raum besagt, daß die einelementigen Teilmengen abgeschlossen sind.
\end{Bem*}

\begin{Satz} \label{FA.1.4} 
Ist $X$ ein fastmetrischer Raum, so ist $X$ als topologischer Raum hausdorffsch.
\end{Satz}

\textit{Beweis.} Seien $x,y \in X$ mit $x \ne y$.

1. Fall: $d(x,y) < \infty$.
Wir setzen $\varepsilon := \frac{1}{2} d(x,y)$.
Dann gilt 
$$ U_{\varepsilon}(x) \cap U_{\varepsilon}(y) = \emptyset, $$
da andernfalls $z \in X$ mit $d(x,z) < \varepsilon$ und $d(y,z) < \varepsilon$, also 
$$ 2 \varepsilon = d(x,y) \le d(x,z) + d(z,y) < 2 \varepsilon, $$ 
existierte. 

2. Fall: $d(x,y) = \infty$. 
Dann gilt für jedes $\varepsilon \in \R_+$
$$ U_{\varepsilon}(x) \cap U_{\varepsilon}(y) = \emptyset, $$
da andernfalls $z \in X$ mit $d(x,z) < \varepsilon$ und $d(y,z) < \varepsilon$, also 
$$ \infty = d(x,y) \le d(x,z) + d(z,y) < 2 \varepsilon, $$ 
existierte. 
\q

\begin{Def}[Teilraumtopologie] \label{FA.1.TT}
Seien $X$ ein topologischer Raum und $Y$ eine Teilmenge von $X$ sowie $G$ eine Teilmenge von $Y$.

$G$ heißt \emph{offen im Teilraum $Y$ von $X$}\index{Menge!offene -- einer Teilmenge} genau dann, wenn eine in $X$ offene Menge $U$ existiert mit $$G = U \cap Y.$$

\begin{Satz*} 
Seien $X$ ein topologischer Raum und $Y$ eine Teilmenge von $X$.
Dann ist 
$$ \boxed{{\rm Top}_X(Y)} := \left\{ G \subset Y \, | \, \mbox{$G$ offen im Teilraum $Y$ von $X$} \right\} $$
eine Topologie für $Y$, die sog.\ \emph{Teilraumtopologie\index{Topologie!Teilraum-} von $Y$ bzgl.\ $X$}.
Sofern keine Verwechselungen auftreten können, schreiben wir auch $\boxed{{\rm Top}(Y)}$ anstelle von ${\rm Top}_X(Y)$.
\end{Satz*}

\textit{Beweis als Übung.} \q

\begin{Bem*}
Wir betrachten eine Teilmenge eines topologischen Raumes stets als topologischen Raum mit der durch den umgebenden Raum induzierten Teilraumtopologie.
\end{Bem*}
\end{Def}

\begin{Satz} \label{FA.1.TT.1}
Seien $X$ ein fastmetrischer Raum und $Y$ eine Teilmenge von $X$.

Dann stimmt die kanonische Topologie von $(Y,d|_{Y \times Y})$ mit der Teilraumtopologie von $Y$ bzgl.\ $X$ überein.
\end{Satz}

\textit{Beweis als Übung.} \q

\begin{Satz} \label{FA.1.TT.2}
Seien $X$ ein hausdorffsch topologischer Raum und $Y \subset X$.

Dann ist der topologische Teilraum $Y$ von $X$ hausdorffsch.
\end{Satz}

\textit{Beweis.} Seien $y, \tilde{y} \in Y \subset X$ mit $y \ne \tilde{y}$.
Da nach Voraussetzung offene Mengen $U, V \subset X$ mit $y \in U$, $ \tilde{y} \in V$ und $U \cap V = \emptyset$ existieren, folgt für die in $Y$ offenen Mengen $U \cap Y$ und $ V \cap Y$:
$y \in U \cap Y$, $\tilde{y} \in V \cap Y$ und 
$$ (U \cap Y) \cap (V \cap Y) = (U \cap V) \cap Y = \emptyset. $$
\q

\begin{Nr}[Konvergenz von Folgen in topologischen Räumen] \index{Konvergenz} \index{Folgen!konvergente} \label{FA.1.5} $\,$
\\[2mm]
\noindent \textbf{Definition 1.}
Seien $X$ ein topologischer Raum, $(x_n)_{n \in \N}$ eine \emph{Folge in $X$} -- d.i.\ per definitionem ein Element von $X^{\N}$ -- und $x \in X$.
Wir definieren dann:

\emph{$(x_n)_{n \in \N}$ konvergiert (in $X$) gegen $x$} $:\Longleftrightarrow \forall_{U \in \U(x,X)} \exists_{n_0 \in \N} \forall_{n \in \N, n \ge n_0} \, x_n \in U$.

\begin{Satz*}
Sind $X$ ein topologischer Raum, $Y \subset X$, $(y_n)_{n \in \N}$ eine Folge in $Y$ und $y \in Y$, so konvergiert $(y_n)_{n \in \N}$ offenbar genau dann in dem topologischen Teilraum $Y$ von $X$ gegen $y$, wenn $(y_n)_{n \in \N}$ in $X$ gegen $y$ konvergiert.
\end{Satz*}

\textit{Beweis als Übung.} \q

\begin{Lemma*}
Sind $X$ ein Hausdorff-Raum und $\left( x_n \right)_{n \in \N}$ eine Folge in $X$, die sowohl gegen $x \in X$ als auch gegen $\tilde{x} \in X$ konvergiert, so gilt $x = \tilde{x}$.
\end{Lemma*}

\textit{Beweis.} Angenommen, $x \ne \tilde{x}$. 
Nach Voraussetzung existieren $U \in \U(x,X)$ und $V \in \U(\tilde{x},X)$ mit $U \cap V = \emptyset$. 
Weiter existieren $n_1, n_2 \in \N$ mit 
$$ \forall_{n \in \N, n \ge n_1} \, x_n \in U ~~ \mbox{ und } ~~ \forall_{n \in \N, n \ge n_2} \, x_n \in V,$$
also folgt für $n_0 := \max \{n_1,n_2\}$: $x_{n_0} \in U \cap V$, Widerspruch! \q

\begin{Bem*}
Das Lemma ist i.a.\ falsch, wenn $X$ nicht hausdorffsch ist.
Ist z.B.\ $X$ ein topologischer Raum, der mit der trivialen Topologie $\{ \emptyset, X \}$ versehen ist, so konvergiert jede Folge in $X$ gegen jedes Element von $X$.
\end{Bem*}

\noindent \textbf{Definition 2.} Sind daher $X$ ein Hausdorff-Raum und $\left( x_n \right)_{n \in \N}$ eine Folge in $X$, die gegen $x \in X$ konvergiert, so ist $x$ mit dieser Eigenschaft (nach dem Lemma) eindeutig bestimmt, und wir setzen $$\boxed{\lim_{n \to \infty} x_n} := x$$ und nennen $x$ den \emph{Grenzwert}\index{Grenzwert} oder \emph{Limes\index{Limes} von $\left( x_n \right)_{n \in \N}$ für $n$ gegen unendlich}.
\end{Nr}

\begin{Satz} \label{FA.1.KA}
Seien $A$ eine abgeschlossene Teilmenge eines topologischen Raumes $X$ und $(x_n)_{n \in \N}$ eine Folge in $A$.

Ist dann $x \in X$ derart, daß $(x_n)_{n \in \N}$ in $X$ gegen $x$ konvergiert, so gilt $x \in A$.
\end{Satz}

\textit{Beweis.} Da $A$ abgeschlossen ist, ist $X \setminus A$ offen. 
Angenommen, $x \in X \setminus A$. 
Dann existiert $n_0 \in \N$ derart, daß $\forall_{n \in \N, \, n \ge n_0} \, x_n \in X \setminus A$, im Widerspruch zur Voraussetzung. \q

\begin{Satz}[Konvergenz von Folgen in fastmetrischen Räumen] \index{Konvergenz} \index{Folgen!konvergente} \label{FA.1.6}
Sei $X$ ein fastmetrischer Raum mit Fastmetrik $d$ (also ist $X$ nach \ref{FA.1.4} hausdorffsch).

Dann gilt für alle Folgen $\left( x_n \right)_{n \in \N}$ in $M$ und alle $x \in X$:
\begin{eqnarray} \label{FA.1.6.S}
\lim_{n \to \infty} x_n = x & \Longleftrightarrow & \forall_{\varepsilon \in \R_+} \exists_{n_0 \in \N} \forall_{n \in \N, n \ge n_0} \, \underbrace{x_n \in U_{\varepsilon}(x)}_{\mbox{$\Leftrightarrow d(x_n,x) < \varepsilon$}}.
\end{eqnarray} 
Die rechte Seite von (\ref{FA.1.6.S}) besagt gerade, daß $\lim_{n \to \infty} d(x_n,x) = 0$ bzgl.\ der Standardmetrik $d_{\R}$ auf $\R$ gilt.
\end{Satz}

\textit{Beweis.} ,,$\Rightarrow$`` ist klar, da jede $\varepsilon$-Umgebung nach \ref{FA.1.2} offen ist.

Zu ,,$\Leftarrow$``: Sei $U \in \U(x,X)$. 
Nach Definition der Topologie von $X$ (vgl.\ \ref{FA.1.2}) existiert dann $\varepsilon \in \R_+$ mit $U_{\varepsilon}(x) \subset U$. 
Nach Voraussetzung gilt daher für fast alle $n \in \N$: $x_n \in U_{\varepsilon}(x) \subset U$. \q

\begin{Bsp} \label{FA.1.7.dsup}
Seien $M$ eine nicht-leere Menge, $(X,d)$ ein fastmetrischer Raum, $(f_n)_{n \in \N}$ eine Folge in $X^M$ und $f \in X^M$.
Dann gilt:

$(f_n)_{n \in \N}$ konvergiert genau dann in $(X^M, d_{\infty})$ gegen $f$, wenn sie im üblichen Sinne \emph{gleichmäßig konvergent gegen $f$} ist\index{Konvergenz!gleichmäßige -- einer Folge von Abbildungen}\index{Folgen!konvergente!gleichmäßig -- von Abbildungen}, d.h.\ per definitionem
$$ \forall_{\varepsilon \in \R_+} \exists_{n_0 \in \N} \forall_{n \in \N, \, n \ge n_0} \forall_{p \in M} \, d(f_n(p),f(p)) < \varepsilon. $$

\textit{Beweis als Übung.} \q
\end{Bsp}

\begin{Bem*}
Wir erinnern daran, daß man eine Folge $(f_n)_{n \in \N}$ in $X^M$, wobei $M$ eine nicht-leere Menge und $X$ ein fastmetrischer Raum sei, \emph{punktweise konvergent gegen $f \in X^M$}\index{Konvergenz!gleichmäßige -- einer Folge von Abbildungen}\index{Folgen!konvergente!gleichmäßig -- von Abbildungen} nennt, wenn gilt
$$ \forall_{p \in M} \forall_{\varepsilon \in \R_+} \exists_{n_0 \in \N} \forall_{n \in \N, \, n \ge n_0} \, d(f_n(p),f(p)) < \varepsilon. $$
\end{Bem*}

\begin{Def}[offener Kern, abgeschlossene Hülle, Rand] \label{FA.1.7} 
Seien $X$ ein topologischer Raum und $A \subset X$.
\begin{itemize}
\item[(i)] Wir definieren den \emph{offenen Kern\index{offener Kern} von $A$ in $X$} als
\begin{equation*}
\boxed{A^{\circ}} := \bigcup_{U \in {\rm Top}(X), \, U \subset A} U.
\end{equation*}

Nach (T2) ist $A^{\circ}$ offen in $X$, also ist $A^{\circ}$ offenbar die größte offene Teilmenge von $X$, die in $Y$ enthalten ist.
\item[(ii)] Die \emph{abgeschlossene Hülle von $A$ in $X$}\index{abgeschlossene Hülle}\index{Hülle!abgeschlossene} ist definiert als
\begin{equation*}
\boxed{\overline{A}} := \bigcap_{A' \subset X {\rm abgeschlossen}, \, A' \supset A} A'.
\end{equation*}

Nach \ref{FA.1.3} (ii) (a) ist $\overline{A}$ abgeschlossen in $X$, also ist $\overline{A}$ offenbar die kleinste abgeschlossene Teilmenge von $X$, die $A$ umfaßt.
\item[(iii)] Der \emph{Rand von $A$ in $X$}\index{Rand} ist definiert als
$$ \boxed{\partial A} := \overline{A} \setminus A^{\circ}. $$
\end{itemize}

\begin{Bsp*}
Sei $X := [-1,1] \cup \{\i\} \subset \C$ mit der durch $d_{\C}$ induzierten Metrik.
Dann gilt
$$ U_1(0) \subsetneqq \overline{U_1(0)} \subsetneqq B_1(0) := \{ x \in X \, | \, d_{\C}(x,0) \le 1 \}, $$
da $U_1(0) = {]}-1,1{[}$, $\overline{U_1(0)} = [-1,1]$ und $B_1(0) = X$.
\end{Bsp*}
\end{Def}

\begin{Satz} \label{FA.1.8}
Seien $X$ ein topologischer Raum und $A \subset X$.
Dann gilt:
\begin{itemize}
\item[(i)] $A^{\circ} = \{ x \in X \, | \, \exists_{U \in \U(x,X)} \, U \subset A \}$.

Einen Punkt $x \in X$ mit $\exists_{U \in \U(x,X)} \, U \subset A$ nennen wir einen \emph{inneren Punkt\index{innerer Punkt} von $A$ in $X$}.
\item[(ii)] $\overline{A} = \{ x \in X \, | \, \forall_{U \in \U(x,X)} \, U \cap A \neq \emptyset \}$. 

Einen Punkt $x \in X$ mit $\forall_{U \in \U(x,X)} \, U \cap A \neq \emptyset$ nennen wir einen \emph{Berührungspunkt von $A$ in $X$}\index{Berührungspunkt}.
\item[(iii)] $\partial A = \{ x \in X \, | \, \forall_{U \in \U(x,X)} \, U \cap A \neq \emptyset \, \wedge \, U \cap (X \setminus A) \neq \emptyset \} = \overline{A} \cap \overline{X \setminus A}$.
Insbesondere ist $\partial A$ abgeschlossen in $X$.

Ein Element von $\partial A$ nennen wir einen \emph{Randpunkt\index{Rand!-punkt} von $A$ in $X$}.
\end{itemize}
\end{Satz}

\textit{Beweis.} Zu (i): ,,$\subset$`` Sei $x \in A^{\circ}$. 
Dann existiert nach \ref{FA.1.7} (i) ein $U \in {\rm Top}(X)$ mit $U \subset A$ und $x \in U$, also $U \in \U(x,X)$ mit $U \subset A$, d.h.\ $x$ ist innerer Punkt von $A$ in $X$.

,,$\supset$`` Sei $x$ ein innerer Punkt von $A$ in $X$. 
Dann existiert $U \in \U(x,X)$ mit $U \subset A$, also $x \in A^{\circ}$.

Zu (ii): ,,$\subset$`` Sei $x \in \overline{A}$.
Angenommen, $x$ ist kein Berührungspunkt von $A$ in $X$, d.h. es existiert $U \in \U(x,X)$ mit $U \cap A = \emptyset$. 
Dann ist $A' := X \setminus U$ eine abgeschlossene Teilmenge von $X$ mit $A' \supset A$, folglich nach \ref{FA.1.7} (ii): $\overline{A} \subset A'$, insbesondere wegen $x \in \overline{A}$: 
$x \in A' = X \setminus U$, im Widerspruch zu $U \in \U(x,X)$.

,,$\supset$`` Sei $x$ Berührungspunkt von $A$ in $X$.
Zu zeigen ist $x \in \overline{A}$, d.h.\ für jede abgeschlossene Teilmenge $A'$ von $X$ mit $A' \supset A$ gilt $x \in A$.
Angenommen, dies ist falsch, d.h.\ es existiert eine abgeschlossene Teilmenge $A'$ von $X$ mit $A' \supset A$ und $x \notin A'$.
Dann ist $U := X \setminus A'$ offen in $X$, und es gilt $x \in U$, also $U \in \U(x,X)$, folglich (da $x$ Berührungspunkt von $A$ in $X$): 
$U \cap A \ne \emptyset$, im Widerspruch zu $U = X \setminus A'$ und $A \subset A'$. 

Zu (iii): Die Behauptung folgt sofort aus (i) und (ii). \q

\begin{Satz} \label{FA.1.8.L}
Seien $A,B$ Teilmengen eines topologischen Raumes $X$
Dann gilt:
\begin{itemize}
\item[(i)] $\left(X \setminus A \right)^{\circ} = X \setminus \overline{A} ~~ \mbox{ und } ~~ \overline{X \setminus A} = X \setminus A^{\circ}$,
\item[(ii)] $A \subset B \Longrightarrow \left( A^{\circ} \subset B^{\circ} \, \wedge \, \overline{A} \subset \overline{B} \right)$,
\item[(iii)] $(A \cap B)^{\circ} = A^{\circ} \cap B^{\circ}$,
\item[(iv)] $(A \cup B)^{\circ} \supset A^{\circ} \cup B^{\circ}$,
\item[(v)] $\overline{A \cap B} \subset \overline{A} \cap \overline{B}$,
\item[(vi)] $\overline{A \cup B} = \overline{A} \cup \overline{B}$.
\end{itemize}
Die Inklusionen in (iv) und (v) sind i.a.\ echt.
\end{Satz}

\textit{Beweis als Übung.} \q

\begin{Satz} \label{FA.1.9}
Seien $X$ ein fastmetrischer Raum, $A \subset X$ und $x \in X$.

Dann ist $x$ genau dann ein Berührungspunkt von $A$ in $X$, wenn eine Folge $(x_n)_{n \in \N}$ in $A$ mit $\lim_{n \to \infty} x_n = x$ existiert.
\end{Satz}

\textit{Beweis.} ,,$\Rightarrow$`` Ist $x$ ein Berührungspunkt von $A$ in $X$, so können wir zu jedem $n \in \N$ ein $x_n \in U_{\frac{1}{n+1}}(x) \cap A$ wählen.
Dann gilt offenbar $\lim_{n \to \infty} x_n = x$.

,,$\Leftarrow$`` Ist $(x_n)_{n \in \N}$ eine Folge in $A$ mit $\lim_{n \to \infty} x_n = x$, so liegen in jeder Umgebung $U \in \U(x,X)$ fast alle Folgenglieder, d.h.\ $U \cap A \neq \emptyset$. \q

\begin{Def}[Cauchyfolgen, vollständige fastmetrische Räume] \label{FA.1.10}
Sei $X$ ein fastmetrischer Raum.
\begin{itemize}
\item[(i)] Eine Folge $(x_n)_{n \in \N}$ in $X$ heißt \emph{Cauchyfolge (in $X$)}\index{Cauchy!-folge}\index{Folgen!Cauchy-} genau dann, wenn gilt
$$ \forall_{\varepsilon \in \R_+} \exists_{n_0 \in \N} \forall_{n,m \in \N, \, n,m \ge n_0} \, d(x_n,x_m) < \varepsilon. $$
\begin{Bsp*}
Konvergiert eine Folge $(x_n)_{n \in \N}$ in $X$ gegen $x \in X$, so ist $(x_n)_{n \in \N}$ eine Cauchyfolge:

Denn zu $\varepsilon \in \R_+$ existiert $n_0 \in \N$ mit $\forall_{n \in \N, n \ge n_0} \, d(x_n,x) < \frac{\varepsilon}{2}$, also gilt für alle $n,m \in \N$ mit $n,m \ge n_0$
$$ d(x_n,x_m) \le d(x_n,x) + d(x, x_m) < \varepsilon. $$
\end{Bsp*}
\item[(ii)] $X$ (bzw.\ $d$) heißt \emph{vollständig}\index{Raum!fastmetrischer!vollständiger}\index{Menge!vollständige}\index{Vollständigkeit}\index{Fastmetrik!vollständige} genau dann, wenn jede Cauchyfolge in $X$ gegen ein Element von $X$ konvergiert.
\end{itemize}
\end{Def}

\begin{Lemma} \label{FA.1.21.K}
Seien $X$ ein fastmetrischer Raum und $A$ eine Teilmenge von $X$.

Dann ist $\overline{A}$ genau dann vollständig, wenn jede Cauchyfolge in $A$ in $X$ konvergiert.
\end{Lemma}

\textit{Beweis.} ,,$\Rightarrow$`` ist trivial.

,,$\Leftarrow$`` Sei $(x_n)_{n \in \N}$ eine Cauchyfolge in $\overline{A}$, d.h.\ eine Folge von Berührungspunkten von $A$.
Dann gilt
$$ \forall_{n \in \N} \, U_{\frac{1}{n+1}}(x_n) \cap A \ne \emptyset, $$
und wir wählen zu jedem $n \in \N$ ein 
\begin{equation} \label{FA.1.21.K.1}
y_n \in U_{\frac{1}{n+1}}(x_n) \cap A.
\end{equation}

Wir behaupten:
\begin{equation} \label{FA.1.21.K.2}
(y_n)_{n \in \N} \mbox{ ist Cauchyfolge.}
\end{equation}

{[} Zu (\ref{FA.1.21.K.2}): Sei $\varepsilon \in \R_+$.
Da $(x_n)_{n \in \N}$ eine Cauchyfolge ist, existiert $n_0 \in \N$ mit
$$ \forall_{n,m \in \N, \, n,m \ge n_0} \, d(x_n, x_m) < \frac{\varepsilon}{3}, $$
und wir können ohne Einschränkung annehmen, daß $\frac{1}{n_0 + 1} < \frac{\varepsilon}{3}$ gilt.
Dann folgt für alle $n,m \ge n_0$
$$ d(y_n, y_m) \le \underbrace{d(y_n, x_n)}_{\stackrel{(\ref{FA.1.21.K.1})}{<} \frac{1}{n+1} < \frac{\varepsilon}{3}} + \underbrace{d(x_n, x_m)}_{< \frac{\varepsilon}{3}} + \underbrace{d(x_m, y_m)}_{\stackrel{(\ref{FA.1.21.K.1})}{<} \frac{1}{m+1} < \frac{\varepsilon}{3}} < \varepsilon, $$
also gilt (\ref{FA.1.21.K.2}). {]}

Aus (\ref{FA.1.21.K.2}) und der Voraussetzung der rechten Seite folgt, daß $(y_n)_{n \in \N}$ gegen ein gewisses $x \in X$ konvergiert.
Da für jedes $n \in \N$ gilt
$$ d(x_n,x) \le \underbrace{d(x_n,y_n)}_{\stackrel{(\ref{FA.1.21.K.1})}{<} \frac{1}{n+1}} + \underbrace{d(y_n,x)}_{\stackrel{n \to \infty}{\longrightarrow} 0} \stackrel{n \to \infty}{\longrightarrow} 0, $$
konvergiert auch $(x_n)_{n \in \N}$ gegen $x$. \q
\A
Ein zweiter Beweis von ,,$\Leftarrow$`` nutzt die Existenz einer Vervollständigung $\widehat{X}$ von $X$ aus, s.u.\ \ref{FA.1.21}:

Es gilt für jedes $x \in \widehat{X}$
\begin{eqnarray*}
\lefteqn{\mbox{$x$ Berührungspunkt von $A$ in $\widehat{X}$}} ~~~ \\
& \stackrel{\ref{FA.1.9}}{\Longleftrightarrow} & \exists_{(y_n)_{n \in \N} \subset A^{\N}} \, \lim_{n \to \infty} y_n = x \\
& \stackrel{\text{Vor. d. r.S. für ,,$\Rightarrow$``}}{\Longleftrightarrow} & \underbrace{\mbox{$x$ Berührungspunkt von $A$ in $X$}}_{\Rightarrow x \in X},
\end{eqnarray*}
also stimmen die abgeschlossenen Hüllen von $A$ in $\widehat{X}$ und von $A$ in $X$ überein.
Da $\widehat{X}$ vollstänig ist, ergibt sich aus dem folgenden Satz \ref{FA.1.11} (i) die Vollständigkeit von $\overline{A}$. 

\begin{Satz} \label{FA.1.11}
Es seien $X$ ein fastmetrischer Raum und $A$ eine Teilmenge von $X$.
\begin{itemize}
\item[(i)] Ist $X$ vollständig und $A$ abgeschlossen in $X$, so ist $A$ (bzgl.\ $d|_{A \times A}$) vollständig.
\item[(ii)] Ist $A$ (bzgl.\ $d|_{A \times A}$) vollständig, so ist $A$ abgeschlossen in $X$.
\end{itemize}
\end{Satz}

\textit{Beweis.} Zu (i): Sei $(x_n)_{n \in \N}$ eine Cauchyfolge in $(A, d|_{A \times A})$, d.h.\ auch eine Cauchyfolge in $X$.
Wegen der Vollständigkeit von $X$ existiert $x \in X$, gegen das $(x_n)_{n \in \N}$ in $X$ konvergiert.
$A$ ist abgeschlossen in $X$, also ergibt Satz \ref{FA.1.9}: $x \in A$.

Zu (ii): Wir zeigen $\overline{A} \subset A$.
Hierzu sei $x \in \overline{A}$.
Dann existiert nach Satz \ref{FA.1.9} eine gegen $x \in \overline{A} \subset X$ konvergente Folge $(x_n)_{n \in \N}$ in $A$.
Folglich ist $(x_n)_{n \in \N}$ eine Cauchyfolge in $(A, d|_{A \times A})$, also nach Voraussetzung konvergent gegen ein gewisses $y \in A$.
Wegen der Eindeutigkeit des Grenzwertes in $X$ gilt $y = x$, d.h.\ $x \in A$. \q

\begin{Bem*}
Die Umkehrung der Aussage (ii) des letzten Satzes gilt i.a.\ nicht, wie jeder nicht-vollständige fastmetrische Raum als Teilmenge von sich zeigt.
\end{Bem*}

\begin{Def}[Isometrien, isometrische fastmetrische Räume] \label{FA.1.Isom}
Seien $X, Y$ fastmetrische Räume, wobei beide Fastmetriken mit $d$ bezeichnet seien, sowie $f \: X \to Y$ eine Abbildung.
\begin{itemize}
\item[(i)] $f$ heißt genau dann \emph{isometrisch}\index{Abbildung!isometrische}, wenn gilt
$$ \forall_{x, \tilde{x} \in X} \, d(f(x),f(\tilde{x})) = d(x,\tilde{x}). $$

\begin{Bem*} 
Isometrische Abbildungen $X \to Y$ sind offenbar injektiv.
\end{Bem*}

\item[(ii)] $f$ heißt genau dann eine \emph{Isometrie}\index{Isometrie!fastmetrischer Räume}, wenn $f$ isometrisch und surjektiv ist.

\item[(iii)] $X$ und $Y$ heißen \emph{isometrisch}\index{Raum!fastmetrischer!isometrischer} genau dann, wenn eine Isometrie $X \to Y$ existiert.
\end{itemize}
\end{Def}

\begin{Satz} \label{FA.1.Isom.1}
Es seien $M$ eine nicht-leere Menge und $(X,d)$ ein fastmetrischer Raum.
Für jedes $x \in X$ bezeichne $\underline{x} \: M \to X$ die konstante Abbildung vom Wert $x$.
Des weiteren sei $ \iota \: X \to X^M$ gegeben durch
$$ \forall_{x \in X} \, \iota(x) := \underline{x}. $$
\begin{itemize}
\item[(i)] $\iota \: (X,d) \to (X^M, d_{\infty})$ ist isometrisch, und $\iota(X)$ ist eine abgeschlossene Teilmenge von $(X^M,d_{\infty})$. 
\item[(ii)] $(X^M, d_{\infty})$ ist genau dann vollständig, wenn $(X,d)$ vollständig ist.
\end{itemize}
\end{Satz}

\textit{Beweis.}  Zu (i): Daß $\iota$ isometrisch ist, ist klar.
Wir zeigen, daß $\iota(X)$ abgeschlossen ist:
Hierzu sei $f \in X^M \setminus \iota(X)$, d.h.\ $f$ ist nicht konstant.
Dann existieren $p,q \in M$ mit $f(p)\ne f(q)$, und wir setzen 
$$ \varepsilon := \frac{1}{2} d(f(p),f(q)) > 0. $$
Ist nun $g \in U_{\varepsilon}^{d_{\infty}}(f)$, so gilt insbesondere
$$ d(f(p),g(q)) < \varepsilon ~~ \wedge ~~ d(f(q),g(q)) < \varepsilon, $$
\pagebreak
also $g(p) \ne g(q)$, denn andernfalls gälte 
$$ 2 \varepsilon = d(f(p),f(q)) \le d(f(p),g(q)) + d(g(q),g(p)) + d(g(p),f(q)) < 2 \varepsilon. $$
Folglich ist $g$ nicht konstant, d.h.\ $g \in X^M \setminus \iota(X)$.
Damit ist gezeigt, daß gilt
$$ U_{\varepsilon}^{d_{\infty}}(f) \subset X^M \setminus \iota(X). $$
Wegen der Beliebigkeit von $f$ ist $X^M \setminus \iota(X)$ offen.

Zu (ii): ,,$ \Rightarrow$`` Ist $(X^M, d_{\infty})$ vollständig, so ist die nach (i) abgeschlossene Teilmenge $\iota(X)$ von $X^M$ nach \ref{FA.1.11} (i) vollständig.
Sei nun $(x_n)_{n \in \N}$ eine Cauchyfolge in $(X,d)$.
Da $\iota$ nach (i) isometrisch ist, ist auch $(\underline{x_n})_{n \in \N}$ eine Cauchyfolge in $(\iota(X),d_{\infty}|_{\iota(X),\iota(X)})$, die somit gegen ein gewisses $\underline{x} \in \iota(X)$ mit $x \in X$ in $(\iota(X),d_{\infty}|_{\iota(X),\iota(X)})$ konvergiert.
Erneut, da $\iota$ isometrisch ist, konvergiert dann auch $(x_n)_{n \in \N}$ gegen $x$ in $(X,d)$.

,,$\Leftarrow$`` Sei $(f_n)_{n \in \N}$ eine Cauchyfolge in $(X^M, d_{\infty})$.
Wir definieren eine Abbildung $f \: M \to X$ wie folgt:
Zu $p \in M$ ist offenbar auch $(f_n(p))_{n \in \N}$ eine Cauchyfolge in $(X,d)$, die nach Voraussetzung konvergiert, und wir setzen
$$ f(p) := \lim_{n \to \infty} f_n(p), $$
d.h.\ $(f_n)_{n \in \N}$ konvergiert im üblichen Sinne punktweise gegen $f$.
 
Sei $\varepsilon \in \R_+$.
Da $(f_n)_{n \in \N}$ eine Cauchyfolge in $(X^M, d_{\infty})$ ist, gilt 
\begin{equation} \label{FA.1.Isom.1.1}
\exists_{n_0 \in \N} \, \forall_{n,m \in \N, \, n,m \ge n_0} \, d_{\infty}(f_n, f_m) < \frac{\varepsilon}{3}.
\end{equation}
Nach Konstruktion von $f$ gilt des weiteren für jedes $p \in M$
\begin{equation} \label{FA.1.Isom.1.2}
\exists_{n(p) \in \N, n(p) \ge n_0} \, \forall_{n \in \N, \, n \ge n(p)} \, d(f_n(p), f(p)) < \frac{\varepsilon}{3}.
\end{equation}
Daher folgt aus (\ref{FA.1.Isom.1.1}), (\ref{FA.1.Isom.1.2}) für jedes $p \in M$ und alle $n \in \N$ mit $n \ge n_0$
$$ d(f_n(p), f(p)) \le d(f_n(p), f_{n(p)}(p)) +  d(f_{n(p)}(p), f(p)) < \frac{2 \varepsilon}{3}. $$
Bildet man nun das Supremum über alle $p \in M$, so folgt für alle $n \in \N$ mit $n \ge n_0$
$$ d_{\infty}(f_n, f) \le \frac{2 \varepsilon}{3} < \varepsilon, $$
d.h.\ $(f_n)_{n \in \N}$ konvergiert in $(X^M,d_{\infty})$ gegen $f$. \q

\begin{Def}[Stetigkeit, Homöomorphismen und homöomorphe topologische Räume] \label{FA.1.St}
Seien $X,Y$ topologische Räume und $f \: X \to Y$ eine Abbildung.
\begin{itemize}
\item[(i)] Sei $x \in X$.

$f$ heißt \emph{stetig in $x$}\index{Stetigkeit} $: \Longleftrightarrow \forall_{U \in \U(f(x),Y)} \exists_{V \in \U(x,X)} \, f(V) \subset U.$

Es gilt:
\begin{eqnarray} 
f \mbox{ stetig in } x & \Longrightarrow & \mbox{Für jede Folge } \left( x_n \right)_{n \in \N}, \mbox{ die in } X \mbox{ gegen } x \mbox{ kon-} \nonumber \\
&& \mbox{vergiert, konvergiert die Folge } \left( f(x_n) \right)_{n \in \N} \mbox{in } Y ~~~~ \label{FA.1.St.S} \\
&& \mbox{gegen } f(x). \nonumber
\end{eqnarray}

{[} Denn sind $U \in \U(f(x),Y)$ beliebig und $V \in \U(x,X)$ wie oben gewählt, so gilt für fast alle $n \in \N$: $x_n \in V$, also auch $f(x_n) \in f(V) \subset U$. {]}

\begin{Bem*}
In (\ref{FA.1.St.S}) gilt i.a.\ nicht ,,$\Leftrightarrow$`` anstelle von ,,$\Rightarrow$``, betrachte z.B.\ die Abbildung $\id_{[0,1]} \: ([0,1],\mathcal{T}) \to ([0,1],\mathfrak{P}([0,1]))$ im Punkte $1$, wobei $\mathcal{T} := \{U \in \mathfrak{P}([0,1]) \, | \, \# ([0,1] \setminus U) \le \# \N \} \cup \{ \emptyset \}$.
Dagegen gilt stets ,,$\Leftrightarrow$``, falls $X$ und $Y$ sogar fastmetrische Räume sind, s.u.\ \ref{FA.1.14}.
\end{Bem*}

\item[(ii)] $f$ heißt \emph{stetig}\index{Abbildung!stetige} $: \Longleftrightarrow \forall_{x \in X} \, f \mbox{ ist stetig in } x.$

$\boxed{\mathcal{C}(X,Y)}$ bezeichne die Menge aller stetigen Abbildungen $X \to Y$.

\item[(iii)] Eine bijektive stetige Abbildung $X \to Y$ mit stetiger Umkehrabbildung heißt ein \emph{Homöomorhismus}\index{Homöomorphismus}.
Existiert eine solche Abbildung, so heißen $X$ und $Y$ zueinander \emph{homöomorph}.
\end{itemize}
\end{Def}

\begin{Satz} \label{FA.1.12} \index{Stetigkeit}
Es seien $X, Y$ topologische Räume und $f \: X \to Y$ eine Abbildung.
Dann gilt:
\begin{eqnarray}
f \mbox{ ist stetig } & \Longleftrightarrow & \forall_{U \in {\rm Top}(Y)} \, \overline{f}^1 (U) \in {\rm Top}(X), \label{FA.1.12.1} \\
f \mbox{ ist stetig } & \Longleftrightarrow & \forall_{A \subset Y \text abgeschlossen} \, \overline{f}^1 (A) \mbox{ abgeschlossen in $X$.} \label{FA.1.12.2}
\end{eqnarray}
\end{Satz} 

\textit{Beweis.} Zu (\ref{FA.1.12.1}): ,,$\Rightarrow$`` Sei $U$ offen in $Y$. 
Wegen der Stetigkeit von $f$ existiert zu jedem $x \in \overline{f}^1(U)$ -- d.h.\ $U \in \U(f(x),Y)$ -- ein $V_x \in \U(x,X)$ mit $f(V_x) \subset U$, d.h.\ $V_x \subset \overline{f}^1(U)$. 
Dann folgt
\begin{equation} \label{FA.1.12.3}
\overline{f}^1 (U) = \bigcup_{x \in \overline{f}^1(U)} V_x,
\end{equation}
und die rechte Seite von (\ref{FA.1.12.3}) ist offen nach (T2).

{[} Zu (\ref{FA.1.12.3}): ,,$\subset$`` gilt wegen $\forall_{x \in \overline{f}^1(U)} \, x \in V_x$, und wegen $\forall_{x \in \overline{f}^1(U)} \, V_x \subset \overline{f}^1(U)$ gilt auch ,,$\supset$``. {]}

,,$\Leftarrow$`` Sei $p \in M$ und $U \in \U(f(p),N)$.

Dann gilt $p \in \overline{f}^1(U)$, und nach Voraussetzung ist $\overline{f}^1(U)$ offen in $M$, also $V:=\overline{f}^1(U) \in \U(p,M)$ und $f(V) = f(\overline{f}^1(U)) \subset U$.

Zu (\ref{FA.1.12.2}): Für jede Teilmenge $A$ von $Y$ gilt 
$$ A \text{ abgeschlossen } \Longleftrightarrow  Y \setminus A \in {\rm Top}(Y) ~~ \mbox{ und } ~~ \overline{f}^1(Y \setminus A) = X \setminus \overline{f}^1(A).$$

Daher folgt (\ref{FA.1.12.2}) aus (\ref{FA.1.12.1}). \q

\begin{Satz} \label{FA.1.13}
Es seien $X, Y, Z$ topologische Räume sowie $f \: Y \to Z$ und $g \: X \to Y$ zwei Abbildungen. Ferner sei $x \in X$.
Dann gilt:
\begin{itemize}
\item[(i)] $g$ stetig in $x$ und $f$ stetig in $g(x)$ $\Longrightarrow$ $f \circ g$ stetig in $x$.
\item[(ii)] $g$ stetig und $f$ stetig $\Longrightarrow$ $f \circ g$ stetig.
\end{itemize}
\end{Satz}

\textit{Beweis.} Zu (i): Zu $U_3 \in \U( (f \circ g) (x), Z)$ existiert zunächst wegen der Stetigkeit von $f$ in $g(x)$ ein $U_2 \in \U( g(x), Y)$ mit $f(U_2) \subset U_3$ und sodann wegen der Stetigkeit von $g$ in $x$ ein $U_1 \in \U(p,M_1)$ mit $g(U_1) \subset U_2$, also folgt
$$ (f \circ g) (U_1) = f(g(U_1)) \subset f(U_2) \subset U_3.$$

Aus der Beliebigkeit von $U_3$ folgt die Behauptung von (i).

(ii) folgt trivial aus (i). \q

\begin{Satz} \label{FA.1.TT.3}
Seien $X,Y$ topologische Räume und $f \: X \to Y$ eine Abbildung.
Dann folgt:
\begin{itemize}
\item[(i)] Für jede Teilmenge $\widetilde{X}$ von $X$ und alle $x \in \widetilde{X}$ gilt:
\newline
$f \: X \to Y$ stetig in $x$ $\Longrightarrow$ $f|_{\widetilde{X}} \: \widetilde{X} \to Y$ stetig in $x$.
\item[(ii)] Für jede Teilmenge $\widetilde{Y}$ von $Y$ mit $f(X) \subset \widetilde{Y}$ und alle $x \in X$ gilt:
\newline
$f \: X \to Y$ stetig in $x$ $\Longleftrightarrow$ $f \: X \to \widetilde{Y}$ stetig in $x$.
\end{itemize}
\end{Satz}

\textit{Beweis als Übung.}

\begin{Bem*}
Die Richtung ,,$\Leftarrow$`` in (i) ist i.a.\ falsch, wie das Beispiel $\widetilde{X} = \emptyset$ zeigt.
\end{Bem*}

\begin{Satz}[Stetige Abbildungen zwischen fastmetrischen Räumen] \index{Stetigkeit} \index{Abbildung!stetige} \label{FA.1.14}
Seien $X$, $Y$ fastmetrische Räume, wobei wir beide Fastmetriken mit $d$ bezeichnen, und $f \: X \to Y$ eine Abbildung sowie $x \in X$.
Dann sind die folgenden drei Aussagen paarweise äquivalent:
\begin{itemize}
\item[(i)] $f$ ist stetig in $x$.
\item[(ii)] $\forall_{\varepsilon \in \R_+} \exists_{\delta \in \R_+} \underbrace{\forall_{\tilde{x} \in X} \, \big( d(x,\tilde{x}) < \delta \Longleftrightarrow d(f(x),f(\tilde{x})) < \varepsilon \big)}_{\mbox{$\Longleftrightarrow f(U_{\delta}(x)) \subset U_{\varepsilon}(f(x))$}}$.
\item[(iii)] Für jede Folge $(x_n)_{n \in \N}$ in $X$ mit $\lim_{n \to \infty} x_n = x$ gilt $\lim_{n \to \infty} f(x_n) = f(x)$.
\end{itemize}
\end{Satz}

\textit{Beweis.} ,,(i) $\Rightarrow$ (iii)`` gilt nach (\ref{FA.1.St.S}). 

Zu ,,(iii) $\Rightarrow$ (ii)``: Angenommen, (ii) ist falsch.
Dann existieren $\varepsilon \in \R_+$ und eine Folge $(x_n)_{n \in \N}$ in $X$ mit 
$$ \forall_{n \in \N} \, \left( d(x_n,x) < \frac{1}{n+1} \, \wedge \, d(f(x_n),f(x)) \ge \varepsilon \right), $$
und dies widerspricht (iii). 

Zu ,,(ii) $\Rightarrow$ (i)``: Sei $U \in  \U(f(x),Y)$.
Nach Definition der Topologie von $Y$ existiert $\varepsilon \in \R_+$ mit $U_{\varepsilon}(f(x)) \subset U$.
Wähle zu $\varepsilon$ ein $\delta \in \R_+$ gemäß (ii).
Dann gilt $U_{\delta}(x) \in \U(x,X)$ und $$\forall_{\tilde{x} \in U_{\delta}(x)} d(\tilde{x},x) < \delta, \mbox{ also nach (ii): } d(f(\tilde{x}),f(x)) < \varepsilon,$$ d.h.\ $\forall_{\tilde{x} \in U_{\delta}(x)} f(\tilde{x}) \in U_{\varepsilon}(f(x)) \subset U$, also $f(U_{\delta}(x)) \subset U$.
Damit ist (i) gezeigt. \q

\begin{Def*}
Die Eigenschaft (iii) nennt man \emph{Folgenstetigkeit von $f$ in $x$}\index{Folgen!-stetigkeit}\index{Stetigkeit!Folgen-}. 
\end{Def*}

\begin{Bsp*}
Isometrische Abbildungen zwischen fastmetrischen Räumen sind nach ,,(ii) $\Rightarrow$ (i)`` des letzten Satzes stetig.
Folglich sind Isometrien zwischen solchen Homöomorphismen.
\end{Bsp*}

\begin{Satz} \label{FA.1.Isom.2}
Seien $M$ ein nicht-leerer topologischer Raum und $(X,d)$ ein fastmetrischer Raum.

Dann ist $\mathcal{C}(M,X)$ abgeschlossen in $(X^M, d_{\infty})$.
\end{Satz}

\textit{Beweis.} Sei $(f_n)_{n \in \N}$ eine Folge in $\mathcal{C}(M,X)$, die in $(X^M, d_{\infty})$ gegen $f \in X^M$ konvergiert, d.h.\ 
\begin{equation} \label{FA.1.Isom.2.1}
\forall_{\varepsilon \in \R_+} \exists_{n_0 \in \N} \forall_{n \in \N, \, n \ge n_0} \, d_{\infty}(f_n,f) < \frac{\varepsilon}{3}.
\end{equation}

Nach Satz \ref{FA.1.9} genügt es zu zeigen, daß $f$ stetig ist:
Für alle $n \in \N$ ist $f_n$ stetig, also gilt
\begin{equation} \label{FA.1.Isom.2.2}
\forall_{p \in M} \forall_{\varepsilon \in \R_+} \exists_{\delta_p \in \R_+} \forall_{q \in U_{\delta_p}(p)} \, d(f_n(p),f_n(q)) < \frac{\varepsilon}{3}.
\end{equation}
Sind nun $p \in M$ beliebig, $\varepsilon \in \R_+$ und $n_0, \delta$ gemäß (\ref{FA.1.Isom.2.1}), (\ref{FA.1.Isom.2.2}) gewählt, so folgt für jedes $q \in U_{\delta}(p)$
$$ d(f(p),f(q)) \le d(f(p),f_{n_0}(p)) + d(f_{n_0}(p),f_{n_0}(q)) + d(f_{n_0}(q),f(q)) < \varepsilon. $$
\q
\A
Aus \ref{FA.1.Isom.1} (ii), dem letzten Satz sowie \ref{FA.1.11} (i) ergibt sich unmittelbar das folgende Korollar.

\begin{Kor} \label{FA.1.15}
Seien $M$ ein nicht-leerer topologischer und $(X,d)$ ein vollständiger fastmetrischer Raum.

Dann ist $(\mathcal{C}(M,X),d_{\infty}|_{\mathcal{C}(M,X) \times \mathcal{C}(M,X)})$ vollständig. \q
\end{Kor}

\begin{Satz}[Fixpunktsatz von \textsc{Banach}] \index{Satz!Fixpunkt- von \textsc{Banach}} \label{FA.1.BFix}
Seien $X$ ein vollständiger metrischer Raum, $A \subset X$ eine nicht-leere abgeschlosssene Teilmenge und $f \: A \to X$ eine Abbildung mit folgenden Eigenschaften:

$f$ ist \emph{kontrahierend}, d.h.\ per definitionem 
\begin{equation} \label{FA.1.BFix.1}
\exists_{C \in {]0,1[}} \forall_{x,y \in A} \, d \left( f(x) , f(y) \right) \le C \, d(x,y),
\end{equation}
also insbes.\ stetig, und eine \emph{Selbstabbildung}, d.h.\ per definitionem
\begin{equation} \label{I.1.BFix.2}
f(A) \subset A.
\end{equation}

Dann gilt:
\begin{itemize}
\item[(i)] Es existiert genau ein $x_* \in A$ mit $f(x_*) = x_*$.
\item[(ii)] Sind $x_0 \in A$ beliebig gewählt und die Folge $(x_n)_{n \in \N}$ rekursiv definiert durch $\forall_{n \in \N} \, x_{n+1} := f(x_n)$, so gilt $\lim_{n \to \infty} x_n = x_*$ und darüber hinaus
\begin{equation} \label{FA.1.BFix.3}
\forall_{n \in \N} \, d(x_n,x_*) \le \frac{C^n}{1-C} \, d(x_0,x_1).
\end{equation}
\end{itemize}
\end{Satz}

\textit{Beweis.} Wir zeigen zunächst:
\begin{equation} \label{FA.1.BFix.4}
(x_n)_{n \in \N} \mbox{ ist Cauchyfolge.}
\end{equation}

{[} Es gilt für alle $n \in \N_+$
\begin{eqnarray*}
d(x_n,x_{n+1}) & = & d \left( f(x_{n-1}), f(x_n) \right) \stackrel{(\ref{FA.1.BFix.1})}{\le} C \, d \left( x_{n-1}, x_n \right) \\
& \stackrel{(\ref{FA.1.BFix.1})}{\le} & C^2 \, d \left( x_{n-2}, x_{n-1} \right) \le \ldots \le C^n \, d(x_0,x_1),
\end{eqnarray*}
\pagebreak
also auch für $m \in \N_+$
\begin{equation} \label{FA.1.BFix.5}
d(x_n,x_{n+m}) \le \sum_{i=1}^{m-1} \underbrace{d(x_{n+i},x_{n+i+1})}_{\le C^{n+i} \, d(x_0,x_1)} \le \underbrace{\frac{C^n}{1-C} \, d(x_0,x_1)}_{\stackrel{n \to \infty}{\longrightarrow} 0},
\end{equation}
und hieraus folgt (\ref{FA.1.BFix.4}). {]}

Da $X$ vollständig ist, folgt aus (\ref{FA.1.BFix.4}) die Existenz von $x_* \in X$ mit 
$$ \lim_{n \to \infty} x_n = x_*,  $$
und weil $(x_n)_{n \in \N}$ eine Folge in der abgeschlossenen Menge $A$ ist, gilt $x_* \in A$ nach Satz \ref{FA.1.KA}.

Die Gültigkeit von (\ref{FA.1.BFix.3}) folgt aus (\ref{FA.1.BFix.5}) durch Bildung des Grenzwertes für $m \to \infty$.

Zum Nachweis des Satzes bleibt daher zu zeigen, daß $x_*$ die in (i) genannte Eigenschaft hat:

Aus (\ref{FA.1.BFix.1}) folgt, daß $f \: A \to X$ stetig ist, also gilt
$$ x_* = \lim_{n \to \infty} x_{n+1} = \lim_{n \to \infty} f( x_n ) = f( x_* ). $$

Sei $x \in A$ beliebig mit $f(x) = x$. Dann ergibt (\ref{FA.1.BFix.1})
$$ d( x_* , x ) = d \left( f(x_*) , f(x) \right) \le \underbrace{C}_{\in {]0,1[}} d( x_* , x ), $$
und dies ist nur im Falle $d( x_* , x ) = 0$, d.h.\ $x = x_*$, möglich. \q

\begin{Def}[(Cauchy-)Äquivalenz von Fastmetriken] \label{FA.1.16}
Seien $X$ eine Menge und $d$, $\tilde{d}$ Fastmetriken auf $X$.
\begin{itemize}
\item[(i)] $d$ und $\tilde{d}$ heißen \emph{äquivalent}\index{Fastmetrik!Äquivalenz von --en} genau dann, wenn der in \ref{FA.1.2} beschriebe Prozeß für $d$ und $\tilde{d}$ dieselbe Topologie liefert.
\item[(ii)] $d$ und $\tilde{d}$ heißen \emph{Cauchy-äquivalent}\index{Cauchy!-äquivalent}\index{Fastmetrik!Cauchy-Äquivalenz von --en} genau dann, wenn jede Cauchyfolge in $(X,d)$ eine Cauchyfolge in $(X,\tilde{d})$ ist und umgekehrt.
\end{itemize}
\end{Def}

\begin{Satz} \label{FA.1.17}
Seien $X$ eine Menge und $d$, $\tilde{d}$ Fastmetriken auf $X$.

Sind $d$ und $\tilde{d}$ Cauchy-äquivalent, so sind $d$ und $\tilde{d}$ äquivalent.
\end{Satz}

\textit{Beweis.} Bezeichnen $\mathcal{T}, \widetilde{\mathcal{T}}$ die durch $d, \tilde{d}$ induzierten Topologien.
Ohne Beschränkung der Allgemeinheit genügt es zu zeigen, daß gilt $\mathcal{T} \subset \widetilde{\mathcal{T}}$.
Hierfür wiederum ist nach (\ref{FA.1.12.1}) zu zeigen, daß
\begin{equation} \label{FA.1.17.S}
\id_X \: (X,\tilde{d}) \longrightarrow (X,d)
\end{equation}
stetig ist.

Sei $x \in X$.
Wir zeigen, daß (\ref{FA.1.17.S}) in $x$ folgenstetig ist.
Sei daher $(x_n)_{n \in \N}$ eine Folge in $X$, die bzgl.\ $\tilde{d}$ gegen $x$ konvergiert.
Dann konvergiert auch $(y_n)_{n \in \N}$, definiert durch
$$ \forall_{n \in \N} \, y_{2n} := x_n \, \wedge \, y_{2n+1} := x, $$
bzgl.\ $\tilde{d}$ gegen $x$, ist also eine Cauchyfolge bzgl. $\tilde{d}$ und somit nach Voraussetzung auch eine Cauchyfolge bzgl. $d$, die die gegen $x$ konvergente Teilfolge $(y_{2n+1})_{n \in \N}$ besitzt.
Hieraus folgt offenbar, daß $(y_n)_{n \in \N}$ und damit auch die Teilfolge $(x_n)_{n \in \N}$ bzgl.\ $d$ gegen $x$ konvergiert. \q

\begin{Bem*}
Durch
$$ \forall_{x,y \in \R} \, d(x,y) := d_{\R}(\arctan(x), \arctan(y)) = | \arctan(x) - \arctan(y) | $$
wird eine Metrik $d$ auf $\R$ definiert, die zu $d_{\R}$ äquivalent aber nicht Cauchy-äquivalent ist. 
Denn $\arctan \: (\R, d) \to \left( {]} - \frac{\pi}{2}, \frac{\pi}{2} {[}, d_{\R} \right)$ ist eine Isometrie und somit ebenso wie $\tan \:  \left( \left] - \frac{\pi}{2}, \frac{\pi}{2} \right[, d_{\R} \right) \to (\R,d_{\R})$ ein Homöoomorphismus, also ist auch $\id_{\R} \: (\R, d) \to (\R,d_{\R})$ ein Homöomorphismus
$$ (\R,d) \stackrel{\arctan}{\longrightarrow} \left( \left] - \frac{\pi}{2}, \frac{\pi}{2} \right[, d_{\R} \right) \stackrel{\tan}{\longrightarrow} (\R,d_{\R}), $$
d.h.\ $d$ und $d_{\R}$ sind äquivalent.
Der zu $(\R, d)$ isometrische Raum $\left( \left] - \frac{\pi}{2}, \frac{\pi}{2} \right[, d_{\R} \right)$ ist aber im Gegensatz zu $(\R,d_{\R})$ nicht vollständig, d.h.\ $d$ und $d_{\R}$ sind nicht Cauchy-äquivalent.
\end{Bem*}

\begin{Satz} \label{FA.1.18}
Seien $(X,d)$ ein fastmetrischer Raum und $\varphi \: [0, \infty] \to [0, 1]$ definiert durch
$$ \forall_{t \in [0, \infty[} \, \varphi(t) := \frac{t}{1+t} ~~ \wedge ~~ \varphi(\infty) := \lim_{t \to \infty} \frac{t}{1+t} = 1. $$

Dann ist $\varphi \circ d$ eine zu $d$ (Cauchy-)äquivalente Metrik auf $X$.
\end{Satz}

\textit{Beweis.} 1.) $\varphi \circ d$ ist eine Metrik:
Daß $\varphi \circ d$ reellwertig ist, ist ebenso wie (M1) und (M2) trivial.
(M3) folgt offenbar daraus, daß für alle $t_1, t_2 \in [0, \infty[$
$$ \varphi(t_1 + t_2) = \frac{t_1 + t_2}{1 + t_1 + t_2} = \frac{t_1}{1 + t_1 + t_2} + \frac{t_2}{1 + t_1 + t_2} \le \frac{t_1}{1 + t_1} + \frac{t_2}{1 + t_2} = \varphi(t_1) + \varphi(t_2) $$
und für alle $t  \in [0, \infty]$
$$ \varphi(t + \infty) = \varphi(\infty) = 1 \le \varphi(t) + 1 = \varphi(t) + \varphi(\infty) $$
gilt.

2.) Sei $(x_n)_{n \in \N}$ eine Cauchyfolge bzgl.\ $d$, d.h.\
\begin{equation} \label{FA.1.18.1}
\forall_{\varepsilon \in \R_+} \exists_{n_0 \in \N} \forall_{n,m \in \N, \, n,m \ge n_0} \, d(x_n,x_m) < \varepsilon.
\end{equation}
Wegen $\forall_{t \in [0,\infty[} \, \varphi(t) = \frac{t}{1+t} \le t$ folgt hieraus
\begin{equation} \label{FA.1.18.2}
\forall_{\varepsilon \in \R_+} \exists_{n_0 \in \N} \forall_{n,m \in \N, \, n,m \ge n_0} \, (\varphi \circ d)(x_n,x_m) < \varepsilon,
\end{equation}
d.h.\ $(x_n)_{n \in \N}$ ist eine Cauchyfolge bzgl.\ $\varphi \circ d$.

3.) Sei $(x_n)_{n \in \N}$ eine Cauchyfolge bzgl.\ $\varphi \circ d$, d.h.\ es gilt (\ref{FA.1.18.2}).
Sei $\varepsilon \in \R_+$. 
Dann existiert gemäß (\ref{FA.1.18.2}) eine Zahl $n_0 \in \N$ derart, daß für alle $n,m \in \N$ mit $n,m \ge n_0$ gilt
\begin{equation} \label{FA.1.18.3}
\forall_{\varepsilon \in \R_+} \exists_{n_0 \in \N} \forall_{n,m \in \N, \, n,m \ge n_0} \, (\varphi \circ d)(x_n,x_m) < \varphi(\varepsilon) \in {]}0,1{[}.
\end{equation}
$\varphi|_{[0, \infty[} \: [0, \infty[ \to [0,1[$ ist differenzierbar mit $\forall_{t \in [0, \infty[} \, \varphi'(t) = \frac{1}{(1+t)^2} > 0$, also ist $\varphi \: [0,\infty] \to [0,1]$ streng monoton wachsend und offenbar bijektiv. 
Daher ist auch $\varphi^{-1} \: [0, 1] \to [0, \infty]$ streng monoton wachsend, und aus (\ref{FA.1.18.3}) folgt die Gültigkeit von (\ref{FA.1.18.1}), d.h.\ $(x_n)_{n \in \N}$ ist eine Cauchyfolge bzgl.\ $d$. \q 

\subsection*{Vervollständigung fastmetrischer Räume} \addcontentsline{toc}{subsection}{Vervollständigung fastmetrischer Räume}

\begin{Def}[Gleichmäßige Stetigkeit] \label{FA.1.19}
Seien $X,Y$ fastmetrische Räume, wobei beide Fastmetriken mit $d$ bezeichnet seien, und $f \: X \to Y$ eine Abbildung.

$f$ heißt \emph{gleichmäßig stetig}\index{Stetigkeit!gleichmäßige}\index{Abbildung!stetige!gleichmäßig} genau dann, wenn gilt
$$ \forall_{\varepsilon \in \R_+} \exists_{\delta \in \R_+} \forall_{x, \tilde{x} \in X} \, \left( d(x, \tilde{x}) < \delta \Longrightarrow d(f(x), f(\tilde{x})) < \varepsilon \right). $$

\begin{Bsp*}
Eine isometrische Abbildung zwischen fastmetrischen Räumen ist gleichmäßig stetig.
\end{Bsp*}
\end{Def}

\begin{Lemma} \label{FA.1.BspSt} 
Sei $(X,d)$ ein metrischer Raum.

Dann ist $d \: ( X \times X, \tilde{d} ) \to (\R,d_{\R})$,
wobei $\tilde{d}$ die Produktmetrik auf $X \times X$ und $d_{\R}$ die übliche Metrik auf $\R$ bezeichne, gleichmäßig stetig.
\end{Lemma}

\textit{Beweis.} Nach \ref{FA.1.L} gilt für alle $(x,\tilde{x}), (y,\tilde{y}) \in X \times X$
$$ | d(x, \tilde{x}) - d(y, \tilde{y}) | \le d(x,y) + d(\tilde{x}, \tilde{y}) = \tilde{d}( (x, \tilde{x}), (y, \tilde{y}) ), $$
und hieraus ergibt sich die Behauptung. \q

\begin{Satz}[Fortsetzungssatz gleichmäßig stetiger Abbbildungen] \label{FA.1.20} \index{Satz!Fortsetzungs-!gleichmäßig stetiger Abbildung\-en} $\,$

\noindent \textbf{Vor.:} Seien $X,Y$ fastmetrische Räume, wobei beide Fastmetriken mit $d$ bezeichnet seien, $A \subset X$ und $f \: A \to Y$ eine Abbildung mit
\begin{equation} \label{FA.1.20.S}
\forall_{x \in \overline{A}} \exists_{U \in \U(x,X)} \, f|_{A \cap U} \mbox{ ist gleichmäßig stetig.}
\end{equation}
Ferner sei $Y$ vollständig.

\noindent \textbf{Beh.:} Es existiert genau eine stetige Abbildung $\overline{f} \: \overline{A} \to Y$ mit $\overline{f}|_A = f$.
\end{Satz}

\textit{Beweis.} 1.) Wenn eine Abbildung $\overline{f}$ wie in der Behauptung existiert, so muß aus Stetigkeitsgründen
\begin{equation} \label{FA.1.20.D}
\forall_{x \in \overline{A}} \forall_{(x_n)_{n \in \N} \in A^{\N}, \, \lim_{n \to \infty} x_n = x}  \, \overline{f}(x) = \lim_{n \to \infty} f( \underbrace{x_n}_{\in A} )
\end{equation}
gelten, d.h.\ $\overline{f}$ ist durch $f$ eindeutig bestimmt.
Beachte, daß nach \ref{FA.1.9} zu jedem $x \in \overline{A}$ eine gegen $x$ konvergente Folge $(x_n)_{n \in \N}$ in $A$ existiert.
Wir definieren $\overline{f}$ durch (\ref{FA.1.20.D}) und haben zu zeigen, daß dies wohldefiniert ist.
Es ist nachzuweisen, daß zu $(x_n)_{n \in \N}$ wie in (\ref{FA.1.20.D}) die Folge $(f(x_n))_{n \in \N}$ in $Y$ konvergiert und der Grenzwert unabhängig von der speziellen Wahl der Folge $(x_n)_{n \in \N}$ ist.

a) Seien $x \in \overline{A}$ und $(x_n)_{n \in \N}$ eine Folge in $A$, die gegen $x$ konvergiert, insbesondere ist $(x_n)_{n \in \N}$ eine Cauchyfolge in $A$, d.h.\
\begin{equation} \label{FA.1.20.1}
\forall_{\delta \in \R_+} \exists_{n_0 \in \N} \forall_{n,m \in \N, \, n,m \ge n_0} \, d(x_n,x_m) < \delta.
\end{equation}
Um zu zeigen, daß $(f(x_n))_{n \in \N}$ in $Y$ konvergiert, genügt es wegen der Vollständigkeit von $Y$ zu zeigen
\begin{equation} \label{FA.1.20.2}
(f(x_n))_{n \in \N} \mbox{ ist Cauchyfolge.}
\end{equation}

{[} Zu (\ref{FA.1.20.2}): Sei $\varepsilon \in \R_+$.
Wegen (\ref{FA.1.20.S}) können wir ein $\delta \in \R_+$ wählen derart, daß gilt
$$ \forall_{a,b \in A} \, \left( d(a,b) < \delta \Longrightarrow d(f(a),f(b)) < \varepsilon \right), $$
also folgt aus (\ref{FA.1.20.1}) die Existenz einer Zahl $n_0 \in \N$ mit
$$ \forall_{n,m \in \N} \, \left( n,m \ge n_0 \Longrightarrow d(f(x_n),f(x_m)) < \varepsilon \right) $$
und (\ref{FA.1.20.2}) ist gezeigt. {]}

b) Seien nun $x \in \overline{A}$ und $(x_n)_{n \in \N}, (\tilde{x}_n)_{n \in \N}$ Folgen in $A$, die beide gegen $x$ konvergieren.
Dann wird durch 
$$ \forall_{n \in \N} \, \tilde{\tilde{x}}_{2n} := x_n \, \wedge \, \tilde{\tilde{x}}_{2n+1} := \tilde{x}_n $$
eine weitere Folge in $A$ definiert, die $x$ als Grenzwert besitzt.
Nach a) konvergiert dann auch $(f(\tilde{\tilde{x}}_n))_{n \in \N}$ in $Y$ und somit jede Teilfolge von $(f(\tilde{\tilde{x}}_n))_{n \in \N}$ gegen denselben Grenzwert.
Insbesondere ergibt sich $\lim_{n \to \infty} f(x_n) = \lim_{n \to \infty} f(\tilde{x}_n)$.

2.) Für das in 1.) definierte $\overline{f}$ gilt $\overline{f}|_A = f$, denn für jedes $x \in A$ besitzt die konstante Folge vom Wert $x$ die konstante Folge vom Wert $f(x)$ als Bildfolge.

3.) Zu zeigen bleibt, daß $\overline{f}$ stetig ist. 
Seien $x \in \overline{A}$ und $\varepsilon \in \R_+$.
Es existiert eine Zahl $\delta \in \R_+$ derart, daß
\begin{equation} \label{FA.1.20.3}
\forall_{a \in A} \left( d(a,x) < 2 \delta \Longrightarrow d(f(a),\overline{f}(x)) < \frac{\varepsilon}{2} \right),
\end{equation}
denn andernfalls gäbe es eine Folge in $A$, die gegen $x$ konvergiert, deren Bilderfolge nicht gegen $f(x)$ konvergiert, im Widerspruch zu 1.).

Sei $\tilde{x} \in \overline{A}$ mit $d(x,\tilde{x}) < \delta$.

Da $\tilde{x}$ ein Berührungspunkt von $A$ ist, existiert eine Folge $(\tilde{x}_n)_{n \in \N}$ in $A$ mit $\lim_{n \to \infty} \tilde{x}_n = \tilde{x}$, d.h.\ insbes.\
\begin{equation} \label{FA.1.20.4}
\exists_{n_1 \in \N} \forall_{n \in \N, \, n \ge n_1} \, d(\tilde{x}_n,\tilde{x}) < \delta, 
\end{equation}
und wegen 1.) (angewandt auf $\tilde{x}$ anstelle von $x$) gilt $\lim_{n \to \infty} f(\tilde{x}_n) = \overline{f}(\tilde{x})$, d.h.\ insbes.\
\begin{equation} \label{FA.1.20.5}
\exists_{n_2 \in \N} \forall_{n \in \N, \, n \ge n_2} \, d(f(\tilde{x}_n),\overline{f}(\tilde{x})) < \frac{\varepsilon}{2}. 
\end{equation}

Nun gilt für $n_0 := \max \{n_1, n_2\}$ wegen der Voraussetzung an $\tilde{x}$ und (\ref{FA.1.20.4})
$$ d(x,\tilde{x}_{n_0}) \le d(x,\tilde{x}) + d(\tilde{x},\tilde{x}_{n_0}) < 2 \delta, $$
d.h.\ nach (\ref{FA.1.20.3})
\begin{equation} \label{FA.1.20.6}
d(f(\tilde{x}_{n_0}),\overline{f}(x)) < \frac{\varepsilon}{2},
\end{equation}
und es ergibt sich aus (\ref{FA.1.20.6}), (\ref{FA.1.20.5})
$$d(\overline{f}(x), \overline{f}(\tilde{x})) \le d(\overline{f}(x), f(\tilde{x}_{n_0})) + d(f(\tilde{x}_{n_0}), \overline{f}(\tilde{x})) < \varepsilon. $$

Mit. 1.), 2.) und 3.) ist der Satz vollständig bewiesen. \q

\begin{Bem*}
Die Bedingung (\ref{FA.1.20.S}) ist natürlich erfüllt, wenn $f$ sogar gleichmäßig stetig ist.
Der letzte Satz gilt dagegen nicht, wenn man (\ref{FA.1.20.S}) durch die schwächere Bedingung der Stetigkeit von $f$ ersetzt.
Z.B.\ läßt sich $\frac{1}{x} \: \R_+ \to \R$ nicht auf $[0,\infty[ \subset \R$ stetig fortsetzen.
\end{Bem*}

\begin{HS}[Vervollständigung (fast-)metrischer Räume] \label{FA.1.21} \index{Vervollständigung!(fast)metrischer Räume} $\,$

\noindent \textbf{Vor.:} Sei $(X,d)$ ein metrischer Raum.

\noindent \textbf{Beh.:} Es existiert ein metrischer Raum $(\widehat{X},\widehat{d})$ mit folgenden Eigenschaften:
\begin{itemize}
\item[(i)] $\widehat{X}$ ist vollständig.
\item[(ii)] $X \subset \widehat{X}$ und $d = \widehat{d}|_{X \times X}$.
\item[(iii)] $\overline{X} = \widehat{X}$, d.h.\ per definitionem $X$ ist \emph{dicht in $\widehat{X}$}\index{Menge!dichte}.
\item[(iv)] $(\widehat{X},\widehat{d})$ mit (i) - (iii) ist bis auf Isometrie eindeutig bestimmt, daher bezeichnet man $(\widehat{X},\widehat{d})$ als ,,die`` \emph{Vervollständigung von $(X,d)$}. 
\end{itemize}
\end{HS}

\begin{Zusatz}
Da ein fastmetrischer Raum $X$ die disjunkte Vereinigung seiner Äquivalenzklassen bzgl.\ der Äquivalenzrelation $\sim$, wobei $\sim$ durch
$$ \forall_{x,y \in X} \, x \sim y :\Longleftrightarrow d(x,y) < \infty $$
definiert sei, welches dann metrische Räume sind, gilt der letzte Hauptsatz auch, wenn man ,,metrisch`` jeweils durch ,,fastmetrisch`` ersetzt.
\end{Zusatz}

\textit{Beweis.} Es sei $\mathfrak{X}$ die Menge aller Cauchyfolgen in $X$.
Auf $\mathfrak{X}$ wird durch
$$ \forall_{\mathfrak{x}=(x_n)_{n \in \N}, \mathfrak{y}=(y_n)_{n \in \N} \in \mathfrak{X}} \, \mathfrak{x} \sim \mathfrak{y} : \Longleftrightarrow (d(x_n,y_n))_{n \in \N} \mbox{ ist Nullfolge} $$
offenbar eine Äquivalenzrelation $\sim$ in $\mathfrak{X}$ definiert.
Wir setzen 
$$ \widehat{X} := \mathfrak{X} / \sim \, = \left\{ [\mathfrak{x}]_{\sim} \, | \, \mathfrak{x} \in \mathfrak{X} \right\} $$
und fassen die offenbar injektive Abbildung
$$ \iota \: X \longrightarrow \widehat{X}, ~~ x \longmapsto [\underline{x}]_{\sim}, $$
wobei $\underline{x}$ für jedes $x \in X$ die konstante Folge vom Wert $x$ bezeichne, als Inklusion auf, d.h.\ $\forall_{x \in X} \, \iota(x) \equiv x$.

Zunächst gilt für alle $\mathfrak{x}=(x_n)_{n \in \N}, \tilde{\mathfrak{x}}=(\tilde{x}_n)_{n \in \N}, \mathfrak{y}=(y_n)_{n \in \N}, \tilde{\mathfrak{y}}=(\tilde{y}_n)_{n \in \N} \in \mathfrak{X}$
\begin{gather}
(d(x_n,y_n))_{n \in \N} \mbox{ konvergiert in $\R$ (bzgl.\ der üblichen Metrik $d_{\R}$),} \label{FA.1.21.1} \\
\mathfrak{x} \sim \tilde{\mathfrak{x}} \, \wedge \, \mathfrak{y} \sim \tilde{\mathfrak{y}} \Longrightarrow \lim_{n \to \infty} d(x_n,y_n) = \lim_{n \to \infty} d(\tilde{x}_n, \tilde{y}_n). \label{FA.1.21.2}
\end{gather}

{[} Zu (\ref{FA.1.21.1}): Da $\R$ vollständig ist, genügt es zu zeigen, daß $(d(x_n,y_n))_{n \in \N}$ eine Cauchyfolge ist.
Sei $\varepsilon \in \R_+$.
Da $\mathfrak{x}, \mathfrak{y}$ Cauchyfolgen sind existiert $n_0 \in \N$ derart, daß für alle $n,m \in \N$ mit $n,m \ge n_0$ gilt
$$ d(x_n,x_m) < \frac{\varepsilon}{2} ~~ \wedge ~~ d(y_n,y_m) <\frac{\varepsilon}{2}, $$
also folgt aus \ref{FA.1.L}
$$ | d(x_n,y_n) - d(x_m,y_m) | \le d(x_n,x_m) + d(y_n,y_m) < \varepsilon. $$

Zu (\ref{FA.1.21.2}): Wegen $\mathfrak{x} \sim \tilde{\mathfrak{x}}$ und $\mathfrak{y} \sim \tilde{\mathfrak{y}}$ gilt
$$ \lim_{n \to \infty} d(x_n,\tilde{x}_n) = 0 ~~ \wedge ~~ \lim_{n \to \infty} d(y_n,\tilde{y}_n) = 0. $$
Erneute Anwendung von \ref{FA.1.L} ergibt für alle $n \in \N$
$$ | d(x_n,y_n) - d(\tilde{x}_n,\tilde{y}_n) | \le d(x_n,\tilde{x}_n) + d(y_n,\tilde{y}_n), $$
also folgt $\lim_{n \to \infty} d(x_n,y_n) = \lim_{n \to \infty} d(\tilde{x}_n, \tilde{y}_n)$. {]}

Wegen (\ref{FA.1.21.1}), (\ref{FA.1.21.2}) wird durch
$$ \forall_{\mathfrak{x}=(x_n)_{n \in \N}, \mathfrak{y}=(y_n)_{n \in \N} \in \mathfrak{X}} \, \widehat{d}([\mathfrak{x}]_{\sim}, [\mathfrak{y}]_{\sim}) := \lim_{n \to \infty} d(x_n,y_n) $$
offenbar eine Metrik auf $\widehat{X}$ definiert.\footnote{Auf $\mathfrak{X}$ wird durch die entsprechende Definition nur eine Halbmetrik definiert. Eine \emph{Halbmetrik}\index{Halb!-metrik}\index{Metrik (siehe auch Fastmetrik)!Halb-} definiert man analog zu einer Metrik, indem man in (D1) nur ,,$\Leftarrow$`` anstelle von ,,$\Leftrightarrow$`` fordert.}
Damit ist (ii) bereits gezeigt.

Zu (i): Sei $(\hat{x}_k)_{k \in \N}$ eine Cauchyfolge in $\widehat{X}$, also gilt
\begin{equation} \label{FA.1.21.3}
\forall_{i \in \N} \exists_{k_i \in \N} \forall_{k,l \in \N, \, k,l \ge k_i} \, \widehat{d}(\hat{x}_k, \hat{x}_l) < \frac{1}{i+1},
\end{equation}
und wir fixieren zu jedem $i \in \N$ ein $k_i$ mit (\ref{FA.1.21.3}) und
\begin{equation} \label{FA.1.21.3.1}
(k_i)_{i \in \N} \mbox{ ist streng monoton wachsende Folge in $\N$.}
\end{equation}
Wir wählen zu $k \in \N$ ein $\mathfrak{x}_k = (x_{k,n})_{n \in \N} \in \mathfrak{X}$ mit $[\mathfrak{x}_k]_{\sim} = \hat{x}_k$.
Dann besagt (\ref{FA.1.21.3})
\begin{equation} \label{FA.1.21.4}
\forall_{i,k,l \in \N, \, k,l \ge k_i} \exists_{n_{i,k,l} \in \N} \forall_{n \in \N, \, n \ge n_{i,k,l}} \, d(x_{k,n}, x_{l,n}) < \frac{1}{i+1}. 
\end{equation}
Wir wählen zu $i,k,l \in \N$ mit $k,l \ge k_i$ ein $n_{i,k,l} \in \N$ gemäß (\ref{FA.1.21.4}) und definieren für alle $i \in \N$ eine Folge $(n_{i,j})_{j \in \N, \, j \ge i}$ durch
$$ \forall_{j \in \N, \, j \ge i} \, n_{i,j} := \max \{ n_{i,k,k_j} \, | \,  k \in \N \, \wedge \, k_i \le k \le k_j \}, $$
beachte, daß aus $j \ge i$ wegen (\ref{FA.1.21.3.1}) folgt $k_j \ge k_i$.
(\ref{FA.1.21.4}) ergibt nun (mit $l = k_j$)
\begin{equation} \label{FA.1.21.5}
\forall_{i,j,k,n \in \N, \, j \ge i, \, k_i \le k \le k_j, \, n \ge n_{i,j}} \, d(x_{k,n}, x_{k_j,n}) < \frac{1}{i+1}. 
\end{equation}

Da $\mathfrak{x}_k = (x_{k,n})_{n \in \N}$ für jedes $k \in \N$ als Element von $\mathfrak{X}$ eine Cauchyfolge in $X$ ist, können wir nach eventueller Vergrößerung der $n_{i,j}$ erreichen, daß gilt
\begin{equation} \label{FA.1.21.6}
\forall_{i,j,k,n \in \N, \, j \ge i, n \ge n_{i,j}}  \, d(x_{k,n}, x_{k,n_{i,j}}) < \frac{1}{i+1}.
\end{equation}
und außerdem nach ggf.\ zusätzlicher Vergrößerung der $n_{j,j}$
\begin{equation} \label{FA.1.21.7}
(n_{j,j})_{j \in \N} \mbox{ ist streng monoton wachsende Folge in $\N$.}
\end{equation}

Wir zeigen:
\begin{equation} \label{FA.1.21.8}
\mathfrak{x}_{\infty} := (x_{k_j,n_{j,j}})_{j \in \N} \in \mathfrak{X}.
\end{equation}

{[} Zu (\ref{FA.1.21.8}): Seien $\varepsilon \in \R_+$ und $j,l \in \N$ mit $j \ge l > \frac{2}{\varepsilon} - 1$.
Dann gilt wegen $n_{j,j} \stackrel{(\ref{FA.1.21.7})}{\ge} n_{l,l}$ nach (\ref{FA.1.21.6}) (wähle dort $i = l, j = l, k = k_j, n = n_{j,j}$)
$$ d(x_{k_j,n_{j,j}}, x_{k_j,n_{l,l}}) < \frac{1}{l+1} $$
und nach (\ref{FA.1.21.5}) (wähle dort $i = l, j = l, k = k_l, n = n_{l,l}$)
$$ d(x_{k_j,n_{l,l}}, x_{k_l,n_{l,l}}) < \frac{1}{l+1}, $$
also auch
$$ d(x_{k_j,n_{j,j}}, x_{k_l,n_{l,l}}) \le d(x_{k_j,n_{j,j}}, x_{k_j,n_{l,l}}) + d(x_{k_j,n_{l,l}}, x_{k_l,n_{l,l}}) < \frac{1}{l+1} + \frac{1}{l+1} < \varepsilon, $$
d.h.\ $\mathfrak{x}_{\infty}$ ist eine Cauchyfolge in $X$. {]}

Zum Beweis von $(i)$ bleibt nachzuweisen, daß $([\mathfrak{x}_k]_{\sim})_{k \in \N}$ in $(\widehat{X}, \widehat{d})$ gegen $[\mathfrak{x}_{\infty}]_{\sim}$ konvergiert, d.h.\ genau
\begin{equation} \label{FA.1.21.9}
\lim_{k \to \infty} \underbrace{\lim_{j \to \infty} d(x_{k,j}, x_{k_j,n_{j,j}})}_{= \widehat{d}([\mathfrak{x}_k]_{\sim}, [\mathfrak{x}_{\infty}]_{\sim})} = 0.
\end{equation}

{[} Zu (\ref{FA.1.21.9}): Seien $\varepsilon \in \R_+$ und $i \in \N$ mit $i > \frac{2}{\varepsilon} - 1$.
Für jedes $k \in \N$ ist $(x_{k,j})_{j \in \N}$ eine Cauchyfolge in $X$, also existiert $j_{i,k} \in \N$ derart, daß gilt
\begin{equation} \label{FA.1.21.10}
\forall_{j,l \in \N, \, j,l \ge j_{i,k}} \, d(x_{k,j}, x_{k,l}) < \frac{1}{i+1}.
\end{equation}
Dann gilt für $j \in \N$ mit $j \ge j_{i,k}$ nach (\ref{FA.1.21.7}) auch $n_{j,j} \ge j_{i,k}$ und somit wegen
(\ref{FA.1.21.10})
$$ d(x_{k,j}, x_{k,n_{j,j}}) < \frac{1}{i+1} $$
und wegen (\ref{FA.1.21.5}), falls zusätzlich $k_i \le k \le k_j$, 
$$ d(x_{k,n_{j,j}}, x_{k_j,n_{j,j}}) < \frac{1}{i+1}, $$
folglich
$$ \forall_{j,k \in \N, \, j \ge j_{i,k}, \, k_i \le k \le k_j} \, d(x_{k,j}, x_{k_j,n_{j,j}}) \le d(x_{k,j}, x_{k,n_{j,j}}) + d(x_{k,n_{j,j}}, x_{k_j,n_{j,j}}) < \frac{2}{i+1}. $$
Hieraus folgt für $j \to \infty$, d.h.\ nach (\ref{FA.1.21.3.1}) auch $k_j \to \infty$,
$$ \forall_{k \in \N, \, k \ge k_i} \, \lim_{j \to \infty} d(x_{k,j}, x_{k_j,n_{j,j}}) \le \frac{2}{i+1} < \varepsilon, $$
und (\ref{FA.1.21.9}) ist bewiesen. {]}

Zu (iii): Seien $\mathfrak{x} = (x_n)_{n \in \N} \in \mathfrak{X}$ und $\varepsilon \in \R_+$ beliebig.
$X \cap U^{\widehat{d}}_{\varepsilon}([\mathfrak{x}]_{\sim}) \ne \emptyset$ ist zu zeigen, d.h.\ die Existenz von $x \in X$ mit
$$ \underbrace{\widehat{d}( [\mathfrak{x}]_{\sim}, [\underline{x}]_{\sim} )}_{= \lim_{n \to \infty} d(x_n,x)} < \varepsilon. $$
Da $(x_n)_{n \in \N}$ eine Cauchyfolge in $X$ ist, existiert $n_0 \in \N$ derart, daß für alle $n,m \in \N$ mit $n,m \ge n_0$ gilt $d(x_n,x_m) < \frac{\varepsilon}{2}$.
Hieraus folgt für $x := x_{n_0}$
$$ \forall_{n \in \N, \, n \ge n_0} \, d(x_n,x) < \frac{\varepsilon}{2}, $$
also auch $\lim_{n \to \infty} d(x_n,x) \le \frac{\varepsilon}{2} < \varepsilon$.

Zu (iv): Sei auch $(\widetilde{X},\widetilde{d})$ ein metrischer Raum, der (i) - (iii) für $(\widetilde{X},\widetilde{d})$ anstelle von $(\widehat{X},\widehat{d})$ erfüllt.
Die Inklusionen $\iota \: (X,d) \to (\widehat{X},\widehat{d})$ und $\tilde{\iota} \: (X,d) \to (\widetilde{X},\widetilde{d})$ sind nach (ii) isometrisch, also Isometrien -- und damit auch Bijektionen -- auf ihr Bild $X = \iota(X) = \tilde{\iota}(X)$.
Dann ist auch
\begin{equation} \label{FA.1.21.11}
f := \tilde{\iota} \circ \left( \iota \: X \to \iota(X) \right)^{-1} \: \underbrace{(X,d)}_{\subset (\widehat{X},\widehat{d})} \longrightarrow (\widetilde{X},\widetilde{d}) \mbox{ eine isometrische Abbildung,}
\end{equation}
d.h.\ insbes.\ eine gleichmäßge stetige Abbildung. 
$(\widetilde{X},\widetilde{d})$ ist vollständig, und es gilt $\overline{X} = \widehat{X}$.
Nach \ref{FA.1.20} läßt sich $f$ daher zu einer stetigen Abbildung
$$ \overline{f} \: (\widehat{X}, \widehat{d}) \longrightarrow (\widetilde{X},\widetilde{d}) $$
fortsetzen.
Wir behaupten:
\begin{equation} \label{FA.1.21.12}
\overline{f} \: (\widehat{X}, \widehat{d}) \longrightarrow (\widetilde{X},\widetilde{d}) \mbox{ ist eine isometrische Abbildung.}
\end{equation}

{[} Zu (\ref{FA.1.21.12}): Seien $\hat{x},\hat{y} \in \widehat{X}$.
Wegen $\widehat{X} = \overline{X}$ existieren Folgen $(x_n)_{n \in \N}, (y_n)_{n \in \N}$ in $X$ mit $\lim_{n \to \infty}^{\widehat{d}} x_n = \hat{x}$ und $\lim_{n \to \infty}^{\widehat{d}} y_n = \hat{y}$.
Aus Stetigkeitsgründen folgt dann $\overline{f}(\hat{x}) = \lim_{n \to \infty}^{\widetilde{d}} \overline{f}(x_n) = \lim_{n \to \infty}^{\widetilde{d}} f(x_n)$, $\overline{f}(\hat{y}) = \lim_{n \to \infty}^{\widetilde{d}} \overline{f}(y_n) = \lim_{n \to \infty}^{\widetilde{d}} f(y_n),$
also mittels \ref{FA.1.BspSt}, angewandt auf $\widetilde{d}, \widehat{d}$,
$$ \widetilde{d}(\overline{f}(\hat{x}),\overline{f}(\hat{y})) = \lim_{n \to \infty} \widetilde{d}(f(x_n),f(y_n)) \stackrel{(\ref{FA.1.21.11})}{=} \lim_{n \to \infty} \widehat{d}(x_n,y_n) = \widehat{d}(\hat{x}, \hat{y}), $$
und (\ref{FA.1.21.12}) ist gezeigt. {]}

Wegen (\ref{FA.1.21.12}) bleibt zu zeigen, daß $\overline{f}$ surjektiv ist:
Durch Vertauschung der Rollen von $\iota$ und $\tilde{\iota}$ erhält man eine stetige Fortsetzung $\overline{g} \: (\widetilde{X}, \widetilde{d}) \longrightarrow (\widehat{X},\widehat{d})$ von
$$ g := \iota \circ \left( \tilde{\iota} \: X \to \tilde{\iota}(X) \right)^{-1} \: \underbrace{(X,d)}_{\subset (\widetilde{X},\widetilde{d})} \longrightarrow (\widehat{X},\widehat{d}). $$
Dann ist $\overline{f} \circ \overline{g} \: (\widetilde{X}, \widetilde{d}) \to (\widetilde{X}, \widetilde{d})$ offenbar ebenso wie $\id_{\widetilde{X}} \: (\widetilde{X}, \widetilde{d}) \to (\widetilde{X}, \widetilde{d})$ eine stetige Fortsetzung der gleichmäßig stetigen Abbildung
$$ \id_X \: \underbrace{(X,d)}_{\subset (\widetilde{X}, \widetilde{d})} \longrightarrow (\widetilde{X}, \widetilde{d}). $$
Da $\overline{X} = \widetilde{X}$ und $(\widetilde{X}, \widetilde{d})$ vollständig ist, ergibt sich aus \ref{FA.1.20} die Eindeutigkeit einer solchen Fortsetzung, d.h.\ $\overline{f} \circ \overline{g} = \id_{\widetilde{X}}$.
Somit muß $\overline{f}$ surjektiv sein. \q

\begin{Bem*}
Der hier gegebene Beweis geht wesentlich auf \textsc{G.\ Cantor} zurück.
Er gab eine Konstruktion zur Vervollständigung der Menge der rationalen Zahlen zu der der reellen Zahlen.
\textsc{F.\ Hausdorff} verallgemeinerte diese Ideen zum hier angegebenen Beweis.
\textsc{H.\ König} hat einen weiteren Beweis gegeben, der mehr ,,funktionalanalytischen Charakter`` hat, siehe Übung \ref{FA.UA.K}.
Letztgenannter Beweis benötigt aber bereits die Vollständigkeit der Menge der reellen Zahlen.
Dies hat den Autor dieser Zeilen bewogen, den oben vorgeführten Beweis zu präsentieren.
\end{Bem*}

\subsection*{Kompakta und der Satz von \textsc{Arzel\`{a}-Ascoli}} \addcontentsline{toc}{subsection}{Kompakta und der Satz von \textsc{Arzel\`{a}-Ascoli}}

\begin{Def}[(Relative) Kompaktheit] \label{FA.1.22}
Seien $X$ ein topologischer Raum und $K$ eine Teimenge von $X$.
\begin{itemize}
\item[(i)] $K$ heißt \emph{kompakt}\index{Kompaktheit} \index{Menge!kompakte} genau dann, wenn jede Überdeckung von $K$ durch offene Teilmengen von $X$ eine endliche Teilüberdeckung besitzt; d.h.\ genauer:
Zu jeder Abbildung $I \to {\rm Top}(X)$, $i \mapsto U_i$, einer beliebigen Menge $I$ mit $\bigcup_{i \in I} U_i \supset K$ existieren endlich viele $i_1, \ldots, i_k \in I$ mit $U_{i_1} \cup \ldots \cup U_{i_k} \supset K$.

Es gilt:
\begin{eqnarray}
\mbox{$K$ ist kompakt} & \Longleftrightarrow & \mbox{Jede Überdeckung von $K$ durch} \nonumber \\
&& \mbox{offene Teilmengen von $K$ besitzt} \label{FA.1.22.S} \\
&& \mbox{eine endliche Teilüberdeckung.} \nonumber
\end{eqnarray}

{[} Zu (\ref{FA.1.22.S}): ,,$\Rightarrow$`` Sei $\left( V_i \right)_{i \in I}$ eine Überdeckung von $K$ durch offene Teilmengen des Teilraumes $K$ von $X$.
Nach Definition der Teilraumtopologie existiert für jedes $i \in I$ eine offene Menge $U_i$ von $X$ mit $V_i = U_i \cap K$.
Dann ist $\left( U_i \right)_{i \in I}$ eine Überdeckung von $K$ durch offene Teilmengen von $X$, welche nach Voraussetzung eine endliche Teilüberdeckung $U_{i_1} \cup \ldots \cup U_{i_k} \supset K$ besitzt, und es gilt folglich auch $V_{i_1} \cup \ldots \cup V_{i_k} \supset K$. 

,,$\Leftarrow$`` Sei $\left( U_i \right)_{i \in I}$ eine Überdeckung von $K$ durch offene Teilmengen von $X$.
Dann ist $\left( U_i \cap K \right)_{i \in I}$ eine Überdeckung von $K$ durch offene Teilmengen des Teilraumes $K$ von $X$, welche nach Voraussetzung eine endliche Teilüberdeckung
$$ (U_{i_1} \cap K) \cup \ldots \cup (U_{i_k} \cap K) \supset K $$
besitzt, d.h.\ $(U_{i_1} \cup \ldots \cup U_{i_k}) \cap K \supset K$.{]}
\item[(ii)] $K$ heißt \emph{relativ kompakt}\index{Kompaktheit!relative} \index{Menge!relativ kompakte} (i.Z. $\boxed{K \subset \subset X}$) genau dann, wenn $\overline{K}$ kompakt ist.
\end{itemize}
\end{Def}

\begin{Bsp} \label{FA.1.22.Bsp}
Wir versehen $\N$ mit der \emph{Topologie der kofiniten Mengen}\index{Topologie!der kofiniten Mengen}, d.h.\ per definitionem, daß die echten Teilmengen genau dann abgeschlossenen sind, wenn sie endlich sind.
Dieser topologische Raum besitzt nur kompakte Teilmengen (und ist somit selbst kompakt).
Insbesondere müssen kompakte Teilmengen i.a.\ nicht abgeschlossen sein.

Seien  nämlich $M \subset \N$ und $(U_i)_{i \in I}$ eine offene Überdeckung von $M$.
Ohne Einschränkung gelte $\forall_{i \in I} \, U_i \ne \emptyset$.
Fixiere $i_0 \in I$.
Dann existieren $k \in \N$ sowie $n_1, \ldots, n_k \in \N$ mit $U_{i_0} = \N \setminus \{ n_1, \ldots, n_k \}$.
Nach eventueller Umnumerierung der $n_i$ finden wir $l \in \{1, \ldots, k\}$ mit $\{n_1, \ldots, n_l\} \subset M$ und $\{n_{l+1}, \ldots, n_k\} \subset \N \setminus M$.
Zu jedem $j \in \{1, \ldots, l\}$ können wir ein $i_j \in I$ mit $n_j \in U_{i_j}$ wählen, und es gilt $\bigcup_{i=0}^l U_{i_j} \supset M$.
\end{Bsp}

\begin{Satz}[Äquivalente Charakterisierung der Kompaktheit] \label{FA.1.22.Durchschnitt} \index{Kompaktheit}
Es seien $X$ ein topologischer Raum und $K$ eine Teilmenge von $X$.
Dann sind die folgenden beiden Aussagen äquivalent:
\begin{itemize}
\item[(i)] $K$ ist kompakt.
\item[(ii)] Jedes zentrierte System in $K$ abgeschlossener Mengen besitzt einen nicht-leeren Schnitt, dabei heißt eine Familie $(A_i)_{i \in I}$ von Teilmengen von $X$, wobei $I$ eine beliebige Menge sei, \emph{zentriertes System}\index{zentriertes System} genau dann, wenn für jede endliche Teilmenge $I_0$ von $I$ gilt: $\D \bigcap_{i \in I_0} A_i \ne \emptyset$.
\end{itemize}
\end{Satz}

\textit{Beweis.} ,,(i) $\Rightarrow$ (ii)`` Seien also $K$ kompakt, $I$ eine beliebige Menge und $(A_i)_{i \in I}$ ein zentriertes System in $K$ abgeschlossener Mengen.
Wäre $\bigcap_{i \in I} A_i = \emptyset$, so folgte $K = K \setminus \left( \bigcap_{i \in I} A_i \right) = \bigcup_{i \in I} \left( K \setminus A_i \right)$, also existierte wegen (i) und der Offenheit von $K \setminus A_i$ in $K$ für jedes $i \in I$ eine endliche Teilmenge $I_0$ von $I$ mit
$$ K = \bigcup_{i \in I_0} \left( K \setminus A_i \right) = K \setminus \left( \bigcap_{i \in I_0} A_i \right), $$
d.h.\ $\bigcap_{i \in I_0} A_i = \emptyset$, Widerspruch!

,,(ii) $\Rightarrow$ (i)`` Gelte (ii), und angenommen, $K$ ist nicht kompakt.
Dann existiert eine Überdeckung $(U_i)_{i \in I}$ durch offene Mengen von $K$, die keine endliche Teilüberdeckung besitzt, und es gilt
$$ \emptyset = K \setminus \left( \bigcup_{i \in I} U_i \right) = \bigcap_{i \in I} \left( K \setminus U_i \right) $$
sowie für jede endliche Teilmenge $I_0$ von $I$
$$ \emptyset \ne K \setminus \left( \bigcup_{i \in I_0} U_i \right) = \bigcap_{i \in I_0} \left( K \setminus U_i \right), $$
im Widerspruch zu (ii), da $K \setminus U_i$ für jedes $i \in I$ abgeschlossen in $K$ ist. \q

\begin{Satz} \label{FA.1.24}
Seien $X$ ein topologischer Raum und $K \subset X$ eine kompakte Teilmenge von $X$.
Dann gilt:
\begin{itemize}
\item[(i)] Ist $A$ eine abgeschlossene Teilmenge von $X$ mit $A \subset K$, so ist $A$ kompakt.
\item[(ii)] Ist $X$ zusätzlich hausdorffsch, so ist $K$ abgeschlossen in $X$.
\end{itemize}
\end{Satz}

\textit{Beweis.} Zu (i): Sei $\left( U_i \right)_{i \in I}$ eine Überdeckung von $A$ durch offene Teilmengen von $X$.
Dann ist $\left( U_i \right)_{i \in I}$ zusammen mit der offenen Menge $X \setminus A$ eine offene Überdeckung von $K$ (sogar von $X$).
Wegen der Kompaktheit von $K$ existieren daher endlich viele $i_1, \ldots, i_k \in I$ mit
$$ K \subset U_{i_1} \cup \ldots \cup U_{i_k} \cup (M \setminus A), $$
also gilt auch (wegen $A \subset K$)
$$ A \subset U_{i_1} \cup \ldots \cup U_{i_k}. $$
Damit ist die Kompaktheit von $A$ gezeigt.

Zu (ii): Wir zeigen
\begin{equation} \label{FA.1.24.S}
\forall_{x \in X \setminus K} \exists_{W_x \in \U(x,X)} \, W_x \subset X \setminus K.
\end{equation}
Aus (\ref{FA.1.24.S}) folgt offenbar $X \setminus K = \bigcup_{x \in X \setminus K} W_x$, und diese Menge ist offen in $X$, d.h.\ es gilt (ii).

Zu (\ref{FA.1.24.S}): Sei $x \in X \setminus K$ fest gewählt.
Da $X$ hausdorffsch ist, existieren zu jedem $y \in K$
$$ U_y \in \U(x,X) \quad \mbox{ und } \quad V_y \in \U(y,X) $$
mit $U_y \cap V_y = \emptyset$.

$\left( V_y \right)_{y \in K}$ ist eine Überdeckung von $K$ durch offene Teilmengen von $X$.
Wegen der Kompaktheit von $K$ existieren daher endlich viele Punkte $y_1, \ldots, y_k \in K$ mit $K \subset \bigcup_{i=1}^k V_{y_i}$.

Dann gilt $W_x := \bigcap_{i=1}^k U_{y_i} \in \U(x,X)$ und $V := \bigcup_{i=1}^k V_{y_i}$ ist eine offene Obermenge von $K$ mit
$$ V \cap W_x = \left( \bigcup_{i=1}^k V_{y_i}  \right) \cap W_x = \bigcup_{i=1}^k ( \underbrace{ V_{y_i} \cap \underbrace{W_x}_{\subset U_{y_i}} }_{\subset V_{y_i} \cap U_{y_i} = \emptyset} ) = \emptyset, $$
also gilt (\ref{FA.1.24.S}). \q

\begin{Kor} \label{FA.1.24.Kor}
Seien $X$ ein Hausdorff-Raum und $K \subset X$ eine kompakte Teilmenge von $X$.

Dann gilt $\forall_{x \in X \setminus K} \exists_{U \in \U(x,X)} \exists_{V \in {\rm Top}(X), \, V \supset K} \, U \cap V = \emptyset$.
\end{Kor}

\textit{Beweis.} Mit $U := W_x$ und $V$ wie im Beweis von (\ref{FA.1.24.S}) ergibt dieser die Behauptung. \q

\begin{Satz} \label{FA.1.23}
Seien $X$ ein topologischer Raum, $I$ eine beliebige Menge und $K_i \subset X$ für jedes $i \in I$ kompakt.
Dann gilt:
\begin{itemize}
\item[(i)] $X \mbox{ hausdorffsch } \Longrightarrow \D \bigcap_{i \in I} K_i$ ist kompakt.
\item[(ii)] $\# I < \infty \Longrightarrow \D \bigcup_{i \in I} K_i \mbox{ kompakt.}$
\end{itemize}
\end{Satz}

\textit{Beweis.} (ii) ist klar, wir beweisen (i):
Da $X$ hausdorffsch ist, ist $K_i$ für jedes $i \in I$ nach \ref{FA.1.24} (ii) abgeschlossen, also ist auch $\bigcap_{i \in I} K_i$ abgeschlossen.
Ferner ist $\bigcap_{i \in I} K_i$ eine Teilmenge des Kompaktums $K_{i_0}$, wobei $i_0 \in I$ beliebig sei, und somit nach \ref{FA.1.24} (i) kompakt. \q

\begin{Satz}[Kompaktheitstreue stetiger Abbildungen] \label{FA.1.25}
Seien $X,Y$ topologische Räume, $K$ eine kompakte Teilmenge von $X$ und $f \: X \to Y$ eine stetige Abbildung.

Dann ist $f(K)$ eine kompakte Teilmenge von $Y$.
\end{Satz}

\textit{Beweis.} Sei $\left( V_i \right)_{i \in I}$ eine Überdeckung von $f(K)$ durch offene Teilmengen von $Y$.
Wegen der Stetigkeit von $f$ ist $U_i := \overline{f}^1(V_i)$ für jedes $i \in I$ eine offene Teilmenge von $X$, und es gilt
$$ \bigcup_{i \in i} U_i = \bigcup_{i \in I} \overline{f}^1(V_i) = \overline{f}^1 \left( \bigcup_{i \in I} V_i \right) \supset \overline{f}^1( f(K) ) \supset K. $$
Daher ist $\left( U_i \right)_{i \in I}$ eine Überdeckung von $K$ durch offene Teilmengen von $X$.
Wegen der Kompaktheit von $K$ existieren $i_1, \ldots, i_k \in I$ mit
$$ U_{i_1} \cup \ldots \cup U_{i_k} \supset K, $$
also gilt auch
$$ V_{i_1} \cup \ldots \cup V_{i_k} \supset f(U_{i_1}) \cup \ldots \cup f(U_{i_k}) = f(U_{i_1} \cup \ldots \cup U_{i_k}) \supset f(K), $$
d.h.\ $f(K)$ ist kompakt. \q

\begin{Satz} \label{FA.1.26} $\,$

\noindent \textbf{Vor.:} Seien $f \: X \to Y$ eine bijektive stetige Abbildung zwischen topologischen Räumen.
Zusätzlich sei $X$ kompakt und $Y$ hausdorffsch.

\noindent \textbf{Beh.:} $f$ ist ein Homöomorphismus.
\end{Satz}

\textit{Beweis.} Zu zeigen ist, daß $f^{-1} \: Y \to X$ stetig ist.
Hierzu sei $A \subset X$ eine beliebige abgeschlossene Teilmenge von $X$.
Nach (\ref{FA.1.12.2}) genügt es zu zeigen, daß dann auch das Urbild von $A$ unter $f^{-1}$ eine abgeschlossene Teilmenge von $Y$ ist.
Dieses Urbild ist gleich
$$ \{ y \in Y \, | \, f^{-1}(y) \in A \} \stackrel{f \text{ bij.}}{=} f(A). $$
Wegen der Kompaktheit von $X$ ist $A$ nach \ref{FA.1.24} (i) eine kompakte Teilmenge von $X$, also ist nach \ref{FA.1.25} $f(A)$ eine kompakte Teilmenge von $Y$.
Aus der Hausdorff-Eigenschaft von $Y$ folgt schließlich mit \ref{FA.1.24} (ii), daß $f(A)$ eine abgeschlossene Teilmenge von $Y$ ist. \q

\begin{Satz}[Satz von \textsc{Dini}] \index{Satz!von \textsc{Dini}} \label{FA.1.38}
Seien $X$ ein kompakter topologischer Raum und $(f_n)_{n \in \N}$ eine monoton wachsende (bzw.\ fallende) Folge in $\mathcal{C}(X,\R)$, wobei $\R$ mit der Standardmetrik $d_{\R}$ versehen sei, die punktweise gegen $f \in \mathcal{C}(X,\R)$ konvergiert.

Dann konvergiert $(f_n)_{n \in \N}$ gleichmäßig gegen $f$.
\end{Satz}

\textit{Beweis.} Wir führen den Beweis nur für denn Fall, daß $(f_n)_{n \in \N}$ monoton wachsend ist.
(Hieraus folgt durch Übergang von $(f_n)_{n \in \N}$ und $f$ zu $(-f_n)_{n \in \N}$ und $-f$ nämlich die Behauptung im anderen Falle.)
Es gilt also
\begin{equation} \label{FA.1.38.1}
\forall_{n,m \in \N \, m \ge n}  \, 0 \le f - f_m \le f - f_n.
\end{equation}

Sei $\varepsilon \in \R_+$.
Für jedes $n \in \N$ ist $f-f_n \: X \to \R$ stetig, also ist
\begin{equation} \label{FA.1.38.2}
U_n := \{ x \in X \, | \, |f - f_n|(x) < \varepsilon \}
\end{equation}
als Urbild von ${]} - \varepsilon, \varepsilon {[}$ unter $f - f_n$ eine in $X$ offene Menge.
Wegen der punktweisen Konvergenz von $(f_n)_{n \in \N}$ gegen $f$ bildet $(U_n)_{n \in \N}$ somit eine offene Überdeckung des Kompaktums $X$.
Folglich existieren endlich viele $n_1, \ldots, n_k \in \N$ mit
$$ X = \bigcup_{i=1}^n U_{n_i}. $$
Wegen (\ref{FA.1.38.1}) gilt aber $\forall_{n \in \N} \, U_n \subset U_{n+1}$, d.h.\
$$ X = U_{n_0} \mbox{ mit } n_0 := \max \{ n_1, \ldots, n_k \}, $$
also nach (\ref{FA.1.38.2})
$$ \forall_{x \in X} \, |f - f_{n_0}|(x) < \varepsilon, $$
woraus sich mit (\ref{FA.1.38.1}) ergibt: $\forall_{m \ge n_0} \forall_{x \in X} |f - f_m|(x) < \varepsilon.$ \q

\begin{Def}[Folgenkompaktheit] \label{FA.1.27}
Seien $X$ ein topologischer Raum und $K$ eine Teilmenge von $X$.

$K$ heißt \emph{folgenkompakt}\index{Kompaktheit!Folgen-}\index{Folgen!-kompaktheit}\index{Menge!folgenkompakte} genau dann, wenn jede Folge in $K$ eine in $K$ konvergente Teilfolge besitzt.
\end{Def}

\begin{Def}[Präkompaktheit] \label{FA.1.28}
Seien $X$ ein fastmetrischer Raum und $K$ eine Teilmenge von $X$.

$K$ heißt \emph{präkompakt}\index{Präkompaktheit}\index{Kompaktheit!Prä-}\index{Menge!präkompakte} :$\Longleftrightarrow$ $\forall_{\varepsilon \in \R_+} \exists_{x_1, \ldots, x_k \in K} \, K \subset \bigcup_{i=1}^k U_{\varepsilon}(x_i)$.
\end{Def}

\begin{Satz} \label{FA.1.29}
Seien $X$ ein fastmetrischer Raum und $K$ eine Teilmenge von $X$.
Dann gilt:
\begin{itemize}
\item[(i)] $K$ kompakt $\Longrightarrow$ $K$ präkompakt.
\item[(ii)] $K$ präkompakt $\Longleftrightarrow$ $\overline{K}$ präkompakt.
\item[(iii)] $K$ relativ kompakt $\Longrightarrow$ $K$ präkompakt.
\end{itemize}
\end{Satz}

\textit{Beweis.} Zu (i): Sei $\varepsilon \in \R_+$.
Dann ist $(U_{\varepsilon}(x))_{x \in K}$ eine offene Überdeckung der nach Voraussetzung kompakten Menge $K$, also existieren $x_1, \ldots, x_k \in K$ mit $K \subset \bigcup_{i=1}^k U_{\varepsilon}(x_i)$. 

Zu (ii): ,,$\Rightarrow$`` Sei $\varepsilon \in \R_+$.
Dann existieren nach Voraussetzung der linken Seite $x_1, \ldots, x_k \in K$ mit $K \subset \bigcup_{i=1}^k U_{\frac{\varepsilon}{2}}(x_i)$, also gilt $K \subset \bigcup_{i=1}^k \overline{U_{\frac{\varepsilon}{2}}(x_i)}$, und $\bigcup_{i=1}^k \overline{U_{\frac{\varepsilon}{2}}(x_i)}$ ist als Vereinigung endlich vieler abgeschlossener Mengen selbst abgeschlossen.
Da $\overline{K}$ die kleinste abgeschlossene Obermenge von $K$ ist, folgt
$$ \overline{K} \subset \bigcup_{i=1}^k \overline{U_{\frac{\varepsilon}{2}}(x_i)} \subset \bigcup_{i=1}^k U_{\varepsilon}(x_i). $$

,,$\Leftarrow$`` Sei $\varepsilon \in \R_+$.
Dann existieren nach Voraussetzung der rechten Seite $\tilde{x}_1, \ldots, \tilde{x}_k \in \overline{K}$ mit $\overline{K} \subset \bigcup_{i=1}^k U_{\frac{\varepsilon}{2}}(\tilde{x}_i)$.
Wähle zu jedem $i \in \{1, \ldots, k\}$ ein 
$$ x_i \in K \mbox{ mit } d(\tilde{x}_i,x_i) < \frac{\varepsilon}{2}. $$
Folglich gilt für alle $x \in U_{\frac{\varepsilon}{2}}(\tilde{x}_i)$
$$ d(x,x_i) < d(x,\tilde{x}_i) + d(\tilde{x}_i,x_i) < \varepsilon, $$
d.h.\ $U_{\frac{\varepsilon}{2}}(\tilde{x}_i) \subset U_{\varepsilon}(x_i)$.
Daher ergibt sich $K \subset \overline{K} \subset \bigcup_{i=1}^k U_{\frac{\varepsilon}{2}}(\tilde{x}_i) \subset \bigcup_{i=1}^k U_{\varepsilon}(x_i)$.

Zu (iii): Wegen (ii) genügt es zu zeigen, daß die nach Voraussetzung kompakte Menge $\overline{K}$ präkompakt ist, und dies gilt nach (i).
\q

\pagebreak
\begin{HS} \label{FA.1.30} \index{Kompaktheit} \index{Menge!kompakte} \index{Folgen!-kompaktheit} \index{Menge!folgenkompakte} \index{Präkompaktheit} \index{Kompaktheit!Prä-} \index{Menge!präkompakte} \index{Menge!volständige} \index{Vollständigkeit}
Seien $X$ ein fastmetrischer Raum und $K$ eine Teilmenge von $X$.
Dann sind die folgenden Aussagen paarweise äquivalent:
\begin{itemize}
\item[(i)] $K$ ist kompakt.
\item[(ii)] $K$ ist folgenkompakt
\item[(iii)] $K$ ist präkompakt und vollständig.
\end{itemize}
\end{HS}

\textit{Beweis.} ,,(i) $\Rightarrow$ (ii)`` Angenommen, $K$ ist nicht folgenkompakt.
Dann existiert eine Folge $\left( x_n \right)_{n \in \N}$ in $K$ derart, daß keine ihrer Teilfolgen gegen ein Element von $K$ konvergiert.
Wir behaupten
\begin{equation} \label{FA.1.30.1}
\forall_{y \in K} \exists_{\varepsilon \in \R_+} \, \# \{ r \in \N \, | \, x_r \in U_{\varepsilon}(y) \} < \infty.
\end{equation}

{[} Zu (\ref{FA.1.30.1}): Angenommen, es existiert $y \in K$ mit
\begin{equation} \label{FA.1.30.2}
\forall_{\varepsilon \in \R_+} \, \# \{ r \in \N \, | \, x_r \in U_{\varepsilon}(y) \} = \infty.
\end{equation}
Wir definieren dann rekursiv eine Teilfolge $\left( i_n \right)_{n \in \N}$ von $\left( n \right)_{n \in \N}$ durch
$$ i_0 := \min \underbrace{\{ r \in \N \, | \, x_r \in U_1(y) \}}_{\stackrel{(\ref{FA.1.30.2})}{\ne} \emptyset} $$
und
$$ \forall_{n \in \N} \, i_{n+1} := \min \underbrace{\{ r \in \N \, | \, r > i_n \wedge x_r \in U_{\frac{1}{n+2}}(y) \}}_{\stackrel{(\ref{FA.1.30.2})}{\ne} \emptyset}. $$
Dann ist $\left( x_{i_n} \right)_{n \in \N}$ eine Teilfolge von $\left( x_n \right)_{n \in \N}$, und wegen $x_{i_n} \in U_{\frac{1}{n+1}}(y)$ für alle $n \in \N$ konvergiert diese Teilfolge gegen $y \in K$ im Widerspruch zu oben. {]}

Zu jedem $y \in K$ können wir gemäß (\ref{FA.1.30.1}) eine Zahl $\varepsilon_y \in \R_+$ wählen mit
\begin{equation} \label{FA.1.30.3}
\# \{ r \in \N \, | \, x_r \in U_{\varepsilon_q}(y) \} < \infty.
\end{equation}
Dann ist $\left( U_{\varepsilon_y}(y) \right)_{y \in K}$ eine offene Überdeckung von $K$.
Wegen der Kompaktheit von $K$ existieren endlich viele Punkte $y_1, \ldots, y_m \in K$ mit
$$ K \subset \bigcup_{j=1}^m U_{\varepsilon_{y_j}}(y_j). $$
Hieraus folgt
$$ \N = \bigcup_{j=1}^m \{ r \in \N \, | \, x_r \in U_{\varepsilon_{y_j}}(y_j) \}, $$
und die rechte Seite ist nach (\ref{FA.1.30.3}) eine endliche Menge, im Widerspruch dazu, daß $\N$ bekanntlich eine unendliche Menge ist.

,,(ii) $\Rightarrow$ (iii)`` 1.) Zur Präkompaktheit: Ohne Beschränkung sei der Allgemeinheit $K \ne \emptyset$. 
Wir bezeichnen mit $\widetilde{\mathfrak{P}}(K) \subset \mathfrak{P}(K)$ die Menge aller endlichen nicht-leeren Teilmengen von $K$.

Angenommen, $K$ ist nicht präkompakt, d.h.\ es existiert $\varepsilon \in \R_+$ mit 
$$ \forall_{M \in \widetilde{\mathfrak{P}}(K)} \, \underbrace{K \nsubset \bigcup_{x \in M} U_{\varepsilon}(x)}_{\Leftrightarrow K \setminus \bigcup_{x \in M} U_{\varepsilon}(x) \ne \emptyset}. $$
Nach dem Auswahlaxiom existiert dann eine Abbildung $\varphi \: \widetilde{\mathfrak{P}}(K) \to K$ mit
\begin{equation} \label{FA.1.30.4}
\forall_{M \in \widetilde{\mathfrak{P}}(K)} \, \varphi(M) \in K \setminus \bigcup_{x \in M} U_{\varepsilon}(x).
\end{equation}

Wir wählen $x_0 \in K$ und definieren rekursiv
\begin{equation} \label{FA.1.30.5}
\forall_{n \in \N} \, x_{n+1} := \varphi(\{x_0, \ldots, x_n\}) \in K.
\end{equation}
Dann ist $\left( x_n \right)_{n \in \N}$ eine Folge in $K$, und es gilt
\begin{equation} \label{FA.1.30.6}
\forall_{j,k \in \N} \, \left( j<k \Longrightarrow d(p_j,p_k) \ge \varepsilon \right).
\end{equation}

{[} Denn für $j<k$ gilt $j \le k-1$ und daher 
$$ x_k \stackrel{(\ref{FA.1.30.5})}{=} \varphi(\{x_0, \ldots, x_{k-1}\}) \stackrel{(\ref{FA.1.30.4})}{\in} \left( K \setminus \bigcup_{i=0}^{k-1} U_{\varepsilon}(x_i) \right) \subset K \setminus U_{\varepsilon}(x_j), $$
also $d(x_k,x_j) \ge \varepsilon$. {]}

Aufgrund der Folgenkompaktheit von $K$ können wir ohne Beschränkung der Allgemeinheit annehmen, daß $x \in K$ mit $\lim_{n \to \infty} x_n = x$ existiert.
Dann existiert ein $n_0 \in \N$ mit 
$$\D \forall_{n \in \N, \, n \ge n_0} \, d(x_n,x) < \frac{\varepsilon}{2}, $$
also auch
$$ \forall_{m,n \in \N, \, n_0 \le n < m} \, d(x_n,x_m) \le d(x_n,x) + d(x,x_m) < \varepsilon, $$
im Widerspruch zu (\ref{FA.1.30.6}).

2.) Zur Vollständigkeit: Sei $\left( x_n \right)_{n \in \N}$ eine Cauchyfolge in $K$.
Da $K$ folgenkompakt ist, existieren eine Teilfolge $\left( x_{i_n} \right)_{n \in \N}$ von $\left( x_n \right)_{n \in \N}$ und $x \in K$ mit
$$ \lim_{n \to \infty} x_{i_n} = x, $$
und wir behaupten, daß auch $\left( x_n \right)_{n \in \N}$ gegen $x$ konvergiert.

Hierzu sei $\varepsilon \in \R_+$.
Da $\left( x_n \right)_{n \in \N}$ eine Cauchyfolge ist, so existiert $n_0 \in \N$ mit
$$ \forall_{n,m \in \N, n,m \ge n_0} \, d(x_n,x_m) < \frac{\varepsilon}{2}. $$
Es existiert weiter $m_0 \in \N$ derart, daß gilt
$$ \forall_{m \in \N, m \ge m_0} \, i_m \ge n_0 \, \wedge \, d(x_{i_m},x) < \frac{\varepsilon}{2}. $$

Somit ergibt sich für jedes $n,m \in \N$ mit $n \ge n_0$ und $m \ge m_0$
$$ d(x_n,x) \le d(x_n,x_{i_m}) + d(x_{i_m},p) < \varepsilon, $$
womit die Behauptung bewiesen ist.

\pagebreak
,,(iii) $\Rightarrow$ (i)`` Ohne Beschränkung der Allgemeinheit sei $K \ne \emptyset$, und angenommen, $K$ ist nicht kompakt.
Dann existiert eine offene Überdeckung $(U_i)_{i \in I}$ von $K$ derart, daß gilt
\begin{equation} \label{FA.1.30.7}
\forall_{I_0 \subset I} \, \left( \# I_0 < \infty \Longrightarrow K \nsubset \bigcup_{i \in I_0} U_i \right).
\end{equation}
Wir definieren rekursiv eine Folge $(K_n)_{n \in \N}$ nicht-leerer Teilmengen von $K$ derart, daß für alle $n \in \N$ 
\begin{gather}
n \in \N_+ \Longrightarrow K_n \subset K_{n-1}, \label{FA.1.30.8} \\
\sup \{ d(x,\tilde{x}) \, | \, x,\tilde{x} \in K_n \} \le \frac{1}{n} (:= \infty, \mbox{ falls } n=0), \label{FA.1.30.9} \\
\forall_{I_0 \subset I} \, \left( \# I_0 < \infty \Longrightarrow K_n \nsubset \bigcup_{i \in I_0} U_i \right) \label{FA.1.30.10}
\end{gather}
gilt, wie folgt:

$K_0 := K$ erfüllt (\ref{FA.1.30.8}) - (\ref{FA.1.30.10}) für $n=0$, denn die erste dieser Aussagen ist für $n=0$ leer, die zweite trivial, und die dritte gilt nach (\ref{FA.1.30.7}).

Seien $n \in \N$ und $K_0, \ldots, K_n$ bereits definiert.
Wegen der Präkompaktheit von $K$ existieren $x_1, \ldots, x_k \in K$ mit
$$ K \subset \bigcup_{j=1}^k U_{\frac{1}{2(n+1)}}(x_j). $$
Setze $K_{n+1,j} := K_n \cap U_{\frac{1}{2(n+1)}}(x_j)$ für $j \in \{1, \ldots, k\}$.
Dann gilt wegen $K_n \subset K$
\begin{equation} \label{FA.1.30.11}
K_n = \bigcup_{j=1}^k K_{n+1,j}
\end{equation}
sowie (\ref{FA.1.30.8}), (\ref{FA.1.30.9}) für $n+1$ mit $K_{n+1,j}$ anstelle von $K_{n+1}$ für alle $j \in \{1, \ldots, k\}$.
Wir zeigen
\begin{equation} \label{FA.1.30.12}
\exists_{j_0 \in \{1, \ldots, k\}} \forall_{I_{j_0} \subset I} \, \left( \# I_{j_0} < \infty \Longrightarrow K_{n+1,j_0} \nsubset \bigcup_{i \in I_{j_0}} U_i \right),
\end{equation}
d.h.\ $K_{n+1} := K_{n+1,j_0}$ erfüllt (\ref{FA.1.30.8}) - (\ref{FA.1.30.10}) für $n+1$.

{[} Angenommen, (\ref{FA.1.30.12}) ist falsch, d.h.\
$$ \forall_{j \in \{1, \ldots, k\}} \exists_{I_j \subset I} \, \left( \# I_j < \infty \, \wedge \, K_{n+1,j} \subset \bigcup_{i \in I_j} U_i \right). $$
Dann ist $I_0 := \bigcup_{j=1}^k I_j$ eine endliche Menge, und es gilt
$$ K_n \stackrel{(\ref{FA.1.30.11})}{=} \bigcup_{j=1}^k K_{n+1,j} \subset \bigcup_{i \in I_0} U_i, $$
im Widerspruch zu (\ref{FA.1.30.10}) für $n$. {]}

Wähle nun zu jedem $n \in \N$ ein $x_n \in K_n$.
Wegen (\ref{FA.1.30.8}), (\ref{FA.1.30.9}) ist $(x_n)_{n \in \N}$ dann offenbar eine Cauchyfolge in $K \subset \bigcup_{i \in I} U_i$.
Da $K$ vollständig ist, existieren $x \in K$ und $i_0 \in I$ derart, daß gilt
\begin{equation} \label{FA.1.30.13}
\lim_{n \to \infty} x_n = x \in U_{i_0}.
\end{equation}
$U_{i_0}$ ist offen, also existiert eine Zahl $\varepsilon \in \R_+$ mit
\begin{equation} \label{FA.1.30.14}
U_{2 \varepsilon}(x) \subset U_{i_0},
\end{equation}
und wegen (\ref{FA.1.30.13}) existiert $n_0 \in \N$ mit
\begin{equation} \label{FA.1.30.15}
\forall_{n \in \N, \, n \ge n_0} \, d(x_n,x) < \varepsilon.
\end{equation}

Sei schließlich $n \in \N_+$ mit $n \ge n_0$ und $\frac{1}{n} < \varepsilon$.
Wir zeigen 
$$ K_n \subset U_{i_0}, $$
womit ein Widerspruch zu (\ref{FA.1.30.10}) hergeleitet ist.
Sei also $y \in K_n$ beliebig.
Dann gilt
$$ d(y,x) \le d(y,x_n) + d(x_n,x) \stackrel{(\ref{FA.1.30.9}), (\ref{FA.1.30.15})}{<} \frac{1}{n} + \varepsilon < 2 \varepsilon, $$
$y \in U_{2 \varepsilon}(x) \stackrel{(\ref{FA.1.30.14})}{\subset} U_{i_0}$. \q

\begin{Kor} \label{FA.1.30.K}
Seien $X$ ein fastmetrischer Raum und $K$ eine Teilmenge von $X$.
Dann gilt:

$K$ präkompakt und $\overline{K}$ vollständig $\Longrightarrow$ $K$ relativ kompakt.
\end{Kor}

\textit{Beweis.} Das Korollar folgt sofort aus \ref{FA.1.29} (ii) und \ref{FA.1.30}. \q

\begin{Satz} \label{FA.1.31} \index{Stetigkeit} \index{Stetigkeit!gleichmäßige} \index{Abbildung!stetige} \index{Abbildung!stetige!gleichmäßig} $\,$ 

\noindent \textbf{Vor.:} Seien $X,Y$ fastmetrische Räume und $f \: X \to Y$ eine stetige Abbildung.
Ferner sei $X$ kompakt.

\noindent \textbf{Beh.:} $f$ stetig $\Longleftrightarrow$ $f$ gleichmäßig stetig.
\end{Satz}

\textit{Beweis.} ,,$\Leftarrow$`` ist trivial.

,,$\Rightarrow$`` Angenommen, $f$ ist nicht gleichmäßig stetig. 
Dann gibt es eine Zahl $\varepsilon \in \R_+$ derart, daß gilt
$$ \forall_{\delta \in \R_+} \exists_{x, \tilde{x} \in X} \, \left( d(x,\tilde{x}) < \delta \, \wedge \, d(f(x),f(\tilde{x})) \ge \varepsilon \right). $$
Hieraus folgt die Existenz von Folgen $(x_n)_{n \in \N}, (\tilde{x}_n)_{n \in \N}$ in $X$ mit
\begin{equation} \label{FA.1.31.1}
\forall_{n \in \N} \, d(x_n,\tilde{x}_n) < \frac{1}{n+1} \, \wedge \, d(f(x_n),f(\tilde{x}_n)) \ge \varepsilon. 
\end{equation}

$X$ ist als Kompaktum folgenkompakt, also besitzt $(x_n)_{n \in \N}$ eine etwa gegen $x \in X$ konvergente Teilfolge $(x_{i_n})_{n \in \N}$.
Aus $\forall_{n \in \N} \, d(x_{i_n},\tilde{x}_{i_n}) < \frac{1}{i_n + 1}$ folgt dann wegen
$$ d(\tilde{x}_{i_n},x) \le d(\tilde{x}_{i_n},x_{i_n}) + d(x_{i_n},x), $$
daß auch $(\tilde{x}_{i_n})_{n \in \N}$ gegen $x$ konvergiert, also ergibt die Stetigkeit von $f$ in $x$
$$ \lim_{n \to \infty} f(x_{i_n}) = f(x) = \lim_{n \to \infty} f(\tilde{x}_{i_n}), $$
d.h.\
$$ \lim_{n \to \infty} \underbrace{d(f(x_{i_n}), f(\tilde{x}_{i_n}))}_{\le d(f(x_{i_n}), f(x)) + d(f(x), f(\tilde{x}_{i_n}))} = 0, $$
im Widerspruch zu (\ref{FA.1.31.1}). \q

\begin{Def} \label{FA.1.32}
Seien $X$ ein fastmetrischer Raum und $A$ eine nicht-leere Teilmenge von $X$.
Wir definieren den \emph{Durchmesser von $A$}\index{Durchmesser} als
$$ \boxed{{\rm diam}(A)} := \sup \{ d(x,y) \, | \, x,y \in A \} \in {[} 0, \infty {]} $$
und nennen $A$ \emph{beschränkt}\index{Menge!beschränkte}, falls ${\rm diam}(A) < \infty$ gilt.
Zusätzlich definieren wir die leere Menge als beschränkt.

\begin{Bem*}
Für eine Teilmenge $A$ eines fastmetrischen Raumes $X$ gilt
$$ A \mbox{ beschränkt} \Longleftrightarrow \exists_{x \in X} \exists_{C \in R_+} \forall_{y \in A} \, d(x,y) \le C, $$
denn ,,$\Rightarrow$`` ist trivial und ,,$\Leftarrow$`` folgt aus der Dreiecksungleichung.
\end{Bem*}  
\end{Def}

\begin{Satz} \label{FA.1.33}
Seien $X$ ein metrischer Raum und $K$ eine Teilmenge von $X$.
Dann gilt:

$K$ kompakt $\Longrightarrow$ $K$ beschränkt.
\end{Satz}

\textit{Beweis.} Wäre $K$ nicht beschränkt, so existierten Folgen $\left( x_n \right)_{n \in \N}$, $\left( y_n \right)_{n \in \N}$ in $K$ mit $\lim_{n \to \infty} d(x_n,y_n) = \infty$.
Da $K$ kompakt und somit folgenkompakt ist, können wir zunächst ohne Beschränkung der Allgemeinheit annehmen, daß $\left( x_n \right)_{n \in \N}$ gegen ein $x \in K$ konvergiert, denn andernfalls ersetzen wir $(x_n)_{n \in \N}$ durch eine in $K$ konvergente Teilfolge.
Sodann können ebenfalls wir annehmen, daß $\left( y_n \right)_{n \in \N}$ gegen ein $y \in K$ konvergiert.
Nun folgt für jedes $n \in \N$
$$ d(x_n,y_n) \le d(x_n,x) + d(x,y) + d(y,y_n) \stackrel{n \to \infty}{\longrightarrow} d(x,y) < \infty, $$
im Widerspruch zu $\lim_{n \to \infty} d(x_n,y_n) = \infty$. \q

\begin{Satz}[von \textsc{Weierstraß}] \index{Satz!von \textsc{Weierstraß}} \label{FA.1.34} $\,$

\noindent \textbf{Vor.:} Seien $X$ ein topologischer Raum, $K$ eine nicht-leere kompakte Teilmenge von $X$ und $f \: X \to \R$ eine stetige Abbildung, wobei wir $\R$ mit der üblichen Metrik $d_{\R}$ versehen.

\noindent \textbf{Beh.:} $f$ nimmt auf $K$ ein Maximum und ein Minimum an.
\end{Satz}

\textit{Beweis.} $K$ ist kompakt und $f$ stetig, also ist nach \ref{FA.1.25} $f(K)$ eine kompakte Teilmenge von $\R$, die dann nach \ref{FA.1.33} und \ref{FA.1.24} (ii) beschränkt und abgeschlossen ist.
Hieraus folgt offenbar die Behauptung. \q

\begin{Def}[Gleichgradige Stetigkeit] \index{Stetigkeit!gleichgradige}\index{Abbildung!gleichgradige stetige --en} \label{FA.1.35}
Seien $X$ ein topologischer Raum, $Y$ ein fastmetrischer Raum und $\mathcal{F} \subset Y^X$.
\begin{itemize}
\item[(i)] Sei $x \in X$.
$\mathcal{F}$ heißt \emph{gleichgradig stetig in $x$} genau dann, wenn gilt
$$ \forall_{\varepsilon \in \R_+} \exists_{U \in \U(x,X)} \forall_{f \in \mathcal{F}} \, f(U) \subset U_{\varepsilon}(f(x)). $$
\item[(ii)] $\mathcal{F}$ heißt \emph{gleichgradig stetig} $:\Longleftrightarrow$ $\forall_{x \in X} \, \mathcal{F}$ gleichgradig stetig in $x$.
\end{itemize}
\end{Def}

\begin{Bsp*}
$\{ \underbrace{x^n}_{=: f_n}|_{[ 0, 1]} \, | \, n \in \N \}$ ist in allen $t \in [0,1[$ gleichgradig stetig, nicht aber im Punkte $1$.

\textit{Beweisskizze.} 1.) Für die erste Aussage genügt es zu zeigen, daß $$\{ f_n|_{[ 0, C]} \, | \, n \in \N \}$$ für alle $C \in {]}0,1{[}$ gleichgradig stetig ist.
Dies wiederum ergibt sich aus
\begin{equation} \label{FA.1.35.1}
L := \sup \{ \| {f_n}' \|_{\infty} \, | \, n \in \N \} < \infty,
\end{equation}
wobei
$$ \| {f_n}' \|_{\infty} := \sup \{ |{f_n}'(t)| \, | \, t \in [0,C] \} = n \, C^{n-1}, $$
denn dann folgt aus dem Mittelwertsatz
$$ \forall_{n \in \N} \forall_{t,t_0 \in [0,C]} \, |f_n(t) - f_n(t_0)| \le L \, |t - t_0|, $$
und für jedes $t_0 \in [0,C]$ sowie $\varepsilon \in \R_+$ gilt
$$ \forall_{n \in \N} \forall_{t \in [0,C]} \, \left( |t - t_0| < \frac{\varepsilon}{L} \Longrightarrow |f_n(t) - f_n(t_0)| < \varepsilon \right). $$
Zum Nachweis von (\ref{FA.1.35.1}) zeigt man, daß die Funktion $g := x \, C^{x-1}|_{{[}0, \infty{[}}$ in $- \frac{1}{\ln(c)}$ ein globales Maximum besitzt.

2.) Die zweite Aussage folgt daraus, daß $(f_n|_{{[}0,1{[}})_{n \in \N}$ punktweise gegen die konstante Funktion vom Wert $0$ konvergiert und $\forall_{n \in \N} \, f_n(1) = 1$ gilt. \q
\end{Bsp*}

\begin{Bem*}
Eine Abbildung $f \: X \to Y$ zwischen fastmetrischen Räumen heißt \emph{Lip\-schitz-stetig}\index{Stetigkeit!Lipschitz-}\index{Lipschitz-stetig}\index{Abbildung!Lipschitz-stetige}, wenn $\exists_{L \in \R} \forall_{x,\tilde{x} \in X} \, d(f(x),f(\tilde{x})) \le L \, d(x,\tilde{x})$ gilt. 
Solche Abbildungen sind offenbar gleichmäßig stetig.
\end{Bem*}

\begin{Satz} \label{FA.1.37}
Seien $X$ ein kompakter topologischer Raum, $Y$ ein fastmetrischer Raum, $(f_n)_{n \in \N}$ eine Folge in $\mathcal{C}(X,Y)$ und $f \in Y^X$.

Dann konvergiert $(f_n)_{n \in \N}$ genau dann gleichmäßig gegen $f$, wenn $(f_n)_{n \in \N}$ punktweise gegen $f$ konvergiert und $\{ f_n \, | \, n \in \N \}$ gleichgradig stetig ist.
\end{Satz}

\textit{Beweis.} ,,$\Rightarrow$`` Trivialerweise ist nur zu zeigen, daß $\{ f_n \, | \, n \in \N \}$ gleichgradig stetig ist.
Seien hierzu $x_0 \in X$ und $\varepsilon \in \R_+$.
Da $\{ f_n \, | \, n \in \N \}$ gleichmäßig gegen $f$ konvergiert, existiert zum einen $n_0 \in \N$ derart, daß gilt
\begin{equation*} \label{FA.1.37.1}
\forall_{n \in \N, \, n \ge n_0} \forall_{x \in X} \, d(f_n(x), f(x)) < \frac{\varepsilon}{3},
\end{equation*}
und zum anderen ist wegen der Stetigkeit von $f_n$ für jedes $n \in \N$ auch $f$ stetig, insbes.\ in $x_0$.
Daher gibt es ein $U \in \U(x_0,X)$ mit
\begin{equation*} \label{FA.1.37.2}
\forall_{x \in U} \, d(f(x), f(x_0)) < \frac{\varepsilon}{3},
\end{equation*}
und es folgt für alle $n \in \N$ mit $n \ge n_0$ sowie jedes $x \in U$
\begin{equation*} \label{FA.1.37.3}
d(f_n(x), f_n(x_0)) \le d(f_n(x),f(x)) + d(f(x),f(x_0)) + d(f(x_0),f_n(x_0)) < \varepsilon.
\end{equation*}
Damit ist gezeigt, daß $\{ f_n \, | \, n \in \N \, \wedge \, n \ge n_0 \}$ in $x_0$ gleichgradig stetig ist.
Die Hinzunahme der endlich vielen in $x_0$ stetigen Abbildungen $f_0, \ldots, f_{n_0 - 1}$ spielt offenbar im Bezug auf die gleichgradige Stetigkeit in $x_0$ keine Rolle.

,,$\Leftarrow$`` Sei $\varepsilon \in \R_+$.
Da $\{ f_n \, | \, n \in \N \}$ gleichgradig stetig ist, existiert zu jedem $x_0 \in X$ ein $U_{x_0} \in \U(x_0,X)$ derart, daß gilt
\begin{equation} \label{FA.1.37.4}
\forall_{n \in \N} \forall_{x \in U_{x_0}} \, d(f_n(x), f_n(x_0)) < \frac{\varepsilon}{3},
\end{equation}
also aufgrund der punktweisen Konvergenz von $(f_n)_{n \in \N}$ gegen $f$ auch
\begin{equation} \label{FA.1.37.5}
\forall_{x \in U_{x_0}} \, d(f(x), f(x_0)) \le \frac{\varepsilon}{3},
\end{equation}
beachte, daß $d$ auf $U_{\varepsilon}(f_n(x_0))$ und $U_{\varepsilon}(f(x))$ eine Metrik und damit stetig ist.

$(U_x)_{x \in X}$ ist eine offene Überdeckung des Kompaktums $X$, also existieren $x_1, \ldots, x_k \in X$ mit 
\begin{equation} \label{FA.1.37.6}
\bigcup_{i=1}^k U_{x_i} = X.
\end{equation}
Für $i \in \{1, \ldots, k\}$ konvergiert $(f_n(x_i))_{n \in \N}$ nach Voraussetzung der rechten Seite gegen $f(x_i)$, und wir können offenbar ein (gemeinsames) $n_0 \in \N$ wählen so, daß gilt
\begin{equation} \label{FA.1.37.7}
\forall_{i \in \{1, \ldots, k\}} \forall_{n \in \N, \, n \ge n_0} \, d(f_n(x_i),f(x_i)) < \frac{\varepsilon}{3}.
\end{equation}

Sei nun $x \in X$ beliebig.
Wegen (\ref{FA.1.37.6}) existiert $i \in \{1, \ldots, k\}$ mit $x \in U_{x_i}$.
Dann gilt für alle $n \in \N$ mit $n \ge n_0$
$$ d(f_n(x), f(x)) \le \underbrace{d(f_n(x), f_n(x_i))}_{\stackrel{(\ref{FA.1.37.4})}{<} \frac{\varepsilon}{3}} + \underbrace{d(f_n(x_i), f(x_i))}_{\stackrel{(\ref{FA.1.37.7})}{<} \frac{\varepsilon}{3}} + \underbrace{d(f(x_i), f(x))}_{\stackrel{(\ref{FA.1.37.5})}{\le} \frac{\varepsilon}{3}} < \varepsilon, $$
d.h.\ die linke Seite ist gezeigt. \q

\begin{HS}[Satz von \textsc{Arzel\`{a}-Ascoli}] \index{Satz!von \textsc{Arzel\`{a}-Ascoli}} \label{FA.1.36} $\,$

\noindent \textbf{Vor.:} Seien $X$ ein kompakter topologischer Raum, $(Y,d)$ ein fastmetrischer Raum und $\mathcal{H}$ eine nicht-leere Teilmenge von $\mathcal{C}(X,Y)$.
Wir versehen $\mathcal{C}(X,Y)$ mit der Fastmetrik $d_{\infty}$. 
(Ist $d$ sogar eine Metrik, so folgt aus \ref{FA.1.34} wegen der Kompaktheit von $X$ und der Stetigkeit von $d$, vgl.\ \ref{FA.1.BspSt}, daß $d_{\infty}$ eine Metrik auf $\mathcal{C}(X,Y)$ ist.)

\noindent \textbf{Beh.:} Die folgenden Aussagen sind paarweise äquivalent:
\begin{itemize}
\item[(i)] $\mathcal{H}$ ist gleichgradig stetig und alle Orbiten sind relativ kompakte Teimengen von $Y$, letzteres bedeutet
\begin{equation} \label{FA.1.36.0}
\forall_{x \in X} \, \mathcal{H}(x) := \{ f(x) \, | \, f \in \mathcal{H} \} \subset \subset Y.
\end{equation}
\item[(ii)] $\mathcal{H} \subset \subset \mathcal{C}(X,Y).$
\item[(iii)] Jede Folge $(f_n)_{n \in \N}$ in $\mathcal{H}$ besitzt eine Teilfolge die gleichmäßig gegen eine stetige Abbildung $f \: X \to Y$ konvergiert.
\end{itemize}
\end{HS}

\pagebreak
\textit{Beweis.} ,,(i) $\Rightarrow$ (ii)`` Wegen \ref{FA.1.30.K} haben wir zu zeigen
\begin{gather}
\mathcal{H} \mbox{ präkompakt,} \label{FA.1.36.1} \\
\overline{\mathcal{H}} \mbox{ vollständig.} \label{FA.1.36.2}
\end{gather}

Zu (\ref{FA.1.36.1}): Sei $\varepsilon \in \R_+$.
Wegen der gleichgradigen Stetigkeit von $\mathcal{H}$ existiert zu jedem $x_0 \in X$ ein $U_{x_0} \in \U(x_0,X)$ mit
\begin{equation} \label{FA.1.36.3}
\forall_{f \in \mathcal{H}} \, f(U_{x_0}) \subset U_{\frac{\varepsilon}{5}}(f(x_0)),
\end{equation}
und $(U_x)_{x \in X}$ ist eine offene Überdeckung des Kompaktums $X$.
Daher gibt es $x_1, \ldots, x_k \in X$ derart, daß gilt
\begin{equation} \label{FA.1.36.4}
X = \bigcup_{i=1}^k U_{x_i}.
\end{equation}
Wir definieren
$$ K := \bigcup_{i=1}^k \mathcal{H}(x_i) \subset \bigcup_{i=1}^k \overline{\mathcal{H}(x_i)}. $$
Dann ist $\overline{K}$ eine abgeschlossene Teilmenge der nach Voraussetzung von (i) kompakten Menge $\bigcup_{i=1}^k \overline{\mathcal{H}(x_i)}$, denn $\bigcup_{i=1}^k \overline{\mathcal{H}(x_i)}$ ist eine abgschlossene Obermenge von $K$ und $\overline{K}$ ist die kleinste solche.
Aus \ref{FA.1.24} (i) folgt die Kompaktheit von $\overline{K}$, d.h.\ $K$ ist relativ kompakt und somit nach \ref{FA.1.29} (iii) präkompakt.
Daher existieren $y_1, \ldots, y_l \in Y$ mit
\begin{equation} \label{FA.1.36.5}
\bigcup_{i=1}^k \mathcal{H}(x_i) = K \subset \bigcup_{j=1}^l U_{\frac{\varepsilon}{5}}(y_j).
\end{equation}

Wir betrachten nun die endliche Menge von Abbildung
$$ \Phi := \{1, \ldots, l \}^{\{1, \ldots, k \}} $$
(der Mächtigkeit $l^k$, falls $\mathcal{H} \ne \emptyset$) und für jedes $\varphi \in \Phi$
\begin{equation} \label{FA.1.36.6}
\mathcal{H}_{\varphi} := \{ f \in \mathcal{H} \, | \, \forall_{i \in \{1, \ldots, k\}} \, f(x_i) \subset U_{\frac{\varepsilon}{5}}(y_{\varphi(i)}) \}. 
\end{equation}
Dann gilt
\begin{equation} \label{FA.1.36.7}
\mathcal{H} = \bigcup_{\varphi \in \Phi} \mathcal{H}_{\varphi},
\end{equation}
denn ,,$\supset$`` ist trivial und ,,$\subset$`` folgt daraus, daß nach (\ref{FA.1.36.5}) für jedes $i \in \{1, \ldots, k\}$ sowie $f \in \mathcal{H}$
$$ f(x_i) \in \bigcup_{j=1}^l U_{\frac{\varepsilon}{5}}(y_j) $$
gilt, also ein $\varphi \in \Phi$ mit 
$$ \forall_{i \in \{1, \ldots, k\}} \, f(x_i) \subset U_{\frac{\varepsilon}{5}}(y_{\varphi(i)}) $$
existiert.
(Definiere $\varphi$ durch $\forall_{i \in \{1, \ldots, k\}} \, \varphi(i) = j$, wobei $f(x_i) \in U_{\frac{\varepsilon}{5}}(y_j)$.)

Wir behaupten
\begin{equation} \label{FA.1.36.8}
\forall_{\varphi \in \Phi} \, \forall_{f,g \in \mathcal{H}_{\varphi}} \, d_{\infty}(f,g) < \varepsilon.
\end{equation}

{[} Zu (\ref{FA.1.36.8}): Seien $\varphi \in \Phi$, $f,g \in \mathcal{H}_{\varphi}$ und $x \in X$ beliebig.
Dann existiert nach (\ref{FA.1.36.4}) ein $i \in \{1, \ldots, k\}$ mit $x \in U_{x_i}$, und es gilt
\begin{eqnarray*}
\lefteqn{d(f(x),g(x))} ~~ \\
& \le & \underbrace{d(f(x),f(x_i))}_{\stackrel{(\ref{FA.1.36.3})}{<} \frac{\varepsilon}{5}} + \underbrace{d(f(x_i), y_{\varphi(i)})}_{\stackrel{(\ref{FA.1.36.6})}{<} \frac{\varepsilon}{5}} + \underbrace{d(y_{\varphi(i)}, g(x_i))}_{\stackrel{(\ref{FA.1.36.6})}{<} \frac{\varepsilon}{5}} + \underbrace{d(g(x_i), g(x))}_{\stackrel{(\ref{FA.1.36.3})}{<} \frac{\varepsilon}{5}} < \frac{4 \varepsilon}{5},
\end{eqnarray*}
also $d_{\infty}(f,g) \le \frac{4 \varepsilon}{5} < \varepsilon$. {]}

Nun ist $\Phi_0 := \{ \varphi \in \Phi \, | \, \mathcal{H}_{\varphi} \ne \emptyset \}$ mit $\Phi$ eine endliche Menge.
Zu jedem $\varphi \in \Phi_0$ existiert ein $f_{\varphi} \in \mathcal{H}_{\varphi}$, und es gilt nach (\ref{FA.1.36.8})
$$ \mathcal{H}_{\varphi} \subset U^{d_{\infty}}_{\varepsilon}(f_{\varphi}), $$
also nach (\ref{FA.1.36.7})
$$ \mathcal{H} \subset \bigcup_{\varphi \in \Phi_0} U^{d_{\infty}}_{\varepsilon}(f_{\varphi}), $$
womit die Präkompaktheit von $\mathcal{H}$ gezeigt ist.

Zu (\ref{FA.1.36.2}): Für jedes $x \in X$ ist die Evaluationsbbildung
$$ E_x \: \mathcal{C}(X,Y) \longrightarrow Y, ~~ f \longmapsto f(x), $$
gleichmäßig setig, denn für alle $\varepsilon \in \R_+$ gilt mit $\delta := \varepsilon$
$$ \forall_{f,g \in \mathcal{C}(X,Y)} \, \left( d_{\infty}(f,g) < \delta \Longrightarrow d(f(x),g(x)) < \varepsilon \right). $$
Hieraus folgt -- auch ohne die Voraussetzungen von (i) --
\begin{equation} \label{FA.1.36.9}
\forall_{x \in X} \, \overline{\mathcal{H}}(x) = \{ f(x) \, | \, f \in \overline{\mathcal{H}} \} = E_x(\overline{\mathcal{H}}) \subset \overline{\mathcal{H}(x)}.
\end{equation}

{[} Zur letzten Inklusion in (\ref{FA.1.36.9}): Seien $x \in X$ und $f \in \overline{\mathcal{H}}$.
Dann existiert eine Folge $(f_n)_{n \in \N}$ aus $\mathcal{H}$, die bzgl.\ $d_{\infty}$ gegen $f$ konvergiert.
Da $E_x$ stetig (sogar gleichmäßig stetig) ist, konvergiert auch die Folge $(f_n(x))_{n \in \N}$ aus $\mathcal{H}(x)$ bzgl.\ $d$ gegen $f(x)$, d.h.\ $f(x) \in \overline{\mathcal{H}(x)}$. {]}

Es sei $(f_n)_{n \in \N}$ eine Cauchyfolge in $\overline{\mathcal{H}}$.
Für jedes $x \in X$ ist $E_x$ gleichmäßig stetig, also ist offenbar auch $(f_n(x))_{n \in \N}$ eine Cauchyfolge in $E_x(\overline{\mathcal{H}}) \stackrel{(\ref{FA.1.36.9})}{\subset} \overline{\mathcal{H}(x)}$.
Nach Voraussetzung von (i) ist $\overline{\mathcal{H}(x)}$ kompakt, also vollständig, d.h.\ $(f_n(x))_{n \in \N}$ konvergiert in $\overline{\mathcal{H}(x)}$ (bzgl.\ $d$).
Wir setzen
$$ f(x) := \lim_{n \to \infty} f_n(x) \in \overline{\mathcal{H}(x)} \subset Y. $$
Wegen der Beliebigkeit von $x \in X$ ist somit eine Abbildung $f \: X \to Y$ definiert.

Wir behaupten
\begin{equation} \label{FA.1.36.10}
(f_n)_{n \in \N} \mbox{ konvergiert in } Y^X \mbox{ (bzgl.\ $d_{\infty}$).}
\end{equation}

{[} Zu (\ref{FA.1.36.10}): Sei $\varepsilon \in \R_+$.
$(f_n)_{n \in \N}$ ist eine Cauchyfolge, also existiert $n_0 \in \N$ derart, daß gilt
$$ \forall_{n,m \in \N, \, n,m \ge n_0} \, d_{\infty}(f_n,f_m) < \frac{\varepsilon}{2}, $$
d.h.\
$$ \forall_{n,m \in \N, \, n,m \ge n_0} \forall_{x \in X} \, d_{\infty}(f_n(x),\underbrace{f_m(x)}_{\stackrel{m \to \infty}{\longrightarrow} f(x)}) < \frac{\varepsilon}{2}, $$
also $\forall_{n \in \N, \, n \ge n_0} \forall_{x \in X} \, d(f_n(x),f(x)) \le \frac{\varepsilon}{2}$.
Somit gilt auch für jedes $n \in \N$ mit $n \ge n_0$: $d_{\infty}(f_n,f) \le \frac{\varepsilon}{2} < \varepsilon$. {]}

Da $(f_n)_{n \in \N}$ eine Folge in $\overline{\mathcal{H}} \subset \mathcal{C}(X,Y)$ ist, folgt aus \ref{FA.1.Isom.2}, daß auch $f \: X \to Y$ stetig ist.
Da $\overline{\mathcal{H}}$ außerdem abgeschlossen in $\mathcal{C}(X,Y)$ ist, ergibt sich aus \ref{FA.1.KA} schließlich $f \in \overline{\mathcal{H}}$, und wir haben (\ref{FA.1.36.2}) bewiesen.

,,(ii) $\Rightarrow$ (i)`` Wir zeigen zunächst, daß $\mathcal{H}$ gleichgradig stetig ist.
Seien hierzu $x \in X$ und $\varepsilon \in \R_+$ beliebig.
Wegen (ii) ist $\mathcal{H}$ nach \ref{FA.1.29} (iii) präkompakt, also existieren $f_1, \ldots, f_k \in \mathcal{H}$ derart, daß gilt
\begin{equation} \label{FA.1.36.11}
\mathcal{H} \subset \bigcup_{i=1}^k U^{d_{\infty}}_{\frac{\varepsilon}{3}}(f_i),
\end{equation}
und wir definieren
\begin{equation} \label{FA.1.36.12}
U := \bigcap_{i=1}^k \overline{f_i}^1 \left( U_{\frac{\varepsilon}{3}}(f_i(x)) \right) \in \U(x,X),
\end{equation}
beachte, daß $f_1, \ldots, f_k$ stetig sind.

Seien nun $f \in \mathcal{H}$ beliebig und $i \in \{1, \ldots, k\}$ mit 
\begin{equation} \label{FA.1.36.13}
f \in U^{d_{\infty}}_{\frac{\varepsilon}{3}}(f_i),
\end{equation}
welches wegen (\ref{FA.1.36.11}) existiert.

Dann gilt für jedes $\tilde{x} \in U$
$$ d(f(\tilde{x}),f(x)) \le \underbrace{d(f(\tilde{x}), f_i(\tilde{x}))}_{\stackrel{(\ref{FA.1.36.13})}{<} \frac{\varepsilon}{3}} + \underbrace{d(f_i(\tilde{x}), f_i(x))}_{\stackrel{(\ref{FA.1.36.12})}{<} \frac{\varepsilon}{3}} + \underbrace{d(f_i(x), f(x))}_{\stackrel{(\ref{FA.1.36.13})}{<} \frac{\varepsilon}{3}} < \varepsilon, $$
d.h.\ $\mathcal{H}$ ist gleichgradig stetig in $x$.

Zum Nachweis von (i) bleibt (\ref{FA.1.36.0}) zu zeigen.
Sei also $x \in X$.
Wie oben erwähnt, ist für jedes $x \in X$ die Evaluationsabbildung $E_x \: \mathcal{C}(X,Y) \to Y$ stetig.
Daher folgt aus der Kompaktheitstreue stetiger Abbildungen und (ii), daß $E_x(\overline{\mathcal{H}})$ kompakt, insbes.\ abgeschlossen (da $Y$ hausdorffsch als fastmetrischer Raum), ist.
Somit ist $\overline{E_x(\mathcal{H})} = \overline{\mathcal{H}(x)}$ eine abgeschlossene Teilmenge des Kompaktums $E_x(\overline{\mathcal{H})}$ und somit kompakt.

,,(ii) $\Rightarrow$ (iii)`` Sei $(f_n)_{n \in \N}$ eine Folge in $\mathcal{H} \subset \overline{\mathcal{H}}$.
Wegen (ii) ist $\overline{\mathcal{H}}$ kompakt, also nach \ref{FA.1.30.K} folgenkompakt in $\mathcal{C}(X,Y)$.
Somit besitzt $(f_n)_{n \in \N}$ eine Teilfolge, die bzgl.\ $d_{\infty}$ in $\mathcal{C}(X,Y)$ konvergiert, und das bedeutet genau (iii).

,,(iii) $\Rightarrow$ (ii)`` Wir zeigen, daß $\overline{\mathcal{H}}$ folgenkompakt, also kompakt, ist.
Sei $(f_n)_{n \in \N}$ eine Folge in $\overline{\mathcal{H}}$, d.h.\ eine Folge von Berührungspunkten von $\mathcal{H}$.
Dann gilt
$$ \forall_{n \in \N} \, U^{d_{\infty}}_{\frac{1}{n+1}}(f_n) \cap \mathcal{H} \ne \emptyset, $$
und wir wählen zu jedem $n \in \N$ ein $g_n \in U^{d_{\infty}}_{\frac{1}{n+1}}(f_n) \cap \mathcal{H}$.
Wegen (iii) besitzt $(g_n)_{n \in \N}$ eine Teilfolge $(g_{i_n})_{n \in \N}$, die gegen ein gewisses $f \in \mathcal{C}(X,Y)$ bzgl.\ $d_{\infty}$ konvergiert, und es gilt sogar $f \in \overline{\mathcal{H}}$, da $(g_{i_n})_{n \in \N}$ eine Folge in $\mathcal{H}$ ist.
Nun folgt für jedes $n \in \N$
$$ d_{\infty}(f_{i_n},f) \le \underbrace{d_{\infty}(f_{i_n},g_{i_n})}_{\le \frac{1}{i_n + 1}} + \underbrace{d_{\infty}(g_{i_n},f)}_{\stackrel{n \to \infty}{\longrightarrow} 0} \stackrel{n \to \infty}{\longrightarrow} 0, $$
d.h.\ $(f_{i_n})_{n \in \N}$ konvergiert in $\overline{\mathcal{H}}$. \q

\begin{Bem*} $\,$
\begin{itemize}
\item[1.)] Der Beweis zeigt, daß die Kompaktheit von $X$ nur für ,,(i) $\Rightarrow$ (ii)`` benötigt wird.
Verzichtet man bei ,,(ii) $\Rightarrow$ (i)`` und ,,(ii) $\Leftrightarrow$ (iii)`` auf die Voraussetzung der Kompaktheit von $X$, so ist $d_{\infty}$ allerdings auch im Falle eines metrischen Raumes $(X,d)$ i.a.\ keine Metrik auf $\mathcal{C}(X,Y)$.
\item[2.)] Wenn $(\mathcal{C}(X,Y), d_{\infty})$ ein metrischer Raum ist, so folgt aus (ii) die Be\-schränkt\-heit von $\mathcal{H}$ bzgl.\ $d_{\infty}$, vgl.\ \ref{FA.1.33}.
\item[3.)] Ist in $Y$ jede beschränkte Menge relativ kompakt\footnote{Dies ist z.B.\ erfüllt, wenn $(Y,d)$ ein kompakter metrischer Raum oder (vgl.\ Kapitel \ref{FAna2}) gleich $(\K^n, d_k)$, wobei $\K \in \{\R, \C\}, n \in \N_+$ sowie $k \in \{1,2,\infty\}$ sei, ist.}, und ist $\mathcal{H}$ bzgl.\ $d_{\infty}$ beschränkt, so gilt (\ref{FA.1.36.0}).
\item[4.)] In der Literatur wird bei der Formulierung des Satzes von \textsc{\textsc{Arzel\`{a}-Ascoli}} häufig zusätzlich gefordert, daß $X$ hausdorffsch ist.
Nach Meinung des Autors hat dies entweder historische Gründe, da z.B.\ in \cite{Hirz} oder \cite{Reck} ein kompakter Raum per definitionem auch hausdorffsch ist.
(Bei \textsc{N.\ Bourbaki} werden kompakte Räume in unserem Sinne \emph{quasikompakt} genannt.)
Oder es liegt daran, daß man den Satz auch allgemeiner für lokal-kompakte Räume $X$ zeigen kann, wenn man $\mathcal{C}(X,Y)$ mit der Topologie der kompakten Konvergenz versieht, worauf wir hier nicht eingehen, vgl.\ z.B.\ \cite[Theroem 47.1]{Munkres}. 
Da in der Literatur aber teilweise unterschiedliche Definitionen der lokalen Kompaktheit gegeben werden, die nur für Hausdorff-Räume äquivalent sind, setzt man die Hausdorff-Ei\-gen\-schaft vermutlich voraus, um Konfusionen zu vermeiden.
In \cite{Munkres} wird auch ein lokal-kompakter Hausdorff-Raum vorausgesetzt, obwohl der Beweis die Haus\-dorff-Ei\-gen\-schaft eigentlich nicht verwendet.\footnote{\textsc{Munkres} verwendet die Hausdorff-Eigenschaft, um die Stetigkeit der Evaluationsabbildung $X \times \mathcal{C}(X,Y) \to Y, \, (x,f) \mapsto f(x),$ zu erhalten. Der Beweis des Satzes erfordert aber nur die Stetigkeit für festes $x \in X$, welche auch allgemein gegeben ist.} 
(Wir wollen hier einen topologischen Raum $X$ \emph{lokal-kompakt}\index{Raum!topologischer!lokal-kompakter} nennen, wenn gilt $\forall_{x \in X} \exists_{U \in \U(x,X)} \, U \subset \subset X$.)
\end{itemize}
\end{Bem*}

\subsection*{Der Satz von \textsc{Baire}} \addcontentsline{toc}{subsection}{Der Satz von \textsc{Baire}}

\begin{Def}[(Nirgends) dichte, generische und magere Mengen] \label{FA.1.39}
Es seien $X$ ein topologischer Raum und $A$ eine Teilmenge von $X$.
\begin{itemize}
\item[(i)] $A$ heißt \emph{dicht (in $X$)}\index{Menge!dichte} $: \Longleftrightarrow \overline{A}= X$.
\item[(ii)] $A$ heißt \emph{nirgends dicht (in $X$)}\index{Menge!nirgends dichte} $: \Longleftrightarrow \overline{A}^{\circ}= \emptyset$.
\item[(iii)] $A$ heißt \emph{generisch (in $X$)}\index{Menge!generische} genau dann, wenn eine Folge $(U_n)_{n \in \N}$ offener sowie dichter Teilmengen von $X$ mit 
$$ \bigcap_{n \in \N} U_n \subset A $$
existiert.

\begin{Bem*}
Der endliche Schnitt von Teilmengen von $X$, die offen bzw.\ dicht sind, ist wieder offen bzw.\ dicht.
\end{Bem*}

\item[(iv)] $A$ heißt \emph{mager (in $X$)}\index{Menge!magere} genau dann, wenn eine Folge $(N_n)_{n \in \N}$ nirgends dichter Teilmengen von $X$ mit 
$$ A = \bigcup_{n \in \N} N_n $$
existiert.

\begin{Bem*}
Die mageren Mengen heißen bei \textsc{R.\ Baire} \emph{Mengen erster Kategorie}, die nicht-mageren heißen \emph{Mengen zweiter Kategorie}.
\end{Bem*}
\end{itemize}

\begin{Bsp*} 
Sei $\K \in \{\R,\C\}$. 
Wir versehen $\K$ mit der Standardmetrik $d_{\K}$.
\begin{itemize}
\item[1.)] Für jedes $q \in \Q$ ist $\R \setminus \{q\}$ offen und dicht in $\R$, also ist die dichte und nicht-offene Teilmenge
$$ \R \setminus \Q = \bigcap_{q \in \Q} ( \R \setminus \{q\} ) $$
von $\R$ generisch -- beachte, daß $\Q$ abzählbar ist.
\item[2.)] $\boxed{S^1} := \{ \sigma \in \C \, | \, |\sigma| = 1 \}$ ist nirgends dicht in $\C$.
\item[3.)] ${[}0,1{]}$ ist weder dicht noch nirgends dicht in $\R$.
\item[4.)] Das \emph{Cantorsche Diskontinuum}\index{Cantorsches Diskontinuum}\index{Menge!Cantor-} $C$ ist eine nirgends dichte überabzählbare Teilmenge von $\R$. 
($C$ entsteht aus $[0,1]$, indem man im ersten Schritt das in der Mitte liegende offene Intervall $]\frac{1}{3}, \frac{2}{3}[$ der Länge $\frac{1}{3}$ entfernt, im zweiten Schritt aus jedem verbleibenden Intervall das in der Mitte liegende offene Intervall der Länge $\frac{1}{9} = \left(\frac{1}{3}\right)^2$ entfernt, im $k$-ten Schritt ($k \in \N_+$)  aus jedem verbliebenen Intervall das in der Mitte liegende offene Intervall der Länge $\left(\frac{1}{3}\right)^k$ entfernt usw.)
\item[5.)] Jede einelementige Teilmenge von $\K$ ist nirgends dicht in $\K$, also sind abzählbare Teilmengen von $\K$ mager.
Insbesondere ist $\Q$ mager in $\R$.
\end{itemize}
\end{Bsp*}
\end{Def}

\begin{Satz} \label{FA.1.40}
Sei $X$ ein topologischer Raum.
Dann gilt:
\begin{itemize}
\item[(i)] $G \subset X$ ist genau dann generisch, wenn eine Folge $(G_n)_{n \in \N}$ von Teilmengen von $X$ derart, daß $(G_n)^{\circ}$ für jedes $n \in \N$ dicht in $X$ ist und $\bigcap_{n \in \N} G_n = G$ gilt, existiert.
\item[(ii)] $G \subset X$ ist genau dann generisch, wenn $X \setminus G$ mager ist.
\item[(iii)] Der abzählbare Schnitt generischer Teilmengen von $X$ ist wieder generisch.
\end{itemize}
\end{Satz}

\textit{Beweis.} Zu (i): ,,$\Rightarrow$`` Sei $G \subset X$ generisch, d.h.\ es existiert eine Folge $(U_n)_{n \in \N}$ offener sowie dichter Teilmengen von $X$ mit $\bigcap_{n \in \N} U_n \subset G$.
Definiere dann für jedes $n \in \N$
$$ G_n := U_n \cup G. $$
Dann gilt zum einen $(G_n)^{\circ} \supset (U_n)^{\circ} = U_n$, d.h.\ wegen der Dichtheit von $U_n$, daß auch $(G_n)^{\circ}$ dicht in $X$ ist, und zum anderen
$$ \bigcap_{n \in \N} G_n = \underbrace{\left( \bigcap_{n \in \N} U_n \right)}_{\subset G} \cup G = G. $$

,,$\Leftarrow$`` Existiere eine Folge $(G_n)_{n \in \N}$ gemäß der rechten Seite.
Dann wird für jedes $n \in \N$ durch $U_n := (G_n)^{\circ}$ eine offene sowie dichte Teilmenge von $X$ definiert, und es gilt $\bigcap_{n \in \N} U_n \subset G$.

Zu (ii): Sei $G \subset X$ generisch, d.h.\ genau nach (i), daß eine Folge $(G_n)_{n \in \N}$ von Teilmengen von $X$ mit
$$ \forall_{n \in \N} \, {(G_n)^{\circ}} = X ~~ \wedge ~~ \bigcap_{n \in \N} G_n = G $$
existiert.
Dies ist gleichbedeutend damit daß für $N_n := X \setminus G_n$, $n \in \N$, gilt
\begin{gather*}
\overline{N_n} = \overline{X \setminus G_n} = X \setminus (G_n)^{\circ}, \\
\overline{N_n}^{\circ} = (X \setminus (G_n)^{\circ})^{\circ} = X \setminus \overline{(G_n)^{\circ}} = X \setminus X = \emptyset, \\
X \setminus G = X \setminus \left( \bigcap_{n \in \N} G_n \right) = \bigcup_{n \in \N} N_n
\end{gather*}
gilt, d.h.\ genau, daß $X \setminus G$ mager ist.

Zu (iii): Sei $G_k$ für jedes $k \in \N$ eine generische Menge.
Dann existiert zu $k$ eine Folge $(U_{k,n})_{n \in \N}$ offener sowie dichter Teilmengen von $X$ mit 
$$ \bigcap_{n \in \N} U_{k,n} \subset G_k, $$
also gilt auch
$$ \bigcap_{k \in \N} \bigcap_{n \in \N} U_{k,n} \subset \bigcap_{k \in \N} G_k. $$
Da $\N \times \N$ abzählbar ist, folgt, daß auch $\bigcap_{k \in \N} G_k$ generisch ist. \q

\begin{Bem*}
,,Generisch`` ist nicht das Gegenteil von ,,mager``.
Betrachte z.B.\ $\N$ mit der Topologie der kofiniten Mengen, vgl.\ \ref{FA.1.22.Bsp}.
In diesem kompakten topologischen Raum ist $\N$ sowohl generisch als auch mager.
(Daß $\N = \bigcup_{n \in \N} \{n\}$ mager ist, liegt an $\forall_{n \in \N} \, \overline{\{n\}}^{\circ} = \emptyset$.)
Analog sieht man ein, daß die Teilmenge $\Q$ des nicht-vollständigen metrischen Raumes $\Q$ (mit der durch $d_{\R}$ induzieren Metrik) ebenfalls sowohl generisch als auch mager ist.
\end{Bem*}

\begin{Satz} \label{FA.1.41}
Sei $X$ ein topologischer Raum.
Dann sind die folgenden Aussagen paarweise äquivalent:
\begin{itemize}
\item[(i)] Jede generische Menge in $X$ ist dicht in $X$.
\item[(ii)] Jede magere Menge in $X$ hat einen leeren offenen Kern.
\item[(iii)] Ist $(A_n)_{n \in \N}$ eine Folge abgeschlossener Mengen von $X$ und hat $\bigcup_{n \in \N} A_n$ einen inneren Punkt, so existiert ein $n_0 \in \N$ mit $(A_{n_0})^{\circ} \ne \emptyset$.
\item[(iv)] Ist $(A_n)_{n \in \N}$ eine Folge abgeschlossener Mengen von $X$ mit $\forall_{n \in \N} \, (A_n)^{\circ} = \emptyset$, so hat auch $\bigcup_{n \in \N} A_n$ keinen inneren Punkt.
\item[(v)] Ist $(U_n)_{n \in \N}$ eine Folge offener sowie dichter Mengen in $X$, so ist auch $\bigcap_{n \in \N} U_n$ dicht in $X$.
\end{itemize}
\end{Satz}

\textit{Beweis.} ,,(i) $\Leftrightarrow$ (ii)`` folgt sofort aus \ref{FA.1.40} (ii), und ,,(iii) $\Leftrightarrow$ (iv)`` ist trivial.
Die Gültigkeit von ,,(iv) $\Leftrightarrow$ (v)`` erhält man, indem man für jedes $n \in \N$ setzt $U_n := X \setminus A_n$ bzw.\ $A_n := X \setminus U_n$.
,,(i) $\Rightarrow$ (v)`` ergibt sich aus \ref{FA.1.40} (iii), und die Aussage ,,(v) $\Rightarrow$ (i)`` ist wieder trivial. \q

\begin{Bem*}
Dichte Teilmengen können einen leeren offenen Kern besitzen.
Betrachte hierzu erneut $\N$ mit der Topologie der kofiniten Mengen.
In diesem Raum gilt $\overline{2 \N} = \N$ und $(2 \N)^{\circ} = \emptyset$.
(Letzteres liegt daran, daß es keine nicht-leeren Teilmenge von $2 \N$ gibt, deren Komplement in $\N$ endlich ist).
\end{Bem*}

\begin{Def} \label{FA.1.42}
Ein topologischer Raum $X$, der eine -- und damit alle -- der Aussagen (i) - (v) des letzten Satzes erfüllt, heißt ein \emph{Baire-Raum}\index{Raum!Baire-} oder \emph{Bairesch}.
\end{Def}

\begin{Kor} \label{FA.1.42.Kor}
Seien $X$ ein nicht-leerer Baire-Raum und $G$ eine generische Teilmenge von $X$.

Dann ist $G$ nicht mager in $X$.
\end{Kor}

\textit{Beweis.} Anderenfalls wäre nach \ref{FA.1.40} (ii), (iii)
$$ G \cap (X \setminus G) = \emptyset $$
generisch in $X$, also -- da $X$ Bairesch -- dicht in $X$, Widerspruch! \q

\begin{Bsp*}
Die in der Bemerkung zu Satz \ref{FA.1.40} angegebenen Räume sind nicht Bairesch.
\end{Bsp*}

Der folgende Hauptsatz -- der Satz von \textsc{Baire} -- rechtfertigt überhaupt erst die Einführung der obigen Begrifflichkeiten.
In einer reichhaltigen Klasse von Räumen -- nämlich der der Baire-Räume -- ist eine generische Teilmenge so etwas wie ,,groß``.
Entsprechend ist sie ,,klein``, wenn sie mager ist.
Daß die Klasse der Baire-Räume viele ,,schöne`` Räume enthält, besagt gerade der Satz von \textsc{Baire}.
Vor der Formulierung benötigen wir eine Definition.

\begin{Def}[Lokale Kompaktheit] \label{FA.1.43}
Ein topologischer Raum $X$ heißt \emph{lokal-kompakt}\index{Kompaktheit!lokale}\index{Raum!topologischer!lokal-kompakter} genau dann, wenn jeder Punkt in $X$ eine relativ kompakte Umgebung besitzt, d.h.\ genau $\forall_{x \in X} \exists_{U \in \U(x,X)} \, U \subset \subset X$.
\end{Def}

\begin{Bsp*} $\,$
\begin{itemize}
\item[1.)] Jeder kompakte topologische Raum ist lokal-kompakt.
\item[2.)] Jeder metrische Raum, in dem jede beschränkte Menge relativ kompakt ist, ist lokal-kompakt.
Aus \ref{FA.2.E.4} und \ref{FA.2.E.11} (s.u.) folgt daher, daß $(\K^n,d_k)$ für $\K \in \{\R,\C\}$, $n \in \N_+$ und $k \in \{1,2, \infty\}$ lokal-kom\-pakt ist.
\end{itemize}
\end{Bsp*}

\begin{HS}[Satz von \textsc{Baire}] \index{Satz!von \textsc{Baire}} \label{FA.1.44}
Jeder vollständige fastmetrische Raum und jeder lokal-kompakte Hausdorff-Raum ist Bairesch.
\end{HS}

Wir bereiten den Beweis des Satzes von \textsc{Baire} durch die folgenden Eregbnisse vor:

\begin{Satz}[Durchschnittssatz von \textsc{Cantor}] \index{Satz!Durchschnitts- von \textsc{Cantor}} \label{FA.1.45} $\,$

\noindent \textbf{Vor.:} Seien $X$ ein vollständiger fastmetrischer Raum und $(A_n)_{n \in \N}$ eine Folge nicht-leerer abgeschlossener Teilmengen von $X$ mit
\begin{equation} \label{FA.1.45.1}
\forall_{n \in \N} \, A_{n+1} \subset A_n 
\end{equation}
sowie
\begin{equation} \label{FA.1.45.2}
\lim_{n \to \infty} {\rm diam}(A_n) = 0. 
\end{equation}

\noindent \textbf{Beh.:} Es existiert $a \in X$ mit $\D \bigcap_{n \in \N} A_n = \{a\}.$
\end{Satz}

\textit{Beweis.} 1.) Wähle zu jedem $n \in \N$ ein $a_n \in A_n$.
Dann ist $(a_n)_{n \in \N}$ eine Cauchyfolge in $X$, denn zu $\varepsilon \in \R_+$ existiert nach (\ref{FA.1.45.2}) ein $n_0 \in \N$ mit ${\rm diam}(A_{n_0}) < \varepsilon$, und es gilt wegen (\ref{FA.1.45.1}) $a_n, a_m \in A_{n_0}$ für alle $n,m \in \N$ mit $n,m \ge n_0$, also $d(a_n,a_m) \le {\rm diam}(A_{n_0}) < \varepsilon$.
Da $X$ vollständig ist, konvergiert $(a_n)_{n \in \N}$ gegen ein gewisses $a \in X$ und wir behaupten: $a \in \bigcap_{n \in \N} A_n$.

Beweis hiervon: Für jedes $n \in \N$ gilt nach (\ref{FA.1.45.1})
$$ \forall_{m \in \N, \, m \ge n} \, \underbrace{a_m}_{\stackrel{m \to \infty}{\longrightarrow} a} \in A_n, $$
also folgt aus der Abgeschlossenheit von $A_n$ und \ref{FA.1.KA}: $a \in A_n$.

2.) Seien $a, b \in \bigcap_{n \in \N} A_n$, also gilt
$$ \forall_{n \in \N} \, d(a,b) \le {\rm diam}(A_n), $$
also nach (\ref{FA.1.45.2}) $d(a,b) = 0$, d.h.\ $a = b$. 

Mit 1.) und 2.) ist der Satz bewiesen. \q

\begin{Satz} \label{FA.1.46}
Es seien $X$ ein Hausdorff-Raum und $(A_n)_{n \in \N}$ eine Folge nicht-leerer kompakter Teilmengen von $X$ mit
\begin{equation} \label{FA.1.46.1}
\forall_{n \in \N} \, A_{n+1} \subset A_n.
\end{equation}

Dann gilt $\D \bigcap_{n \in \N} A_n \ne \emptyset$.
\end{Satz}

\textit{Beweis.} Da $X$ hausdorffsch ist, ist $(A_n)_{n \in \N}$ eine Folge nicht-leerer abgeschlossener Teilmengen von $X$.
Für jedes $n \in \N$ gilt nach (\ref{FA.1.46.1}) $A_n \subset A_0$ und
$$ A_0 \setminus A_n = \underbrace{\underbrace{(X \setminus A_n)}_{\text{offen in $X$}} \cap A_0}_{\text{offen in $A_0$}}, $$
d.h.\ $(A_n)_{n \in \N}$ eine Folge nicht-leerer in $A_0$ abgeschlossener Teilmengen von $A_0$.
$A_0$ ist kompakt, und aus (\ref{FA.1.46.1}) folgt offenbar, daß $(A_n)_{n \in \N}$ ein zentriertes System ist.
Daher ergibt sich die Behauptung aus \ref{FA.1.22.Durchschnitt} ,,(i) $\Rightarrow$ (ii)``. \q

\begin{Satz} \label{FA.1.47}
Seien $X$ ein lokal-kompakter Hausdorff-Raum und $K$ eine kompakte sowie $U$ eine offene Teilmenge von $X$ mit $K \subset U$.

Dann gilt $\forall_{x \in K} \exists_{V \in \U(x,X)} \, K \subset V \subset \overline{V} \subset U \, \wedge \, \overline{V} \mbox{ kompakt}$. 
\end{Satz}

\textit{Beweis.} Sei $x \in K$.
Wir behaupten zunächst, daß wir ohne Beschränkung der Allgemeinheit annehmen können, daß gilt $U \subset \subset X$.

{[} Nachweis hiervon: Aufgrund der lokalen Kompaktheit von $X$ existiert zu jedem $y \in K$ eine relativ kompakte Umgebung $U_y \in \U(y,X)$.
Wegen der Kompaktheit von $K$ besitzt die offene Überdeckung $\left( U_y \right)_{y \in K}$ von $K$ eine endliche Teilüberdeckung $\left( U_{y_i} \right)_{i \in \{1, \ldots, n\}}$, also gilt
$\widetilde{U} := U \cap \left( \bigcup_{i=1}^n U_{y_i} \right) \in \U(x,X)$, $K \subset \widetilde{U}$ und $\overline{\widetilde{U}} \subset \overline{U} \cap \overline{\left( \bigcup_{i=1}^n U_{y_i} \right)} \subset \overline{\bigcup_{i=1}^n U_{y_i}} = \bigcup_{i=1}^n \overline{U_{x_i}}$, wobei $\bigcup_{i=1}^n \overline{U_{x_i}}$ kompakt ist, d.h.\ nach \ref{FA.1.24} (i): $\widetilde{U} \subset \subset X$. {]}

Nun ist die abgeschlossene Menge $\partial U \stackrel{\ref{FA.1.8} (iii)}{=} \overline{U} \cap \overline{X \setminus U}$ als Teilmenge des Kompaktums $\overline{U}$ selbst kompakt, und wegen $\partial U \subset \overline{X \setminus U} = X \setminus U$ sowie $K \subset U$ gilt $\partial U \cap K = \emptyset$.
Da $X$ hausdorffsch ist, folgt aus \ref{FA.1.24.Kor} offenbar die Existenz offener Mengen $\widetilde{U},\widetilde{V}$ mit $\widetilde{U} \cap \widetilde{V} = \emptyset$, $\partial U \subset \widetilde{U}$ und $K \subset \widetilde{V}$, also für $V := \widetilde{V} \cap U \in \U(x,X)$ auch $K \subset V$ und $\widetilde{U} \cap V = \emptyset$.
Somit gilt $V \subset U \cap (X \setminus \widetilde{U}) \subset \overline{U \cap (X \setminus \widetilde{U})} \subset \overline{U} \cap \overline{X \setminus \widetilde{U}}$, d.h.\ $\overline{V} \subset \overline{U} \cap \overline{X \setminus \widetilde{U}}$.
Wegen $\underbrace{X \setminus U}_{= \overline{X \setminus U}} \cap ( \overline{U} \cap \overline{X \setminus \widetilde{U}} ) \stackrel{\ref{FA.1.8} (iii)}{=} \underbrace{\partial U}_{\subset \widetilde{U}} \cap \underbrace{\overline{X \setminus \widetilde{U}}}_{= X \setminus \widetilde{U}} = \emptyset$ gilt daher $\overline{V} \subset U$.
Außerdem ist $\overline{V}$ als abgeschlossene Teilmenge der kompakten Menge $\overline{U}$ kompakt. \q

\begin{Kor} \label{FA.1.48}
Sei $X$ ein lokal-kompakter Hausdorff-Raum.

Dann gilt $\forall_{x \in X} \forall_{U \in \U(x,X)} \exists_{V \in \U(x,X)} \, \overline{V} \subset U \, \wedge \, \overline{V} \mbox{ kompakt}$. \q
\end{Kor}

\begin{Bem} \label{FA.1.48.B}
Wir haben soeben bewiesen, daß eine offene Teilmenge eines lokal-kompakten Hausdorff-Raumes ebenfalls lokal-kompakt (und hausdorffsch) ist.
Man kann leicht zeigen, daß jeder lokal-kompakte Hausdorff-Raum $X$, der nicht bereits kompakt ist, eine offene dichte Teilmenge seiner \emph{(Alexandroff) Einpunkt-Kompaktifizierung} $\boxed{X^{\infty}} := X \cupdot \{\infty\}$, wobei $\infty$ ein ,,unendlich ferner Punkt``, der nicht zu $X$ gehört, sei, versehen mit der Topologie
$$ {\rm Top}(X^{\infty}) := {\rm Top}(X) \cup \{ (X \setminus K) \cup \{\infty\} \, | \, K \subset X \mbox{ kompakt} \}, $$
ist. 
$X^{\infty}$ ist ebenfalls ein Hausdorff-Raum, da $X \setminus K$ für kompaktes $K$ offen in $X$ ist, weil $X$ als hausdorffsch vorausgesetzt ist.
Daher folgt, daß die lokal-kompakten Hausdorff-Räume genau die offenen Teilmengen der kompakten Hausdorff-Räume sind.
\end{Bem}

\textit{Beweis des Satzes von \textsc{Baire}.} Seien $X$ ein (ohne Beschränkung der Allgemeinheit nicht-leerer)
\begin{itemize}
\item[(A)] fastmetrischer Raum bzw.\ ein
\item[(B)] lokal-kompakter Hausdorff-Raum
\end{itemize}
und $G \subset X$ generisch. 
Es existiert also eine Folge $(U_n)_{n \in \N}$ offener sowie dichter Teilmengen von $X$ mit $\bigcap_{n \in \N} U_n \subset G$. 

Zu zeigen ist $\overline{G} = X$, d.h.\ genau $\forall_{x \in X} \forall_{U \in \U(x,X)} \, U \cap G \ne \emptyset$ bzw.\ 
$$ \forall_{U \in {\rm Top}(X) \setminus \{ \emptyset \}} \, U \cap G \ne \emptyset. $$
Hierfür genügt es nachzuweisen, daß gilt
$$ \forall_{U \in {\rm Top}(X) \setminus \{ \emptyset \}} \, U \cap \left( \bigcap_{n \in \N} U_n \right) \ne \emptyset. $$

Sei daher $U$ eine beliebige nicht-leere offene Menge.
Wir konstruieren rekursiv eine Folge $(A_n)_{n \in \N}$ nicht-leerer Teilmengen von $U$ derart, daß für alle $n \in \N$ gilt $A_{n+1} \subset A_n$, $(A_n)^{\circ} \ne \emptyset$, $A_n \subset U_n$ und
\begin{itemize}
\item[(A)] $A_n$ ist abgeschlossen mit ${\rm diam}(A_n) \le \frac{1}{n+1}$ bzw.\
\item[(B)] $A_n$ ist kompakt.
\end{itemize}

$U_0$ ist dicht in $X$, also existiert $x_0 \in U \cap U_0$.
Wegen der Offenheit von $U, U_0$ gilt $U \cap U_0 \in \U(x_0,X)$.
Im Falle (A) leistet $A_0 := \{ x \in X \, | \, d(x,x_0) \le r \}$ für hinreichend kleines $r \in {]}0, \frac{1}{2}{]}$ das Gewünschte.
Im Falle (B) existiert nach \ref{FA.1.48} eine relativ kompakte Umgebung $V_0 \in \U(x_0,X)$ mit $\overline{V_0} \subset U \cap U_0$, und wir können $A_0 := \overline{V_0}$ setzen.

Sei $n \in \N$ und $A_n$ bereits konstruiert.
$U_{n+1}$ ist dicht in $X$, also existiert $x_{n+1} \in (A_n)^{\circ} \cap U_{n+1}$, d.h. $(A_n)^{\circ} \cap U_{n+1} \in \U(x_{n+1},X)$, da $U_{n+1}$ offen ist.
Im Falle (A) existiert $r \in {]}0,\frac{1}{2 (n+1)}{]}$ derart, daß $A_{n+1} := \{ x \in X \, | \, d(x,x_{n+1}) \le r \}$ die gewünschten Eigenschaften hat.
Im Falle (B) existiert wieder nach \ref{FA.1.48} eine relativ kompaktes $V_{n+1} \in \U(x_{n+1},X)$ mit $\overline{V_{n+1}} \subset (A_n)^{\circ} \cap U_{n+1}$, und wir können $A_{n+1} := \overline{V_{n+1}}$ setzen.

Nun folgt im Falle (A) aus dem Cantorschen Durchschnittssatz \ref{FA.1.45} und im Falle (B) aus \ref{FA.1.46} die Existenz von $a \in \bigcap_{n \in \N} A_n$.
Wegen $\forall_{n \in \N} \, A_n \subset U_n$ gilt dann auch $a \in \bigcap_{n \in \N} U_n$.
Da $U_n$ für jedes $n \in \N$ ferner eine Teilmenge von $U$ ist, ergibt sich $a \in U \cap \left( \bigcap_{n \in \N} U_n \right)$. \q

\begin{Bsp*} 
Wir betrachten $\R$ wieder mit der Standardmetrik $d_{\R}$.
Nach dem Satz von \textsc{Baire} ist dieser Raum Bairesch, also ist ${[}0,1{]}$ weder generisch noch mager in $\R$.
(Wir sehen also, daß auch in Baire-Räumen ,,generisch`` nicht das Gegenteil von ,,mager`` ist.) 

Des weiteren folgt aus \ref{FA.1.42.Kor}, daß die oben bereits als generisch erkannte Teilmenge $\R \setminus \Q$ von $\R$ nicht mager ist.
Nach \ref{FA.1.40} (ii) ist daher auch die magere Teilmenge $\Q$ von $\R$ nicht generisch.
\end{Bsp*}

Für topologische Anwendungen des Satzes von \textsc{Baire} sei auf \cite[S.\ 297 ff.]{Munkres} verwiesen.
Funktionalanalytische Folgerungen werden wir im nächsten Kapitel sehen.

\subsection*{Übungsaufgaben}

\begin{UA}
Beweise den Satz in \ref{FA.1.TT}.
\end{UA}

\begin{UA} \label{FA.UA.A}
Seien $X$ ein topologischer Raum und $Y$ eine offene (bzw.\ abgeschlossene) Teilmenge von $X$.

Zeige, daß jede offene (bzw.\ abgeschlossene) Teilmenge von $Y$ auch offen (bzw.\ abgeschlossen) in $X$ ist.
\end{UA}

\begin{UA}
Beweise Satz \ref{FA.1.TT.1}.
\end{UA}

\begin{UA}
Beweise den Satz in \ref{FA.1.5}.
\end{UA}

\begin{UA}
Beweise die in Beispiel \ref{FA.1.7.dsup} genannte Aussage.
\end{UA}

\begin{UA}
Beweise Satz \ref{FA.1.8.L}.
\end{UA}

\begin{UA}
Seien $X, Y$ topologische Räume.

Zeige, daß eine Abbildung $f \: X \to Y$ genau dann stetig ist, wenn für jede Teilmenge $A$ von $X$ gilt $f(\overline{A}) \subset \overline{f(A)}$, d.h.\ $\overline{A} \subset \overline{f}^1(\overline{f(A)})$.
\end{UA}

\begin{UA}
Beweise Satz \ref{FA.1.TT.3}.
\end{UA}

\begin{UA}
Es seien $X, Y$ topologische Räume, $I$ eine Menge und $A_i$ für alle $i \in I$ eine abgeschlossene Teilmenge von $X$ derart, daß gilt $\D \bigcup_{i \in I} A_i = X$.
Ferner sei $f_i \: A_i \to Y$ für jedes $i \in I$ eine stetige Abbildung mit
$$ \forall_{i,j \in I} \, f_i|_{A_i \cap A_j} = f_j|_{A_i \cap A_j}. $$

Zeige, daß genau eine stetige Abbildung $f \: X \to Y$ mit $\forall_{i \in I} \, f|_{A_i} = f_i$ existiert.
\end{UA}
Tip: Aufgabe \ref{FA.UA.A}.

\begin{UA}[Existenz- und Eindeutigkeitssatz von \textsc{Picard-Lindelöf}]\index{Satz!Existenz- und Eindeutigkeits- von \textsc{Picard-Lindelöf}}

Es seien $a,b \in \R$ mit $a<b$, $X$ ein endlich-dimensionaler normierter $\R$-Vektorraum und $f \in \mathcal{C}({[}a,b{]} \times X, \R)$ sowie $L \in \R_+$ derart, daß gilt
$$ \forall_{x \in {[}a,b{]}} \forall_{y_1, y_2 \in X} \, |f(x,y_1) - f(x,y_2)| \le L \, \| y_1 - y_2\|. $$

Zeige, daß zu jedem $y_0 \in X$ genau eine Lösung $y \in \mathcal{C}({[}a,b{]}, X)$ der \emph{Anfangswertaufgabe}
$$ \boxed{y' = f(x,y), ~~ y(a) = y_0,} $$
existiert.
\end{UA}
Tip: Durch
$$ \forall_{g,h \in \mathcal{C}({[}a,b{]}, X)} \, d(g,h) := \sup \{ \exp (-2 L (t-a)) \, \|g(t) - h(t)\| \, | \, t \in {[}a,b{]} \} $$
wird eine Metrik $d$ auf $\mathcal{C}({[}a,b{]}, X)$ definiert bzgl.\ derer $\mathcal{C}({[}a,b{]}, X)$ vollständig ist.
Für jedes $y \in \mathcal{C}({[}a,b{]}, X)$ gilt des weiteren $T_y \in \mathcal{C}({[}a,b{]}, X)$, wobei
$$  T_y := y_0 + \int_a^x f(t,y(t)) \, \d t, $$
sei, und
$$ T \: \mathcal{C}({[}a,b{]}, X) \longrightarrow \mathcal{C}({[}a,b{]}, X), ~~ y \longmapsto T_y, $$
ist kontrahierend bzgl.\ $d$ (mit ,,Kontrahierungsfaktor`` $\frac{1}{2}$).

\begin{UA}[Vervollständigung metrischer Räume nach \textsc{H.\ König}] \label{FA.UA.K} 
Es sei $(X,d)$ ein metrischer Raum.
Wir definieren eine Teilmenge
$$ Z := \{ f \in \R^X \, | \, \forall_{x,y \in X} \, f(x) - f(y) \le d(x,y) \le f(x) + f(y) \}. $$
von $\R^X$.

Zeige, daß für jedes $a \in X$
$$ f_a \: X \longrightarrow \R, ~~ x \longmapsto d(x,a), $$
ein Element von $Z$ ist, so daß wir eine wohldefinierte Abbildung
$$ j \: X \longrightarrow Z, ~~ a \longmapsto f_a, $$
erhalten.

Zeige weiterhin
\begin{gather*}
\forall_{f,g \in Z} \forall_{x,y \in X} \, f(x) - f(y) \le g(x) + g(y), \\
\forall_{f,g \in Z} \forall_{x,y \in X} \, |f(x) - g(x)| \le f(x) + g(y)
\end{gather*}
und folgere, daß $d_{\infty}|_{Z \times Z}$ eine Metrik auf $Z$ ist.

Beweise nun, daß $Z$ eine bzgl.\ $d_{\infty}$ abgeschlossene Teilmenge von $\R^X$ ist, und daher ist $(Z,d_{\infty}|_{Z \times Z})$ ein vollständiger metrischer Raum mit der Eigenschaft
$$ \forall_{a \in X} \forall_{f \in Z} \, d_{\infty}(f_a,f) = f(a). $$

Folgere hieraus, daß 
$$ j \: (X,d) \longrightarrow (\mathcal{H}, d_{\infty}|_{Z \times Z}) $$
eine isometrische Abbildung ist.
Daher (!) ist $j(X)$ eine dichte Teilmenge des vollständigen metrischen Raumes $(\overline{j(X)}, d_{\infty}|_{\overline{j(X)} \times \overline{j(X)}})$, wobei $\overline{j(X)}$ die abgeschlossene Hülle von $j(X)$ in $({\R}^X, d_{\infty})$ bezeichne.

Identifiziere nun $X \equiv j(X)$ und setze
$$ \widehat{X} := X \cup \left( \overline{j(X)} \setminus j(X) \right). $$
Definiere des weiteren für alle $x,y \in \widehat{X}$
$$ \hat{d}(x,y) = \left\{
\begin{array}{cl}
d(x,y), & \mbox{falls $x,y \in X$,} \\ y(x), & \mbox{falls $x \in X$ und $y \in \overline{j(X)} \setminus j(X)$,} \\ d_{\infty}(x,y), & \mbox{falls $x,y \in \overline{j(X)} \setminus j(X)$.}
\end{array} \right. $$

Beweise schließlich, daß $(\widehat{X},\hat{d})$ eine Vervollständigung von $(X,d)$ ist.
\end{UA}

\begin{UA} \label{FA.disj abgeschl und komp} 
Seien $X$ ein fastmetrischer Raum und $A, B$ nicht-leere disjunkte Teilmengen von $X$ derart, daß $A$ kompakt und $B$ abgeschlossen ist.

Zeige, daß dann gilt $d(A,B) := \inf \{ d(a,b) \, | \, a \in A \, \wedge \, b \in B \} > 0$.
\end{UA}

\begin{UA}
Wir wollen in dieser Aufgabe ein nicht-triviales Beispiel (d.h.\ eine echte Teilmenge eines fastmetrischen Raumes) vorstellen, die abgeschlossen aber nicht vollständig ist.
Seien hierzu $X := \R[x] \cap \mathcal{C}([0,1],\R) \subset (\R^{[0,1]}, d_{\infty})$ und $Y := \{ p \in X \, | \, p(0) \ge 1 \}$.

Der Leser überzeuge sich davon, daß $Y$ in $X$ abgeschlossen aber nicht vollständig ist.
Er überlege des weiteren, ob $X$ vollständig sein kann.
\end{UA}

Tip: Die Evaluationsabbildung ist stetig.
Außerdem ist $\sum_{k=0}^n \frac{x^k}{k!}|_{[0,1]}$ für jedes $n \in \N$ ein Element von $Y$, das für $n \to \infty$ punktweise gegen $\e^x|_{[0,1]}$ konvergiert, und es gilt z.B.\ der Satz von \textsc{Dini}.

\begin{UA}
Es seien $X$ ein topologischer Raum, $Y$ ein metrischer Raum und $\mathcal{H} \subset \mathcal{C}(X,Y)$.
Ferner sei jede in $Y$ beschränkte Menge relativ kompakt in $Y$.

Zeige, daß
$$ A := \{ x \in X \, | \, \mathcal{H}(x) \subset \subset Y \} $$
sowohl offen als auch abgeschlossen in $X$ ist.

Ist $X$ sogar zusammenhängend, d.h.\ per definitionem, daß $\emptyset$ und $X$ die einzigen Teilmengen von $X$ sind, die sowohl offen als auch abgechlossen sind, so folgt daher
$$ \left( \exists_{x \in X} \, \mathcal{H}(x) \subset \subset X \right) \Longleftrightarrow \left( \forall_{x \in X} \, \mathcal{H}(x) \subset \subset X \right). $$
\end{UA}
Tip: $\overline{A} \subset A^{\circ}$.

\begin{UA}[Zum Beweis des Satzes von \textsc{Montel} der \emph{Funktionentheorie}]
Es sei $G \subset \C$ ein Gebiet, d.h.\ per definitionem, daß $G \subset \C$ offen und zusammenhängend ist.
Ferner sei $C \in \R_+$.

Zeige, daß die Menge der holomorphen Funktionen $f \: G \to \C$ mit $|f| \le C$ gleichgradig stetig ist.
\end{UA}

\begin{Bem*}
Der Satz von \textsc{Montel}\index{Satz!von \textsc{Montel}} besagt, daß eine lokal beschränkte Funktionenfolge auf einem Gebiet definierter holomorpher Funktionen eine lokal gleichmäßig konvergente Teilfolge besitzt.
\end{Bem*}

\begin{UA}
Versehe $\Q$ mit der durch die Standardmetrik $d_{\R}$ induzierten Metrik.

Zeige, daß jede Teilmenge von $\Q$ generisch ist.
\end{UA}

\begin{UA}
Beweise, daß jede generische und jede nicht-leere offene Teilmenge eines Baire-Raumes als Teilraum auch wieder ein Bairescher Raum ist.
\end{UA}

\cleardoublepage
\section{Normierte Räume und Algebren} \label{FAna2}
Sei von nun an stets $\K \in \{\R,\C\}$.

\subsection*{Grundlagen} \addcontentsline{toc}{subsection}{Grundlagen}

\begin{Def}[Norm, normierter Vektorraum, normierte Algebra] \label{FA.2.1}
Sei $X$ ein $\K$-Vektorraum.
\begin{itemize}
\item[(i)] Eine Abbildung $\| \ldots \| \: \, X \to {[}0, \infty{[}$ heißt \emph{Norm auf $X$}\index{Norm} (es sei hier bereits auf \ref{FA.N.1} (i) unten hingewiesen), genau dann, wenn gilt
\begin{eqnarray*}
({\rm N}1) & \forall_{x \in X} \, \| x \| = 0 \Longleftrightarrow x = 0, & \mbox{(Positiv-Definitheit)},\\
({\rm N}2) & \forall_{\lambda \in \K} \forall_{x \in X} \, \| \lambda \, x \| = | \lambda | \, \| x \|,& \mbox{(Homogenität)},\\
({\rm N}3) & \forall_{x,y \in X} \, \| x + y \| \le \| x \| + \| y \|, & \mbox{(Dreiecksungleichung).}\index{Ungleichung!Dreiecks-}
\end{eqnarray*}

Gelten nur (N2) und (N3), so spricht man von einer \emph{Halbnorm}\index{Halb!-norm}\index{Norm!Halb-}.
\item[(ii)] Ist $\| \ldots \| \: X \to {[}0, \infty{[}$ eine Norm auf $X$, so heißt das Paar $(X, \| \ldots \|)$, für das wir auch kurz $X$ schreiben, ein \emph{normierter $\K$-Vektorraum}\index{Raum!normierter Vektor-}\index{Vektorraum!normierter}.
\item[(iii)] Sei $X$ sogar ein normierter $\K$-Vektorraum.
Besitzt $X$ zusätzlich die Struktur einer \emph{assoziativen $\K$-Algebra mit Einselement $e$},%
\footnote{\textbf{Definition} (Assoziative Algebra, Unteralgebra, Ideal)\textbf{.}\index{Algebra}
Sei $A$ ein $\K$-Vektorraum.
\begin{itemize}
\item[(i)] $A$ heißt \emph{assoziative $\K$-Algebra}, wenn eine $\K$-bilineare Abbildung, eine sog.\ \emph{Multiplikation}\index{Multiplikation},
$$ A \times A \longrightarrow A, ~~ (a,b) \longmapsto a \cdot b, $$
wobei wir für $a \cdot b$ auch oft $a \, b$ schreiben, mit
$$ \forall_{a,b,c \in A} \, a \, (b \, c) = (a \, b) \, c $$
existiert.
Bei unserer Notation bindet die Multiplikation stärker als die Addition.
\item[(ii)] $A$ wie in (i) heißt \emph{assoziative $\K$-Algebra mit Einselement}, wenn $e \in A$ mit
$$ \forall_{a \in A} \, e \, a = a \, e = a $$
existiert.
$e$ ist also eindeutig bestimmt, und es gilt $e = 0 \Rightarrow A = \{0\}$.
\item[(iii)] Ist $A$ eine assoziative $\K$-Algebra, so heißt ein $\K$-Untervektorraum $B$ von $A$ eine \emph{$\K$-Un\-ter\-al\-ge\-bra von $A$}, wenn gilt
$$ \forall_{b,\tilde{b} \in B} \, b \, \tilde{b} \in B. $$
Gilt darüber hinaus $A \, B \subset B$ und $B \, A \subset B$, wobei wir
$$ \boxed{A \, B} := \{ a \, b \, | a \in A \, \wedge \, b \in B \} ~~ \mbox{ bzw.\ } ~~ \boxed{B \, A} := \{ b \, a \, | b \in B \, \wedge \, a \in A \}, $$
setzen, so heißt $B$ ein \emph{$\K$-Ideal von $A$}.\index{Ideal}
Im Falle $B \subsetneqq A$ heißt $B$ \emph{echtes $\K$-Ideal}.

\begin{Bem*}
Besitzt $A$ ein Einselement, so ist ein $\K$-Ideal offenbar genau dann echt, wenn es das Einselement nicht enthält.
\end{Bem*}
\end{itemize}}
so heißt $X$ eine \emph{normierte $\K$-Algebra}\index{Algebra!normierte}, wenn gilt
\begin{eqnarray*}
({\rm N}4) && \|e\| = 1 \, \wedge \, \forall_{x,y \in X} \, \|x \, y \| \le \|x\| \, \|y\|. 
\end{eqnarray*}

Gilt anstelle von (N4) nur $\forall_{x,y \in X} \, \|x \, y \| \le \|x\| \, \|y\|$, so heißt $\| \dots \|$ \emph{submultiplikativ}.

\pagebreak

\begin{Bem*}
Eine normierte $\K$-Algebra besitzt für uns also per definitionem immer ein Einselement.
Dies ist aber keine wesentliche Einschränkung, da sich jede assoziative $\K$-Algebra $A$, die kein Einselement hat, durch \emph{Adjunktion der Eins} wie folgt zu einer assoziativen $\K$-Algebra $A_e$ mit Eins erweitern läßt, die $A$ als Unteralgebra enthält:
Man setzt $A_e := A \times \K$ und definiert für alle $\lambda \in \K$ sowie $(a,\mu), (b,\nu) \in A_e$
\begin{gather*}
\lambda \, (a,\mu) := (\lambda \, a, \lambda \, \mu), \\
(a,\mu) + (b,\nu) := (a+b,x+y), \\
(a,\mu) \, (b,\nu) := (a \, b + \mu \, b + \nu \, a, \mu \, \nu).
\end{gather*}
Diese Operationen machen $A_e$ zu einer assoziativen $\K$-Algebra mit Einselement $(0,1)$, die $A$ durch die Identifikation $A \equiv A \times \{0\}$ als Unteralgebra enthält.
Ist auf $A$ bereits eine submultiplikative Norm $\| \ldots \|$ gegeben, so wird diese durch
$$ \forall_{(a,\mu) \in A_e} \, \| (a,\mu) \|_e := \|a\| + |\mu| $$
zu einer Norm $\| \ldots \|_e$ auf $A_e$ erweitert, die $A_e$ zu einer normierten $\K$-Al\-ge\-bra macht. 
\end{Bem*}
\end{itemize}

\begin{Bsp*} $\,$
\begin{itemize}
\item[1.)] Seien $n \in \N_+$ sowie $p \in \{1,2\}$.
Dann werden durch
$$ \forall_{x = (x_1, \ldots, x_n) \in \K^n} \, \boxed{\|x\|_p} := \left( \sum_{i=1}^n |x_i|^p \right)^{\frac{1}{p}} $$
und
$$ \forall_{x = (x_1, \ldots, x_n) \in \K^n} \, \|x\|_{\infty} := \max \{ |x_i| ~ | ~ i \in \{1, \ldots, n\} \} $$
bekanntlich Normen auf $\K^n$ definiert, siehe ggf.\ Beispiel \ref{FA.5.6.Bsp} 1.) unten.

Im Falle $n=1$ erhält man also hier jeweils die übliche Betragsfunktion $| \ldots | \: \K \to {[}0, \infty{[}$, und $\K$ versehen mit der Multiplikation in $\K$ ist eine normierte $\K$-Algebra. Wir betrachten $\K$ stets als diese normierte $\K$-Algebra.
\item[2.)] Sind $X_1, \ldots, X_n$ normierte $\K$-Vektorräume, wobei sämtliche Normen mit $\| \ldots \|$ bezeichnet seien, so wird durch
$$ \forall_{(x_1, \ldots, x_n) \in \bigtimes_{i=1}^n X_i} \, \boxed{\|x\|_{\infty}} := \max \{ \|x_i\| ~ | ~ i \in \{1, \ldots, n\} \} $$
eine Norm, die sog.\ \emph{Maximumsnorm auf $\bigtimes_{i=1}^n X_i$}\index{Norm!Maximums- auf $\bigtimes_{i=1}^n X_i$}, definiert.
In dem Falle, daß alle $X_i$ gleich $\K$ mit der Betragsfunktion sind, ist die Definition mit 1.) verträglich. 
\item[3.)] Sei $M$ ein nicht-leerer kompakter topologischer Raum.
Auf $\mathcal{C}(M,\K)$ wird dann durch
\begin{equation} \label{FA.2.1.S}
\forall_{f \in \mathcal{C}(M,\K)} \, \boxed{\|f\|_{\infty}} := \sup \{ |f(p)| \, | \, p \in M \} = \max \{ |f(p)| \, | \, p \in M \}
\end{equation}

eine Norm definiert, die sog. \emph{Maximumsnorm auf $\mathcal{C}(M,\K)$}\index{Norm!Maximums- auf $\mathcal{C}(M,\K)$}.
(Daß das Maximum in (\ref{FA.2.1.S}) tatsächlich existiert, liegt an der Stetigkeit von $|f|$ und der Kompaktheit von $M$.)
Im Falle $n \in \N_+$ und $M = \{1, \ldots, n\} \equiv n$ ist auch diese Definition mit obigem verträglich.
\item[4.)] Seien $a,b \in \R$ mit $a<b$.
Auf dem $\R$-Vektorraum $\mathcal{C}([a,b], \R)$ wird durch
$$ \forall_{f \in \mathcal{C}([a,b], \R)} \, \boxed{\|f\|_1} := \int_a^b |f(t)| \, \d t $$
eine Norm definiert, wobei wir hier (zunächst) das Riemann-Integral verwenden, vgl.\ hierzu ggf.\ z.B.\ \cite[§ 6]{Henke} oder \cite[Kapitel 8]{ElAna}.
\end{itemize}
\end{Bsp*}
\end{Def}

\begin{Nr} \label{FA.2.2}
Sei $(X, \| \ldots \|)$ ein normierter $\K$-Vektorraum.
\begin{itemize}
\item[(i)] Durch
$$ \forall_{x,y \in X} \, d(x,y) := \| x - y \| $$
wird offenbar eine Metrik auf $X$ definiert.
Wir betrachten $X$ stets als metrischen Raum mit dieser induzierten Metrik.
Damit ist $X$ auch ein topologischer Raum, dessen Topologie ${\rm Top}(X,\|\ldots\|)$ wir die \emph{Normtopologie}\index{Norm!-topologie}\index{Topologie!Norm-} nennen.

Im Falle $X = \K$ ergibt unsere Konvention, daß wir $\K$ stets als $\K$-Algebra mit der Betragsfunktion betrachten also auch, daß wir $\K$ stets als metrischen Raum mit der Standardmetrik betrachten.

\begin{Lemma*}
Seien $X_1, \ldots, X_n$ normierte $\K$-Vektorräume, wobei sämtliche Normen mit $\| \ldots \|$ bezeichnet seien.

Eine Folge $(x_{1,k}, \ldots, x_{n,k})_{k \in \N}$ in $\bigtimes_{i=1}^n X_i$  konvergiert genau dann gegen $(x_1,\ldots,x_n) \in \bigtimes_{i=1}^n X_i$ bzgl.\ $\|\ldots \|_{\infty}$, wenn für jedes $i \in \{1, \ldots, n\}$ die Folge $(x_{i,k})_{k \in \N}$ in $X_i$ gegen $x_i$ konvergiert. 
\end{Lemma*}

\textit{Beweisskizze.} ,,$\Rightarrow$`` ist klar.
Für ,,$\Leftarrow$`` setzt man zu $\varepsilon \in \R_+$ 
$$ n_0 := \max \{ n_i \, | \, i \in \{1, \ldots, n\} \}, $$
wobei für alle $i \in \{1, \ldots,n\}$ gilt $\forall_{n \in \N, \, n \ge n_i} \| x_{i,n} - x_i \| < \varepsilon$. \q
\item[(ii)] Für alle $\varepsilon \in \R_+$ und alle $x \in X$ gilt
\begin{gather*}
U_{\varepsilon}(x) = \boxed{x + U_{\varepsilon}(0)} := \{ x + y \, | \, y \in U_{\varepsilon}(0) \}, \\
U_{\varepsilon}(x) = \boxed{\varepsilon \, \cdot \, U_{1}(x)} := \{ \varepsilon \, y \, | \, y \in U_{\varepsilon}(0) \}, \\
\overline{U_{\varepsilon}(x)} = \boxed{B_{\varepsilon}(x)} := \{ y \in X \, | \, \|x-y\| \le \varepsilon \}.
\end{gather*}

{[} ,,$\supset$`` in der letzten Gleichung liegt daran, daß zu jedem $y \in B_{\varepsilon}(x)$ die Folge $(\frac{n}{n+1} \, y)_{n \in \N}$ gegen $y$ konvergiert; ,,$\subset$`` ist trivial. {]}

$B_{\varepsilon}(x)$ heißt \emph{(abgeschlossene) Vollkugel vom Radius $\varepsilon$ um $x$}.
\item[(iii)] Für alle $x,y \in X$ gilt
$$ \|x\| - \|y\| \le | \| x \| - \| y \| | \le \| x - y \|, $$
denn die erste Ungleichung ist trivial und die zweite folgt aus der Dreiecksungleichung.

Folglich ist $\| \ldots \| \: X \to \R$ eine gleichmäßig stetige Funktion.
\item[(iv)] Wir versehen $\K \times X$ und $X \times X$ jeweils mit der Maximumsnorm.
Dann ist die Addition $X \times X \to X$ gleichmäßig stetig und die skalare Multiplikation $\K \times X \to X$ stetig.
Ist $X$ sogar eine normierte $\K$-Algebra, so ist auch die Multiplikation $X \times X \to X$ stetig.

\textit{Beweis.} 1.) Seien $x,\tilde{x},y,\tilde{y} \in X$.
Dann gilt
\begin{equation*}
\| (x + y) - (\tilde{x} + \tilde{y}) \| = \| (x - \tilde{x}) + (y - \tilde{y}) \| \le \| x - \tilde{x} \| + \| y - \tilde{y} \| \le 2 \, \max \{ \|x - \tilde{x}\|, \|y - \tilde{y}\| \}.
\end{equation*}

2.) Seien $\lambda \in \K$, $x \in V$, $(\lambda_n)_{n \in \N}$ eine Folge in $\K$ mit $\lim_{n \to \infty} \lambda_n = \lambda$ und $(x_n)_{n \in \N}$ eine Folge in $X$ mit $\lim_{n \to \infty} x_n = x$, beachte erneut das Lemma in (i). 
Dann gilt für jedes $n \in \N$
\begin{eqnarray*}
\| \lambda_n \, x_n - \lambda \, x \| & = & \| \lambda_n \, ( x_n - x ) + (\lambda_n - \lambda) \, x \| \\
& \le & \| \lambda_n \, ( x_n - x ) \| +  \| (\lambda_n - \lambda) \, x \| \\
& = & \underbrace{ \underbrace{| \lambda_n |}_{\stackrel{n \to \infty}{\longrightarrow} \lambda} \, \underbrace{\| x_n - x \|}_{\stackrel{n \to \infty}{\longrightarrow} 0} + \underbrace{| \lambda_n - \lambda |}_{\stackrel{n \to \infty}{\longrightarrow} 0} \, \| x \| }_{\mbox{$\stackrel{n \to \infty}{\longrightarrow} 0$}}.
\end{eqnarray*}

3.) Daß im Falle einer normierten $\K$-Algebra auch die Multiplikation stetig ist, sieht man analog zu 2.) ein. \q
\end{itemize}
\end{Nr}

\begin{Def}[Banachraum, -algebra] \label{FA.2.3}
Sei $X$ ein normierter $\K$-Vektorraum.
Ist $X$ als metrischer Raum vollständig, so nennen wir $X$ einen \emph{$\K$-Banach\-raum}\index{Raum!Banach-}\index{Banach!-raum}.
Falls $X$ sogar eine normierte $\K$-Algebra ist, heißt $X$ im Falle der Vollständigkeit eine \emph{$\K$-Banachalgebra}\index{Algebra!Banach-}\index{Banach!-algebra}. 

\begin{Bem*}
Ist $A$ eine $\K$-Algebra ohne Einselement, auf der eine submultiplikative Norm $\| \ldots \|$ gegeben ist, bzgl.\ derer $A$ als metrischer Raum vollständig ist, so ist $(A_e,\| \ldots \|_e$) wie in der Bemerkung zu \ref{FA.2.1} (iii) eine $\K$-Banach\-algebra. 
\end{Bem*}

\begin{Bsp*} $\,$
\begin{itemize}
\item[1.)] Sei $M$ nicht-leerer ein kompakter topologischer Raum.
$(\mathcal{C}(M,\K), \| \ldots \|_{\infty})$ ist dann eine $\K$-Banachalgebra, denn $\| \ldots \|_{\infty}$ induziert auf $\mathcal{C}(M,\K)$ offenbar die Supremumsfastmetrik $d_{\infty}|_{\mathcal{C}(M,\K) \times \mathcal{C}(M,\K)}$, die in diesem Falle natürlich eine Metrik ist, und $\K$ ist ein vollständiger metrischer Raum, also ist nach \ref{FA.1.15} auch $(\mathcal{C}(M,X),d_{\infty}|_{\mathcal{C}(M,X) \times \mathcal{C}(M,X)})$ vollständig.
\item[2.)] $(\mathcal{C}([-1,1],\R), \| \ldots \|_1)$ ist kein $\R$-Banachraum.

\textit{Beweis.} Betrachte für jedes $n \in \N$, die stetige Funktion $f_n \: [-1,1] \to \R$, die durch
$$ \forall_{t \in \R} \, f_n(t) := \left\{ \begin{array}{cc} 0, & -1 \le t < 0, \\ (n+1) t, & 0 \le t \le \frac{1}{n+1}, \\ 1, & \frac{1}{n+1} < t \le 1 \end{array} \right. $$
definiert ist.
Die Funktionenfolge $(f_n)_{n \in \N}$ konvergiert punktweise gegen $f \: [-1,1] \to \R$, gegeben durch
$$ \forall_{t \in \R} \, f(t) := \left\{ \begin{array}{cc} 0, & -1 \le t \le 0, \\ 1, & 0 < t \le 1. \end{array} \right. $$
Es gilt für jedes $n \in \N$
$$ \|f_n - f\|_1 =  \int_0^{\frac{1}{n+1}} |(n+1)t - 1| \, \d t = \frac{1}{n+1} - \frac{1}{2 (n+1)} = \frac{1}{2 (n+1)} \stackrel{n \to \infty}{\longrightarrow} 0,$$
$f$ ist aber nicht stetig.
\end{itemize}
\end{Bsp*}
\end{Def}

\begin{Def}[(Absolute) Konvergenz von Reihen]\index{Konvergenz!von Reihen} \index{Reihe!(absolut) konvergente} \label{FA.2.Reihe.0}
Seien $X$ ein normierter $\K$-Vektor\-raum und $(x_k)_{k \in \N}$ eine Folge in $X$.
\begin{itemize}
\item[(i)] Die \emph{Reihe} $\D \boxed{\sum_{k=0}^{\infty} x_k}$ ist per definitionem gleich der \emph{Folge der Partialsummen $(\sum_{k=0}^n x_k)_{n \in \N}$.}
Die Folge $(x_k)_{k \in \N}$ heißt \emph{Gliederfolge der Reihe $\sum_{k=0}^{\infty} x_k$}.
\item[(ii)] Im Falle der Konvergenz der Reihe $\sum_{k=0}^{\infty} x_k = (\sum_{k=0}^n x_k)_{n \in \N}$ in $X$, so bezeichnet man den Grenzwert $\lim_{n \to \infty} \sum_{k=0}^n x_k$ ebenfalls mit $\D \boxed{\sum_{k=0}^{\infty} x_k}$, also mit demselben Symbol wie die Reihe.

\begin{Lemma*}
Konvergiert $\sum_{k=0}^{\infty} x^k$ in $X$, so ist $(x_k)_{k \in \N}$ eine Nullfolge (in dem $\K$-Vektorraum $X$).
\end{Lemma*}

\textit{Beweis.} Für jedes $n \in \N$ sei $y_n := \sum_{k=0}^n x_k$.
Dann folgt aus der Voraussetzung $\lim_{n \to \infty} \underbrace{y_{n+1} - y_n}_{= x_{n+1}} = 0$, also auch $\lim_{n \to \infty} x_n = 0$. \q 

\item[(iii)] Im Falle der Konvergenz der Reihe $\sum_{k=0}^{\infty} \|x_k\|$ in $\R$, heißt die Reihe $\sum_{k=0}^{\infty} x_k$ \emph{absolut konvergent}.
\end{itemize}
\end{Def}

\begin{Satz} \label{FA.2.Reihe.1} $\,$

\noindent \textbf{Vor.:} Sei $X$ ein normierter $\K$-Vektorraum.

\noindent \textbf{Beh.:} $X$ ist genau dann ein $\K$-Banachraum, wenn jede absolut konvergente Reihe in $X$ konvergiert.
\end{Satz}

\textit{Beweis.} ,,$\Rightarrow$`` Seien $\sum_{k=0}^{\infty} x_k$ eine absolut konvergente Reihe in $X$ und $\varepsilon \in \R_+$.
Aus der Konvergenz von $\sum_{k=0}^{\infty} \|x_k\|$ in $\R$ folgt, daß $\sum_{k=0}^{\infty} \|x_k\|$ eine Cauchyfolge in $\R$ ist, d.h.\ es existiert $n_0 \in \N$ derart, daß für alle $n,m \in \N$ mit $m > n \ge n_0$ gilt
$$ \sum_{k=n+1}^m \|x_k\| = \sum_{k=0}^m \|x_k\| - \sum_{k=0}^n \|x_k\| < \varepsilon, $$
also auch
$$ \left\| \sum_{k=0}^m x_k - \sum_{k=0}^n x_k \right\| = \left\| \sum_{k=n+1}^m x_k \right\| \le \sum_{k=n+1}^m \|x_k\| < \varepsilon. $$
Daher ist $(\sum_{k=0}^n x_k)_{n \in \N}$ eine Cauchyfolge in dem $\K$-Banachrraum $X$, also konvergent. 

,,$\Leftarrow$`` Sei $(y_k)_{k \in \N}$ eine Cauchyfolge in $X$.
Dann existiert zu jedem $k \in \N$ eine $i_k \in \N$ mit
$$ \forall_{i,j \in \N, \, i,j \ge i_k} \, \| y_i - y_j \| < \frac{1}{2^k}, $$
und $(y_{i_k})_{k \in \N}$ ist eine Folge in $X$ mit
$$ \forall_{k \in \N} \, \| y_{i_k} - y_{i_{k+1}} \| < \frac{1}{2^k}. $$

Für $k \in \N$ sei $x_k := y_{i_{k+1}} - y_{i_k}$, also gilt $\sum_{k=0}^{\infty} \| y_k \| \le 2$.
Nach Voraussetzung der rechten Seite existiert dann $x \in X$ mit
$$ \lim_{n \to \infty} \underbrace{\sum_{k=0}^n x_k}_{= y_{i_{n+1}} - y_{i_0}} = x, $$
und es folgt $\lim_{n \to \infty} y_{i_n} = x + y_{i_0}$.

Wir zeigen, daß dann auch $\lim_{k \to \infty} y_k = x + y_{i_0}$ gilt:
Sei nämlich $\varepsilon \in \R_+$.
Da $(y_k)_{k \in \N}$ eine Cauchyfolge ist, existiert $k_0 \in \N$ mit 
$$\forall_{k,l \in \N, \, k,l \ge k_0} \, \| y_k - y_l \| < \frac{\varepsilon}{2},$$
und wegen $\lim_{k \to \infty} y_{i_k} = x + y_{i_0}$ existiert $k_1 \in \N$ mit $k_1 \ge k_0$ und
$$\forall_{k \in \N, \, k \ge k_1} \, \| y_{i_k} - ( x + y_{i_0} ) \| < \frac{\varepsilon}{2}.$$
Daher folgt für für jedes $k \in \N$ mit $k \ge k_1 (\ge k_0)$
$$ \| y_k - (x + y_{i_0}) \| \le \| y_k - y_{i_{k_1}} \| + \| y_{i_{k_1}} - (x + y_{i_0} ) \| < \varepsilon.$$
\q

\begin{Satz}[Neumannsche Reihen] \label{FA.2.Reihe.2} \index{Reihe!Neumannsche}
Es seien $A$ eine $\K$-Banachalgebra mit Einselement $e$ und $x \in A$. Wie üblich setzen wir $x^0 := e$.
Dann gilt:
\begin{itemize}
\item[(i)] Die Folge $(\sqrt[k]{\|x^k\|})_{k \in \N_+}$ konvergiert in $\R$.
$$ r_{\sigma}(x) := \lim_{k \to \infty} \sqrt[k]{\|x^k\|} $$
heißt \emph{Spektralradius an der Stelle $x$}.
\item[(ii)] $r_{\sigma}(x) \le \|x\|$.
\item[(iii)] Konvergiert die \emph{Neumannsche Reihe $\sum_{k=0}^{\infty} x^k$ (an der Stelle $x$)} in $A$, so besitzt $e-x$ ein inverses Element in $A$.%
\footnote{\textbf{Definition.} In jeder assoziativen $\K$-Algebra $A$ mit Einselement $e$ besitzt $a \in A$ per definitionem ein \emph{inverses Element $b \in A$}\index{Algebra!inverses Element in einer --}, wenn gilt $a \, b = b \, a = e$.
$b$ mit dieser Eigenschaft ist dann eindeutig bestimmt und wird mit $\boxed{a^{-1}}$ bezeichnet. 
In diesem Falle heißt $a$ \emph{invertierbar}.}%
Es gilt dann $(e-x)^{-1} = \sum_{k=0}^{\infty} x^k$.
\item[(iv)] $r_{\sigma}(x) < 1 \Longrightarrow \sum_{k=0}^{\infty} x^k \mbox{ konvergiert in $A$}$.
\item[(v)] $\|x\| < 1 \Longrightarrow \forall_{n \in \N} \, \|(e-x)^{-1} - \sum_{k=0}^n x^k\| \le \frac{\|x\|^{n+1}}{1 - \|x\|}.$
\end{itemize}
\end{Satz}

Wir bereiten den Beweis des Satzes durch das folgende Lemma vor:

\begin{Lemma} \label{FA.2.Reihe.L}
Es sei $(a_k)_{k \in \N}$ eine Folge nicht-negativer reeller Zahlen derart, daß gilt
\begin{equation} \label{FA.2.Reihe.L.1}
\forall_ {k,l \in \N} \, a_{k+l} \le a_k \, a_l.
\end{equation}

Dann konvergiert $(\sqrt[k]{a_k})_{k \in \N_+}$ und es gilt
$$ \lim_{k \to \infty} \sqrt[k]{a_k} = \inf \{ \sqrt[k]{a_k} \, | k \in \N_+ \}. $$
\end{Lemma}

\textit{Beweis.} Seien $\lambda := \inf \{ \sqrt[k]{a_k} \, | k \in \N_+ \} \in {[}0, \infty{[}$ und $\varepsilon \in \R_+$.
Dann ist $\lambda + \varepsilon$ keine untere Schranke von $\{ \sqrt[k]{a_k} \, | k \in \N_+ \}$, also können wir $l \in \N_+$ mit $\sqrt[l]{a_l} < \lambda + \varepsilon$ fixieren, d.h.\
\begin{equation} \label{FA.2.Reihe.L.2}
a_l < (\lambda + \varepsilon)^l.
\end{equation}

Sei $k \in \N_+$.
Dann existieren $p_k, q_k \in \N$ mit
\begin{gather}
k = p_k \, l + q_k, \label{FA.2.Reihe.L.3} \\
q_k \in \underbrace{\{0, \ldots , l-1 \}}_{= \{0\}, \text{ falls } l=1} , \label{FA.2.Reihe.L.4}
\end{gather}
also folgt mit 
\begin{equation} \label{FA.2.Reihe.L.5}
C := \max \{ a_0, \ldots, a_{l-1} \} \in {[}0, \infty {[}, 
\end{equation}
daß gilt
\begin{equation} \label{FA.2.Reihe.L.6}
a_k = a_{p_k \, l + q_k} \stackrel{(\ref{FA.2.Reihe.L.1})}{\le} (a_l)^{p_k} \, \overbrace{a_{q_k}}^{\text{evtl.\ } = 0} \stackrel{(\ref{FA.2.Reihe.L.2})}{\le} (\lambda + \varepsilon)^{l \, p_k} \, a_{q_k} \stackrel{(\ref{FA.2.Reihe.L.4}), (\ref{FA.2.Reihe.L.5})}{\le} (\lambda + \varepsilon)^{l \, p_k} \, C. 
\end{equation}

Für alle $k \in \N_+$ mit $k \ge l$ gilt $p_k \ne 0$, also wegen $l \stackrel{(\ref{FA.2.Reihe.L.3})}{=} \frac{k - q_k}{p_k}$ und (\ref{FA.2.Reihe.L.6})
$$ a_k \le (\lambda + \varepsilon)^{k - q_k} \, C, $$
d.h.\
\begin{equation} \label{FA.2.Reihe.L.7}
\forall_{k \in \N, \, k \ge l} \, \lambda \le \sqrt[k]{a_k} \le (\lambda + \varepsilon)^{1 - \frac{q_k}{k}} \, C^{\frac{1}{k}} = (\lambda + \varepsilon) \frac{C^{\frac{1}{k}}}{(\lambda + \varepsilon)^{\frac{q_k}{k}}}.
\end{equation}

Im Falle $C=0$ folgt aus $0 \le \lambda$ und (\ref{FA.2.Reihe.L.7}) sofort $\lambda = \lim_{k \to \infty} \sqrt[k]{a_k}$.

Sei daher im folgenden $C>0$.
Wir begründen unten, daß gilt
\begin{equation} \label{FA.2.Reihe.L.8}
\lim_{k \to \infty} \frac{C^{\frac{1}{k}}}{(\lambda + \varepsilon)^{\frac{q_k}{k}}} = 1.
\end{equation}
Dann folgt die Existenz von $k_0 \in \N$ mit
$$ \forall_{k \in \N, \, k \ge k_0} \, \frac{C^{\frac{1}{k}}}{(\lambda + \varepsilon)^{\frac{q_k}{k}}} < \underbrace{\frac{\lambda + 2 \varepsilon}{\lambda + \varepsilon}}_{> 1}, $$
also aus (\ref{FA.2.Reihe.L.7}) für jedes $k \in \N$ mit $k \ge l$ und $k \ge k_0$
$$ \lambda \le \sqrt[k]{a_k} < \lambda + 2 \varepsilon, $$
also aus der Beliebigkeit von $\varepsilon \in \R_+$, daß gilt $\lambda = \lim_{k \to \infty} \sqrt[k]{a_k}$.
\pagebreak

Zu zeigen bleibt (\ref{FA.2.Reihe.L.8}):
Wegen $C>0$ gilt $\lim_{k \to \infty} \sqrt[k]{C} = 1$.
Des weiteren folgt aus (\ref{FA.2.Reihe.L.4}) $\lim_{k \to \infty} \frac{q_k}{k} = 0$, also $\lim_{k \to \infty} (\lambda + \varepsilon)^{\frac{q_k}{k}} = (\lambda + \varepsilon)^0 = 1$,
insgesamt somit (\ref{FA.2.Reihe.L.8}). \q
\A
\textit{Beweis des Satzes.} (i) und (ii) folgen sofort aus dem Lemma, angewandt auf $a_k := \|x^k\|$ für alle $k \in \N$. Beachte bei (i), daß $X$ eine normierte $\K$-Algebra ist, und bei (ii), daß $\|x\| \in \{ \sqrt[k]{\|x_k\|} \, | \, k \in \N_+ \}$ gilt.

Zu (iii): Konvergiere also $\sum_{k=0}^{\infty} x^k$.
Dann gilt für jedes $n \in \N$ 
$$ (e - x) \sum_{k=0}^n x^k = \sum_{k=0}^n x^k - \sum_{k=1}^{n+1} x^k = e - x^{n+1} \stackrel{n \to \infty}{\longrightarrow} e $$ 
und analog
$$ \left( \sum_{k=0}^n x^k \right) (e - x) \stackrel{n \to \infty}{\longrightarrow} e, $$
beachte das Lemma in \ref{FA.2.Reihe.0} (ii).
Hieraus folgt (iii).

Zu (iv): Aus dem Wurzelkriterium der Analysis I folgt mittels der Voraussetzung von (iv) und
$$ \lim_{k \to \infty} \sqrt[k]{\|x^k\|} = r_{\sigma}(x) < 1 $$
sowie \ref{FA.2.Reihe.1} die Behauptung, da $X$ eine $\K$-Banachalgebra ist.

Zu (v): Aus $\|x\| < 1$ ergibt sich via (iv), (ii) und (iii):
$(e-x)^{-1} = \sum_{k=0}^{\infty} x^k$.
Es gilt außerdem für jedes $n \in \N$
\begin{eqnarray*}
\left\| \sum_{k=0}^{\infty} x^k - \sum_{k=0}^n x^k \right\| & = & \left\| \sum_{k=n+1}^{\infty} x^k \right\| \le \sum_{k=n+1}^{\infty} \|x\|^k = \sum_{k=0}^{\infty} \|x\|^k - \sum_{k=0}^n \|x\|^k
\\ & = & \frac{1}{1-\|x\|} - \frac{1-\|x\|^{n+1}}{1-\|x\|} = \frac{\|x\|^{n+1}}{1-\|x\|},
\end{eqnarray*}
also folgt (v). \q

\begin{Satz} \label{FA.2.4}
Seien $X$ ein normierter $\K$-Vektorraum (bzw.\ eine normierte $\K$-Al\-ge\-bra mit Einselement $e$) und $Y \subset X$ ein $\K$-Untervektorraum (bzw.\ eine $\K$-Unteralgebra).
Dann gilt:
\begin{itemize}
\item[(i)] $Y$ ist in kanonischer Weise ebenfalls ein normierter $\K$-Vektorraum (bzw.\ eine normierte $\K$-Algebra, falls $e \in Y$).
\item[(ii)] $\overline{Y} \subset X$ ist ebenfalls ein $\K$-Untervektorraum (bzw.\ eine $\K$-Unter\-algebra; und sogar ein echtes $\K$-Ideal, falls $Y$ ein solches ist).
\end{itemize}
\end{Satz}

\textit{Beweis.} (i) ist trivial.

Zu (ii): Seien $\lambda \in \K$ und $y_1, y_2 \in \overline{Y}$.
Es existieren also Folgen $(y_{i,n})_{n \in \N}$ in $Y$ mit $\lim_{n \to \infty} y_{i,n} = y_i$ für $i \in \{1,2\}$.
Dann gilt für jedes $n \in \N$ sowohl $\lambda \, y_{1,n} \in Y$ also auch $y_{1,n} + y_{2,n} \in Y$, also folgt mit \ref{FA.2.2} (iv)
$$ \lambda \, y_1 = \lim_{n \to \infty} \underbrace{\lambda \, y_{1,n}}_{\in Y} ~~ \mbox{ und } ~~ y_1 + y_2 = \lim_{n \to \infty} \underbrace{(y_{1,n} + y_{2,n})}_{\in Y}, $$
d.h.\ $\lambda \, y_1 \in \overline{Y}$ und $y_1 + y_2 \in \overline{Y}$.
\pagebreak

Im Falle einer normierten $\K$-Algebra gilt zusätzlich $y_{1,n} \cdot y_{2,n} \in Y$ für jedes $n \in \N$,
also wiederum aus Stetigkeitsgründen $y_1 \cdot y_2 = \lim_{n \to \infty} (y_{1,n} \cdot y_{2,n}) \in \overline{Y}$. 
Im folgenden sei $Y$ ein $\K$-Ideal. Man argumentiert analog, daß für jedes $x \in X$ gilt $x \cdot y_2 \in \overline{Y}$, $y_1 \cdot x \in \overline{Y}$.
Zu zeigen bleibt daher, daß aus $e \notin Y$ auch $e \notin \overline{Y}$ folgt:
Gälte $e \in \overline{Y}$, so existierte $y \in Y$ mit $\|e - y\| < 1$, d.h.\ nach \ref{FA.2.Reihe.2} (v), daß $y = e - (e - y)$ invertierbar wäre, Widerspruch!%
\footnote{\textbf{Lemma.}\label{FA.2.4.L} Seien $A$ eine assoziative $\K$-Algebra mit Einselement $e$ und $\mathfrak{b}$ ein $\K$-Ideal in $A$.

Dann ist $\mathfrak{b}$ genau dann ein echtes $\K$-Ideal, wenn $\mathfrak{b}$ keine invertierbaren Elemente enthält.

\textit{Beweis.} ,,$\Leftarrow$`` ist trivial.

,,$\Rightarrow$`` Angenommen, es existieren $b \in \mathfrak{b}$ und $a \in A$ mit $b \, a = a \, b = e$.
Dann folgt sofort, da $\mathfrak{b}$ ein $\K$-Ideal ist, daß gilt $e \in \mathfrak{b}$, im Widerspruch dazu, daß $\mathfrak{b}$ ein echtes $\K$-Ideal ist. \q}
\q

\begin{Satz} \label{FA.2.5}
Seien $X$ ein normierter $\K$-Vektorraum und $Y$ ein abgeschlossener $\K$-Untervektorraum von $X$.
Dann gilt:
\begin{itemize}
\item[(i)] Durch
$$ \forall_{\xi \in X/Y} \, \boxed{\| \xi \|_{\sim}} := \inf \{ \|x\| \in \R \, | \, x \in \xi \} = d(\{0\},Y)$$
wird eine Norm auf $X/Y$ definiert.%
\footnote{$X/Y$ ist (auch für nicht-abgeschlossenes $Y$) per definitionem der $\K$-Vektorraum der Äquivalenzklassen bzgl.\ der Äquivalenzrelation $\sim$, welche durch
$$ \forall_{x,\tilde{x} \in X} \, x \sim \tilde{x} : \Longleftrightarrow x - \tilde{x} \in Y $$
gegeben ist.
Es gilt dann für jedes $x \in X$: $[x]_{\sim} = \boxed{x + Y} := \{ x + y \, | \, y \in Y \}$, d.h.\ per definitionem, daß $[x]_{\sim}$ die Menge aller zu $Y$ \emph{parallelen} Mengen, die $x$ enthalten, ist.

In einem beliebigen metrischen Raum $X$ definiert man für je zwei nicht-leere Teilmengen $Y_1,Y_2$ von $X$ den \emph{Abstand von $Y_1$ zu $Y_2$}\index{Abstand zweier Teilmengen} als $\boxed{d(Y_1, Y_2)} := \inf \{ d(y_1,y_2) \, | \, y_1 \in Y_1 \, \wedge \, y_2 \in Y_2 \}$.}

Wir betrachten $X/Y$ im folgenden stets mit dieser induzierten Norm, für die wir häufig auch wieder $\| \ldots \|$ anstelle von $\| \ldots \|_{\sim}$ schreiben.
\item[(ii)] Der kanonische Homomorphismus $\pi \: X \to X/Y$ ist gleichmäßig stetig und eine \emph{offene Abbildung}.\index{Abbildung!offene}
Letzteres heißt per definitionem, daß $\pi$ offene Mengen auf offene Mengen abbildet.
\item[(iii)] Ist $X$ ein $\K$-Banachraum, so ist auch $X/Y$ ein solcher.
\item[(iv)] Sind zusätzlich $X$ eine normierte $\K$-Algebra und $Y$ ein echtes $\K$-Ideal von $X$, so ist $X/Y$ eine normierte $\K$-Algebra.%
\footnote{Sind $A$ eine assoziative $\K$-Algebra und $\mathfrak{b}$ ein $\K$-Ideal von $A$, so ist eine kanonische Multiplikation in $A/\mathfrak{b}$ durch
$$ \forall_{a, \tilde{a} \in A} \, (a + \mathfrak{b}) \, (\tilde{a} + \mathfrak{b})  = (a \, \tilde{a}) + \mathfrak{b} $$
gegeben, die $A/\mathfrak{b}$ zu einer $\K$-Algebra macht. Daß dies wohldefiniert ist, liegt daran, daß $\mathfrak{b}$ ein $\K$-Ideal von $X$ ist:
Gilt nämlich sowohl $a + \mathfrak{b} = a' + \mathfrak{b}$ als auch $\tilde{a} + \mathfrak{b} = \tilde{a}' + \mathfrak{b}$, so existieren $b,\tilde{b} \in \mathfrak{b}$ mit $a - a' = b$ und $\tilde{a} - \tilde{a}' = \tilde{b}$, also folgt
$$ a \, \tilde{a} = (a' + b) \, (\tilde{a}' + \tilde{b}) = a' \, \tilde{a}' + \underbrace{a' \, \tilde{b}}_{\in \mathfrak{b}} + \underbrace{b \, \tilde{a}'}_{\in \mathfrak{b}} + \underbrace{b \, \tilde{b}}_{\in \mathfrak{b}}, $$ 
d.h.\ $(a \, \tilde{a}) + \mathfrak{b} = (a' \, \tilde{a}') + \mathfrak{b}$.
Bezeichnet $\pi \: A \to A/\mathfrak{b}$ die kanonische Projektion, so gilt demnach $\forall_{a, \tilde{a} \in A} \, \pi(a) \, \pi(\tilde{a}) = \pi(a \, \tilde{a})$.
Besitzt $A$ ein Einselment $e$, so besitzt $A/\mathfrak{b}$ das Einselement $e + \mathfrak{b}$ und es gilt $\pi(e) = e + Y$.}
\end{itemize}
\end{Satz}

\textit{Beweis.} Zu (i): Sei $\xi \in X/Y$.
$\|\xi\|_{\sim} \ge 0$ ist trivial.

Gilt $\|\xi\|_{\sim} = 0$ und ist $x \in \xi$, letzteres heißt $\xi = x + Y$, so existiert offenbar eine Folge $(y_n)_{n \in \N}$ in $Y$, die gegen $x$ konvergiert.
Dann ergibt die Abgeschlossenheit von $Y$: $x = \lim_{n \to \infty} y_n \in Y$, also $\xi = Y = 0_{X/Y}$.

Ferner gilt für $\lambda \in \K$
$$ \| \lambda \, \xi \|_{\sim} =  \inf \{ \| \lambda \, x  \| \, | \, x \in \xi \} = |\lambda| \, \inf \{ \|x\| \, | \, x \in \xi \}  = |\lambda| \, \|\xi\|_{\sim} $$
sowie für jedes $\xi_1, \xi_2, \xi_3 \in X/Y$
\begin{eqnarray*}
\| \xi_1 + \xi_2 \|_{\sim} & = & \inf \{ \| x_1 + x_2 \| \, | \, \forall_{i \in \{1,2\}} \, x_i \in \xi_i \} \\
& \le & \inf\{ \| x_1 + x_3 \| + \| x_3 + x_2 \| \, | \, \forall_{i \in \{1,2,3\}} \, x_i \in \xi_i \} \\
& \le & \inf\{ \| x_1 + x_3 \| \, | \, \forall_{i \in \{1,3\}} \, x_i \in \xi_i \} + \inf \{ \| x_3 + x_2 \| \, | \, \forall_{i \in \{2,3\}} \, x_i \in \xi_i \} \\
& = & \| \xi_1 + \xi_3 \|_{\sim} + \| \xi_3 + \xi_2 \|_{\sim}.
\end{eqnarray*}

Zu (ii): Daß $\pi$ gleichmäßig stetig ist, folgt aus
\begin{equation} \label{FA.2.5.0}
\forall_{x, \tilde{x} \in X} \, \| \pi(x) - \pi(\tilde{x}) \| \stackrel{\pi \text{ linear}}{=} \| \pi(x - \tilde{x}) \|_{\sim} \stackrel{\text{$x-\tilde{x} \in \pi(x-\tilde{x})$}}{\le} \| x - \tilde{x} \|.
\end{equation}

Wir zeigen, daß für alle $x \in X$ und $\varepsilon \in \R_+$ gilt
\begin{equation} \label{FA.2.5.S}
\pi \left( U_{\varepsilon}(x) \right) = U_{\varepsilon}(\pi(x)).
\end{equation}

{[} Zum Nachweis von (\ref{FA.2.5.S}) ist wegen 
\begin{gather*}
\pi \left( U_{\varepsilon}(x) \right) \stackrel{\ref{FA.2.2} (ii)}{=} \pi \left( x + U_{\varepsilon}(0) \right) \stackrel{\pi \text{ linear}}{=} \pi(x) + \pi \left( U_{\varepsilon}(0_{X/Y}) \right), \\
U_{\varepsilon}( \pi(x) ) \stackrel{\ref{FA.2.2} (ii)}{=} \pi(x) + U_{\varepsilon}(\pi(0)) 
\end{gather*}
zu zeigen, daß $\pi(U_{\varepsilon}(0_{X/Y})) = U_{\varepsilon}(\pi(0))$ gilt.

Beweis hiervon: Zu ,, $\subset$`` sei $x_0 \in U_{\varepsilon}(0)$.
Dann gilt
\begin{equation} \label{FA.2.5.1}
\| \pi(x_0) - \pi(0) \|_{\sim} \stackrel{\text{$\pi$ linear}}{=} \|\pi(x_0)\|_{\sim} \stackrel{(\ref{FA.2.5.0})}{\le} \|x_0\| < \varepsilon, 
\end{equation}
also $\pi(x_0) \in U_{\varepsilon}(\pi(0))$. 

Zu ,,$\supset$`` sei $\xi \in U_{\varepsilon}(\pi(0))$, also existiert $x_0 \in \xi$ mit $\|x_0\| < \varepsilon$.
Es ergibt sich
$$ \| \underbrace{\xi}_{= \pi(x_0)} - \, \pi(0)  \|_{\sim} \stackrel{(\ref{FA.2.5.1})}{<} \varepsilon, $$
d.h.\ $\xi = \pi(x_0) \in U_{\varepsilon}(\pi(0))$.   {]}

Aus (\ref{FA.2.5.S}) folgt direkt die Offenheit der Abbildung $\pi \: X \to X/Y$.

Zu (iii): Sei $(\xi_n)_{n \in \N}$ eine Cauchyfolge in $X/Y$, d.h.\ insbesondere
$$ \forall_{n \in \N} \exists_{k_0 \in \N} \forall_{k,l \in \N, \, k,l \ge k_0} \, \|\xi_k - \xi_l\|_{\sim} < \frac{1}{2^n}. $$
Wir können daher eine Teilfolge $(\xi_{i_n})_{n \in \N}$ von $(\xi_n)_{n \in \N}$ mit
$$ \forall_{n \in \N} \, \|\xi_{i_{n+1}} - \xi_{i_n}\|_{\sim} < \frac{1}{2^n} $$
wählen.
Zu jedem $n \in \N$ existiert dann $x_n \in \xi_{i_n}$, d.h.\ $\pi(x_n) = \xi_{i_n}$, mit
$$  \|x_{n+1} - x_n\| < \frac{1}{2^n}, $$
und es folgt für alle $n_0,k \in \N$ 
$$ \|x_{n_0 + k + 1} - x_{n_0}\| \le  \|x_{n_0 + k + 1}  - x_{n_0 + k}\| + \ldots + \|x_{n_0 + 1} - x_{n_0} \| < \underbrace{\frac{1}{2^{n_0 + k}} + \ldots + \frac{1}{2^{no}}}_{\stackrel{n_0 \to \infty}{\longrightarrow}0}. $$
Dann ist offenbar $(x_n)_{n \in \N}$ eine Cauchyfolge in $X$, die nach Voraussetzung von (iii) gegen ein gewisses $x \in X$ konvergiert.
Folglich konvergiert wegen der Stetigkeit von $\pi$ auch $(\xi_{i_n})_{n \in \N}$ gegen $\pi(x)$. 

Zu (iv): Bezeichne $e$ das Einselement von $X$.
Trivialerweise gilt $\| e + \mathfrak{b} \|_{\sim} \le 1$.
Aus $\| e + Y \|_{\sim} < 1$ folgte die Existenz von $y \in Y$ mit $\|e - y\| < 1$, d.h.\ wie oben nach \ref{FA.2.Reihe.2} (v), daß $y = e - (e - y)$ invertierbar wäre, im Widerspruch zu dem Lemma in der Fußnote auf Seite \pageref{FA.2.4.L}.

Seien $\xi_1, \xi_2 \in X/Y$.
Dann gilt
\begin{eqnarray*}
\| \xi_1 \, \xi_2 \|_{\sim} & = & \inf \{ \| x_1 \, x_2 \| \, | \, \forall_{i \in \{1,2\}} \, x_i \in \xi_i \} \\
& \le & \inf\{ \| x_1 \| \, \| x_2 \| \, | \, \forall_{i \in \{1,2\}} \, x_i \in \xi_i \} \\
& \le & \inf\{ \| x_1 \| \, | \, x_1 \in \xi_1 \} \cdot \inf \{ \| x_2 \| \, | \,  x_2 \in \xi_2 \} \\
& = & \| \xi_1 \|_{\sim} \, \| \xi_2 \|_{\sim},
\end{eqnarray*}
womit der Satz bewiesen ist. \q

\begin{Bem*}
Aus (ii) folgt, daß für jede Teilmenge $U \subset X/Y$ gilt
$$ U \mbox{ offen } \Longleftrightarrow \overline{\pi}^1(U) \mbox{ offen,} $$
d.h.\ $\| \ldots \|_{\sim}$ induziert die Quotiententopologie auf $X/Y$.\footnote{Ist $X$ ein topologischer Raum und $\sim$ eine Äquivalenzrelation auf $X$, so definiert man die \emph{Quotiententopologie für $X/\sim \, =: X_{\sim}$}\index{Topologie!Quotienten-} durch
$$ \forall_{U \subset X_{\sim}} \, \mbox{$U$ offen in  $X_{\sim}$  :$\Longleftrightarrow$ $\overline{\pi}^1(U)$ offen in  $X$}, $$
wobei $\pi \: X \to X_{\sim}$ die kanonische Projektion bezeichne.
Somit ist eine Abbildung $f \: X_{\sim} \to M$, wobei $M$ ein topologischer Raum sei, genau dann stetig, wenn $f \circ \pi \: X \to M$ stetig ist.}
\end{Bem*}

\subsection*{Beschränkte Operatoren} \addcontentsline{toc}{subsection}{Beschränkte Operatoren}

\begin{Def}[(Stetige bzw.\ beschränkte) Operatoren] \label{FA.2.6}
Seien $X,Y$ $\K$-Vek\-tor\-räume.
\begin{itemize}
\item[(i)] Den $\K$-Vektorraum aller $\K$-linearen Abbildungen $X \to Y$ bezeichnen wir mit $\boxed{{\rm Hom}_{\K}(X,Y)}$.
Ein Element von ${\rm Hom}_{\K}(V,W)$ nennen wir einen \emph{Operator von $X$ nach $Y$}\index{Operator}.
Im Falle $Y=\K$ spricht man auch von einem \emph{Funktional von $X$}\index{Funktional}.
\end{itemize}
Im folgenden seien $X,Y$ sogar normiert, wobei wir beide Normen mit $\| \ldots \|$ bezeichnen.
\begin{itemize}
\item[(ii)] Ein \emph{stetiger Operator von $X$ nach $Y$}\index{Operator!stetiger} ist per definitionem eine stetiges Element von ${\rm Hom}_{\K}(X,Y)$.
Den $\K$-Vektorraum aller solchen Abbildungen bezeichnen wir mit $\boxed{\mathcal{L}_{\K}(X,Y)}$.
(Daß dies tatsächlich ein $\K$-Vektorraum ist, folgt daraus, daß für alle $\lambda, \mu \in \K$ und $T,S \in \mathcal{L}_{\K}(X,Y)$
$$ \forall_{x \in X} \, \| (\lambda \, T + \mu \, S) (x) \| \le |\lambda| \, \|T(x)\| + |\mu| \, \|S(x)\| $$
gilt, d.h.\ $\lambda \, T + \mu \, S$ ist mit $T,S$ stetig.)
\item[(iii)] Ein Operator $T \: X \to Y$ heißt \emph{beschränkt}\index{Operator!beschränkter} genau dann, wenn gilt
$$ \exists_{C \in \R} \forall_{x \in X} \, \|T(x)\| \le C \, \|x\|. $$
\end{itemize}
\end{Def}

\begin{Satz} \label{FA.2.7}
Seien $X,Y$ normierte $\K$-Vektorräume. (Beide Normen seien mit $\| \ldots \|$ bezeichnet.)
Für jedes $T \in {\rm Hom}_{\K}(X,Y)$ sind die folgenden Aussagen paarweise äquivalent:
\begin{itemize}
\item[(i)] $T$ ist Lipschitz-stetig, also insbesondere gleichmäßig stetig.
\item[(ii)] $T$ ist stetig, d.h.\ genau $T \in \mathcal{L}_{\K}(X,Y)$.
\item[(iii)] $T$ ist stetig in $0$.
\item[(iv)] $T$ ist ein beschränkter Operator.
\end{itemize}
\end{Satz}

\textit{Beweis.} ,,(i) $\Rightarrow$ (ii)`` und ,,(ii) $\Rightarrow$ (iii)`` sind trivial.

,,(iii) $\Rightarrow$ (iv)`` Angenommen, (iv) ist falsch. 
Dann existiert eine Folge $(x_n)_{n \in \N_+}$ in $X$ mit $\forall_{n \in \N_+} \, \| T(x_n) \| > n^2 \| x_n \|$, also $x_n \neq 0$. 
Folglich ist $\left( \frac{x_n}{n \, \| x_n \|} \right)_{n \in \N_+}$ eine Nullfolge mit (beachte, daß $T$ $\K$-linear ist)
$$\lim_{n \to \infty} \left\| T \left( \frac{x_n}{n \, \| x_n \|} \right) \right\| = \lim_{n \to \infty} \frac{\| T(x_n) \|}{n \, \| x_n \|} = \infty,$$ 
im Widerspruch zur Stetigkeit von $T$ in $0 = T(0)$.

,,(iv) $\Rightarrow$ (i)`` Es seien $x,\tilde{x} \in X$ beliebig.
Dann folgt aus der $\K$-Linearität von $T$ und (iv) die Existenz von $C \in \R$ derart, daß gilt
$$ \| T(x) - T(\tilde{x}) \| = \| T(x - \tilde{x}) \| \le C \| x - \tilde{x} \|, $$
d.h.\ $T$ ist Lipschitz-stetig. \q

\begin{Satz} \label{FA.2.H} $\,$

\noindent \textbf{Vor.:} Seien $X,Y$ normierte $\K$-Vektorräume, $T \in \mathcal{L}_{\K}(X,Y)$ und $W$ ein abgeschlossener $\K$-Untervektorraum von $X$ mit $W \subset {\rm Kern} \, T$, also ist auch $X/W$ ein normierter $\K$-Vekrorraum und die kanonische Abbildung $\pi \: X \to X/W$ stetig und offen, vgl. \ref{FA.2.5}.
Beachte, daß $W = {\rm Kern} \, T$ möglich ist, da ${\rm Kern} \, T$ als stetiges Urbild einer abgeschlossenen Menge selbst abgeschlossen in $X$ ist.
\pagebreak

\noindent \textbf{Beh.:}
\begin{itemize}
\item[(i)] Es existiert genau ein $S \in \mathcal{L}_{\K}(X/W,Y)$ mit $T = S \circ \pi$, insbes.\ kommutiert das folgende Diagramm:
%
%\begin{diagram}
%X            & \rTo^{T}  & Y \\
%\dOnto^{\pi} & \ruTo>{S} &   \\
%X/W          &           &   \\
%\end{diagram}
%
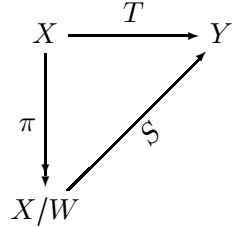
\begin{center}
\begin{tikzpicture}[node distance=2cm, every arrow/.style={thick,->}]
    % Knoten definieren
    \node (X) {$X$};
    \node (Y) [right of=X] {$Y$};
    \node (X/W) [below of=X] {$X/W$};
    
    % Pfeile zeichnen
    \draw[->] (X) -- node[above] {$T$} (Y);
    \draw[-{>>}] (X) -- node[left] {$\pi$} (X/W);
    \draw[->] (X/W) -- node[right] {$S$} (Y);
\end{tikzpicture}
\end{center}
\item[(ii)](Homomorphiesatz) \index{Satz!Homomorphie-}
\newline
Im Falle $W = {\rm Kern} \, T$ ist $S$ wie in (i) injektiv.
\end{itemize}
\end{Satz}

\textit{Beweisskizze.} Zu (i): Wegen der Surjektivität von $\pi$ ist die Eindeutigkeit von $S$ mit $T = S \circ \pi$ sofort klar.
Man kann $S$ daher nur folgendermaßen definieren 
$$ \forall_{\xi \in X/W} \, S(\xi) := T(x), \mbox{ wobei } x \in X \mbox{ mit } \pi(x) = \xi \mbox{ sei}. $$
Wir haben zu zeigen, daß dies wohldefiniert ist.
Seien daher $\xi \in X/W$ und $x_1, x_2 \in X$ mit $\pi(x_1) = \pi(x_2) = \xi$.
Dann gilt $x_1 - x_2 \in W$, also wegen $W \subset {\rm Kern} \, T$ auch $T(x_1) = T(x_2)$.

Die $\K$-Linearität von $S$ ergibt sich sofort aus der Definition von $S$.
Zu zeigen bleibt die Stetigkeit von $S$.
Hierzu sei $U \subset Y$ eine offene Menge.
Wegen der Stetigkeit von $T$ und der Offenheit von $\pi$ ist dann auch $\pi(\overline{T}^1(U))$ offen in $X/W$.
Wir zeigen
\begin{equation} \label{FA.2.H.1}
\overline{S}^1(U) = \pi(\overline{T}^1(U)).
\end{equation}
Dann folgt aus der Beliebigkeit von $U$, daß $S$ stetig ist.

{[} Zu (\ref{FA.2.H.1}): Es gilt wegen der Surjektivität von $\pi$
\begin{eqnarray*}
\xi \in \overline{S}^1(U) & \Longleftrightarrow & S(\xi) \in U \, \wedge \, \exists_{x \in X} \, \pi(x) = \xi \Longleftrightarrow \exists_{x \in X} \, T(x) \in U \, \wedge \, \pi(x) = \xi \\
& \Longleftrightarrow & \exists_{x \in \overline{T}^1(U)} \, \pi(x) = \xi \Longleftrightarrow \xi \in \pi(\overline{T}^1(U))
\end{eqnarray*}
für jedes $\xi \in X/{\rm Kern} \, S$. {]}

Zu (ii): Seien $\xi_1 ,\xi_2 \in X/{\rm Kern} \, T$ mit $S(\xi_1) = S(\xi_2)$.
Es existieren $x_1,x_2 \in X$ mit $\pi(x_1) = \xi_1$ sowie $\pi(x_2) = \xi_2$.
Dann gilt $T(x_1) = S(\xi_1) = S(\xi_2) = T(x_2)$, d.h.\ $x_1 - x_2 \in {\rm Kern} \, T$, also folgt $\pi(x_1 - x_2) = 0$ und somit $\xi_1 = \pi(x_1) = \pi(x_2) = \xi_2$. \q

\begin{Bem*}
Seien $X,Y$ normierte $\K$-Vektorräume und $T \in \mathcal{L}_{\K}(X,Y)$.
Dann ist $S \: X/{\rm Kern} \, T \to T(X)$ wie in (ii) ein stetiger $\K$-Vektorraum-Isomorphismus aber i.a.\ kein Homöomorphismus.
Der Satz vom inversen Operator (siehe \ref{FA.2.37} unten) ergibt allerdings die Stetigkeit von $S^{-1} \: T(X) \to X/{\rm Kern} \, T$, falls $X$ und $T(X)$ vollständig sind.
\end{Bem*}

\begin{Satz}[über die Operatornorm] \label{FA.2.8}
Seien $X,Y$ normierte $\K$-Vektorräume, wobei beide Normen mit $\| \ldots \|$ bezeichnet seien.
\begin{itemize}
\item[(i)] Indem man für jedes $T \in \mathcal{L}_{\K}(X,Y)$
$$ \boxed{\|T\|_{\mathcal{L}_{\K}(X,Y)}} := \sup \left\{ \frac{\|T(x)\|}{\|x\|} \, | \, x \in X \setminus \{0\} \right\} ( := 0, \mbox{ falls } X = \{0\}) $$
setzt, wird eine Norm, die sog.\ \emph{Operatornorm\index{Operatornorm}\index{Norm!Operator-} auf $\mathcal{L}_{\K}(V,W)$}, definiert.

Wir betrachten $\mathcal{L}_{\K}(X,Y)$ im folgenden stets als durch die Operatornorm, die wir meistens auch wieder mit $\| \ldots \|$ bezeichnen, normierten $\K$-Vek\-tor\-raum.
\item[(ii)] Für jedes $T \in \mathcal{L}_{\K}(X,Y)$ gilt
\begin{equation} \label{FA.2.8.Stern}
\| T \| = \sup \{ \|T(x)\| \, | \, x \in X \, \wedge \, \| x \| \le 1 \} = \underbrace{\sup \{ \|T(x)\| \, | \, x \in X  \, \wedge \, \| x \| = 1 \}}_{:= 0, \text{ falls } X = \{0\}}
\end{equation}
sowie
\begin{equation} \label{FA.2.8.Doppelstern}
\forall_{x \in X} \, \| T(x) \| \le \|T\| \, \|x\|.
\end{equation}

Darüber hinaus ist $\|T\|$ die kleinste reelle Zahl mit der Eigenschaft (\ref{FA.2.8.Doppelstern}).
\end{itemize}
\end{Satz}

\textit{Beweis.} Ohne Einschränkung gelte $\dim_{\K} X > 0$.

Zu (i): \ref{FA.2.7} ,,(i) $\Rightarrow$ (iii)`` gibt die Wohldefiniertheit von $\|T\|$ als reelle Zahl für $T \in \mathcal{L}_{\K}(X,Y)$.

Beweis, daß $\| \ldots \|$ eine Norm für $\mathcal{L}_{\K}(X,Y)$ ist:

(N1) ist trivial. Weiter gilt für alle $T,S \in \mathcal{L}_{\K}(X,Y)$, $\lambda \in \K$ und $x \in X \setminus \{0\}$
\begin{gather*}
\D \frac{\| (\lambda \, T)(x) \|}{\|x\|} = |\lambda| \frac{\|T(x)\|}{\|x\|}, \\
\D \frac{\|(T+S)(x)\|}{\|x\|} \le \frac{\|T(x)\|}{\|x\|} + \frac{\|S(x)\|}{\|x\|} \le \|T\| + \|S\|, \\
\end{gather*}
also folgen (N2) und (N3).

Zu (ii): Sei $T \in \mathcal{L}_{\K}(V,W)$.
Für alle $x \in X \setminus \{0\}$ mit $\|x\| \le 1$ gilt 
$$ \|T(x)\| \le \frac{\|T(x)\|}{\|x\|} \le \|T\|, $$
weshalb 
$$ \sup \{ \|T(x)\| \, | \, x \in X  \, \wedge \, \| x \| = 1 \} \le \sup \{ \|T(x)\| \, | \, x \in X  \, \wedge \, \| x \| \le 1 \} \le \|T\| $$
folgt.

Andererseits gilt für $x \in X \setminus \{0\}$
$$ \frac{\|T(x)\|}{\|x\|} = \left\| T \left( \frac{x}{\|x\|} \right) \right\| \le \sup \{ \|T(\tilde{x})\| \, | \, \tilde{x} \in X  \, \wedge \, \| \tilde{x} \| = 1 \},$$
also folgt auch $\|T\| \le \sup \{ \|T(\tilde{x})\| \, | \, \tilde{x} \in X  \, \wedge \, \| \tilde{x} \| = 1 \}$.
Damit ist (\ref{FA.2.8.Stern}) gezeigt.

Die übrigen Aussagen in (ii) folgen sofort aus der Definition in (i). \q

\begin{Bem*}
Sind $X,Y$ normierte $\K$-Vektorräume, so spielt die Fastmetrik $d_{\infty}$, die ${\rm Hom}_{\K}(V,W)$ als $\K$-Untervektorraum des $\K$-Vektorraumes $Y^X$ erbt, keine Rolle, da für alle Abbildungen $T,S \in {\rm Hom}_{\K}(X,Y)$ mit $T \ne S$
$$ d_{\infty}(T,S) = \sup \{ \| T(x) - S(x) \| \, | \, x \in X \} = \infty $$
gilt; denn aus $T \ne S$ folgt die Existenz von $x \in X$ mit $T(x) \ne S(x)$, und es gilt für alle $x \in X$ und $\lambda \in \R$
$$ \| T(\lambda \, x) - S(\lambda \, x)\| = |\lambda| \, \underbrace{\|T(x) - S(x)\|}_{> 0} \stackrel{\lambda \to \infty}{\longrightarrow} \infty, $$
also $d_{\infty}(T,S) = \infty$.
Aus diesem Grunde ist auch der Begriff der ,,gleichmäßigen Konvergenz`` für Operatoren nicht relevant.
\end{Bem*}

\begin{Satz} \label{FA.2.9}
Seien $X,Y,Z$ normierte $\K$-Vektorräume.

Dann gilt für alle $T \in \mathcal{L}_{\K}(X,Y)$ und $S \in \mathcal{L}_{\K}(Y,Z)$
$$ S \circ T \in \mathcal{L}_{\K}(X,Z) ~~ \wedge ~~ \|S \circ T\| \le \|S\| \, \|T\|. $$
\end{Satz}

\textit{Beweis.} Für jedes $x \in X$ gilt nach (\ref{FA.2.8.Doppelstern})
$$ \|(S \circ T)(x)\| \le \|S\| \, \|T(x)\| \le \underbrace{\|S\| \, \|T\|}_{=: C} \, \|x\|,$$ 
d.h.\ $S \circ T$ ist ein beschränkter Operator, und aus der letzten Aussage in \ref{FA.2.8} (ii) folgt
$$ \|S \circ T\| \le C = \|S\| \, \|T\|. $$
\q

\begin{Satz} \label{FA.2.10}
Sei $X$ ein normierter $\K$-Vektorraum.

Dann ist $\mathcal{L}_{\K}(X,X)$ bzgl.\ der Multiplikation
$$ \mathcal{L}_{\K}(X,X) \times \mathcal{L}_{\K}(X,X) \longrightarrow \mathcal{L}_{\K}(X,X), ~~ (S,T) \longmapsto S \circ T, $$
eine normierte $\K$-Algebra mit Einselement $\id_X$.
\end{Satz}

\textit{Beweis.} Wegen des vorherigen Satzes ist nur zu zeigen, daß $\| \id_X \| = 1$ gilt.
Dies folgt aus (\ref{FA.2.8.Stern}). \q

\begin{Satz} \label{FA.2.11}
Es seien $X,Y$ normierte $\K$-Vektorräume und $Y$ sogar ein $\K$-Banach\-raum.

Dann ist $\mathcal{L}_{\K}(X,Y)$ ein $\K$-Banachraum.
\end{Satz}

\textit{Beweis.} Sei $(T_n)_{n \in \N}$ eine Cauchyfolge in $\mathcal{L}_{\K}(X,Y)$.

1.) Sei $x \in X$.
Dann gilt für alle $n,m \in \N$
$$ \|T_n(x) - T_m(x)\| = \| (T_n-T_m)(x) \| \stackrel{(\ref{FA.2.8.Doppelstern})}{\le} \|T_n - T_m\| \, \|x\|, $$
also ist offenbar $(T_n(x))_{n \in \N}$ eine Cauchyfolge in $Y$, da $(T_n)_{n \in \N}$ eine Cauchyfolge in $\mathcal{L}_{\K}(X,Y)$ ist.
$Y$ ist als vollständig vorausgesetzt, also können wir definieren
$$ T(x) := \lim_{n \to \infty} T_n(x). $$

2.) Wir behaupten, daß für die durch 1.) definierte Abbildung $T \: X \to Y$ gilt
\begin{gather}
\mbox{$T$ ist $\K$-linear,} \label{FA.2.11.1} \\
\mbox{$T$ ist stetig,} \label{FA.2.11.2} \\
\lim_{n \to \infty} T_n = T, \label{FA.2.11.3}
\end{gather}
womit der Satz bewiesen ist.

Zu (\ref{FA.2.11.1}): Seien $x, \tilde{x} \in X$ sowie $\lambda, \mu \in \K$.
Dann gilt
\begin{eqnarray*}
T(\lambda \, x + \mu \, \tilde{x}) & = & \lim_{n \to \infty} (T_n(\lambda \, x + \mu \, \tilde{x})) = \lim_{n \to \infty} (\lambda \, T_n(x) + \mu \, T_n(\tilde{x})) \\
& \stackrel{\ref{FA.2.2} (iv)}{=} & \lambda \, \lim_{n \to \infty} T_n(x) + \mu \, \lim_{n \to \infty} T_n(\tilde{x}) = \lambda \, T(x) + \mu \, T(\tilde{x}).
\end{eqnarray*}

Zu (\ref{FA.2.11.2}): Für alle $n,m \in \N$ gilt nach \ref{FA.2.2} (iii) $| \|T_n\| - \|T_m\| | \le \| T_n - T_m \|$, also ist $(\|T_n\|)_{n \in \N}$ eine Cauchyfolge in $\R$, da $(T_n)_{n \in \N}$ eine Cauchyfolge in $\mathcal{L}_{\K}(X,Y)$ ist.
Aus der Vollständigkeit von $\R$ folgt die Konvergenz von $(\|T_n\|)_{n \in \N}$ in $\R$, d.h.\ insbes.\
$$ \exists_{C \in \R} \forall_{n \in \N} \, \|T_n\| \le C, $$
also gilt auch für jedes $x \in X$
$$ \| T(x) \| = \| \lim_{n \to \infty} T_n(x) \| \stackrel{\ref{FA.2.2}(iii)}{=} \lim_{n \to \infty} \underbrace{\| T_n(x) \|}_{\le \|T_n\| \, \|x\|} \le C \, \|x\|, $$
d.h.\ $T$ ist ein beschränkter Operator.

Zu (\ref{FA.2.11.3}): Für jedes $x \in X$ mit $\|x\| \le 1$ und jedes $n \in \N$ gilt
\begin{eqnarray*}
\|(T_n-T)(x)\| & = & \|T_n(x) - T(x)\| = \|T_n(x) - \lim_{m \to \infty} T_m(x) \| \\
& \stackrel{\ref{FA.2.2} (ii)}{=} & \lim_{m \to \infty} \|T_n(x)-T_m(x)\| = \lim_{m \to \infty} \underbrace{\|(T_n - T_m)(x)\|}_{\stackrel{(\ref{FA.2.8.Doppelstern})}{\le} \|T_n - T_m\| \, \underbrace{\|x\|}_{\le 1}} \\
& \le & \lim_{m \to \infty} \|T_n - T_m \|.
\end{eqnarray*}
Daher genügt es zum Nachweis von (\ref{FA.2.11.3}) wegen (\ref{FA.2.8.Stern}) zu zeigen, daß für hinreichend großes $n \in \N$ gilt $\lim_{m \to \infty} \|T_n - T_m \| = 0$.
Sei hierzu $\varepsilon \in \R_+$.
Da $(T_n)_{n \in \N}$ eine Cauchyfolge ist, existiert $n_0 \in \N$ mit
$$ \forall_{n,m \in \N, \, n \ge n_0} \, \|T_n - T_m\| < \frac{\varepsilon}{2}, $$
also gilt für alle $n \ge n_0$: $\lim_{m \to \infty} \|T_n - T_m\| < \varepsilon$. \q

\begin{Bem*}
Da die Operatornorm nicht die Supremumsfastmetrik induziert, vgl.\ die Bemerkung nach \ref{FA.2.8}, kann der letzte Satz nicht auf \ref{FA.1.Isom.1} (ii) zurückgeführt werden.
Man hat übrigens auch keine zu \ref{FA.1.Isom.1} analoge kanonische Abbildung von $Y$ nach $\mathcal{L}_{\K}(X,Y)$, da konstante Abbildungen, die nicht die Nullabbildung sind, nicht linear sind.
\end{Bem*}

Aus \ref{FA.2.10} und \ref{FA.2.11} ergibt sich sofort das folgende Korollar.

\begin{Kor} \label{FA.2.11.K1}
Ist $X$ ein $\K$-Banachraum, so ist $\mathcal{L}_{\K}(X,X)$ eine $\K$-Banach\-algebra mit Einselement $\id_X$. \q
\end{Kor}

\begin{Satz} \label{FA.2.10.M} 
Sei $A$ eine normierte $\K$-Algebra mit Einselement $e$.
\begin{itemize}
\item[(i)] Für alle $a \in A$ sind die \emph{Linksmultiplikation mit $a$}\index{Multiplikation!Links-}
$$ \boxed{L_a \:  A \longrightarrow A, ~~ x \longmapsto a \, x,} $$
und die \emph{Rechtsmultiplikation mit $a$}\index{Multiplikation!Rechts-}
$$ \boxed{R_a \: A \longrightarrow A, ~~ x \longmapsto x \, a,} $$
beschränkte Operatoren, und es gilt $\|L_a\| = \|R_a\| = \|a\|$.
\item[(ii)] Die Abbildung $\boxed{L \: A \to \mathcal{L}_{\K}(A,A)}$, die durch 
$$ \forall_{a \in A} \, L \, a := L_a $$
definiert ist, ist ein isometrischer $\K$-Algebra-Homomorphismus%
\footnote{\textbf{Definition.} Sind $A,B$ assoziative $\K$-Algebren mit Einselement, wobei wir beide Einselemente mit $e$ bezeichnen, so heißt $\varphi \in {\rm Hom}_{\K}(A,B)$ ein \emph{$\K$-Algebra-Homomorphismus}\index{Algebra!-Homomorphismus} genau dann, wenn gilt $\varphi(e) = e \, \wedge \, \forall_{a, \tilde{a} \in A} \, \varphi(a \, \tilde{a}) = \varphi(a) \, \varphi(\tilde{a})$.

Sind $A$ und $B$ normiert und gilt $\forall_{a \in A} \, \|\varphi(a)\| = \|a\|$, so heißt $\varphi$ ein \emph{isometrischer $\K$-Algebra-Ho\-mo\-mor\-phis\-mus}\index{Algebra!-Homomorphismus!isometrischer}\index{Abbildung!isometrische} und, falls $\varphi$ zusätzlich surjektiv ist, eine \emph{$\K$-Algebra-Isometrie}.\index{Algebra!-Isometrie}\index{Isometrie!normierter Algebren}
Wegen der $\K$-Li\-ne\-a\-ri\-tät von $\varphi$ ist dies mit den Definitionen des letzten Kapitels konform.}
, also inbesondere injektiv und gleichmäßig stetig.
\item[(iii)] $L(A)$ ist eine $\K$-Unteralgebra von $\mathcal{L}_{\K}(A,A)$ mit Einselement $L_e = \id_X$.
\end{itemize}
\end{Satz}

\textit{Beweis.} Zu (i): Sei $a \in A$.

$L_a$ ist offenbar $\K$-linear, und für jedes $x \in X$
\begin{gather*}
\|L_a(x)\| = \| a  \, x \| \le \|a\| \, \|x\|, \\
\|L_a(e)\| = \| a \|,
\end{gather*}
also folgt $L_a \in \mathcal{L}_{\K}(A,A)$ und $\|L_a\| = \|a\|$, beachte (\ref{FA.2.8.Stern}).

Die Behauptung für $R_a$ sieht man analog ein.

Zu (ii): Da für alle $\lambda \in \K$, $a,b \in A$ sowie $x \in A$
\begin{gather*}
L_{\lambda \, a}(x) = \lambda \, (a \, x) = \lambda \, L_a(x) = (\lambda \, L_a)(x), \\
L_{a+b}(x) = (a+b) \, x = a \, x + b \, x = L_a(x) + L_b(x) = (L_a + L_b)(x), \\
L_{e}(x) = e \, x = x = \id_A(x), \\
L_{a \, b}(x) = (a \, b) \, x = a \, (b \, x) = L_a(L_b(x)) = (L_a \circ L_b)(x)
\end{gather*}
gilt, folgt (ii) aus (i). 

Zu (iii): Die Behauptung ergibt sich sofort aus (ii). \q

\begin{Satz} \label{FA.2.Reihe.3}
Sei $A$ eine $\K$-Banachalgebra mit Einselement $e$.

Dann ist die (multilplikative) Gruppe der invertierbaren Elemente von $A$
$$ \boxed{G(A)} := \{ x \in A \, | \, \exists_{y \in A} \, x \, y = y \, x = e \} $$
offen in $A$ und die \emph{Inversenabbildung}
$$ \boxed{\inv \: G(A) \longrightarrow G(A), ~~ x \longmapsto x^{-1},} $$ 
ist ein Homöomorphismus.
\end{Satz}

\textit{Beweis.} 1.) Wir zeigen zunächst für jedes $a \in G(A)$, daß die Linksmultiplikation mit $a$
$$ L_a \: A \longrightarrow A, ~~ x \longmapsto a \, x, $$
ein Homöomorphismus ist:

a) $L_a$ ist nach \ref{FA.2.10.M} (i) ein beschränkter Operator, d.h. genau stetig.

b) Wegen der Beliebigkeit von $a \in G(A)$ ist auch $L_{a^{-1}} \: A \to A$ stetig, also folgt aus
$$ L_a \circ L_{a^{-1}} = L_{a^{-1}} \circ L_a = \id_A $$
die Bijektivität von $L_a$ mit $(L_a)^{-1} = L_{a^{-1}}$, d.h.\ $L_a$ ist ein Homöomorphismus.

2.) Analog zu 1.) sieht man ein, daß für jedes $a \in G(A)$ die Rechtsmultiplikation mit $a$ 
$$ R_a \: A \longrightarrow A, ~~ x \longmapsto x \, a, $$
ein Homöomorphismus ist.

3.) Wir zeigen weiterhin, daß $G(A)$ offen ist:

Hierfür zeigen wir zunächst
\begin{gather} 
U_1(e) \subset G(A), \label{FA.2.Reihe.3.1} \\
\forall_{a \in G(A)} \, L_a(U_1(e)) \subset G(A). \label{FA.2.Reihe.3.2}
\end{gather}

{[} Zu (\ref{FA.2.Reihe.3.1}): Sei $x \in U_1(e)$, d.h.\ $\|e-x\| < 1$.
Aus \ref{FA.2.Reihe.2} (ii), (iv), (iii) folgt dann auch $x = e-(e-x) \in G(A)$. 

Zu (\ref{FA.2.Reihe.3.2}): Seien $a \in G(A)$ und $x \in U_1(e) \stackrel{(\ref{FA.2.Reihe.3.1})}{\subset} G(A)$.
Dann gilt natürlich auch $L_a(x) = a \, x \in G(A)$ (mit $(a \, x)^{-1} = x^{-1} \, a^{-1}$). {]}

Sei nun $a \in G(A)$ beliebig.
Als Homöomorphismus ist $L_a$ eine offene Abbildung, also folgt aus (\ref{FA.2.Reihe.3.2}), daß $a = L_a(e)$ ein innerer Punkt von $G(A)$ ist.

4.) Wir zeigen schließlich, daß $\inv \: G(A) \to G(A)$ ein Homöomorphismus ist:

a) Sei $(x_n)_{n \in \N}$ eine Folge in $G(A)$ mit $\lim_{n \to \infty} x_n = e$, d.h.\ $\lim_{n \to \infty} e - x_n =  0$, und wir können ohne Einschränkung annehmen, daß gilt $\forall_{n \in \N} \, \|e - x_n\| < 1$.
Dann folgt wieder aus \ref{FA.2.Reihe.2} (ii), (iv), (iii) $x_n = e-(e-x_n) \in G(A)$ sowie des weiteren ${x_n}^{-1} = \sum_{k=0}^{\infty} (e-x_n)^k$ für jedes $n \in \N$, und es gilt
\begin{eqnarray*}
\| \inv(x_n) - \inv(e) \| & = & \left\| \left( \sum_{k=0}^{\infty} (e - x_n)^k \right) - e \right\| = \left\| \sum_{k=1}^{\infty} (e - x_n)^k \right\| \\
& \le & \sum_{k=1}^{\infty} \|e - x_n\|^k  = \frac{1}{1- \|e - x_n\|} - 1 \stackrel{n \to \infty}{\longrightarrow} 0,
\end{eqnarray*}
d.h.\ $\inv$ ist stetig in $e$.
\pagebreak

b) Für jedes $a \in G(A)$ gilt offenbar
$$ \inv = L_{a^{-1}} \circ \inv \circ R_{a^{-1}}. $$
$R_{a^{-1}}$ ist nach 2.) stetig in $a$, J ist nach a) stetig in $e = R_{a^{-1}}(a)$ und $L_{a^{-1}}$ ist nach 1.) stetig in $e = \inv(e)$, also ist auch $\inv$ stetig in $a$.
Aus der Beliebigkeit von $a \in G(A)$ folgt die Stetigkeit von $\inv$.

c) Da $\inv$ offenbar selbstinvers ist, ergibt b), daß $\inv$ ein Homöomorphismus ist. \q

\begin{Satz}[Fortsetzungssatz beschränkter Operatoren] \label{FA.2.12} \index{Satz!Fortsetzungs-!beschränkter Operatoren} $\,$

\noindent \textbf{Vor.:} Es seien $X,Y$ normierte $\K$-Vek\-tor\-räume und $Y$ sogar ein $\K$-Banach\-raum.
Ferner seien $W$ ein $\K$-Untervektorraum von $X$ sowie $T \in \mathcal{L}_{\K}(W,Y)$.

\noindent \textbf{Beh.:} Es existiert genau ein $\overline{T} \in \mathcal{L}_{\K}(\overline{W},Y)$ mit $\overline{T}|_{W} = T$.
Des weiteren gilt $\|\overline{T}\| = \|T\|$.
\end{Satz}

\textit{Beweis.} $T$ ist nach \ref{FA.2.7} ,,(ii) $\Rightarrow$ (i)`` gleichmäßig stetig, also existiert gemäß \ref{FA.1.20} ein $\overline{T} \in \mathcal{C}(\overline{W},Y)$ mit $\overline{T}|_W = T$.
Zu zeigen bleiben die $\K$-Linearität von $\overline{T}$ und $\|\overline{T}\| = \|T\|$.

Zur $\K$-Linearität: Seien $\lambda \in \K$ und $x,\tilde{x} \in \overline{W}$.
Dann existieren Folgen $(w_n)_{n \in \N}$ sowie $(\tilde{w}_n)_{n \in \N}$ in $W$ mit $\lim_{n \to \infty} w_n = x$ sowie $\lim_{n \to \infty} \tilde{w}_n = \tilde{x}$, und es folgt aus Stetigkeitsgründen und der $\K$-Linearität von $T$
\begin{gather*}
\overline{T}(\lambda \, x) = \overline{T}( \lim_{n \to \infty} \lambda \, w_n) = T( \lim_{n \to \infty} \lambda \, w_n ) = \lambda \, T( \lim_{n \to \infty} w_n ) = \lambda \, \overline{T}( \lim_{n \to \infty} w_n ) = \lambda \, \overline{T}(x), \\
\overline{T}(x + \tilde{x}) = T( \lim_{n \to \infty} (w_n + \tilde{w}_n) ) = \lim_{n \to \infty} T(w_n) + \lim_{n \to \infty} T(\tilde{w}_n) = \overline{T}(x) + \overline{T}(\tilde{x}).
\end{gather*}

Zu $\|\overline{T}\| = \|T\|$:
Zu $x \in X$ existiert wieder eine Folge $(w_n)_{n \in \N}$ in $W$ mit $\lim_{n \to \infty} w_n = x$, und es ergibt sich erneut aus Stetigkeisgründen
$$ \| \overline{T} \, x \| = \| \lim_{n \to \infty} T \, w_n \| = \lim_{n \to \infty} \| T \, w_n \| \le \|T\| \, \lim_{n \to \infty} \|w_n\| = \|T\| \, \|x\|, $$
d.h.\ nach der letzten Aussage in \ref{FA.2.8} (ii): $\|\overline{T}\| \le \|T\|$.
Des weiteren gilt 
\begin{eqnarray*}
\| T \| & = & \sup \{ \|T \, w\| \, | \, w \in W \, \wedge \, \|w\| \le 1 \} = \sup \{ \|\overline{T} \, w\| \, | \, w \in W \, \wedge \, \|w\| \le 1 \} \\
& \le & \sup \{ \|\overline{T} \, x\| \, | \, x \in \overline{W} \, \wedge \, \|x\| \le 1 \} = \|\overline{T}\|,
\end{eqnarray*}
also folgt $\|\overline{T}\| = \|T\|$. \q

\begin{Kor} \label{FA.2.12.K}
Sind $A,B$ normierte $\K$-Alge\-bren, wobei wir beide Einselemente mit $e$ bezeichnen, $C$ eine $\K$-Unteralgebra von $A$ mit $e \in C$ und $T \: C \to B$ ein stetiger $\K$-Algebra-Homomorphismus, so existiert genau ein stetiger $\K$-Algebra-Homomorphismus $\overline{T} \: \overline{C} \to B$ mit $\overline{T}|_C = T$. 
Des weiteren gilt $\|\overline{T}\| = \|T\|$.
\end{Kor}

\textit{Beweis.} Da $\overline{T}(e) = T(e) = e$ klar ist, müssen wir wegen des letzten Satzes nur 
$$ \forall_{a,\tilde{a} \in \overline{C}} \, \overline{T}(a \, \tilde{a}) = \overline{T}(a) \, \overline{T}(\tilde{a}) $$
nachweisen:
Seien also $a,\tilde{a} \in \overline{C}$ und $(c_n)_{n \in \N}$ sowie $(\tilde{c}_n)_{n \in \N}$ in $C$ mit $\lim_{n \to \infty} c_n = a$ sowie $\lim_{n \to \infty} \tilde{c}_n = \tilde{c}$.
Dann folgt analog zum Beweis der $\K$-Linearität im Beweis des letzten Satzes
$$ \overline{T}(a \, \tilde{a}) = T( \lim_{n \to \infty} (c_n \, \tilde{c}_n) ) = \lim_{n \to \infty} T(c_n) \cdot \lim_{n \to \infty} T(\tilde{c}_n) = \overline{T}(a) \, \overline{T}(\tilde{a}). $$
\q 

\begin{HS}[Vervollständigung normierter $\K$-Vektorräume] \label{FA.2.13} \index{Vervollständigung!normierter!Vektorräume} $\,$

\noindent \textbf{Vor.:} Sei $(X, \| \ldots \|)$ ein normierter $\K$-Vektorraum.

\noindent \textbf{Beh.:} Es existiert ein $\K$-Banachraum $(\widehat{X},\widehat{\| \ldots \|})$ mit folgenden Eigenschaften:
\begin{itemize}
\item[(i)] $X$ ist $\K$-Untervektorraum von $\widehat{X}$ mit $\overline{X} = \widehat{X}$ und $\| \ldots \| = \widehat{\| \ldots \|}|_X$.
\item[(ii)] Der $\K$-Banachraum $(\widehat{X},\widehat{\| \ldots \|})$ mit (i) ist bis auf Isometrie normierter $\K$-Vektorräume%
\footnote{\textbf{Definition.} Eine $\K$-lineare Abbildung $T \: X \to Y$ zwischen normierten $\K$-Vektorräumen heißt genau dann \emph{isometrisch}\index{Abbildung!isometrische}, wenn $\forall_{x \in X} \, \|T(x)\| = \|x\|$ gilt. 
Ist $T$ zusätzlich surjektiv, so heißt $T$ eine \emph{$\K$-Vektorraum-Isometrie}.\index{Isometrie!Vektorraum-}
Wegen der $\K$-Linearität von $T$ ist dies mit den Definitionen des letzten Kapitels konform.}%
eindeutig bestimmt, daher nennt man $(\widehat{X},\widehat{\| \ldots \|})$ ,,die`` \emph{Vervollständigung von $(X,\| \ldots \|)$}. 
\end{itemize}
\end{HS}

\textit{Beweis.} Zu (i): Obwohl die Behauptung eine einfache Folgerung aus \ref{FA.2.TB.6} sein wird, zeigen wir, wie man sie aus \ref{FA.1.21} herleitet.
Hiernach besitzt $X$ mit der induzierten Metrik, die durch
$$ \forall_{x,y \in X} \, d(x,y) = \|x-y\| $$
gegeben ist, eine bis auf Isometrie metrischer Räume eindeutig bestimmte Vervollständigung $(\widehat{X},\widehat{d})$, die $X$ als dichte Teilmenge enthält.
Wir konstruieren auf $\widehat{X}$ die Struktur eines normierten $\K$-Vektorraumes:

Wir betrachten auf $\widehat{X} \times \widehat{X}$ die Produktmetrik, vgl.\ \ref{FA.1.1} Beispiel 4.), die wir hier mit $\widehat{d} + \widehat{d}$ bezeichnen, und behaupten
\begin{equation} \label{FA.2.13.1}
\overline{X \times X} = \widehat{X} \times \widehat{X}.
\end{equation}

{[} Zu (\ref{FA.2.13.1}): ,,$\subset$`` ist trivial.

,,$\supset$`` Sei $(x_*,y_*) \in \widehat{X} \times \widehat{X}$. 
Wegen $\overline{X} = \widehat{X}$ existieren Folgen $(x_n)_{n \in \N}$ und $(y_n)_{n \in \N}$ in $X$ mit $\lim_{n \to \infty} x_n = x_*$ sowie $\lim_{n \to \infty} y_n = y_*$, und es gilt für jedes $n \in \N$
$$ (\widehat{d} + \widehat{d}) ( (x_n,y_n), (x_*,y_*) ) = \widehat{d}(x_n,x_*) + \widehat{d}(y_n,y_*) \stackrel{n \to \infty}{\longrightarrow} 0, $$
also nach \ref{FA.1.9}: $(x_*,y_*) \in \overline{X \times X}$. {]}

Nun ist $A \: X \times X \to X \subset \widehat{X}$, definiert durch
$$ \forall_{x,y \in X} \, A(x,y) := x + y, $$
wegen
$$ \forall_{(x,y), (\tilde{x},\tilde{y}) \in X \times X} \| (x + y) - (\tilde{x} + \tilde{y}) \| = \| (x - \tilde{x}) + (y - \tilde{y}) \| \le \| x - \tilde{x} \| + \| y - \tilde{y} \| $$
eine gleichmäßig stetige Abbildung.
Aus der Vollständigkeit von $\widehat{X}$ und (\ref{FA.2.13.1}) folgt mittels \ref{FA.1.20} die Existenz einer eindeutig bestimmten stetigen Fortsetzung $\overline{A} \: \widehat{X} \times \widehat{X} \to \widehat{X}$ von $A$.
Wir erweitern die Addition von $X$ nun auf $\widehat{X}$, indem wir
$$ \forall_{x_*, y_* \in \widehat{X}} \, x_* + y_* := \overline{A}(x_*,y_*) $$
setzen.

Des weiteren wird für jedes $\lambda \in \K$ durch
$$ \forall_{x \in X} \, S_{\lambda}(x) := \lambda \, x $$
offenbar ein beschränkter Operator $S_{\lambda} \: X \to X \subset \widehat{X}$ definiert, der nach \ref{FA.2.12} zu einem beschränkten Operator $\overline{S_{\lambda}} \: \widehat{X} \to \widehat{X}$ fortsetzbar ist.
Wir erweitern die skalare Multiplikation von $X$ auf $\widehat{X}$ via
$$ \forall_{\lambda \in \K} \forall_{x_* \in \widehat{X}} \, \lambda \cdot x_* := \overline{S_{\lambda}}(x_*). $$

Nun weisen wir nach:
\begin{equation} \label{FA.2.13.2}
(\widehat{X},+, \, \cdot) \mbox{ ist ein $\K$-Vektorraum.}
\end{equation}

{[} Zu (\ref{FA.2.13.2}): Es ist zu zeigen, daß für alle $\lambda, \mu \in \K$ sowie $x_*, y_*, z_* \in \widehat{X}$ folgende Eigenschaften gelten: 
\begin{gather}
x_* + ( y_* + z_* ) = ( x_* + y_* ) + z_*, \label{FA.2.13.2.1} \\
0_X + x_* = x_* \, \wedge \, \, x_* + (-1) \, x_* = 0_X, \label{FA.2.13.2.2} \\
x_* + y_* = y_* + x_*, \label{FA.2.13.2.3} \\
1 \, x_* = x_*, \label{FA.2.13.2.4} \\
(\lambda \, \mu) \, x_* = \lambda \, ( \mu \, x_*), \label{FA.2.13.2.5} \\
\lambda \, ( x_* + y_* ) = \lambda \, x_* + \lambda \, y_*, \label{FA.2.13.2.6} \\
(\lambda + \mu) \, x_* = \lambda \, x_* + \mu \, x_*. \label{FA.2.13.2.7}
\end{gather}
Sämtliche dieser Aussagen lassen sich analog zum Beweis von \ref{FA.2.12} mittels Stetigkeitsargumenten nachweisen.
(\ref{FA.2.13.2.4}) sowie (\ref{FA.2.13.2.6}) ergeben sich bereits aus der $\K$-Linearität der Abbildung $\overline{S_{\lambda}}$.
Wir zeigen jetzt exemplarisch (\ref{FA.2.13.2.3}) und (\ref{FA.2.13.2.7}):

Seien $(x_n)_{n \in \N}$ sowie $(y_n)_{n \in \N}$ Folgen in $X$, die gegen $x_*$ bzw.\ $y_*$ konvergieren.
Dann gilt aus Stetigkeitsgründen
$$ x_* + y_* = \lim_{n \to \infty} A(x_n,y_n) = A( \lim_{n \to \infty} y_n, \lim_{n \to \infty} x_n) = y_* + x_*. $$

Sind zusätzlich $\lambda, \mu \in \K$, so folgt des weiteren
\begin{eqnarray*}
(\lambda + \mu) \, x_* & = & \overline{S_{\lambda + \mu}} (\lim_{n \to \infty} x_n) = \lim_{n \to \infty} S_{\lambda + \mu} (x_n) = \lim_{n \to \infty} (S_{\lambda} \circ S_{\mu}) (x_n) \\
& = & \lim_{n \to \infty} (\lambda \, \mu) (x_n) = \lim_{n \to \infty} (\lambda \, x_n + \mu \, x_n) = \lim_{n \to \infty} S_{\lambda}(x_n) + \lim_{n \to \infty} S_{\mu}(x_n) \\
& = & \overline{S_\lambda}(x_*) + \overline{S_{\mu}}(y_*) = \lambda \, x_* + \mu \, y_*,
\end{eqnarray*}
womit die Behauptung gezeigt ist. {]}

Wir weisen nun die Existenz einer kanonischen Norm nach und nutzen aus, daß sich die gleichmäßig stetige Abbildung $\| \ldots \| \: X \to \R$, vgl.\ \ref{FA.2.2} (iii), nach \ref{FA.1.20} zu einer einer stetigen Abbildung $\widehat{\| \ldots \|} \: \widehat{X} \to \R$ erweitern läßt.

Seien $x_* \in \widehat{X}$ und $(x_n)_{n \in \N}$ eine Folge in $X$ mit $\lim_{n \to \infty} x_n = x_*$.
Dann ist aus Stetigkeitsgründen direkt klar, daß $\widehat{\| x_* \|} = \lim_{n \to \infty} \|x_n\| \ge 0$ folgt.
Zu zeigen bleiben (N1) - (N3):

Zu (N1): 
Es gilt
$$ x_* = 0 = \lim_{n \to \infty} 0 \stackrel{\widehat{\| \ldots \|} \, \text{stetig}}{\Longrightarrow} \widehat{\| x_* \|} = \lim_{n \to \infty} \| 0 \| = 0 $$
und
$$ \widehat{\| x_* \|} = 0 \Longleftrightarrow \widehat{\| \lim_{n \to \infty} x_n \|} = 0 \stackrel{\widehat{ \| \ldots \|} \, \text{stetig}}{\Longrightarrow} \underbrace{\lim_{n \to \infty} \|x_n\|}_{= \|x\|} = 0. $$

Zu (N2): Sei $\lambda \in \K$.
Dann gilt
$$ \widehat{\| \lambda \, x_* \|} = \| \lambda \, \lim_{n \to \infty} x_n \| = |\lambda| \, \| \lim_{n \to \infty} x_n \| = |\lambda| \, \widehat{\| x_* \|} $$

Zu (N3): Es seien $x_{*i} \in \widehat{X}$ und $(x_{i,n})_{n \in \N}$ Folgen in $X$ mit $\lim_{n \to \infty} x_{i,n} = x_{*i}$ für $i \in \{1,2\}$.
Es ergibt sich
$$ \widehat{\|x_1 + x_2\|} = \lim_{n \to \infty}\| x_{1,n} + x_{2,n}\| \le \lim_{n \to \infty} \|x_{1,n} \| + \lim_{n \to \infty} \|x_{2,n} \| = \widehat{\|x_1\|} + \widehat{\|x_2\|}. $$

Zu (ii): Der Beweis wird analog zu dem von \ref{FA.1.21} (iv) geführt. 
Im jetzigen Falle sind die Inklusionen isometrische Abbildungen normierter $\K$-Vektorräume, und eine stetige Abbildung läßt sich nach \ref{FA.2.12} stetig auf eine Vervollständigung fortsetzen.
Die genaue Ausführung des Beweises sei dem Leser als Übung überlassen.
\q

\begin{Kor}[Vervollständigung normierter $\K$-Algebren] \label{FA.2.14} \index{Vervollständigung!normierter!Algebren} $\,$

\noindent \textbf{Vor.:} Sei $(A, \| \ldots \|)$ eine normierte $\K$-Algebra.

\noindent \textbf{Beh.:} Es existiert eine $\K$-Banachalgebra $(\widehat{A},\widehat{\| \ldots \|})$ mit folgenden Eigenschaften:
\begin{itemize}
\item[(i)] $A$ ist $\K$-Unteralgebra von $\widehat{A}$ mit $\overline{A} = \widehat{A}$ und $\| \ldots \| = \widehat{\| \ldots \|}|_A$.
\item[(ii)] Die $\K$-Banachalgebra $(\widehat{A},\widehat{\| \ldots \|})$ mit (i) ist bis auf Isometrie normierter $\K$-Algebren eindeutig bestimmt, daher nennt man $(\widehat{A},\widehat{\| \ldots \|})$ ,,die`` \emph{Vervollständigung von $(A,\| \ldots \|)$}. 
\end{itemize}
\end{Kor}

\textit{Beweis.} (ii) zeigt man analog zu \ref{FA.2.13} (ii) unter Verwendung von \ref{FA.2.12.K}.

Zu (i): Sei $X$ der $A$ zugrundeliegende normierte $\K$-Vektorraum.
Wegen des letzten Hauptsatzes existiert ein $\K$-Banachraum $\widehat{X}$, der $X$ als normierten $\K$-Vektorraum vervollständigt.
Nach \ref{FA.2.11.K1} ist $\mathcal{L}_{\K}(\widehat{X},\widehat{X})$ mit der Operatornorm eine $\K$-Banach\-algebra, die $A$ unter der Identifikation $A \equiv L(A)$, vgl.\ \ref{FA.2.10.M}, als $\K$-Unteralgebra enthält, und die Einschränkung der Operatornorm auf $A \equiv L(A)$ stimmt mit der Norm von $A \equiv L(A)$ überein.
Wegen \ref{FA.2.4} (ii) ist $\overline{L(A)}$ ebenfalls eine $\K$-Unteralgebra der $\K$-Banachalgebra $\mathcal{L}_{\K}(\widehat{X},\widehat{X})$, also nach \ref{FA.1.11} (i) bzgl.\ der Operatornorm vollständig.
Offenbar ist $\overline{L(A)}$ eine Vervollständigung von $A \equiv L(A)$. \q 

\subsection*{Endlich-dimensionale Räume} \addcontentsline{toc}{subsection}{Endlich-dimensionale Räume}

\begin{Def}[Äquivalenz von Normen] \label{FA.E.1} 
Seien $X$ ein $\K$-Vektorraum und $\| \ldots \|$ sowie $\| \ldots \|_*$ Normen auf $X$. 

$\| \ldots \|, \| \ldots \|_*$ heißen \emph{äquivalent}\index{Norm! Äquivalenz von --en} (i.Z. $\| \ldots \| \sim \| \ldots \|_*$) genau dann, wenn gilt ${\rm Top}(X,\|\ldots\|) = {\rm Top}(X,\|\ldots\|_*)$, vgl.\ \ref{FA.2.2} (i).
\end{Def}

\begin{Bem*}
Sei $X$ ein $\K$-Vektorraum.
Dann ist $\sim$ offenbar eine Äquivalenzrelation in der Menge aller Normen auf $X$.
\end{Bem*}

\begin{Satz} \label{FA.2.E.2}
Seien $X$ ein $\K$-Vektorraumm und $\| \ldots\|$ sowie $\| \ldots \|_*$ Normen auf $X$.
Dann sind die folgenden Aussagen paarweise äquivalent:

\begin{itemize}
\item[(i)] $\| \ldots \| \sim \| \ldots \|_*$.
\item[(ii)] $\id_X \: (X, \| \ldots\|) \to (X, \| \ldots \|_*)$ ist ein Homöomorphismus.
\item[(iii)] Es existieren $C,D \in \R_+$ mit $\forall_{x \in X} \, C \, \|x\| \le \|x\|_* \le D \, \|x\|$.
\end{itemize}
\end{Satz}

\textit{Beweis als Übung.} \q

\begin{Kor} \label{FA.2.E.3} 
Seien $X$ ein $\K$-Vektorraum und $\| \ldots\|, \| \ldots \|_*$ äquivalente Normen auf $X$.

Dann sind die durch $\| \ldots\|$ und $\| \ldots \|_*$ induzierten Metriken Cauchy-äqui\-va\-lent; insbes.\ ist $(X, \| \ldots\|)$ genau dann ein $\K$-Banachraum, wenn $(X, \| \ldots\|_*)$ ein solcher ist. \q
\end{Kor}

\begin{Satz} \label{FA.2.E.4}
Sei $X$ ein endlich-dimensionaler $\K$-Vektorraum.

Dann sind je zwei Normen auf $X$ äquivalent.
\end{Satz}

Wir bereiten den Beweis des Satzes durch ein Lemma vor, dessen Formulierung zunächst folgende Definition benötigt.

\begin{Def}[Beschränkte Mengen und Folgen] \label{FA.2.E.5} \index{Menge!beschränkte} \index{Folgen!beschränkte}
Sei $X$ ein normierter $\K$-Vektor\-raum.
\begin{itemize}
\item[(i)] Nach der Bemerkung zu \ref{FA.1.32} ist eine Teilmenge $B$ von $X$ genau dann \emph{beschränkt}, wenn gilt $\exists_{C \in \R_+} \forall_{x \in B} \, \|x\| \le C$.
\item[(ii)] Eine Folge in $X$ heißt per definitionem \emph{beschränkt}, wenn die Menge ihrer Folgenglieder beschränkt ist.
\end{itemize}
\end{Def}

Das folgende Lemma zeige der Leser als Übung.
Beachte, daß der Fall $n=1$ nach dem Satz von \textsc{Bolzano-Weierstraß} der Analysis I bereits klar ist.

\begin{Lemma}[Satz von \textsc{Bolzano-Weierstraß} in $(\R^n,\| \ldots \|_{\infty})$] \label{FA.2.E.6}
Jede beschränkte Folge in $(\R^n,\| \ldots \|_{\infty})$ besitzt eine konvergente Teilfolge. \q
\end{Lemma} 

\textit{Beweis des Satzes.} Zunächst ist jeder normierte $\C$-Vektorraum in kanonischer Weise ein normierter $\R$-Vektorraum, dessen Topologie mit der des $\C$-Vektorraumes übereinstimmt.
Daher können wir ohne Einschränkung $\K = \R$ annehmen.
Des weiteren genügt es offenbar den Fall $n := \dim_{\R} X \in \N_+$ zu betrachten.
Wir wählen eine Basis $\{b_1, \ldots, b_n\}$ von $X$.
Bezeichnet dann $\{e_1, \ldots, e_n\}$ die kanonische Basis von $\R^n$, so ist
$$ \R^n \longrightarrow X, ~~ \sum_{i=1}^n \lambda_i \, e_i \longmapsto \sum_{i=1}^n \lambda_i \, b_i, $$
ein $\R$-Vektorraum-Isomorphismus, unter dem den Normen auf $\R^n$ umkehrbar eindeutig die Normen auf $X$ entsprechen, wobei äquivalenten Normen auf $\R^n$ äquivalente Normen auf $X$ entsprechen.
Daher können wir weiterhin annehmen, daß $X = \R^n$ gilt.

Sei nun $\| \ldots \|$ eine Norm auf $\R^n$.
Da $\sim$ eine Äquivalentrelation ist, genügt es zu zeigen, daß $\| \ldots \|$ zu $\| \ldots \|_{\infty}$ äquivalent ist, d.h.\ nach \ref{FA.2.E.2} ,,(iii) $\Rightarrow$ (i)``
$$ \exists_{C,D \in \R_+} \forall_{x \in \R^n} \, C \, \|x\|_{\infty} \le \|x\| \le D \, \|x\|_{\infty}. $$

Beweis hiervon:
1.) Für alle $x = \sum_{i=1}^n \lambda_i \, e_i \in \R^n$ gilt mit $D := \sum_{i=1}^n \| e_i \| \in \R_+$
$$\|x\| \le \sum_{i=1}^n |\lambda_i| \, \|e_i\| \le \|x\|_{\infty} \sum_{i=1}^n \| e_i \| =  D \, \| x \|_{\infty}.$$

2.) Wir setzen
\begin{equation} \label{FA.2.E.4.1}
C := \inf \{ \|x\| \, | \, x \in \R^n \, \wedge \, \| x \|_{\infty} = 1 \}
\end{equation}
und zeigen zunächst, daß gilt
\begin{equation} \label{FA.2.E.4.2}
C > 0.
\end{equation}

{[} Zu (\ref{FA.2.E.4.2}): Nach Definition von $C$ existiert offenbar eine Folge $\left( x_k \right)_{k \in \N}$ in $\R^n$ mit
\begin{equation} \label{FA.2.E.4.3}
\forall_{k \in \N} \, \| x_k \|_{\infty} = 1 ~ \mbox{ und } ~ C = \lim_{k \to \infty} \| x_k \|.
\end{equation}
Als beschränkte Folge in $\R^n$ bzgl.\ $\| \ldots \|_{\infty}$ besitzt $\left( x_k \right)_{k \in \N}$ nach \ref{FA.2.E.6} eine in $(\R^n,\| \ldots \|_{\infty})$ konvergente Teilfolge.
Wir können ohne Beschränkung der Allgemeinheit annehmen, daß $x \in \R^n$ mit
\begin{equation} \label{FA.2.E.4.4}
\lim_{k \to \infty} x_k = x \mbox{ in } (\R^n, \| \ldots \|_{\infty})
\end{equation}
existiert.
Da wir 1.) bereits gezeigt haben, folgt dann auch
\begin{equation} \label{FA.2.E.4.5}
\lim_{k \to \infty} x_k = x \mbox{ in } (\R^n, \| \ldots \|).
\end{equation}
Nach \ref{FA.2.2} (iii) sind $\| \ldots \|_{\infty} \: (\R^n, \| \ldots \|_{\infty}) \to \R$ und $\| \ldots \| \: (\R^n, \| \ldots \|) \to \R$ stetig, also folgt zunächst $\|x\|_{\infty} = \lim_{k \to \infty} \|x_k\|_{\infty} \stackrel{(\ref{FA.2.E.4.3})}{=} 1$, insbesondere also $x \ne 0$, und sodann aus (\ref{FA.2.E.4.5})
$$ 0 < \| x \| = \lim_{k \to \infty} \| x_k \| \stackrel{(\ref{FA.2.E.4.3})}{=} C. $$
Damit ist (\ref{FA.2.E.4.2}) gezeigt. {]}

Sei nun $x \in \R^n \setminus \{0\}$ beliebig.
Dann folgt aus (\ref{FA.2.E.4.1}) wegen $\|\frac{x}{\|x\|_{\infty}}\|_{\infty} = 1$
$$ C \le \left\| \frac{x}{\|x\|_{\infty}} \right\| = \frac{\|x\|}{\|x\|_{\infty}}, $$
also ist die Behauptung wegen (\ref{FA.2.E.4.2}) offenbar gezeigt. \q

\begin{Satz} \label{FA.2.E.7}
Sei $X$ ein endlich-dimensionaler normierter $\K$-Vektorraum.

Dann ist $X$ ein $\K$-Banachraum.
\end{Satz}

\textit{Beweis.} Wir bemerken zunächst, daß die Begriffe ,,Cauchyfolge`` und ,,Konvergenz`` von der Norm nach \ref{FA.2.E.2} ,,(i) $\Leftrightarrow$ (iii)`` nur bis auf Äquivalenz abhängen.
Dieselbe Begründung wie im Beweis des letzten Satzes zeigt nun, daß es genügt den Fall $X = \R^n$ mit $n \in \N_+$ zu betrachten.
Der letzte Satz \ref{FA.2.E.4} zeigt des weiteren, daß wir annehmen können, daß die Norm die Maximumsnorm auf $\R^n$ ist.
Zeige unter Ausnutzung der Vollständigkeit von $\R$ als Übung, daß $(\R^n,\| \ldots \|_{\infty})$ ein $\R$-Banachraum ist.
Damit ist der Satz dann bewiesen. \q 

\begin{Kor} \label{FA.2.E.8}
Seien $X$ ein normierter $\K$-Vektorraum und $Y$ ein endlich-dimensionaler $\K$-Untervektorraum von $X$.
Dann gilt:
\begin{itemize}
\item[(i)] $Y$ ist abgeschlossen in $X$.
\item[(ii)] Ist $Z$ ein abgeschlossener $\K$-Untervektorraum von $X$, so ist auch $Y+Z$ ein solcher.%
\footnote{\textbf{Definition.} Sind $X$ ein beliebiger $\K$-Vektorraum und $Y_1, Y_2$ $\K$-Untervektorräume von $X$, so heißt der $\K$-Untervektorraum $\boxed{Y_1 + Y_2} := \{ y_1 + y_2 \, | \, y_1 \in Y_1 \, \wedge y_2 \in Y_2 \}$ von $X$ die \emph{Summe von $Y_1$ und $Y_2$}. 
Diese Summe heißt \emph{direkt} (i.Z. \boxed{Y_1 \oplus Y_2}) genau dann, wenn gilt $Y_1 \cap Y_2 = \{0\}$.
Letzteres ist offenbar gleichbedeutend damit, daß zu jedem $y \in Y_1 + Y_2$ eindeutig bestimmte $y_1 \in Y_1$ und $y_2 \in Y_2$ mit $y = y_1 + y_2$ existieren.}
\end{itemize}
\end{Kor}

\textit{Beweis.} (i) ergibt sich sofort aus \ref{FA.2.E.7} und \ref{FA.1.11}. 

Zu (ii): Nach \ref{FA.2.5} (i), (ii) ist der kanonische Homomorphismus $\pi \: X \to X/Z$ eine stetige Abbildung zwischen normierten $\K$-Vektorräumen.
$\pi(Y)$ ist ein end\-lich-dimensionaler $\K$-Untervektorraum von $X/Z$, also nach (i) abgeschlossen in $X/Z$.
Wegen der Stetigkeit von $\pi$ ist dann auch $Y+Z = \overline{\pi}^1(\pi(U))$ abgeschlossen in $X$. \q

\begin{Satz} \label{FA.2.E.9}
Seien $X$ ein $\K$-Banachraum, $I$ eine höchstens abzählbare Menge und $\mathfrak{B} = \{ b_i \, | \, i \in I \}$ eine Basis von $X$, vgl.\ \ref{FA.A.1} und beachte \ref{FA.A.2}.

Dann ist $\mathfrak{B}$ endlich, d.h.\ $\dim_{\K} X < \infty$.
\end{Satz}

\textit{Beweis.} Angenommen, es gilt $\# I = \# \N$.
Sei dann ohne Beschränkung der Allgemeinheit $I = \N$.
Nach \ref{FA.2.E.8} (i) ist für jedes $i \in \N$
$$ A_i := {\rm Span}_{\K} \{b_1, \ldots, b_{i+1}\} $$
ein abgeschlossener $\K$-Untervektorraum des $\K$-Banachraumes $X$, wobei letzterer nach dem Satz von Baire \ref{FA.1.44} Bairesch ist, vgl.\ \ref{FA.1.42}.
Wegen $\bigcup_{i \in \N} A_i = X$ und $\dim_{\K} X > 0$ besitzt $\bigcup_{i \in \N} A_i$ einen inneren Punkt, also existiert $i_0 \in \N$ derart, daß $(A_{i_0})^{\circ} \ne \emptyset$, und wir zeigen
\begin{equation} \label{FA.2.E.9.1}
A_{i_0} = X,
\end{equation}
womit dann ein Widerspruch hergeleitet ist.

Zu (\ref{FA.2.E.9.1}): Da ,,$\subset$`` ist trivial ist, haben wir nur ,,$\supset$`` zu zeigen.
Wie oben eingesehen existiert $a \in (A_{i_0})^{\circ}$, d.h.\ auch, daß wir eine Zahl $\varepsilon \in \R_+$ mit
$$ U_{\varepsilon}(a) \subset A_{i_0} $$
finden können.
Sei nun $x \in X$ mit $x \ne a$ beliebig.
Dann gilt
$$ y := a + \frac{\varepsilon}{2} \, \frac{x-a}{\|x-a\|} \in U_{\varepsilon}(a) \subset A_{i_0}, $$
also $x = \underbrace{a}_{\in A_{i_0}} + \frac{2 \|x - a\|}{\varepsilon} \, \underbrace{(y-a)}_{\in A_{i_0}} \in A_{i_0}$. \q

\begin{Satz} \label{FA.2.E.10}
Seien $X$ ein endlich-dimensionaler normierter $\K$-Vektorraum und $Y$ ein normierter $\K$-Vektorraum.

Dann ist jede $\K$-lineare Abbildung $T \: X \to Y$ stetig.
\end{Satz}

\textit{Beweis.} Ohne Einschränkung sei $n := \dim_{\K} X \in \N_+$.
Wir wählen eine Basis $\mathfrak{B} := \{b_1, \ldots, b_n\}$ von $X$ und betrachten auf $X$ die \emph{Maximumsnorm bzgl.\ $\mathfrak{B}$}, die durch
$$ \forall_{x = \sum_{i=1}^n \lambda_i b_i \in X} \, \|x\|_{\infty}^{\mathfrak{B}} := \max \{ |\lambda_i| \, | \, i \in \{1, \ldots, n\} \} $$
gegeben ist.
Es genügt zu zeigen, daß $T \: X \to Y$ bzgl.\ dieser Norm auf $X$ stetig ist, und dies ist klar, da für alle $x = \sum_{i=1}^n \lambda_i \, b_i \in X$ mit $C := \sum_{i=1}^n \|T(b_i)\| \in \R_+$ 
$$ \|T(x)\| = \sum_{i=1}^n |\lambda_i| \, \| T(b_i) \| \le C \, \|x\|_{\infty}^{\mathfrak{B}} $$
gilt. \q

\begin{Kor} \label{FA.2.E.10.K}
Seien $X$ ein endlich-dimensionaler $\K$-Vektorraum und $r \in \N_+$.
Ferner seien auf $X^r$ eine Norm gegeben und $Y$ ein normierter $\K$-Vektorraum.

Dann ist jede $r$-fach $\K$-multilineare Abbildung $\varphi \: X^r \to Y$ stetig.
\end{Kor}

\textit{Beweisskizze.} Seien $n := \dim_{\K} X \in \N_+$ und $\mathfrak{B} = \{b_1, \ldots, b_n\}$ eine Basis von $X$.  
Wähle eine Norm $\| \ldots \|$ auf $X$. 
Bezeichnet dann $\{b_1^*, \ldots, b_n^*\}$ die zu $\mathfrak{B}$ duale Basis von ${\rm Hom}_{\K}(X,\K)$, d.h.\ $\forall_{i,j \in \{1, \ldots, n\}} \, b_i^*(b_j) = \delta_{ij}$, so sind $b_1^*, \ldots, b_n^*$ nach dem letzten Satz stetig.
Betrachte auf $X^r$ die durch $\| \ldots \|$ induzierte Maximumsnorm.
Dann ist für jedes $k \in \{1, \ldots, r\}$ die Projektion $\pi_k \: X^r \to X$ auf die $k$-te Komponente stetig, und es gilt
$$ \varphi = \sum_{i_1,\ldots,i_r=1}^n (b_{i_1}^* \circ \pi_1) \, \cdots \, (b_{i_r}^* \circ \pi_r) \, \varphi(b_{i_1}, \ldots, b_{i_r}), $$
also ist auch $\varphi$ stetig. \q

\begin{HS}[Charakterisierung endlich-dimensionaler normierter Vektorräume] \label{FA.2.E.11} $\,$

\noindent \textbf{Vor.:} Sei $X$ ein normierter $\K$-Vektorraum.

\noindent \textbf{Beh.:} Die folgenden drei Aussagen sind paarweise äquivalent:
\begin{itemize}
\item[(i)] $\dim_{\K} X < \infty$.
\item[(ii)] In $X$ gilt der \emph{Satz von \textsc{Heine-Borel}}\index{Satz!von \textsc{Heine-Borel}}, welcher besagt, daß für alle Teilmengen $K$ von $X$ gilt
\begin{equation} \label{FA.2.E.11.0}
K \mbox{ kompakt } \Longleftrightarrow  \mbox{ $K$ beschränkt und abgeschlossen}.
\end{equation}
\item[(iii)] Die Einheitsvollkugel $B_1(0)$ ist kompakt.
\end{itemize}
\end{HS}
\pagebreak

\textit{Beweis.} ,,(i) $\Rightarrow$ (ii)`` Sei $X$ endlich-dimensional und sei $K$ eine Teilmenge von $X$.
Zu zeigen ist (\ref{FA.2.E.11.0}):

,,$\Rightarrow$`` ist klar nach \ref{FA.1.24} (ii) und \ref{FA.1.33}.

,,$\Leftarrow$`` Wie oben können wir annehmen, daß gilt $X = (\R^n,\|\ldots\|_{\infty})$ mit $n \in \N_+$.
Sei $K$ beschränkt und abgeschlossen.
Nach \ref{FA.1.30} ,,(i) $\Leftrightarrow$ (ii)`` genügt es zu zeigen, daß $K$ folgenkompakt ist.
Ist eine Folge in $K$ gegeben, so besitzt diese wegen \ref{FA.2.E.6} eine in $(\R^n,\|\ldots\|_{\infty})$ konvergente Teilfolge, deren Folgenglieder in $K$ liegen.
Da $K$ abgeschlossen ist, ist der Grenzwert der Teilfolge nach \ref{FA.1.KA} ein Element von $K$.
Damit ist die Folgenkompaktheit von $K$ gezeigt.

,,(ii) $\Rightarrow$ (iii)`` $B_1(0)$ ist beschränkt und abgeschlossen, also nach (ii) kompakt.

,,(iii) $\Rightarrow$ (i)`` Sei $B_1(0)$ kompakt.
Nach \ref{FA.1.29} (i) ist $B_1(0)$ präkompakt, d.h.\ es existieren insbes.\ $y_1, \ldots, y_k \in B_1(0)$ mit 
\begin{equation} \label{FA.2.E.11.1}
B_1(0) \subset \bigcup_{i=1}^k U_{\frac{1}{2}}(y_i).
\end{equation}

Wir definieren einen endlich-dimensionalen $\K$-Untervektorraum $Y$ von $X$ durch
$$ Y := {\rm Span}_{\K} \{y_1, \ldots, y_k\} $$
und behaupten
\begin{equation} \label{FA.2.E.11.2}
Y = X,
\end{equation}
womit (i) bewiesen ist.

Angenommen, (\ref{FA.2.E.11.2}) ist falsch, d.h.\ $Y \subsetneqq X$, also existiert $\xi \in X/Y \setminus \{0\}$.
Da $Y$ endlich-dimensional ist, ist $Y$ nach \ref{FA.2.E.8} (i) abgeschlossen in $X$.
Wegen \ref{FA.2.5} (i) ist dann $X/Y$ in kanonischer Weise ein normierter $\K$-Vektorraum, und es gilt
$$ \alpha := \|\xi\| = \inf \{ \|x\| \, | \, x \in \xi \} \in \R_+. $$
Daher existiert $x \in \xi$, d.h.\ $\xi = x + Y$, mit
\begin{equation} \label{FA.2.E.11.3}
0 < \alpha \le \|x\| < \frac{3}{2} \alpha.
\end{equation}
Es gilt $\frac{x}{\|x\|} \in B_1(0)$, also finden wir nach (\ref{FA.2.E.11.1}) ein $i \in \{1, \ldots, k\}$ derart, daß gilt
\begin{equation} \label{FA.2.E.11.4}
\frac{x}{\|x\|} \in U_{\frac{1}{2}}(y_i).
\end{equation}
Nun folgt zunächst aus der Definition von $Y$
$$ \|x\| \left( \frac{x}{\|x\|} - y_i \right) = x - \|x\| \, y_i \in x + Y = \xi $$
und sodann aus der Definition von $\alpha$ sowie (\ref{FA.2.E.11.3}), (\ref{FA.2.E.11.4})
$$ \alpha \le \|x\| \left\| \frac{x}{\|x\|} - y_i \right\| < \frac{3}{4} \alpha, $$
Widerspruch. \q 

\subsection*{Existenz von Normen} \addcontentsline{toc}{subsection}{Existenz von Normen}

Wir wollen eine Charakterisierung dafür angeben, daß ein $\K$-Vektorraum eine Norm besitzt, letztere hatten wir in \ref{FA.2.1} (i) als eine Abbildung $\| \ldots \| \: X \to {[}0, \infty{[}$ mit (N1) - (N3) definiert.
Tatsächlich folgt aus (N2) und (N3) bereits die Nicht-Negativität:

\begin{Lemma} \label{FA.N.1}
Seien $X$ ein $\K$-Vektorraum und $p \: X \to \R$ eine Funktion mit
$$ \forall_{x,y \in X} \forall_{\lambda \in \K} \, p(\lambda \, x) = |\lambda| \, p(x) ~~ \mbox{ und } ~~ p(x+y) \le p(x) + p(y). $$
\begin{itemize}
\item[(i)] Es folgt $p \ge 0$, d.h.\ $p$ ist eine Halbnorm.
\item[(ii)] Gilt zusätzlich
$$ \forall_{x \in X} \, p(x) = 0 \, \Longrightarrow \, x=0, $$
so ist $p$ eine Norm.
\end{itemize}
\end{Lemma}

\textit{Beweis.} Für alle $x,y \in X$ gilt $p(x) = p(x-y + y) \le p(x-y) + p(y)$, d.h.\ $p(x) - p(y) \le p(x-y)$ und ebenso $p(y) - p(x) \le p(y-x)$.
Aus $p(x-y) = p(y-x)$ folgt daher $|p(x) - p(y)| \le p(x-y)$, also für $y=0$: $p(x) \ge 0$.
Damit ist (i) gezeigt, und (ii) folgt trivial. \q 

\begin{Def}[Absorbierende, konvexe und ausgewogene Mengen sowie das Minkowski-Funktional einer absorbierenden Menge] \label{FA.N.2}
Seien $X$ ein $\K$-Vektorraum und $A,C$ Teilmengen von $X$.
Wir definieren dann:

\begin{itemize}
\item[(i)] $A$ \emph{absorbierend}\index{Menge!absorbierende} $: \Longleftrightarrow$ $X = \bigcup_{\tau \in \R_+} \tau \, A$.

\begin{Bem*} $\,$
\begin{itemize}
\item[1.)] Eine absorbierende Menge enthält notwendigerweise das Nullelement.
\item[2.)] $A$ ist z.B.\ absorbierend, wenn eine Norm auf $X$ existiert bzgl.\ derer $0 \in A^{\circ}$ gilt.
\end{itemize}
\end{Bem*}
\item[(ii)] Ist $A$ eine absorbierende Menge, so heißt die Funktion $p_A \: X \to \R$, die durch
$$ \forall_{x \in X} \, \boxed{p_A(x)} := \inf \{ \tau \in \R_+ \, | \, x \in \tau \, A \} $$
definiert ist, das \emph{Minkowski-Funktional von $A$}\index{Minkowski!-Funktional}.
\item[(iii)] Für alle $x,y \in X$ bezeichnen wir die \emph{Verbindungsstrecke von x und y} mit
$$ \boxed{[x,y]} := \{ \underbrace{(1-t) \, x + t \, y}_{= x + t \, (y-x)} \, | \, t \in [0,1] \}. $$
\item[(iv)] $C$ \emph{konvex}\index{Menge!konvexe} $: \Longleftrightarrow$ $\forall_{x,y \in C} \, [x,y] \subset C$.

\begin{Bsp*} $\,$
\begin{itemize}
\item[1.)] $\emptyset$, jede Verbindungsstrecke und jeder $\K$-Untervektorraum von $X$ ist konvex.
Insbesondere ist $X$ konvex.
\item[2.)] Ist $C \subset X$ konvex, so ist für jedes $x_0 \in X$ auch $x_0 + C$ konvex, denn für alle $x,y \in C$ und jedes $t \in [0,1]$ gilt
$$ (1-t) \, (x_0 + x) + t \, (x_0 + y) = x_0 + \underbrace{(1-t) \, x + t \, y}_{\in C}. $$
\item[3.)] Seien $C, \widetilde{C}$ konvexe Teilmengen von $X$ und $\lambda \in \K$.
Dann sind auch $C + \widetilde{C}$ und $\lambda \, C$ konvexe Teilmengen von $X$.%
\footnote{\textbf{Definition.} Für beliebige Teilmengen $A,B$ eines $\K$-Vektorraumes und $\lambda \in \K$ setzt man $\boxed{A+B} := \{a+b \, | \, a \in A \, \wedge \, b \in B\}$ und $\boxed{\lambda \, A} := \{ \lambda \, a \, | \, a \in A \}$. 
(Beachte, daß $2 \, A \subsetneqq A + A$ gelten kann.)
Gilt $A = \{a\}$, so schreibt man auch $\boxed{a+B}$ anstelle von $\{a\} + B$.} %
Sind nämlich $x,y \in C$, $\tilde{x},\tilde{y} \in \widetilde{C}$ und $t \in [0,1]$, so gilt
\begin{gather*}
(1-t) \, (x+\tilde{x}) + t \, (y+\tilde{y}) = (\overbrace{(1-t) \, x + t \, y}^{\in C}) + (\overbrace{(1-t) \, \tilde{x} + t \, \tilde{y}}^{\in \widetilde{C}}), \\
(1-t) \, (\lambda \, x) + t \, (\lambda \, x) = \lambda \, (\underbrace{(1-t) \, x + t \, y}_{\in C}).
\end{gather*}
\item[4.)] Der Schnitt beliebig vieler konvexer Teilmengen von $X$ ist wieder konvex. (Daher existiert zu $A$ eine kleinste konvexe Menge $\boxed{C(A)}$, die $A$ enthält, die sog.\ \emph{konvexe Hülle von $A$}\index{konvexe Hülle}\index{Hülle!konvexe}.)
\item[5.)] Im Falle eines normierten $\K$-Vektorraumes $X$ ist mit $C$ auch $\overline{C}$ eine konvexe Teilmenge von $X$, da zu $x,y \in \overline{C}$ Folgen $(x_n)_{n \in \N}$ und $(y_n)_{n \in \N}$ in $C$ mit $\lim_{n \to \infty} x_n = x$ sowie $\lim_{n \to \infty} y_n = y$ existieren, d.h.\ für jedes $t \in [0,1]$ gilt
$$ (1-t) \, x + t \, y = \lim_{n \to \infty} \underbrace{(1-t) \, x + t \, y}_{\in C} \in \overline{C}. $$
\item[6.)] Im Falle eines normierten $\K$-Vektorraumes ist mit $C$ auch $C^{\circ}$ konvex, vgl.\ \ref{FA.N.9.L}.
\end{itemize}
\end{Bsp*}
\item[(v)] $A$ \emph{ausgewogen}\index{Menge!ausgewogene} $: \Longleftrightarrow$ $\forall_{\lambda \in \K, \, |\lambda| \le 1} \, \lambda \, A \subset A$.
\end{itemize}
\end{Def}

\begin{Def}[Sublineare Funktionen] \label{FA.2.16}
Sei $X$ ein $\K$-Vektorraum.
Eine Funktion $p \: X \to \R$ heißt \emph{sublinear}\index{Funktion!sublineare} genau dann, wenn für alle $t \in \R$ mit $t \ge 0$ und alle $x,y \in X$ gilt
$$ p( t \, x ) = t \, p(x) ~~ \wedge ~~ p(x + y) \le p(x) + p(y). $$

\begin{Bsp*}
Lineare Abbildungen $X \to \K$ und Halbnormen auf $X$ sind sublinear.
\end{Bsp*}
\end{Def}

\begin{Satz} \label{FA.N.3}
Seien $X$ ein $\K$-Vektorraum, $p \: X \to \R$ eine sublineare Funktion und $\alpha \in \R_+$.
Dann gilt:
\begin{itemize}
\item[(i)] $\overline{p}^1({]}- \infty, \alpha{[})$ und $\overline{p}^1({]}- \infty, \alpha{]})$ sind absorbierend sowie konvex.
Die Konvexität ist sogar für $\alpha \in \R$ gegeben.
\item[(ii)] Ist $p \: X \to \R$ sogar eine Halbnorm, so sind $\overline{p}^1({[}0, \alpha{[})$ und $\overline{p}^1({[}0, \alpha{]})$ ausgewogen.
\end{itemize}
\end{Satz}

\textit{Beweis.} Zu (i): 1.) Sei $x \in X$.
Dann existiert wegen $\alpha \in \R_+$ eine Zahl $\tau \in \R_+$ mit $p(x) < \tau \cdot \alpha$.
Aus der Sublinearität von $p$ folgt $p(\frac{x}{\tau}) < \alpha$, also ergibt sich $\frac{x}{\tau} \in \overline{p}^1({]}- \infty, \alpha{[}) \subset \overline{p}^1({]}- \infty, \alpha{]})$.
Somit gilt $x \in \tau \, \overline{p}^1({]}- \infty, \alpha{[}) \subset \tau \, \overline{p}^1({]}- \infty, \alpha{]})$.

2.) Seien $\alpha \in \R$ und $x,y \in X$ mit $p(x)<\alpha$ und $p(y)<\alpha$.
Dann gilt für jedes $t \in {[}0,1{[}$
$$ p( (1-t) \, x + t \, y ) \le p((1-t) \, x) + p(t \, x) = (1-t) \, p(x) + t \, p(y) < \alpha, $$
womit gezeigt ist, daß $\overline{p}^1({]}- \infty, \alpha{[})$ konvex ist.
Die Aussage für $\overline{p}^1({]}- \infty, \alpha{]})$ sieht man analog ein.

Zu (ii): Für jedes $\lambda \in \K$ mit $|\lambda| \le 1$ und alle $x \in X$ gilt
$$ p( \lambda \, x ) = |\lambda| \, p(x) \le p(x), $$
also folgt aus $p(x) < \alpha$ bzw.\ $p(x) \le \alpha$ auch $p(\lambda \, x) < \alpha$ bzw.\ $p(\lambda \, x) \le \alpha$. \q

\begin{Bem*}
Der letzte Satz zeigt, daß im Falle eines normierten $\K$-Vektor\-raum\-es die Mengen $U_{\varepsilon}(X)$ und $B_{\varepsilon}(x)$ für jedes $x \in X$ und $\varepsilon \in \R_+$ absorbierend, konvex und ausgewogen sind.
\end{Bem*}

\begin{Satz} \index{Minkowski!-Funktional} \label{FA.N.4}
Seien $X$ ein $\K$-Vektorraum und $C$ eine absorbierende konvexe Teilmenge von $X$.
Dann folgt:
\begin{itemize}
\item[(i)] Das Minkowski-Funktional $p_C$ ist eine sublineare Funktion mit $p_C \ge 0$.
\item[(ii)] $\overline{p_C}^1({[}0,1{[}) \subset C$.
\item[(iii)] Ist zusätzlich
\begin{equation} \label{FA.2.23.2}
\forall_{x \in C} \exists_{\tau \in {]}1, \infty{[}} \, \tau \, x \in C
\end{equation}
erfüllt, so gilt in (ii) Gleichheit.
\item[(iv)] Unter der Voraussetzung
\begin{equation}
\forall_{\sigma \in \K, \, |\sigma|=1} \, \sigma \, C \subset C \label{FA.2.23.3}
\end{equation}
ist $p_C$ eine Halbnorm.
\item[(v)] Gilt zusätzlich zu (\ref{FA.2.23.3}) auch
\begin{equation}
\bigcap_{\tau \in \R_+} \, \tau \, C = \{0\}, \label{FA.2.23.4}
\end{equation}
so ist $p_C$ eine Norm.
\end{itemize}
\end{Satz}

\textit{Beweis.} Zu (i): $p_C \ge 0$ ist trivial. 
Zu zeigen ist also die Sublinearität von $p_C$:

1.) $p_C(0) = 0$ ist klar.
Seien $x \in X$ und $t \in \R_+$.
Dann gilt
\begin{eqnarray*}
p_C( t \, x ) & = & \inf \{ \tau \in \R_+ \, | \, t \, x \in \tau \, C \} = \inf \{ t \, \tilde{\tau} \in \R_+ \, | \, x \in \tilde{\tau} \, C \} \\
& = & t \, \inf \{ \tilde{\tau} \in \R_+ \, | \, x \in \tilde{\tau} \, C \} = t \, p_C(x).
\end{eqnarray*}

2.) Seien $x,y \in X$ und $\varepsilon \in \R_+$ beliebig.
Es existieren $\tau, \sigma \in \R_+$ mit $x \in \tau \, C$ und $p_C(x) \le \tau < p_C(x) + \varepsilon$ sowie $y \in \sigma \, C$ und $p_C(y) \le \sigma < p_C(y) + \varepsilon$.
Insbes.\ gilt $\frac{x}{\tau}, \frac{y}{\sigma} \in C$.
Dann folgt wegen $\frac{\tau}{\tau + \sigma}, \frac{\sigma}{\tau + \sigma} \in [0,1]$ mit $\frac{\tau}{\tau + \sigma} + \frac{\sigma}{\tau + \sigma} = 1$
$$ \frac{x+y}{\tau + \sigma} = \frac{\tau}{\tau + \sigma} \, \frac{x}{\tau} + \frac{\sigma}{\tau + \sigma} \, \frac{y}{\sigma} \in C $$
-- beachte, daß $C$ konvex ist --, also $x+y \in (\tau + \sigma) \, C$ bzw.\ $p_C(x+y) \le \tau + \sigma$.
Es ergibt sich
$$ p_C(x+y) \le \tau + \sigma < p_C(x) + p_C(y) + 2 \varepsilon, $$
also aus der Beliebigkeit von $\varepsilon \in \R_+$: $p_C(x + y) \le p_C(x) + p_C(y)$. 

Zu (ii): Sei $x \in X$ mit $p_C(x) < 1$.
Dann existiert $\tau \in {]}p_C(x),1{[}$ mit $x \in \tau \, C$, d.h.\ $\frac{x}{\tau} \in C$, wobei $\frac{1}{\tau} > 1$.
Da $C$ konvex ist, folgt $[0, \frac{x}{\tau}] \subset C$, insbes.\ $x \in C$.

Zu (iii): Seien $x \in C$ und $\tau \in {]}1, \infty{[}$ mit $\tau \,x \in C$.
Dann folgt $x \in \frac{1}{\tau} \, C$, wobei $0 < \frac{1}{\tau} < 1$, also auch $p_C(x) \le \frac{1}{\tau} < 1$. 

Zu (iv): Wir folgern aus (\ref{FA.2.23.3}), daß sogar gilt
\begin{equation} \label{FA.N.4.3}
\forall_{\sigma \in \K, \, |\sigma| = 1} \, \sigma \, C = C. 
\end{equation}

{[} Zu (\ref{FA.N.4.3}): Seien $x \in C$ und $\sigma \in \K$ mit $|\sigma| = 1$.
Dann gilt auch $\frac{1}{|\sigma|} = 1$, also $\frac{x}{\sigma} \in A$ nach (\ref{FA.2.23.3}), d.h.\ $x \in \sigma \, A$. {]}

Nun ergibt sich für jedes $x \in X$ und jedes $\lambda \in \K \setminus \{0\}$
$$ p_C( \lambda \, x) = p_C \left( |\lambda| \, \frac{\lambda}{|\lambda|} \, x \right) \stackrel{(i)} {=} |\lambda| \, p_C \left( \frac{\lambda}{|\lambda|} \, x \right) \stackrel{(\ref{FA.N.4.3})}{=} |\lambda| \, p_C(x). $$
Aus (i) und \ref{FA.N.1} (i) folgt daher (ii).

Zu (v): Da für jedes $x \in X$ gilt 
$$ p_C(x) = 0 \, \Longrightarrow \, \forall_{\tau \in \R_+} \, x \in \tau \, C \, \stackrel{(\ref{FA.2.23.4})}{\Longrightarrow} \, x = 0, $$
folgt (iii) aus \ref{FA.N.1} (ii).
\q

\begin{Bem*}
Seien $X$ ein $\K$-Vektorraum und $C$ eine Teilmenge von $X$.
\begin{itemize}
\item[1.)] Es ist leicht einzusehen, daß (ii) unter den Voraussetzungen des Lemmas durch $\overline{p_C}^1({[}0,1{[}) \subset C \subset \overline{p_C}^1({[}0,1{]})$ verschärft werden kann.
Auch $\overline{p_C}^1({[}0,1{[})$ und $\overline{p_C}^1({[}0,1{]})$ sind dann absorbierend sowie konvex, und es gilt $p_{\overline{p_C}^1({[}0,1{[})} = p_C = p_{\overline{p_C}^1({[}0,1{]})}$, vgl.\ \ref{FA.N.6}.
\item[2.)] (\ref{FA.2.23.2}) ist z.B.\ erfüllt, wenn eine Norm auf $X$ existiert, bzgl.\ derer $C$ offen ist, denn dann gilt für $x \in C$ zum einen $C \in \U(x,X)$ und zum anderen $\lim_{n \to \infty} \frac{n+1}{n} \, x = x$.
\item[3.)] (\ref{FA.2.23.4}) ist z.B.\ erfüllt, wenn eine Norm auf $X$ existiert, bzgl.\ derer $C$ beschränkt ist.
\end{itemize}
\end{Bem*}

\begin{Satz} \label{FA.N.5}
Sei $X$ ein $\K$-Vektorraum.

Dann entsprechen die Halbnormen auf $X$ den Minkowski-Funktionalen absorbierender konvexer ausgewogener Teilmengen von $X$.
\end{Satz}

Wir bereiten den Beweis des Satzes durch ein Lemma vor.

\begin{Lemma} \label{FA.N.6}
Seien $X$ ein $\K$-Vektorraum, $p \: X \to \R$ eine sublineare Funktion mit $p \ge 0$ und
$$ C \in \left\{ \overline{p}^1({[}0,1{[}), \, \overline{p}^1({[}0,1{]}) \right\}, $$
insbes.\ ist $C$ nach \ref{FA.N.3} (i) absorbierend sowie konvex.

Dann gilt $p = p_C$.
\end{Lemma}

\textit{Beweis.} Sei $x \in X$.
Zum einen gilt für jedes $\tau \in \R_+$ mit $\tau > p(x) (\ge 0)$: $p(\frac{x}{\tau}) < 1$, d.h.\ $\frac{x}{\tau} \in C$, also $x \in \tau \, C$ und somit $p(x) \le \tau$.
Daher folgt $p(x) \le p_{C}(x)$.
Für jedes $\tau \in \R_+$ mit $\tau < p(x)$ gilt zum anderen $p(\frac{x}{\tau}) > 1$, d.h.\ $\frac{x}{\tau} \notin C$, also $x \notin \tau \, C$ und somit $p(x) \ge \tau$.
Daher folgt auch $p(x) \ge p_{C}(x)$. \q
\A
\textit{Beweis des Satzes.} Jedes Minkowski-Funktional einer absorbierenden konvexen ausgewogenen Teilmenge von $X$ ist nach \ref{FA.N.4} (iv) eine Halbnorm auf $X$.
Sei $p$ eine Halbnorm auf $X$.
Dann ist $C := \overline{p}^1({[}0,1{[})$ nach \ref{FA.N.3} eine absorbierende konvexe ausgewogene Teilmenge von $X$ und aus Lemma \ref{FA.N.6} folgt $p = p_C$. \q

\begin{Satz} \label{FA.N.E}
Zu jeder absorbierenden konvexen Teilmenge $C$ von $X$, die (\ref{FA.2.23.2}), d.h.\ $\forall_{x \in C} \exists_{\tau \in {]}1, \infty{[}} \, \tau \, x \in C$, sowie (\ref{FA.2.23.3}), d.h.\ $\forall_{\sigma \in \K, \, |\sigma|=1} \, \sigma \, C \subset C$, und (\ref{FA.2.23.4}), d.h.\ $\bigcap_{\tau \in \R_+} \tau \, C = \emptyset$, genügt, existiert genau eine Norm -- nämlich das Minkowski-Funktional $p_C$ -- derart, daß  $U_1(0) = C$ gilt.
Im Falle ihrer Existenz läßt sich darüber hinaus jede Norm so erzeugen.
\end{Satz}

\textit{Beweis.} 1.) Sei $C$ eine Menge wie im im Satz.
Daß das Min\-kow\-ski-Funktional $p_C$ von $C$ dann eine Norm ist, haben wir in \ref{FA.N.4} eingesehen.
Wegen (\ref{FA.2.23.2}) gilt dann nach \ref{FA.N.4} (iii) auch $U_1^{p_C}(0) = C$.

Sei $\|\ldots\| \: X \to \R$ eine weitere Norm auf $X$ mit $U^{\|\ldots\|}_1(0) = C$.
Wir behaupten
\begin{equation} \label{FA.N.E.1}
B^{p_C}_1(0) = B^{\|\ldots\|}_1(0).
\end{equation}

{[} Zu (\ref{FA.N.E.1}): ,,$\subset$`` Sei $x \in X$ mit $p_C(x) \le 1$.
Dann gilt $p_C(\frac{n}{n+1} \, x) < 1$ für alle $n \in \N$, also $\frac{n}{n+1} \, x \in C$, d.h.\ wegen $U^{\|\ldots\|}_1(0) = C$
$$ \underbrace{\frac{n}{n+1}}_{\stackrel{n \to \infty}{\longrightarrow} 1} \|x\| = \left\| \frac{n}{n+1} \, x \right\| < 1. $$
Somit folgt $\|x\| \le 1$.

,,$\supset$`` sieht man analog zu ,,$\subset$`` ein, indem man die Rollen von $p_C$ und $\|\ldots\|$ vertauscht. {]}

Aus (\ref{FA.N.E.1}) folgt wegen $U_1^{p_C}(0) = C = U^{\|\ldots\|}_1(0)$
\begin{equation*} \label{FA.N.E.2}
\{ x \in X \, | \, p_C(x) = 1 \} = \{ x \in X \, | \, \|x\| = 1 \} =: S,
\end{equation*}
also gilt 
\begin{equation} \label{FA.N.E.3}
p_C|_S = \|\ldots\||_S.
\end{equation}
Hieraus ergibt sich 
\begin{equation} \label{FA.N.E.4}
p_C|_C = \|\ldots\||_C.
\end{equation}

{[} Zu (\ref{FA.N.E.4}): Sei $x \in C = U^{\|\ldots\|}_1(0)$ beliebig.
Ohne Einschränkung gelte $x \ne 0$.
Aus (\ref{FA.N.E.3}) folgt $p_C(\frac{x}{\|x\|}) = \|\frac{x}{\|x\|}\|$.
Weil $p_C$ und $\|\ldots\|$ Normen sind, gilt dann auch $p_C(x) = \|x\|$. {]}

Sei schließlich $\tilde{x} \in X$.
Da $C$ absorbierend ist, existieren $\tau \in \R_+$ und $x \in C$ mit $\tilde{x} = \tau \, x$.
Dann folgt $p_C(\tilde{x}) = \tau \, p_C(x) \stackrel{(\ref{FA.N.E.4})}{=} \tau \, \|x\| = \|\tilde{x}\|$, also gilt $\|\ldots\| = p_C$.

2.) Sei $\|\ldots\| \: X \to \R$ eine Norm auf $X$.
Wir setzen
$$ C := U_1^{\|\ldots\|}(0). $$
Wegen \ref{FA.N.6} bleibt zu zeigen, daß $C$ die Eigenschaften (\ref{FA.2.23.2}) - (\ref{FA.2.23.4}) erfüllt.

(\ref{FA.2.23.2}) ist klar, da $C$ bzgl.\ $\|\ldots\|$ offen ist -- beachte, daß für jedes $x \in C$ gilt $C \in U(x,X)$ und $\lim_{n \to \infty} \frac{n+1}{n} \, x = x$ --;
(\ref{FA.2.23.3}) ist erfüllt, da $C$ nach \ref{FA.N.3} (ii) sogar ausgewogen ist; 
und (\ref{FA.2.23.4}) gilt, da $C$ bzgl.\ $\|\ldots\|$ beschränkt ist. \q

\subsection*{\textsc{Hahn-Banach} Sätze} \addcontentsline{toc}{subsection}{\textsc{Hahn-Banach} Sätze}

Wir beweisen in \ref{FA.2.17} eine rein algebraische Fassung des Fortsetzungssatzes von \textsc{Hahn-Banach}.
Unter stärkeren Prämissen erhält man auch stärkere Konklusionen, die wir im Anschluß aus \ref{FA.2.17} folgern.
Die für die Theorie der normierten Vektorräume interessante Version des Fortsetzungssatzes wird in \ref{FA.2.19} hergeleitet.

\begin{Def}[Algebraischer und topologischer Dualraum] \label{FA.2.15}
Sei $X$ ein $\K$-Vek\-tor\-raum.
\begin{itemize}
\item[(i)] Wie üblich bezeichnen wir mit 
$$ \boxed{X^*} := {\rm Hom}_{\K}(X,\K) $$
den \emph{algebraischen Dualraum von $X$}\index{Dualraum!algebraischer}.
Seine Elemente sind die Funktionale von $X$, die wir auch \emph{Linearformen von $X$} nennen.

\begin{Bsp*}
Sei $n \in \N_+$.
Im Falle $X = \K^n$ ist
$$ \K^n \longrightarrow (\K^n)^*, ~~ (a_1, \ldots, a_n) \longmapsto \sum_{i=1}^n a_i \, x_i, $$
wobei $x_i \: \K^n \to \K$ für $i \in \{1, \ldots, n\}$ hier die Projektion auf die $i$-te Komponente bezeichne, ein $\K$-Vektorraum-Isomorphismus.
\end{Bsp*} 
\item[(ii)] Ist $X$ zusätzlich normiert, so ist der \emph{topologische Dualraum von $X$}\index{Dualraum!topologischer} per definitionem
$$ \boxed{X'} := \mathcal{L}_{\K}(X,\K), $$
welcher nach \ref{FA.2.11} ein $\K$-Banachraum ist.

\begin{Bsp*}
$T \in \mathcal{C}([-1,1],\R)^*$ sei durch $\forall_{f \in \mathcal{C}([-1,1],\R)} \, T(f) = f(0)$ definiert.
\begin{itemize}
\item[1.)] $T \: (\mathcal{C}([-1,1],\R), \|\ldots\|_{\infty}) \to \R$ ist stetig, denn für $f \in \mathcal{C}([-1,1],\R)$ gilt $T(f) = f(0) \le \|f\|_{\infty}$.
\item[2.)] $T \: (\mathcal{C}([-1,1],\R), \|\ldots\|_1) \to \R$ ist nicht stetig, denn für jedes $k \in \N_+$ ist durch
$$ \forall_{t \in [-1,1]} \, f_k(t) := 
\left\{ \begin{array}{cc}
0, & t \in [-1,-\frac{1}{k}] \cup [\frac{1}{k},1], \\ 
k^2 \, t + k, & t \in {]}-\frac{1}{k},0{[}, \\
-k^2 \, t + k, & t \in {[}0,\frac{1}{k}{[}
\end{array} \right. $$
eine stetige Funktion $f_k \: [-1,1] \to \R$ mit $|T(f_k)| = |f_k(0)| = k$ und $\|f_k\|_1 = \int_{-1}^1 |f_k(t)| \, \d t = 1$ gegeben.
\end{itemize}
\end{Bsp*}
\end{itemize}
\end{Def}

Beim Beweis des oben angekündigten folgenden Satzes werden wir das Zornsche Lemma und die damit verbundenen Begriffe verwenden, vgl.\ Anhang \ref{FAnaA}.

\begin{HS}[Algebraischer Fortsetzungssatz von \textsc{Hahn-Banach} in der Version für $\R$-Vek\-tor\-räume] \index{Satz!von \textsc{Hahn-Banach}!Fortsetzungs-} \index{Satz!Fortsetzungs-!von \textsc{Hahn-Banach}} \label{FA.2.17} $\,$

\noindent \textbf{Vor.:} Es seien $X$ ein $\R$-Vektorraum, $p \: X \to \R$ eine sublineare Funktion und $Y$ ein $\R$-Untervektorraum von $X$.

\noindent \textbf{Beh.:} Zu jedem $f \in Y^*$ mit $f \le p|_Y$ existiert ein $F \in X^*$ mit $F|_Y = f$ und $F \le p$.
\end{HS}

\textit{Beweis.} Sei $f \in Y^*$ mit $f \le p|_Y$. 
Auf der Menge
$$ \mathcal{X} := \{ g \in {Y_g}^* \, | \, \mbox{$Y_g$ $\R$-Untervektorraum von $X \, \wedge \, Y \subset Y_g \, \wedge \, g|_Y = f \, \wedge \, g \le p|_{Y_g}$} \} $$
wird durch
$$ \forall_{(g_1 \: Y_{g_1} \to \R), (g_2 \: Y_{g_2} \to \R) \in \mathcal{X}} \, g_1 \preceq g_2 : \Longleftrightarrow Y_{g_1} \subset Y_{g_2} \, \wedge \, g_2|_{Y_{g_1}} = g_1 $$
offenbar eine Ordnung definiert.

Wegen $f \in \mathcal{X}$ gilt $\mathcal{X} \ne \emptyset$, und wir behaupten:
\begin{equation} \label{FA.2.17.1}
\begin{array}{l}
\mbox{Jede bzgl.\ $\preceq$ totalgeordnete nicht-leere Teilmenge von $\mathcal{X}$ besitzt eine} \\
\mbox{obere Schranke in $\mathcal{X}$.}
\end{array}
\end{equation}

{[} Zu (\ref{FA.2.17.1}): Sei $M$ eine bzgl.\ $\preceq$ totalgeordnete nicht-leere Teilmenge von $\mathcal{X}$.
Dann ist
\begin{equation} \label{FA.2.17.2}
Z := \bigcup_{(g \: Y_g \to \R) \in M} Y_g 
\end{equation}
ein $\R$-Untervektorraum von $X$ mit $Y \subset Z$.

Beweis hiervon:
Seien $z_1,z_2 \in Y$ und $\lambda \in \R$.
Dann existieren für $i \in \{1,2\}$ Elemente $g_i \: Y_{g_i} \to \R$ von $M$ mit $z_i \in Y_i := Y_{g_i}$.
Da $M$ bzgl.\ $\preceq$ totalgeordnet ist, können wir ohne Beschränkung der Allgemeinheit $Y_1 \subset Y_2$ annehmen, also gilt $z_1, z_2 \in Y_2$.
$Y_2$ ist ein $\R$-Untervektorraum von $X$, daher folgt $z_1 + z_2, \lambda \, z_1 \in Y_2 \subset Z$, d.h.\ $Z$ ist ein $\R$-Untervektorraum von $X$.
$Z \subset Y$ ist trivial.

Wir behaupten weiter, daß durch
\begin{equation} \label{FA.2.17.3}
\tilde{g} \: Z \longrightarrow \R, ~~ z \longmapsto g(z), \mbox { wobei } (g \: Y_g \longrightarrow \R) \in M \mbox{ mit } z \in Y_g \mbox{ beliebig,}
\end{equation}
eine Linearform auf $X$ wohldefiniert ist.

Beweis der Wohldefiniertheit von $\tilde{g}$:
Seien $z \in Z$ und $g_1 \: Y_{g_1} \to \R$ sowie $g_2 \: Y_{g_2} \to \R$ Elemente von $M$ mit $z \in Y_1 \cap Y_2$, wobei wir wieder $Y_i := Y_{g_i}$ für $i \in \{1,2\}$ setzen.
Da $M$ bzgl.\ $\preceq$ totalgeordnet ist, können wir ohne Einschränkung $g_1 \preceq g_2$ annehmen, also gilt $z \in Y_1 \subset Y_2$ und $g_2(z) = g_1(z)$.
Damit ist die Wohldefiniertheit von $\tilde{g}$ gezeigt.

Beweis der $\R$-Linearität von $\tilde{g}$: Seien $z_1,z_2 \in Y$ und $\lambda \in \R$.
Es existieren zu $i \in \{1,2\}$ wieder Elemente $g_i \: Y_{g_i} \to \R$ von $M$ mit $z_i \in Y_i := Y_{g_i}$, und wir können erneut ohne Beschränkung der Allgemeinheit annehmen $g_1 \preceq g_2$.
Dann gilt
\begin{gather*}
\tilde{g}( \lambda \, z_1 ) = g_1( \lambda \, z_1 ) = \lambda \, g_1(z_1) = \lambda \, \tilde{g}(z_1), \\
\tilde{g}( z_1 + z_2 ) = g_2( z_1 + z_2 ) = g_2(z_1) + g_2(z_2) = \tilde{g}(z_1) + \tilde{g}(z_2).
\end{gather*}

Damit haben wir gezeigt: $\tilde{g} \in Z^*$.
Daß $\tilde{g}$ zusätzlich ein Element von $\mathcal{X}$ und eine obere Schranke von $M$ ist, ergibt sich aus (\ref{FA.2.17.2}) und (\ref{FA.2.17.3}).
Damit haben (\ref{FA.2.17.1}) bewiesen. {]}

Als einelementige Teilmenge von $\mathcal{X}$ ist $\{ f \: Y \to \R \}$ eine bzgl.\ $\preceq$ totalgeordnete Teilmenge von $\mathcal{X}$.
Daher folgt aus (\ref{FA.2.17.1}) und der verschärften Version des Lemmas von \textsc{Zorn} \ref{FA.A.0} die Existenz eines $F \in \mathcal{X}$ derart, daß
\begin{equation} \label{FA.2.17.M}
\mbox{$F \: Y_F \longrightarrow \R$ maximales Element von $\mathcal{X}$ bzgl.\ $\preceq$}
\end{equation}
und obere Schranke von $\{ f \: Y \to \R \}$ ist, d.h.\ insbes.\ $Y_F$ ist ein $\R$-Unter\-vek\-tor\-raum von $X$ mit $Y \subset Y_F$ sowie
\begin{equation} \label{FA.2.17.X}
F \in (Y_F)^* ~~ \wedge ~~ F|_Y = f ~~ \wedge ~~ F \le p|_{Y_F}.
\end{equation}
Zum Nachweis des Hauptsatzes bleibt daher zu zeigen:
$$ Y_F = X. $$

Beweis hiervon:
Angenommen, es existiert $a \in X \setminus Y_F$.
Dann ist 
\begin{equation} \label{FA.2.17.4}
\mbox{$\widetilde{Y} := Y_F \oplus \R \, a$ ein $\R$-Untervektorraum von $X$ mit $Y \subsetneqq \widetilde{Y}$.}
\end{equation}

{[} Zu (\ref{FA.2.17.4}): $Y \subsetneqq \widetilde{Y}$ ist wegen $a \in X \setminus Y$ klar, und da für alle $y_1, y_2 \in Y_F$ sowie $t_1, t_2 \in \R$ gilt
$$ y_1 + t_1 \, a = y_2 + t_2 \, a \Longrightarrow \underbrace{y_1 - y_2}_{\in Y_F} = \underbrace{(t_2 - t_1) \, a}_{\in \R a} \stackrel{a \in X \setminus Y_F}{\Longrightarrow} t_2 - t_1 = 0 \, \wedge y_1 - y_2 = 0, $$
ist $Y_F + \R \, a$ sogar eine direkte Summe von $\R$-Vektorräumen. {]}

Wir zeigen nun
\begin{equation} \label{FA.2.17.5}
\forall_{y_0, y \in Y_F} \, F(y) - p(y-a) \le p(y_0 + a) - F(y_0).
\end{equation}

{[} Zu (\ref{FA.2.17.5}): Seien $y_0, y \in Y_F$ beliebig.
Dann folgt aus der $\R$-Linearität von $F$, $F \le p|_{Y_F}$ und der Sublinearität von $p$
\begin{eqnarray*}
F(y) + F(y_0) & = & F(y + y_0) \le p(y + y_0) = p( (y - a) + (y_0 + a) ) \\
& \le & p(y - a) + p(y_0 + a), 
\end{eqnarray*}
also $F(y) - p(y-a) \le p(y_0 + a) - F(y_0)$. {]}

Wegen (\ref{FA.2.17.5}) gilt
\begin{equation} \label{FA.2.17.s}
s := \sup \{ F(y) - p(y-a) \, | \, y \in Y_F \} < \infty,
\end{equation}
und wir können $\widetilde{F} \: \widetilde{Y} \to \R$ durch
\begin{equation} \label{FA.2.17.F}
\forall_{y + t \, a \in \widetilde{Y} = Y_F \oplus \R a} \, \widetilde{F}(y + t \, a) := F(y) + t \, s
\end{equation}
definieren.
Dann folgt
\begin{gather} 
\widetilde{F} \in \widetilde{Y}^*, \label{FA.2.17.6} \\
\widetilde{F}|_{Y_F} = F. \label{FA.2.17.7}
\end{gather}

{[} (\ref{FA.2.17.7}) ist trivial, und (\ref{FA.2.17.6}) ergibt sich daraus, daß für alle $y_1, y_2 \in Y_F$, $t_1, t_2 \in \R$ und $\lambda \in \R$
\begin{eqnarray*}
\widetilde{F}( \lambda \, (y_1 + t_1 \, a) ) & = & \widetilde{F}( \lambda \, y_1 + (\lambda \, t_1) \, a ) = F( \lambda \, y_1 ) + (\lambda \, t_1) \, s \\
& = & \lambda \, F(y_1) + \lambda \, (t_1 \, s) = \lambda \, (F(y_1) + t_1 \, s) = \lambda \, \widetilde{F}(y_1 + t_1 \, a)
\end{eqnarray*}
sowie
\begin{eqnarray*}
\widetilde{F}( (y_1 + t_1 \, a) + (y_2 + t_2 \, a) ) & = & \widetilde{F}( (y_1 + y_2) + (t_1 + t_2) \, a ) \\
& = & F( y_1 + y_2 ) + (t_1 + t_2) \, s  \\
& = & F(y_1) + F(y_2) + t_1 \, s + t_2 \, s \\
& = & ( F(y_1) + t_1 \, s ) + ( F(y_2) + t_2 \, s ) \\
& = & \widetilde{F}(y_1 + t_1 \, a) + \widetilde{F}(y_2 + t_2 \, a),
\end{eqnarray*}
gilt. {]}

Schließlich weisen wir nach
\begin{equation} \label{FA.2.17.8}
\widetilde{F} \le p|_{\widetilde{Y}}.
\end{equation}

{[} Zu (\ref{FA.2.17.8}): Aus (\ref{FA.2.17.s}) folgt
\begin{equation} \label{FA.2.17.9}
\forall_{y_0 \in Y_F} \, F(y_0) - s \le p(y_0 - a),
\end{equation}
und aus (\ref{FA.2.17.s}), (\ref{FA.2.17.5}) folgt
\begin{equation} \label{FA.2.17.10}
\forall_{y_0 \in Y_F} \, F(y_0) + s \le p(y_0 + a).
\end{equation}

Sei $y + t \, a \in \widetilde{Y} = Y_F \oplus \R \, a$.

1. Fall: $t=0$.
Dann gilt
$$ \widetilde{F}(y + t \, a) = \widetilde{F}(y) \stackrel{(\ref{FA.2.17.7})}{=} F(y) \stackrel{\text{Vor.}}{\le} p(y) = p(y + t \, a). $$

2. Fall: $t \in \R_+$.
Dann gilt
\begin{eqnarray*}
\widetilde{F}(y + t \, a) & \stackrel{(\ref{FA.2.17.6})}{=} & t \, \widetilde{F}(\frac{y}{t} + a) \stackrel{(\ref{FA.2.17.F})}{=} t \, (F(\frac{y}{t}) + s) \\
& \stackrel{(\ref{FA.2.17.10})}{\le} & t \, p(\frac{y}{t} + a) \stackrel{\text{Vor.}}{=} p(y + t \, a).
\end{eqnarray*}

3. Fall: $t \in \R_-$.
Dann gilt
\begin{eqnarray*}
\widetilde{F}(y + t \, a) & \stackrel{(\ref{FA.2.17.6})}{=} & -t \, \widetilde{F}(-\frac{y}{t} - a) \stackrel{(\ref{FA.2.17.F})}{=} -t \, (F(-\frac{y}{t}) - s) \\
& \stackrel{(\ref{FA.2.17.9})}{\le} & -t \, p(-\frac{y}{t} - a) \stackrel{\text{Vor.}}{=} p(y + t \, a).
\end{eqnarray*}

Damit ist (\ref{FA.2.17.8}) gezeigt. {]}

Aus (\ref{FA.2.17.6}), (\ref{FA.2.17.4}), $\widetilde{F}|_Y = f$ -- wegen (\ref{FA.2.17.7}), (\ref{FA.2.17.X}) -- und (\ref{FA.2.17.8}) folgt $\widetilde{F} \in \mathcal{X}$.
Nach (\ref{FA.2.17.4}) und (\ref{FA.2.17.7}) gilt ferner $\widetilde{F} \succ F$, im Widerspruch zu (\ref{FA.2.17.M}). \q

\begin{Satz}[Algebraischer Fortsetzungssatz von \textsc{Hahn-Banach} in der Version für $\K$-Vek\-tor\-räume] \index{Satz!von \textsc{Hahn-Banach}!Fortsetzungs-} \index{Satz!Fortsetzungs-!von \textsc{Hahn-Banach}} \label{FA.2.18} $\,$

\noindent \textbf{Vor.:} Es seien $X$ ein $\K$-Vektorraum, $p \: X \to \R$ eine Halbnorm und $Y$ ein $\K$-Un\-ter\-vek\-tor\-raum von $X$.

\noindent \textbf{Beh.:} Zu jedem $f \in Y^*$ mit $|f| \le p|_Y$ existiert ein $F \in X^*$ mit $F|_Y = f$ und $|F| \le p$.
\end{Satz}

\textit{Beweis.} Sei $f \in Y^*$ mit $|f| \le p|_Y$.

1. Fall: $\K = \R$.
Dann folgt aus $f \le |f| \le p|_Y$ und \ref{FA.2.17} die Existenz von $F \in X^*$ mit $F|_Y = f$ und $F \le p$.
Da $p$ eine Halbnorm ist, gilt zusätzlich 
$$ \forall_{x \in X} \, -F(x) = F(-x) \le p(-x) = p(x), $$
also $|F| \le p$.

2. Fall: $\K = \C$.
$X$ sowie $Y$ sind in kanonischer Weise $\R$-Vektorräume, die wir mit $X_{\R}$ sowie $Y_{\R}$ bezeichnen.
$p$ ist dann auch eine Halbnorm für $X_{\R}$, und $Y_{\R}$ ist ein $\R$-Untervektorraum von $X_{\R}$.
Es gilt offenbar $f_1 := {\rm Re} f \in {Y_{\R}}^* = {\rm Hom}_{\R}(Y_{\R},\R)$ und $|f_1| \le |f| \le p|_{Y_{\R}}$.
Nach dem 1. Fall existiert $F_1 \in {X_{\R}}^* = {\rm Hom}_{\R}(X_{\R},\R)$ mit $F_1|_{Y_{\R}} = f_1$ und $|F_1| \le p$.
Der Leser zeige als Übung, daß durch
$$ \forall_{x \in X} \, F(x) := F_1(x) - \i \, F_1( \i \, x ) $$
ein Element $F$ von $X^* = {\rm Hom}_{\C}(X,\C)$ mit $F|_Y = f$ gegeben ist.\footnote{Verwende, daß die Abbildung $$X^* \longrightarrow (X_{\R})^*, ~~ F \longmapsto {\rm Re} \, F,$$ eine $\R$-lineare Bijektion mit $\R$-linearer Umkehrabbildung $$(X_{\R})^* \longrightarrow X^*, ~~ F_1 \longmapsto \left( x \mapsto F_1(x) - \i \, F_1( \i \, x ) \right),$$ bei dessen Nachweis $\forall_{z \in \C} \, {\rm Re}(\i \, z) = - {\rm Im} (z)$ eingeht, ist.}
Des weiteren existiert zu jedem $x \in X$ ein $t \in {[} 0, 2 \pi {[}$ mit $\e^{\i t} \, F(x) = |F(x)| \in {[}0, \infty{[}$, d.h.\ es gilt
$$ |F(x)| = F(|e^{\i \, t}| \, x) = |({\rm Re} \, F) (|e^{\i \, t}| \, x)| = |e^{\i \, t}| \, |F_1(x)| = |F_1(x)| \le p(x), $$
womit die Behauptung auch im 2. Falle bewiesen ist. \q

\begin{Satz}[Topologischer Fortsetzungssatz von \textsc{Hahn-Banach} für normierte $\K$-Vek\-tor\-räume] \index{Satz!Fortsetzungs-!von \textsc{Hahn-Banach}} \label{FA.2.19} $\,$

\noindent \textbf{Vor.:} Seien $X$ ein normierter $\K$-Vektorraum und $Y$ ein $\K$-Un\-ter\-vek\-tor\-raum von $X$, der dann in kanonischer Weise ebenfalls ein normierter $\K$-Vek\-tor\-raum ist.

\noindent \textbf{Beh.:} Zu jedem $f \in Y'$ existiert ein $F \in X'$ mit $F|_Y = f$ und $\|F\| = \|f\|$.
\end{Satz}

\textit{Beweis.} Sei $f \in Y'$. 
Dann wird durch
$$ \forall_{x \in X} \, p(x) := \|f\| \, \|x\| $$
offenbar eine Halbnorm auf $X$ definiert, und es gilt -- beachte, daß $f$ ein beschränkter Operator ist --
$$ \forall_{y \in Y} \, |f(y)| \stackrel{(\ref{FA.2.8.Doppelstern})}{\le} \|f\| \, \|y\| = p(y), $$
also folgt aus \ref{FA.2.18} die Existenz von $F \in X^*$ mit $F|_Y = f$ sowie $|F| \le p = \|f\| \, \| \ldots \|$.
Letzteres bedeutet, daß auch $F$ ein beschränkter Operator ist, d.h.\ $F \in X'$, und
wegen (\ref{FA.2.8.Stern}) weiterhin $\| F \| \le \|f\|$. 
Da trivialerweise 
$$ \sup \{ |\underbrace{f(y)}_{= F(y)}| \, | \, y \in Y \, \wedge \|y \| \le 1 \} \le \sup \{ |F(x)| \, | \, x \in X \, \wedge \|x\| \le 1 \} $$
gilt, folgt aus (\ref{FA.2.8.Stern}) außerdem $\|f\| \le \|F\|$. \q

\begin{Kor} \label{FA.2.20}
Sei $X$ ein normierter $\K$-Vektorraum. 
Dann gilt:
\begin{itemize}
\item[(i)] Ist $Y$ ein $\K$-Un\-ter\-vek\-tor\-raum von $X$, so existiert zu jedem $x_0 \in X$ mit positivem Abstand zu $Y$, d.h.\
\begin{equation}
\delta := d(\{x_0\},Y) = \inf \{ \| x_0 - y \| \, | \, y \in Y \} > 0,%
\footnote{Es sei darauf hingewiesen, daß aus der Kompaktheit von $\{x_0\}$ und Übung \ref{FA.disj abgeschl und komp} $$\delta = d(\{x_0\}, Y) > 0$$ folgt, falls $Y$ ein abgeschlossener $\K$-Untervektorraum von $X$ ist.}
\label{FA.2.20.0}
\end{equation}
ein $F \in X'$ mit $F(x_0) = \delta$, $F|_Y = 0$ und $\|F\| = 1$.
\item[(ii)] Zu jedem $x_0 \in X \setminus \{0\}$ existiert ein $F \in X'$ mit $F(x_0) = \|x_0\|$ und $\|F\| = 1$.
\item[(iii)] Sind $x_1, x_2 \in X$ mit $x_1 \ne x_2$, so existiert $F \in X'$ mit $F(x_1) \ne F(x_2)$.
Man nennt $X'$ daher auch \emph{punktetrennend}.\index{Punktetrennungseigenschaft!des topolpologischen Dualraumes}\index{Trennungs!-eigenschaft von Punkten}
\end{itemize}
\end{Kor}

\textit{Beweis.} Zu (i): Wir machen eine ähnliche Konstruktion wie im Beweis von \ref{FA.2.17} und setzen $\widetilde{Y} := Y \oplus \K \, x_0$.
Dann wird durch
$$ \forall_{y \in Y} \forall_{\lambda \in \K} \, f(y + \lambda \, x_0) := \lambda \, \delta $$
offenbar eine Linearform $f \: \widetilde{Y} \to \K$ mit $f(x_0) = \delta$ und $f|_Y = 0$ definiert.
Wir werden zeigen, daß gilt
\begin{equation}
f \in \widetilde{Y}' ~~ \wedge ~~ \|f\| = 1. \label{FA.2.20.1}
\end{equation}
Daher folgt aus \ref{FA.2.19} die Existenz von $F \in X'$ mit $F|_{\widetilde{Y}} = f$ und $\|F\| = \|f\|$, d.h.\ $F$ hat die gewünschten Eigenschaften.

{[} Zu (\ref{FA.2.20.1}): Für alle $y \in Y$ und $\lambda \in \K$ gilt
$$ |f(y + \lambda \, x_0)| \stackrel{\text{Def.}}{=} |\lambda| \, \delta \le \| y + \lambda \, x_0 \|, $$
denn für $\lambda = 0$ ist dies trivial, und für $\lambda \ne 0$ ergibt sich
$$ \| y + \lambda \, x_0 \| = |\lambda| \, \left\| x_0 - \frac{-y}{\lambda} \right\| \stackrel{(\ref{FA.2.20.0})}{\ge} |\lambda| \, \delta. $$
Also ist $f$ stetig mit $\|f\| \le 1$; zu zeigen bleibt $\|f\| \ge 1$.

Hierzu sei $\varepsilon \in \R_+$.
Aus (\ref{FA.2.20.0}) folgt die Existenz eines $y \in Y$ mit
$$ \delta + \varepsilon \, \delta \ge \| x_0 - y \|, $$
also gilt
$$ 0 < \delta \stackrel{(\ref{FA.2.20.0})}{\le} \| x_0 - y \| \le \delta + \varepsilon \, \delta \stackrel{\text{Def.}}{=} |f(y - x_0)| + \varepsilon \, \delta \le  \|f\| \, \| y - x_0 \| + \varepsilon \, \delta, $$
d.h.\
$$ 1 \le \|f\| + \frac{\varepsilon \, \delta}{\|x_0 - y\|} \le \|f\| + \varepsilon, $$
und aus der Beliebigkeit von $\varepsilon \in \R_+$ ergibt sich $1 \le \|f\|$. {]}

Zu (ii): Die Behauptung ergibt sich sofort, indem man (i) auf $Y := \{0\}$ anwendet, denn dann gilt $\delta = \|x_0\|$.

Zu (iii): Sei $x_0 := x_2 - x_1 \in X \setminus \{0\}$.
Dann existiert nach (ii) $F \in X'$ mit $0 \ne F(x_0) = F(x_2 - x_1) = F(x_2) - F(x_1)$, also $F(x_1) \ne F(x_2)$. \q

\begin{Bem*} $\,$
\begin{itemize}
\item[1.)] Das Korollar zeigt, daß für einem normierten $\K$-Vektorraum die Menge der stetigen Linearformen ,,reichhaltig`` ist.
Wir haben oben bereits gesehen, daß im Falle endlicher Dimension sogar jede Linearform stetig ist.
\item[2.)] Ist $X$ ein normierter $\K$-Vektorraum mit $\dim_{\K} X = \infty$, so existiert stets eine unstetige Linearform:
Nach \ref{FA.A.2} existiert nämich eine Basis $\{ b_i \, | \, i \in I \}$ von $X$, wobei $I$ eine Menge mit $\# I \ge \# \N$ sein muß.
Ohne Beschränkung der Allgemeinheit gelte $\forall_{i \in I} \, \|b_i\| = 1$.
Wir können eine abzählbare Teilmenge $\{ b_{i_k} \, | \, k \in \N \}$ von $I$ wählen und eine Linearform $F \: X \to \K$ durch
$$ \forall_{x = \sum_{i \in I} \lambda_i \, b_i \text{ mit eindeutig bestimmten } \lambda_i \in \K, \text{ wobei fast alle } \lambda_i = 0 \text{ sind}} \, F(x) = \sum_{k=0}^{\infty} k \, {\lambda}_{i_k} $$
definieren.
Dann gilt $\forall_{k \in \N} \, |F(b_{i_k})| = k$, d.h.\ $F$ kann kein beschränkter Operator sein. 
\end{itemize}
\end{Bem*}

Unser nächstes Ziel ist der Beweis der Trennungssätze von \textsc{Hahn-Banach}, welche besagen, daß sich disjunkte konvexe Mengen in einem nun zu erklärenden Sinne trennen lassen.
 
Ist $X$ ein (ggf.\ normierter) $\C$-Vektorraum so bezeichnen wir mit $\boxed{X_{\R}}$ den kanonisch aus $X$ erhaltenen (normierten) $\R$-Vektorraum, d.i.\ seine sog.\ \emph{Reellifizierung}\index{Reellifizierung}.
Als Mengen (bzw.\ ggf.\ Abbildungen) stimmen $X$ und $X_{\R}$ (bzw.\ ihre Normen) natürlich überein.
Um die Definitionen und Sätze einheitlich formulieren zu können, vereinbaren wir im Falle eines (ggf.\ normierten) $\R$-Vektor\-raumes $X$ zusätzlich:
$$ X_{\R} := X ~~ \mbox{ und } ~~  \forall_{F \in X^*} \, {\rm Re} \, F := F.  $$

\begin{Def}[Hyperebene] \index{Hyperebene} \label{FA.2.21.H}
Es sei $X$ ein $\K$-Vektorraum.
\begin{itemize}
\item[(i)] Seien $F \in X^* \setminus \{0\}$ und $\alpha \in \K$.

Wir setzen
$$ \boxed{H_{F,\alpha}} := \overline{F}^1(\{\alpha\}) = \{ x \in X \, | \, F(x) = \alpha \}. $$
Diese Menge ist nicht-leer, da $F$ offenbar surjektiv ist.

\begin{Bem*} $\,$
\begin{itemize}
\item[1.)] $H_{F,\alpha}$ ist genau dann ein $\K$-Untervektorraum von $X$, wenn $\alpha = 0$ gilt.
In diesem Fall ergibt der Homomorphiesatz der linearen Algebra $\dim_{\K} (V/H_{F,0}) = 1$.
\item[2.)] Ist $Y$ ein $\K$-Untervektorraum von $X$ mit $\dim_{\K} X/Y = 1$, so existiert nach \ref{FA.A.2} eine Basis von $Y$, die durch Hinzunahme eines Elementes $x_0 \in X \setminus Y$ zu einer Basis von $X$ wird, und die Projektion $X \to \K \, x_0$ liefert mittels der Identifikation $\K \, x_0 \equiv \K$ ein Element von $X^* \setminus \{0\}$, welches $Y$ als Kern besitzt.
\item[3.)] Für jedes $a \in H_{F,\alpha}$ gilt $H_{F,\alpha} = a + H_{F,0}$.
\end{itemize}
\end{Bem*}

\item[(ii)] Ist $H$ eine Teilmenge von $X$, so heißt $H$ eine \emph{{$\R$}-Hyperebene} genau dann, wenn $F \in (X_{\R})^* \setminus \{0\}$ und $\alpha \in \R$ mit $H = H_{F,\alpha}$ existieren.

\begin{Bem*}
Für jedes $F \in X^* \setminus \{0\}$ gilt ${\rm Re} \, F \in (X_{\R})^* \setminus \{0\}$.
\end{Bem*}
\end{itemize}
\end{Def}

\begin{Bem*}
Sind $X$ ein normierter $\K$-Vektorraum und $F \in X'$, so gilt auch ${\rm Re} \, F \in (X_{\R})'$, denn aus der Existenz eines $C \in \R$ mit $\forall_{x \in X} \, |F(x)| \le C \, \|x\|$ folgt, daß für jedes $x \in X_{\R}$ gilt: $|({\rm Re} \, F)(x)| \le |F(x)| \le C \, \|x\|$.
\end{Bem*} 

\begin{Lemma} \label{FA.2.21}
Es seien $X$ ein normierter $\K$-Vektorraum und $F \in X^* \setminus \{0\}$.
Dann gilt:
\begin{itemize}
\item[(i)] $F \in X' \Longleftrightarrow \forall_{\alpha \in \K} \, H_{F,\alpha} \mbox{ ist abgeschlossen in } X$.
\item[(ii)] $F \notin X' \Longleftrightarrow \forall_{\alpha \in \K} \, H_{F,\alpha} \mbox{ ist dicht in } X$.
\end{itemize}
\end{Lemma}

\textit{Beweis.} Sei $\alpha \in \K$ beliebig.
Es genügt offenbar zu zeigen, daß die folgenden Aussagen paarweise äquivalent sind:
\begin{gather}
F \in X'. \label{FA.2.21.1} \\
H_{F,\alpha} \mbox{ ist abgeschlossen in } X. \label{FA.2.21.2} \\
\overline{H_{F,\alpha}} \ne X. \label{FA.2.21.3}
\end{gather}

,,(\ref{FA.2.21.1}) $\Rightarrow$ (\ref{FA.2.21.2})`` ist wegen (\ref{FA.1.12.2}) klar.

,,(\ref{FA.2.21.2}) $\Rightarrow$ (\ref{FA.2.21.3})`` Aus der Surjektivität von $F$ folgt $H_{F,\alpha} \ne X$.
Hieraus und aus (\ref{FA.2.21.2}) ergibt sich (\ref{FA.2.21.3}).

,,(\ref{FA.2.21.3}) $\Rightarrow$ (\ref{FA.2.21.1})`` Wegen (\ref{FA.2.21.3}) können wir $x_0 \in X \setminus \overline{H_{F,\alpha}}$ fixieren.
Letzteres ist eine offene Menge, also existiert $\varepsilon \in \R_+$ mit $U_{\varepsilon}(x_0) \subset X \setminus \overline{H_{F,\alpha}}$, d.h.\
\begin{equation}
U_{\varepsilon}(x_0) \cap \overline{H_{F,\alpha}} = \emptyset. \label{FA.2.21.4}
\end{equation}
Wir setzen 
$$ C := \frac{F(x_0) - \alpha}{\varepsilon}. $$
Zum Nachweis von (\ref{FA.2.21.1}) genügt es zu zeigen, daß gilt
\begin{equation}
\forall_{x \in X} \, |F(x)| \le C \, \|x\|. \label{FA.2.21.5}
\end{equation}

Beweis hiervon:
Sei $x \in X$ beliebig.
Im Falle $F(x) = 0$ ist die Behauptung klar, gelte daher $F(x) \ne 0$.
Dann folgt offenbar $y := x_0 - \frac{F(x_0) - \alpha}{F(x)} \, x \in H_{F,\alpha}$, also nach (\ref{FA.2.21.4}): $\| y - x_0 \| \ge \varepsilon$.
Dies bedeutet genau (\ref{FA.2.21.5}). \q

\begin{Def} \label{FA.2.T} \index{Hyperebene!trennende -- zweier Teimengen}
Seien $X$ ein $\K$-Vektorraum sowie $F \in (X_{\R})^* \setminus \{0\}$ und $\alpha \in \R$, also ist $H := H_{F,\alpha}$ eine $\R$-Hyperebene von $X$.
Ferner seien $A,B$ nicht-leere Teilmengen von $X$.

Per definitionem \emph{trennt $H$ die Mengen $A$ und $B$}, wenn
\begin{equation}
\sup \{ F(a) \, | \, a \in A \} \le \alpha \le \inf \{ F(b) \, | \, b \in B \} \label{FA.2.T.1}
\end{equation}
bzw.\
\begin{equation}
\sup \{ F(b) \, | \, b \in B \} \le \alpha \le \inf \{ F(a) \, | \, a \in A \} \label{FA.2.T.2}
\end{equation}
gilt. 
Die Trennung heißt \emph{strikt}, wenn eine der Ungleichungen in (\ref{FA.2.T.1}) bzw.\ (\ref{FA.2.T.2}) echt ist.

\begin{Bem*}
Aus (\ref{FA.2.T.1}) folgt $A \subset \overline{F}^1( {]} - \infty, \alpha {]} )$ und $B \subset \overline{F}^1( {[} \alpha, \infty {[} )$, entsprechendes folgt aus (\ref{FA.2.T.2}).
Dies rechtfertigt den Begriff der ,,Trennung``, obwohl $A \cap B$ Elemente von $H = \overline{F}^1(\{\alpha\})$ enthalten kann.
Im Falle einer strikten Trennung gilt $A \cap B = \emptyset$.
\end{Bem*}
\end{Def}

\begin{HS}[Trennungssätze von \textsc{Hahn-Banach}] \index{Satz!von \textsc{Hahn-Banach}!Trennungs-}\index{Trennungs!-sätze von \textsc{Hahn-Banach}} \label{FA.2.24} $\,$

\noindent \textbf{Vor.:} Seien $X$ ein normierter $\K$-Vektorraum und $A,B$ nicht-leere disjunkte Teilmengen von $X$.
Ferner seien $A,B$ konvex.

\noindent \textbf{Beh.:}
\begin{itemize}
\item[(i)] Ist $A$ offen, so existieren $F \in X' \setminus \{0\}$ und $\alpha \in \R$ derart, daß gilt
$$ \forall_{a \in A} \forall_{b \in B} \, ({\rm Re} \, F)(a) < \alpha \le ({\rm Re} \, F)(b), $$
insbes.\ trennt die abgeschlossene $\R$-Hyperebene $H_{{\rm Re} \, F, \alpha}$ die konvexen Mengen $A$ und $B$.
\item[(ii)] Gilt $\delta := d(A,B) \in \R_+$, so existieren $F \in X' \setminus \{0\}$ und $\alpha, \beta \in \R$ mit 
$$ \forall_{a \in A} \forall_{b \in B} \, ({\rm Re} \, F)(a) \le \alpha < \beta \le ({\rm Re} \, F)(b), $$
insbes.\ trennt die abgeschlossene $\R$-Hyperebene $H_{{\rm Re} \, F, \alpha}$ die konvexen Mengen $A$ und $B$ strikt.
\end{itemize}
\end{HS}

Wir bereiten den Beweis des Hauptsatzes durch das folgende Lemma vor.

\begin{Lemma} \label{FA.2.24.L}
Seien $X$ ein normierter $\K$-Vektorraum und $C$ eine nicht-leere konvexe Teilmenge von $X$ sowie $\varepsilon \in \R_+$.

Dann ist auch die \emph{$\varepsilon$-Umgebung $U_{\varepsilon}(C) := \bigcup_{a \in C} U_{\varepsilon}(a)$ von $C$} eine konvexe Teilmenge von $X$.
\end{Lemma}

\textit{Beweis.} Seien $x,y \in U_{\varepsilon}(C)$.
Dann existieren $a,b \in C$ mit $x \in U_{\varepsilon}(a)$ sowie $y \in U_{\varepsilon}(b)$ und für jedes $t \in [0,1]$ gilt $(1-t) \, a + t \, b \in C$.
Nun folgt für $t \in [0,1]$
$$ \| (1-t) \, x + t \, y - ( \underbrace{(1-t) \, a + t \, b}_{\in C} ) \| \le (1-t) \, \underbrace{\| x-a\|}_{< \varepsilon} + t \, \underbrace{\|y-b\|}_{< \varepsilon} < \varepsilon, $$
also $(1-t) \, x + t \, y \in U_{\varepsilon}(C)$. \q
\A
\textit{Beweis des Hauptsatzes.} 1. Fall: $\K = \R$.

Zu (i): Sei also $A$ offen.
Es seien $a_0 \in A$ und $b_0 \in B$.
Wir setzen $x_0 := b_0 - a_0 \ne 0$ sowie
$$ C := x_0 + (A - B). $$
Wir behaupten:
\begin{gather}
C \in \U(0,X) \mbox{ ist konvex,} \label{FA.2.24.1} \\
x_0 \notin C. \label{FA.2.24.2}
\end{gather}

{[} Zu (\ref{FA.2.24.1}): Daß $C$ konvex ist, ergibt sich sofort aus den Beispielen 2.) und 3.) zu \ref{FA.N.2} (iv).
Weiterhin ist $0 = x_0 + a_0 - b_0 \in C$.

Für jedes $y \in X$ ist $T_y \: X \to X$, definiert durch
$$ \forall_{x \in X} \, T_y(x) := y + x, $$
offenbar eine bijektive stetige Abbildung mit Umkehrabbildung $T_{-y}$.
Insbes.\ ist $T_y$ ein Homöomorphismus, also eine offene Abbildung.
Wegen
$$ C = \bigcup_{b \in B} T_{x_0 - b}(A) $$
ist $C$ daher mit $A$ offen.

Zu (\ref{FA.2.24.2}): Wäre $x_0 \in C$, so existierten $a \in A$ und $b \in B$ mit $x_0 = x_0 + a - b$, also $a=b$, im Widerspruch zu $A \cap B = \emptyset$. {]}

Wegen (\ref{FA.2.24.1}) können wir das Minkowosi-Funktional $p_C \: X \to \R$ von $C$, vgl.\ \ref{FA.N.2} (i) inkl.\ Bemerkung 2.), (ii) und \ref{FA.N.4} (i) - (iii) inkl.\ Bemerkung 2.), betrachten, und aus (\ref{FA.2.24.2}) folgt $p_C(x_0) \ge 1$.
Wir setzen $Y := \R \, x_0$ und zeigen
\begin{equation}
f \: Y \longrightarrow \R, ~ t \, x_o \longmapsto t, \mbox{ ist eine Linearform mit } f \le p_C|_Y. \label{FA.2.24.3}
\end{equation}

{[} Zu (\ref{FA.2.24.3}): Die $\R$-Linearität ist trivial.
Des weiteren gilt für $t \in \R_+$
$$ f(t \, x_0) = t \, \stackrel{p_C(x_0) \ge 1}{\le} t \, p_C(x_0) = p_C(t \, x_0), $$
sowie für $t \in \R_-$
$$ f(t \, x_0) = t < 0 \stackrel{\ref{FA.N.4} (i)}{\le}  p_C(t \, x_0) $$
und außerdem $f(0) = 0 \stackrel{\ref{FA.N.4} (i)}{=} p_C(0)$. {]} 

\pagebreak
Aus dem Fortsetzungssatz von \textsc{Hahn-Banach} \ref{FA.2.17} ergibt sich wegen (\ref{FA.2.24.3}) die Existenz eines $F \in X^*$ mit $F|_Y = f$ und $F \le p_C$, also gilt offenbar $F \ne 0$ und nach \ref{FA.N.4} (iii) inkl.\ Bemerkung 2.) außerdem $F|_C \le p_C|_C < 1$.
Hieraus folgt, daß $H := H_{F,1}$ eine $\R$-Hyperebene von $X$ mit
$$ H \cap C = \emptyset $$
ist. 
Da $C$ nach (\ref{FA.2.24.1}), (\ref{FA.2.24.2}) eine nicht-leere offene Teilmenge von $X$ ist, folgt hieraus mittels \ref{FA.1.8} (ii): $\overline{H} \ne X$, d.h.\ wegen \ref{FA.2.21} (ii), daß $F$ stetig ist.
Es gilt also $F \in X'$.

Zum Nachweis von (i) genügt es nun zu zeigen 
\begin{gather}
\forall_{a \in A} \forall_{b \in B} \, F(a) < F(b), \label{FA.2.24.4} \\
\mbox{$F$ ist eine offene Abbildung}, \label{FA.2.24.5}
\end{gather}
denn (\ref{FA.2.24.4}) ergibt $\alpha := \inf \{ F(b) \, | \, b \in B \} \in \R$ sowie $\forall_{a \in A} \forall_{b \in B} \, F(a) \le \alpha \le F(b)$, und hieraus sowie (\ref{FA.2.24.5}) folgt zusätzlich, da $A$ offen ist, $\forall_{a \in A} \, F(a) < \alpha$.

{[} Zu (\ref{FA.2.24.4}): Seien $a \in A$ und $b \in B$, d.h.\ $x_0 + a - b \in C$, also wegen $F|_C < 1$ und $F(x_0) = f(x_0) = f(1 \, x_0) = 1$
$$ 1 > F(x_0 + a - b) = F(x_0) + F(a) - F(b) = F(a) - F(b). $$
Somit gilt $F(a) < F(b)$.

Zu (\ref{FA.2.24.5}): Da $F$ stetig ist, ist $Z := {\rm Kern} \, F$ ein abgeschlossener $\R$-Unter\-vek\-tor\-raum des normierten $\R$-Vektorraumes $X$.
Nach \ref{FA.2.5} (i), (ii) ist dann auch $X/Z$ ein normierter $\R$-Vektorraumraum und $\pi \: X \to X/Z$ eine offene Abbildung.
Außerdem ist $F$ wegen $F \ne 0$ surjektiv, daher folgt aus dem Homomorphiesatz, daß ein $\R$-Vektorraum-Isomorphismus $T \: X/Z \to \R$ existiert so, daß das folgende Diagramm kommutiert:
%
%\begin{diagram}
%X          & \rTo^{F}    & \R \\
%\dOnto^{\pi} & \ruTo>{T} &    \\
%X/Z        &             &    \\
%\end{diagram}
%
\begin{center}
\begin{tikzpicture}[node distance=2cm, every arrow/.style={thick,->}]
    % Knoten definieren
    \node (X) {$X$};
    \node (R) [right of=X] {$\R$};
    \node (X/Z) [below of=X] {$X/Z$};
    
    % Pfeile zeichnen
    \draw[-{>>}] (X) -- node[above] {$T$} (R);
    \draw[-{>>}] (X) -- node[left] {$\pi$} (X/Z);
    \draw[>-{>>}] (X/Z) -- node[right] {$S$} (R);
\end{tikzpicture}
\end{center}
Nun sind $\R$-lineare Abbildungen zwischen endlich-dimensionalen normierten $\R$-Vektorräumen stetig, vgl.\ \ref{FA.2.E.10}, also ist $T$ als Abbildung zwischen den normierten $\R$-Vektorräumen $X/Y$ und $\R$ ein Homöomorphismus, d.h.\ insbosondere eine offene Abbildung.
Somit ist auch $F = T \circ \pi$ eine offene Abbildung, und (\ref{FA.2.24.5}) ist gezeigt. {]}

Zu (ii): Sei nun $\delta = d(A,B) \in \R_+$.
Mit $A$ ist nach Lemma \ref{FA.2.24.L} auch die offene Menge $U_{\delta}(A) := \bigcup_{a \in A} U_{\delta}(a)$ konvex, und es gilt $U_{\delta}(A) \cap B = \emptyset$.
Aus dem Beweis von (i), insbes.\ (\ref{FA.2.24.4}), (\ref{FA.2.24.5}) -- angewandt auf $U_{\delta}(A)$ anstelle von $A$ --, folgt die Existenz eines $F \in X' \setminus \{0\}$ mit
\begin{equation}
\forall_{\tilde{a} \in U_{\delta}(A)} \forall_{b \in B} \, F(\tilde{a}) < \underbrace{\inf \{ F(\tilde{b}) \, | \, \tilde{b} \in B \}}_{=: \beta \in \R} \le F(b). \label{FA.2.24.6}
\end{equation}
Wir setzen
$$ \alpha := \sup \{ F(a) \, | \, a \in A \} \in \R $$
und haben $\alpha < \beta$ zu zeigen:
$F \in X' \setminus \{0\}$ ist eine offene Abbildung, vgl.\ (\ref{FA.2.24.5}), also ist $F(U_{\delta}(0))$ ein offenes Intervall, das $0$ enthält, d.h.\ $\rho := \sup F(U_{\delta}(0)) \in \R_+$.
Es genügt nachzuweisen, daß 
$$ \alpha + \rho \le \beta. $$

Beweis hiervon:
Sei $\varepsilon \in \R_+$.
Nach Definition von $\alpha, \rho$ existieren $a \in A$ und $x \in U_{\delta}(0)$ mit
$$ F(a) \ge \alpha - \varepsilon ~~ \wedge ~~ F(x) \ge \rho - \varepsilon, $$
also folgt mittels (\ref{FA.2.24.6}) aus $a+x \in U_{\varepsilon}(A)$
$$ \beta > F(a+x) = F(a) + F(x) \ge \alpha + \rho - 2 \varepsilon. $$
Die Beliebigkeit von $\varepsilon \in \R_+$ ergibt $\beta \ge \alpha + \rho$.

2. Fall: $\K = \C$.
Analog zum Beweis von \ref{FA.2.18} wendet man den 1. Fall auf $X_{\R}$ an, erhält dann $F_1 \in (X_{\R})' \setminus \{0\}$ mit den entsprechenden reellen Eigenschaften und zeigt, daß $F \: X \to \C$, definiert durch
$$ \forall_{x \in X} \, F(x) := F_1(x) - \i \, F_1( \i \, x ), $$
die Behauptung erfüllt. \q

\subsection*{Transponierte Operatoren und der topologische Bidualraum} \addcontentsline{toc}{subsection}{Transponierte Operatoren und der topologische Bidualraum}

\begin{Def} \label{FA.2.TB.1}
Seien $X,Y$ normierte $\K$-Vektorräume.
Dann induziert jedes $T \in \mathcal{L}_{\K}(X,Y)$ ein $\boxed{T' \: Y' \to X'}$, den sog.\ \emph{zu $T$ transponierten Operator}\index{Operator!transponierter}, der durch
$$ \forall_{F \in Y'} \, T'\, F := F \circ T $$
gegeben ist.
\end{Def}

\begin{Satz} \label{FA.2.TB.2}
Seien $X,Y$ normierte $\K$-Vektorräume und $T,S \in \mathcal{L}_{\K}(X,Y)$ sowie $\lambda \in \K$.
Dann gilt:
\begin{itemize}
\item[(i)] $T \in \mathcal{L}_{\K}(Y',X')$ und $\|T'\| = \|T\|$.
\item[(ii)] $(T + S)' = T' + S'$ und $(\lambda \, T)' = \lambda \, T'$.
\end{itemize}
\end{Satz}

\textit{Beweis.} Zu (i): Seien $F,G \in Y'$ und $\lambda \in \K$.
Dann gilt
\begin{gather*}
T' (F+G) = (F+G) \circ T = (F \circ T) + (G \circ T) = T' \, F + T' \, G, \\
T' (\lambda \, F) = (\lambda F) \circ T = \lambda \, (F \circ T) = \lambda \, (T'\, F), \\
\|T' \, F\| = \|F \circ T\| \le \|F\| \, \|T\| = \|T\| \, \|F\|,
\end{gather*}
also folgt $T \in \mathcal{L}_{\K}(Y',X')$ und $\|T'\| \le \|T\|$.

Zu zeigen bleibt $\|T\| \le \|T'\|$.
Hierzu sei $x \in X$ mit $T(x) \ne 0$.
Nach \ref{FA.2.20} (ii) existiert $F \in Y'$ mit $F(T(x)) = \|T(x)\| \in \R_+$ und $\|F\| = 1$, also ergibt sich
$$ \|T(x)\| = |F(T(x))| \le \|F \circ T\| \, \|x\| = \|T' \, F\|  \, \|x\| \le \|T'\| \, \|F\| \, \|x\| = \|T'\| \, \|x\|. $$
Da $\|T(x)\| \le \|T'\| \, \|x\|$ natürlich auch im Falle $T(x) = 0$ gilt, folgt $\|T\| \le \|T'\|$.

Zu (ii): Für alle $F \in Y'$ und $x \in X$ gilt
\begin{eqnarray*}
((T+S)' \, F)(x) & = & (F \circ (T+S))(x) = F((T+S)(x)) = F(T(x) + S(x)) \\
& = & F(T(x)) + F(S(x)) =  (F \circ T)(x) + (F \circ S)(x) \\
& = & (T' \, F)(x) + (S' \, F)(x)
\end{eqnarray*}
und
\begin{eqnarray*}
((\lambda \, T')(F))(x) & = & (F \circ (\lambda \, T))(x) = F(\lambda \, T(x)) = \lambda \, F(T(x)) = \lambda \, (F \circ T)(x) \\
& = & \lambda \, (T' \, F)(x), 
\end{eqnarray*}
womit der Satz bewiesen ist. \q

\begin{Satz} \label{FA.2.TB.3}
Die Zuordnung
$$ X \longmapsto X', ~~ T \longmapsto T', $$
definiert einen \emph{kontravarianten Funktor $\boxed{\,'}$ von der Kategorie der normierten $\K$-Vektorräume mit ihren stetigen linearen Abbildungen in die Kategorie der $\K$-Banachräume mit ihren stetigen linearen Abbildungen}. 
D.h.\ per definitionem, daß für alle normierten $\K$-Vektorräume $X,Y,Z$ und alle $T \in \mathcal{L}_{\K}(X,Y)$ sowie $S \in \mathcal{L}_{\K}(Y,Z)$ gilt:
\begin{itemize}
\item[(i)] $X'$ ist ein $\K$-Banachraum.
\item[(ii)] $T' \in \mathcal{L}_{\K}(Y',X')$.
\item[(iii)] $(\id_X)' = \id_{X'}$.
\item[(iv)] $(S \circ T)' = T' \circ S'$.
\end{itemize}
\end{Satz}

\textit{Beweis.} (i) sowie (ii) haben wir bereits eingesehen, und (iii) ist klar.

Zu (iv): Für alle $F \in Z'$ gilt
$$ (S \circ T)' \, F = F \circ (S \circ T) = (F \circ S) \circ T = T'(S'(F)) = (T' \circ S') \, F, $$
d.h.\ (iv) ist gezeigt. \q

\begin{Bem*}
In \cite[Satz 7.10]{Hirz} wird behauptet, der Funktor $'$ sei exakt, d.h., daß er kurze exakte Sequenzen auf kurze exakte Sequenzen abbilde.%
\footnote{Sind $n \in \N$, $X_0, \ldots, X_{n+1}$ $\K$-Vektorräume und $T_i \: X_i \to X_{i+1}$ $\K$-lineare Abbildungen für $i \in \{0, \ldots, n\}$, so heißt die Sequenz
$$ X_0 \overset{T_0}{\longrightarrow} X_1 \overset{T_1}{\longrightarrow} \ldots \overset{T_n}{\longrightarrow} T_{n+1} $$
\emph{exakt} genau dann, wenn gilt $\forall_{i \in \{0, \ldots, n-1\}} \, T_i(X_i) = \overline{T_{i+1}}^1(\{0\})$.
Im Falle $n=2$ spricht man auch von einer \emph{kurzen exakten Sequenz}.\index{exakte Sequenz}}
Diese Behauptung ist falsch.%
\footnote{Mit den Bezeichnungen wie in \cite[S.\ 36 f.]{Hirz} wird beim ,,Beweis`` von \cite[Satz 7.10]{Hirz} $Y/\alpha(X)$ mit $\beta(Y)$ identifiziert. 
$Y/\alpha(X)$ ist jedoch i.a.\ nicht homöomorph zu $\beta(Y)$.
Ist zusätzlich $Y$ ein $\K$-Banachraum und $\beta(Y)$ vollständig, so liefert der Satz vom inversen Operator (s.u.\ \ref{FA.2.37}) diese Homöomorphie.}
Seien nämlich z.B.\ $X := (\mathcal{C}([0,1],\R), \|\ldots\|_{\infty})$ und $T \: X \to X$ definiert durch $\forall_{f \in \mathcal{C}([0,1],\R)} \, T \, f := \int_0^x f(t) \, \d t$.
Dann ist $T \in \mathcal{L}_{\R}(X,X)$ injektiv, d.h.\ $\{0\} \rightarrow X \overset{T}{\rightarrow} X$ ist exakt.
$T' \: X' \to X'$ ist aber nicht surjektiv, denn $F \in X'$, definiert durch $\forall_{f \in \mathcal{C}([0,1],\R)} \, F(f) := f(1)$, liegt nicht im Bild von $T'$.
Daher ist $X' \overset{T'}{\rightarrow} X' \rightarrow \{0\}$ nicht exakt.
\end{Bem*}

\begin{Kor} \label{FA.2.TB.3.K}
Seien $X,Y$ normierte $\K$-Vektorräume und $T \in \mathcal{L}_{\K}(X,Y)$.
Dann gilt:
\begin{itemize}
\item[(i)] $T \: X \longrightarrow Y$ $\K$-Vektorraum-Isomorphismus
\newline $\Longrightarrow$ $T' \: Y' \longrightarrow X'$ $\K$-Vektorraum-Isomorphismus.
\item[(ii)] $T \: X \longrightarrow Y$ $\K$-Vektorraum-Isometrie
\newline $\Longrightarrow$ $T' \: Y' \longrightarrow X'$ $\K$-Vektorraum-Isometrie.
\end{itemize}
\end{Kor}

\textit{Beweis.} (i) folgt sofort aus \ref{FA.2.TB.3} (iii) und (iv).

Zu (ii): Im Falle $X = \{0\}$ ist die Behauptung klar.
Sei daher $\dim_{\K} X \ne 0$ und $T \: X \to Y$ eine $\K$-Vektorraum-Isometrie.
Als surjektive isometrische Abbildung bildet $T$ die Einheitskugel $\{ x \in X \, | \, \|x\|=1 \}$ in $X$ surjektiv auf die Einheitskugel $\{ y \in Y \, | \, \|y\|=1 \}$ in $Y$ ab, also folgt für jedes $F \in Y'$: $\|T' \, F\| = \|F \circ T\| = \|F\|$. \q

\begin{Def}[Topologischer Bidualraum] \label{FA.2.TB.4}
Sei $X$ ein normierter $\K$-Vektor\-raum.
\begin{itemize}
\item[(i)] Der $\K$-Banachraum
$$ \boxed{X''} := (X')' = \mathcal{L}_{\K}(\mathcal{L}_{\K}(X,\K),\K) $$
heißt der \emph{topologische Bidualraum von $X$}.
\item[(ii)] Sind $Y$ ein weiterer normierter $\K$-Vektorraum und $T \in \mathcal{L}_{\K}(X,Y)$, so setzen wir
$$ \boxed{T''} := (T')' \in \mathcal{L}_{\K}(X'',Y''). $$
Es gilt also $\forall_{\Phi \in X''} \forall_{F \in X'} \, (T'' \Phi)(F) = (\Phi \circ T')(F) = \Phi(T'(F)) = \Phi(F \circ T)$.
\end{itemize}
\end{Def}

\begin{Lemma} \label{FA.2.TB.5}
Seien $X$ ein normierter $\K$-Vektorraum und $x \in X$.

Dann ist die \emph{Evaluationsabbildung von $x$}
$$ \boxed{E_x \: X' \longrightarrow \K, ~~ F \longmapsto F(x),} $$
$\K$-linear und stetig, d.h.\ $E_x \in X''$, mit $\|E_x\| = \|x\|$.
\end{Lemma}

\textit{Beweis.} Die $\K$-Linearität von $E_x$ ist trivial.

Für jedes $F \in X'$ gilt
$$ |E_x(F)| = |F(x)| \le \|F\| \, \|x\| =  \|x\| \, \|F\|, $$
also folgt die Stetigkeit von $E_x \: X' \to \K$ und $\|E_x\| \le \|x\|$.

Aus der letzten Ungleichung ergibt sich sofort, daß in ihr im Falle $x=0$ Gleichheit gelten muß.
Im Falle $x \ne 0$ existiert nach \ref{FA.2.20} (ii) ein Funktional $F \in X'$ mit $F(x) = \|x\| \in \R_+$ und $\|F\| = 1$, also folgt
$$ |E_x(F)| = |F(x)| = \|x\| = \|x\| \, \|F\|, $$
d.h.\ $\|E_x\| \ge \|x\|$. \q

\begin{Satz} \label{FA.2.TB.6}
Sei $X$ ein normierter $\K$-Vektorraum.

Dann ist 
$$ \boxed{i_X \: X \longrightarrow X'', ~~ x \longmapsto E_x,} $$
eine $\K$-lineare isometrische -- also insbes.\ stetige und injektive -- Abbildung.
\end{Satz}

\textit{Beweis.} 1.) Aus \ref{FA.2.TB.5} folgt $\forall_{x \in X} \, E_x \in X''$, d.h.\ $i_X$ ist wohldefiniert.

2.) Die $\K$-Linearität von $i_X$ ist trivial.

3.) Wegen \ref{FA.2.TB.5} gilt außerdem $\forall_{x \in X} \, \|i_X(x)\| = \|E_x\| = \|x\|$.
Daher ist $i_X$ wegen 2.) eine isometrische Abbildung. \q

\begin{Kor} \label{FA.2.TB.7}
Sei $X$ ein normierter $\K$-Vektorraum.

Dann ist $\overline{i_X(X)} \subset X''$ eine Vervollständigung des zu $X$ $\K$-linear isometrischen $\K$-Vektor\-raumes $i_X(X)$.
\end{Kor}

\textit{Beweis.} $X''$ ist vollständig und $\overline{i_X(X)}$ ist abgeschlossen in $X''$.
Daher folgt die Behauptung aus \ref{FA.1.11} (i). \q

\begin{Def} \label{FA.2.TB.8}
Ein normierter $\K$-Vektorraum heißt \emph{reflexiv}\index{Raum!reflexiver Vektor-} genau dann, wenn $i_X \: X \to X''$ surjektiv, d.h.\ genau eine $\K$-Vektorraum-Isometrie, ist.

\begin{Bsp*}
Jeder endlich-dimensionale $\K$-Vektorraum ist reflexiv.
\end{Bsp*}
\end{Def}

\begin{Satz} \label{FA.2.TB.9}
Es seien $X,Y$ normierte $\K$-Vektorräume und $T \in \mathcal{L}_{\K}(X,Y)$.

Dann gilt $\|T''\| = \|T\|$ und $T'' \circ i_X = i_Y \circ T$.
Letzteres bedeutet, daß das folgende Diagramm kommutiert:
%
%\begin{diagram}
%X             & \rTo^{T}   & Y             \\
%\dEmbed^{i_X} &            & \dEmbed_{i_Y} \\
%X''           & \rTo^{T''} & Y''           \\
%\end{diagram}
%
\begin{center}
\begin{tikzpicture}[node distance=2cm, every arrow/.style={thick,->}]
    % Knoten definieren
    \node (X) {$X$};
    \node (Y) [right of=X] {$Y$};
    \node (X'') [below of=X] {$X''$};
    \node (Y'') [right of=X''] {Y''};
    
    % Pfeile zeichnen
    \draw[->] (X) -- node[above] {$T$} (Y);
    \draw[>->] (X) -- node[left] {$i_X$} (X'');
    \draw[>->] (Y) -- node[right] {$i_Y$} (Y'');
    \draw[->] (X'') -- node[above] {$T''$} (Y'');
\end{tikzpicture}
\end{center}
\end{Satz}

\textit{Beweis.} $\|T''\| = \|T'\| = \|T\|$ folgt aus \ref{FA.2.TB.2} (i), und für alle $x \in X$ und $F \in X'$ gilt
\begin{eqnarray*}
(T''(i_X(x)))(F) & = & (T'' \, E_x)(F) = E_x(F \circ T) = F(T(x))) = (E_{T(x)})(F) \\
& = &(i_Y(T(x)))(F),
\end{eqnarray*}
also $(T'' \circ i_X)(x) = (i_Y \circ T)(x)$. \q   

\subsection*{Die Prinzipien der gleichmäßigen Beschränktheit und der offenen Abbildung} \addcontentsline{toc}{subsection}{Die Prinzipien der gleichmäßigen Beschränktheit und der offenen Abbildung}

Der folgende Satz wird in der Literatur auch als \emph{das Prinzip der gleichmäßigen Beschränktheit} bezeichnet.

\begin{HS}[Satz von \textsc{Banach-Steinhaus}] \index{Satz!von \textsc{Banach-Steinhaus}} \index{Prinzip!der gleichmäßigen Beschränktheit} \label{FA.2.25} $\,$

\noindent \textbf{Vor.:} Es seien $X,Y$ normierte $\K$-Vektorräume, $\mathcal{H} \subset \mathcal{L}_{\K}(X,Y)$ und 
$$ B := \{ x \in X \, | \, \underbrace{\mathcal{H}(x)}_{= \{ T(x) \, | \, T \in \mathcal{H} \}} \mbox{ beschränkt} \}. $$

\noindent \textbf{Beh.:} $B$ nicht mager in $X$ $\Longrightarrow$ $\mathcal{H}$ beschränkt in $\mathcal{L}_{\K}(X,Y)$.
\end{HS}

\begin{Bem*}
Sind $X,Y$ normierte $\K$-Vektorräume, so ist $\mathcal{H} \subset \mathcal{L}_{\K}(X,Y)$ genau dann beschränkt, wenn $C \in \R_+$ mit
$$ \forall_{T \in \mathcal{H}} \forall_{x \in X} \, \|T(x)\| \le C $$
existiert, d.h.\ $\mathcal{H}$ ist \emph{gleichmäßig} beschränkt.
\end{Bem*}

\textit{Beweis.} Wir zeigen die zur Behauptung äquivalente Aussage
$$ \sup \{ \|T\| \, | \, T \in \mathcal{H} \} = \infty \Longrightarrow \underbrace{B \mbox{ mager in } X}_{\Leftrightarrow X \setminus B \text{ generisch in } X}. $$

Für jedes $n \in \N$ definieren wir eine offene Menge $U_n$ durch
$$ U_n := \{ x \in X \, | \, \exists_{T \in \mathcal{H}} \, \|T(x)\| > n \} = \bigcup_{T \in \mathcal{H}} \{ x \in X \, | \, \|T(x)\| > n \} $$
und behaupten
\begin{equation} \label{FA.2.25.1}
U_n \mbox{ ist dicht in } X.
\end{equation}

{[} Zu (\ref{FA.2.25.1}): Angenommen, $U_n$ ist nicht dicht in $X$, d.h.\ es existieren $x_0 \in X$ und $\varepsilon \in \R_+$ mit
$$ U_{\varepsilon}(x_0) \cap U_n = \emptyset. $$
Dann folgt $T \tilde{x} \le n$ für jedes $\tilde{x} \in U_{\varepsilon}(x_0)$ sowie $T \in \mathcal{H}$, also für jedes $x \in U_{\varepsilon}(0)$ sowie $T \in \mathcal{H}$
$$ \|T(x)\| \le \|T(x+x_0)\| + \|T(x_0)\| \le 2n, $$
und für jedes $T \in \mathcal{H}$ sowie $x \in X \setminus \{0\}$ gilt
$$ \|T(x)\| = \frac{2 \|x\|}{\varepsilon} \, \left\| T \left( \frac{\varepsilon \, x}{2 \|x\|} \right) \right\| \le \frac{2 \|x\|}{\varepsilon} \, 2n = \frac{4n}{\varepsilon} \, \|x\|. $$
Wegen $\|T(0)\| \le \frac{4n}{\varepsilon} \, \|0\|$ bedeutet die letzte Ungleichung
$$ \forall_{T \in H} \, \|T\| \le \frac{4n}{\varepsilon}, $$
im Widerspruch zu $\sup \{ \|T\| \, | \, T \in \mathcal{H} \} = \infty$. {]}

Wir haben damit bewiesen -- vgl.\ \ref{FA.1.39} --, daß
$$ \bigcap_{n \in \N} U_n = \{ x \in X \, | \, \mathcal{H}(x) \mbox{ unbeschränkt} \} = X \setminus B $$
generisch ist. \q

\begin{Satz} \label{FA.2.26} $\,$

\noindent \textbf{Vor.:} Es seien $X$ ein $\K$-Banachraum, $Y$ ein normierter $\K$-Vektorraum und $\mathcal{H} \subset \mathcal{L}_{\K}(X,Y)$ sowie $B := \{ x \in X \, | \, \mathcal{H}(x) \mbox{ beschränkt} \}$.

\noindent \textbf{Beh.:} $B^{\circ} \ne \emptyset$ $\Longrightarrow$ $\mathcal{H}$ beschränkt in $\mathcal{L}_{\K}(X,Y)$.%
\footnote{$B^{\circ} \ne \emptyset$ ist z.B.\ erfüllt, wenn $\mathcal{H}(x)$ für jedes $x \in X$ beschränkt ist.}
\end{Satz}

\textit{Beweis.} $X$ ist vollständig, also Bairesch nach dem Satz von Baire \ref{FA.1.44}.
Daher folgt -- vgl.\ \ref{FA.1.42}, \ref{FA.1.41} (ii) --, daß $B$ nicht mager ist.
Somit ergibt sich die Behauptung aus dem Satz von \textsc{Banach-Steinhaus} \ref{FA.2.25}. \q

\begin{Kor} \label{FA.2.27}
Seien $X$ ein $\K$-Banachraum, $Y$ ein normierter $\K$-Vektorraum und $(T_n)_{n \in \N}$ eine Folge in $\mathcal{L}_{\K}(X,Y)$, die punktweise gegen $T \in Y^X$ konvergiert, d.h.\
\begin{equation} \label{FA.2.27.0}
\forall_{x \in X} \, \lim_{n \to \infty} T_n(x) = T(x).
\end{equation}

Dann gilt $T \in \mathcal{L}_{\K}(X,Y)$ und $\|T\| \le \liminf_{n \to \infty} \|T_n\| < \infty$.
\end{Kor}

\textit{Beweis.} Die $\K$-Linearität von $T$ ist trivial.
Zum Nachweis der Stetigkeit von $T$ wenden wir \ref{FA.2.26} auf
$$ \mathcal{H} := \{ T_n \, | \, n \in \N \} \subset \mathcal{L}_{\K}(X,Y) $$
an.
Wegen (\ref{FA.2.27.0}) gilt $\{ x \in X \, | \, T_n(x) \mbox{ beschränkt} \} = X$, also folgt die Existenz einer reellen Zahl $C \in \R_+$ derart, daß $\forall_{n \in \N} \, \|T_n\| \le C$ gilt.
Hieraus folgt zum einen 
\begin{equation} \label{FA.2.27.1}
\liminf_{n \to \infty} \|T_n\| \le C
\end{equation}
und wegen der Stetigkeit von $\| \ldots \| \: Y \to \R$ zum anderen für jedes $x \in X$
$$ \|T(x)\| = \lim_{n \to \infty} \underbrace{\|T_n(x)\|}_{\le \|T_n\| \, \|x\| \le C \, \|x\|} \le C \, \|x\|, $$
also ist $T$ stetig, und es gilt des weiteren
\begin{eqnarray*}
\lim_{k \to \infty} \inf \{ \overbrace{\|T_k(x)\|}^{\le \|T_k\| \, \|x\|} \, \| \, n \in \N \, \wedge \, n \ge k \} & \le & \lim_{k \to \infty} \|x\| \, \inf \{ \|T_k\| \, | \, n \in \N \, \wedge \, n \ge k \} \\
& = & \|x\| \, \lim_{k \to \infty} \inf \{ \|T_k\| \, | \, n \in \N \, \wedge \, n \ge k \},
\end{eqnarray*}
also
$$ \|T(x)\| = \lim_{n \to \infty} \|T_n(x)\| = \liminf_{n \to \infty} \|T_n(x)\| \le \|x\| \, \liminf_{n \to \infty} \|T_n\|, $$
d.h.\ $\|T\| \le \liminf_{n \to \infty} \|T_n\| \stackrel{(\ref{FA.2.27.1})}{<} \infty$. \q

\begin{Kor} \label{FA.2.28}
Seien $X$ ein normierter $\K$-Vektorraum und $B$ eine Teilmenge von $X$ derart, daß gilt
$$ \forall_{F \in X'} \, F(B) \mbox{ beschränkte Teilmenge von } \K. $$

Dann ist $B$ eine beschränkte Teilmenge von $X$.
\end{Kor}

\textit{Beweis.} Wir wenden \ref{FA.2.26} auf
$$ \mathcal{H} :=  i_X(B) = \{ E_b \, | \, b \in B \} \subset X'' = \mathcal{L}_{\K}(X',\K) $$
an.
Beachte, daß $X'$ ein $\K$-Banachraum ist.
Es gilt $\mathcal{H}(F) = F(B)$ für jedes $F \in X'$, also nach Voraussetzung $\{ F \in X' \, | \, \mathcal{H}(F) \mbox{ beschränkt} \} = X$.
Somit folgt die Beschränktheit von $\mathcal{H} = i_X(B)$.
Da $i_X \: X \to X''$ eine isometrische Abbildung ist, ist auch $B$ beschränkt. \q
\pagebreak

Obwohl wir den Begriff der ,,Offenheit einer Abbildung`` bereits verwendet haben, geben wir der Vollständigkeit halber die Definition erneut an.

\begin{Def}[Offene Abbildungen] \index{Abbildung!offene} \label{FA.2.29}
Seien $X,Y$ topologische Räume und $f \in Y^X$.
Dann heißt $f$ eine \emph{offene Abbildung}, wenn für jedes $U \in {\rm Top}(X)$ gilt $f(U) \in {\rm Top}(Y)$.
\end{Def}

\begin{Satz} \label{FA.2.30}
Es seien $X,Y$ normierte $\K$-Vektorräume und $T \in {\rm Hom}_{\K}(X,Y)$.
Dann sind die folgenden drei Aussagen paarweise äquivalent:
\begin{itemize}
\item[(i)] $T$ ist eine offene Abbildung.
\item[(ii)] $\exists_{C \in \R_+} \forall_{y \in Y} \exists_{x \in X} \, T(x) = y \, \wedge \, C \, \|x\| \le \|y\|$.
\item[(iii)] $\exists_{\delta \in \R_+} \, U_{\delta}(0_Y) \subset T(U_1(0_X))$.
\end{itemize}
\end{Satz}

\textit{Beweis.} ,,(i) $\Rightarrow$ (iii)`` ist trivial.

,,(iii) $\Rightarrow$ (i)`` Seien $U \subset X$ offen und $y_0 \in T(U)$.
Dann existieren $x_0 \in U$ und $\varepsilon \in \R_+$ mit $T(x_0) = y_0$ sowie $U_{\varepsilon}(x_0) \subset U$.
Ist $\delta \in \R_+$ wie in (iii), so gilt für alle $y \in U_{\varepsilon \, \delta}(y_0)$
$$ \frac{1}{\varepsilon}(y-y_0) \in U_{\delta}(0_Y) \subset T(U_1(0_X)), $$
also existiert $x \in U_1(0_X)$ mit $T(x) = \frac{1}{\varepsilon}(y-y_0)$, d.h.\
$$ y = y_0 + \varepsilon \, T(x) = T(x_0) + \varepsilon \, T(x) = T(x_0 + \varepsilon \, x) \in T(\underbrace{U_{\varepsilon}(x_0)}_{\subset U}) \subset T(U). $$
Damit ist $U_{\varepsilon \, \delta}(y_0) \subset T(U)$ gezeigt.
Wegen der Beliebigkeit von $y_0 \in T(U)$ folgt die Offenheit von $T(U)$ in $Y$.

,,(iii) $\Rightarrow$ (ii)`` Sei $\delta$ gemäß (iii) gewählt.
Wir setzen $C := \frac{\delta}{2} \in \R_+$.
Für $y = 0_Y$ ist (ii) trivial.
Sei daher $y \in Y \setminus \{0\}$.
Dann gilt 
$$ \frac{\delta}{2} \, \frac{y}{\|y\|} \in U_{\delta}(0_Y) \subset T(U_1(0_X)), $$
also existiert $\tilde{x} \in U_1(0_X)$ mit $T(\tilde{x}) = \frac{\delta}{2} \, \frac{y}{\|y\|}$, also gilt mit $x := \frac{2 \|y\|}{\delta} \tilde{x}$
$$ T(x) = y \, \wedge \, \|x\| = \frac{2 \|y\|}{\delta} \|\tilde{x}\| < \frac{1}{C} \, \|y\|. $$

,,(ii) $\Rightarrow$ (iii)`` Seien $C \in \R_+$ wie in (ii) und $\delta := C$.
Dann existiert zu jedem $y \in U_{\delta}(0_Y)$ gemäß (ii) ein $x \in X$ mit 
$$ T(x) = y ~~ \wedge ~~ \delta \, \|x\| \le \|y\| < \delta, $$
also gilt $y = T(x)$ und $\|x\| < 1$. \q

\begin{Kor} \label{FA.2.31}
Es seien $X,Y$ normierte $\K$-Vektorräume und $T \in {\rm Hom}_{\K}(X,Y)$.

Ist $T$ eine offene Abbildung, so ist $T$ surjektiv.
\end{Kor}

\textit{Beweis.} Die Behauptung folgt trivial aus \ref{FA.2.30} ,,(i) $\Rightarrow$ (ii)``. \q
\pagebreak

\begin{Def}[Graph eines Operators] \label{FA.2.32} 
Seien $X,Y$ normierte $\K$-Vektor\-räume, $W$ ein $\K$-Unter\-vektor\-raum von $X$ und $T \: W \to Y$ eine $\K$-lineare Abbildung.
Dann heißt der $\K$-Untervektorraum
$$ \boxed{\Gamma_T} := \{ (x, T(x)) \in W \times Y \, | \, x \in W \} $$
von $X \times Y$ der \emph{Graph von $T$}\index{Graph}.

Sprechen wir im folgenden davon, daß \emph{$\Gamma_T$ abgeschossen in $X \times Y$} ist, so meinen wir dies bzgl.\ der Maximumsnorm auf $X \times Y$.

\begin{Bem*}
In der funktionalanalytischen Literatur wird $T$ wie oben gelegentlich \emph{abgeschlosser Operator} genannt, wenn ihr Graph in $X \times Y$ abgeschlossen ist.
Um Verwechselungen zu vermeiden, verwenden wir diese Terminologie nicht, da man eine Abbildung zwischen topologischen Räumen grundsätzlich \emph{abgeschlossen}\index{Abbildung!abgeschlossene} nennt, wenn sie abgeschlossene Mengen auf solche abbildet.
\end{Bem*}
\end{Def}

\begin{Lemma} \label{FA.2.33}
Seien $X,Y$ normierte $\K$-Vektorräume, $W$ ein $\K$-Unter\-vektor\-raum von $X$ und $T \: W \to Y$ eine $\K$-lineare Abbildung.
Dann sind die folgenden Aussagen äquivalent:
\begin{itemize}
\item[(i)] $\Gamma_T$ ist abgeschlossen in $X \times Y$.
\item[(ii)] Ist $(x_n)_{n \in \N}$ eine Folge in $W$ derart, daß $x \in X$ und $y \in Y$ mit $\lim_{n \to \infty} x_n = x$ und $\lim_{n \to \infty} T(x_n) = y$ existieren, so gilt $x \in W$ und $T(x) = y$.
\end{itemize}
\end{Lemma}

\textit{Beweis.} ,,(i) $\Rightarrow$ (ii)`` Sind $(x_n)_{n \in \N}, x, y$ wie in (ii), so ist $(x_n, T(x_n))_{n \in \N}$ eine Folge in der abgeschlossenen Menge $\Gamma_T$, die gegen $(x,y) \in X \times Y$ konvergiert.
Dann folgt $(x,y) \in \Gamma_T$, d.h.\ $x \in W$ und $T(x) = y$.

,,(ii) $\Rightarrow$ (i)`` Ist $(x_n, T(x_n))_{n \in \N}$ eine Folge in $\Gamma_T$, die gegen $(x,y) \in X \times Y$ konvergiert, so folgt offenbar aus (ii), daß gilt $(x,y) \in \Gamma_T$. \q 

\begin{Satz} \label{FA.2.34}
Seien $X,Y$ normierte $\K$-Vektor\-räume, $W$ ein $\K$-Unter\-vektor\-raum von $X$ und $T \: W \to Y$ eine $\K$-lineare Abbildung.
\begin{itemize}
\item[(i)] $\Gamma_T$ abgeschlossen in $X \times Y$ $\Longrightarrow$ ${\rm Kern} \, T$ abgeschlossen in $X$.
\item[(ii)] $W$ abgeschlossen in $X$ und $T \in \mathcal{L}_{\K}(W,Y)$ $\Longrightarrow$ $\Gamma_T$ abgeschlossen in $X \times Y$. 
\end{itemize}
\end{Satz}

\textit{Beweis.} Zu (i): Seien $\Gamma_T$ abgeschlossen und $(x_n)_{n \in\N}$ eine Folge in ${\rm Kern} \, T$ sowie $x \in X$ mit $\lim_{n \to \infty} x_n = x$.
Es folgt $\lim_{n \to \infty} (x_n,T(x_n)) = (x,0)$, also ergibt die Abgeschlossenheit von $\Gamma_T$ zusammen mit \ref{FA.2.33} ,,(i) $\Rightarrow$ (ii)`` sowohl $x \in W$ als auch $T(x) = 0$, d.h.\ $x \in {\rm Kern} \, T$.

Zu (ii): Seien $W$ abgeschlossen in $X$, $T \in \mathcal{L}_{\K}(W,Y)$ und $(x_n)_{n \in \N}$ eine Folge in $W$, die gegen $x \in X$ konvergiert.
Aus der Abgeschlossenheit von $W$ folgt $x \in W$ und aus der Stetigkeit von $T$, daß $\lim_{n \to \infty} T(x_n) = T(x)$ gilt, also liefert \ref{FA.2.33} ,,(ii) $\Rightarrow$ (i)`` die Abgeschlosssenheit von $\Gamma_T$. \q

\begin{Bsp} \label{FA.2.32.B}
Der $\R$-Banachraum $(\mathcal{C}([0,1], \R), \|\ldots\|_{\infty})$ besitzt den $\R$-Unter\-vektor\-raum $\mathcal{C}^1([0,1], \R)$ der stetig differenzierbaren Funktionen $[0,1] \to \R$, und der $\R$-lineare \emph{Differentialoperator}
$$ \frac{\d}{\d x} \: \mathcal{C}^1([0,1], \R) \longrightarrow  \mathcal{C}([0,1], \R) $$
hat einen abgeschlossenen Graphen, ist aber nicht stetig.
\pagebreak

{[} Verwende beim Nachweis der Abgeschlossenheit des Graphen des Differentialoperators \ref{FA.2.34} ,,(ii) $\Rightarrow$ (i)`` und beim Nachweis der Unstetigkeit, daß für jedes $n \in \N_+$ 
$$ \frac{\d}{\d x} \, x^n|_{[0,1]} = n \, x^{n-1}|_{[0,1]} $$
gilt. {]}
\end{Bsp}

\begin{Satz} \label{2.35.Z}
Seien $X,Y$ normierte $\K$-Vektor\-räume und $T \in \mathcal{L}_{\K}(X,Y)$ injektiv.

Dann ist $\Gamma_{T^{-1} \: T(X) \to X}$ abgeschlossen in $Y \times X$.
\end{Satz}

\textit{Beweis.} Seien $(y_n)_{n \in \N}$ eine Folge in $T(X)$ und $y \in Y$ mit $\lim_{n \to \infty} y_n = y$.
Dann existiert wegen der Injektivität von $T$ eine eindeutig bestimmte Folge $(x_n)_{n \in \N}$ in $X$ mit $\forall_{n \in \N} \, T(x_n) = y_n$, also $\forall_{n \in \N} \, T^{-1}(y_n) = x_n$.
Sei ferner $x \in X$ mit $\lim_{n \to \infty} x_n = x$.
Da $T$ stetig ist, gilt 
$$ \underbrace{\lim_{n \to \infty} \underbrace{T(x_n)}_{= y_n}}_{= y} = T(x), $$
d.h.\ $y \in T(X)$ und $T^{-1}(y) = x$.
Aus \ref{FA.2.33} ,,(ii) $\Rightarrow$ (i)`` ergibt sich daher die Behauptung. \q

\begin{HS}[Das Prinzip der offenen Abbildung] \index{Prinzip!der offenen Abbildung} \label{FA.2.35} $\,$

\noindent \textbf{Vor.:} Seien $X$ ein $\K$-Banachraum, $W$ ein $\K$-Untervektorraum von $X$, $Y$ ein normierter $\K$-Vektorraum und $T \: W \to Y$ eine $\K$-lineare Abbildung mit abgeschlossenem Graphen.
Ferner sei $T(W)$ nicht mager in $Y$.

\noindent \textbf{Beh.:} $T \: W \to Y$ ist eine offene -- insbes.\ surjektive -- Abbildung.

\begin{Zusatz} 
Zusätzlich zu Voraussetzungen des Hauptsatzes sei $T$ injektiv. 

Dann ist $T$ surjektiv und $T^{-1} \: Y \to X$ eine stetige Abbildung.
\end{Zusatz}
\end{HS}

\textit{Beweis.} Für jedes $r \in \R_+$ definieren wir eine abgeschlossen Teilmenge von $Y$ durch
\begin{equation} \label{FA.2.35.0}
A_r := \overline{T(U_r^W(0))}
\end{equation}
und zeigen nacheinander
\begin{gather}
\forall_{r \in \R_+} \, (A_r)^{\circ} \in \U(0,Y), \label{FA.2.35.1} \\
\forall_{r \in \R_+} \, A_r \subset T(U_{3 r}^W(0)). \label{FA.2.35.2}
\end{gather}
Hieraus folgt die Existenz von $\delta \in \R_+$ mit  $U_{\delta}^Y(0) \subset (A_{\frac{1}{3}})^{\circ} \subset T(U_1^W(0))$, also ist der Hauptsatz wegen \ref{FA.2.30} ,,(iii) $\Rightarrow$ (i)`` dann bewiesen.

Den \underline{Nachweis von (\ref{FA.2.35.1})} bereiten wir durch die folgende Aussage vor:
\begin{equation} \label{FA.2.35.1-}
\forall_{r \in \R_+} \, (A_r)^{\circ} \ne \emptyset.
\end{equation}

Zu (\ref{FA.2.35.1-}): Es gilt $W = \bigcup_{n \in \N_+} U_n^W(0)$, also
\begin{equation} \label{FA.2.35.3}
T(W) = \bigcup_{n \in \N_+} \, T(U_n^W(0)) \stackrel{(\ref{FA.2.35.0})}{\subset} \bigcup_{n \in \N_+} A_n.
\end{equation}
Wir behaupten
\begin{equation} \label{FA.2.35.4}
\bigcup_{n \in \N_+} \, A_n \mbox{ ist nicht mager in } Y.
\end{equation}

{[} Zu (\ref{FA.2.35.4}): $T(W)$ ist nicht mager in $Y$, d.h.\ es existiert keine Folge nirgends dichter Teilmengen von $Y$, deren Vereinigung gleich $T(W)$ ist.
Wäre $\bigcup_{n \in \N_+} A_n$ mager in $Y$, so existierte eine Folge nirgends dichter Teilmengen $(N_k)_{k \in \N}$ von $Y$ mit $\bigcup_{k \in \N} N_k = \bigcup_{n \in \N_+} A_n$.
Dann folgte für jedes $k \in \N$
$$ \left( \overline{T(W) \cap N_k} \right)^{\circ} = \overline{T(W)}^{\circ} \cap \underbrace{\overline{N_k}^{\circ}}_{= \emptyset} = \emptyset $$
und 
$$ \bigcup_{k \in \N} T(W) \cap N_k = T(W) \cap \left( \bigcup_{k \in \N} N_k \right) = T(W) \cap \left( \bigcup_{n \in \N_+} A_n \right) \stackrel{(\ref{FA.2.35.3})}{=} T(W), $$ 
Widerspruch! {]}

Aus (\ref{FA.2.35.4}) ergibt sich
\begin{equation} \label{FA.2.35.5}
\exists_{n_0 \in \N_+} \, (A_{n_0})^{\circ} \ne \emptyset.
\end{equation}

{[} Zu (\ref{FA.2.35.5}): Angenommen, es gilt $\forall_{n \in \N_+} \, (A_n)^{\circ} = \emptyset$.
Wegen der Abgeschlossenheit von $A_n$ in $Y$ für $n \in \N$ sind dann alle $A_n$ nirgends dicht in $Y$, also ist $\bigcup_{n \in \N} A_n$ mager in $Y$, im Widerspruch zu (\ref{FA.2.35.4}). {]} 

Seien nun $r \in \R_+$ und $n_0 \in \N_+$ wie in (\ref{FA.2.35.5}).
Wir betrachten für $t := \frac{r}{n_0} \in \R_+$ stetige Abbildungen
$$ R^W_t \: W \longrightarrow W, ~~ x \longmapsto t \, x, $$
(mit stetiger Umkehrabbildung $R^W_{\frac{1}{{t}}} \: W \to W$) und
$$ R^Y_t \: Y \longrightarrow Y, ~~ y \longmapsto t \, y, $$
mit stetiger Umkehrabbildung $R^Y_{\frac{1}{{t}}} \: Y \to Y$.
Aus der $\K$-Linearität von $T$ folgt
$$ T \circ R^W_t = R^Y_t \circ T, $$
also wegen $U^W_r(0) = R^W_t(U^W_{n_0}(0))$
$$ T(U^W_r(0)) = T(R^W_t(U^W_{n_0}(0))) = R^Y_t(T(U^W_{n_0}(0))). $$
Da $R^Y_t \: Y \to Y$ ein Homöomorphismus ist, ergibt sich hieraus
$$ A_r \stackrel{(\ref{FA.2.35.0})}{=} \overline{T(U^W_r(0))} = R^Y_t \left( \overline{T(U^W_{n_0}(0))} \right) \stackrel{(\ref{FA.2.35.0})}{=} R^Y_t(A_{n_0}), $$
also erneut wegen der Homöomorphie von $R^Y_t \: Y \to Y$
$$ (A_r)^{\circ} = R^Y_t((A_{n_0})^{\circ}) \stackrel{(\ref{FA.2.35.5})}{\ne} \emptyset, $$
d.h.\ (\ref{FA.2.35.1-}) ist gezeigt.

Zu (\ref{FA.2.35.1}): Sei $r \in \R_+$.
(\ref{FA.2.35.1-}), angewandt auf $\frac{r}{2}$ anstelle von $r$, ergibt
$$ (A_{\frac{r}{2}})^{\circ} \ne \emptyset, $$
also existieren $y_0 \in (A_{\frac{r}{2}})^{\circ}$ und $\varepsilon \in \R_+$ mit
\begin{equation} \label{FA.2.35.6}
U_{\varepsilon}^Y(y_0) \subset (A_{\frac{r}{2}})^{\circ} \subset A_{\frac{r}{2}}.
\end{equation}
Zum Nachweis von (\ref{FA.2.35.1}) genügt es zu zeigen, daß gilt
$$ U_{\varepsilon}^Y(0) \subset \overline{T(U_r^W(0))} \stackrel{(\ref{FA.2.35.0})}{=} A_r. $$

Beweis hiervon:
Sei $y \in U_{\varepsilon}^Y(0)$.
Dann gilt 
$$ y_0 + y \in U_{\varepsilon}^Y(y_0) \stackrel{(\ref{FA.2.35.6})}{\subset} A_{\frac{r}{2}} \stackrel{(\ref{FA.2.35.0})}{=} \overline{T(U_{\frac{r}{2}}^W(0))}, $$
also existieren Folgen $(x_n)_{n \in \N}, (\tilde{x}_n)_{n \in \N}$ in $U_{\frac{r}{2}}^W(0)$ mit
$$ \lim_{n \to \infty} T(x_n) = y_0 ~~ \mbox{ und } ~~ \lim_{n \to \infty} T(\tilde{x}_n) = y_0 + y. $$
Es folgt $\tilde{x}_n - x_n \in U_r^W(0)$ sowie
$$ y = (y_0 + y) - y_0 = \lim_{n \to \infty} (T(\tilde{x}_n) - T(x_n)) = \lim_{n \to \infty} (T(\underbrace{\tilde{x}_n - x_n}_{\in U_r^W(0)})) \in \overline{T(U_r^W(0))}, $$
womit (\ref{FA.2.35.1}) bewiesen ist.

Zum \underline{Nachweis von (\ref{FA.2.35.2})} zeigen wir zunächst
\begin{equation} \label{FA.2.35.2-}
\forall_{r \in \R_+} \forall_{y \in A_r} \exists_{x \in U^W_r(0)} \exists_{\tilde{y} \in U^Y_{\frac{r}{2}}(0) \cap (A_{\frac{r}{2}})^{\circ}} \, y = T(x) + \tilde{y}. 
\end{equation}

Zu (\ref{FA.2.35.2-}): Seien $r \in \R_+$ und $y \in A_r$.
Aus (\ref{FA.2.35.1}), angewandt auf $\frac{r}{2}$ anstelle von $r$, folgt $0 \in (A_{\frac{r}{2}})^{\circ}$, also
\begin{equation} \label{FA.2.35.7}
y \in \{y\} - (A_{\frac{r}{2}})^{\circ}.
\end{equation}
Da die Translation $y - \ldots \: Y \to Y$ offenbar ein Homöomorphismus ist, gilt 
$$ \{y\} - (A_{\frac{r}{2}})^{\circ} = (\{y\} - A_{\frac{r}{2}})^{\circ}, $$
also ergibt (\ref{FA.2.35.7})
$$ U := U^Y_{\frac{r}{2}}(y) \cap (\{y\} - A_{\frac{r}{2}})^{\circ} \in \U(y,Y). $$
Wegen $y \in A_r \stackrel{(\ref{FA.2.35.0})}{=} \overline{T(U^W_r(0))}$ folgt hieraus mittels \ref{FA.1.8} (ii), angewandt auf $T(U^W_r(0))$ anstelle von $Y$,
$$ U \cap T(U^W_r(0)) \ne \emptyset, $$
d.h.\ es existiert ein $x \in U^W_r(0)$ derart, daß gilt
$$ T(x) \in U = U^Y_{\frac{r}{2}}(y) \cap (\{y\} - A_{\frac{r}{2}})^{\circ}. $$
Dies bedeutet offenbar
$$ \tilde{y} := y - T(x) \in U^Y_{\frac{r}{2}}(0) \cap (A_{\frac{r}{2}})^{\circ}, $$
womit (\ref{FA.2.35.2-}) gezeigt ist.

Zu (\ref{FA.2.35.2}): Seien $r \in \R_+$ sowie $y \in A_r$.
Wir setzen $y_0 := y \in A_{\frac{r}{2^0}}$. 
Dann können für jedes $n \in \N$ gemäß (\ref{FA.2.35.2-}) durch
$$ y_n = T(x_{n+1}) + y_{n+1}, \mbox{ wobei } x_{n+1} \in U^W_{\frac{r}{2^n}}(0) \mbox{ und } y_{n+1} \in U^Y_{\frac{r}{2^{n+1}}}(0) \cap (A_{\frac{r}{2^{n+1}}})^{\circ} \mbox{ seien,} $$
rekursiv Folgen $(y_n)_{n \in \N}$ in $A_r$ sowie $(x_n)_{n \in \N_+}$ in $U^W_r(0)$ mit
\begin{gather}
\forall_{m \in \N_+} \, y_0 = T(x_1) + y_1 = T(x_1) + T(x_2) + y_2 = \ldots = \sum_{n=1}^m T(x_n) + y_m, \label{FA.2.35.8} \\
\forall_{m \in \N_+} \, \sum_{n=1}^m \|x_n\| \le r \, \sum_{n=0}^{m-1} \frac{1}{2^n} \stackrel{m \to \infty}{\longrightarrow} 2 r, \label{FA.2.35.9} \\ 
\lim_{m \to \infty} y_m = 0 \label{FA.2.35.10}
\end{gather}
definieren.
Nach Voraussetzung ist $X$ ein $\K$-Banachraum, also folgt aus \ref{FA.2.Reihe.1} und (\ref{FA.2.35.9}) die Existenz von $x \in \overline{U^X_{2 r}(0)}$ mit
$$ \sum_{n=1}^{\infty} x_n = x. $$
Weil $W$ nicht vollständig sein muß, folgt aus $(x_n)_{n \in \N_+} \subset U^W_r(0) \subset W$ a priori nicht die Gültigkeit von $x \in W$.
Die $\K$-Linearität von $T$, (\ref{FA.2.35.8}) und (\ref{FA.2.35.10}) ergeben jedoch
$$ T \left( \sum_{n=1}^{\infty} x_n \right) = y_0, $$
und wegen $\forall_{n \in \N_+} \, x_n \in W$ gilt $\forall_{m \in \N_+} \, \sum_{n=1}^m x_n \in W$, also folgt aus der Tatsache, daß $\Gamma_T$ nach Voraussetzung abgeschlossen ist, und \ref{FA.2.33} ,,(i) $\Rightarrow$ (ii)`` doch $x \in W$, d.h.\ $x \in \overline{U_{2 r}^X(0)} \cap W \subset U_{3 r}^W(0)$, sowie $T(x) = y_0 = y$.
Damit ist auch (\ref{FA.2.35.2}) gezeigt. 

Der Zusatz ist klar. \q

\begin{Bsp} \label{FA.2.35.B}
Wir betrachten $X := (\{ f \in \mathcal{C}([0,1],\R) \, | \, f(0) = 0 \}, \| \ldots \|_{\infty})$ und $Y := (\{ f \in \mathcal{C}^1([0,1],\R) \, | \, f(0) = 0 \}, \| \ldots \|_{\infty})$ sowie $T \: X \to Y$ definiert durch
$$ \forall_{f \in X} \, T \,f := \int_0^x f(t) \, \d t. $$
$T \: X \to Y$ ein $\R$-Vektorraum-Isomorphismus, dessen Graph folglich abgeschlossen ist.
Die Umkehrabbildung $\frac{\d}{\d x}|_Y \: Y \to X$ ist jedoch nicht stetig! (Das in \ref{FA.2.32.B} genannte Argument zeigt auch dies.)
Da $X$ ein $\R$-Banachraum ist, folgt aus dem Zusatz, daß $T(X) =  Y$ mager in $Y$ ist, d.h.\ $Y$ ist nicht Bairesch -- also auch kein $\R$-Banach\-raum.
\end{Bsp}

\begin{Satz} \label{FA.2.36} 
Seien $X,Y$ $\K$-Banachräume und $T \in \mathcal{L}_{\K}(X,Y)$ eine surjektive Abbildung.

Dann ist $T$ eine offene Abbildung.
\end{Satz}

\textit{Beweis.} Zum einen ist $X$ abgeschlossen in $X$, also folgt aus \ref{FA.2.34} (ii) die Abgeschlossenheit von $\Gamma_T$.
Zum anderen ist $Y$ als $\K$-Banachraum ein Baire-Raum, also ist die nach \ref{FA.1.40} generische Teilmenge $T(X) = Y$ von $Y$ wegen \ref{FA.1.42.Kor} nicht mager in $Y$.
Daher folgt aus \ref{FA.2.35} die Offenheit von $T$. \q

\begin{Kor}[Satz vom inversen Operator] \index{Satz!vom inversen Operator} \label{FA.2.37}
Seien $X,Y$ $\K$-Banachräume und $T \in \mathcal{L}_{\K}(X,Y)$ eine bijektive Abbildung.

Dann ist $T^{-1} \: Y \to X$ eine stetige Abbildung, also ein Homöomorphismus. \q
\end{Kor}

\begin{Kor} \label{FA.2.38}
Seien $X$ ein $\K$-Vektorraum und $\| \ldots \|, \| \ldots \|_*$ zwei Normen für $X$ derart, daß $(X,\| \ldots \|)$ sowie $(X,\| \ldots \|_*)$ $\K$-Banachräume sind.
Ferner existiere $D \in \R_+$ mit $\| \ldots \|_* \le D \, \| \ldots \|$.

Dann existiert $C \in \R_+$ mit $C \, \| \ldots \| \le \| \ldots \|_*$, also sind $\| \ldots \|$ und $\| \ldots \|_*$ äquivalent.
\end{Kor}

\textit{Beweis.} Nach Voraussetzung ist $\id_X \: (X,\|\ldots\|) \to (X,\|\ldots\|_*)$ eine stetige Abbildung, die die Voraussetzung des letzten Korollares erfüllt.
Daher ist auch $\id_X \: (X,\|\ldots\|_*) \to (X,\|\ldots\|)$ eine stetige Abbildung, d.h.\ es existiert $\widetilde{C} \in \R$ mit $\|\ldots\| \le \widetilde{C} \, \| \ldots \|_*$.
Offenbar gilt $\widetilde{C} > 0$, d.h.\ $C := \frac{1}{\widetilde{C}}$ leistet das Gewünschte. \q

\begin{Satz}[Satz vom abgeschlossen Graphen] \index{Satz!vom abgeschlossenen Graphen} \label{FA.2.39}
Seien $X,Y$ $\K$-Banachräume und $T \: X \to Y$ eine $\K$-lineare Abbildung mit abgeschlossenem Graphen.

Dann ist $T \: X \to Y$ stetig.
\end{Satz}

\textit{Beweis.} Wir versehen $X \times Y$ wieder mit der Maximumsnorm, womit $X \times Y$ mit $X$ und $Y$ offenbar ebenfalls ein $\K$-Banachraum ist, der $\Gamma_T$ nach Voraussetzung als abgeschlossenen $\K$-Untervektorraum besitzt.
Daher ist $\Gamma_T$ zusammen mit der Maximumsnorm ein $\K$-Banachraum.
$\pi_1 \: \Gamma_T \to X, \, (x,T(x)) \mapsto x,$ ist eine $\K$-lineare stetige bijektive Abbildung, also nach dem Satz vom inversen Operator ein Homöomorphismus.
Da auch $\pi_2 \: \Gamma_T \to Y, \, (x,T(x)) \mapsto T(x),$ stetig ist, ergibt sich die Stetigkeit von $T = \pi_2 \circ {\pi_1}^{-1}$. \q

\subsection*{Anhang: Normierbare topologische Vektorräume} \addcontentsline{toc}{subsection}{Anhang: Normierbare topologische Vektorräume}

Für die Existenz einer Norm auf einem $\K$-Vektorraum ist es offenbar notwendig, daß es sich um einen topologischen $\K$-Vektorraum im folgenden Sinne handelt.

\begin{Def}[Topologische Vektorräume] \label{FA.N.7}
Ein \emph{topologischer $\K$-Vektorraum}\index{Raum!topologischer Vektor-}\index{Vektorraum!topologischer} ist ein topologischer Raum $X$, der eine Struktur als $\K$-Vektoraum besitzt derart, daß die Addition und die skalare Multiplikation stetig sind, d.h.\
\begin{gather*}
\forall_{(x_1,x_2) \in X \times X} \forall_{U \in \U(x_1+x_2,X)} \exists_{V_1 \in \U(x_1,X)} \exists_{V_2 \in \U(x_2,X)} \, V_1 + V_2 \subset U, \\
\forall_{\lambda \in \K} \forall_{x \in X} \forall_{U \in \U(\lambda \, x, X)} \exists_{\varepsilon \in \R_+} \exists_{V \in \U(x,X)} \forall_{\mu \in \U_{\varepsilon}(0,\K)} \, \mu \, V \subset U.
\end{gather*}
\end{Def}

\begin{Bem*}
Im Gegensatz zur Definition in \cite{Kad} oder \cite{Rud} braucht ein topologischer $\K$-Vektorraum für uns nicht hausdorffsch zu sein.
Z.B.\ ist jeder nicht-nulldimensionale $\K$-Vek\-tor\-raum $X$ versehen mit der trivialen Topologie $\{ \emptyset, X \}$ ein nicht-hausdorffscher topologischer $\K$-Vektorraum.
\end{Bem*}

Wir wollen uns zum Abschluß dieses Kapitels mit der Frage beschäftigen, ob ein topologischer $\K$-Vektorraum \emph{normierbar}\index{Raum!topologischer Vektor-!normierbarer}\index{Vektorraum!topologischer!normierbarer} ist.
Dies soll bedeuten, daß eine Norm existiert so, daß die entsprechende Normtopologie mit der gegebenen Topologie übereinstimmt.

Der folgende Satz ergibt sich unmittelbar aus der letzten Definition.

\begin{Satz} \label{FA.N.8}
Sei $X$ ein topologischer $\K$-Vektorraum.
\begin{itemize}
\item[(i)] Für jedes $a \in X$ ist die \emph{Translation $T_a$ mit $a$}
$$ \boxed{T_a \: X \longrightarrow X, ~~ x \longmapsto a + x,} $$
ein Homöomorphismus mit Umkehrabbildung $T_{-a}$.
\item[(ii)] Für jedes $\lambda \in \K \setminus \{0\}$ ist die \emph{skalare Multiplikation mit $\lambda$}
$$ \boxed{S_{\lambda} \: X \longrightarrow X, ~~ x \longmapsto \lambda \, x}, $$
ein Homöomorphismus mit Umkehrabbildung $S_{\lambda^{-1}}$. \q
\end{itemize}
\end{Satz}

\begin{Bem*}
Sei $X$ ein topologischer $\K$-Vektorraum.
\begin{itemize}
\item[1.)] Wegen (i) gilt
$$ \forall_{a \in X} \, \U(a,X) = a + \U(0,X) := \{ a + U \, | \, U \in \U(0,X) \}, $$
d.h.\ die Topologie von $X$ ist durch $\U(0,X)$ vollständig bestimmt.
\item[2.)] Für jede offene Teilmenge $U$ von $X$ und jedes $\lambda \in \K \setminus \{0\}$ ist $\lambda \, U$ nach (ii) offen in $X$.
\end{itemize}
\end{Bem*}

\begin{Satz} \label{FA.N.H}
Sei $X$ ein topologischer $\K$-Vektorraum.
Dann sind die folgenden Aussagen paarweise äquivalent:
\begin{itemize}
\item[(i)] $X$ ist hausdorffsch.
\item[(ii)] $\forall_{x \in X \setminus \{0\}} \exists_{U \in \U(0,X)} \, x \notin U$.
\item[(iii)] $\{0\}$ ist abgeschlossen in $X$.
\item[(iv)] $\D \bigcap_{U \in \U(0,X)} U = \{0\}$.
\end{itemize}
\end{Satz}

\textit{Beweis.} ,,(i) $\Rightarrow$ (ii)`` Zu $x \in X \setminus \{0\}$ existieren $U \in \U(0,X)$ und $V \in \U(x,X)$ mit $U \cap V = \emptyset$, also $x \notin U$.

,,(ii) $\Rightarrow$ (i)`` Seien $x,y \in X$ mit $x \ne y$, also $x-y \ne 0$.
Dann gibt es aufgrund von (ii) eine Umgebung $U \in \U(0,X)$ derart, daß $x-y \notin U$.
Wegen der Stetigkeit von $\ldots - \ldots \: X \times X \to X$ existiert offenbar $V \in \U(0,X)$ mit $V - V \subset U$, und wir zeigen $(x+V) \cap (y+V) = \emptyset$, womit (i) bewiesen ist:
Existierten nämlich $v_1,v_2 \in V$ mit $x + v_1 = y + v_2$, so folgte $x-y = v_1 - v_2 \in V - V \subset U$, Widerspruch!

,,(i) $\Rightarrow$ (iii)`` ist nach \ref{FA.1.Punkt} klar.

,,(iii) $\Rightarrow$ (iv)`` folgt sofort aus $\bigcap_{U \in \U(0,X)} U = \overline{\{0\}}$.

,,(iv) $\Rightarrow$ (i)`` Wäre $X$ nicht hausdorffsch, so existierte wegen der bereits bewiesenen Aussage ,,(ii) $\Rightarrow$ (i)`` ein $x \in X \setminus \{0\}$ mit $x \in \bigcap_{U \in \U(0,X)} U$, im Widerspruch zu (iv). \q

\begin{Bem*}
Aus ,,(ii) $\Rightarrow$ (i)`` ergibt sich, daß ein topologischer $\K$-Vek\-tor\-raum mit dem ersten Trennungsaxiom\index{Trennungs!-axiome} auch das zweite erfüllt.
\end{Bem*}

\begin{Satz} \label{FA.N.9}
Seien $X$ ein topologischer $\K$-Vektorraum und $U \in \U(0,X)$.
Dann gilt:
\begin{itemize}
\item[(i)] $U$ enthält eine ausgewogene Umgebung von $0$.
\item[(ii)] Ist $U$ konvex, so enthält $U$ eine ausgewogene konvexe Umgebung von $0$.
\end{itemize}
\end{Satz}

Wir bereiten den Beweis des Satzes durch das folgende Lemma vor.

\begin{Lemma} \label{FA.N.9.L}
Seien $X$ ein topologischer $\K$-Vektorraum und $C$ eine konvexe Teilmenge von $X$.

Dann ist auch $C^{\circ}$ eine konvexe Teilmenge von $X$.
\end{Lemma}

\textit{Beweis.} Seien $x,y \in C^{\circ}$.
Dann existiert offenbar $U \in \U(0,X)$ mit $x+U \subset C$ sowie $y+U \subset C$ und für jedes $t \in [0,1]$ sowie $u \in U$ gilt
$$ (1+t) \, x + t \, y + u = (1+t) \, x + t \, y + (1+t) \, u + t \, u = (1+t) \, (\underbrace{x+u}_{\in C}) + t \, (\underbrace{y+u}_{\in C}) \in C, $$
d.h.\ $((1+t) \, x + t \, y) + U \subset C$, also gilt $(1+t) \, x + t \, y \in C^{\circ}$. \q
\A
\textit{Beweis des Satzes.} Zu (i): Wegen der Stetigkeit von $S_0$ existieren $\varepsilon \in \R_+$ und $V \in \U(0,X)$ derart, daß für alle $\mu \in U_{\varepsilon}(0,\K)$ gilt $\mu \, V \subset U$.
Offenbar ist $W := \bigcup_{\mu \in U_{\varepsilon}(0,\K)} \mu \, V \in \U(0,X)$ eine ausgewogene Teilmenge von $X$ mit $W \subset U$.

Zu (ii): Sei also $U$ konvex.
Wir setzen
\begin{equation} \label{FA.N.9.0}
C := \bigcap_{\sigma \in \K, \, |\sigma|=1} \sigma \, U, 
\end{equation}
also ist $C$ nach den Beispielen 3.) und 4.) zu \ref{FA.N.2} (iv) konvex.

$W \in \U(0,X)$ sei eine nach (i) existierende ausgewogene Teilmenge von $U$.
Dann gilt für jedes $\sigma \in \K$ mit $|\sigma|=1$
$$ \frac{1}{\sigma} \, W \subset W \subset U, $$
also $W \subset \sigma \, U$.
Aus (\ref{FA.N.9.0}) folgt $W \subset C$ und wegen $W \in \U(0,X)$ somit auch $C^{\circ} \in \U(0,X)$.
Es genügt nun zu zeigen, daß $C^{\circ}$ konvex und ausgewogen ist, wobei wir die Konvexität von $C^{\circ}$ in Lemma \ref{FA.N.9.L} bewiesen haben.

Zur Ausgewogenheit von $C^{\circ}$: Wir behaupten
\begin{equation} \label{FA.N.9.1}
\forall_{\lambda \in \K, \, |\lambda| \le 1} \, \lambda \, C = \bigcap_{\sigma \in \K, \, |\sigma|=1} |\lambda| \, \sigma \, U.
\end{equation}

{[} Zu (\ref{FA.N.9.1}): Sei $\lambda \in \K$ mit $|\lambda| \le 1$.
Wegen (\ref{FA.N.9.0}) können wir ohne Einschränkung $\lambda \ne 0$ annehmen.
Dann gilt 
$$ \lambda \, C = |\lambda| \, \frac{\lambda}{|\lambda|} \, C \stackrel{(\ref{FA.N.9.0})}{=} \bigcap_{\sigma \in \K, \, |\sigma|=1} |\lambda| \, \frac{\lambda}{|\lambda|} \, \sigma \, U = \bigcap_{\tilde{\sigma} \in \K, \, |\tilde{\sigma}|=1} |\lambda| \, \tilde{\sigma} \, U, $$
womit (\ref{FA.N.9.1}) gezeigt ist. {]}
\pagebreak

Für jedes $\sigma \in \K$ (mit $|\sigma| = 1$) ist $\sigma \, U$ nach \ref{FA.N.2} (iv) Beispiel 3.) eine konvexe Teilmenge von $X$, die zusätzlich $0$ enthält, also gilt 
$$ \forall_{\lambda \in \K, \, |\lambda| \le 1} \, |\lambda| \, \sigma \, U  \subset \sigma \, U. $$ 
Hieraus und aus (\ref{FA.N.9.1}), (\ref{FA.N.9.0}) folgt $\lambda \, C \subset C$ für $\lambda \in \K$ mit $|\lambda| \le 1$. \q

\begin{Def}[Beschränkte Mengen] \label{FA.N.10}
Seien $X$ ein topologischer $\K$-Vektor\-raum und $B$ eine Teilmenge von $X$.

$B$ heißt \emph{beschränkt}\index{Menge!beschränkte} genau dann, wenn zu jeder Umgebung $U \in \U(0,X)$ eine Zahl $C \in \R_+$ mit $\forall_{t \in \R_+, \, t > C} \, B \subset t \, U$ existiert.

\begin{Bem*}
Im Falle eines normierten $\K$-Vektorraumes $X$ stimmt diese Definition einer beschränkten Teilmenge von $X$ offenbar mit der in \ref{FA.2.E.5} (i) gegebenen überein.
\end{Bem*}
\end{Def}

\begin{Satz} \label{FA.N.11}
Seien $X$ ein topologischer $\K$-Vektorraum und $U \in \U(0,X)$.
\begin{itemize}
\item[(i)] Ist $(t_n)_{n \in \N}$ eine streng monoton wachsende Folge positiver reeller Zahlen mit $\lim_{n \to \infty} t_n = \infty$, so gilt $X = \bigcup_{n \in \N} t_n \, U$.
Insbesondere ist $U$ absorbierend.
\item[(ii)] Ist $U$ beschränkt und $(t_n)_{n \in \N}$ eine streng monoton fallende Folge in $\R$ mit $\lim_{n \to \infty} t_n = 0$, so existiert zu jedem $\widetilde{U} \in \U(0,X)$ ein $n_0 \in \N$ mit $\forall_{n \in \N, \, n \ge n_0} \, t_n \, U \subset \widetilde{U}$.
\end{itemize}
\end{Satz}

\textit{Beweis.} Zu (i): Sei $x \in X$.
Da $X$ ein topologischer $\K$-Vektorraum ist, ist die Abbildung $M_x \: \K \to X, \, \lambda \mapsto \lambda \, x,$ stetig.
Daher gilt $\overline{M_x}^1(U) \in \U(0,\K)$ und es existiert $n_0 \in \N$ mit $\forall_{n \in \N, \, n \ge n_0} \, \frac{1}{t_n} \in \overline{M_x}^1(U)$, d.h.\ $x \in t_n \, U$.

Zu (ii): Da $U$ beschränkt ist, existiert $C \in \R_+$ derart, daß $\forall_{t \in \R, \, t > C} \, U \subset t \, \widetilde{U}$ gilt.
Insbesondere existiert $n_0 \in \N$ mit $\forall_{n \in \N, \, n \ge n_0} \, U \subset \frac{1}{t_n} \, \widetilde{U}$, d.h.\ $t_n \, U \subset \widetilde{U}$. \q 

\begin{Def}[Lokal-beschränkte sowie lokal-konvexe Vektorräume] \label{FA.N.12} 
Sei $X$ ein topologischer $\K$-Vek\-tor\-raum.
Wir definieren dann:
\begin{itemize}
\item[(i)] $X$ \emph{lokal-beschränkt}\index{Raum!lokal-beschränkter Vektor-}\index{Vektorraum!lokal-beschränkter} $: \Longleftrightarrow$ $\forall_{x \in X} \forall_{U \in \U(x,X)} \exists_{V \in \U(x,X)} \, V \subset U \, \wedge \, V \mbox{ beschränkt}$.

\begin{Bem*}
Die Definition der lokalen Beschränktheit erfolgt in der Literatur durch \ref{FA.N.13} (i).
Der Autor dieser Seiten hat obige Form gewählt, um den typischen Aspekt einer Lokalitätsbedingung hervorzuheben.
\end{Bem*}
\item[(ii)] $X$ \emph{lokal-konvex}\index{Raum!lokal-konvexer Vektor-}\index{Vektorraum!lokal-konvexer} $: \Longleftrightarrow$ $\forall_{x \in X} \forall_{U \in \U(x,X)} \exists_{V \in \U(x,X)} \, V \subset U \, \wedge \, V \mbox{ konvex}$.
\end{itemize}
\end{Def}

\begin{Satz} \label{FA.N.13}
Sei $X$ ein topologischer $\K$-Vek\-tor\-raum.
Dann gilt:
\begin{itemize}
\item[(i)] $X$ ist genau dann lokal-beschränkt, wenn $0$ eine beschränkte Umgebung besitzt.
\item[(ii)] $X$ ist genau dann lokal-konvex, wenn gilt
$$ \forall_{U \in \U(0,X)} \exists_{V \in \U(0,X)} \, V \subset U \, \wedge \, V \mbox{ konvex}. $$
\item[(iii)] $X$ ist genau dann lokal-beschränkt und lokal-konvex, wenn eine beschränkte konvexe Umgebung von $0$ existiert.
\end{itemize}
\end{Satz}

Wir bereiten den Beweis der Satzes durch ein Lemma vor.

\begin{Lemma} \label{FA.N.14}
Seien $X$ ein topologischer $\K$-Vek\-tor\-raum und $B_1,B_2$ beschränkte Teilmengen von $X$.

Dann ist $B_1+B_2$ beschränkt.
\end{Lemma}

\textit{Beweis.}
Sei $U \in \U(0,X)$ beliebig.
Aus der Stetigkeit der Addition folgt offenbar die Existenz von $V \in \U(0,X)$ mit $V + V \subset U$.
Da $B_1,B_2$ beschränkt sind, existieren des weiteren $C_1,C_2 \in \R_+$ mit $\forall_{i \in \{1,2\}} \forall_{t \in \R_+, \, t > C_i} \, B_i \subset t \, V$, also folgt für jedes $t \in \R$ mit $t > \max \{C_1,C_2\}$: $B_1 + B_2 \subset t \, V + t \, V = t \, (V+V) \subset t \, U$. \q
\A
\textit{Beweisskizze des Satzes.} (i) ,,$\Rightarrow$`` sowie (ii) ,,$\Rightarrow$`` sind trivial.

(ii) ,,$\Leftarrow$`` folgt aus \ref{FA.N.8} (i), beachte, daß mit $V \in U(0,X)$ auch $x+V$ für jedes $x \in X$ konvex ist.

(iii) ,,$\Rightarrow$`` folgt aus (i), (ii), da jede Teilmenge einer beschränkten Menge ebenfalls beschränkt ist.

(iii) ,,$\Leftarrow$`` folgt aus (i), (ii) und \ref{FA.N.11} (ii), weil die Multiplikation einer konvexen Menge mit einem Skalar ebenfalls konvex ist.

Zu (i) ,,$\Leftarrow$``: Wir begründen zunächst
\begin{equation} \label{FA.N.13.1}
\forall_{U \in \U(0,X)} \exists_{V \in \U(0,X)} \, V \subset U \, \wedge \, V \mbox{ beschränkt}.
\end{equation}

{[} Zu (\ref{FA.N.13.1}): Seien $\widetilde{U} \in \U(0,X)$ die nach Voraussetzung existierende beschränkte Umgebung von $0$ und $U \in \U(0,X)$.
Dann ist $V := U \cap \widetilde{U} \in U(0,X)$ eine beschränkte Teilmenge von $U$. {]}

Seien $x \in X$ und $x + U \in x + \U(0,X) \stackrel{\ref{FA.N.8} (i)}{=} \U(x,X)$, wobei $U \in \U(0,X)$.
Wegen (\ref{FA.N.13.1}) existiert eine beschränkte Teilmenge $V \in \U(0,X)$ von $U$.
Es genügt offenbar zu zeigen, daß $x + V$ beschränkt ist, und dies ist nach Lemma \ref{FA.N.14} klar, da $\{x\}$ wegen \ref{FA.N.11} (i) beschränkt ist. \q

\begin{HS}[Charakterisierung normierbarer Vektorräume] \index{Raum!topologischer Vektor-!normierbarer}\index{Vektorraum!topologischer!normierbarer} \label{FA.N.15} $\,$

\noindent \textbf{Vor.:} Sei $X$ ein $\K$-Vektorraum.

\noindent \textbf{Beh.:} $X$ besitzt genau dann eine Norm, wenn eine hausdorffsche Topologie für $X$ existiert, die $X$ zu einem lokal-beschränkten lokal-konvexen $\K$-Vek\-tor\-raum macht.
Genauer gilt:
\begin{itemize}
\item[(i)] Bzgl.\ jeder Norm auf $X$ ist die Normtopologie hausdorffsch.
Sie macht $X$ zu einem lokal-beschränkten lokal-konvexen $\K$-Vek\-tor\-raum.
\item[(ii)] Ist $X$ ein hausdorffscher lokal-beschränkter lokal-konvexer $\K$-Vek\-tor\-raum  mit Topologie ${\rm Top}(X)$, so existiert eine Norm $\| \ldots \|$ auf $X$ derart, daß die entsprechende Normtopologie ${\rm Top}(X, \| \ldots \|)$ mit ${\rm Top}(X)$ übereinstimmt.
\end{itemize}
\end{HS}

\textit{Beweis.} Zu (i): Besitzt $X$ eine Norm, so ist die Normtopologie eine hausdorffsche Topologie für $X$.
Wegen \ref{FA.N.13} (iii) ist daher zu zeigen, daß eine beschränkte konvexe Umgebung von $0$ existiert. 
Dies ist aber klar; beachte die Bemerkung in \ref{FA.N.10} und \ref{FA.N.3} (i).

Zu (ii): Sei $X$ ein hausdorffscher lokal-beschränkter lokal-konvexer $\K$-Vek\-tor\-raum mit Topologie ${\rm Top}(X)$.
Nach \ref{FA.N.13} (iii) existiert eine konvexe beschränkte Umgebung von $0$, die wiederum nach \ref{FA.N.9} (ii) eine ausgewogene konvexe Teilmenge $C \in \U(0,X)$ enthält, welche natürlich ebenfalls beschränkt und nach \ref{FA.N.11} (i) absorbierend ist.
Wenn wir begründen, daß (\ref{FA.2.23.2}) - (\ref{FA.2.23.4}) gelten, folgt aus \ref{FA.N.E}, daß $\| \ldots \| := p_C$ eine Norm auf $X$ mit
\begin{equation} \label{FA.N.15.1}
U_1(0) = C
\end{equation}
ist, wobei $p_C$ das Minkowsi-Funtional von $C$ bezeichne.

{[} Zu (\ref{FA.2.23.2}):
Sei $x \in C$.
Dann gilt wegen der Offenheit von $C$ auch $C \in \U(x,X)$.
Da $X$ ein topologischer $\K$-Vektorraum ist, ist $M_x \: \K \to X, \, \lambda \mapsto \lambda \, x,$ stetig, also folgt die Existenz eines $n_0 \in \N$ derart, daß $\forall_{n \in \N, \, n \ge n_0} \, \frac{n+1}{n} \, x \in C$ gilt.

Zu (\ref{FA.2.23.3}): Dies folgt sofort, da $C$ ausgewogen ist.

Zu (\ref{FA.2.23.4}):
Da $C$ beschränkt ist, existiert zu jedem $U \in \U(0,X)$ eine Zahl $\tau \in \R_+$ mit $C \subset \frac{1}{\tau} \, U$, d.h.\ $\tau \, C \subset U$ und somit $\bigcap_{\tau \in \R_+} \tau \, C \subset U$.
Weil $X$ hausdorffsch ist, ergibt \ref{FA.N.H} ,,(i) $\Rightarrow$ (iv)`` daher $\bigcap_{\tau \in \R_+} \tau \, C \subset \bigcap_{U \in \U(0,X)} U = \emptyset$. {]}

Es bleibt zu zeigen, daß gilt 
$$ {\rm Top}(X) = {\rm Top}(X, \| \ldots\|). $$ 

Beweis hiervon: Für alle $x \in X$ und $\varepsilon \in \R_+$ gilt
$$ U_{\varepsilon}(x) \stackrel{\ref{FA.2.2} (ii)}{=} x + \varepsilon \, U_1(0) \stackrel{(\ref{FA.N.15.1})}{=} x + \varepsilon \, C = (T_x \circ S_{\varepsilon})(C). $$
Wegen $C \in {\rm Top}(X)$ und \ref{FA.N.8} folgt hieraus $U_{\varepsilon}(x) \in {\rm Top}(X)$, woraus sich offenbar ,,$\supset$`` ergibt.
Zum Nachweis von ,,$\subset$`` genügt es wegen \ref{FA.N.8} Bemerkung 1.) zu zeigen, daß gilt
$$ \forall_{U \in \U(0,X)} \exists_{V \in \U(0,(X,\| \ldots \|))} \, V \subset U. $$
Dies ist aber klar, da wegen der Beschränktheit von $C \in \U(0,X)$ und Satz \ref{FA.N.11} (ii) zu jedem $U \in \U(0,X)$ eine Zahl $\varepsilon \in \R_+$ mit $\varepsilon \, C \subset U$ existiert und 
$$ \varepsilon \, C = \varepsilon \, U_1(0) = U_{\varepsilon}(0) \in \U(0,(X,\| \ldots \|)) $$ 
gilt.  
\q

\subsection*{Übungsaufgaben}

\begin{UA}
Sei $A$ eine assoziative $\K$-Algebra mit Einselement.

Zeige, daß ein $\K$-Ideal von $A$ genau dann echt ist, wenn es das Einselement nicht enthält.
\end{UA}

\begin{UA}
$\ell_{\K}^1$ sei der $\K$-Vektorraum aller $\K$-wertigen Folgen derart, daß die Summe der Folgenglieder absolut konvergent ist, zusammen mit der Norm, die durch
$$ \forall_{(x_n)_{n \in \N} \in \ell_{\K}^1} \, \left\| (x_n)_{n \in \N} \right\|_1 := \sum_{n = 0}^{\infty} |x_n|. $$
gegeben ist.

Zeige, daß $\ell_{\K}^1$ ein $\K$-Banachraum ist.
\end{UA}

\begin{UA} \label{FA.absbedkonv}
Es seien $X$ ein $\K$-Banachraum und $\sum_{k=0}^{\infty} x_k$ eine Reihe in $X$, die gegen $x \in X$ konvergiert.

Beweise die folgenden Aussagen:
\begin{itemize}
\item[(i)] Ist $\sum_{k=0}^{\infty} x_k$ sogar absolut konvergent und $i \: \N \to \N$ eine bijektive Abbildung, so konvergiert $\sum_{k=0}^{\infty} x_{i(k)}$ absolut, und es gilt $\sum_{k=0}^{\infty} x_k = \sum_{k=0}^{\infty} x_{i(k)}$.
\item[(ii)] Gelte zusätzlich $n := \dim_{\K} X \in \N$.
Ist $\sum_{k=0}^{\infty} x_k$ nicht absolut konvergent, so existiert eine bijektive Abbildung $i \: \N \to \N$ derart, daß $\sum_{k=0}^{\infty} x_{i(k)}$ gegen ein Element von $X \setminus \{x\}$ oder gar nicht konvergiert.
\end{itemize}
\end{UA}

Tip: Der Beweis von (i) verläuft völlig analog zum Beweis im Falle $X = \R$, vgl.\ hierzu z.B.\ \cite[Hauptsatz 4.37 (i)]{ElAna}.
Zum Nachweis von (ii) überlege man sich zunächst, daß es genügt $X = \R^n$ mit der euklidischen Norm $\|\ldots\|_2$ zu betrachten, um sodann den \emph{Riemannschen Umordnungssatz}, den man ggf.\ z.B.\ in \cite[Hauptsatz 4.37 (ii)]{ElAna} findet, auf eine nicht absolut konvergente Komponente der $\R^n$-wertigen Reihe $\sum_{k=0}^{\infty} x_k$ anwenden zu können.

\begin{Bem*}
Seien $X$ ein normierter $\K$-Vektorraum und $\sum_{k=0}^{\infty} x_k$ eine Reihe in $X$.
Dann heißt $\sum_{k=0}^{\infty} x_k$ genau dann \emph{unbedingt konvergent}\index{Reihe!unbedingt konvergente}, wenn ein $x \in X$ existiert, so daß $\sum_{k=0}^{\infty} x_{i(k)}$ für jede bijektive Abbildung $i \: \N \to \N$ gegen $x$ konvergiert.
\textsc{Dvoretzky} und \textsc{Rogers}' \cite[Theorem 1]{DR} besagt, daß absolute und unbedingte Konvergenz von Reihen im Falle eines $\K$-Banachraumes $X$ genau dann äquivalent sind, wenn $X$ endlich-dimensional ist, d.h.\ also nach der letzten Aufgabe, daß es in unendlich-dimensionalen $\K$-Banachräumen stets eine Reihe gibt, die unbedingt aber nicht absolut konvergent ist.
\end{Bem*}

\begin{UA} \label{FA.absbedkonvBsp}
Wir werden unten in Kapitel \ref{FAna5} sehen, daß der $\K$-Vektorraum $\ell_{\K}^2 := \{ (x_k)_{k \in \N} \in \K^{\N} \, | \, \sum_{k=0}^{\infty} |x_k|^2 < \infty \}$ zusammen mit der Norm $\| \ldots \|_2$, die durch
$$ \forall_{(x_k)_{k \in \N} \in \ell_{\K}^2} \, \| (x_k)_{k \in \N} \|_2 := \sqrt{\sum_{k=0}^{\infty} |x_k|^2} $$
gegeben ist, ein $\K$-Banachraum ist.
Für jedes $k \in \N$ sei $e_k \in \ell^2_{\K}$ die Folge mit $\forall_{j \in \N} \, e_{k,j} = \delta_{kj}$.

Zeige, daß die Reihe $\sum_{k=0}^{\infty} \frac{1}{k+1} \, e_k$ in $\ell^2_{\K}$ unbedingt aber nicht absolut konvergiert.
\end{UA}

\begin{UA}
Beweise, daß die in der Fußnote auf Seite \pageref{FA.2.6} definierte Quotiententopologie tatsächlich eine Topologie ist. 
\end{UA}

\begin{UA}
Führe den Beweis der Aussage (ii) des Hauptsatzes \ref{FA.2.13} in allen Einzelheiten aus.
\end{UA}

\begin{UA}
Beweise Satz \ref{FA.2.E.2}.
\end{UA}

\begin{UA}[Satz von \textsc{Bolzano-Weierstraß}] $\,$
\begin{itemize}
\item[(i)] Beweise Lemma \ref{FA.2.E.6}, das beim Beweis des Satzes \ref{FA.2.E.4} benötigt wird.
\item[(ii)] Sei $X$ ein endlich-di\-men\-sionaler normierter $\K$-Vektorraum.
\\
Zeige, daß jede beschränkte Folge in $X$ eine konvergente Teilfolge besitzt.
\end{itemize}
\end{UA}

\begin{UA}\label{FA.UA.C(A)}
Es seien $X$ ein $\K$-Vektorraum und $A$ eine Teilmenge von $X$.
$C(A)$ bezeichne die konvexe Hülle von $A$, vgl.\ Beispiel 4.) zu \ref{FA.N.2} (iv).

\begin{itemize}
\item[(i)] $\D C(A) = \left\{ \left. \sum_{i=1}^{n+1} \lambda_i \, x_i \, \right| \, n \in \N \, \wedge \, \forall_{i \in \{1, \ldots, n+1\}} \, \left( \lambda_i \in {[} 0, 1 {]} \, \wedge \, x_i \in A \right) \, \wedge \sum_{i=1}^{n+1} \lambda_i = 1 \right\}$.
\item[(ii)] Ist $n := \dim_{\R} {\rm Span}_{\R} A \in \N$, so gilt
$$ C(A) = \left\{ \left. \sum_{i=1}^{n+1} \lambda_i \, x_i \, \right| \, \forall_{i \in \{1, \ldots, n+1\}} \, \left( \lambda_i \in {[} 0, 1 {]} \, \wedge \, x_i \in A \right) \, \wedge \sum_{i=1}^{n+1} \lambda_i = 1 \right\}. $$
\item[(iii)] Ist $X$ zusätzlich normiert und $A$ kompakt mit $\dim_{\R} {\rm Span}_{\R} A < \infty$, so ist $C(A)$ kompakt.
\end{itemize}
\end{UA}

\begin{Def*}
Es seien $X$ ein fastmetrischer Raum, $A$ eine nicht-leere Teilmenge von $X$ und $\varepsilon \in \R_+$.

Dann heißt
$$ \boxed{U_{\varepsilon}(A)} := \bigcup_{a \in A} U_{\varepsilon}(a) $$
die \emph{$\varepsilon$-Umgebung von $A$}.\index{Umgebung!$\varepsilon$-}
\end{Def*}

\begin{UA}
Sei $X$ ein normierter $\K$-Vektorraum, $A$ eine nicht-leere Teilmenge von $X$ und $\varepsilon \in \R_+$.

\begin{itemize}
\item[(i)] $A$ konvex $\Longrightarrow$ $U_{\varepsilon}(A)$ konvex.
\item[(ii)] $C(U_{\varepsilon}(A)) = U_{\varepsilon}(C(A))$.
\item[(iii)] $A$ präkompakt $\Longrightarrow$ $C(A)$ präkompakt.
\item[(iv)](Satz von \textsc{Mazur}).\index{Satz!von \textsc{Mazur}}\\
Ist $X$ ein $\K$-Banachrraum und $A$ eine relativ kompakte Teilmenge von $X$, so ist $C(A)$ ebenfalls eine relativ kompakte Teilmenge von $X$.
\end{itemize}
\end{UA}

Tip zu (iii): Aufgabe \ref{FA.UA.C(A)} (iii).

\begin{UA}
Seien $X$ ein $\K$-Vektorraum und $C$ eine absorbierende konvexe Teilmenge von $X$.
$p_C$ bezeichne das Minkowski-Funktional von $C$.

Zeige $C \subset \overline{p_C}^1({[}0,1{]})$.
\end{UA}

\begin{UA} $\,$
\begin{itemize}
\item[(i)] Beweise die Aussage der Fußnote auf Seite \pageref{FA.2.18}. (Dort ist $X$ ein $\C$-Vek\-tor\-raum!)
\item[(ii)] Sei $X$ ein normierter $\C$-Vektorraum.
\\
Zeige, daß eine $\R$-Vektorraum-Isometrie $(X')_{\R} \to (X_{\R})'$ existiert.
\end{itemize}
\end{UA}

\begin{UA}
Führe den Beweis des Trennungssatzes von \textsc{Hahn-Banach} \ref{FA.2.24} (ii) im Falle $\K = \C$ in allen Einzelheiten aus.
\end{UA}

\begin{UA} \label{FA.UA.E}
Es bezeichne $\mathcal{C}^1([0,1],\R)$ den $\R$-Vektorraum aller stetig differenzierbaren Funktionen ${[}0,1{]} \to \R$.
\begin{itemize}
\item[(i)] Der Integraloperator
$$ T \: (\mathcal{C}([0,1],\R), \| \ldots \|_{\infty}) \longrightarrow (\mathcal{C}^1([0,1],\R), \| \ldots \|_{\infty}), $$
definiert durch
$$ \forall_{f \in \mathcal{C}([0,1],\R)} \, T \, f := \int_0^x f(t) \, \d t, $$
ist $\R$-linear, injektiv und stetig.
Das Bild von $\mathcal{C}([0,1],\R)$ unter $T$ ist ein abgeschlossener Teilraum von $(\mathcal{C}^1([0,1],\R), \| \ldots \|_{\infty})$.
\item[(ii)] Der Differentialoperator
$$ \frac{\d}{\d x} \: (\mathcal{C}^1([0,1],\R), \| \ldots \|_{\infty}) \longrightarrow (\mathcal{C}([0,1],\R), \| \ldots \|_{\infty}) $$
ist $\R$-linear, hat einen abgeschossenen Graphen und ist nicht stetig.
Ferner ist das Funktional
$$ \left. \frac{\d}{\d x} \right|_{x=1} \: (\mathcal{C}^1([0,1],\R), \| \ldots \|_{\infty}) \longrightarrow \R $$
ebenfalls unstetig.
\item[(iii)] $T' \: (\mathcal{C}^1([0,1],\R), \| \ldots \|_{\infty})' \to (\mathcal{C}([0,1],\R), \| \ldots \|_{\infty})'$ ist nicht surjektiv.
\item[(iv)] $(\mathcal{C}^1([0,1], \R), \|\ldots\|_{\infty})$ ist kein $\R$-Banachraum.
\end{itemize}
\end{UA}

\begin{UA}
Es seien $X, Z$ normierte $\K$-Vektorräume, $Y$ ein $\K$-Banach\-raum, $T \in \L_{\K}(X,Y)$, $S \in \L_{\K}(Y,Z)$ und $S(Y)$ vollständig sowie
$$ X \stackrel{T}{\longrightarrow} Y \stackrel{S}{\longrightarrow} Z $$
eine exakte Sequenz.

Zeige, daß dann auch
$$ Z' \stackrel{S'}{\longrightarrow} Y' \stackrel{T'}{\longrightarrow} X' $$
eine exakte Sequenz ist.
\end{UA}

\begin{UA}
Seien $X,Y$ normierte $\K$-Vektorräume und $T \in \L_{\K}(X,Y)$.
\begin{itemize}
\item[(i)] $X, T(X)$ vollständig und $T$ injektiv $\Longrightarrow$ $T'$ surjektiv.
\item[(ii)] $Y$ vollständig und $T$ surjektiv $\Longrightarrow$ $T'$ injektiv.
\item[(iii)] Ist $X$ vollständig, so ist $T'(X')$ abgeschlossen in $Y'$, und es existiert ein kanonischer stetiger $\K$-Vektorraum-Isomorphismus
$$ \underbrace{X'/T'(Y')}_{= \, {\rm Cokern} \, (T')} \longrightarrow ({\rm Kern} \, T)' $$
mit stetiger Umkehrabbildung.
\item[(iv)] Ist $Y$ vollständig und $T(X)$ abgeschlossen in $Y$, so existiert ein kanonischer stetiger $\K$-Vektorraum-Isomorphismus
$$ \underbrace{(Y/T(X))'}_{= \, ({\rm Cokern} \, T)'} \longrightarrow {\rm Kern} \, (T') $$
mit stetiger Unkehrabbildung.
\end{itemize}
\end{UA}

\begin{UA}
Es seien $X,Y,Z$ $\K$-Banachräume und $T \: X \to Y$ eine $\K$-li\-ne\-a\-re Abbildung sowie $S \: Y \to Z$ eine injektive $\K$-lineare Abbildung.
Ferner seien $S \: Y \to Z$ und $S \circ T \: X \to Z$ stetig.

Zeige, daß $T \: X \to Y$ stetig ist.
\end{UA}

\cleardoublepage
\section{Räume stetiger Abbildungen} \label{FAna3}

\begin{Def} \label{FA.3.1}
Es seien $M$ eine nicht-leere Menge und $Y$ ein normierter $\K$-Vektorraum (bzw.\ eine normierte $\K$-Algebra).
\begin{itemize}
\item[(i)] Für jedes $f \in Y^M$ setzen wir
$$ \boxed{\|f\|_{\infty}} := \sup \{ \|f(p)\| \, | \, p \in M \} \in [0, \infty]. $$
\item[(ii)] Eine Abbildung $f \in Y^M$ heißt \emph{beschränkt}\index{Abbildung!beschränkte}, wenn gilt $\|f\|_{\infty} < \infty$.
Den $\K$-Untervektorraum (bzw.\ die $\K$-Unteralgebra) von $Y^M$ der beschränkten Abbildungen $M \to Y$ bezeichnen wir mit $\boxed{\mathcal{B}(M,Y)}$.
\item[(iii)] Sei $M$ sogar ein topologischer Raum.
Wir bezeichnen dann den $\K$-Unter\-vek\-tor\-raum (bzw.\ die $\K$-Unteralgebra) von $Y^M$ der beschränkten stetigen Abbildungen $M \to Y$ mit $\boxed{\mathcal{C}_b(M,Y)} = \mathcal{C}(M,Y) \cap \mathcal{B}(M,Y)$.
\end{itemize}
\end{Def}

\begin{Satz} \label{FA.3.2}
Es seien $M$ eine nicht-leere Menge und $Y$ ein normierter $\K$-Vektor\-raum.
Dann gilt:
\begin{itemize}
\item[(i)] $\| \ldots \|_{\infty}$ ist eine Norm auf $\mathcal{B}(M,Y)$, sie heißt die \emph{Supremumsnorm auf $\mathcal{B}(M,Y)$}\index{Norm!Supremums- auf $\mathcal{B}(M,Y)$} und induziert die Metrik $d_{\infty}$ auf $\mathcal{B}(M,Y) \subset Y^M$, vgl.\ Beispiel 3.) zu \ref{FA.1.1}.
\item[(ii)] $\mathcal{B}(M,Y)$ ist eine abgeschlossene Teilmenge des fastmetrischen Raumes $(Y^M,d_{\infty})$.
\item[(iii)] Ist $Y$ ein $\K$-Banachraum, so ist auch $(\mathcal{B}(M,Y),\| \ldots \|_{\infty})$ ein $\K$-Banach\-raum.
\item[(iv)] Sei $Y$ sogar eine normierte $\K$-Algebra $Y$ mit Einselement $e$.
Dann ist auch $(\mathcal{B}(M,Y),\| \ldots \|_{\infty})$ eine normierte $\K$-Al\-ge\-bra, die die konstante Abbildung vom Wert $e$ als Einselement besitzt.

Insbes.\ ist $(\mathcal{B}(M,Y),\| \ldots \|_{\infty})$ im Falle einer $\K$-Banachalgebra $Y$ eine ebensolche.
\end{itemize}
\end{Satz}

\textit{Beweis.} (i) ist klar.

Zu (ii): Seien $f \in Y^M$ und $(f_n)_{n \in \N}$ eine Folge in $\mathcal{B}(M,Y)$, die bzgl.\ $d_{\infty}$ gegen $f$ konvergiert.
Dann existiert $n_0 \in \N$ mit $\|f_{n_0} - f\|_{\infty} < 1$, also gilt
$$ \|f\|_{\infty} \le \|f - f_{n_0}\|_{\infty} + \|f_{n_0}\|_{\infty} < 1 + \|f_{n_0}\|_{\infty} < \infty, $$
d.h.\ $f \in \mathcal{B}(M,Y)$.
Damit ist gezeigt, daß jeder Berührungspunkt von $\mathcal{B}(M,Y)$ in $(Y^M,d_{\infty})$ zu $\mathcal{B}(M,Y)$ gehört, also gilt (ii).

(iii) folgt aus \ref{FA.1.Isom.1}, (ii) sowie \ref{FA.1.11} (i), und (iv) ist mit (iii) klar. \q

\begin{Satz} \label{FA.3.3}
Es seien $X$ ein nicht-leerer topologischer Raum und $Y$ ein normierter $\K$-Vektor\-raum.
Dann gilt: 
\begin{itemize}
\item[(i)] $\mathcal{C}_b(X,Y)$ ist eine abgeschlossene Teilmenge des fastmetrischen Raumes $(\mathcal{C}(X,Y),d_{\infty})$.
\item[(ii)] $\mathcal{C}_b(X,Y)$ ist abgeschlossener $\K$-Untervektorraum von $(\mathcal{B}(X,Y),\| \ldots \|_{\infty})$.
\item[(iii)] Ist $Y$ ein $\K$-Banachraum, so ist auch $(\mathcal{C}_b(X,Y),\| \ldots \|_{\infty})$ ein $\K$-Banach\-raum.
\item[(iv)] Sei $Y$ sogar eine normierte $\K$-Algebra $Y$ mit Einselement $e$.
Dann ist auch $(\mathcal{C}_b(X,Y),\| \ldots \|_{\infty})$ eine normierte $\K$-Al\-ge\-bra, die die konstante Abbildung vom Wert $e$ als Einselement besitzt.

Insbes.\ ist $(\mathcal{C}_b(X,Y),\| \ldots \|_{\infty})$ im Falle einer $\K$-Banachalgebra $Y$ eine ebensolche.
\end{itemize}
\end{Satz}

\textit{Beweis.} (i) folgt aus \ref{FA.1.Isom.2} und \ref{FA.3.2} (ii).

Zu (ii): Daß $\mathcal{C}_b(X,Y)$ ein $\K$-Untervektorraum von $\mathcal{B}(X,Y)$ ist, ist klar.
Es seien $f \in \mathcal{B}(X,Y)$ und $(f_n)_{n \in \N}$ eine Folge in $\mathcal{C}_b(X,Y)$, die bzgl.\ $\| \ldots \|_{\infty}$ gegen $f$ konvergiert.
Dann konvergiert die Folge $(f_n)_{n \in \N}$ stetiger Abbildungen auch in $(Y^X,d_{\infty})$ gegen $f$, also folgt aus \ref{FA.1.Isom.2} und \ref{FA.1.KA} die Stetigkeit von $f$. 

(iii) und (iv) folgen aus (ii) und \ref{FA.3.2} (iii), (iv), beachte \ref{FA.1.11} (i). \q

\begin{Def} \label{FA.3.4}
Seien $X$ ein nicht-leerer topologischer Raum und $Y$ ein normierter $\K$-Vektorraum (bzw.\ eine normierte $\K$-Algebra).
\begin{itemize}
\item[(i)] Eine Abbildung $f \: X \to Y$ heißt \emph{im Unendlichen verschwindend}\index{Abbildung!im Unendlichen ver\-schwin\-den\-de} genau dann, wenn zu jedem $\varepsilon \in \R_+$ eine kompakte Teilmenge $K$ von $X$ mit $\forall_{x \in X \setminus K} \, \|f(x)\| < \varepsilon$ existiert.

Den $\K$-Untervektorraum (bzw.\ die $\K$-Unteralgebra) von $Y^X$ der stetigen Abbildungen $X \to Y$, die im Unendlichen verschwinden, bezeichnen wir mit $\boxed{\mathcal{C}_0(X,Y)}$.
\item[(ii)] Für jede Abbildung $f \: X \to Y$ heißt $\boxed{{\rm Tr}(f)} := \overline{\{ x \in X \, | \, f(x) \ne 0\}}$ der \emph{Träger von $f$}\index{Träger}.

Mit $\boxed{\mathcal{C}_c(X,Y)}$ bezeichnen wir den $\K$-Untervektorraum (bzw.\ die $\K$-Un\-ter\-al\-ge\-bra) von $Y^X$ der stetigen Abbildungen $X \to Y$ mit kompaktem Träger.\footnote{Beachte, daß für alle $f,g \in \mathcal{C}_c(X,Y)$ gilt $\forall_{x \in X} \, \left( x \notin {\rm Tr}(f) \cup {\rm Tr}(g) \Rightarrow f(x) + g(x) = 0 \right)$, also auch $\{ x \in X \, | \, (f+g)(x) \ne 0 \} \subset {\rm Tr}(f) \cup {\rm Tr}(g)$ und somit ${\rm Tr}(f+g) \subset {\rm Tr}(f) \cup {\rm Tr}(g)$.}
\end{itemize}
\end{Def}

\begin{Satz} \label{FA.3.5}
Seien $X$ ein nicht-leerer topologischer Raum und $Y$ ein normierter $\K$-Vektorraum.
\begin{itemize}
\item[(i)] $\mathcal{C}_0(X,Y)$ ist abgeschlossener $\K$-Untervektorraum von $(\mathcal{C}_b(X,Y),\| \ldots \|_{\infty})$.
\item[(ii)] Ist $Y$ ein $\K$-Banachraum, so ist auch $(\mathcal{C}_0(X,Y),\| \ldots \|_{\infty})$ ein $\K$-Banach\-raum.
\item[(iii)] Sei $X$ zusätzlich als kompakt vorausgesetzt.

Ist $Y$ eine normierte $\K$-Algebra, so ist auch $(\mathcal{C}_0(X,Y),\| \ldots \|_{\infty})$ eine normierte $\K$-Al\-gebra, die die konstante Abbildung vom Wert $e$ als Einselement besitzt, wobei $e$ das Einselement von $Y$ bezeichne.

Insbes.\ ist dann $(\mathcal{C}_0(X,Y),\| \ldots \|_{\infty})$ im Falle einer $\K$-Banachalgebra $Y$ eine ebensolche.
\item[(iv)] $\mathcal{C}_c(X,Y) \subset \mathcal{C}_0(X,Y)$.
\end{itemize}
\end{Satz}

\textit{Beweis.} Zu (i): 1.) Zu jedem $f \in \mathcal{C}_0(X,Y)$ existiert eine kompakte Teilmenge $K$ von $X$ mit $\forall_{x \in X \setminus K} \, \|f(x)\| < 1$.
Da $K$ kompakt und $\|f\|$ stetig ist, existiert weiter $C \in \R$ mit $\|f|_K\| \le C$, also gilt $\|f\|_{\infty} \le \max \{1,C\}$, d.h.\ $f \in \mathcal{C}_b(X,Y)$.

2.) Seien $f \in \mathcal{C}_b(X,Y)$ und  $(f_n)_{n \in \N}$ eine Folge in $\mathcal{C}_0(X,Y)$, die bzgl.\ $\| \ldots \|_{\infty}$ gegen $f$ konvergiert.
Ferner sei $\varepsilon \in \R_+$.
Dann existiert $n_0 \in \N$ derart, daß gilt $\forall_{n \in \N, \, n \ge n_0} \, \|f_n -f\|_{\infty} < \frac{\varepsilon}{2}$.
Wegen $f_{n_0} \in \mathcal{C}_0(X,Y)$ gibt es ein Kompaktum $K \subset X$ mit $\forall_{x \in X \setminus K} \, \|f_{n_0}(x)\| < \frac{\varepsilon}{2}$, also folgt für jedes $x \in X \setminus K$
$$ \|f(x)\| \le \|f_{n_0}(x)\| + \|f_{n_0}(x) - f(x)\| \le \|f_{n_0}(x)\| + \|f_{n_0} - f\|_{\infty} < \frac{\varepsilon}{2} + \frac{\varepsilon}{2} = \varepsilon, $$
also gilt $f \in \mathcal{C}_0(X,Y)$.

Aus 1.) und 2.) folgt (i).

(ii) und (iii) folgen aus (i) und \ref{FA.3.3} (iii), (iv), beachte \ref{FA.1.11} (i).
(iv) ist trivial. \q

\begin{Bem*} 
Sind $X \ne \emptyset$ ein nicht-kompakter topologischer Raum und $Y$ eine normierte $\K$-Algebra, so besitzen die assoziativen $\K$-Algebren $\mathcal{C}_0(X,Y)$ und $\mathcal{C}_c(X,Y)$ kein Einselement.
\end{Bem*}

\begin{Bsp} \label{FA.3.B}
Wir betrachten $\N$ wie üblich als topologischen Teilraum von $\R$ und setzen
$$ \boxed{\mathbf{c_0}_{\K}} := (\mathcal{C}_0(\N,\K),\|\ldots\|_{\infty}) ~~ \mbox{ und } ~~ \boxed{\ell^{\infty}_{\K}} := (\mathcal{B}(\N,\K),\|\ldots\|_{\infty}). $$
$\mathbf{c_0}_{\K}$ bzw.\ $\ell^{\infty}_{\K}$ ist also ein $\K$-Banachraum bzw.\ eine $\K$-Banachalgebra, dessen zugrundeligender $\K$-Vektorraum aus den Null- bzw.\ beschränkten Folgen in $\K$ besteht.
Ferner sei \boxed{\ell^1_{\K}} der $\K$-Banachraum mit zugrundeliegendem $\K$-Vektorraum
$$ | \ell^1_{\K} | = \left\{ (x_n)_{n \in \N} \, | \, \sum_{n=0}^{\infty} x_n \mbox{ ist absolut konvergent } \right\} $$
und der Norm $\boxed{\| \ldots \|_1}$, die durch
$$ \forall_{(x_n)_{n \in \N} \in | \ell^1_{\K} |} \, \| (x_n)_{n \in \N} \|_1 := \sum_{n=0}^{\infty} |x_n| $$
gegeben ist.
Dann sind
$$ T \: \ell^{\infty}_{\K} \longrightarrow (\ell^1_{\K})', ~~ (x_n)_{n \in \N} \longmapsto \left( (\tilde{x}_n)_{n \in \N} \mapsto \sum_{n=0}^{\infty} x_n \, \tilde{x}_n \right) $$
und
$$ \widetilde{T} \: \ell^1_{\K} \longrightarrow (\mathbf{c_0}_{\K})', ~~ (x_n)_{n \in \N} \longmapsto \left( (\tilde{x}_n)_{n \in \N} \mapsto \sum_{n=0}^{\infty} x_n \, \tilde{x}_n \right) $$
$\K$-Vektorraum-Isometrien.
Bezeichnet $e_n \in \mathbf{c_0}_{\K}$ für jedes $n \in \N$ die Folge mit $\forall_{k \in \N} \, e_{n,k} = \delta_{nk}$, so gilt für die $\K$-Vektorraum-Isometrie $S := \widetilde{T}^{-1} \: (\mathbf{c_0}_{\K})' \to \ell^1_{\K}$
$$ \forall_{F \in (\mathbf{c_0}_{\K})'} \, S \, F = (F((e_{n,k})_{k \in \N}))_{n \in \N}, $$
und nach \ref{FA.2.TB.3.K} (ii) ist auch $S' \: (\ell^1_{\K})' \to (\mathbf{c_0}_{\K})''$ und somit $S' \circ T \: \ell^{\infty}_{\K} \to (\mathbf{c_0}_{\K})''$ eine $\K$-Vektorraum-Isometrie.
Des weiteren kommutiert das folgende Diagramm:
%
%\begin{diagram}
%\mathbf{c_0}_{\K}  & \rEmbed^{i_{\mathbf{c_0}_{\K}}} & (\mathbf{c_0}_{\K})'' \\
%\dEmbed            & \ruEmbedonto_{S' \circ T}       &                       \\
%\ell^{\infty}_{\K} &                                 &                       \\
%\end{diagram}
%
\begin{center}
\begin{tikzpicture}[node distance=2cm, every arrow/.style={thick,->}]
    % Knoten definieren
    \node (X) {$\mathbf{c_0}_{\K}$};
    \node (Y) [right of=X] {$(\mathbf{c_0}_{\K})''$};
    \node (Z) [below of=X] {$\ell^{\infty}_{\K}$};
    
    % Pfeile zeichnen
    \draw[>->] (X) -- node[above] {$i_{\mathbf{c_0}_{\K}}$} (Y);
    \draw[>->] (X) -- (Z);
    \draw[>-{>>}] (Z) -- node[right] {$S' \circ T$} (Y);
    \end{tikzpicture}
\end{center}
Folglich ist $\mathbf{c_0}_{\K}$ nicht reflexiv.
\end{Bsp}

\begin{Satz} \label{FA.3.6}
Seien $X$ ein nicht-leerer lokal-kompakter Hausdorffraum und $Y$ ein normierter $\K$-Vektorraum.

Dann liegt $\mathcal{C}_c(X,Y)$ dicht in $(\mathcal{C}_0(X,Y), \| \ldots \|_{\infty})$.
\end{Satz}

\textit{Beweis.} Seien $f \in \mathcal{C}_0(X,Y)$ und $n \in \N$.
Dann existiert eine kompakte Teilmenge $\widetilde{K}_n$ von $X$ derart, daß gilt
$$ \forall_{x \in X \setminus \widetilde{K}_n} \, \|f(x)\| < \frac{1}{n+1}. $$
Wegen der Stetigkeit von $\|f\| \: X \to \R$ ist
$$ K_n := \overline{\|f\|}^1({[} \frac{1}{n+1}, \infty {[}) $$
eine abgeschlossene Teilmenge des Kompaktums $\widetilde{K}_n$, also selbst kompakt.
Aus dem Urysohnschen Lemma \ref{FA.B.1} folgt die Existenz von $\varphi_n \in \mathcal{C}_c(X,[0,1])$ mit $\varphi_n|_{K_n} = 1$.
Dann gilt
$$ {\rm Tr}(\varphi_n f) = \overline{ \{x \in X \, | \, \varphi_n(x) \ne 0 \} \cap \{x \in X \, | \, f(x) \ne 0 \} } \subset {\rm Tr}(\varphi_n) \cap {\rm Tr}(f) \subset {\rm Tr}(\varphi_n), $$
also ist ${\rm Tr}(\varphi_n f)$ als abgeschlossene Teilmenge eines Kompaktums selbst kompakt.
Somit folgt $\varphi_n f \in \mathcal{C}_c(X,Y)$.
Wegen $|1-\varphi_n| \le 1$, $\forall_{x \in K_n} \, \varphi_n(x) = 1$ und $\forall_{x \in X \setminus K_n} \, \|f(x)\| < \frac{1}{n+1}$ ergibt sich
$$ \| \varphi_n f - f \|_{\infty} = \|1 - \varphi_n\|_{\infty} \, \|f\|_{\infty} \le \frac{1}{n+1} \stackrel{n \to \infty}{\longrightarrow} 0, $$
womit die Behauptung gezeigt ist. \q

\begin{Bem*}
Sind $X$ ein nicht-leerer kompkater toplologischer Raum und $Y$ ein normierter $\K$-Vektorraum, so gilt
$$ \mathcal{C}_c(X,Y) = \mathcal{C}_0(X,Y) = \mathcal{C}_b(X,Y) = \mathcal{C}(X,Y). $$
\end{Bem*}

\begin{HS}[Satz von \textsc{Stone-Weierstraß}] \index{Satz!von \textsc{Stone-Weierstraß}} \label{FA.3.7} $\,$

\noindent \textbf{Vor.:} Seien $X$ ein kompakter topologischer Raum und $A$ eine $\K$-Unteralgebra von $\mathcal{C}(X,\K)$, die nicht notwendig $1_X$ enthalten muß, derart, daß gilt
\begin{itemize}
\item[(a)] $A$ ist \emph{punktetrennend}, d.h.\ $\forall_{x,y \in X, \, x \ne y} \exists_{f \in A} \, f(x) \ne f(y)$,
\item[(b)] $\forall_{x \in X} \exists_{f \in A} \, f(x) \ne 0$,
\item[(c)] $\forall_{f \in A} \, \overline{f} = {\rm Re} f - \i \, {\rm Im} f \in A$.\footnote{Im Falle $\K = \R$ gilt per definitionem ${\rm Im} \, f = 0$ für alle $f \in A$, also ist (c) dann stets erfüllt.}
\end{itemize} 

\noindent \textbf{Beh.:} $A$ liegt dicht in $(\mathcal{C}(X,\K), \| \ldots \|_{\infty})$.
\end{HS}

Wir bereiten den Beweis des Hauptsatzes durch folgendes Lemma vor.

\begin{Lemma} \label{FA.3.8}
Durch
\begin{equation} \label{FA.3.8.S}
p_0 := 0 ~~ \mbox{ und } ~~ \forall_{n \in \N} \, p_{n+1} := p_n + \frac{x^2 - {p_n}^2}{2}
\end{equation}
wird eine Folge $(p_n)_{n \in \N}$ von Polynomfunktionen $\R \to \R$ definiert derart, daß $(p_n|_{[-1,1]})_{n \in \N}$ gleichmäßig gegen $|x|_{[-1,1]}|$ konvergiert.
\end{Lemma}

\textit{Beweis.} Wir zeigen zunächst
\begin{equation} \label{FA.3.8.1}
\forall_{n \in \N} \forall_{t \in [-1,1]} \, 0 \le p_n(t) \le |t|.
\end{equation}

{[} Zu (\ref{FA.3.8.1}): Für $n=0$ ist die Behauptung trivial.
Sei daher $n \in \N$ und gelte die Behauptung für $n$.
Dann folgt für jedes $t \in [-1,1]$
$$ |t| - p_{n+1}(t) = |t| - p_n(t) - \frac{t^2 + p_n(t)^2}{2} = \underbrace{(|t|-p_n(t))}_{\ge 0} \, \underbrace{( 1 - \frac{1}{2}(|t| + \underbrace{p_n(t)}_{\le |t|}) )} _{\ge 1 - |t| \, \ge 0} \ge 0 $$
und 
$$ p_{n+1}(t) = \underbrace{p_n(t)}_{\ge 0} + \frac{1}{2} \underbrace{(t^2 - \underbrace{p_n(t)}_{\le t^2}))}_{\ge 0} \ge 0, $$
also ist die Behauptung auch für $n+1$ gezeigt. {]}

Als nächstes beweisen wir
\begin{equation} \label{FA.3.8.2}
\forall_{n \in \N} \forall_{t \in [-1,1]} \, p_n(t) \le p_{n+1}(t).
\end{equation}

{[} Zu (\ref{FA.3.8.2}): Für jedes $n \in \N$ und jedes $t \in [-1,1]$ gilt
$$ (\ref{FA.3.8.1}) \Longrightarrow p_n(t)^2 \le t^2 \Longrightarrow p_n(t) \le p_n + \frac{t^2 - t^2}{2} \Longrightarrow p_n(t) \le p_{n+1}(t), $$
d.h.\ es gilt (\ref{FA.3.8.2}). {]}

Aus (\ref{FA.3.8.1}), (\ref{FA.3.8.2}) und dem Satz von \textsc{Bolzano-Weierstraß} folgt die Existenz einer Funktion $f \: [-1,1] \to \R$ mit
\begin{equation} \label{FA.3.8.3}
\forall_{t \in [-1,1]} \, f(t) = \lim_{n \to \infty} p_n(t) \ge 0,
\end{equation}
und nach (\ref{FA.3.8.S}) gilt $f(t) = f(t) + \frac{t^2 - f(t)^2}{2}$, also $f(t) = t^2$, für jedes $t \in [0,1]$.
Wegen (\ref{FA.3.8.3}) bedeutet dies
$$ f = |x|_{[-1,1]}| \in \mathcal{C}([-1,1],\R). $$
Somit folgt das Lemma aus (\ref{FA.3.8.2}), (\ref{FA.3.8.3}) und dem Satz von \textsc{Dini} \ref{FA.1.38}. \q
\A
\textit{Beweis des Hauptsatzes.} Wir betrachten $\mathcal{C}(X,\K)$ in diesem Beweis stets als durch die Supremumsnorm $\|\ldots\|_{\infty}$ normierte $\K$-(Banach-)Algebra. 

1. Fall: $\K = \R$.
Nach \ref{FA.2.4} (ii) ist die abgeschlossene Hülle $\overline{A}$ von $A$ ebenfalls eine $\R$-Unteralgebra von $\mathcal{C}(X,\R)$.
Wir behaupten
\begin{equation} \label{FA.3.7.1}
\forall_{f,g \in \overline{A}} \,\, \max(f,g), \min(f,g) \in \overline{A}.
\end{equation}

{[} Zu (\ref{FA.3.7.1}): Da für alle $f,g \in \R^X$ gilt
$$ \max(f,g) = \frac{1}{2}(f+g+|f-g|), ~~ \min(f,g) = \frac{1}{2}(f+g-|f-g|), $$
ist nur zu zeigen, daß mit $f$ auch $|f|$ ein Element von $\overline{A}$ ist.

Hierzu sei $f \in \overline{A}$.
Ohne Beschränkung der Allgemeinheit gelte 
$$ C := \|f\|_{\infty} \in \R_+, $$
denn andernfalls gilt $f = 0$ und $|f| = f \in \overline{A}$.
Wir zeigen
\begin{equation} \label{FA.3.7.2}
\forall_{\varepsilon \in \R_+} \, \overline{A} \cap U_{\varepsilon}(|f|) \ne \emptyset,
\end{equation}
denn dann folgt $f \in \overline{\overline{A}} = \overline{A}$.
Sei $\tilde{f} := \frac{f}{C} \: X \to [-1,1]$.
Wegen Lemma \ref{FA.3.8} existiert eine Polynomfunktion $p \: \R \to \R$ derart, daß gilt $p(0) = 0$ und
$$ \forall_{t \in [-1,1]} \, | |t| - p(t) | < \frac{\varepsilon}{2 C}. $$
Aus $p(0) = 0$ ergibt sich $\forall_{t \in [-1,1]} \, p(t) = \sum_{k=1}^n a_k t^k$ mit $a_1, \ldots, a_n \in \R$, also folgt
$$ \forall_{x \in X} \, \frac{1}{C} \left| |f(x)| - \sum_{k=1}^n a_k \frac{f(x)^k}{C^{k-1}} \right| = | |\tilde{f}(x)| - p(\tilde{f}(x)) | < \frac{\varepsilon}{2 C} $$
d.h.\ $\| |f| - \sum_{k=1}^n a_k \frac{f^k}{C^{k-1}} \|_{\infty} \le \frac{\varepsilon}{2} < \varepsilon$.
Mit $f$ ist auch $\sum_{k=1}^n a_k \frac{f^k}{C^{k-1}}$ ein Element von $\overline{A}$, also ist (\ref{FA.3.7.2}) gezeigt. {]}

Seien nun $f \in \mathcal{C}(X,\R)$ und $\varepsilon \in \R_+$.
Zum Nachweis des 1. Falles genügt es zu zeigen, daß gilt
\begin{equation} \label{FA.3.7.3}
\exists_{g \in \overline{A}} \,\, \|f - g\|_{\infty} < \varepsilon,
\end{equation}
d.h.\ $\overline{A} \cap U_{\varepsilon}(f) \ne \emptyset$, denn dann folgt aus der Beliebigkeit von $\varepsilon \in \R_+$: $f \in \overline{\overline{A}} = \overline{A}$.

Um (\ref{FA.3.7.3}) zu zeigen, beweisen wir nacheinander 
\begin{gather}
\forall_{a,b \in X} \exists_{g_{a,b} \in \overline{A}} \,\, g_{a,b}(a) = f(a) \, \wedge \, g_{a,b}(b) = f(b), \label{FA.3.7.4} \\
\forall_{a \in X} \exists_{g_a \in \overline{A}} \,\, g_a(a) = f(a) \, \wedge \, f - \varepsilon < g_a. \label{FA.3.7.5}
\end{gather}

{[} Zu (\ref{FA.3.7.4}): Seien $a,b \in X$.
Im Falle $a=b$ existiert nach Voraussetzung (b) eine Funktion $g \in A$ mit $g(b) \ne 0$ und
$$ g_{b,b} := \frac{f(b)}{g(b)} \, g \in A $$
leistet das Gewünschte.
Im Falle $a \ne b$ existiert nach Voraussetzung (a) eine Funktion $h \in A$ mit $h(a) \ne h(b)$.
Ohne Einschränkung gelte $h(a) \ne 0$.
Dann folgt
$$ \tilde{h} := \frac{(h(b) - h) \, h}{(h(b) - h(a)) \, h(a)} = \frac{h(b)}{(h(b) - h(a)) \, h(a)} - \frac{h^2}{(h(b) - h(a)) \, h(a)} \in A $$
und $\tilde{h}(a) = 1$ sowie $\tilde{h}(b) = 0$.
Daher gilt auch
$$ g_{a,b} := (f(a) - g_{b,b}(a)) \, \tilde{h} + g_{b,b} \in A, $$
wobei $g_{b,b} \in A$ mit $g_{b,b}(b) = f(b)$ wie oben sei, $g_{a,b}(a)  = f(a)$ und $g_{a,b}(b) = f(b)$.

Zu (\ref{FA.3.7.5}): Sei $a \in X$.
Für jedes $b \in X$ gilt
$$ V_b := \{ x \in X \, | \, f(x) - \varepsilon < g_{a,b}(x) \} \in \U(b,X), $$
wobei $g_{a,b}$ wie in (\ref{FA.3.7.4}) sei.
Dann ist $(V_b)_{b \in X}$ eine offene Überdeckung von $X$, die wegen der Kompaktheit von $X$ eine endliche Teilüberdeckung $(V_{b_i})_{i \in \{1, \ldots, n\}}$ besitzt.
Aus (\ref{FA.3.7.1}) folgt
$$ g_a := \max ( g_{a,b_i}, i \in \{1, \ldots, n\} ) \in \overline{A}, $$
und es gilt nach (\ref{FA.3.7.4}): $g_a(a) = f(a)$.
Des weiteren existiert zu jedem $x \in X$ ein $i \in \{1, \ldots, n\}$ mit $x \in V_{b_i}$, also ergibt die Definition von $V_{b_i}$
$$ f(x) - \varepsilon < g_{a,b_i}(x) \le g_a(x). $$
Damit ist (\ref{FA.3.7.5}) gezeigt. {]}

Wir zeigen nun (\ref{FA.3.7.3}):
Für jedes $a \in X$ gilt
$$ U_a := \{ x \in X \, | \, g_a(x) < f(x) + \varepsilon \} \in \U(a,X), $$
wobei $g_a$ wie in (\ref{FA.3.7.5}) sei.
Dann ist $(U_a)_{a \in X}$ eine offene Überdeckung von $X$, die wegen der Kompaktheit von $X$ eine endliche Teilüberdeckung $(U_{a_i})_{i \in \{1, \ldots, m\}}$ besitzt.
Aus (\ref{FA.3.7.1}) folgt
$$ g := \min ( g_{a_i}, i \in \{1, \ldots, m\} ) \in \overline{A}, $$
und es gilt nach (\ref{FA.3.7.5}): $f - \varepsilon < g$.
Des weiteren existiert zu jedem $x \in X$ ein $i \in \{1, \ldots, m\}$ mit $x \in U_{a_i}$, also ergibt die Definition von $U_{a_i}$
$$ g(x) \le g_a(x) < f(x) + \varepsilon. $$
Somit gilt $f - \varepsilon < g < f + \varepsilon$, also folgt aus der Kompaktheit von $X$ und der Stetigkeit von $f-g$: $\|f-g\|_{\infty} < \varepsilon$.
Damit ist (\ref{FA.3.7.3}), also auch der 1. Fall, bewiesen.

2. Fall: $\K = \C$.
Wir setzen
$$ A_{\R} := A \cap \mathcal{C}(X,\R). $$
Dann gilt nach Voraussetzung (c) für jedes $f \in A$
$$ {\rm Re} \, f = \frac{1}{2}(f + \overline{f}) \in A_{\R} ~~ \mbox{ sowie } ~~ {\rm Im} \, f = \frac{1}{2 \i}(f - \overline{f}) \in A_{\R}, $$
und die $\R$-Algebra $A_{\R}$ erfüllt die Voraussetzungen (a) und (b):
Sind nämlich $a,b \in X$ mit $a \ne b$, so existiert nach Voraussetzung (a) eine Funktion $f \in A$ mit $f(a) \ne f(b)$, d.h.\ ${\rm Re} f(a) \ne {\rm Re} f(b)$ oder ${\rm Im} f(a) \ne {\rm Im} f(b)$, und es gilt ${\rm Re} f, {\rm Im} f \in A$; oder ist $a \in X$, so existiert nach Voraussetzung (b) eine Funktion $f \in A$ mit $f(a) \ne 0$, d.h.\ ${\rm Re} f(a) \ne 0$ oder ${\rm Im} \, f(a) \ne 0$, und es gilt ${\rm Re} \, f, {\rm Im} \, f \in A$.
Daher folgt aus dem 1. Fall
$$ \overline{A_{\R}} = \mathcal{C}(X,\R). $$
Somit existieren zu $f \in \mathcal{C}(X,\C)$ und $\varepsilon \in \R_+$ Funktionen $g,h \in A_{\R}$ mit
$$ \| ({\rm Re} \, f) - g \|_{\infty} < \frac{\varepsilon}{2} ~~ \mbox{ und } ~~ \| ({\rm Im} \, f) - h \|_{\infty} < \frac{\varepsilon}{2}, $$
also gilt für jedes $x \in X$
$$ |f(x) - (g + \i \, h)(x)| = \sqrt{(({\rm Re} \, f) - g)^2(x) + (({\rm Im} \, f) - h)^2(x)} \le \sqrt{\frac{2 \varepsilon^2}{4}} = \frac{\varepsilon}{\sqrt{2}} < \varepsilon, $$
d.h.\ $\| f - (g + \i \, h) \|_{\infty} < \varepsilon$.
Aus $g,h \in A_{\R} \subset A$ folgt ferner $g + \i h \in A$, also ist die Behauptung auch im 2. Falle gezeigt. \q

\begin{Kor}[Approximationssatz von \textsc{Weierstraß}] \index{Satz!Approximations-!von \textsc{Weierstraß}} \label{FA.3.9}
Es seien $n \in \N_+$, $K$ eine nicht-leere kompakte Teilmenge von $\R^n$ und $f \in \mathcal{C}(K,\R)$.

Dann läßt sich $f$ auf $K$ gleichmäßig durch Polynomfunktionen $\R^n \to \R$ approximieren.\q
\end{Kor}

\begin{Kor} \label{FA.3.9.F}
Es seien $n \in \N_+$ und $K$ eine nicht-leere kompakte Teilmenge von $\R^n$.

Dann besitzt $(\mathcal{C}(K,\R), \| \ldots \|_{\infty})$ eine abzählbare dichte Teilmenge, nämlich die Menge aller Polynomfunktionen $K \to \R$ mit rationalen Koeffizienten. \q
\end{Kor}

\begin{Def*} \index{Raum!topologischer!separabler}
Ein topologischer Raum, der eine höchstens abzählbare dichte Teilmenge besitzt, heißt \emph{separabel}.
\end{Def*}

\begin{Satz}[Satz von \textsc{Stone-Weierstraß} für lokal-kompakte Räume] \index{Satz!von \textsc{Stone-Weierstraß}} \label{FA.3.10} $\,$

\noindent \textbf{Vor.:} Seien $X$ ein lokal-kompakter topologischer Raum, der nicht bereits kompakt ist, und $A$ eine $\K$-Un\-ter\-al\-ge\-bra von $\mathcal{C}_0(X,\K)$ derart, daß gilt
\begin{itemize}
\item[(a)] $A$ ist \emph{punktetrennend}, d.h.\ $\forall_{x,y \in X, \, x \ne y} \exists_{f \in A} \, f(x) \ne f(y)$,
\item[(b)] $\forall_{x \in X} \exists_{f \in A} \, f(x) \ne 0$,
\item[(c)] $\forall_{f \in A} \, \overline{f} = {\rm Re} f - \i \, {\rm Im} f \in A$.
\end{itemize} 

\noindent \textbf{Beh.:} $A$ liegt dicht in $(\mathcal{C}_0(X,\K), \| \ldots \|_{\infty})$.
\end{Satz}

\textit{Beweisskizze.} Sei $\widehat{X} := X^{\infty}$ ,,die`` Einpunkt-Kompaktifizierung von $X$ wie in \ref{FA.1.48.B}.
Offenbar kann jedes $f \in \mathcal{C}_0(X,\K)$ durch
$$ \hat{f}|_X := f ~~ \mbox{ und } ~~ \hat{f}(\infty) := 0 $$
zu einer stetigen Abbildung $\hat{f} \in \mathcal{C}(\widehat{X},\K)$ fortgesetzt werden.
Dann ist
$$ \widehat{A} := \{ \lambda 1_{\widehat{X}} + \hat{f}  \, | \, \lambda \in \K \, \wedge \, f \in \mathcal{C}_0(X,\K) \} $$
eine $\K$-Unteralgebra von $\mathcal{C}(\widehat{X},\K)$, die die Voraussetzungen des klassischen Satzes von \textsc{Stone-Weierstraß} \ref{FA.3.7} erfüllt -- beim Nachweis der Eigenschaft (a) für $\widehat{A}$ gehen sowohl Eigenschaft (a) als auch Eigenschaft (b) für $A$ ein.
Folglich existieren insbesondere zu jedem $f \in \mathcal{C}_0(X,\K)$ und jedem $\varepsilon \in \R_+$ eine Zahl $\lambda \in \K$ und eine Funktion $g \in \mathcal{C}_0(X,\K)$ derart, daß gilt $\|\hat{f} - \lambda 1_{\widehat{X}} - \hat{g} \|_{\infty} < \frac{\varepsilon}{2}$.
Hieraus und aus $\hat{f}(\infty) = \hat{g}(\infty) = 0$ ergibt sich zunächst $|\lambda| < \frac{\varepsilon}{2}$ und sodann
$$ \|f - g\|_{\infty} = \|\hat{f} - \hat{g} \|_{\infty} \le \|\hat{f} - \lambda 1_{\widehat{X}} - \hat{g} \|_{\infty} + |\lambda| < \varepsilon, $$
womit der Satz bewiesen ist. \q

\begin{Bem} \label{FA.3.11}
Wir wollen abschließend noch erwähnen, daß für alle end\-lich-dimensionalen normierten $\R$-Vektorräume $V$, $W$ und jede offene Teilmenge $U$ von $V$ durch
\begin{gather*}
\forall_{k \in \N_+} \, \boxed{\mathcal{C}^k(U,W)} := \{ f \in W^U \, | \, \mbox{$f$ ist $k$-mal stetig differenzierbar} \}, \\
\boxed{\mathcal{C}^{\infty}(U,W)} := \bigcap_{k \in \N_+} \mathcal{C}^k(U,W)
\end{gather*}
$\R$-Untervektorräume von $\mathcal{C}(U,W)$ definiert werden.
Im Falle $V=\R$ können diese Definitionen auch für abgeschlossene bzw.\ halb-offene Intervalle $I$ von $\R$ anstelle von $U$ gegeben werden.

Der Leser mache sich klar, daß z.B.\ $(f_n)_{n \in \N}$, definiert durch
$$ \forall_{n \in \N} \forall_{t \in {]}-1,1{[}} \, f_n(t) := \sqrt{ \frac{1}{n+1} + t^2 }, $$
eine Folge in $\mathcal{C}^1({]}-1,1{[},\R)$ ist, die bzgl.\ der Supremumsnorm gegen die Betragsfunktion in $\mathcal{C}({]}-1,1{[},\R)$, welche in Null nicht differenzierbar ist, konvergiert.
$\mathcal{C}^1({]}-1,1{[},\R)$ ist also nicht abgeschlossen in $(\mathcal{C}({]}-1,1{[},\R), \| \ldots \|_{\infty})$.
\end{Bem}

\subsection*{Übungsaufgaben}

\begin{UA} \label{FA.ÜA.BV}
Es seien $a,b \in \R$ mit $a<b$.

Für eine Funktion $f \: {[}a,b{]} \to \K$ ist die \emph{Variation $\boxed{V_a^b(f)}$ von $f$} per definitionem das Supremum aller reellen Zahlen
$$ \sum_{i=1}^n |f(a_i) - f(a_{i-1})|, $$
wobei $n \in \N_+$ und $a = a_0 < a_1 < \ldots < a_n = b$ eine Zerlegung von ${[}a,b{]}$ seien.
\pagebreak

Wir setzen
$$ \boxed{{\rm BV}_{\K}([a,b])} := \{ f \in \K^{{[}a,b{]}} \, | \, V_a^b(f) < \infty \} $$
und nennen ein Element dieser Menge eine \emph{Funktion beschränkter Variation}\index{Funktion!beschränkter Variation} oder \emph{rektifizierbar}.

Zeige, daß ${\rm BV}_{\K}([a,b])$ ein $\K$-Vektorraum ist und durch
$$ \forall_{f \in {\rm BV}_{\K}([a,b])} \, \|f\|_{\rm BV} := |f(a)| + V_a^b(f) $$
eine Norm auf ${\rm BV}_{\K}([a,b])$ definiert wird, die ${\rm BV}_{\K}([a,b])$ zu einem $\K$-Banach\-raum macht.
\end{UA}

Tip: Verwende, daß $(\mathcal{B}([a,b],\K), \| \ldots \|_{\infty})$ nach \ref{FA.3.2} (iii) ein $\K$-Banachraum ist.

\begin{UA}
Es seien $a,b \in \R$ mit $a<b$.

Eine Funktion $f \: {[}a,b{]} \to \R$ heißt \emph{Treppenfunktion auf ${[}a,b{]}$}\index{Treppenfunktion}\index{Funktion!Treppen-} genau dann, wenn eine Zerlegung $a = a_0 < a_1 < \ldots < a_n = b$ von ${[}a,b{]}$ und $\alpha_1, \ldots, \alpha_n \in \R$ mit
$$ \forall_{i \in \{1, \ldots, n\}} \, f|_{{]}a_{i-1}, a_i{[}} = \alpha_i $$
existieren.
Für ein solches $f$ ist
$$ T \, f := \sum_{i=1}^n (a_i - a_{i-1}) \, \\alpha_i $$
offenbar wohldefiniert.

Zeige:
\begin{itemize}
\item[(i)] Die Menge $\boxed{\mathcal{S}([a,b])}$ aller Treppenfunktionen auf ${[}a,b{]}$ bildet einen abgeschlossenen $\R$-Untervektorraum von $\mathcal{B}({[}a,b{]},\R)$ zusammen mit der Supremumsnorm, und
$$ T \: \mathcal{S}([a,b]) \longrightarrow \R, ~~ f \longmapsto T \, f, $$
ist ein gleichmäßig stetiges Funktional.
Daher existiert eine eindeutig bestimmte stetige Fortsetzung $\overline{T} \: \overline{\mathcal{S}([a,b])} \to \R$ von $T$.

\begin{Bem*}
Ein Element von $\overline{\mathcal{S}([a,b])}$ nennt man eine \emph{Regelfunktion auf ${[}a,b{]}$}\index{Funktion!Regel-}.
Ist $f \in \overline{\mathcal{S}([a,b])}$, so heißt $\overline{T} \, f$ das \emph{Regelintegral von $f$ über ${[}a,b{]}$}\index{Integral!Regel-}.
\end{Bem*}

\item[(ii)] $f \in \overline{\mathcal{S}([a,b])}$
\newline $\Longrightarrow$ $f$ Riemann-integrierbar über ${[}a,b{]}$ und $\overline{T} \, f = \int_a^b f(t) \, \d t$.\footnote{Hier bezeichne das Integral auf der rechten Seite (noch) das Riemann-Integral.} 

\item[(iii)] $\mathcal{C}({[}a,b{]},\R) \subset \overline{\mathcal{S}([a,b])}$.
\end{itemize}
\end{UA}

\begin{UA} \label{FA.3.B.U}
In dieser Aufgabe sollen Einzelheiten ausgeführt werden, die in Beispiel \ref{FA.3.B} nur erwähnt wurden.
Beweise, daß die Abbildungen
$$ T \: \ell^{\infty}_{\K} \longrightarrow (\ell^1_{\K})', ~~ (x_n)_{n \in \N} \longmapsto \left( (\tilde{x}_n)_{n \in \N} \mapsto \sum_{n=0}^{\infty} x_n \, \tilde{x}_n \right) $$
und
\begin{equation*} \label{FA.3.B.U.S}
\widetilde{T} \: \ell^1_{\K} \longrightarrow (\mathbf{c_0}_{\K})', ~~ (x_n)_{n \in \N} \longmapsto \left( (\tilde{x}_n)_{n \in \N} \mapsto \sum_{n=0}^{\infty} x_n \, \tilde{x}_n \right)
\end{equation*}
wohldefiniert und tatsächlich $\K$-Vektorraum-Isometrien sind.
Beweise hierbei auch, daß
$$ S \: (\mathbf{c_0}_{\K})' \longrightarrow \ell^1_{\K}, ~~ F \longmapsto (F(e_n))_{n \in \N}, $$
wobei $e_n \in \mathbf{c_0}_{\K}$ für jedes $n \in \N$ durch $\forall_{k \in \N} \, e_{n,k} := \delta_{nk}$ definiert sei, eine wohldefinierte $\K$-lineare Abbildung mit $S = \widetilde{T}^{-1}$ ist.
Zeige schließlich $S' \circ T|_{\mathbf{c_0}_{\K}} = i_{\mathbf{c_0}_{\K}}$.
\end{UA}

\begin{UA} \label{FA.3.c}
$\boxed{\mathbf{c}_{\K}}$ bezeichne den $\K$-Vektorraum aller $\K$-wertigen Folgen, die in $\K$ konvergieren, versehen mit der Norm, die durch
$$ \forall_{(x_n)_{n \in \N} \in \mathbf{c}_{\K}} \, \| (x_n)_{n \in \N} \|_{\infty} = \sup\{ |x_n| \, | \, n \in \N \} $$
gegeben ist.\footnote{\textbf{Warnung.} Diese allgemein übliche Notation läßt fälschlicherweise vermuten, daß mit der $\mathbf{c}_{\K}$ zugrundeliegenden Menge $\mathcal{C}(\N,\K) = \K^{\N}$ gemeint ist.}

Beweise die folgenden beiden Aussagen:
\begin{itemize}
\item[(i)] $\mathbf{c}_{\K}$ ist ein $\K$-Banachraum.
\item[(ii)] Die Abbildung 
$$ \ell^1_{\K} \longrightarrow (\mathbf{c}_{\K})', ~~ (x_n)_{n \in \N} \longmapsto \left( (\tilde{x}_n)_{n \in \N} \mapsto x_0 \, \lim_{n \to \infty} \tilde{x}_n + \sum_{n=1}^{\infty} x_{n+1} \, \tilde{x}_n \right) $$
ist wohldefiniert und eine $\K$-Vektorraum-Isometrie.
\end{itemize}
\end{UA}

\begin{Bem*} 
Die Aufgaben \ref{FA.3.B.U} und \ref{FA.3.c} (ii) besagen nicht, daß $\mathbf{c}_{\K}$ und $\mathbf{c_0}_{\K}$ denselben topologischen Dualraum haben.
\end{Bem*}

\begin{UA}[Zweiter Approximationssatz von \textsc{Weierstraß}] \index{Satz!Approximations-!Zweiter -- von \textsc{Weierstraß}} \label{FA.3.W2}
Sei $f \in \mathcal{C}(\R,\K)$ eine $2 \pi$-periodische Funktion.

Beweise, daß eine Folge $2 \pi$-periodischer trigonometrischer Polynomfunktionen\footnote{Eine \emph{$2 \pi$-periodische trigonometrische Polynomfunktion} ist per definitionem eine Funktion $$ \sum_{k=-n}^n c_k \, \e^{\i k x} = \frac{a_0}{2} + \sum_{k=1}^n ( a_k \cos(k x) + b_k \sin(k x) ) \: \R \longrightarrow \K $$ mit $n \in \N$, $\forall_{k \in \{-n, \ldots, n\}} \, c_k \in \C$ sowie $\forall_{k \in \{0, \ldots, n\}} \, a_k := c_k + c_{-k} \, \wedge \, b_k := \i \, (c_k - c_{-k})$ und im Falle $\K = \R$ des weiteren $\forall_{k \in \{0, \ldots, n\}} \, c_{-k} = \overline{c_k}$.} existiert, die gleichmäßig gegen $f$ konvergiert.
\end{UA}

\cleardoublepage
\section{Lebesguesche Integrationstheorie} \label{FAna4}

In diesem und im nächsten Kapitel sei stets $n \in \N_+$.
Wir versehen $\K^n$ des weiteren mit seiner Normtopologie.
Ist im folgenden die Rede von einer Norm, so induziert diese also eine Metrik und die Normtopologie.

\subsection*{Quadermaße} \addcontentsline{toc}{subsection}{Quadermaße}

\begin{Def}[Intervalle, Quader, parkettierbare Mengen] \label{FA.4.1} $\,$
\begin{itemize}
\item[(i)] Eine Teilmenge $J$ von $\R^n$ heißt \emph{Intervall von $\R^n$}\index{Intervall} genau dann, wenn Intervalle $J_1, \ldots, J_n$ von $\R$ mit $J = \bigtimes_{i=1}^n J_i$ existieren.

Mit $\boxed{\mathfrak{I}_n}$ bezeichnen wir die Menge aller Intervalle von $\R^n$.
\begin{Bem*}
Die leere Menge und die einpunktigen Teilmengen von $\R^n$ sind Intervalle.
\end{Bem*}
$J = \bigtimes_{i=1}^n J_i \in \mathfrak{I}_n$, wobei $J_1, \ldots, J_n \in \mathfrak{I}_1$ seien, heißt \emph{entartet} genau dann, wenn $\exists_{i \in \{1, \ldots, n\}} \, \# J_i \le 1$ gilt.
\item[(ii)] Eine Teilmenge $Q$ von $\R^n$ heißt \emph{Quader des $\R^n$}\index{Quader} genau dann, wenn $Q$ ein beschränktes Intervall von $\R^n$ ist.

Mit $\boxed{\mathfrak{Q}_n}$ bezeichnen wir die Menge aller Quader des $\R^n$.
\item[(iii)] Eine Teilmenge $P$ von $\R^n$ heißt \emph{parkettierbar}\index{Menge!parket\-tier\-ba\-re} genau dann, wenn wenn paarweise disjunkte $Q_1, \ldots, Q_k \in \mathfrak{Q}_n$ existieren derart, daß gilt $P = \bigcupdot_{i=1}^k Q_i$.
In diesem Falle heißt $\bigcupdot_{i=1}^k Q_i$ eine \emph{Quaderdarstellung von $P$}.
 
Mit $\boxed{\mathfrak{P}_n}$ bezeichnen wir die Menge der parkettierbaren Teilmengen von $\R^n$.
\end{itemize}
\end{Def}

\begin{Bsp} \label{FA.4.1.B}
Seien $J \in \mathfrak{I}_n$ und $J_1, \ldots, J_n \in \mathfrak{I}_n$ mit $J = \bigtimes_{i=1}^n J_i$.
Dann gilt
\begin{itemize}
\item[(i)] $J$ beschränkt $\Longleftrightarrow$ $\forall_{i \in \{1, \ldots, n\}} \, J_i$ beschränkt,
\item[(ii)] $J$ offen in $\R^n$ $\Longleftrightarrow$ $\forall_{i \in \{1, \ldots, n\}} \, J_i$ offen in $\R$,
\item[(iii)] $J$ abgeschlossen in $\R^n$ $\Longleftrightarrow$ $\forall_{i \in \{1, \ldots, n\}} \, J_i$ abgeschlossen in $\R$,
\item[(iv)] $J$ kompakt $\Longleftrightarrow$ $\forall_{i \in \{1, \ldots, n\}} \, J_i$ kompakt.
\end{itemize}

{[} (i), (ii) ,,$\Leftarrow$`` sowie (iii) ,,$\Leftarrow$`` sind trivial, und (iv) folgt aus (i), (iii).

Zu (ii) ,,$\Rightarrow$``: Zu jedem $x = (x_1,\ldots,x_n) \in J$ existiert nach Voraussetzung eine Zahl $\varepsilon \in \R_+$ derart, daß gilt $U_{\varepsilon}^{\|\ldots\|_{\infty}}(x) \subset J$, also gilt ${]}x_i - \varepsilon, x_i + \varepsilon{[} \subset J_i$ für alle $i \in \{1,\ldots,n\}$.
Daher folgt die Behauptung aus der Beliebigkeit von $p_i \in J_i$. 

Zu (iii) ,,$\Rightarrow$``: Ersetze im Beweis von (ii) ,,$\Rightarrow$`` jeweils $J$ durch $\R^n \setminus J$ sowie $J_i$ durch $\R \setminus J_i$. {]}
\end{Bsp}

\begin{Satz} \label{FA.4.2}
$\mathfrak{P}_n$ ist ein \emph{Mengenkörper}, d.h.\ per definitionem
$$ \forall_{P,P' \in \mathfrak{P}_n} \, P \cap P', P \cup P', P \setminus P' \in \mathfrak{P}_n. $$
\end{Satz}

Wir bereiten den Beweis von \ref{FA.4.2} durch das folgende Lemma vor.
\pagebreak

\begin{Lemma} \label{FA.4.3}
Seien $Q, \widetilde{Q} \in \mathfrak{Q}_n$ mit $Q \subset \widetilde{Q}$.

Dann existieren paarweise disjunkte $Q_0, Q_1, \ldots, Q_{2n} \in \mathfrak{Q}_n$ mit
$$ Q_0 = Q, ~~ \bigcupdot_{i=0}^{2n} Q_i = \widetilde{Q} ~~ \mbox{ und } ~~ \forall_{k \in \{0, \ldots, 2n\}} \, \bigcupdot_{i=0}^k Q_i \in \mathfrak{Q}_n. $$
\end{Lemma}

\textit{Beweis.} Im Falle $Q = \emptyset$ gilt die Behauptung mit $Q_0 = \ldots = Q_{2n-1} = \emptyset$, $Q_{2n} = \widetilde{Q}$.
Sei daher $Q \ne \emptyset$.
Wir beweisen das Lemma nun durch vollständige Induktion über $n \in \N_+$:

Für $n=1$ gilt die Behauptung mit $Q_0 = Q$, $Q_1 = \{ t \in \widetilde{Q} \setminus Q \, | \, t \le \inf Q \}$ und $Q_2 = \{ t \in \widetilde{Q} \setminus Q \, | \, t \ge \sup Q \}$.

Gelte die Behauptung für $(n-1) \in \N_+$ und seien $J_1, \ldots, J_n, \widetilde{J}_1, \ldots, \widetilde{J}_n \in \mathfrak{Q}_1$ mit $Q = \bigtimes_{i=1}^n J_i$ und $\widetilde{Q} = \bigtimes_{i=1}^n \widetilde{J}_i$.  
Wir setzen $Q' := \bigtimes_{i=1}^{n-1} J_i$ und $\widetilde{Q}' = \bigtimes_{i=1}^{n-1} \widetilde{J}_i$.
Dann gilt $Q', \widetilde{Q}' \in \mathfrak{Q}_{n-1}$ sowie $\emptyset \ne Q' \subset \widetilde{Q}'$, also existieren nach Induktionsvoraussetzung $Q_0', Q_1' \ldots, Q_{2n-2}' \in \mathfrak{Q}_{n-1}$ mit
$$ Q_0' = Q', ~~ \bigcupdot_{i=0}^{2n-2} Q_i' = \widetilde{Q}' ~~ \mbox{ und } ~~ \forall_{k \in \{0, \ldots, 2n-2\}} \, \bigcupdot_{i=0}^k Q_i' \in \mathfrak{Q}_n. $$
Wegen $J_n, \widetilde{J}_n \in \mathfrak{Q}_1$ sowie $\emptyset \ne J_n \subset \widetilde{J}_n$ existieren nach dem bereits bewiesenen Fall für $n=1$ paarweise disjunkte $I_0, I_1, I_2 \in \mathfrak{Q}_1$ mit $I_0 = J_n$, $I_0 \cup I_1 \cup I_2 = \widetilde{J}_n$ und $I_0 \cup I_1 \in \mathfrak{Q}_n$.
Dann werden durch
$$ \forall_{i \in \{0, \ldots, 2n-2\}} \, Q_i := Q_i' \times J_n, ~~ Q_{2n-1} := \widetilde{Q}' \times I_1, ~~ Q_{2n} := \widetilde{Q}' \times I_2 $$
paarweise disjunkte Elemente von $\mathfrak{Q}_n$ definiert, die offenbar die Behauptung für $n$ erfüllen. \q
\A
\textit{Beweis des Satzes.} Seien $P,P' \in \mathfrak{P}_n$ und $\bigcupdot_{i=1}^k Q_i$ bzw.\ $\bigcupdot_{j=1}^l Q_j'$ Quaderdarstellungen von $P$ bzw.\ $P'$.

1.) Es gilt $P \cap P' = \left( \bigcupdot_{i=1}^k Q_i \right) \cap \left( \bigcupdot_{j=1}^l Q_j' \right) = \bigcupdot_{(i,j) \in \{1, \ldots, k\} \times \{1, \ldots, l\}} Q_i \cap Q_j'$ und $\forall_{i,j \in \{1, \ldots, k\} \times \{1, \ldots, l\}} \, Q_i \cap Q_j' \in \mathfrak{Q}_n$.

2.) Seien $Q, Q' \in \mathfrak{Q}_n$.
Wir zeigen
\begin{gather}
Q \setminus Q' \in \mathfrak{P}_n, \label{FA.4.2.1} \\
P \setminus Q' \in \mathfrak{P}_n. \label{FA.4.2.2}
\end{gather}

{[} Zu (\ref{FA.4.2.1}): Aus $Q,Q' \in \mathfrak{Q}_n$ folgt zunächst $Q \cap Q' \in \mathfrak{Q}_n$ und $Q \cap Q' \subset Q$.
Aus Lemma \ref{FA.4.3} folgt daher die Existenz paarweise disjunkter $Q_0, \ldots, Q_{2n} \in \mathfrak{Q}_n$ mit $Q_0 = Q \cap Q'$ sowie $\bigcupdot_{i=0}^{2n} Q_i = Q$, also gilt $Q \setminus Q' = \bigcupdot_{i=1}^{2n} Q_i \in \mathfrak{P}_n$.

Zu (\ref{FA.4.2.2}): Für jedes $i \in \{1, \ldots, k\}$ ist $Q_i \setminus Q'$ nach (\ref{FA.4.2.1}) disjunkte Vereinigung endlich vieler Quader, also ist auch $P \setminus Q' = \bigcupdot_{i=1}^k (Q_i \setminus Q')$ disjunkte Vereinigung endlich vieler Quader, d.h.\ $P \setminus Q' \in \mathfrak{P}_n$. {]}

$l$-malige Anwendung von (\ref{FA.4.2.2}) ergibt nun sukzessiv
\begin{gather*}
P \setminus Q_1' \in \mathfrak{P}_n, \\
P \setminus (Q_1' \cupdot Q_2') = (P \setminus Q_1') \setminus Q_2'  \in \mathfrak{P}_n, \\
\vdots \\
P \setminus P' = P \setminus \left( \bigcupdot_{i=1}^l Q_i' \right) = \left( \ldots \left( (P \setminus Q_1') \setminus Q_2' \right) \ldots \right) \setminus Q_l' \in \mathfrak{P}_n.
\end{gather*}
\pagebreak

3.) $P \cup P'$ ist eine beschränkte Teilmenge von $\R^n$.
Daher existiert $Q \in \mathfrak{Q}_n$ mit $P \cup P' \subset Q$.
Wegen 2.) und 1.) gilt $Q \setminus (P \cup P') = (Q \setminus P) \cap (Q \setminus P') \in \mathfrak{Q}_n$, also nach 2.) auch $P \cup P' = Q \setminus ( Q \setminus (P \cup P') ) \in \mathfrak{P}_n$.

Mit 1.), 2.) und 3.) ist der Satz vollständig bewiesen. \q

\begin{Def}[Quadermaß] \label{FA.4.4}
Ein \emph{Quadermaß auf $\R^n$}\index{Quader!-maß} ist per definitionem eine Funktion 
$$ \varphi \: \mathfrak{Q}_n \longrightarrow \R $$
derart, daß gilt:
\begin{itemize}
\item[(a)] $\varphi$ ist \emph{additiv}, d.h.\ per definitionem
$$ \forall_{Q, Q', Q'' \in \mathfrak{Q}_n} \, \left( Q = Q' \cupdot Q'' \Longrightarrow \varphi(Q) = \varphi(Q') + \varphi(Q'') \right), $$
insbesondere gilt also $\varphi(\emptyset) = 0$.
\item[(b)] $\varphi$ ist \emph{monoton}, d.h.\ per definitionem
$$ \forall_{Q, Q' \in \mathfrak{Q}_n} \, \left( Q \subset Q' \Longrightarrow \varphi(Q) \le \varphi(Q') \right), $$
insbesondere folgt aus (a): $\varphi \ge 0$.
\item[(c)] $\varphi$ ist \emph{regulär}, d.h.\ per definitionem
$$\forall_{Q \in \mathfrak{Q}_n} \forall_{\varepsilon \in \R_+} \exists_{Q' \in \mathfrak{Q}_n} \, Q \subset Q', \, Q' \mbox{ offen in $\R^n$ und } \varphi(Q') \le \varphi(Q) + \varepsilon. $$
\end{itemize}
\begin{Bem*}
Die Regularität eines Quadermaßes stellt eine Beziehung zur Normtopologie des $\R^n$ her.
\end{Bem*}
\end{Def}

Wir möchten nun (in \ref{FA.4.4.B}) Beispiele für Quadermaße geben.
Beim Nachweis davon, daß es sich dabei tatsächlich um solche handelt, wird das folgende Lemma benötigt.

\begin{Lemma} \label{FA.4.4.L}
Seien $Q, Q', Q'' \in \mathfrak{Q}_n$ mit $Q = Q' \cupdot Q''$.
Für $i \in \{1, \ldots, n\}$ seien ferner $J_i, J_i', J_i'' \in \mathfrak{Q}_1$ derart, daß gilt $Q = \bigtimes_{i=1}^n J_i$, $Q' = \bigtimes_{i=1}^n J_i'$  und $Q'' = \bigtimes_{i=1}^n J_i''$.

Dann existiert $i_0 \in \{1, \ldots, n\}$ mit $J_{i_0} = J_{i_0}' \cupdot J_{i_0}''$ und $\forall_{i \in \{1, \ldots,n\}} \, J_i' = J_i'' = J_i$.
\end{Lemma}

\textit{Beweis.} Wir beweisen das Lemma durch vollständige Induktion über $n \in \N_+$.

Im Falle $n=1$ ist nichts zu zeigen.

Gelte die Behauptung für $(n-1) \in \N_+$.

1. Fall: $J_n' = J_n''$, also $J_n = J_n' = J_n''$.
Wir setzen $\widetilde{Q}' := \bigtimes_{i=1}^{n-1} J_i', \widetilde{Q}'' := \bigtimes_{i=1}^{n-1} J_i''$ und $\widetilde{Q} := \widetilde{Q}' \cupdot \widetilde{Q}''$.
Dann gilt $\widetilde{Q}, \widetilde{Q}', \widetilde{Q}'' \in \mathfrak{Q}_{n-1}$ sowie $\widetilde{Q} = \widetilde{Q}' \cupdot \widetilde{Q}''$.
Daher folgt die Behauptung aus der Induktionsvoraussetzung.

2. Fall: $J_n' \ne J_n''$, also $J_n = J_n' \cupdot J_n''$.
Angenommen, es existiert $i_0 \in \{1, \ldots, n-1\}$ mit $J_{i_0}' \ne J_{i_0}''$, d.h.\ $J_{i_0} = J_{i_0}' \cupdot J_{i_0}''$.
Ohne Beschränkung der Allgemeinheit sei $i_0 = 1$.
Dann ist $J_1' \times  J_2 \times \cdots \times J_{n-1} \times J_n''$ eine Teilmenge von $Q$, die sowohl mit $\widetilde{Q}'$ als auch mit $\widetilde{Q}''$ einen leeren Schnitt hat, Widerspruch! \q  
\pagebreak

\begin{Bsp} \label{FA.4.4.B} $\,$
\begin{itemize}
\item[1.)] Das \emph{eindimensionale Volumen} ist als das Quadermaß $\boxed{\mu_1}$ auf $\R$, welches gegeben ist durch
$$ \mu_1(\emptyset) = 0 ~~ \wedge ~~ \forall_{Q \in \mathfrak{Q}_1 \setminus \{ \emptyset \}} \, \mu_1(Q) := \sup Q - \inf Q, $$
definiert.

{[} Daß $\mu_1$ additiv und monoton ist, ist klar.
Zur Regularität seien $Q \in \mathfrak{Q}_1$, $a := \inf Q, b:= \sup Q$ und $\varepsilon \in \R_+$.
Für jedes $t \in \R_+$ ist $Q_t := {]}a - t, b + t{[}$ eine offene Obermenge von $Q$ mit $Q_t \in \mathfrak{Q}_1$ und $\mu_1(Q_t) = b-a + 2t = \mu_1(Q) + 2t$, also gilt $\mu_1(Q_t) \le \mu_1(Q) + \varepsilon$ für hinreichend kleines $t \in \R_+$. {]}
\item[2.)] Seien $k,l \in \N_+$ mit $k+l=n$, $\varphi_1$ ein Quadermaß auf $\R^k$ und $\varphi_2$ ein Quadermaß auf $\R^l$.
Dann wird durch
$$ \forall_{Q_1 \in \mathfrak{Q}_k} \forall_{Q_2 \in \mathfrak{Q}_l} \, (\varphi_1 \times \varphi_2)(Q_1 \times Q_2) := \varphi_1(Q_1) \cdot \varphi(Q_2) $$
ein Quadermaß $\boxed{\varphi_1 \times \varphi_2}$ auf $\R^n$, das sog.\ \emph{Produktquadermaß von $\varphi$ und $\psi$ auf $\R^n$}\index{Quader!-maß!Produkt-}, definiert.

{[} Zur Additivität: Seien $Q, Q', Q'' \in \mathfrak{Q}_n$ mit $Q = Q' \cupdot Q''$.
Für $i \in \{1, \ldots, n\}$ seien ferner $J_i, J_i', J_i'' \in \mathfrak{Q}_1$ derart, daß gilt $Q = \bigtimes_{i=1}^n J_i$, $Q' = \bigtimes_{i=1}^n J_i'$  und $Q'' = \bigtimes_{i=1}^n J_i''$.
Nach dem obigen Lemma \ref{FA.4.4.L} existiert dann $i_0 \in \{1, \ldots, n\}$ mit $J_{i_0} = J_{i_0}' \cupdot J_{i_0}''$ und $\forall_{i \in \{1, \ldots,n\}} \, J_i' = J_i'' = J_i$.
Ohne Beschränkung der Allgemeinheit gelte $i_0 = 1$.
Dann sind $Q_1' := \bigtimes_{i=1}^k J_i' = J_1' \times \left( \bigtimes_{i=2}^k J_i \right)$ und $Q_1'' := \bigtimes_{i=1}^k J_i'' = J_1'' \times \left( \bigtimes_{i=2}^k J_i \right)$ disjunkte Elemente von $\mathfrak{Q}_k$, d.h.\ 
$$ \varphi_1(Q_1' \cupdot Q_1'') = \varphi_1(Q_1') + \varphi_1(Q_1''), $$
sowie $Q_2' := \bigtimes_{i=k+1}^n J_i' = \bigtimes_{i=k+1}^n J_i'' = \bigtimes_{i=k+1}^n J_i \in \mathfrak{Q}_l$ mit $Q' = Q_1' \times Q_2'$ und $Q'' = Q_1'' \times Q_2'$.
Daher folgt $Q = ( Q_1' \times Q_2' ) \cupdot ( Q_1'' \times Q_2' ) = ( Q_1' \cupdot Q_1'' ) \times Q_2'$ und
\begin{eqnarray*}
(\varphi_1 \times \varphi_2)(Q) & = & \varphi_1( Q_1' \cupdot Q_1'' ) \cdot \varphi_2(Q_2') \\
& = & \varphi_1(Q_1') \cdot \varphi_2(Q_2') + \varphi_1(Q_1'') \cdot \varphi_2(Q_2') \\
& = & (\varphi_1 \times \varphi_2)(Q') + (\varphi_1 \times \varphi_2)(Q'').
\end{eqnarray*}

Zur Monotonie: Seien $Q, Q' \in \mathfrak{Q}_n$ mit $Q \subset Q'$.
Seien ferner $Q_1,Q_1' \in \mathfrak{Q}_k$ sowie $Q_2,Q_2' \in \mathfrak{Q}_l$ mit $Q = Q_1 \times Q_2$ und $Q' = Q_1' \times Q_2'$. 
Dann gilt für $i \in \{1,2\}$: $Q_i \subset Q_i'$ und $\varphi_i(Q_i) \le \varphi_(Q_i')$.
Daher folgt 
$$ (\varphi_1 \times \varphi_2)(Q) = \varphi_1(Q_1) \cdot \varphi_2(Q_2) \le \varphi_1(Q_1') \cdot \varphi_2(Q_2') = (\varphi_1 \times \varphi_2)(Q'). $$

Zur Regularität: Seien $Q \in \mathfrak{Q}_n$, $Q_1 \in \mathfrak{Q}_k$ und $Q_2 \in \mathfrak{Q}_l$ mit $Q = Q_1 \times Q_2$ sowie $\varepsilon \in \R_+$.
Dann existiert für $i \in \{1,2\}$ zu jedem $t \in \R_+$ ein offener Quader $Q_{i,t}$ mit $Q_i \subset Q_{i,t}$ sowie $\varphi_i(Q_{i,t}) \le \varphi_i(Q_i) + t$, und es folgt
\begin{eqnarray*}
(\varphi_1 \times \varphi_2)(Q_{1,t} \times Q_{2,t}) & \le & (\varphi_1(Q_1) + t) \cdot (\varphi_2(Q_2) + t) \\
& = & (\varphi_1 \times \varphi_2)(Q) + t \, \varphi_1(Q_1) + t \, \varphi_2(Q_2) + t^2,
\end{eqnarray*}
d.h.\ es gilt $(\varphi_1 \times \varphi_2)(Q_{1,t} \times Q_{2,t}) \le (\varphi_1 \times \varphi_2)(Q) + \varepsilon$ für hinreichend kleines $t \in \R_+$. {]}
\item[3.)] Das Quadermaß $\boxed{\mu_n} := \bigtimes_{i=1}^n \mu_1$ heißt das \emph{$n$-dimensionale Volumen}.
\item[4.)] Seien $M$ eine Teilmenge von $\R^n$, die keinen Häufungspunkt in $\R^n$ besitzt\footnote{\textbf{Definition.} Seien $X$ ein topologischer Raum, $M$ eine Teilmenge von $X$ und $x \in X$. Dann heißt $x$ \emph{Häufungspunkt von $M$ in $X$}, wenn gilt $\forall_{U \in \U(x,X)} \, U \cap (M \setminus \{x\}) \ne \emptyset$.}, -- insbes.\ ist $M$ höchstens abzählbar\footnote{Denn jede beschränkte unendliche Teilmenge von $\R^n$ besitzt wegen des Satzes von \textsc{Bolzano-Weierstraß} einen Häufungspunkt (vgl.\ z.B.\ \cite[Satz 9.19 (ii)]{ElAna}), und $M \subset \R^n$ ist disjunkte Vereinigung abzählbar vieler beschränkter Mengen.} -- und  $m \: M \to \R$ eine positivwertige Funktion, eine sog.\ \emph{diskrete Massenverteilung}.
Dann wird das \emph{Quadermaß \boxed{\varphi_m} der diskreten Massenverteilung $m$}\index{Quader!-maß!einer diskreten Massenverteilung} (für $M = \{0\}, m(0)=1$ nennt man dieses das \emph{Diracmaß}\index{Quader!-maß!Diracsches} und bezeichnet es mit $\boxed{\delta}$) gegeben durch
$$ \forall_{Q \in \mathfrak{Q}_n} \, \varphi_m(Q) := \sum_{p \in Q \cap M} m(p). $$
Der Satz von \textsc{Bolzano-Weierstraß} ergibt die Endlichkeit der Summe.

{[} Additivität und Monotonie von $\varphi_m$ sind trivial.
Die Regularität beweisen wir, indem wir sogar zeigen, daß zu jedem Quader $Q \in \mathfrak{Q}_n$ ein offener Quader $Q' \in \mathfrak{Q}_n$ mit $\varphi_m(Q') = \varphi_m(Q)$ existiert.
Seien also $Q \in \mathfrak{Q}_n$, $J_1, \ldots, J_n \in \mathfrak{Q}_1$ mit $Q = \bigtimes_{i=1}^n J_i$ und $a_i := \sup J_i$ sowie $b_i := \inf J_i$ für $i \in \{1, \ldots, n\}$.
Da $M$ keinen Häufungspunkt in $\R^n$ besitzt, existiert zu jedem $x \in \partial Q = \overline{Q} \setminus Q^{\circ} \stackrel{\ref{FA.4.1.B}}{=} \bigtimes_{i=1}^n J_1 \times \ldots \times J_{i-1} \times \{a_i,b_i\} \times J_{i+1} \times \ldots \times J_n$ eine Zahl $\varepsilon_x \in \R_+$ mit $U^{\|\ldots\|_{\infty}}_{\varepsilon_x}(x) \cap (M \setminus \{x\}) = \emptyset$.
Dann ist $\bigcup_{x \in \partial Q} U^{\|\ldots\|_{\infty}}_{\varepsilon_x}(x)$ eine offene Überdeckung des Kompaktums $\partial Q$.
Daher existieren $x_1, \ldots, x_k \in \partial Q$ mit
$$ \partial Q \subset \bigcup_{j=1}^k U^{\|\ldots\|_{\infty}}_{\varepsilon_{x_j}}(x_j) ~~ \mbox{ und } \forall_{j \in \{1, \ldots, k\}} \, U^{\|\ldots\|_{\infty}}_{\varepsilon_{x_j}}(x_j) \cap M \subset \{x_j\} \subset \partial Q. $$
Wir setzen nun 
$$ \varepsilon := \min \{ \varepsilon_{x_1}, \ldots, \varepsilon_{x_k} \} \in \R_+ $$
und für jedes $i \in \{1, \ldots, n\}$
$$ J_i' := \left\{ \begin{array}{cl} {]} a_i, b_i {[}, & \mbox{ falls } J_i = {]} a_i, b_i {[}, \\
                                     {]} a_i, b_i + \varepsilon {[}, & \mbox{ falls } J_i = {]} a_i, b_i {]}, \\
                                     {]} a_i - \varepsilon, b_i {[}, & \mbox{ falls } J_i = {[} a_i, b_i {[}, \\
                                     {]} a_i - \varepsilon, b_i + \varepsilon {[}, & \mbox{ falls } J_i = {[} a_i, b_i {]}.
           \end{array} \right. $$
Dann ist $Q' := \bigtimes_{i=1}^n J_i' \in \mathfrak{Q}_n$ eine offene Obermenge von $Q$, die keinen Punkt von $M$ enthält, der nicht bereits zu $Q$ gehört, d.h.\ $\varphi_m(Q') = \varphi_m(Q)$. {]}
\item[5.)] Es sei $g \: \R \to \R$ eine monoton wachsende Funktion.
Wir schreiben wie üblich $g(c-) := \lim_{t \nearrow c} g(c)$ sowie $g(c+) := \lim_{t \searrow c} g(c)$ für jedes $c \in \R$ und definieren $\varphi_g \: \mathfrak{Q}_1 \to \R$ durch
$$ \forall_{Q \in \mathfrak{Q}_1} \, \varphi_g(Q) := \left\{ \begin{array}{cl} 0, & \mbox{ falls } Q = \emptyset, \\
                                                                               g(b-) - g(a+), & \mbox{ falls } Q = {]}a,b{[} \mbox{ mit } a,b \in \R, a < b, \\
                                                                               g(b+) - g(a+), & \mbox{ falls } Q = {]}a,b{]} \mbox{ mit } a,b \in \R, a < b, \\
                                                                               g(b-) - g(a-), & \mbox{ falls } Q = {[}a,b{[} \mbox{ mit } a,b \in \R, a < b, \\
                                                                               g(b+) - g(a-), & \mbox{ falls } Q = {[}a,b{]} \mbox{ mit } a,b \in \R, a \le b. \\
                                                     \end{array} \right. $$
Der Leser zeige als Übung, daß hierdurch ein Quadermaß $\boxed{\varphi_g}$ auf $\R$, ein sog.\ \emph{Stieltjessches Maß}\index{Quader!-maß!Stieltjessches}, definiert wird, welches im Falle $g = \id_{\R}$ offenbar mit $\mu_1$ übereinstimmt.
\end{itemize}
\end{Bsp}

\begin{Satz} \label{FA.4.5}
Sei $\varphi \: \mathfrak{Q}_n \to \R$ ein Quadermaß.

Dann läßt sich $\varphi$ auf genau eine Weise fortsetzen zu einer Abbildung
$$ \varphi \: \mathfrak{P}_n \longrightarrow \R $$
derart, daß gilt:
\begin{itemize}
\item[(a)] $\varphi$ ist \emph{additiv}, d.h.\ per definitionem
$$ \forall_{P, P', P'' \in \mathfrak{P}_n} \, \left( P = P' \cupdot P'' \Longrightarrow \varphi(P) = \varphi(P') + \varphi(P'') \right), $$
insbesondere gilt also $\varphi(\emptyset) = 0$.
\item[(b)] $\varphi$ ist \emph{monoton}, d.h.\ per definitionem
$$ \forall_{P, P' \in \mathfrak{P}_n} \, \left( P \subset P' \Longrightarrow \varphi(P) \le \varphi(P') \right), $$
insbesondere folgt aus (a): $\varphi \ge 0$.
\item[(c)] $\varphi$ ist \emph{regulär}, d.h.\ per definitionem
$$\forall_{P \in \mathfrak{P}_n} \forall_{\varepsilon \in \R_+} \exists_{P' \in \mathfrak{P}_n} \, P \subset P', \, P' \mbox{ offen in $\R^n$ und } \varphi(P') \le \varphi(P) + \varepsilon. $$
\end{itemize}

\begin{Zusatz}
Ist $P \in \mathfrak{P}_n$ und ist $P = \bigcupdot_{i=1}^k Q_i$ eine Quaderdarstellung von $P$, so gilt $\varphi(P) = \sum_{i=1}^k \varphi(Q_i)$.
\end{Zusatz}
\end{Satz}

Wir bereiten den Beweis des Satzes durch das folgende Lemma vor.

\begin{Lemma} \label{FA.4.6}
Seien $\varphi \: \mathfrak{Q}_n \to \R$ ein Quadermaß auf $\R^n$ und $Q_1, \ldots, Q_k \in \mathfrak{Q}_n$ paarweise disjunkt derart, daß gilt
$$ Q := \bigcupdot_{i=1}^k Q_i \in \mathfrak{Q}_n. $$

Dann folgt $\varphi(Q) = \sum_{i=1}^k \varphi(Q_i)$.
\end{Lemma}

\textit{Beweis.} Wir führen den Beweis durch vollständige Induktion über $k \in \N_+$.

Der Fall $k=1$ ist trivial.

Gelte die Behauptung für $(k-1) \in \N_+$ und seien $Q_1, \ldots, Q_n \in \mathfrak{Q}_n$ paarweise disjunkt mit $Q := \bigcupdot_{i=1}^k Q_i \in \mathfrak{Q}_n$.
Wegen $Q_k \subset Q$ existieren nach Lemma \ref{FA.4.3} paarweise disjunkte $Q_1', \ldots, Q_{2n}' \in \mathfrak{Q}_n$ mit
$$ Q_0' = Q_k, ~~ \bigcupdot_{i=0}^{2n} Q_i = Q ~~ \mbox{ und } ~~ \forall_{m \in \{0, \ldots, 2n\}} \, \bigcupdot_{i=0}^m Q_i' \in \mathfrak{Q}_n. $$
Dann folgt für jeden Quader $\widetilde{Q} \in \mathfrak{Q}_n$ mit $\widetilde{Q} \subset Q$:
$$ \widetilde{Q} = \bigcupdot_{i=1}^{2n} (Q_i' \cap \widetilde{Q}) \mbox{ ist disjunkte Vereinigung von Quadern} $$
und
$$ \forall_{m \in \{0, \ldots, 2n\}} \, \bigcupdot_{i=0}^m (Q_i' \cap \widetilde{Q}) \in \mathfrak{Q}_n. $$
Hieraus ergibt sich unter der Benutzung der Additivität von $\varphi \: \mathfrak{Q}_n \to \R$
\begin{equation} \label{FA.4.6.1}
\forall_{\widetilde{Q} \in \mathfrak{Q}_n, \, \widetilde{Q} \subset Q} \, \varphi(\widetilde{Q}) = \sum_{i=0}^{2n} \varphi(Q_i' \cap \widetilde{Q}).
\end{equation}

{[} Zu (\ref{FA.4.6.1}): Wir zeigen durch endliche Induktion über $m \in \{0, \ldots, 2n\}$
$$ \forall_{m \in \{0, \ldots, 2n\}} \, \varphi \left( \bigcupdot_{i=0}^m (Q_i' \cap \widetilde{Q}) \right) = \sum_{i=0}^m \varphi(Q_i \cap \widetilde{Q}). $$

Der Fall $m=0$ ist trivial.

Gelte die Behauptung für $(m-1) \in \{0, \ldots, 2n-1\}$.
Dann folgt aus der Additivität von $\varphi \: \mathfrak{Q}_n \to \R$ und der Induktionsvoraussetzung
\begin{eqnarray*}
\varphi \left( \bigcupdot_{i=0}^m (Q_i' \cap \widetilde{Q}) \right) & = & \varphi \left( \left( \bigcupdot_{i=0}^{m-1} (Q_i' \cap \widetilde{Q}) \right) \cup (Q_m' \cap \widetilde{Q}) \right) \\
& = & \varphi \left( \bigcupdot_{i=0}^{m-1} (Q_i' \cap \widetilde{Q}) \right) + \varphi (Q_m' \cap \widetilde{Q}) \\
& = & \sum_{i=0}^{m-1} \varphi (Q_i' \cap \widetilde{Q}) + \varphi (Q_m' \cap \widetilde{Q}) = \sum_{i=0}^m \varphi (Q_i' \cap \widetilde{Q}).
\end{eqnarray*}

Für $m=2n$ ergibt sich gerade (\ref{FA.4.6.1}). {]}

Aus (\ref{FA.4.6.1}) folgt für $\widetilde{Q} = Q$
\begin{equation} \label{FA.4.6.2}
\varphi(Q) = \sum_{i=0}^{2n} \varphi(Q_i') = \varphi(Q_k) + \sum_{i=0}^{2n} \varphi(Q_i')
\end{equation}
und für $\widetilde{Q} = Q_j$ mit $j \in \{1, \ldots, k-1\}$ wegen $Q_0' \cap Q_j = Q_k \cap Q_j = \emptyset$
\begin{equation} \label{FA.4.6.3}
\varphi(Q_j) = \sum_{i=0}^{2n} \varphi(Q_i' \cap Q_j) = \sum_{i=1}^{2n} \varphi(Q_i' \cap Q_j).
\end{equation}

Für jedes $i \in \{1, \ldots, 2n\}$ gilt
$$ Q_i' \subset Q \setminus Q_0' = Q \setminus Q_k = \bigcupdot_{j=1}^{k-1} Q_j, $$
also $Q_i' = \bigcupdot_{j=1}^{k-1} (Q_i' \cap Q_j),$ folglich nach Induktionsvoraussetzung
$$ \varphi(Q_i') = \sum_{j=1}^{k-1} \varphi(Q_i' \cap Q_j). $$
Hieraus folgt schließlich
\begin{eqnarray*}
\varphi(Q) & \stackrel{(\ref{FA.4.6.2})}{=} & \varphi(Q_k) + \sum_{i=1}^{2n} \sum_{j=1}^{k-1} \varphi(Q_i' \cap Q_j) \\
& \stackrel{(\ref{FA.4.6.3})}{=} & \varphi(Q_k) + \sum_{j=1}^{k-1} \varphi(Q_j) = \sum_{i=1}^k \varphi(Q_i),
\end{eqnarray*}
womit das Lemma bewiesen ist. \q
\pagebreak

\textit{Beweis des Satzes.} Sei $\varphi \: \mathfrak{Q}_n \to \R$ das gegebene Quadermaß.
Falls überhaupt eine Fortsetzung $\tilde{\varphi} \: \mathfrak{P}_n \to \R$ von $\varphi$ mit (a), (b) und (b) existiert, so gilt für diese wegen (a)
\begin{equation} \label{FA.4.5.4}
\forall_{P = \bigcupdot_{i=1}^k Q_i \in \mathfrak{P}_n \text{ mit } Q_i \in \mathfrak{Q}_n} \, \tilde{\varphi}(P) = \sum_{i=1}^k \varphi(Q_i).
\end{equation}
Damit ist die Eindeutigkeit von $\tilde{\varphi}$ bereits klar.
Wir haben zu überlegen, daß durch (\ref{FA.4.5.4}) tatsächlich eine Abbildung $\tilde{\varphi} \: \mathfrak{P}_n \to \R$ definiert wird und daß diese Abbildung additiv, monoton sowie regulär ist.

Zur Wohldefiniertheit: Sei $P \in \mathfrak{P}_n$ und seien $\bigcupdot_{i=1}^k Q_i$ sowie $\bigcupdot_{j=1}^l Q_j'$ zwei Quaderdarstellungen von $P$.
Dann folgt
$$ \forall_{i \in \{1, \ldots, k\}} \, Q_i = \bigcupdot_{j=1}^l (\underbrace{Q_i \cap Q_j'}_{\in \mathfrak{Q}_n}) ~~~ \mbox{ und } ~~~ \forall_{j \in \{1, \ldots, l\}} \, Q_j' = \bigcupdot_{i=1}^k (\underbrace{Q_i \cap Q_j'}_{\in \mathfrak{Q}_n}), $$
also nach Lemma \ref{FA.4.6}
$$ \forall_{i \in \{1, \ldots, k\}} \, \varphi(Q_i) = \sum_{j=1}^l \varphi(Q_i \cap Q_j') ~~~ \mbox{ und } ~~~ \forall_{j \in \{1, \ldots, l\}} \, \varphi(Q_j') = \sum_{i=1}^k \varphi(Q_i \cap Q_j'). $$
Daher gilt auch $\sum_{i=1}^k \varphi(Q_i) = \sum_{i=1}^k \sum_{j=1}^l \varphi(Q_i \cap Q_j') = \sum_{j=1}^l \varphi(Q_j')$.

Die Additivität folgt sofort aus (\ref{FA.4.5.4}).

Zur Monotonie: Seien $P,P' \in \mathfrak{P}_n$ mit $P \subset P'$.
Dann folgt
$$ P' = P \cupdot (P' \setminus P) \mbox{ und } P \in \mathfrak{P}_n, (P' \setminus P) \stackrel{\ref{FA.4.2}}{\in} \mathfrak{P}_n, $$
also gilt aufgrund der Additivität von $\tilde{\varphi}$
$$ \tilde{\varphi}(P') = \tilde{\varphi}(P) + \tilde{\varphi}(P' \setminus P) \stackrel{(\ref{FA.4.5.4})}{\ge} \tilde{\varphi}(P). $$

Zur Regularität: Seien $P \in \mathfrak{P}_n$ und $\varepsilon \in \R_+$.
Sei ferner $\bigcupdot_{i=1}^k Q_i$ eine Quaderdarstellung von $P$.
Für jedes $i \in \{1, \ldots, k\}$ existiert wegen der Regularität von $\varphi \: \mathfrak{Q}_n \to \R$ ein offener Quader $Q_i' \in \mathfrak{Q}_n$ mit $Q_i \subset Q_i'$ und $\varphi(Q_i') \le \varphi(Q_i) + \frac{\varepsilon}{k}$.
Dann gilt nach \ref{FA.4.2}: $P' := \bigcup_{i=1}^k Q_i' \in \mathfrak{Q}_n$.
Ferner ist $P'$ offene Obermenge von $P$ und die Additivität sowie die Monotonie von $\tilde{\varphi} \: \mathfrak{P}_n \to \R$ ergeben
\begin{eqnarray*}
\tilde{\varphi}(P') & = & \tilde{\varphi} \left( \bigcup_{i=1}^k Q_i' \right) = \tilde{\varphi} \bigg( \bigcupdot_{i=1}^k \bigg( \overbrace{Q_i' \setminus \bigcup_{j=1}^{i-1} Q_j'}^{\stackrel{\ref{FA.4.2}}{\in} \mathfrak{P}_n} \bigg) \bigg) = \sum_{i=1}^k \tilde{\varphi} \bigg( Q_i' \setminus \bigcup_{j=1}^{i-1} Q_j' \bigg) \\
& \le & \sum_{i=1}^k \underbrace{\tilde{\varphi}(Q_i')}_{= \varphi(Q_i')} \le \sum_{i=1}^k \varphi(Q_i) + \varepsilon = \tilde{\varphi}(P) + \varepsilon.
\end{eqnarray*}

Damit ist \ref{FA.4.5} bewiesen. \q

\begin{Bem*}
Wir werden künftig jedes gegebene Quadermaß $\varphi \: \mathfrak{Q}_n \to \R$ als seine nach dem letzten Satz eindeutig bestimmte Fortsetzung $\varphi \: \mathfrak{P}_n \to \R$ betrachten.
\end{Bem*}

\begin{Kor} \label{FA.4.5.K}
Seien $\varphi \: \mathfrak{P}_n \to \R$ ein Quadermaß und $P,P' \in \mathfrak{P}_n$.
Dann gilt:
\begin{itemize}
\item[(i)] $P \subset P' \Longrightarrow \varphi(P') = \varphi(P' \setminus P) + \varphi(P)$.
\item[(ii)] $\varphi(P \cup P') = \varphi(P) + \varphi(P') - \varphi(P \cap P')$.
\end{itemize}
\end{Kor}

\textit{Beweis.} (i) haben wir beim Nachweis der Monotonie im vorangegangenen Beweis gezeigt, und (ii) folgt aus $P \cup P' = (P \setminus (P \cap P')) \cupdot (P' \setminus (P \cap P')) \cupdot (P \cap P')$, (i) und der Additivität von $\varphi \: \mathfrak{P}_n \to \R$. \q

\subsection*{Treppenfunktionen} \addcontentsline{toc}{subsection}{Treppenfunktionen}

\begin{Def}[Treppenfunktionen] \label{FA.4.7}
Es seien $J \in \mathfrak{I}_n$ ein nicht-leeres Intervall von $\R^n$ und $T \: J \to \R$ eine Funktion.

$T \: J \to \R$ heißt \emph{Treppenfunktion auf $J$}\index{Treppenfunktion}\index{Funktion!Treppen-} genau dann, wenn gilt $\# T(J) < \infty$ und $\forall_{\alpha \in \R \setminus \{0\}} \, \overline{T}^1(\{\alpha\}) \in \mathfrak{P}_n$.

Die Menge aller Treppenfunktionen auf $J$ bezeichnen wir mit $\boxed{\mathcal{S}(J)}$.
\end{Def}

\begin{Satz} \label{FA.4.8}
Sei $J \in \mathfrak{I}_n$ ein nicht-leeres Intervall von $\R^n$.
Dann gilt:
\begin{itemize}
\item[(i)] $T \: J \to \R$ ist genau dann eine Treppenfunktion auf $J$, wenn $k \in \N_+$ und paarweise disjunkte Quader $Q_1, \ldots, Q_k \in \mathfrak{Q}_n$, die in $J$ enthalten sind, sowie $\alpha_1, \ldots, \alpha_k \in \R$ mit $T = \sum_{i=1}^k \alpha_i \, \chi_{Q_i}^J$ existieren.\footnote{Sind $\widetilde{M}$ eine Menge und $M$ eine Teilmenge von $\widetilde{M}$, so bezeichnen wir mit $\boxed{\chi_M^{\widetilde{M}} \: \widetilde{M} \to \R}$ die \emph{charakteristische Funktion von $M$ in $\widetilde{M}$}, die auf $M$ konstant vom Wert $1$ und auf $\widetilde{M} \setminus M$ konstant vom Wert $0$ ist. Im Falle $\widetilde{M} = \R^n$ schreiben wir $\boxed{\chi_M}$ für $\chi_M^{\widetilde{M}}$.}
\item[(ii)] Sind $r \in \N_+$ und $T_1, \ldots, T_r \in \mathcal{S}(J)$, so existieren $k \in \N_+$ und paarweise disjunkte Quader $Q_1, \ldots, Q_k \in \mathfrak{Q}_n$, die in $J$ enthalten sind, derart, daß gilt
$$ \forall_{\rho \in \{1, \ldots, r\}} \, T_{\rho} \in {\rm Span}_{\R} \{ \chi_{Q_1}^J, \ldots, \chi_{Q_k}^J \}.$$
\item[(iii)] Sind $r \in \N_+$, $T_1, \ldots, T_r \in \mathcal{S}(J)$ und $\psi \: \R^r \to \R$ eine Funktion mit $\psi(0) = 0$, so gilt $\psi \circ (T_1, \ldots, T_r) \in \mathcal{S}(J)$.
\item[(iv)] Sind $r \in \N_+$ und $T_1, \ldots, T_r \in \mathcal{S}(J)$, so gilt
$$ T_1 + \ldots + T_r, \, T_1 \cdot \ldots \cdot T_r, \, \sup (T_1, \ldots, T_r), \, \inf ( T_1, \ldots, T_r ) \in \mathcal{S}(J). $$
\item[(v)] Sind $T \in \mathcal{S}(J)$ und $\lambda \in \R$, so gilt
$$ \lambda \, T, \, |T|, \, T^+ = \sup ( T, 0 ), \, T^- = \sup ( -T, 0 ) \in \mathcal{S}(J).\footnote{Für jede reellwertige Funktion $f$ heißt $\boxed{f^+} := \sup(f,0) \ge 0$ der \emph{Positivteil von $f$} und $\boxed{f^-} := \sup(-f,0) = - \inf(f,0) \ge 0$ der \emph{Negativteil von $f$}. Dann gilt also $f = f^+ - f^-$ sowie $|f| = f^+ + f^-$.} $$   
\end{itemize}
\end{Satz}

Auch den Beweis dieses Satzes bereiten wir durch ein Lemma vor.

\begin{Lemma} \label{FA.4.8.L}
Es seien $k,l \in \N_+$.
Ferner seien sowohl $Q_1', \ldots, Q_k' \in \mathfrak{Q}_n$ als auch $Q_1'', \ldots, Q_l'' \in \mathfrak{Q}_n$ paarweise disjunkte Quader.

Dann existieren $r \in \N_+$ und paarweise disjunkte Quader $Q_1, \ldots, Q_r \in \mathfrak{Q}_n$ derart, daß jeder der Quader $Q_1', \ldots, Q_k', Q_1'', \ldots, Q_l''$ Vereinigung gewisser der $Q_1, \ldots, Q_r$ ist.
\end{Lemma}

\textit{Beweisskizze.} Schreibe jeweils
$$ Q_i' \setminus \left( \bigcupdot_{j=1}^l Q_j'' \right) = Q_i' \setminus \left( \left( \bigcup_{\substack{\iota = 1 \\ \iota \ne i}}^k Q_j' \right) \cup \left( \bigcup_{j=1}^l Q_j'' \right) \right) \stackrel{\ref{FA.4.2}}{\in} \mathfrak{P}_n $$
für jedes $i \in \{1, \ldots, k\}$ und
$$ Q_j'' \setminus \left( \bigcupdot_{i=1}^k Q_i' \right) = Q_j'' \setminus \left( \left( \bigcup_{i=1}^k Q_i' \right) \cup \left( \bigcup_{\substack{\iota = 1 \\ \iota \ne j}}^l Q_j'' \right) \right) \stackrel{\ref{FA.4.2}}{\in} \mathfrak{P}_n $$
für jedes $j \in \{1, \ldots, l\}$ als disjunkte Vereinigung von endlich vielen Quadern.
Diese bilden dann zusammen mit den Quadern $Q_i' \cap Q_j''$, $(i,j) \in \{1, \ldots, k\} \times \{1, \ldots, l\}$, die endlich vielen Quader $Q_1, \ldots, Q_r$ mit der gewünschten Eigenschaft. \q
\A
\textit{Beweis des Satzes.} Zu (i): ,,$\Rightarrow$`` Im Falle $T=0$ gilt $T = 0 \cdot \chi_{\emptyset}^J$.
Sei im folgenden $T \ne 0$.
Dann existieren $l \in \N_+$ und paarweise verschiedene $\beta_1, \ldots, \beta_l \in \R \setminus \{0\}$ derart, daß gilt $T(J) \setminus \{0\} = \{\beta_1, \ldots, \beta_l\}$ sowie $\forall_{j \in \{1, \ldots, l\}} \, \overline{T}^1(\{\beta_j\}) \in \mathfrak{P}_n \setminus \{ \emptyset \}$.
Folglich existieren zu jedem $j \in \{1, \ldots, l\}$ ein $m_j \in \N_+$ sowie paarweise disjunkte Quader $Q^j_1, \ldots, Q^j_{m_j} \in \mathfrak{Q}_n$, die in $J$ enthalten sind, mit $\overline{T}^1(\{\beta_j\}) = \bigcupdot_{i=1}^{m_j} Q_i$.
Dann gilt $T = \sum_{j=1}^l \sum_{i=1}^{m_j} \beta_j \, \chi_{Q^j_i}^J$.
Ferner sind $\underbrace{Q^1_1, \ldots, Q^1_{m_1}}_{\text{paarweise disjunkt}}, \ldots, \underbrace{Q^l_1, \ldots, Q^l_{m_1}}_{\text{paarweise disjunkt}}$ paarweise disjunkt, da $\beta_1, \ldots, \beta_l$ paarweise verschieden sind.

,,$\Leftarrow$`` Es gilt $T(J) \subset \{0, \alpha_1, \ldots, \alpha_k\}$ und für jedes $\alpha \in \R \setminus \{0\}$
$$ \overline{T}^1(\{\alpha\}) = \left\{ \begin{array}{cl} \emptyset \in \mathfrak{P}_n, & \mbox{ falls } \alpha \notin \{\alpha_1, \ldots, \alpha_k\}, \\
                                                          \bigcupdot_{i=1}^k Q_i \in \mathfrak{P}_n, & \mbox{ falls } \alpha \in \{\alpha_1, \ldots, \alpha_k\}.
                                \end{array} \right. $$

Zu (ii): Wir beweisen die Aussage durch vollständige Induktion über $r \in \N_+$.

Der Fall $r=1$ gilt nach (i).

Gelte die Behauptung für $(r-1) \in \N_+$.
Dann existieren nach Induktionsvoraussetzung $k \in \N_+$ und paarweise disjunkte $Q_1', \ldots, Q_k' \in \mathfrak{Q}_n$, die in $J$ enthalten sind, mit
$$ \forall_{\rho \in \{1, \ldots, r-1\}} \, T_{\rho} \in {\rm Span}_{\R} \{ \chi_{Q_1'}^J, \ldots, \chi_{Q_k'}^J \}.$$
Wegen (i) gibt es des weiteren paarweise disjunkter Quader $Q_1'',\ldots, Q_l'' \in \mathfrak{Q}_n$ mit $T_r \in {\rm Span}_{\R} \{ \chi_{Q_1''}^J, \ldots, \chi_{Q_l''}^J \}$, also finden wir nach Lemma \ref{FA.4.8.L} offenbar paarweise disjunkte Quader $Q_1, \ldots, Q_r \in \mathfrak{Q}_n$, die in $J$ enthalten sind, derart, daß jeder der Quader $Q_1', \ldots, Q_k', Q_1'', \ldots, Q_l''$ Vereinigung gewisser der $Q_1, \ldots, Q_r$ ist.
Dann gilt offenbar
$$ \forall_{\rho \in \{1, \ldots, r\}} \, T_{\rho} \in {\rm Span}_{\R} \{ \chi_{Q_1}^J, \ldots, \chi_{Q_r}^J \},$$
d.h.\ die Behauptung ist auch für $r$ gezeigt.
\pagebreak

Zu (iii): Seien $Q_1, \ldots, Q_k$ gemäß (ii) gewählt.
Dann ist $\psi \circ (T_1, \ldots, T_r)$ auf jedem der Quader $Q_1, \ldots, Q_k$ konstant.
Wegen $\psi(0) = 0$ ist $\psi \circ (T_1, \ldots, T_r)$ außerdem auf $J \setminus (\bigcupdot_{i=1}^k Q_i)$ gleich Null, also folgt aus (i) $\psi \circ (T_1, \ldots, T_r) \in \mathcal{S}(J)$. 

(iv) und (v) folgen sofort aus (iii). \q    

\begin{Def}[Das Lebesgue-Integral für Treppenfunktionen] \label{FA.4.10}
Es seien $\varphi$ ein Quadermaß auf $\R^n$, $J \in \mathfrak{I}_n$ ein nicht-leeres Intervall von $\R^n$ und $T \in \mathcal{S}(J)$.
Wir definieren dann das \emph{$\varphi$-Integral\index{Integral} von $T$ über $J$} durch
$$ \boxed{\int_J T \, \d \varphi} := \sum_{\alpha \in T(J) \setminus \{0\}} \alpha \cdot \varphi( \overline{T}^1(\{\alpha\}) ) = \sum_{\alpha \in \R \setminus \{0\}} \alpha \cdot \varphi( \overline{T}^1(\{\alpha\}) ). $$
Beachte, daß obige Summen de facto endlich sind und die leere Summe als Null definiert ist.
\end{Def}

\begin{Satz} \label{FA.4.11}
Es seien $\varphi$ ein Quadermaß auf $\R^n$, $J \in \mathfrak{I}_n$ ein nicht-leeres Intervall von $\R^n$ und $T \in \mathcal{S}(J)$.
Seien ferner $k \in \N_+$ und $Q_1, \ldots, Q_k \in \mathfrak{Q}_n$ paarweise disjunkte Teilmengen von $J$ sowie $\alpha_1, \ldots, \alpha_k \in \R$ mit $T = \sum_{i=1}^k \alpha_i \, \chi_{Q_i}^J$.
Dann gilt 
$$ \int_J T \, \d \varphi = \sum_{i=1}^k \alpha_i \, \varphi(Q_i). $$
\end{Satz}

\textit{Beweis.} Mit $I_{\alpha} := \{ \kappa \in \{1, \ldots,k\} \, | \, \alpha_{\kappa} = \alpha \}$ für jedes $\alpha \in T(J) \setminus \{0\}$ gilt $\overline{T}^1(\{\alpha\}) = \bigcupdot_{i \in I_{\alpha}} Q_i$, also nach dem Zusatz zu \ref{FA.4.5} $\varphi(\overline{T}^1(\{\alpha\})) = \sum_{i \in I_{\alpha}} \varphi(Q_i)$ und somit $\alpha \cdot \varphi(\overline{T}^1(\{\alpha\})) = \sum_{i \in I_{\alpha}} \alpha_i \, \varphi(Q_i)$.
Hieraus folgt
$$ \int_J T \, \d \varphi = \sum_{\alpha \in T(J) \setminus \{0\}} \alpha \cdot \varphi( \overline{T}^1(\{\alpha\}) ) = \sum_{\substack{i=1 \\ \alpha_i \ne 0}}^k  \alpha_i \, \varphi(Q_i) = \sum_{i=1}^k \alpha_i \, \varphi(Q_i). $$
\q

\begin{Def} \label{FA.4.fhut}
Seien $M \subset \R^n$ und $f \:M \to \K$ eine Funktion.
Wir definieren dann $\boxed{\hat{f} \: \R^n \to \K}$ durch
$$ \forall_{x \in \R^n} \, \hat{f}(x) := \left\{ \begin{array}{cl} f(x), & \mbox{falls } x \in M, \\ 0, & \mbox{falls } x \in \R^n \setminus M. \end{array} \right. $$
\end{Def}

\begin{HS} \label{FA.4.12} $\,$

\noindent \textbf{Vor.:} Seien $\varphi$ ein Quadermaß auf $\R^n$, $J \in \mathfrak{I}_n$ ein nicht-leeres Intervall von $\R^n$ sowie $T, \widetilde{T} \in \mathcal{S}(J)$ und $\lambda \in \R$.

\noindent \textbf{Beh.:} Es gilt:
\begin{itemize}
\item[(i)] $\D T + \widetilde{T} \in \mathcal{S}(J) \, \wedge \, \int_J (T + \widetilde{T}) \, \d \varphi = \int_J T \, \d \varphi + \int_J \widetilde{T} \, \d \varphi$.
\item[(ii)] $\D \lambda \, T \in \mathcal{S}(J) \, \wedge \, \int_J \lambda \, T \, \d \varphi = \lambda \int_J T \, \d \varphi$.
\item[(iii)] $\D T \le \widetilde{T} \Longrightarrow \int_J T \, \d \varphi \le \int_J \widetilde{T} \, \d \varphi$.
\item[(iv)] $\D \left| \int_J T \, \d \varphi \right| \le \int_J |T| \, \d \varphi$.
\item[(v)] Ist $I \in \mathfrak{I}_n$ eine nicht-leere Teilmenge von $J$, so folgt $T|_I \in \mathcal{S}(I)$, $\chi_I \, T \in \mathcal{S}(J)$ und $\D \int_I T|_I \, \d \varphi = \int_J \chi_I \, T \, \d \varphi$.
\item[(vi)] $\chi_J \, \widehat{T} \in \mathcal{S}(\R^n)$ und $\D \int_{\R^n} \chi_J \, \widehat{T} \, \d \varphi = \int_J T \, \d \varphi$.
\end{itemize}
\end{HS}

\textit{Beweis.} Zu (i), (ii): $\lambda \, T, T + \widetilde{T} \in \mathcal{S}(J)$ ist bereits klar.
Nach \ref{FA.4.8} (ii) existieren paarweise disjunkte Quader $Q_1, \ldots, Q_k \in \mathfrak{Q}_n$, die in $J$ enthalten sind, und $\alpha_1, \ldots, \alpha_k, \tilde{\alpha}_1, \ldots, \tilde{\alpha}_k \in \R$ mit $T = \sum_{i=1}^k \alpha_i \, \chi_{Q_i}^J$ sowie $\widetilde{T} = \sum_{i=1}^k \tilde{\alpha}_i \, \chi_{Q_i}^J$.
Dann gilt $T + \widetilde{T} = \sum_{i=1}^k (\alpha_i + \tilde{\alpha}_i) \, \chi_{Q_i}^J$, $\lambda \, T = \sum_{i=1}^k \lambda \, \alpha_i \, \chi_{Q_i}^J$ und aus \ref{FA.4.11} folgt
\begin{eqnarray*}
\int_J (T + \widetilde{T}) \, \d \varphi & = & \sum_{i=1}^k (\alpha_i + \tilde{\alpha}_i) \, \varphi(Q_i) = \sum_{i=1}^k \alpha_i \, \varphi(Q_i) + \sum_{i=1}^k \tilde{\alpha}_i \, \varphi(Q_i) \\
& = & \int_J T \, \d \varphi + \int_J \widetilde{T} \, \d \varphi
\end{eqnarray*}
sowie
$$ \int_J \lambda \, T \, \d \varphi = \sum_{i=1}^k \lambda \, \alpha_i \, \varphi(Q_i) = \lambda \sum_{i=1}^k \alpha_i \, \varphi(Q_i) = \lambda \int_J T \, \d \varphi. $$

Zu (iii): Aus $T \le \tilde{T}$ folgt $0 \le \tilde{T} - T$, also aus der Definition des Integrales sowie (i), (ii)
$$ 0 \le \int_J (\tilde{T} - T) \, \d \varphi = \int_J \tilde{T} \, \d \varphi - \int_J T \, \d \varphi, $$
d.h.\ $\int_J T \, \d \varphi \le \int_J \tilde{T} \, \d \varphi$.

Zu (iv): Es gilt $-|T| \le T \le |T|$ und $\pm |T| \stackrel{\ref{FA.4.8} (v)}{\in} \mathcal{S}(J)$, also folgt aus (iii), (ii)
$$ - \int_J |T| \, \d \varphi \le \int_J T \, \d \varphi \le \int_J |T| \, \d \varphi $$
und somit $\left| \int_J T \, \d \varphi \right| \le \int_J |T| \, \d \varphi$. 

Zu (v): Sei $I$ wie in (v).
Zunächst gilt
\begin{equation} \label{FA.4.12.1}
\forall_{\alpha \in \R \setminus \{0\}} \, \overline{T_I}^1(\{\alpha\}) = \overline{(\chi_I \, T)}^1(\{\alpha\}).
\end{equation}

{[} Zu (\ref{FA.4.12.2}): Sei $\alpha \in \R \setminus \{0\}$.

,,$\subset$`` Ist $x \in \overline{T_I}^1(\{\alpha\})$, so gilt $x \in I$ und $T(x) = \alpha$, also auch $\chi_I(x) = 1$ und $(\chi_I \, T)(x) = \alpha$, d.h.\ $x \in \overline{(\chi_I \, T)}^1(\{\alpha\})$.

,,$\supset$`` Sei $x \in \overline{(\chi_I \, T)}^1(\{\alpha\})$, d.h.\ $\chi_I(x) \, T(x) = \alpha$.
Wegen $\alpha \ne 0$ gilt $\chi_I(x) \ne 0$, also $\chi_I(x) = 1$ und somit $x \in I$ sowie $T(x) = \alpha$.
Dies bedeutet $x \in \overline{T_I}^1(\{\alpha\})$. {]} 

Wir behaupten
\begin{equation} \label{FA.4.12.2}
\chi_I^J \, T \in \mathcal{S}(J).
\end{equation}

{[} Zu (\ref{FA.4.12.2}): Seien paarweise disjunkte Quader $Q_1, \ldots, Q_k \in \mathfrak{Q}_n$, die in $J$ enthalten sind, und $\alpha_1, \ldots, \alpha_k \in \R$ mit $T = \sum_{i=1}^k \alpha_i \, \chi_{Q_i}^J$ gemäß \ref{FA.4.8} (i) gewählt. 
Dann gilt offenbar $\chi_I^J \, T = \sum_{i=1}^k \alpha_i \, \chi_{Q_i \cap I}^J$, wobei $Q_1 \cap I, \ldots, Q_k \cap I \in \mathfrak{Q}_n$ paarweise disjunkte Teilmengen von $J$ sind, also folgt (\ref{FA.4.12.2}) aus \ref{FA.4.8} (i). {]}

Da nach Definition des Produktes von Funktionen $\chi_I \, T = \chi_I^J \, T$ gilt, ergibt (\ref{FA.4.12.2})
\begin{equation} \label{FA.4.12.S}
\chi_I \, T \in \mathcal{S}(J).
\end{equation}
Aus (\ref{FA.4.12.1}) und (\ref{FA.4.12.S}) folgt für jedes $\alpha \in \R \setminus \{0\}$: $\overline{T_I}^1(\{\alpha\}) \in \mathfrak{P}_n$.
Hieraus und aus $\# T|_I(I) \le \# T(J) < \infty$ -- wegen $I \subset J$, d.h.\ $T|_I(I) \subset T(J)$ -- ergibt sich $T|_I \in \mathcal{S}(I)$ und
$$ \int_I T|_I \, \d \varphi = \sum_{\alpha \in \R \setminus \{0\}} \alpha \, \varphi(\overline{T_I}^1(\{\alpha\})) \stackrel{(\ref{FA.4.12.1})}{=} \sum_{\alpha \in \R \setminus \{0\}} \alpha \, \varphi(\overline{(\chi_I \, T)}^1(\{\alpha\})) = \int_J \chi_I \, T \, \d \varphi. $$

Zu (vi): Offenbar gilt $\widehat{T} \in \mathcal{S}(\R^n)$ sowie
$$ T = \widehat{T}|_J ~~ \mbox{ und } ~~ \widehat{T} = \chi_J \, \widehat{T}. $$
Daher folgt aus (v), indem man dort $I$ durch $J$, $J$ durch $\R^n$ und $T$ durch $\widehat{T}$ ersetzt, $\chi_J \, \widehat{T} \in \mathcal{S}(\R^n)$ sowie $\int_J T \, \d \varphi = \int_J \widehat{T}|_J \, \d \varphi = \int_{\R^n} \chi_J \, \widehat{T} \, \d \varphi.$
\q

\subsection*{Nullmengen} \addcontentsline{toc}{subsection}{Nullmengen}

Um die obige kanonische Definition des Integrales für Treppenfunktionen auf eine möglichst große Menge zu erweitern, führen wir nun den Begriff der \emph{Nullmenge} ein.
Nullmengen werden bei der Entwicklung der Integrationstheorie dadurch eine Rolle spielen, daß sie keine Rolle spielen; das soll heißen, daß die Abänderung des Integranden auf einer Nullmenge den Wert des Integrales nicht ändert. 

\begin{Def}[Nullmengen] \label{FA.4.13}
Seien $\varphi$ ein Quadermaß auf $\R^n$ und $N$ eine Teilmenge von $\R^n$.

$N$ heißt genau dann eine \emph{$\varphi$-Nullmenge}\index{Menge!Null-}\index{Nullmenge}, wenn zu jedem $\varepsilon \in \R_+$ eine Folge $(Q_i)_{i \in \N}$ in $\mathfrak{Q}_n$ mit $N \subset \bigcup_{i=0}^{\infty} Q_i$ und $\sum_{i=0}^{\infty} \varphi(Q_i) < \varepsilon$ existiert.
\end{Def}

\begin{Bem*}
Nullmengen müssen nicht parkettierbar sein.
\end{Bem*}

\begin{Satz} \label{FA.4.13.S}
Es seien $\varphi$ ein Quadermaß auf $\R^n$ und $N, \widetilde{N}$ sowie $N_i$ für jedes $i \in \N$ Teilmengen von $\R^n$.
Dann gilt:
\begin{itemize}
\item[(i)] $N$ $\varphi$-Nullmenge $\wedge$ $\widetilde{N} \subset N$ $\Longrightarrow$ $\widetilde{N}$ $\varphi$-Nullmenge.
\item[(ii)] $\forall_{i \in \N} \, N_i$ $\varphi$-Nullmenge $\Longrightarrow$ $\D \bigcup_{i=0}^{\infty} N_i$ $\varphi$-Nullmenge.
\item[(iii)] $N$ ist genau dann eine $\varphi$-Nullmenge, wenn zu jedem $\varepsilon \in \R_+$ eine Folge offener Quader $(Q_i)_{i \in \N}$ in $\mathfrak{Q}_n$ mit $\D N \subset \bigcup_{i=0}^{\infty} Q_i$ und $\D \sum_{i=0}^{\infty} \varphi(Q_i) < \varepsilon$ existiert.
\end{itemize}
\end{Satz}

\textit{Beweis.} (i) ist trivial.

Zu (ii): Sei $\varepsilon \in \R_+$.
Dann existiert zu jedem $i \in \N$ eine Folge $(Q^i_j)_{j \in \N}$ in $\mathfrak{Q}_n$ mit $N_i \subset \bigcup_{j=0}^{\infty} Q^i_j$ und $\sum_{j=0}^{\infty} \varphi(Q^i_j) < \frac{\varepsilon}{2^{i+1}}$.
Sei dann $(Q_{\iota})_{\iota \in \N}$ eine Umnumerierung von $(Q^i_j)_{(i,j) \in \N^2}$ -- beachte, daß $\N^2$ abzählbar ist.
Dann gilt $\bigcup_{i=0}^{\infty} N_i \subset \bigcup_{\iota = 0}^{\infty} Q_{\iota}$ und
$$ \sum_{\iota = 0}^{\infty} \varphi(Q_{\iota}) = \sum_{i=0}^{\infty} \sum_{j=0}^{\infty} \varphi(Q^i_j) \le \sum_{i=0}^{\infty} \frac{\varepsilon}{2^{i+1}} = \varepsilon. $$

Zu (iii): ,,$\Leftarrow$`` Sei $\varepsilon \in \R_+$.
Dann existiert eine Folge $(\widetilde{Q}_i)_{i \in \N}$ mit $N \subset \bigcup_{i=0}^{\infty} \widetilde{Q}_i$ und $\sum_{i=0}^{\infty} \varphi(\widetilde{Q}_i) < \frac{\varepsilon}{2}$.
Wegen der Regularität von $\varphi$ existiert zu jedem $i \in \N$ ein offener Quader $Q_i \in \mathfrak{Q}_n$ mit $\widetilde{Q}_i \subset Q_i$ und $\varphi(Q_i) \le \varphi(\widetilde{Q}_i) + \frac{\varepsilon}{2^{i+2}}$.
Dann folgt $N \subset \bigcup_{i=0}^{\infty} Q_i$ und
$$ \sum_{i=0}^{\infty} \varphi(Q_i) \le \sum_{i=0}^{\infty} \varphi(\widetilde{Q}_i) + \sum_{i=0}^{\infty} \frac{\varepsilon}{2^{i+2}} < \frac{2 \varepsilon}{2} = \varepsilon. $$

,,$\Rightarrow$`` ist trivial. \q

\begin{Bsp} \label{FA.4.13.B} $\,$
\begin{itemize}
\item[1.)] Ist $P \in \mathfrak{P}_n$ mit $\varphi(P) = 0$, so ist $P$ eine $\varphi$-Nullmenge.
\item[2.)] Jedes entartete Intervall von $\R^n$ ist eine $\mu_n$-Nullmenge.
\item[3.)] Seien $M$ eine Teilmenge von $\R^n$, die keinen Häufungspunkt in $\R^n$ besitzt, und $m \: M \to \R_+$ eine Funktion.
Dann ist $\{p\}$ für jedes $p \in M$ keine $\varphi_m$-Nullmenge.
Im Gegensatz dazu ist $\R^n \setminus M$ eine $\varphi_m$-Nullmenge -- beachte $\forall_{Q \in \mathfrak{Q}_n} \, Q \setminus M = Q \setminus (Q \cap M) \stackrel{\ref{FA.4.2}}{\in}\mathfrak{P}_n$, da $Q \cap M$ die endliche Vereinigung (einelementiger) Quader ist.
\item[4.)] Das \emph{Cantorsche Diskontinuum}\index{Cantorsches Diskontinuum}\index{Menge!Cantor-} $C$ ist eine überabzählbare kompakte $\mu_1$-Nullmenge.
(Sind $C_0 := [0,1]$ und $\forall_{i \in \N_+} \, C_i := \frac{1}{3} C_{i-1} \cup (\frac{2}{3} + \frac{1}{3} C_{i-1})$, d.h.\ $C_i$ ist disjunkte Vereinigung von $2^i$ disjunkten abgeschlossenen Intervallen der Länge $\frac{1}{3^i}$, so gilt $C = \bigcap_{i \in \N} C_i$.)
\end{itemize}
\end{Bsp}

\begin{Def} \label{FA.4.fü} 
Es seien $\varphi$ ein Quadermaß auf $\R^n$ und $M \subset \R^n$.
\begin{itemize}
\item[(i)] Sei $H$ eine einstellige Aussageform mit Einsetzungsklasse $\Omega$.
Wir definieren:
\emph{$H$ gilt $\varphi$-fast überall auf $M$}\index{fast überall} genau dann, wenn $\{x \in M \, | \, \neg H(x) \}$ eine $\varphi$-Nullmenge ist, insbes.\ sind alle Punkte von $M$ mit Ausnahme einer $\varphi$-Nullmenge in $\Omega$ enthalten.
\item[(ii)] Seien $f, g$ auf Teilmengen von $\R^n$ definierte $\K$-wertige Funktionen.
Wir schreiben $\boxed{f =_{\varphi} g \mbox{\emph{ auf $M$}}}$ (bzw.\ $\boxed{f \, {\ne}_{\varphi} \, g \mbox{\emph{ auf $M$}}}$) genau dann, wenn $f=g$ (bzw.\ $f \ne g$) $\varphi$-fast-überall auf $M$ gilt.
$\boxed{f <_{\varphi} g \mbox{\emph{ auf $M$}}}$, $\boxed{f \le_{\varphi} g \mbox{\emph{ auf $M$}}}$, $\boxed{f >_{\varphi} g \mbox{\emph{ auf $M$}}}$ und $\boxed{f \ge_{\varphi} g \mbox{\emph{ auf $M$}}}$ sind im Falle $\K = \R$ analog zu verstehen.
Wir werden bei der Notation häufig auf den Zusatz ,,auf $M$`` verzichten, wenn der Bezug zu $M$ klar ist.
\item[(iii)] Sei $f_i$ für jedes $i \in \N$ eine auf einer Teilmenge von $\R^n$ definierte $\K$-wertige Funktion.
$(f_i)_{i \in \N}$ heißt \emph{$\varphi$-konvergent auf $M$} genau dann, wenn eine $\varphi$-Nullmenge $N \subset \R^n$ existiert derart, daß für $x \in M \setminus N$ gilt:
$f_i$ ist in $x$ für jedes $i \in \N$ definiert und $\D \lim_{i \to \infty} f_i(x) \mbox{ existiert in $\K$}$.
\end{itemize}
\end{Def}

\begin{Lemma}[Lemma A von \textsc{Riesz}] \index{Lemma!von \textsc{Riesz}} \label{FA.4.14} $\,$

\noindent \textbf{Vor.:} Seien $\varphi$ ein Quadermaß auf $\R^n$, $J \in \mathfrak{I}_n$ ein nicht-leeres Intervall von $\R^n$ und $(T_i)_{i \in \N}$ eine Folge in $\mathcal{S}(J)$ mit
\begin{equation} \label{FA.4.14.1}
\forall_{i \in \N} \, T_i \ge T_{i+1} \ge 0.
\end{equation}
Ferner existiere eine $\varphi$-Nullmenge $N \subset \R^n$ derart, daß gilt
\begin{equation} \label{FA.4.14.2}
\forall_{x \in J \setminus N} \, \lim_{i \to \infty} T_i(x) = 0.
\end{equation}

\noindent \textbf{Beh.:} $\D \lim_{i \to \infty} \int_J T_i \, \d \varphi = 0$.
\end{Lemma}

\textit{Beweis.} 1. Fall: $J = \R^n$.
Aus (\ref{FA.4.14.1}) und \ref{FA.4.12} (iii) folgt
$$ \forall_{i \in \N} \, \int_J T_i \, \d \varphi \ge \int_J T_{i+1} \, \d \varphi \ge 0, $$
also existiert $\lim_{i \to \infty} \int_J T_i \, \d \varphi$ in ${[}0, \infty{[}$.
Zu zeigen ist daher
\begin{equation} \label{FA.4.14.3}
\forall_{\varepsilon \in \R_+} \exists_{m \in \N_+} \, \int_J T_m \, \d \varphi < \varepsilon.
\end{equation}

Sei $\varepsilon \in \R_+$.
Als reellwertige Funktion mit nur endlich vielen Werten besitzt $T_0$ ein Maximum.
$C \in \R_+$ sei eine positive reelle Zahl größer als dieses Maximum, d.h.\ nach (\ref{FA.4.14.1})
\begin{equation} \label{FA.4.14.4}
\forall_{i \in \N} \forall_{x \in J} \, T_i(x) < C.
\end{equation}
Da $N$ eine $\varphi$-Nullmenge ist, existiert nach \ref{FA.4.13.S} (iii), in dessen Beweis die Regularität von $\varphi$ entscheidend einging, eine Folge $(Q_i)_{i \in \N}$ offener Quader des $\R^n$ mit
\begin{equation} \label{FA.4.14.5}
N \subset \bigcup_{i=0}^{\infty} Q_i ~~ \mbox{ und } ~~ \sum_{i=0}^{\infty} \varphi(Q_i) < \frac{\varepsilon}{3 C}.
\end{equation}

Wegen $T_0 \in \mathcal{S}(J)$ existieren $r \in \N_+$ und $Q_1^*, \ldots, Q_r^* \in \mathfrak{Q}_n$ derart, daß gilt $\overline{T_0}^1(\R \setminus \{0\}) = \bigcup_{\rho = 1}^r Q_{\rho}^*$, und wir betrachten die kompakte Menge
$$ P := \bigcup_{\rho = 1}^r \overline{Q_{\rho}^*} \in \mathfrak{P}_n. $$
Nach (\ref{FA.4.14.1}) gilt dann
\begin{equation} \label{FA.4.14.6}
\forall_{x \in J \setminus P} \forall_{i \in \N} \, T_i(x) = 0.
\end{equation}

Für jedes $i \in \N$ folgt aus $T_i \in \mathcal{S}(J)$
\begin{eqnarray*}
P_i & := & \overline{T_i}^1 \left( \left[ 0, \frac{\varepsilon}{3 \, (\varphi(P) + 1)} \right[ \right) \cap P = \bigg( \overline{T_i}^1(\{0\}) \cup \overbrace{\overline{T_i}^1 \left( \left] 0, \frac{\varepsilon}{3 \, (\varphi(P) + 1)} \right[ \right)}^{\stackrel{(\ref{FA.4.14.6})}{\subset} P} \bigg) \cap P \\
& = &  \left( \overline{T_i}^1(\{0\}) \cap P \right) \cup \overline{T_i}^1 \left( \left] 0, \frac{\varepsilon}{3 \, (\varphi(P) + 1)} \right[ \right) \\
& = & \big( \underbrace{P}_{\in \mathfrak{P}_n} \setminus \underbrace{\overline{T_i}^1(\R \setminus \{0\})}_{\in \mathfrak{P}_n} \big) \cup \underbrace{\overline{T_i}^1 \left( \left] 0, \frac{\varepsilon}{3 \, (\varphi(P) + 1)} \right[ \right)}_{\in \mathfrak{P}_n} \stackrel{\ref{FA.4.2}}{\in} \mathfrak{P}_n,
\end{eqnarray*}
und es gilt
\begin{equation} \label{FA.4.14.7}
\forall_{x \in J = \R^n} \, \left( x \in P_i \Longleftrightarrow \left( x \in P \, \wedge \, T_i(x) < \frac{\varepsilon}{3 \, (\varphi(P) + 1)} \right) \right).
\end{equation}
Wegen der Regularität von $\varphi \: \mathfrak{P}_n \to \R$ existiert eine Folge $(\widetilde{P}_i)_{i \in \N}$ offener parkettierbarer Teilmengen des $\R^n$ mit
\begin{equation} \label{FA.4.14.8}
\forall_{i \in \N} \, P_i \subset \widetilde{P}_i \, \wedge \, \varphi( \widetilde{P}_i ) \le \varphi(P_i) + \frac{\varepsilon}{3 C} \cdot \frac{1}{2^{i+1}},
\end{equation}
also gilt auch -- da $\widetilde{P}_i = P_i \cupdot (\widetilde{P}_i \setminus P_i)$ und somit $\varphi(\widetilde{P}_i) = \varphi(P_i) + \varphi(\widetilde{P}_i \setminus P_i)$ --
\begin{gather}
\forall_{i \in \N} \, \varphi(\widetilde{P}_i \setminus P_i) \le \frac{\varepsilon}{3 C} \cdot \frac{1}{2^{i+1}}, \nonumber \\
\sum_{i=0}^{\infty} \varphi(\widetilde{P}_i \setminus P_i) \le \frac{\varepsilon}{3 C}. \label{FA.4.14.9}
\end{gather}

Wir behaupten
\begin{equation} \label{FA.4.14.10}
P \subset \left( \bigcup_{i=0}^{\infty} Q_i \right) \cup \left( \bigcup_{i=0}^{\infty} \widetilde{P}_i \right)
\end{equation}

{[} Zu (\ref{FA.4.14.10}): Sei $x \in P$.
Im Falle $x \in N$ gilt $x \in \left( \bigcup_{i=0}^{\infty} Q_i \right) \cup \left( \bigcup_{i=0}^{\infty} \widetilde{P}_i \right)$ nach (\ref{FA.4.14.5}).
Sei daher $x \in P \setminus N \subset \R^n = J$.
Aus (\ref{FA.4.14.2}) folgt dann $\lim_{i \to \infty} T_i(x) = 0$, also existiert $i_0 \in \N$ mit $T_{i_0}(x) < \frac{\varepsilon}{3 \, (\varphi(P) + 1)}$, d.h.\ $x \stackrel{(\ref{FA.4.14.7})}{\in} P_{i_0} \stackrel{(\ref{FA.4.14.8})}{\subset} \widetilde{P}_{i_0}$ und somit $x \in \left( \bigcup_{i=0}^{\infty} Q_i \right) \cup \left( \bigcup_{i=0}^{\infty} \widetilde{P}_i \right)$. {]}

Wegen der Kompaktheit von $P$ und der Offenheit der $Q_i$ sowie der $\widetilde{P}_i$ existiert nach (\ref{FA.4.14.10}) nun $m \in \N$ mit
\begin{equation} \label{FA.4.14.11}
P \subset \bigcup_{i=0}^m \left( Q_i \cup \widetilde{P}_i \right) = \underbrace{\bigcup_{i=0}^m \left( Q_i \cup \left( \widetilde{P}_i \setminus P_i \right) \right)}_{=: P' \in \mathfrak{P}_n} \cup \underbrace{\bigcup_{i=0}^m P_i}_{=: P'' \in \mathfrak{P}_n}.
\end{equation}
Wir setzen
\begin{equation} \label{FA.4.14.12}
T' := C \cdot \underbrace{\chi_{P'}}_{= \chi_{P'}^J} ~~ \mbox{ und } ~~ T'' := \frac{\varepsilon}{3 \, (\varphi(P) + 1)} \cdot \underbrace{\chi_{P''}}_{= \chi_{P''}^J}
\end{equation}
-- beachte $J = \R^n$ --, also gilt $T', T'' \in \mathcal{S}(J)$ mit $T' \ge 0$ sowie $T'' \ge 0$, und behaupten
\begin{equation} \label{FA.4.14.13}
T_m \le T' + T''.
\end{equation}

{[} Zu (\ref{FA.4.14.13}): Sei $x \in J = \R^n$.
Im Falle $x \in J \setminus P$ folgt aus (\ref{FA.4.14.6}), daß gilt $T_m(x) = 0$, also $T_m(x) \le T'(x) + T''(x)$.
Im Falle $x \in P \setminus P'$ ergibt (\ref{FA.4.14.11}), daß gilt $x \in P''$, also existiert $i \in \{0, \ldots, m\}$ mit $x \in P_i$ und es folgt 
$$ T_m(x) \stackrel{(\ref{FA.4.14.1})}{\le} T_i(x) \stackrel{(\ref{FA.4.14.7})}{<} \frac{\varepsilon}{3 \, (\varphi(P) + 1)} = T''(x) \le T'(x) + T''(x). $$
Im Falle $x \in P'$ gilt $T_m(x) \stackrel{(\ref{FA.4.14.4})}{\le} C = T'(x) \le T'(x) + T''(x)$. {]}

Aus (\ref{FA.4.14.11}) und (\ref{FA.4.14.7}) ergibt sich $P'' \subset P$ , also folgt aus der Monotonie von $\varphi \: \mathfrak{P}_n \to \R$, daß gilt 
\begin{equation} \label{FA.4.14.14}
\varphi(P'') \le \varphi(P) < \varphi(P) + 1.
\end{equation}
Somit ergibt sich schließlich
\begin{eqnarray*}
\int_J T_m \, \d \varphi & \stackrel{(\ref{FA.4.14.13})}{\le} & \int_J (T' + T'') \, \d \varphi = \int_J T' \, \d \varphi + \int_J T'' \, \d \varphi \\
& \stackrel{(\ref{FA.4.14.12})}{=} & C \cdot \varphi(P') + \frac{\varepsilon}{3 \, (\varphi(P) + 1)} \cdot \varphi(P'') \\
& \stackrel{(\ref{FA.4.14.14})}{<} & C \cdot \varphi(P') + \frac{\varepsilon}{3 \, (\varphi(P) + 1)} \cdot (\varphi(P) + 1) \\
& \stackrel{(\ref{FA.4.14.11})}{\le} & C \sum_{i=1}^m \left( \varphi(Q_i) + \varphi(\widetilde{P}_i \setminus P_i) \right) + \frac{\varepsilon}{3} \stackrel{(\ref{FA.4.14.5}), (\ref{FA.4.14.9})}{<} \varepsilon,
\end{eqnarray*}
womit (\ref{FA.4.14.3}) bewiesen ist. 

2. Fall: $J \in \mathfrak{I}_n$ beliebig. 
Dann erfüllt die Folge $( \widetilde{T}_i )_{i \in \N}$, definiert durch
$$ \forall_{i \in \N} \, \widetilde{T}_i := \chi_J \, \widehat{T_i}, $$
die Voraussetzungen des Lemmas für $\R^n$ anstelle von $J$, also ergibt der 1. Fall
$$ \lim_{i \to \infty} \int_J T_i \,\d \varphi \stackrel{\ref{FA.4.12} (vi)}{=} \lim_{i \to \infty} \int_{\R^n} \widetilde{T}_i \, \d \varphi = 0, $$
d.h.\ die Behauptung gilt auch im 2. Fall. 
\q

\begin{Satz}[Hauptsatz A von \textsc{Riesz}] \index{Hauptsatz A von \textsc{Riesz}} \label{FA.4.16} $\,$

\noindent \textbf{Vor.:} Seien $\varphi$ ein Quadermaß auf $\R^n$, $J \in \mathfrak{I}_n$ ein nicht-leeres Intervall von $\R^n$ und $(T_i)_{i \in \N}$ sowie $(\widetilde{T}_i)_{i \in \N}$ zwei Folgen in $\mathcal{S}(J)$ mit
\begin{equation} \label{FA.4.16.1}
\forall_{i \in \N} \, \left( T_i \le T_{i+1} \, \wedge \, \widetilde{T}_i \le \widetilde{T}_{i+1} \right).
\end{equation}
Ferner seien $C, \widetilde{C} \in \R$ mit
\begin{equation} \label{FA.4.16.2}
\forall_{i \in \N} \, \left( \int_J T_i \, \d \varphi \le C \, \wedge \, \int_J \widetilde{T}_i \, \d \varphi \le \widetilde{C} \right),
\end{equation}
und es existiere eine $\varphi$-Nullmenge $N$ derart, daß gilt
\begin{equation} \label{FA.4.16.3}
\forall_{x \in J \setminus N} \, \left( \lim_{i \to \infty} T_i(x), \lim_{i \to \infty} \widetilde{T}_i(x) \mbox{ existieren in } \R \, \wedge \, \lim_{i \to \infty} T_i(x) \le \lim_{i \to \infty} \widetilde{T}_i(x) \right). 
\end{equation}

\noindent \textbf{Beh.:} Für die wegen (\ref{FA.4.16.1}), \ref{FA.4.12} (iii) und (\ref{FA.4.16.2}) existierenden folgenden Limites gilt 
$$ \lim_{i \to \infty} \int_J T_i \, \d \varphi \le \lim_{i \to \infty} \int_J \widetilde{T}_i \, \d \varphi. $$
\end{Satz}

\textit{Beweis}. Sei zunächst $k \in \N$ fest gewählt.
Dann gilt nach \ref{FA.4.8} (iv), (v)
\begin{equation} \label{FA.4.16.S}
\forall_{i \in \N} \, T^*_i := (T_k - \widetilde{T_i})^+ = \sup(T_k - \widetilde{T}_i, 0) \in \mathcal{S}(J),
\end{equation}
also wegen (\ref{FA.4.16.1}) auch
\begin{equation} \label{FA.4.16.4}
\forall_{i \in \N} \, T^*_i \ge T^*_{i+1} \ge 0.
\end{equation}
Wir zeigen
\begin{equation} \label{FA.4.16.5}
\forall_{x \in J \setminus N} \, \lim_{i \to \infty} T^*_i(x) = 0.
\end{equation}

{[} Zu (\ref{FA.4.16.5}): Sei $x \in J \setminus N$.
Dann folgt aus (\ref{FA.4.16.1}), (\ref{FA.4.16.3})
$$ T_k(x) \le \lim_{i \to \infty} T_i(x) \le \lim_{i \to \infty} \widetilde{T}_i(x), $$
also 
$$ \lim_{i \to \infty} (T_k(x) - \widetilde{T}_i(x)) \le 0. $$
Die Stetigkeit von $(\id_{\R})^+ = \sup (\ldots,0) \: \R \to \R$ und $\forall_{a \in \R, \ a \le 0} \, (\id_{\R})^+(a) = 0$ ergibt
$$ \lim_{i \to \infty} T^*_i(x) = (\id_{\R})^+ \left( \lim_{i \to \infty} (T_k(x) - \widetilde{T}_i(x)) \right) = 0, $$
d.h.\ es gilt (\ref{FA.4.16.5}). {]}

Aus (\ref{FA.4.16.4}), (\ref{FA.4.16.5}) und Lemma A (d.i.\ \ref{FA.4.14}) folgt
\begin{equation} \label{FA.4.16.6}
\lim_{i \to \infty} \int_J T^*_i \, \d \varphi = 0.
\end{equation}

Nun gilt für $i \in \N$
\begin{gather*}
\int_J T_k \, \d \varphi - \int_J \widetilde{T}_i \, \d \varphi \stackrel{\ref{FA.4.12} (i), (iii), (\ref{FA.4.16.S})}{\le} \int_J T^*_i \, \d \varphi, \\
\int_J T_k \, \d \varphi \le \int_J T^*_i \, \d \varphi + \int_J \widetilde{T}_i \, \d \varphi,
\end{gather*}
also nach (\ref{FA.4.16.6})
$$ \int_J T_k \, \d \varphi \le \lim_{i \to \infty} \int_J \widetilde{T}_i \, \d \varphi. $$

Da die letzte Ungleichung für beliebiges $k \in \N$ gezeigt wurde, folgt die Behauptung durch Grenzwertbildung für $k \to \infty$. \q
\A
Beim Aufbau der Theorie wird das folgende Lemma erst beim Beweis des Grenzwertsatzes von \textsc{Levi} benötigt.
\begin{Lemma}[Lemma B von \textsc{Riesz}] \index{Lemma!von \textsc{Riesz}} \label{FA.4.15} $\,$

\noindent \textbf{Vor.:} Seien $\varphi$ ein Quadermaß auf $\R^n$, $J \in \mathfrak{I}_n$ ein nicht-leeres Intervall von $\R^n$ und $(T_i)_{i \in \N}$ eine Folge in $\mathcal{S}(J)$ mit
\begin{equation} \label{FA.4.15.1}
\forall_{i \in \N} \, T_i \le T_{i+1}.
\end{equation}
Ferner existiere eine Zahl $C \in \R$ derart, daß gilt
\begin{equation} \label{FA.4.15.2}
\forall_{i \in \N} \, \int_J T_i \, \d \varphi \le C.
\end{equation}

\noindent \textbf{Beh.:} $N := \{x \in J \, | \, \lim_{i \to \infty} T_i(x) = \infty\}$ ist eine $\varphi$-Nullmenge.
\end{Lemma}

\textit{Beweis.} Indem wir von $T_i$ zu $T_i - T_0$ übergehen, können wir annehmen, daß gilt $\forall_{i \in \N} \, T_i \ge 0$.
Des weiteren gelte ohne Einschränkung $C > 0$.

Sei nun $\varepsilon \in \R_+$.
Wir setzen 
\begin{equation} \label{FA.4.15.3}
\forall_{i \in \N} \, P_i := \overline{T_i}^1 \left( \left]  \frac{C}{\varepsilon}, \infty \right[ \right) \stackrel{\ref{FA.4.7}, \ref{FA.4.2}}{\in} \mathfrak{P}_n,
\end{equation}
also gilt
\begin{equation} \label{FA.4.15.4}
\forall_{i \in \N} P_i \stackrel{(\ref{FA.4.15.1})}{\subset} P_{i+1} ~~ \mbox{ sowie } ~~ N \subset \bigcup_{i = 0}^{\infty} P_i,
\end{equation}
und behaupten
\begin{equation} \label{FA.4.15.5}
\forall_{i \in \N} \, \varphi(P_i) \le \varepsilon.
\end{equation}

{[} Zu (\ref{FA.4.15.5}): Wegen $\forall_{i \in \N} \, T_i \stackrel{(\ref{FA.4.15.1})}{\ge} T_0 \ge 0$ und (\ref{FA.4.15.3}) folgt $\forall_{i \in \N} \, T_i \ge \frac{C}{\varepsilon} \, \chi_{P_i}^J \in \mathcal{S}(J)$.
Da für jedes $i \in \N$ gilt
$$ C \stackrel{(\ref{FA.4.15.2})}{\ge} \int_J T_i \,\d \varphi \ge \int_J \frac{C}{\varepsilon} \, \chi_{P_i}^J \, \d \varphi = \frac{C}{\varepsilon} \, \varphi(P_i), $$ 
ergibt sich (\ref{FA.4.15.5}). {]}

Aus (\ref{FA.4.15.4}), (\ref{FA.4.15.3}) folgt $P_{i+1} = \overbrace{P_i}^{\in \mathfrak{P}_n} \cupdot \overbrace{(P_{i+1} \setminus P_i)}^{\stackrel{\ref{FA.4.2}}{\in} \mathfrak{P}_n}$ für jedes $i \in \N$.
Daher existieren offenbar eine Folge $(Q_j)_{j \in \N}$ paarweise disjunkter Quader des $\R^n$ und eine streng monoton wachsende Abbildung $h \: \N \to \N$ mit
\begin{equation} \label{FA.4.15.6}
\forall_{i \in \N} \, P_i = \bigcupdot_{j=1}^{h(i)} Q_j, \footnote{Idee: $P_0 = Q_0 \cupdot \ldots \cupdot Q_{h(0)}$ und $\forall_{i \in \N} \, P_{i+1} \setminus P_i = Q_{h(i)+1} \cupdot \ldots \cupdot Q_{h(i+1)}$.}
\end{equation}
d.h.\ nach (\ref{FA.4.15.5})
\begin{equation} \label{FA.4.15.7}
\forall_{i \in \N} \, \sum_{j=1}^{h(i)} \varphi(P_i) = \varphi(P_i) \le \varepsilon.
\end{equation}
Schließlich gilt $N \stackrel{(\ref{FA.4.15.4}), (\ref{FA.4.15.6})}{\subset} \bigcup_{j = 1}^{\infty} Q_j$ und $\sum_{j=1}^{\infty} \varphi(Q_j) \stackrel{(\ref{FA.4.15.7})}{\le} \varepsilon$. \q

\begin{Bem*}
Sind die Voraussetzungen des Lemmas B erfüllt, so besagt die Behauptung $\lim_{i \to \infty} T_i \in \mathcal{S}_{\nearrow}(J,\varphi)$, vgl.\ die Definition in \ref{FA.4.17} (i).
\end{Bem*}

\subsection*{Lebesgue-integrierbare Funktionen} \addcontentsline{toc}{subsection}{Lebesgue-integrierbare Funktionen}

Ausgehend vom Integrale für Treppenfunktionen erweitern wir die Definition nun mehrfach, um die Lebesgue-integrierbarbaren Funktionen zu erhalten.

\begin{Def} \label{FA.4.17}
Seien $\varphi$ ein Quadermaß auf $\R^n$ und $J \in \mathfrak{I}_n$ ein nicht-leeres Intervall von $\R^n$.
\begin{itemize}
\item[(i)] Sei $\boxed{\mathcal{S}_{\nearrow}(J,\varphi)}$ die Menge aller auf Teilmengen von $\R^n$ definierten reellwertigen Funktionen $f$, für welche $\left( (T_i)_{i \in \N}, C, N \right) \in \mathcal{S}(J)^{\N} \times \R \times \mathfrak{P}(\R^n)$ mit den folgenden Eigenschaften existiert:
\begin{itemize}
\item[1.)] $\forall_{i \in \N} \, T_i \le T_{i+1}$.
\item[2.)] $\D \forall_{i \in \N} \, \int_J T_i \, \d \varphi \le C$.
\item[3.)] $N$ ist $\varphi$-Nullmenge und
$$ \forall_{x \in J \setminus N} \, \lim_{i \to \infty} T_i(x) = f(x), $$
insbes.\ gehören alle Punkte von $J$ mit Ausnahme einer $\varphi$-Nullmenge zum Definitionsbereich von $f$.
\end{itemize}
$\left( (T_i)_{i \in \N}, C, N \right)$ heißt dann ein \emph{Levitripel\index{Levitripel} für $f$ auf $J$ bzgl.\ $\varphi$}.\footnote{In \cite{Dom} wird ein Levitripel in unserem Sinne als \emph{Levifolge} bezeichnet. Da letzterer Begriff in der Literatur, vgl.\ z.B.\ \cite{Reck}, aber auch für Folgen, die die Voraussetzungen des Grenzwertsatzes von \textsc{Levi} (s.u.\ \ref{FA.4.25}) erfüllen, verwendet wird, hat sich der Autor zur hier gegebenen Namenswahl entschieden.}
\item[(ii)] Für $f \in \mathcal{S}_{\nearrow}(J,\varphi)$ definieren wir das \emph{$\varphi$-Integral\index{Integral} von $f$ über $J$} durch
$$ \boxed{\int_J f \, \d \varphi} := \lim_{i \to \infty} \int_J T_i \, \d \varphi, $$
wobei $\left( (T_i)_{i \in \N}, C, N \right)$ ein beliebiges Levitripel für $f$ auf $J$ bzgl.\ $\varphi$ sei.

{[} Zur Wohldefiniertheit seien $\left( (T_i)_{i \in \N}, C, N \right)$, $\left( (\widetilde{T}_i)_{i \in \N}, \widetilde{C},\widetilde{N} \right)$ zwei Levitripel für $f$ auf $J$.
Dann ist $N^* := N \cup \widetilde{N}$ eine $\varphi$-Nullmenge, und es gilt
$$ \forall_{x \in J \setminus N^*} \, \lim_{i \to \infty} T_i(x) = \lim_{i \to \infty} \widetilde{T}_i(x) \le \lim_{i \to \infty} \widetilde{T}_i(x), $$
also nach \ref{FA.4.16}
$$ \lim_{i \to \infty} \int_J T_i \, \d \varphi \le \lim_{i \to \infty} \int_J \widetilde{T}_i \, \d \varphi. $$
Indem man die Rollen von $(T_i)_{i \in \N}$ und $(\widetilde{T}_i)_{i \in \N}$ vertauscht, erhält man analog
$$ \lim_{i \to \infty} \int_J \widetilde{T}_i \, \d \varphi \le \lim_{i \to \infty} \int_J T_i \, \d \varphi, $$
also gilt in der letzten Ungleichung sogar Gleichheit. {]}
\end{itemize}

\begin{Bem*} 
Es gilt $\mathcal{S}(J) \subset \mathcal{S}_{\nearrow}(J,\varphi)$, und im Falle $T \in \mathcal{S}(J)$ stimmt die frühere Definition des $\varphi$-In\-te\-gra\-les von $T$ über $J$ (in \ref{FA.4.10}) mit der neu gegebenen Definition überein, da man hier die konstante Folge vom Wert $T$ nutzen kann.
\end{Bem*}
\end{Def}

\begin{Satz} \label{FA.4.18} $\,$

\noindent \textbf{Vor.:} Seien $\varphi$ ein Quadermaß auf $\R^n$, $J \in \mathfrak{I}_n$ ein nicht-leeres Intervall von $\R^n$ sowie $f, \tilde{f} \in \mathcal{S}_{\nearrow}(J,\varphi)$ und $\lambda \in \R$.

\noindent \textbf{Beh.:} Es gilt:
\begin{itemize}
\item[(i)] $\D f + \tilde{f} \in \mathcal{S}_{\nearrow}(J,\varphi) \, \wedge \, \int_J (f + \tilde{f}) \, \d \varphi = \int_J f \, \d \varphi + \int_J \tilde{f} \, \d \varphi$.
\item[(ii)] $\D \lambda \ge 0 \Longrightarrow \left( \lambda \, f \in \mathcal{S}_{\nearrow}(J,\varphi) \, \wedge \, \int_J \lambda \, f \, \d \varphi = \lambda \int_J f \, \d \varphi \right)$.
\item[(iii)] Existiert eine $\varphi$-Nullmenge $N^*$ mit $\forall_{x \in J \setminus N^*} \, f(x) \le \tilde{f}(x)$ -- insbesondere sind $f$ und $\tilde{f}$ in $x$ definiert --, so folgt $\D \int_J f \, \d \varphi \le \int_J \tilde{f} \, \d \varphi$.
\item[(iv)] $\sup(f,\tilde{f}), \, f^+ \in \mathcal{S}_{\nearrow}(J,\varphi)$.\footnote{$\sup(f,\tilde{f})$ kann ggf.\ nur auf $\emptyset$ definiert sein!}
\item[(v)] Sind $f^*$ eine auf einer Teilmenge von $\R^n$ definierte reellwertige Funktion und $N^*$ eine $\varphi$-Nullmenge mit $\forall_{x \in J \setminus N^*} \, f^*(x) = f(x)$ -- insbesondere sind $f^*$ und $f$ in $x$ definiert --, so folgt $f^* \in \mathcal{S}_{\nearrow}(J,\varphi)$ und $\D \int_J f^* \, \d \varphi = \int_J f \, \d \varphi$. 
\item[(vi)] Ist $I \in \mathfrak{I}_n$ eine nicht-leere Teilmenge von $J$, so folgt $f \in \mathcal{S}_{\nearrow}(I,\varphi)$ sowie $\chi_I \, f \in \mathcal{S}_{\nearrow}(J,\varphi)$ und $\D \int_I f \, \d \varphi = \int_J \chi_I \, f \, \d \varphi$.
\item[(vii)] $\chi_J \, \hat{f} \in \mathcal{S}_{\nearrow}(\R^n,\varphi)$ und $\D \int_{\R^n} \chi_J \, \hat{f} \, \d \varphi = \int_J f \, \d \varphi$.
\end{itemize}
\end{Satz}

\textit{Beweis.} Seien $\left( (T_i)_{i \in \N}, C, N \right)$ ein Levitripel für $f$ auf $J$ und $\left( (\widetilde{T}_i)_{i \in \N}, \widetilde{C},\widetilde{N} \right)$ ein Levitripel für $\tilde{f}$ auf $J$, jeweils bzgl.\ $\varphi$.

Zu (i): Wegen \ref{FA.4.12} (i) und \ref{FA.4.13.S} (ii) ist $\left( (T_i + \widetilde{T})_{i \in \N}, C + \widetilde{C}, N \cup \widetilde{N} \right)$ ein Levitripel für $f + \tilde{f}$ auf $J$ bzgl.\ $\varphi$, d.h.\ $f + \tilde{f} \in \mathcal{S}_{\nearrow}(J,\varphi)$, und es gilt
\begin{eqnarray*}
\int_J (f + \tilde{f}) \, \d \varphi & = & \lim_{i \to \infty} \int_J (T_i + \widetilde{T}_i) \, \d \varphi = \lim_{i \to \infty} \int_J T_i \, \d \varphi + \int_J \widetilde{T}_i \, \d \varphi \\
& = & \lim_{i \to \infty} \int_J T_i \, \d \varphi + \lim_{i \to \infty} \int_J \widetilde{T}_i \, \d \varphi = \int_J f \, \d \varphi + \int_J \tilde{f} \, \d \varphi.
\end{eqnarray*}

Zu (ii): Sei $\lambda \ge 0$. 
Wegen \ref{FA.4.12} (ii) ist dann $\left( (\lambda \, T_i)_{i \in \N}, \lambda \, C, N \right)$ ein Levitripel für $\lambda \, f$ auf $J$ bzgl.\ $\varphi$, d.h.\ $\lambda \, f \in \mathcal{S}_{\nearrow}(J,\varphi)$, und es gilt
\begin{eqnarray*}
\int_J \lambda \, f \, \d \varphi & = & \lim_{i \to \infty} \int_J \lambda \, T_i \, \d \varphi = \lim_{i \to \infty} \lambda \int_J T_i \, \d \varphi = \lambda \, \lim_{i \to \infty} \int_J T_i \, \d \varphi \\
& = & \lambda \int_J f \, \d \varphi.
\end{eqnarray*}

Zu (iii): Sei $N^*$ wie in (iii).
Dann ist $N \cup \widetilde{N} \cup N^*$ nach \ref{FA.4.13.S} (ii) eine $\varphi$-Nullmenge, und es gilt (wegen $\forall_{i \in \N} \, T_i \le T_{i+1} \, \wedge \, \widetilde{T}_i \le \widetilde{T}_{i+1}$)
$$ \forall_{x \in J \setminus (N \cup \widetilde{N} \cup N^*)} \, \lim_{i \to \infty} T_i(x) \le \lim_{i \to \infty} \widetilde{T}_i(x). $$
Daher folgt aus \ref{FA.4.16}: $\int_J f \, \d \varphi \le \int_J \tilde{f} \, \d \varphi$.

Zu (iv): Zunächst gilt für jedes $i \in \N$
\begin{gather*}
\sup( T_i, \widetilde{T}_i ) \stackrel{\ref{FA.4.8} (iv)}{\in} \mathcal{S}(J), \\
\sup( \underbrace{T_i}_{\le T_{i+1}}, \underbrace{\widetilde{T}_i}_{\le \widetilde{T}_{i+1}} ) \le \sup( T_{i+1}, \widetilde{T}_{i+1} ).
\end{gather*}

Weiterhin gilt
$$ \forall_{\alpha,\beta,\tilde{\alpha},\tilde{\beta} \in \R} \, \sup\{\alpha + \beta, \tilde{\alpha} + \tilde{\beta}\} \le \sup\{\alpha,\tilde{\alpha}\} + \sup\{\beta,\tilde{\beta}\}, $$
denn ist ohne Einschränkung $\sup\{\alpha + \beta, \tilde{\alpha} + \tilde{\beta}\} = \alpha + \beta$, so folgt sofort obige Ungleichung.
Hieraus ergibt sich für jedes $i \in \N$
\pagebreak
\begin{eqnarray*}
\sup ( T_i, \widetilde{T}_i ) & = & \sup ( T_0 + (T_i - T_0), \widetilde{T}_0 + (\widetilde{T}_i - \widetilde{T}_0) ) \\
& \le & \sup ( T_0, \widetilde{T}_0 ) + \sup ( \underbrace{T_i - T_0}_{\ge 0}, \underbrace{\widetilde{T}_i - \widetilde{T}_0}_{\ge 0} ) \\
& \le & \sup ( T_0, \widetilde{T}_0 ) + (T_i - T_0) + (\widetilde{T}_i - \widetilde{T}_0),
\end{eqnarray*}
also wegen \ref{FA.4.12} (iii) und $\int_J T_i \, \d \varphi \le C, \int_J \widetilde{T}_i \, \d \varphi \le \widetilde{C}$
$$ \int_J \sup ( T_i, \widetilde{T}_i ) \, \d \varphi \le \int_J (\sup ( T_0, \widetilde{T}_0 ) - T_0 - \widetilde{T}_0) \, \d \varphi + C + \widetilde{C} =: C^* \in \R. $$

Schließlich zeigen wir 
\begin{equation} \label{FA.4.18.1}
\forall_{x \in J \setminus (N \cup \widetilde{N})} \, \lim_{i \to \infty} \sup(T_i, \widetilde{T}_i)(x) = \sup(f, \tilde{f})(x),
\end{equation}
d.h.\ insgesamt, daß $\left( \sup(T_i, \widetilde{T}_i)_{i \in \N}, C^*, N \cup \widetilde{N} \right)$ ein Levitripel für $\sup(f, \tilde{f})$ auf $J$ ist -- beachte, daß $N \cup \widetilde{N}$ eine $\varphi$-Nullmenge ist --, m.a.W. $\sup(f, \tilde{f}) \in \mathcal{S}_{\nearrow}(J,\varphi)$.
Daß dann auch $f^+ = \sup(f,0) \in \mathcal{S}_{\nearrow}(J,\varphi)$ gilt, folgt aus $0 \in \mathcal{S}(J) \subset \mathcal{S}_{\nearrow}(J,\varphi)$.

{[} Zu (\ref{FA.4.18.1}): Sei $x \in J \setminus (N \cup \widetilde{N})$.
Ohne Einschränkung gelte $f(x) \ge \tilde{f}(x)$.
Dann gilt wegen $T_i(x) \le f(x), \widetilde{T}_i(x) \le \tilde{f}(x)$ für jedes $i \in \N$ auch
$$ \underbrace{T_i(x)}_{\stackrel{i \to \infty}{\longrightarrow} f(x)} \le \sup\{ T_i(x), \widetilde{T}_i(x) \} \le \sup\{ f(x), \tilde{f}(x) \} = f(x), $$
d.h.\ $ \lim_{i \to \infty} \sup\{ T_i(x), \widetilde{T}_i(x) \} = \sup\{ f(x), \tilde{f}(x) \}$. {]}

Zu (v): Sind $f^*$ und $N^*$ wie in (v), so ist $\left( (T_i)_{i \in \N}, C, N^* \cup N \right)$ wegen
$$  \forall_{x \in J \setminus (N^* \cup N)} \, \lim_{i \to \infty} T_i(x) = f(x) = f^*(x) $$
offenbar ein Levitripel für $f^*$ auf $J$. 

Zu (vi): Sei $I$ wie in (vi).
Dann gilt: 
\begin{equation} \label{FA.4.18.2}
\left( (T_i|_I)_{i \in \N}, C + \int_J (\chi_I \, T_0 - T_0) \, \d \varphi , N \right) \mbox{ ist ein Levitripel für $f$ auf $I$ bzgl.\ $\varphi$,}
\end{equation} 
beachte $\chi_I \, \, T_0 \stackrel{\ref{FA.4.12} (v)}{\in} \mathcal{S}(I)$.

{[} Denn für jedes $i \in \N$ gilt $T_i|_I \le T_{i+1}|_I$ und weiter
\begin{equation} \label{FA.4.18.3}
\chi_I \, T_i = \chi_I \, T_0 + \chi_I \, (\underbrace{T_i - T_0}_{\ge 0}) \le \chi_I \, T_0 + (T_i - T_0) = T_i + (\chi_I \, T_0 - T_0),
\end{equation}
also
\begin{equation} \label{FA.4.18.4}
\begin{array}{lcl}
\D \int_I T_i \, \d \varphi & \stackrel{\ref{FA.4.12} (v)}{=} & \D \int_J \chi_I \, T_i \, \d \varphi \stackrel{(\ref{FA.4.18.3})}{\le} \int_J T_i \, \d \varphi +  \int_J (\chi_I \, T_0 - T_0) \,d \varphi \\
& \le & \D C + \int_J (\chi_I \, T_0 - T_0) \,d \varphi.
\end{array}
\end{equation}
Damit ist (\ref{FA.4.18.2}) offenbar bewiesen. {]}

Für jedes $i \in \N$ gilt auch $\chi_I \, T_i \le \chi_I \, T_{i+1}$, weshalb zusammen mit (\ref{FA.4.18.4}) offenbar folgt:
\begin{equation} \label{FA.4.18.5}
\left( (\chi_I \, T_i|_I)_{i \in \N}, C + \int_J (\chi_I \, T_0 - T_0) \, \d \varphi , N \right) \mbox{ ist Levitripel für $\chi_I \, f$ auf $I$ bzgl.\ $\varphi$.}
\end{equation} 

Schließlich gilt
$$ \int_I f \, \d \varphi \stackrel{(\ref{FA.4.18.2})}{=} \lim_{i \to \infty} \int_I T_i|_I \, \d \varphi \stackrel{\ref{FA.4.12} (v)}{=} \lim_{i \to \infty} \int_J \chi_I \, T_i \, \d \varphi \stackrel{(\ref{FA.4.18.5})}{=} \int_J \chi_I \, f \, \d \varphi, $$
womit zusammen mit (\ref{FA.4.18.2}) und (\ref{FA.4.18.5}) die Aussage (vi) bewiesen ist.

Zu (vii): Wegen \ref{FA.4.12} (vi) ist $\left( (\chi_J \, \widehat{T_i})_{i \in \N}, C, N \right)$ ein Levitripel für $\chi_J \, \hat{f}$ auf $\R^n$ bzgl.\ $\varphi$, d.h.\ $\chi_J \, \hat{f} \in \mathcal{S}_{\nearrow}(\R^n,\varphi)$, und es gilt 
$$ \int_{\R^n} \chi_J \, \hat{f} \, \d \varphi = \lim_{i \to \infty} \int_{\R^n} \chi_J \, \widehat{T_i} \, \d \varphi \stackrel{\ref{FA.4.12} (vi)}{=} \lim_{i \to \infty} \int_J T_i \, \d \varphi = \int_J f \, \d \varphi. $$ 
Damit ist auch (vii) bewiesen.
\q

\begin{Bem*}
Die zu (ii) analoge Aussage mit der Prämisse $\lambda < 0$ ist i.a.\ falsch.
Z.B.\ gilt $f := \sum_{i=0}^{\infty} \chi_{{]}0,2^{-i}{]}} \in \mathcal{S}_{\nearrow}(\R^n,\mu_1)$, aber $-f \notin \mathcal{S}_{\nearrow}(\R^n,\mu_1)$, da jede Treppenfunktion nach unten beschränkt ist und somit auf einem Intervall ${]}0,2^{-i}{]}$ für gewisses $i \in \N$ größer als $-f$ ist, d.h.\ es kann kein Levitripel für $-f$ auf $\R^n$ bzgl.\ $\mu_1$ existieren.
Damit ist gezeigt, daß $\mathcal{S}_{\nearrow}(\R^n,\mu_1)$ kein $\R$-Vektorraum ist.
\end{Bem*}

\begin{Def} \label{FA.4.19}
Seien $\varphi$ ein Quadermaß auf $\R^n$ und $J \in \mathfrak{I}_n$ ein nicht-leeres Intervall von $\R^n$.
\begin{itemize}
\item[(i)] Sei $\boxed{\L_{\R}(J,\varphi)}$ die Menge aller auf Teilmengen von $\R^n$ definierten reellwertigen Funktionen $f$, für welche $(g,h,N) \in \mathcal{S}_{\nearrow}(J,\varphi) \times \mathcal{S}_{\nearrow}(J,\varphi) \times \mathfrak{P}(\R^n)$ mit den folgenden Eigenschaften existiert:
\begin{itemize}
\item[1.)] $N$ ist eine $\varphi$-Nullmenge.
\item[2.)] $\forall_{x \in J \setminus N} \, f(x) = g(x) - h(x)$,

insbes.\ gehören alle Punkte von $J$ mit Ausnahme einer $\varphi$-Nullmenge zum Definitionsbereich von $f$.
\end{itemize}
$(g,h,N)$ heißt dann ein \emph{Riesztripel\index{Riesztripel} für $f$ auf $J$ bzgl.\ $\varphi$}.\footnote{In \cite{Dom} wird ein Rieztripel in unserem Sinne als \emph{Rieszpaar} bezeichnet.}

Die Elemente von $\L_{\R}(J,\varphi)$ nennen wir \emph{(reellwertig) $\varphi$-integrierbar über $J$} (oder \emph{(reellwertig) $\varphi$-sum\-mier\-bar über $J$}).\index{Funktion!integrierbare}
\item[(ii)] Sei $f \in \L_{\R}(J,\varphi)$. 
Wir definieren dann das \emph{$\varphi$-Integral\index{Integral} von $f$ über $J$} durch
$$ \boxed{\int_J f \, \d \varphi} :=  \int_J g \, \d \varphi - \int_J h \, \d \varphi, $$
wobei $(g,h,N)$ ein beliebiges Riesztripel für $f$ auf $J$ bzgl.\ $\varphi$ sei.

{[} Zur Wohldefiniertheit seien $(g,h,N)$ sowie $(\tilde{g},\tilde{h},\widetilde{N})$ zwei Riesztripel für $f$ auf $J$ bzgl.\ $\varphi$.
Dann ist $N^* := N \cup \widetilde{N}$ eine $\varphi$-Nullmenge mit
$$ \forall_{x \in J \setminus N^*} \, ( g,h,\tilde{g},\tilde{h} \mbox{ sind in $x$ definiert und } \underbrace{g(x) - h(x) = \tilde{g}(x) - \tilde{h}(x)}_{\Leftrightarrow  \, (g + \tilde{h})(x) = (\tilde{g} + h)(x)} ). $$
Wegen $g + \tilde{h}, \tilde{g} + h \stackrel{\ref{FA.4.18} (i)}{\in} \mathcal{S}_{\nearrow}(J,\varphi)$ folgt hieraus nach \ref{FA.4.18} (iii) (oder (v))
$$ \int_J (g + \tilde{h}) \, \d \varphi \le \int_J (\tilde{g} + h) \, \d \varphi, ~~ \int_J (\tilde{g} + h) \, \d \varphi \le \int_J (g + \tilde{h}) \, \d \varphi, $$
also mittels \ref{FA.4.18} (i): $\D \int_J g \, \d \varphi - \int_J h \, \d \varphi = \int_J \tilde{g} \, \d \varphi - \int_J \tilde{h} \, \d \varphi$. {]} 
\end{itemize}

\begin{Bem*} 
Es gilt $\mathcal{S}_{\nearrow}(J,\varphi) \subset \L_{\R}(J,\varphi)$, und im Falle $f \in \mathcal{S}_{\nearrow}(J,\varphi)$ stimmt die frühere Definition des $\varphi$-In\-te\-gra\-les von $f$ über $J$ mit der neu gegebenen Definition überein -- beachte, daß $f = f - \, 0$ gilt.
Der nächste Hauptsatz zeigt daher auch, daß $\L_{\R}(J,\varphi)$ die kanonische Erweiterung von $\mathcal{S}_{\nearrow}(J,\varphi)$ zu einem $\R$-Vektorraum und $\int_J \ldots \, \d \varphi$ ein Funktional von $\L_{\R}(J,\varphi)$ ist. 
\end{Bem*}
\end{Def}

\begin{HS} \label{FA.4.20} $\,$

\noindent \textbf{Vor.:} Seien $\varphi$ ein Quadermaß auf $\R^n$, $J \in \mathfrak{I}_n$ ein nicht-leeres Intervall von $\R^n$ sowie $f, \tilde{f} \in \L_{\R}(J,\varphi)$ und $\lambda \in \R$.

\noindent \textbf{Beh.:} Es gilt:
\begin{itemize}
\item[(i)] $\D f + \tilde{f} \in \L_{\R}(J,\varphi) \, \wedge \, \int_J (f + \tilde{f}) \, \d \varphi = \int_J f \, \d \varphi + \int_J \tilde{f} \, \d \varphi$.
\item[(ii)] $\D \lambda \, f \in \L_{\R}(J,\varphi) \, \wedge \, \int_J \lambda \, f \, \d \varphi = \lambda \int_J f \, \d \varphi$.
\item[(iii)] Existiert eine $\varphi$-Nullmenge $N^*$ mit $\forall_{x \in J \setminus N^*} \, f(x) \le \tilde{f}(x)$ -- insbesondere sind $f$ und $\tilde{f}$ in $x$ definiert --, so folgt $\D \int_J f \, \d \varphi \le \int_J \tilde{f} \, \d \varphi$.
\item[(iv)] $f^+, f^-, |f| \in \L_{\R}(J,\varphi)$ und $\D \left| \int_J f \, \d \varphi \right| \le \int_J |f| \, \d \varphi$.
\item[(v)] $\sup(f,\tilde{f}), \, \inf(f,\tilde{f}) \in \L_{\R}(J,\varphi)$. 
\item[(vi)] Sind $f^*$ eine auf einer Teilmenge von $\R^n$ definierte reellwertige Funktion und $N^*$ eine $\varphi$-Nullmenge derart, daß gilt $\forall_{x \in J \setminus N^*} \, f^*(x) = f(x)$ -- d.h.\ insbesondere, daß $f^*$ und $f$ in $x$ definiert sind --, so folgt $f^* \in \L_{\R}(J,\varphi)$ und $\D \int_J f^* \, \d \varphi = \int_J f \, \d \varphi$.
\item[(vii)] Ist $I \in \mathfrak{I}_n$ mit $\emptyset \ne I \subset J$, so folgt $f \in \L_{\R}(I,\varphi)$, $\chi_I \, f \in \L_{\R}(J,\varphi)$, und es gilt $\D \int_I f \, \d \varphi = \int_J \chi_I \, f \, \d \varphi$.
\item[(viii)] $\chi_J \, \hat{f} \in \L_{\R}(\R^n,\varphi)$ und $\D \int_{\R^n} \chi_J \, \hat{f} \, \d \varphi = \int_J f \, \d \varphi$.
\end{itemize}
\end{HS}

\textit{Beweis.} Seien $(g,h,N), (\tilde{g},\tilde{h},\widetilde{N})$ Riesztripel für $f$ bzw.\ $\tilde{f}$ auf $J$ bzgl.\ $\varphi$. 

Zu (i): $(g + \tilde{g}, h + \tilde{h}, N \cup \widetilde{N})$ ist offenbar ein Riesztripel für $f + \tilde{f}$ auf $J$ bzgl.\ $\varphi$, und es gilt
\pagebreak
\begin{eqnarray*}
\int_J (f + \tilde{f}) \, \d \varphi & = & \int_J (g + \tilde{g}) \, \d \varphi - \int_J (h + \tilde{h}) \, \d \varphi \\
& \stackrel{\ref{FA.4.18} (i)}{=} & \int_J g \, \d \varphi + \int_J \tilde{g} \, \d \varphi - \int_J h \, \d \varphi - \int_J \tilde{h} \, \d \varphi \\
& = & \left( \int_J g \, \d \varphi - \int_J h \, \d \varphi \right) + \left( \int_J \tilde{g} \, \d \varphi - \int_J \tilde{h} \, \d \varphi \right) \\
& = & \int_J f \, \d \varphi + \int_J \tilde{f} \, \d \varphi.
\end{eqnarray*}

Zu (ii): Aus \ref{FA.4.18} (ii) folgt die Behauptung im Falle $\lambda \ge 0$, da $(\lambda \, g, \lambda \, h, N)$ ein Riesztripel für $f$ auf $J$ bzgl.\ $\varphi$ ist, und im Falle $\lambda < 0$, da $((-\lambda) \, h, (-\lambda) \, g, N)$ ein Riesztripel für $f$ auf $J$ bzgl.\ $\varphi$ ist.

Zu (iii): Dies beweist man analog zum Nachweis der Wohldefiniertheit in \ref{FA.4.19} (ii).

Zu (iv): Aus 
\begin{equation} \label{FA.4.20.0}
\forall_{\alpha, \beta \in \R} \, (\alpha - \beta)^+ = \sup \{ \alpha - \beta, \underbrace{0}_{= \beta - \beta} \} = \sup \{ \alpha, \beta \} - \beta
\end{equation}
folgt
$$ \forall_{x \in J \setminus N} \, f^+(x) = (g-h)^+(x) = \sup(g,h)(x) - g(x), $$
also ist $(\sup(g,h), g, N)$ nach \ref{FA.4.18} (iv) ein Riesztripel für $f^+$ auf $J$, d.h.\ 
$$ f^+ \in \L_{\R}(J,\varphi) $$
und somit nach (ii), (i) ebenfalls 
$$ f^- = \sup(-f,0) = (-f)^+, \, \, |f| = f^+ + f^- \in \L_{\R}(J,\varphi). $$ 
Ferner gilt wegen $\forall_{x \in J \setminus N} \, -|f|(x) \le f(x) \le |f|(x)$ nach (ii), (iii)
$$ - \int_J |f| \, \d \varphi = \int_J -|f| \, \d \varphi \le \int_J f \, \d \varphi \le \int_J |f| \, \d \varphi, $$
also auch $\left| \int_J f \, \d \varphi \right| \le \int_J |f| \, \d \varphi$. 

Zu (v): Wegen (\ref{FA.4.20.0}) ergibt sich mittels (i), (ii), (iv) nacheinander
\begin{gather*}
\sup(f, \tilde{f}) = f + (\tilde{f} - f)^+ \in \L_{\R}(J,\varphi), \\
\inf(f, \tilde{f}) = - \sup(-f, - \tilde{f}) \in \L_{\R}(J,\varphi).
\end{gather*}

Zu (vi): $(g,h,N^* \cup N)$ ist unter den Voraussetzungen von (vi) ein Riesztripel für $f^*$ auf $J$ bzgl.\ $\varphi$.
Daher folgt die Behauptung.

Zu (vii): Infolge von \ref{FA.4.18} (vi) ist $(g,h,N)$ unter der Voraussetung von (vii) auch ein Riesztripel für $f$ auf $I$ bzgl.\ $\varphi$.
Wegen \ref{FA.4.18} (vi) ist des weiteren $(\chi_I \, g, \chi_I \, h,N)$ ein Riesztripel für $\chi_I \, f$ auf $J$ bzgl.\ $\varphi$.
Daher gilt
\begin{eqnarray*}
\int_I f \, \d \varphi & = & \int_I g \, \d \varphi - \int_I h \, \d \varphi \\
& \stackrel{\ref{FA.4.18} (vi)}{=} & \int_J \chi_I \, g \, \d \varphi - \int_J \chi_I \, h \, \d \varphi = \int_J \chi_I \, f \, \d \varphi.
\end{eqnarray*}

Zu (viii): Nach \ref{FA.4.18} (vii) ist $(\chi_J \, \hat{g}, \chi_J \, \hat{h}, N)$ offenbar ein Riesztripel für $\chi_J \, \hat{f}$, also gilt $\chi_J \, \hat{f} \in \L_{\R}(\R^n,\varphi)$ und
\begin{eqnarray*}
\int_{\R^n} \chi_J \, \hat{f} \, \d \varphi & = & \int_{\R^n} \chi_J \, \hat{g} \, \d \varphi - \int_{\R^n} \chi_J \, \hat{h} \, \d \varphi \\
& \stackrel{\ref{FA.4.18} (vii)}{=} & \int_J g \, \d \varphi - \int_J h \, \d \varphi = \int_J f \, \d \varphi.
\end{eqnarray*}
Damit ist der Hauptsatz vollständig bewiesen. \q

\begin{Def} \label{FA.4.21}
Es seien $\varphi$ ein Quadermaß auf $\R^n$ und $M$ eine beliebige Teilmenge von $\R^n$.
\begin{itemize}
\item[(i)] Wir definieren $\boxed{\L_{\R}(M,\varphi)}$ als die Menge aller auf Teilmengen von $\R^n$ definierten reellwertigen Funktionen $f$, für die $\chi_M \, \hat{f} \in \mathcal{L}_{\R}(\R^n,\varphi)$ gilt.
Die Elemente von $\L_{\R}(M,\varphi)$ nennen wir \emph{(reellwertig) $\varphi$-integrierbar über $M$} (oder \emph{(reellwertig) $\varphi$-sum\-mier\-bar über $M$)}.\index{Funktion!integrierbare}

\begin{Bem*} 
Bezeichnet $D$ den Definitionsbereich von $f \in \mathcal{L}_{\R}(M,\varphi)$, so ist $M \setminus D$ ist eine $\varphi$-Nullmenge, d.h.\ alle Punkte von $M$ mit Ausnahme einer $\varphi$-Nullmenge gehören zu $D$.

{[} Denn aus $\chi_M \, \hat{f} \in \mathcal{L}_{\R}(\R^n,\varphi)$ im Sinne von \ref{FA.4.19} (i) folgt die Existenz von $g,h \in \mathcal{S}_{\nearrow}(\R^n,\varphi)$ mit $\chi_M \, \hat{f} =_{\varphi} g - h$ auf $\R^n$, also existieren offenbar eine $\varphi$-Nullmenge $N$ und zwei Folgen $(T_i)_{i \in \N}, (S_i)_{i \in \N}$ in $\mathcal{S}(\R^n)$ mit 
$$ \forall_{x \in M \setminus N} \, f(x) = \lim_{i \to \infty} (T_i(x) - S_i(x)), $$
d.h.\ insbes.\ $M \setminus N \subset D$, und es ergibt sich $M \setminus D \subset M \cap N \subset N$. {]}
\end{Bem*}
\item[(ii)] Sei $f \in \L_{\R}(M,\varphi)$. 
Wir definieren dann das \emph{$\varphi$-Integral\index{Integral} von $f$ über $M$} durch
$$ \boxed{\int_M f \, \d \varphi} :=  \int_{\R^n} \chi_M \, \hat{f} \, \d \varphi. $$
\end{itemize}

\begin{Bem*} 
Die hier gegebenen Definitionen sind im Falle $\emptyset \ne M \in \mathfrak{I}_n$ mit denen in \ref{FA.4.19} konsistent:
Im folgenden sei $J$ ein nicht-leeres Intervall von $\R^n$.

1.) Sei $f \in \L_{\R}(J,\varphi)$ im Sinne von \ref{FA.4.19} (i).
Dann ergibt \ref{FA.4.20} (viii), daß im Sinne von \ref{FA.4.19} (i) gilt $\chi_J \, \hat{f} \in \mathcal{L}_{\R}(\R^n,\varphi)$, d.h.\ $f \in \mathcal{L}_{\R}(J,\varphi)$ im Sinne von (i).

2.) Sei umgekehrt $f \in \mathcal{L}_{\R}(J,\varphi)$ im Sinne von (i), also $\chi_J \, \hat{f} \in \mathcal{L}_{\R}(\R^n,\varphi)$ im Sinne von \ref{FA.4.19} (i).
Aus \ref{FA.4.20} (vii) folgt dann $\chi_J \, \hat{f} \in \mathcal{L}_{\R}(J,\varphi)$ im Sinne von \ref{FA.4.19} (i).\linebreak
Bezeichnet $D$ den Definitionsbereich von $f$, so ist nach der Bemerkung in (i) $J \setminus D$ des weiteren eine $\varphi$-Nullmenge, und es gilt $J \setminus (J \setminus D) = J \cap D$,
$$ \forall_{x \in J \cap D} \, (\chi_J \, \hat{f})(x) = 1 \, f(x) = f(x). $$
\ref{FA.4.20} (vi) besagt daher $f \in \mathcal{L}_{\R}(J,\varphi)$ im Sinne von \ref{FA.4.19} (i).

3.) Für $f \in \mathcal{L}_{\R}(J,\varphi)$ gilt nach \ref{FA.4.20} (viii): 
$\int_{\R^n} \chi_J \, \hat{f} \, \d \varphi = \int_J f \, \d \varphi$ im Sinne von \ref{FA.4.19} (ii).
\end{Bem*}
\end{Def} 

\begin{HS} \label{FA.4.22} $\,$

\noindent \textbf{Vor.:} Seien $\varphi$ ein Quadermaß auf $\R^n$ und $M$ eine Teilmenge von $\R^n$. 
Ferner seien $f, \tilde{f} \in \L_{\R}(M,\varphi)$ und $\lambda \in \R$.

\noindent \textbf{Beh.:} Es gilt:
\begin{itemize}
\item[(i)] $\D f + \tilde{f} \in \L_{\R}(M,\varphi) \, \wedge \, \int_M (f + \tilde{f}) \, \d \varphi = \int_M f \, \d \varphi + \int_M \tilde{f} \, \d \varphi$.
\item[(ii)] $\lambda \, f \in \L_{\R}(M,\varphi) \, \wedge \, \int_M \lambda \, f \, \d \varphi = \lambda \int_M f \, \d \varphi$.
\item[(iii)] Existiert eine $\varphi$-Nullmenge $N^*$ mit $\forall_{x \in M \setminus N^*} \, f(x) \le \tilde{f}(x)$ -- insbesondere sind $f$ und $\tilde{f}$ in $x$ definiert --, so folgt $\D \int_M f \, \d \varphi \le \int_M \tilde{f} \, \d \varphi$.
\item[(iv)] $f^+, f^-, |f| \in \L_{\R}(M,\varphi)$ und $\D \left| \int_M f \, \d \varphi \right| \le \int_M |f| \, \d \varphi$.
\item[(v)] $\sup(f,\tilde{f}), \, \inf(f,\tilde{f}) \in \L_{\R}(M,\varphi)$. 
\item[(vi)] Sind $f^*$ eine auf einer Teilmenge von $\R^n$ definierte reellwertige Funktion und $N^*$ eine $\varphi$-Nullmenge derart, daß gilt $\forall_{x \in M \setminus N^*} \, f^*(x) = f(x)$ -- d.h.\ insbesondere, daß $f^*$ und $f$ in $x$ definiert sind --, so folgt $f^* \in \L_{\R}(M,\varphi)$ und $\D \int_M f^* \, \d \varphi = \int_M f \, \d \varphi$.
\end{itemize}
\end{HS}

\textit{Beweis.} Bezeichne $D$ bzw.\ $\widetilde{D}$ die Definitionsbereiche von $f$ bzw.\ $g$.
Per definitionem ist $D \cap \widetilde{D}$ der Definitionsbereich von $f+g$. 
Mit $M \setminus D$ und $M \setminus \widetilde{D}$ ist auch $M \setminus (D \cap \widetilde{D}) = (M \setminus D) \cup (M \setminus \widetilde{D})$ eine $\varphi$-Nullmenge, und es gilt für jedes $x \in M \cap (D \cap \widetilde{D})$
\begin{eqnarray*}
(\chi_M \, (\widehat{f+\tilde{f}}))(x) & = & (\widehat{f+\tilde{f}})(x) = (f+\tilde{f})(x) = f(x) + \tilde{f}(x) = \hat{f}(x) + \hat{\tilde{f}}(x) \\
& = & (\chi_M \, \hat{f})(x) + (\chi_M \, \hat{\tilde{f}})(x),
\end{eqnarray*}
sowie $(\chi_M \, (\widehat{f+\tilde{f}}))|_{\R^n \setminus M} = 0 = (\chi_M \, \hat{f} + \chi_M \, \hat{\tilde{f}})|_{\R^n \setminus M}$, also 
$$ \chi_M \, (\widehat{f+\tilde{f}}) =_{\varphi} \chi_M \, \hat{f} + \chi_M \, \hat{\tilde{f}} \mbox{ auf } \R^n. $$
Daher folgt (i) aus \ref{FA.4.20} (i) und \ref{FA.4.20} (vi).

(ii), die erste Aussage von (iv) und (v) folgen aus den analog zu oben einsehbaren $\varphi$-Gleichungen
\begin{gather*}
\chi_M \, (\widehat{\lambda \, f}) =_{\varphi} \lambda \, \chi_M \, \hat{f} \mbox{ auf } \R^n, \\
\chi_M \, \hat{f}^{\pm} =_{\varphi} (\chi_M \, \hat{f})^{\pm} \mbox{ auf } \R^n, ~~ \chi_M \, \widehat{|f|} =_{\varphi} | \chi_M \, \hat{f} | \mbox{ auf } \R^n, \\
\chi_M \, \widehat{\sup (f,\overset{\sim}{f})} =_{\varphi} \chi_M \, \sup(\hat{f}, \hat{\tilde{f}}) \mbox{ auf } \R^n, ~~ \chi_M \, \widehat{\inf(f,\overset{\sim}{f})} =_{\varphi} \chi_M \, \inf(\hat{f}, \hat{\tilde{f}}) \mbox{ auf } \R^n
\end{gather*}
und den entsprechenden Aussagen in \ref{FA.4.20} sowie \ref{FA.4.20} (vi). 
Die zweite Aussage von (iv) ergibt sich dann aus der zweiten Aussage in \ref{FA.4.20} (iv).

Unter den Voraussetzungen von (iii) bzw.\ (vi) gilt $\chi_M \, \hat{f} \le_{\varphi} \chi_M \, \hat{\tilde{f}} \mbox{ auf } \R^n$ bzw.\ $\chi_M \, \hat{f^*} =_{\varphi} \chi_M \, \hat{f} \mbox{ auf } \R^n$.
Also gilt die Behauptung von (iii) bzw.\ (vi) nach \ref{FA.4.20} (iii) und \ref{FA.4.20} (vi) bzw.\ \ref{FA.4.20} (vi). \q

\begin{Bem*}
Ersetzt man in \ref{FA.4.20} (vii) das Intervall $I$ durch eine beliebige Teilmenge $M$ des Intervalles $J$, so folgt aus $f \in \L(J,\varphi)$ i.a.\ weder $f \in \L(M,\varphi)$ noch $\chi_M \, f \in \L(J,\varphi)$ -- wähle z.B.\ $\varphi = \mu_1$, $M$ wie unten in \ref{FA.4.49}, $J = [0,1]$ und $f = 1_{\R^n}$.
\end{Bem*}

\begin{Def} \label{FA.4.23}
Es seien $\varphi$ ein Quadermaß auf $\R^n$ und $M$ eine beliebige Teilmenge von $\R^n$.
\begin{itemize}
\item[(i)] Wir definieren $\boxed{\L_{\C}(M,\varphi)}$ als die Menge aller auf Teilmengen von $\R^n$ definierten komplexwertigen Funktionen $f$ mit ${\rm Re} \, f, {\rm Im} \, f \in \L_{\R}(M,\varphi)$.
Die Elemente von $\L_{\C}(M,\varphi)$ nennen wir \emph{$\varphi$-integrierbar} (oder \emph{$\varphi$-sum\-mier\-bar}) \emph{über $M$}.\index{Funktion!integrierbare}
\item[(ii)] Sei $f \in \L_{\C}(M,\varphi)$. 
Wir definieren dann das \emph{$\varphi$-Integral\index{Integral} von $f$ über $M$} durch
$$ \boxed{\int_M f \, \d \varphi} :=  \int_M {\rm Re} \, f \; \d \varphi +  \i \int_M {\rm Im} \, f \; \d \varphi. $$
\end{itemize}

\begin{Bem*} 
Die hier gegebenen Definitionen sind im Falle $f(M) \subset \R$ mit denen in \ref{FA.4.21} konsistent.
\end{Bem*}
\end{Def}

\begin{HS} \label{FA.4.24} $\,$

\noindent \textbf{Vor.:} Seien $\varphi$ ein Quadermaß auf $\R^n$ und $M$ eine Teilmenge von $\R^n$. 
Ferner seien $f, \tilde{f} \in \L_{\K}(M,\varphi)$ und $\lambda \in \K$.

\noindent \textbf{Beh.:} Es gilt:
\begin{itemize}
\item[(i)] $\D f + \tilde{f} \in \L_{\K}(M,\varphi) \, \wedge \, \int_M (f + \tilde{f}) \, \d \varphi = \int_M f \, \d \varphi + \int_M \tilde{f} \, \d \varphi$.
\item[(ii)] $\D \lambda \, f \in \L_{\K}(M,\varphi) \, \wedge \, \int_M \lambda \, f \, \d \varphi = \lambda \int_M f \, \d \varphi$.
\item[(iii)] $|f| \in \L_{\R}(M,\varphi)$ und $\D \left| \int_M f \, \d \varphi \right| \le \int_M |f| \, \d \varphi$.
\item[(iv)] Sind $f^*$ eine auf einer Teilmenge von $\R^n$ definierte $\K$-wertige Funktion und $N^*$ eine $\varphi$-Nullmenge derart, daß gilt $\forall_{x \in M \setminus N^*} \, f^*(x) = f(x)$ -- d.h.\ insbesondere, daß $f^*$ und $f$ in $x$ definiert sind --, so folgt $f^* \in \L_{\K}(M,\varphi)$ und $\D \int_M f^* \, \d \varphi = \int_M f \, \d \varphi$.
\end{itemize}
\end{HS}

\textit{Beweis.} (i), (ii), (iii) -- im Falle $\K = \R$ -- und (iv) folgen sofort aus \ref{FA.4.23} (i), (ii), (iv) und (vi).
Der Beweis der ersten Aussage von (iii) im Falle $\K=\C$ erfolgt unten auf Seite \pageref{Beweis FA.4.24 (iii)} f.\ mittels des Grenzwertsatzes von \textsc{Lebesgue} \ref{FA.4.34}.

Zur zweiten Aussage von (iii) im Falle $\K = \C$:
Es existiert eine Zahl $t \in \R$ mit
$$ \int_M f \, \d \varphi = \left| \int_M f \, \d \varphi \right| \, \e^{\i t}, $$
also ergibt (ii)
$$ \underbrace{\left| \int_M f \, \d \varphi \right|}_{\in \R_+ \cup \{0\}} = \int_M \e^{-\i t} \, f \, \d \varphi = \int_M {\rm Re} (\e^{-\i t} \, f) \, \d \varphi + \i \, \int_M {\rm Im} (\e^{-\i t} \, f) \, \d \varphi. $$
Hieraus folgt zunächst $\int_M {\rm Im} (\e^{-\i t} \, f) \, \d \varphi = 0$ und sodann wegen der Gültigkeit von (iii) im Falle $\K = \R$ und (ii)
\begin{eqnarray*}
\left| \int_M f \, \d \varphi \right| & = & \left| \int_M {\rm Re} (\e^{-\i t} \, f) \, \d \varphi \right| \le \int_M \underbrace{\left| {\rm Re} (\e^{-\i t} \, f) \right|}_{\le \left| \e^{-\i t} \, f \right|} \, \d \varphi \\
& \le & \left| \e^{-\i t} \right| \, \int_M |f| \, \d \varphi = \int_M |f| \, \d \varphi.
\end{eqnarray*}
\q

\begin{Bem*}
Wir werden die erste Aussage von (iii) im Falle $f \notin \L_{\R}(M,\varphi)$ erst nach ihrem Beweis auf Seite \pageref{Beweis FA.4.24 (iii)} f.\ verwenden!
\end{Bem*}
\subsection*{Grenzwertsätze und deren Folgerungen} \addcontentsline{toc}{subsection}{Grenzwertsätze und deren Folgerungen}

Unter Grenzwertsätzen, die in einem gewissen Sinne das Ziel der Lebesgueschen Integrationstheorie sind, versteht man Sätze, die hinreichende Kriterien dafür angeben, daß die Integration des Grenzwertes einer Funktionenfolge mit der Grenzwertbildung der Integrale vertauschbar ist.
Wir erinnern daran, daß der Definitionsbereich des Grenzwertes einer $\K$-wertigen Funktionenfolge $(f_i)_{i \in \N}$ gleich der Menge der Elemente des Schnittes der Definitionsbereiche der $f_i$, in denen $(f_i)_{i \in \N}$ in $\K$ konvergiert, ist.
Diese Menge kann die leere Menge sein.

Der nun folgende Grenzwertsatz von \textsc{Levi} zeigt auch, daß eine zu \ref{FA.4.17} analoge Konstruktion, die wir durchgeführt haben, um das damals bereits bekannte Integral für Treppenfunktionen zu erweitern, zu keiner größeren Menge integrierbarer Funktionen führt.
Dieser Prozeß ist also abgeschlossen.

\begin{HS}[Grenzwertsatz von \textsc{Levi}] \index{Satz!Grenzwert-!von \textsc{Levi}} \label{FA.4.25} $\,$

\noindent \textbf{Vor.:} Es seien $\varphi$ ein Quadermaß auf $\R^n$, $M \subset \R^n$ und $(f_i)_{i \in \N}$ eine Folge in $\mathcal{L}_{\R}(M,\varphi)$.
Ferner existierten eine $\varphi$-Nullmenge $N \subset \R^n$ sowie eine reelle Zahl $C \in \R$ mit
\begin{gather}
\forall_{i \in \N} \forall_{x \in M \setminus N} \, f_i(x) \le f_{i+1}(x), \label{FA.4.25.1} \\
\forall_{i \in \N} \, \int_{M} f_i \, \d \varphi \le C. \label{FA.4.25.2}
\end{gather}

\noindent \textbf{Beh.:} $(f_i)_{i \in \N}$ ist $\varphi$-konvergent auf $M$, $\D \lim_{i \to \infty} f_i \in \mathcal{L}_{\R}(M,\varphi)$ und 
$$ \int_M \left( \lim_{i \to \infty} f_i \right) \, \d \varphi = \lim_{i \to \infty} \int_M f_i \, \d \varphi. $$
\end{HS}

Wir bereiten den Beweis des Hauptsatzes durch zwei Lemmata vor.

\begin{Lemma} \label{FA.4.26} $\,$

\noindent \textbf{Vor.:} Seien $\varphi$ ein Quadermaß auf $\R^n$, $J \in \mathfrak{I}_n$ mit $J \ne \emptyset$ und $f \in \mathcal{L}_{\R}(J,\varphi)$.

\noindent \textbf{Beh.:} Zu jedem $\varepsilon \in \R_+$ existieren $g,h \in \mathcal{S}_{\nearrow}(J,\varphi)$ derart, daß gilt $f =_{\varphi} g - h$ auf $J$,\linebreak $h \ge_{\varphi} 0$ auf $J$ und $\D \int_J h \, \d \varphi \le \varepsilon$.
\end{Lemma}

\textit{Beweis.} Aus $f \in \L_{\R}(J,\varphi)$ folgt zunächst die Existenz von $\tilde{g},\tilde{h} \in \mathcal{S}_{\nearrow}(J,\varphi)$ mit $f =_{\varphi} \tilde{g} - \tilde{h}$.
Also existieren auch monoton wachsende Folgen $(S_i)_{i \in \N}$ bzw.\ $(T_i)_{i \in \N}$ von Treppenfunktionen auf $J$ mit beschränkten $\varphi$-Integralfolgen und $\lim_{i \to \infty} S_i =_{\varphi} \tilde{g}$ sowie $\lim_{i \to \infty} T_i =_{\varphi} \tilde{h}$.
Dann gilt per definitionem 
$$ \int_J \tilde{h} \, \d \varphi = \lim_{i \to \infty} \int_J T_i \, \d \varphi. $$
Ist daher $\varepsilon \in \R_+$ vorgegeben, so existiert $k \in \N_+$ mit
\begin{equation} \label{FA.4.26.1}
\int_J ( \tilde{h} - T_k ) \, \d \varphi = \int_J \tilde{h} \, \d \varphi - \int_J T_k  \, \d \varphi \le \varepsilon.
\end{equation}
Dann sind  auch $(S_i - T_k)_{i \in \N}$ und $(T_i - T_k)_{i \in \N}$ monoton wachsende Folgen von Treppenfunktionen auf $J$ mit beschränkten $\varphi$-Integralfolgen, und es folgt
$$ \lim_{i \to \infty} (S_i - T_k) = \tilde{g} - T_k =: g ~~ \wedge ~~ \lim_{i \to \infty} (T_i - T_k) = \tilde{h} - T_k =: h, $$
also $g,h \in \mathcal{S}_{\nearrow}(J,\varphi)$ sowie $f =_{\varphi} \tilde{g} - \tilde{h} = g - h$.
Wegen $\forall_{i \in \N, \, i \ge k} \, T_i \ge T_k$ gilt $h \ge_{\varphi} 0$.
Des weiteren ergibt (\ref{FA.4.26.1}): $\int_J h \, \d \varphi = \int_J ( \tilde{h} - T_k) \, \d \varphi \le \varepsilon$. \q

\begin{Lemma} \label{FA.4.27}
Es gilt der zu Hauptsatz \ref{FA.4.25} analoge Satz, den man erhält, indem man in seiner Voraussetzung und seiner Behauptung jeweils $\L_{\R}(M,\varphi)$ durch $\mathcal{S}_{\nearrow}(J,\varphi)$ ersetzt, wobei $J \in \mathfrak{I}_n$ mit $J \ne \emptyset$ sei.
\end{Lemma}

\textit{Beweis.} Sei die Voraussetzung von \ref{FA.4.25} mit $\mathcal{S}_{\nearrow}(J,\varphi)$ anstelle von $\L_{\R}(M,\varphi)$ erfüllt.
Für jedes $i \in \N$ folgt aus $f_i \in \mathcal{S}_{\nearrow}(J,\varphi)$ die Existenz einer Folge $(T_{i,j})_{j \in \N}$ in $\mathcal{S}(J)$ mit
\begin{gather}
\forall_{j \in \N} \, T_{i,j} \le T_{i,j+1}, \label{FA.4.25.3} \\
\lim_{j \to \infty} T_{i,j} =_{\varphi} f_i. \label{FA.4.25.4}
\end{gather}
Dann wird nach \ref{FA.4.8} (iv) durch
$$ \forall_{k \in \N} \, T_k := \sup( \underbrace{T_{0,0}, \ldots, T_{0,k}}_{\stackrel{(\ref{FA.4.25.3}), (\ref{FA.4.25.4})}{\le_{\varphi}} f_0}, \underbrace{T_{1,0}, \ldots, T_{1,k}}_{\stackrel{(\ref{FA.4.25.3}), (\ref{FA.4.25.4})}{\le_{\varphi}} f_1}, \ldots, \underbrace{T_{k,0}, \ldots, T_{k,k}}_{\stackrel{(\ref{FA.4.25.3}), (\ref{FA.4.25.4})}{\le_{\varphi}} f_k} ) $$
eine Folge $(T_k)_{k \in \N}$ in $\mathcal{S}(J)$ definiert, und es gilt:
\begin{gather}
\forall_{k \in \N} \, T_k \le_{\varphi} f_k, \label{FA.4.25.5} \\
\forall_{k \in \N} \, T_k \le_{\varphi} T_{k+1}, \label{FA.4.25.6} \\
\exists_{C \in \R} \forall_{k \in \N} \, \int_{J} T_k \, \d \varphi \le C. \label{FA.4.25.7}
\end{gather}

{[} (\ref{FA.4.25.5}) folgt aus (\ref{FA.4.25.3}), (\ref{FA.4.25.4}) und (\ref{FA.4.25.1}), die Aussage (\ref{FA.4.25.6}) ergibt sich sofort aus der Definition, und (\ref{FA.4.25.7}) folgt aus (\ref{FA.4.25.5}) und (\ref{FA.4.25.2}). {]}

Nach (\ref{FA.4.25.6}), (\ref{FA.4.25.7}) und Lemma B von \textsc{Riesz} (d.i. \ref{FA.4.15}) ist daher
$$ N^* := \left\{ x \in J \, | \, \lim_{k \to \infty} T_k(x) = \infty \right\} $$
eine $\varphi$-Nullmenge, d.h.\ $(T_k)_{k \in \N}$ ist $\varphi$-konvergent auf $J$ und
\begin{gather}
\lim_{k \to \infty} T_k \in \mathcal{S}_{\nearrow}(J,\varphi), \label{FA.4.25.8} \\
\int_{J} \left( \lim_{k \to \infty} T_k \right) \, \d \varphi = \lim_{k \to \infty} \int_{J} T_k \, \d \varphi. \label{FA.4.25.9}
\end{gather}

Wegen (\ref{FA.4.25.1}), (\ref{FA.4.25.4}) und \ref{FA.4.13.S} (ii) existiert eine $\varphi$-Nullmenge $\widetilde{N}$ mit
\begin{gather}
\forall_{i \in \N} \forall_{x \in J \setminus \widetilde{N}} \, f_i(x) \le f_{i+1}(x), \label{FA.4.25.10} \\
\forall_{i \in \N} \forall_{x \in J \setminus \widetilde{N}} \, \lim_{j \to \infty} T_{i,j}(x) = f_i(x). \label{FA.4.25.11}
\end{gather}
Wir behaupten:
\begin{equation} \label{FA.4.25.12}
\forall_{k \in \N} \forall_{x \in J \setminus (N^* \cup \widetilde{N})} \, T_k(x) \le f_k(x) \le \lim_{j \to \infty} T_j(x).
\end{equation}

{[} Zu (\ref{FA.4.25.12}): Seien $k \in \N$ und $x \in J \setminus (N^* \cup \widetilde{N})$.
Dann gilt
\begin{eqnarray*}
T_k(x) & = & \sup\{ \underbrace{T_{0,0}(x), \ldots, T_{0,k}(x)}_{\stackrel{(\ref{FA.4.25.3}), (\ref{FA.4.25.11})}{\le} f_0(x)}, \ldots, \underbrace{T_{k,0}(x), \ldots, T_{k,k}(x)}_{\stackrel{(\ref{FA.4.25.3}), (\ref{FA.4.25.11})}{\le} f_k(x)} \} \\
& \stackrel{(\ref{FA.4.25.10})}{\le} & f_k(x) \stackrel{(\ref{FA.4.25.11})}{=} \lim_{j \to \infty} T_{k,j}(x) \le \lim_{j \to \infty} T_j(x),
\end{eqnarray*}
beachte bei der letzten Ungleichung, daß für $j \ge k$ gilt $T_{k,j}(x) \le T_j(x)$. {]}

Aus (\ref{FA.4.25.12}) folgt
$$ \forall_{x \in J \setminus (N^* \cup \widetilde{N})} \, \lim_{k \to \infty} f_k(x) = \lim_{k \to \infty} T_k(x), $$
also ist $(f_k)_{k \in \N}$ $\varphi$-konvergent auf $J$ und $\lim_{k \to \infty} f_k =_{\varphi} \lim_{k \to \infty} T_k$.
Daher ergibt (\ref{FA.4.25.8})
\begin{equation} \label{FA.4.25.13}
\lim_{k \to \infty} f_k \in \mathcal{S}_{\nearrow}(J,\varphi) ~~ \wedge ~~ \int_{J} \left( \lim_{k \to \infty} f_k \right) \, \d \varphi = \int_{J} \left( \lim_{k \to \infty} T_k \right) \, \d \varphi.
\end{equation}
Des weiteren folgt aus (\ref{FA.4.25.12}) für jedes $k \in \N$
$$ \int_{J} T_k \, \d \varphi \le \int_{J} f_k \, \d \varphi \le \int_{J} \left( \lim_{j \to \infty} T_j \right) \, \d \varphi, $$
wobei die linke Seite nach (\ref{FA.4.25.9}) für $k \to \infty$ gegen $\int_{J} (\lim_{k \to \infty} T_k) \, \d \varphi$ konvergiert, also
$$ \int_{J} \left( \lim_{k \to \infty} f_k \right) \, \d \varphi \stackrel{(\ref{FA.4.25.13})}{=} \int_{J} \left( \lim_{k \to \infty} T_k \right) \, \d \varphi = \lim_{k \to \infty} \int_{J} f_k \, \d \varphi. $$
Damit ist das Lemma bewiesen. \q
\A
\textit{Beweis des Grenzwertsatzes von \textsc{Levi}.} 1. Fall: $M = \R^n$.
Sei also $(f_i)_{i \in \N}$ wie in der Voraussetzung von \ref{FA.4.25} (mit $M=\R^n$).
Dann wird durch
$$ \tilde{f}_0 := f_0 \mbox{ und } \forall_{i \in \N_+} \, \tilde{f}_i := f_i - f_{i-1} $$
eine Folge $(\tilde{f}_i)_{i \in \N}$ in $\L_{\R}(\R^n,\varphi)$ mit
$$ \forall_{i \in \N} \, f_i =_{\varphi} \sum_{j=0}^i \tilde{f}_j $$
definiert.
Nach Lemma \ref{FA.4.26} existieren Folgen $(\tilde{g}_i)_{i \in \N}$, $(\tilde{h}_i)_{i \in \N}$ in $\mathcal{S}_{\nearrow}(\R^n,\varphi)$ mit
\begin{gather}
\forall_{i \in \N} \, \tilde{f}_i =_{\varphi} \tilde{g}_i - \tilde{h}_i, \label{FA.4.25.14} \\
\forall_{i \in \N} \, \tilde{h}_i \, {\ge}_{\varphi} \, 0, \label{FA.4.25.15} \\
\forall_{i \in \N} \, \int_{\R^n} \tilde{h}_i \, \d \varphi \le \frac{1}{2^i}. \label{FA.4.25.16}
\end{gather}

Wir setzen für jedes $i \in \N$
$$ \forall_{i \in \N} \, g_i := \sum_{j=0}^i \tilde{g}_j, \, h_i := \sum_{j=0}^i \tilde{h}_j \stackrel{\ref{FA.4.18} (i)}{\in} \mathcal{S}_{\nearrow}(\R^n,\varphi) $$
und behaupten
\begin{gather}
\forall_{i \in \N} \, f_i =_{\varphi} g_i - h_i, \label{FA.4.25.17} \\
\forall_{i \in \N} \, (g_i \le_{\varphi} g_{i+1} ~~ \wedge ~~ h_i \le_{\varphi} h_{i+1}), \label{FA.4.25.18} \\
\forall_{i \in \N} \, \left( \int_{\R^n} g_i \, \d \varphi \le C + 2 ~~ \wedge ~~ \int_{\R^n} h_i \, \d \varphi \le 2 \right), \label{FA.4.25.19}
\end{gather}
wobei $C$ wie in der Voraussetzung von \ref{FA.4.25} sei.

{[} Sei $i \in \N$.

Zu (\ref{FA.4.25.17}): $f_i =_{\varphi} \sum_{j=0}^i \tilde{f}_j \stackrel{(\ref{FA.4.25.14})}{=_{\varphi}} \sum_{j=0}^i \tilde{g}_j - \sum_{j=0}^i \tilde{h}_j = g_i - h_i$.

Zu (\ref{FA.4.25.18}): $h_{i+1} - h_i =_{\varphi} \tilde{h}_{i+1} \stackrel{(\ref{FA.4.25.15})}{{\ge}_{\varphi}} 0$, also auch
\begin{eqnarray*}
g_{i+1} - g_i & =_{\varphi} & \tilde{g}_{i+1} \stackrel{(\ref{FA.4.25.14})}{=_{\varphi}} \tilde{f}_{i+1} + \tilde{h}_{i+1} \\
& = & \underbrace{f_{i+1} - f_i}_{\stackrel{(\ref{FA.4.25.1})}{{\ge}_{\varphi}} 0}  + \underbrace{\tilde{h}_{i+1}}_{{\ge}_{\varphi} 0} {\ge}_{\varphi} \, 0.
\end{eqnarray*}

Zu (\ref{FA.4.25.19}): $\int_{\R^n} h_i \, \d \varphi = \sum_{j=0}^i \int_{\R^n} \tilde{h}_j \, \d \varphi \stackrel{(\ref{FA.4.25.16})}{{\le}_{\varphi}} \sum_{j=0}^i \frac{1}{2^i} \le 2$, also wegen (\ref{FA.4.25.2}) auch $\int_{\R^n} g_i \, \d \varphi \stackrel{(\ref{FA.4.25.17})}{=} \int_{\R^n} f_i \, \d \varphi + \int_{\R^n} h_i \, \d \varphi \le C + 2$. {]}

Aus (\ref{FA.4.25.18}), (\ref{FA.4.25.19}) folgt mittels Lemma \ref{FA.4.27} die $\varphi$-Konvergenz von $(g_i)_{i \in \N}$, $(h_i)_{i \in \N}$ auf $\R^n$ und
\begin{gather}
\lim_{i \to \infty} g_i, \lim_{i \to \infty} h_i \in \mathcal{S}_{\nearrow}(\R^n,\varphi), \label{FA.4.25.20} \\
\int_{\R^n} \left( \lim_{i \to \infty} g_i \right) \, \d \varphi = \lim_{i \to \infty} \int_{\R^n} g_i \, \d \varphi ~~ \wedge ~~ \int_{\R^n} \left( \lim_{i \to \infty} h_i \right) \, \d \varphi = \lim_{i \to \infty} \int_{\R^n} h_i \, \d \varphi. \label{FA.4.25.21}
\end{gather}
(\ref{FA.4.25.17}) und (\ref{FA.4.25.20}) ergeben nun, daß auch $(f_i)_{i \in \N}$ $\varphi$-konvergent ist sowie 
$$ \lim_{i \to \infty} f_i =_{\varphi} \lim_{i \to \infty} g_i - \lim_{i \to \infty} h_i, $$
also nach (\ref{FA.4.25.20})
$$ \lim_{i \to \infty} f_i \in \L_{\R}(\R^n,\varphi) $$
und
\begin{eqnarray*}
\int_{\R^n} \left( \lim_{i \to \infty} f_i \right) \, \d \varphi & = & \int_{\R^n} \left( \lim_{i \to \infty} g_i \right) \, \d \varphi - \int_{\R^n} \left( \lim_{i \to \infty} h_i \right) \, \d \varphi \\
& \stackrel{(\ref{FA.4.25.21})}{=} & \lim_{i \to \infty} \int_{\R^n} g_i \, \d \varphi - \lim_{i \to \infty} \int_{\R^n} h_i \, \d \varphi = \lim_{i \to \infty} \int_{\R^n} (g_i - h_i) \, \d \varphi \\
& \stackrel{(\ref{FA.4.25.17})}{=} & \lim_{i \to \infty} \int_{\R^n} f_i \, \d \varphi.
\end{eqnarray*}

2. Fall: $M$ beliebige Teilmenge von $\R^n$.
Sei nun $(f_i)_{i \in \N}$ wie in der Voraussetzung von \ref{FA.4.25}.
Dann wird durch
\begin{equation} \label{FA.4.25.S}
\forall_{i \in \N} \, \tilde{f}_i := \chi_M \, \widehat{f_i}
\end{equation}
eine Folge $(\tilde{f}_i)_{i \in \N}$ in $\L_{\R}(\R^n,\varphi)$ definiert, die die Voraussetzung des Satzes von \textsc{Levi} im Falle $M = \R^n$ erfüllt.
Daher folgt aus dem 1. Fall die $\varphi$-Konvergenz von $(\tilde{f}_i)_{i \in \N}$ auf $\R^n$ und
\begin{gather}
\lim_{i \to \infty} \tilde{f}_i \in \mathcal{L}_{\R}(\R^n,\varphi), \label{FA.4.25.22} \\
\int_{\R^n} \left( \lim_{i \to \infty} \tilde{f}_i \right) \, \d \varphi = \lim_{i \to \infty} \int_{\R^n} \tilde{f}_i \, \d \varphi. \label{FA.4.25.23}
\end{gather}
Bezeichne $D$ den Definitionsbereich von $\lim_{i \to \infty} \tilde{f}_i$.
Dann ist also $(\R^n \setminus D) \cup N$, wobei $N$ wie in der Voraussetzung von \ref{FA.4.25} sei, nach (\ref{FA.4.25.22}) eine $\varphi$-Nullmenge, und es gilt für jedes $x \in \R^n \setminus ((\R^n \setminus D) \cup N) = D \setminus N$
\begin{gather}
x \in M \Longrightarrow \lim_{i \to \infty} \tilde{f}_i(x) \stackrel{(\ref{FA.4.25.S})}{=} \lim_{i \to \infty} \widehat{f_i}(x) \stackrel{x \notin N}{=} \lim_{i \to \infty} f_i(x), \label{FA.4.25.24}
\end{gather}
\begin{eqnarray*}
\left( \lim_{i \to \infty} \tilde{f}_i \right)(x) & \stackrel{(\ref{FA.4.25.S})}{=} & \lim_{i \to \infty} \chi_M(x) \,\widehat{f_i}(x) \stackrel{(\ref{FA.4.25.24})}{=} \chi_M(x) \left( \widehat{\lim_{i \to \infty} f_i(x)} \right) \\
& = & \left( \chi_M \left( \widehat{\lim_{i \to \infty} f_i} \right) \right)(x),
\end{eqnarray*}
d.h.\ $(f_i)_{i \in \N}$ ist $\varphi$-konvergent auf $M$ und
\begin{equation} \label{FA.4.25.25}
\lim_{i \to \infty} \tilde{f}_i =_{\varphi} \chi_M \left( \widehat{\lim_{i \to \infty} f_i} \right).
\end{equation}
Daher folgt aus (\ref{FA.4.25.22}) und \ref{FA.4.20} (vi)
$$ \chi_M \left( \widehat{\lim_{i \to \infty} f_i} \right) \in \L_{\R}(\R^n,\varphi), $$
d.h.\ per definitionem $\lim_{i \to \infty} f_i \in \L_{\R}(M,\varphi)$ und
$$ \int_{M} \left( \lim_{i \to \infty} f_i \right) \, \d \varphi \stackrel{(\ref{FA.4.25.25})}{=} \int_{\R^n} \left( \lim_{i \to \infty} \tilde{f}_i \right) \, \d \varphi \stackrel{(\ref{FA.4.25.23})}{=} \lim_{i \to \infty} \int_{\R^n} \tilde{f}_i \, \d \varphi \stackrel{(\ref{FA.4.25.S})}{=} \lim_{i \to \infty} \int_M f_i \, \d \varphi, $$
womit der Grenzwertsatz von \textsc{Levi} vollständig bewiesen ist. \q

\pagebreak
Wir erinnern daran, daß das Supremum bzw.\ Infimum einer reellen Funktionenfolge $(f_i)_{i \in N}$ auf der Menge der Elemente des Schnittes der Definitionsbereiche der $f_i$, in denen das Supremum bzw.\ Infimum der $f_i$ in $\R$ existiert, definiert ist.
Diese Menge kann die leere Menge sein.

\begin{Satz} \label{FA.4.32} $\,$

\noindent \textbf{Vor.:} Seien $\varphi$ ein Quadermaß auf $\R^n$, $M$ eine Teilmenge von $\R^n$ und $(f_i)_{i \in \N}$ eine Folge in $\L_{\R}(M,\varphi)$.

\noindent \textbf{Beh.:} Dann gilt:
\begin{itemize}
\item[(i)] Existieren eine Funktion $f \in \L_{\R}(M,\varphi)$ sowie eine $\varphi$-Nullmenge $N$ derart, daß gilt $\forall_{x \in M \setminus N} \forall_{i \in \N} \, f_i(x) \le f(x)$, so folgt $\sup( f_0, f_1, \ldots ) \in \L_{\R}(M,\varphi)$.
\item[(ii)] Existieren eine Funktion $f \in \L_{\R}(M,\varphi)$ sowie eine $\varphi$-Nullmenge $N$ derart, daß gilt $\forall_{x \in M \setminus N} \forall_{i \in \N} \, f_i(x) \ge f(x)$, so folgt $\inf( f_0, f_1, \ldots ) \in \L_{\R}(M,\varphi)$.
\end{itemize}
\end{Satz}

\textit{Beweis.} Zu (i): Für jedes $i \in \N$ definieren wir
$$ g_i := \sup(f_0, \ldots, f_i) \in \L_{\R}(M,\varphi), $$
beachte \ref{FA.4.22} (v) und $g_{i+1} = \sup(g_i,f_{i+1})$.
Dann gilt für jedes $i \in \N$
$$ \forall_{x \in M \setminus N} \, g_i(x) \le g_{i+1}(x) \le f(x), $$
also nach \ref{FA.4.22} (iii)
$$ \int_M g_i \, \d \varphi \le \int_M f \, \d \varphi. $$
Daher folgt aus dem Grenzwertsatz von \textsc{Levi} \ref{FA.4.25} die $\varphi$-Konvergenz von $(g_i)_{i \in \N}$ auf $M$ und $\lim_{i \to \infty} g_i \in \L_{\R}(M,\varphi)$.
Zum Nachweis von (i) genügt es daher zu zeigen, daß gilt
\begin{equation} \label{FA.4.32.S}
\forall_{x \in M \setminus N} \, \lim_{i \to \infty} g_i(x) = \sup \{ f_i(x) \, | \, i \in \N \}.
\end{equation}

{[} Zu (\ref{FA.4.32.S}): Sei $x \in M \setminus N$, d.h.\ $\forall_{i \in \N} \, f_i(x) \le f_{i+1}(x)$.
Dann folgt für jedes $i \in \N$
$$ g_i(x) \le g_{i+1}(x) \le \sup \{ f_i(x) \, | \, i \in \N \} \le f(x), $$
also konvergiert $(g_i(x))_{i \in \N}$ in $\R$ und $\lim_{i \to \infty} g_i(x) \le \sup \{ f_i(x) \, | \, i \in \N \}$.

Angenommen, es gilt $\lim_{i \to \infty} g_i(x) < \sup \{ f_i(x) \, | \, i \in \N \}$.
Dann ist $\lim_{i \to \infty} g_i(x)$ keine obere Schranke von $\{ f_i(x) \, | \, i \in \N \}$, also existiert $j \in \N$ mit
$$ \lim_{i \to \infty} g_i(x) < f_j(x) \le \sup \{ f_0(x), \ldots, f_j(x) \} = g_j(x) \le \lim_{i \to \infty} g_i(x), $$
Widerspruch! {]}

Zu (ii): Durch Anwendung von (i) auf die Folge $(-f_i)_{i \in \N}$ in $\L_{\R}(M,\varphi)$ erhalten wir $\sup(-f_0,-f_1,\ldots) \in \L_{\R}(M,\varphi)$.
Dann folgt aus Hauptsatz \ref{FA.4.22} (ii) auch $\inf(f_0,f_1,\ldots) = - \sup(-f_0,-f_1,\ldots) \in \L_{\R}(M,\varphi)$. \q 

\begin{Satz}[Lemma von \textsc{Fatou}] \index{Lemma!von \textsc{Fatou}} \label{FA.4.33} $\,$

\noindent \textbf{Vor.:} Seien $\varphi$ ein Quadermaß auf $\R^n$, $M$ eine Teilmenge von $\R^n$ und 
\begin{equation} \label{FA.4.33.0}
\mbox{$(f_i)_{i \in \N}$ eine auf $M$ $\varphi$-konvergente Folge in $\L_{\R}(M,\varphi)$.}
\end{equation}
Ferner existierten eine $\varphi$-Nullmenge $N$ und eine Zahl $C \in \R$ mit
\begin{gather}
\forall_{i \in \N} \forall_{x \in M \setminus N} \, f_i(x) \ge 0, \label{FA.4.33.2} \\
\forall_{i \in \N} \, \int_M f_i \, \d \varphi \le C. \label{FA.4.33.1}
\end{gather}

\noindent \textbf{Beh.:} $\D \lim_{i \to \infty} f_i \in \L_{\R}(M,\varphi)$ und $\D \int_M \left( \lim_{i \to \infty} f_i \right) \, \d \varphi \le \liminf_{i \to \infty} \int_M f_i \, \d \varphi \le C$. 
\end{Satz}

\textit{Beweis.} Aus (\ref{FA.4.33.2}) und \ref{FA.4.32} (ii) folgt, daß durch
\begin{equation} \label{FA.4.33.3}
\forall_{i \in \N} \, g_i := \inf (f_i,f_{i+1},\ldots)
\end{equation}
eine Folge $(f_i)_{i \in \N}$ in $\L_{\R}(M,\varphi)$ definiert wird.
Daher existiert eine $\varphi$-Nullmenge $\widetilde{N}$ so, daß gilt
\begin{gather}
\forall_{i \in \N} \forall_{x \in M \setminus \widetilde{N}} \, g_i(x) \le g_{i+1}(x), \label{FA.4.33.4} \\
\forall_{i \in \N} \forall_{x \in M \setminus \widetilde{N}} \, g_i(x) \le f_i(x), \label{FA.4.33.5}
\end{gather}
also auch nach (\ref{FA.4.33.5}), (\ref{FA.4.33.1})
\begin{equation} \label{FA.4.33.6}
\forall_{i \in \N} \, \int_M g_i \, \d \varphi \le \int_M f_i \, \d \varphi \le C.
\end{equation}

Aus (\ref{FA.4.33.4}), (\ref{FA.4.33.6}) und dem Grenzwertsatz von \textsc{Levi} folgt die $\varphi$-Konvergenz von $(g_i)_{i \in \N}$ auf $M$ und
\begin{equation} \label{FA.4.33.7}
\lim_{i \to \infty} g_i \in \L_{\R}(M,\varphi) \mbox{ sowie } \int_M \left( \lim_{i \to \infty} g_i \right) \, \d \varphi = \lim_{i \to \infty} \int_M g_i \, \d \varphi
\end{equation}
Wegen (\ref{FA.4.33.0}) und (\ref{FA.4.33.3}) existiert daher eine $\varphi$-Nullmenge $N^*$ mit
$$ \forall_{x \in M \setminus N^*} \, \lim_{i \to \infty} f_i(x) = \liminf_{i \to \infty} f_i(x) = \lim_{i \to \infty} g_i(x), $$
also ergibt (\ref{FA.4.33.7}): $\lim_{i \to \infty} f_i \in \L_{\R}(M,\varphi)$ sowie
\begin{eqnarray*}
\int_M \left( \lim_{i \to \infty} f_i \right) \, \d \varphi & = & \int_M \left( \lim_{i \to \infty} g_i \right) \, \d \varphi \stackrel{(\ref{FA.4.33.7})}{=} \lim_{i \to \infty} \int_M g_i \, \d \varphi \\
& \stackrel{(\ref{FA.4.33.6})}{\le} & \liminf_{i \to \infty} \int_M f_i \, \d \varphi \stackrel{(\ref{FA.4.33.6})}{\le} C,
\end{eqnarray*}
womit das Lemma von \textsc{Fatou} bewiesen ist. \q

\begin{HS}[Grenzwertsatz von \textsc{Lebesgue}] \index{Satz!Grenzwert-!von \textsc{Lebesgue}} \label{FA.4.34} $\,$

\noindent \textbf{Vor.:} Seien $\varphi$ ein Quadermaß auf $\R^n$, $M$ eine Teilmenge von $\R^n$ und 
\begin{equation} \label{FA.4.34.1}
\mbox{$(f_i)_{i \in \N}$ eine auf $M$ $\varphi$-konvergente Folge in $\L_{\K}(M,\varphi)$.}
\end{equation}
Ferner existierten eine $\varphi$-Nullmenge $N$ und $g \in \L_{\R}(M,\varphi)$
\begin{equation} \label{FA.4.34.2}
\forall_{i \in \N} \forall_{x \in M \setminus N} \, |f_i(x)| \le g(x).
\end{equation}

\noindent \textbf{Beh.:} $\D \lim_{i \to \infty} f_i \in \L_{\K}(M,\varphi)$ und $\D \int_M \left( \lim_{i \to \infty} f_i \right) \, \d \varphi = \lim_{i \to \infty} \int_M f_i \, \d \varphi$.

\begin{Bem*}
Sei $\K = \R$.
\begin{itemize}
\item[1.)] Existieren $g_1, g_2 \in \L_{\R}(M,\varphi)$ mit
$$ \forall_{i \in \N} \forall_{x \in M \setminus N} \, g_1(x) \le f_i(x) \le g_2(x), $$
so erfüllt $g := \sup(-g_1,g_2)$ die Bedingung (\ref{FA.4.34.2}).
\item[2.)] Gegenüber dem Grenzwertsatz von \textsc{Levi} wird die Funktionenfolge $(f_i)_{i \in \N}$ nicht außerhalb einer Nullmenge als monoton wachsend vorausgesetzt.
Dafür muß beim Grenzwertsatz von \textsc{Lebesgue} die $\varphi$-Konvergenz der Funktionenfolge $(f_i)_{i \in \N}$ gefordert werden, die beim Grenzwertsat von \textsc{Levi} gefolgert werden kann.
\end{itemize}
\end{Bem*}
\end{HS}

\textit{Beweis des Hauptsatzes.} 1. Fall: $\K = \R$.
Dann besagt (\ref{FA.4.34.2})
$$ \forall_{i \in \N} \forall_{x \in M \setminus N} \, -g(x) \le \pm f_i(x) \le g(x). $$
Daher ist $(\pm f_i + g)_{i \in \N}$ wegen (\ref{FA.4.34.1}) eine auf $M$ $\varphi$-konvergente Folge in $\L_{\R}(M,\varphi)$ mit
\begin{gather*}
\forall_{i \in \N} \forall_{x \in M \setminus N} \, (\pm f_i + g)(x) \ge 0, \\
\forall_{i \in \N} \, \int_M (\pm f_i + g) \, \d \varphi \le 2 \int_M g \, \d \varphi =: C.
\end{gather*}
Somit ergibt das Lemma von \textsc{Fatou} \ref{FA.4.33}
$$ \lim_{i \to \infty} (\pm f_i + g) \in \L_{\R}(M,\varphi) \mbox{ und } \int_M \left( \lim_{i \to \infty} (\pm f_i + g) \right) \, \d \varphi \le \liminf_{i \to \infty} \int_M (\pm f_i + g) \, \d \varphi, $$
also gilt offenbar auch
$$ \lim_{i \to \infty} \pm f_i \in \L_{\R}(M,\varphi) \mbox{ und } \int_M \left( \lim_{i \to \infty} \pm f_i \right) \, \d \varphi \le \liminf_{i \to \infty} \int_M \pm f_i \, \d \varphi. $$
Hieraus folgt 
\begin{eqnarray*}
\limsup_{i \to \infty} \int_M f_i \, \d \varphi & = & - \liminf_{i \to \infty} \int_M - f_i \, \d \varphi \le - \int_M \left( \lim_{i \to \infty} - f_i \right) \, \d \varphi \\
& = & \int_M \left( \lim_{i \to \infty} f_i \right) \, \d \varphi \le \liminf_{i \to \infty} \int_M f_i \, \d \varphi,
\end{eqnarray*}
womit der Grenzwertsatz von \textsc{Lebesgue} im Falle $\K = \R$ gezeigt ist.

2. Fall: $\K = \C$.
Da die Voraussetzungen des Hauptsatzes auch für ${\rm Re} \, f_i$ bzw.\ ${\rm Im} \, f_i$ anstelle von $f_i$ erfüllt sind, folgt die Behauptung aus dem 1. Fall. \q
\A
Wir können nun den oben angekündigten Beweis der ersten Aussage von \ref{FA.4.24} (iii) im Falle $\K = \C$ führen.
Diesen bereiten wir durch drei Lemmata vor.

\begin{Lemma} \label{FA.4.35} $\,$

\noindent \textbf{Vor.:} Seien $\varphi$ ein Quadermaß auf $\R^n$, $J \in \mathfrak{I}_n$ ein nicht-leeres Intervall von $\R^n$ und $f \in \mathcal{L}_{\R}(J,\varphi)$.

\noindent \textbf{Beh.:} Es existieren $g,h \in \mathcal{S}_{\nearrow}(J,\varphi)$ mit $f =_{\varphi} g - h$ auf $J$ sowie $g \ge 0$ und $h \ge 0$.
\end{Lemma}

\textit{Beweis.} Per definitionem existieren $\tilde{g}, \tilde{h} \in \mathcal{S}_{\nearrow}(J,\varphi)$ so, daß $f =_{\varphi} \tilde{g} - \tilde{h}$ auf $J$.
Daher folgt aus $\tilde{g} = \tilde{g}^+ - \tilde{g}^-$, $\tilde{h} = \tilde{h}^+ - \tilde{h}^-$
$$ f =_{\varphi} (\tilde{g}^+ - \tilde{g}^-) - (\tilde{h}^+ - \tilde{h}^-) = (\tilde{g}^+ + \tilde{h}^-) - (\tilde{g}^- + \tilde{h}^+), $$
also leisten $g := \tilde{g}^+ + \tilde{h}^-$ und $h := \tilde{g}^- + \tilde{h}^+$ offenbar das Gewünschte. \q

\begin{Lemma} \label{FA.4.36}
Seien $t_1, t_2, s_1, s_2$ nicht-negative reelle Zahlen.

Dann gilt $| (t_1 - t_2) + \i \, (s_1 - s_2)| \le t_1 + t_2 + s_1 + s_2$.
\end{Lemma}

\textit{Beweis.} Für alle $t,s \in \R$ mit $t,s \ge 0$ gilt $t^2 + s^2 \le t^2 + 2 \, t \, s + s^2 = (t+s)^2$, d.h.\
\begin{equation} \label{FA.4.36.S}
\forall_{t,s \in \R, \, t,s \ge 0} \, \sqrt{t^2 + s^2} \le t+s,
\end{equation}
also folgt wegen $t_1, t_2, s_1, s_2 \ge 0$
\begin{eqnarray*}
|(t_1 - t_2) + \i \, (s_1 - s_2)| & = & \sqrt{(t_1 - t_2)^2 + (s_1 - s_2)^2} \\
& = & \sqrt{ {t_1}^2 - 2 \, t_1 \, t_2 + {t_2}^2 + {s_1}^2 - 2 \, s_1 \, s_2 + {s_2}^2 } \\
& \le & \sqrt{ {t_1}^2 + 2 \, t_1 \, t_2 + {t_2}^2 + {s_1}^2 + 2 \, s_1 \, s_2 + {s_2}^2 } \\
& = & \sqrt{(t_1 + t_2)^2 + (s_1 + s_2)^2} \stackrel{(\ref{FA.4.36.S})} {\le} t_1 + t_2 + s_1 + s_2 
\end{eqnarray*}
\q

\begin{Lemma} \label{FA.4.31}
Seien $\varphi$ ein Quadermaß auf $\R^n$, $J \in \mathfrak{I}_n$ ein nicht-leeres Intervall von $\R^n$ und $T,S \in \mathcal{S}(J)$.

Dann gilt $|T + \i \, S| \in \L_{\R}(J,\varphi)$.
\end{Lemma}

\textit{Beweis.} Wegen \ref{FA.4.8} (ii) existieren $k \in \N_+$, $\alpha_1, \ldots, \alpha_k$, $\beta_1, \ldots, \beta_k \in \R$ und paarweise disjunkte Quader $Q_1, \ldots, Q_k$ mit $T = \sum_{i=1}^k \alpha_i \, \chi_{Q_i}^J$ und $S = \sum_{i=1}^k \beta_i \, \chi_{Q_i}^J$.
Dann gilt
$$ \bigg| \, | T + \i \, S | \, \bigg| = \left| \sum_{i=1}^k (\alpha_i + \i \, \beta) \, \chi_{Q_i}^J \right| \le \sum_{i=1}^k \left| \alpha_i + \i \, \beta_i \right| \, \chi_{Q_i}^J, $$
wobei $| T + \i \, S |$ als konstante -- also konvergente -- Folge in $\L_{\R}(J,\varphi)$ aufgefaßt werden kann und die rechte Seite der letzten Ungleichung als Treppenfunktion $\varphi$-integrierbar ist.
Der Grenzwertsatz von \textsc{Lebesgue} liefert daher die $\varphi$-Integrierbarkeit von $| T + \i \, S |$.
\q
\A
\textit{Beweis der ersten Aussage von \ref{FA.4.24} (iii) im Falle $\K = \C$.}\label{Beweis FA.4.24 (iii)}
Seien also $\varphi$ ein Quadermaß auf $\R^n$, $M \subset \R^n$ und $f \in \L_{\C}(M,\varphi)$.
Wir wollen zeigen
$$ |f| \in \L_{\R}(M,\varphi). $$

1. Fall: $J := M \in \mathfrak{I}_n$ mit $J \ne \emptyset$.
Per definitionem sind ${\rm Re} \, f$ und ${\rm Im} \, f$ Elemente von $\L_{\R}(J,\varphi)$, also folgt aus \ref{FA.4.35} die Existenz nicht-negativer Funktionen $g_1, h_1, g_2, h_2 \in \mathcal{S}_{\nearrow}(J,\varphi)$ mit
\begin{equation} \label{FA.4.24.1}
{\rm Re} \, f =_{\varphi} g_1 - h_1 ~~ \mbox{ und } ~~ {\rm Im} \, f =_{\varphi} g_2 - h_2.
\end{equation}
Nach Definition von $\mathcal{S}_{\nearrow}(J,\varphi)$ existieren dann weiter monoton wachsende Folgen $(T_{1,i})_{i \in \N}, (T_{2,i})_{i \in \N}, (S_{1,i})_{i \in \N}, (S_{2,i})_{i \in \N}$ in $\mathcal{S}(J)$ mit
\begin{equation} \label{FA.4.24.2}
\lim_{i \to \infty} T_{1,i} =_{\varphi} g_1, ~~ \lim_{i \to \infty} T_{2,i} =_{\varphi} h_1, ~~ \lim_{i \to \infty} S_{1,i} =_{\varphi} g_2, ~~ \lim_{i \to \infty} S_{2,i} =_{\varphi} h_2,
\end{equation}
also gilt wegen  (\ref{FA.4.24.1})
\begin{equation} \label{FA.4.24.3}
f =_{\varphi} \lim_{i \to \infty} \left( ( T_{1,i} - T_{2,i} ) + \i \, ( S_{1,i} - S_{2,i} ) \right)
\end{equation}
und aus Stetigkeitsgründen auch
\begin{equation} \label{FA.4.24.4}
|f| =_{\varphi} \lim_{i \to \infty} \underbrace{\left| ( T_{1,i} - T_{2,i} ) + \i \, ( S_{1,i} - S_{2,i} ) \right|}_{\stackrel{\ref{FA.4.31}}{\in}\L_{\R}(J,\varphi)}.
\end{equation}

Da $g_1, h_1, g_2, h_2$ nicht-negativ sind, können wir durch Übergang von $(T_{1,i})_{i \in \N}$, $(T_{2,i})_{i \in \N}, (S_{1,i})_{i \in \N}$ und $(S_{2,i})_{i \in \N}$ zu ${(T_{1,i}}^+)_{i \in \N}, ({T_{2,i}}^+)_{i \in \N}, ({S_{1,i}}^+)_{i \in \N}$ und $({S_{2,i}}^+)_{i \in \N}$ ohne Beschränkung der Allgemeniheit annehmen, daß auch die Folgen $(T_{1,i})_{i \in \N}$, $(T_{2,i})_{i \in \N}, (S_{1,i})_{i \in \N}$ und $(S_{2,i})_{i \in \N}$ nicht-negativ sind.
Somit ergibt \ref{FA.4.36}
\begin{equation} \label{FA.4.24.5}
\forall_{i \in \N} \, \left| ( T_{1,i} - T_{2,i} ) + \i \, ( S_{1,i} S_{2,i} ) \right| \le  T_{1,i} + T_{2,i} + S_{1,i} + S_{2,i}
\end{equation}
Hiearaus und der Tatsache, daß $(T_{1,i})_{i \in \N}$, $(T_{2,i})_{i \in \N}, (S_{1,i})_{i \in \N}$ und $(S_{2,i})_{i \in \N}$ monoton wachsend sind, folgt wegen (\ref{FA.4.24.2})
\begin{equation} \label{FA.4.24.6}
\forall_{i \in \N} \, \left| ( T_{1,i} - T_{2,i} ) + \i ( S_{1,i} - S_{2,i} ) \right| \le_{\varphi}  g_1 + h_1 + g_2 + h_2  \in \mathcal{S}_{\nearrow}(J,\varphi) \subset \L(J,\varphi).
\end{equation}

Aus (\ref{FA.4.24.4}), (\ref{FA.4.24.6}), $\big| |f| \big| = |f|$ und dem Grenzwertsatz von \textsc{Lebesgue} folgt $|f| \in \L_{\R}(J,\varphi)$, d.h.\ die Behauptung ist im ersten Falle gezeigt.

2. Fall: $M$ beliebige Teilmenge von $\R^n$.
Sei also $f \in \L_{\C}(M,\varphi)$, d.h.\ genau $\chi_M \, \widehat{{\rm Re} \, f}, \, \chi_M \, \widehat{{\rm Im} \, f} \in \L_{\R}(\R^n,\varphi)$ bzw.\ $\chi_M \, \widehat{{\rm Re} \, f} + \i \, \chi_M \, \widehat{{\rm Im} \, f} \in \L_{\C}(\R^n,\varphi)$. 
Zu zeigen ist $\chi_M \, \widehat{|f|} \in \L_{\R}(\R^n,\varphi)$.
Wegen
$$ \chi_M \, \widehat{|f|} = \chi_M \, |\hat{f}| = \chi_M \left| \widehat{{\rm Re} \, f} + \i \, \widehat{{\rm Im} \, f} \right| = \left| \chi_M \, \widehat{{\rm Re} \, f} + \i \, \chi_M \, \widehat{{\rm Im} \, f} \right| $$
folgt dies aus dem ersten Fall. \q

\begin{Bem*}
Seien $\varphi$ ein Quadermaß auf $\R^n$ und $M$ eine Teilmenge von $\R^n$.
Wir werden im folgenden nicht mehr explizit darauf hinweisen, daß mit einer $\varphi$-integrierbaren Funktion über $M$ auch ihr Betrag $\varphi$-integrierbar über $M$ ist.
\end{Bem*}

\begin{Satz} \label{FA.4.37} $\,$

\noindent \textbf{Vor.:} Seien $\varphi$ ein Quadermaß auf $\R^n$, $J \in \mathfrak{I}_n$ ein nicht-leeres Intervall von $\R^n$ und $f \in \mathcal{L}_{\K}(J,\varphi)$.

\noindent \textbf{Beh.:} Zu jedem $\varepsilon \in \R_+$ existieren $T,S \in \mathcal{S}(J)$ mit $\D \int_J |f - (T + \i \, S)| ~ \d \varphi < \varepsilon$.

\begin{Zusatz}
Im Falle $\K = \R$ kann $S = 0$ gewählt werden.
\end{Zusatz}

\begin{Bem*}
Wegen \ref{FA.4.8} (ii) existieren $k \in \N_+$, $\zeta_1, \ldots, \zeta_k \in \C$ und paarweise disjunkte Quader $Q_1, \ldots, Q_k$ mit $T + \i \, S = \sum_{i=1}^k \zeta_i \, \chi_{Q_i}^J$.
\end{Bem*}
\end{Satz}

\textit{Beweis.} Wir verwenden die Bezeichnungen und Eigenschaften des 1. Falles des zuletzt geführten Beweises.
Dann ergibt (\ref{FA.4.24.3})
\begin{equation} \label{FA.4.37.1}
\lim_{i \to \infty} \big| \big( ( T_{1,i} - T_{2,i} ) + \i \, ( S_{1,i} - S_{2,i} ) \big) - f \big| =_{\varphi} 0.
\end{equation}

Da $(T_{1,i})_{i \in \N}, (T_{2,i})_{i \in \N}, (S_{1,i})_{i \in \N}, (S_{2,i})_{i \in \N}$ Folgen in $\mathcal{S}(J)$ sind, gilt außerdem
\begin{equation} \label{FA.4.37.2}
\forall_{i \in \N} \, \big| \big( ( T_{1,i} - T_{2,i} ) + \i \, ( S_{1,i} - S_{2,i} ) \big) - f \big| \in \L_{\R}(J,\varphi).
\end{equation}

Des weiteren folgt aus (\ref{FA.4.24.6})
\begin{equation} \label{FA.4.37.3}
\forall_{i \in \N} \, \big| \big( ( T_{1,i} - T_{2,i} ) + \i \, ( S_{1,i} - S_{2,i} ) \big) - f \big| \le_{\varphi} g_1 + h_1 + g_2 + h_2 + |f| \in \L_{\R}(J,\varphi).
\end{equation}

(\ref{FA.4.37.1}) - (\ref{FA.4.37.3}) und der Grenzwertsatz von \textsc{Lebesgue} liefern 
$$ \lim_{t \to \infty} \big| \big( ( T_{1,i} - T_{2,i} ) + \i \, ( S_{1,i} - S_{2,i} ) \big) - f \big| \in \L_{\R}(J,\varphi) $$
sowie
$$ \lim_{i \to \infty} \int_J \big| \big( ( T_{1,i} - T_{2,i} ) + \i \, ( S_{1,i} - S_{2,i} ) \big) - f \big| = 0, $$
woraus die Behauptung offenbar für hinreichend großes $i \in \N$ folgt.

Der Zusatz ist unter Berücksichtigung von (\ref{FA.4.24.1}) klar. \q

\begin{Satz} \label{FA.4.30} $\,$

\noindent \textbf{Vor.:} Seien $\varphi$ ein Quadermaß auf $\R^n$, $M$ eine Teilmenge von $\R^n$ und $f$ eine auf einer Teilmenge von $\R^n$ definierte $\K$-wertige Funktion.

\noindent \textbf{Beh.:} Die folgenden Aussagen sind äquivalent.
\begin{itemize}
\item[(i)] $f \in \L_{\K}(M,\varphi)$ und $\D \int_M |f| \, \d \varphi = 0$.
\item[(ii)] Es existiert  eine $\varphi$-Nullmenge $N$ mit $\forall_{x \in M \setminus N} \, f(x) = 0$ -- insbes.\ ist $f$ auf $M \setminus N$ definiert.
\end{itemize}
\end{Satz}

\textit{Beweis.} ,,(i) $\Rightarrow$ (ii)`` Nach Voraussetzung ist $(i \, |f|)_{i \in \N}$ eine außerhalb von $N$ mo\-no\-ton wachsende Folge in $\L_{\R}(M,\varphi)$ mit durch $0$ beschränkter $\varphi$-In\-te\-gral\-folge.
Daher folgt aus dem Grenzwertsatz von \textsc{Levi} ihre $\varphi$-Konvergenz auf $M$, d.h.\ es existiert eine $\varphi$-Nullmenge $N$ derart, daß gilt
$$ \forall_{x \in M \setminus N} \, \underbrace{|f(x)|}_{\ge 0} \, \underbrace{\lim_{i \to \infty} i}_{= \infty} = \lim_{i \to \infty} i \, |f(x)| \in \R, $$
also $\forall_{x \in M \setminus N} \, f(x) = 0$.

 ,,(ii) $\Rightarrow$ (i)`` ist klar. 
 \q

\begin{Kor} \label{FA.4.30.K}
Seien $\varphi$ ein Quadermaß auf $\R^n$ und $M$ eine Teilmenge von $\R^n$.

Dann gilt: $M$ $\varphi$-Nullmenge $\Longleftrightarrow$ $\D \left( \chi_M \in \L_{\R}(\R^n,\varphi) \, \wedge \, \int_{\R^n} \chi_M \, \d \varphi = 0 \right)$.
\end{Kor}

\textit{Beweis.} ,,$\Rightarrow$`` Ist $M$ eine $\varphi$-Nullmenge, so gilt $\chi_M =_{\varphi} 0$.
Aus \ref{FA.4.30} ,,(ii) $\Rightarrow$ (i)`` (dort $f = \chi_M$ und $M = \R^n$) folgt die rechte Seite der obigen Äquivalenz.

,,$\Leftarrow$`` Aus der rechten Seite der o.g.\ Äquivalenz und \ref{FA.4.30} ,,(i) $\Leftarrow$ (ii)`` folgt die Existenz einer $\varphi$-Nullmenge $N$ mit $\forall_{x \in M \setminus N} \, \chi_M(x) = 0$, d.h.\ genau $M \subset N$.
Daher ist $M$ mit $N$ eine $\varphi$-Nullmenge. \q

\begin{Satz} \label{FA.4.61} $\,$

\noindent \textbf{Vor.:} Seien $\varphi$ ein Quadermaß auf $\R^n$, $M \subset \R^n$ und $(M_i)_{i \in \N}$ eine Folge von Teilmengen von $M$ derart, daß gilt
\begin{gather}
\forall_{i \in \N} \, M_i \subset M_{i+1} \subset M \mbox{ und} \label{FA.4.61.1} \\
M \setminus \bigcup_{i=0}^{\infty} M_i \mbox{ ist eine $\varphi$-Nullmenge.} \label{FA.4.61.2}
\end{gather}
Ferner seien $D \subset \R^n$ und $f \in {\R}^D$ mit
\begin{equation} \label{FA.4.61.3}
\forall_{i \in \N} \, f \in \L_{\R}(M_i,\varphi).
\end{equation}

\noindent \textbf{Beh.:} Es gilt:
\begin{itemize}
\item[(i)] $\D f \in \L_{\R}(M,\varphi) \Longrightarrow \int_M f \, \d \varphi = \lim_{i \to \infty} \int_{M_i} f \, \d \varphi$.
\item[(ii)] $\D \left( f \ge_{\varphi} 0 \, \wedge \, C := \lim_{i \to \infty} \int_{M_i} f \, \d \varphi \in \R \right) \Longrightarrow f \in \L_{\R}(M,\varphi)$.
\end{itemize}
\end{Satz}

\textit{Beweis.} Aus (\ref{FA.4.61.1}), (\ref{FA.4.61.2}) folgt offenbar
\begin{equation} \label{FA.4.61.4}
\lim_{i \to \infty} \chi_{M_i} \, \hat{f} = \chi_M \, \hat{f}
\end{equation}
und wegen (\ref{FA.4.61.3}) gilt
\begin{equation} \label{FA.4.61.5}
\forall_{i \in \N} \; \chi_{M_i} \, \hat{f} \in \L_{\R}(\R^n,\varphi).
\end{equation}

Zu (i): Aus $f \in \L_{\R}(M,\varphi)$ folgt $|f| \in \L_{\R}(M,\varphi)$, d.h.\ $\chi_M \, \widehat{|f|} \in \L_{\R}(\R^n,\varphi)$.
Daher gilt nach (\ref{FA.4.61.1})
\begin{equation} \label{FA.4.61.6}
\forall_{i \in \N} \, \underbrace{- \chi_M \, \widehat{|f|}}_{\in \L_{\R}(\R^n,\varphi)} \le \chi_{M_i} \, \hat{f} \le \underbrace{\chi_M \, \widehat{|f|}}_{\in \L_{\R}(\R^n,\varphi)}.
\end{equation}
Aus (\ref{FA.4.61.5}), (\ref{FA.4.61.4}), (\ref{FA.4.61.6}) und dem Grenzwertsatz von \textsc{Lebesgue} folgt die Behauptung von (i).

Zu (ii): Aus der Voraussetzung von (ii) folgt mittels (\ref{FA.4.61.1})
\begin{equation} \label{FA.4.61.7}
\forall_{i \in \N} \; \chi_{M_i} \, \hat{f} \le_{\varphi} \chi_{M_{i+1}} \, \hat{f},
\end{equation}
also auch
$$ \forall_{i \in \N} \, \int_{\R^n} \chi_{M_i} \, \hat{f} \d \varphi \le \int_{\R^n} \chi_{M_{i+1}} \, \hat{f} \, \d \varphi $$
und folglich
\begin{equation} \label{FA.4.61.8}
\forall_{i \in \N} \, \int_{\R^n} \chi_{M_i} \, \hat{f} \d \varphi \le C.
\end{equation}
Aus (\ref{FA.4.61.5}), (\ref{FA.4.61.7}), (\ref{FA.4.61.8}), (\ref{FA.4.61.4}) und dem Grenzwertsatz von \textsc{Levi} ergibt sich die Behauptung von (ii). \q

\begin{HS}[Charakterisierung des Lebesgue-Integrales der auf $\R^n$ fast überall-definierten reellwertigen Funktionen] \label{FA.4.28} $\,$

\noindent \textbf{Vor.:} Sei $\varphi$ ein Quadermaß auf $\R^n$.

\noindent \textbf{Beh.:} Es existieren genau ein $\R$-Un\-ter\-vek\-tor\-raum $\L(\varphi)$ des $\R$-Vek\-tor\-raumes $\{ f \in \R^D \, | \, D  \subset \R^n \, \wedge \, (\R^n \setminus D) \mbox{ $\varphi$-Nullmenge} \}$ der sog.\ \emph{auf $\R^n$ $\varphi$-definierten reellwertigen Funktionen} und genau ein Funktional $I \: \L(\varphi) \to \R$ derart, daß folgende Bedingungen erfüllt sind:
\begin{itemize}
\item[1.)] $\forall_{Q \in \mathfrak{Q}_n} \, \chi_Q \in \L(\varphi) \, \wedge \, I(\chi_Q) = \varphi(Q)$.
\item[2.)] $\forall_{f,g \in \L(\varphi)} \, \left( f \le_{\varphi} g \Longrightarrow I(f) \le I(g) \right)$.
\item[3.)] Ist $f$ eine auf einer Teilmenge von $\R^n$ definierte reellwertige Funktion und existiert $g \in \L(\varphi)$ mit $f =_{\varphi} g$, so gilt auch $f \in \L(\varphi)$.
\item[4.)] Zu jeder Folge $(f_i)_{i \in \R}$ in $\L(\varphi)$ mit $\forall_{i \in \N} \, f_i \le_{\varphi} f_{i+1}$ und $\exists_{C \in \R} \forall_{i \in \N} \, I(f_i) \le C$ existiert eine Funktion $f \in \L(\varphi)$ mit $\D f =_{\varphi} \lim_{i \to \infty} f_i$ und $\D I(f) = \lim_{i \to \infty} I(f_i)$.
\item[5.)] Sei $X$ ein $\R$-Untervektorraum von $\L(\varphi)$ mit folgenden Eigenschaften:
\begin{itemize} 
\item[(a)] $\forall_{Q \in \mathfrak{Q}_n} \, \chi_Q \in X$.
\item[(b)] Sind $f \in L(\varphi)$ und $(f_i)_{i \in \R}$ eine Folge in $X$ mit $\forall_{i \in \N} \, 0 \le f_i \le_{\varphi} f_{i+1}$, $\exists_{C \in \R} \forall_{i \in \N} \, I(f_i) \le C$ und $\D f =_{\varphi} \lim_{i \to \infty} f_i$, so gilt $f \in X$.
\end{itemize}

Dann folgt $X = \L(\varphi)$.
\end{itemize}
\end{HS}

Zum Nachweis des Hauptsatzes beweisen wir zunächst das folgende Lemma.

\begin{Lemma} \label{FA.4.29} $\,$

\noindent \textbf{Vor.:} Es seien $\varphi$ ein Quadermaß auf $\R^n$ sowie $X$ ein $\R$-Untervektorraum von $\L_{\R}(\R^n,\varphi)$, der die folgende Eigenschaften erfüllt:
\begin{itemize} 
\item[(a)] $\forall_{Q \in \mathfrak{Q}_n} \, \chi_Q \in X$.
\item[(b)] Sind $f \in \L_{\R}(\R^n,\varphi)$ und $(f_i)_{i \in \R}$ eine Folge in $X$ mit $\forall_{i \in \N} \, 0 \le f_i \le_{\varphi} f_{i+1}$, $\D \exists_{C \in \R} \forall_{i \in \N} \, \int_{\R^n} f_i \, \d \varphi \le C$ und $\D f =_{\varphi} \lim_{i \to \infty} f_i$, so gilt $f \in X$.
\end{itemize}

\noindent \textbf{Beh.:} $X = \L_{\R}(\R^n,\varphi)$.
\end{Lemma}

\textit{Beweis.} Wegen (a) folgt aus \ref{FA.4.8} (i) sofort
\begin{equation} \label{FA.4.29.1}
\mathcal{S}(\R^n) \subset X.
\end{equation}

Zum Nachweis des Lemmas genügt es zu zeigen
\begin{equation} \label{FA.4.29.2}
\mathcal{S}_{\nearrow}(\R^n,\varphi) \subset X,
\end{equation}
denn zu $f \in \L_{\R}(\R^n,\varphi)$ existieren $g,h \in \mathcal{S}_{\nearrow}(\R^n,\varphi)$ mit $f =_{\varphi} g - h \stackrel{(\ref{FA.4.29.2})}{\in} X$.

Zu (\ref{FA.4.29.2}): Sei $f \in \mathcal{S}_{\nearrow}(\R^n,\varphi)$.
Dann existiert eine monoton wachsende Folge $(T_i)_{i \in \N}$ in $\mathcal{S}(\R^n) \stackrel{(\ref{FA.4.29.1})}{\subset} X$ mit beschränkter $\varphi$-Integralfolge und $f =_{\varphi} \lim_{i \to \infty} T_i$.
Somit ist $(T_i + T_0)_{i \in \N}$ eine nicht-negative monoton wachsende Folge in $X$ mit beschränkter $\varphi$-Integralfolge und $\lim_{i \to \infty} (T_i + T_0) =_{\varphi} f + T_0$.
Aus (b) folgt dann $f + T_0 \in X$, also auch $f \in X$. \q
\A
\textit{Beweis des Hauptsatzes.} Wegen $\forall_{Q \in \mathfrak{Q}_n} \, \chi_Q \in \L_{\R}(\R^n,\varphi) \wedge \int_{\R^n} \chi_Q \, \d \varphi = \varphi(Q)$ und \ref{FA.4.20} (i) - (iii), (vi) sowie \ref{FA.4.25} und \ref{FA.4.29} ist die Existenz klar.
Zum Nachweis der Eindeutigkeit genügt es daher zu zeigen, daß für $\L(\varphi)$ und $I$ wie in der Behauptung gilt 
$$ \L(\varphi) = \L_{\R}(\R^n,\varphi) \mbox{ und } I = \int_{\R^n} \ldots \, \d \varphi. $$

Beweis hiervon: Wegen 1.), \ref{FA.4.8} (i) und \ref{FA.4.11} gilt 
\begin{equation} \label{FA.4.28.1}
\mathcal{S}(\R^n) \subset \L(\varphi) \mbox{ und } \forall_{T \in \mathcal{S}(\R^n)} \, I(T) = \int_{\R^n} T \, \d \varphi.
\end{equation}
Hieraus folgern wir
\begin{gather}
\mathcal{S}_{\nearrow}(\R^n,\varphi) \subset \L(\varphi) \mbox{ und } \forall_{f \in \mathcal{S}_{\nearrow}(\R^n,\varphi)} \, I(f) = \int_{\R^n} f \, \d \varphi, \label{FA.4.28.2} \\
\L_{\R}(\R^n,\varphi) \subset \L(\varphi) \mbox{ und } \forall_{f \in \L(\R^n,\varphi)} \, I(f) = \int_{\R^n} f \, \d \varphi. \label{FA.4.28.3}
\end{gather}

{[} Zu (\ref{FA.4.28.2}): Sei $f \in \mathcal{S}_{\nearrow}(\R^n,\varphi)$.
Dann existiert eine monoton wachsende Folge $(T_i)_{i \in \N}$ in $\mathcal{S}(\R^n)$ mit beschränkter $\varphi$-Integralfolge und $f =_{\varphi} \lim_{i \to \infty} T_i$.
Wegen (\ref{FA.4.28.1}) gilt
$$ \forall_{i \in \N} \, T_i \in \L(\varphi) \, \wedge \, I(T_i) = \int_{\R^n} T_i \, \d \varphi. $$
Daher folgt aus 4.) die Existenz von $g \in \L(\varphi)$ mit $g =_{\varphi} \lim_{i \to \infty} T_i = f$ und 
$$ I(g) = \lim_{i \to \infty} I(T_i) = \lim_{i \to \infty} \int_{\R^n} T_i \, \d \varphi \stackrel{\text{Def.}}{=} \int_{\R^n} f \, \d \varphi. $$
Nach 3.), 2.) gilt dann $f \in \L(\varphi)$ und $I(g) = I(f)$, also auch $I(f) = \int_{\R^n} f \, \d \varphi$.

Zu (\ref{FA.4.28.3}): Sei $f \in \L_{\R}(\R^n,\varphi)$.
Dann existieren $g,h \in \mathcal{S}_{\nearrow}(\R^n,\varphi)$ mit $f =_{\varphi} g - h$.
Wegen (\ref{FA.4.28.2}) gilt $g-h \in \L(\varphi)$, also ergibt 3.), daß auch $f \in \L(\varphi)$ sowie
$$ I(f) = I(g-h) = I(g) - I(h) \stackrel{(\ref{FA.4.28.2})}{=} \int_{\R^n} g \, \d \varphi - \int_{\R^n} h \, \d \varphi \stackrel{\text{Def.}}{=} \int_{\R^n} f \, \d \varphi $$
gilt. {]}

Wir beweisen schließlich
\begin{equation} \label{FA.4.28.4}
\L_{\R}(\R^n,\varphi) = \L(\varphi),
\end{equation}
womit der Hauptsatz wegen (\ref{FA.4.28.3}) bewiesen ist.

{[} Zu (\ref{FA.4.28.4}): $\L_{\R}(\R^n,\varphi)$ ist nach (\ref{FA.4.28.3}) ein $\R$-Untervektorraum von $\L(\varphi)$.
Es genügt zu zeigen , daß $\L_{\R}(\R^n,\varphi)$ mutatis mutandis die Eigenschaften (a) und (b) in 5.) erfüllt.
(a) ist trivial.

Zu (b): Seien $f \in L(\varphi)$, $(f_i)_{i \in \R}$ eine Folge in $\L_{\R}(\R^n,\varphi)$ und $C \in \R$ derart, daß $\forall_{i \in \N} \, 0 \le f_i \le_{\varphi} f_{i+1} \, \wedge \, I(f_i) \le C$ sowie $f =_{\varphi} \lim_{i \to \infty} f_i$ gilt.
Dann ist $(f_i)_{i \in \R}$ auch eine Folge in $\L(\varphi)$, also ergibt sich $f \in \L(\varphi)$ aus 4.).{]} \q

\subsection*{Riemann- bzw.\ Stieltjes-Integral und absolut konvergente Reihen} \addcontentsline{toc}{subsection}{Riemann- bzw.\ Stieltjes-Integral und absolut konvergente Reihen}

Der Autor dieser Zeilen hat sich entschieden, hier nicht auf die Konstruktion des Riemann-Integrales einzugehen, verwendet jedoch die z.B.\ in \cite[§ 6]{Henke} oder \cite[Kapitel 8]{ElAna} gegebene Einführung und die in letzterer Quelle definierten Bezeichnungen und Resultate; beachte allerdings \ref{FA.4.40} unten.

\begin{HS}[Charakterisierung eigentlicher Riemann-Integrierbarkeit] \index{Integral!Riemann-} \label{FA.4.38} $\,$

\noindent \textbf{Vor.:} Seien $a,b \in \R$ mit $a<b$ und $f \: [a,b] \to \R$ eine beschränkte Funktion.

\noindent \textbf{Beh.:} Die folgenden Aussagen (i) und (ii) sind äquivalent:
\begin{itemize}
\item[(i)] $f$ ist Rie\-mann-in\-te\-grier\-bar über $[a,b]$.
\item[(ii)] Die Menge der Unstetigkeitsstellen von $f$ ist eine $\mu_1$-Null\-menge.
\end{itemize}

Des weiteren folgt aus (i) -- oder äquivalent (ii) --:
\begin{itemize}
\item[(iii)] $f$ ist $\mu_1$-integrierbar über $[a,b]$ und $\D \int_{[a,b]} f \, \d \mu_1 = \int_a^b f(t) \, \d t$.
\end{itemize}
\end{HS}

\textit{Beweis.} Ohne Einschränkung können wir $a=0$ und $b=1$ annehmen.
Für $k \in \N_+$ seien
\begin{equation} \label{FA.4.38.0}
a_{k,0} := \frac{0}{2^k}, \, a_{k,1} := \frac{1}{2^k}, \, \ldots, a_{k,2^k} := \frac{2^k}{2^k} = 1.
\end{equation}

Dann ist für jedes $k \in \N_+$ durch
\begin{equation} \label{FA.4.38.1}
\SZk : 0 = a_{k,0} < a_{k,1} < \ldots < a_{k,2^k} = 1
\end{equation}
eine Zerlegung $\SZk$ von $[0,1]$ gegeben, und $\mathfrak{Z}_{k+1}$ ist eine Verfeinerung von $\SZk$.
Da $(\SZk)_{k \in \N}$ des weiteren offenbar eine ausgezeichente Folge von Zerlegungen von $[0,1]$ ist, vgl.\ \cite[8.12 sowie Definition 8.1 (iii)]{ElAna}, liefert \cite[8.12]{ElAna}
\begin{equation} \label{FA.4.38.2}
\UIf = \lim_{k \to \infty} \SUS(f,\mathfrak{Z}_k) ~~ \mbox{ und } ~~ \OIf = \lim_{k \to \infty} \SOS(f,\mathfrak{Z}_k).
\end{equation}

Wir setzen für $k \in \N_+$
\begin{gather}
m_{k,1} := \inf f([a_{k,0}, a_{k,1}]), \, \ldots, \, m_{k,2^k} := \inf f([a_{k,2^k-1}, a_{k,2^k}]), \label{FA.4.38.3} \\
M_{k,1} := \sup f([a_{k,0}, a_{k,1}]), \, \ldots, \, M_{k,2^k} := \sup f([a_{k,2^k-1}, a_{k,2^k}]), \label{FA.4.38.4}
\end{gather}
d.h.\ nach (\ref{FA.4.38.1}), (\ref{FA.4.38.3}), (\ref{FA.4.38.4}) sowie \cite[Definition 8.2]{ElAna}
\begin{gather}
\SUS(f,\mathfrak{Z}_k) = \sum_{i=1}^{2^k} m_{k,i} \, (a_{k,i} - a_{k,i-1}), \label{FA.4.38.5} \\
\SOS(f,\mathfrak{Z}_k) = \sum_{i=1}^{2^k} M_{k,i} \, (a_{k,i} - a_{k,i-1}). \label{FA.4.38.6}
\end{gather}

Setzen wir des weiteren für $k \in \N_+$
\begin{equation} \label{FA.4.38.7}
Q_{k,1} := [0, a_{k,1}], \, Q_{k,2} := {]} a_{k,1}, a_{k,2} {]}, \, \ldots, \, Q_{k,2^k} := {]} a_{k,2^k-1}, a_{k,2^k} {]},
\end{equation}
so werden nach (\ref{FA.4.38.1}), (\ref{FA.4.38.3}), (\ref{FA.4.38.4}) -- da $\mathfrak{Z}_{k+1}$ Verfeinerung von $\SZk$ ist -- durch
\begin{equation} \label{FA.4.38.8}
T_k := \sum_{i=1}^{2^k} m_{k,i} \, \chi_{Q_{k,i}}^{[0,1]} ~~  \mbox{ und } ~~ S_k := \sum_{i=1}^{k} M_{k,i} \, \chi_{Q_{k,i}}^{[0,1]} 
\end{equation}
monoton wachsende Folgen $(T_k)_{k \in \N_+}$ sowie $(- S_k)_{k \in \N_+}$ in $\mathcal{S}([0,1])$ mit
\begin{equation} \label{FA.4.38.le}
\forall_{k \in \N_+} \, T_k \le f \le S_k
\end{equation}
definiert, und es gilt nach (\ref{FA.4.38.5}), (\ref{FA.4.38.6}) sowie \cite[Definition 8.4]{ElAna}
\begin{gather}
\forall_{k \in \N_+} \, \int\limits_{[0,1]} T_k \, \d \mu_1 = \SUS(f,\mathfrak{Z}_k) \le \UIf \in \R, \label{FA.4.38.9} \\ 
\forall_{k \in \N_+} \, \int\limits_{[0,1]} - S_k \, \d \mu_1 = - \SOS(f,\mathfrak{Z}_k) \le - \OIf \in \R,  \label{FA.4.38.10}
\end{gather}
somit ergibt der Grenzwertsatz von \textsc{Levi}
\begin{gather}
g_1 := \lim_{k \to \infty} T_k \in \L_{\R}([0,1], \mu_1), \label{FA.4.38.11} \\
\int\limits_{[0,1]} g_1 \, \d \mu_1 = \lim_{k \to \infty} \int\limits_{[0,1]} T_k \, \d \mu_1 \stackrel{(\ref{FA.4.38.9}), (\ref{FA.4.38.2})}{=} \UIf, \label{FA.4.38.12}\\
\tilde{g}_2 := \lim_{k \to \infty} - S_k \in \L_{\R}([0,1], \mu_1),  \mbox{ also auch } g_2 := \lim_{k \to \infty} S_k \in \L_{\R}([0,1], \mu_1), \label{FA.4.38.13} \\
\int\limits_{[0,1]} g_2 \, \d \mu_1 = \lim_{k \to \infty} \int\limits_{[0,1]} S_k \, \d \mu_1 \stackrel{(\ref{FA.4.38.10}), (\ref{FA.4.38.2})}{=} \OIf, \label{FA.4.38.14}
\end{gather}
d.h.\ nach (\ref{FA.4.38.14}) sowie (\ref{FA.4.38.12})
\begin{equation} \label{FA.4.38.15}
\int\limits_{[0,1]} (g_2 - g_1) \, \d \mu_1 = \OIf - \UIf.
\end{equation}
Aus (\ref{FA.4.38.3}), (\ref{FA.4.38.4}), (\ref{FA.4.38.7}), (\ref{FA.4.38.8}), (\ref{FA.4.38.11}) und (\ref{FA.4.38.13}) folgt des weiteren
\begin{equation} \label{FA.4.38.16}
g_1 = \lim_{k \to \infty} T_k \le_{\mu_1} f \le_{\mu_1} \lim_{k \to \infty} S_k = g_2.
\end{equation}

Wir setzen nun
\begin{gather}
R := \bigcup_{i=1}^{\infty} \{a_{k,0}, \ldots, a_{k,2^k}\}, \label{FA.4.38.17} \\
C := \{ t \in [0,1] \, | \, f \mbox{ stetig in } t \} \label{FA.4.38.18}
\end{gather}
und behaupten
\begin{equation} \label{FA.4.38.19}
\big( \{ t \in [0,1] \, | \, g_1(t) = g_2(t) \} \cap \left( [0,1] \setminus R \right) \big) \subset C \subset \{ t \in [0,1] \, | \, g_1(t) = g_2(t) \}.
\end{equation}

{[} Zur ersten Inklusion in (\ref{FA.4.38.19}):
Sei $t_0 \in [0,1]$ mit $g_1(t_0) = g_2(t_0)$ -- d.h.\ in (\ref{FA.4.38.16}) gilt an der Stelle $t_0$ jeweils Gleichheit -- derart, daß $t_0$ nicht in der Menge der ,,Randpunkte`` $R = \left\{ a_{k,i} \, | \, k \in \N_+ \, \wedge \, i \in \left\{ \frac{0}{2^k}, \ldots, \frac{2^k}{2^k} \right\} \right\}$ enthalten ist, und es sei $\varepsilon \in \R_+$.
Wegen der Gleichheit in (\ref{FA.4.38.16}) an der Stelle $t_0$ folgt die Existenz einer Zahl $\kappa \in \N_+$ so, daß gilt
\begin{equation} \label{FA.4.38.20}
S_{\kappa}(t_0) - T_{\kappa}(t_0) < \varepsilon.
\end{equation}
Da $t_0 \notin R$ gilt, können wir wegen (\ref{FA.4.38.0}) eine Zahl $\iota \in \{ a_{\kappa,0}, \ldots,  a_{\kappa_0,2^{\kappa}-1} \}$ mit
\begin{equation*} \label{FA.4.38.21}
t_0 \in {]} a_{\kappa,\iota}, a_{\kappa,\iota + 1} {[}
\end{equation*}
finden, und es gilt sogar für jedes $t \in {]} a_{\kappa,\iota}, a_{\kappa,\iota + 1} {[}$ im Falle $f(t) \ge f(t_0)$
\begin{eqnarray*}
|f(t) - f(t_0)| & = & f(t) - f(t_0) \\
& \le & \sup \left\{ f(\tau) \, | \, \tau \in {]} a_{\kappa,\iota}, a_{\kappa,\iota+1} {[} \right\} \\
&& - \inf \left\{ f(\tau_0) \, | \, \tau_0 \in {]} a_{\kappa,\iota}, a_{\kappa,\iota+1} {[} \right\}
\end{eqnarray*}
sowie im Falle $f(t) < f(t_0)$
\begin{eqnarray*}
|f(t) - f(t_0)| & = & f(t_0) - f(t) \\
& \le & \sup \left\{ f(\tau_0) \, | \, \tau_0 \in {]} a_{\kappa,\iota}, a_{\kappa,\iota+1} {[} \right\} \\
&& - \inf \left\{ f(\tau) \, | \, \tau \in {]} a_{\kappa,\iota}, a_{\kappa,\iota+1} {[} \right\},
\end{eqnarray*}
d.h.\ nach (\ref{FA.4.38.3}), (\ref{FA.4.38.4}) stets für $t \in {]} a_{\kappa,\iota}, a_{\kappa,\iota + 1} {[}$
\begin{equation*} \label{FA.4.38.22}
|f(t) - f(t_0)| \le M_{\kappa,\iota} - m_{\kappa,\iota},
\end{equation*}
also wegen $t, t_0 \in {]} a_{\kappa,\iota}, a_{\kappa,\iota + 1} {[}$ und (\ref{FA.4.38.8}) sowie (\ref{FA.4.38.20})
\begin{equation*} \label{FA.4.38.23}
|f(t) - f(t_0)| \le S_{\kappa}(t_0) - T_{\kappa}(t_0) < \varepsilon.
\end{equation*}
Damit ist die Stetigkeit von $f$ in $t_0$ natürlich gezeigt.

Zur zweiten Inklusion in (\ref{FA.4.38.19}):
Seien $f$ stetig in $t_0 \in [0,1]$ und $\varepsilon \in \R_+$.
Dann existiert $\delta \in \R_+$ derart, daß gilt
\begin{equation} \label{FA.4.38.24}
\forall_{t \in [0,1]} \, \big( t \in \underbrace{{[} t_0 - \delta, t_0 + \delta {]}}_{\subset {]} t_0 - 2 \delta, t_0 + 2 \delta {[}} \Longrightarrow |f(t) - f(t_0)| < \frac{\varepsilon}{2} \big).
\end{equation}
Sei nun $\k_0 \in \N_+$ mit $(\frac{1}{2})^{k_0} < \delta$.
Dann kann nach (\ref{FA.4.38.0}), (\ref{FA.4.38.7}) zu jedem $k \in \N_+$ mit $k \ge k_0$ ein $i_k \in \{1, \ldots, 2^k\}$ mit
\begin{equation*} \label{FA.4.38.25}
t_0 \in Q_{k,i_k} \subset \left[ \frac{i_k - 1}{2^k}, \frac{i_k}{2^k} \right] \subset [ t_0 - \delta, t_0 + \delta ]
\end{equation*}
gewählt werden.
Es folgt mittels (\ref{FA.4.38.le}) sowie (\ref{FA.4.38.3}), (\ref{FA.4.38.4}) und des weiteren (\ref{FA.4.38.24})
\begin{eqnarray*}
0 & \le & S_k(t_0) - T_k(t_0) = M_{k,i_k} - m_{k,i_k} \\
& = & \sup \left\{ f(\tau) - f(t_0) \,| \, \tau \in \left[ \frac{i_k - 1}{2^k}, \frac{i_k}{2^k} \right] \right\} - \inf \left\{ f(\tau) - f(t_0) \, | \, \tau \in \left[ \frac{i_k - 1}{2^k}, \frac{i_k}{2^k} \right] \right\} \\
& \le & \frac{\varepsilon}{2} + \frac{\varepsilon}{2} = \varepsilon.
\end{eqnarray*}
Die Beliebigkeit von $\varepsilon \in \R_+$ ergibt $S_k(t_0) = T_k(t_0)$, also folgt aus (\ref{FA.4.38.11}) und (\ref{FA.4.38.13}): $g_1(t_0) = g_2(t_0)$.
Damit ist (\ref{FA.4.38.19}) vollständig gezeigt. {]}

,,(i) $\Rightarrow$ (ii), (iii)`` 
Sei also $f$ Riemann-integrierbar über $[0,1]$, d.h.\ per definitionem
\begin{equation} \label{FA.4.38.26}
\OIf = \UIf = \If.
\end{equation}
Dann folgt aus (\ref{FA.4.38.15})
\begin{equation*} \label{FA.4.38.27}
\int_{[0,1]} (g_2 - g_1) \, \d \mu_1 = 0,
\end{equation*}
also aus (\ref{FA.4.38.16}) sowie \ref{FA.4.30} ,,(i) $\Rightarrow$ (ii)``
\begin{equation*} \label{FA.4.38.28}
g_2 - g_1 =_{\mu_1} 0,
\end{equation*}
und dies impliziert erneut wegen (\ref{FA.4.38.16})
\begin{equation} \label{FA.4.38.29}
g_1 =_{\mu_1} f =_{\mu_1} g_2,
\end{equation}
m.a.W.\ ist $\{ t \in [0,1] \, | \, g_1(t) \ne g_2(t) \}$ eine $\mu_1$-Nullmenge.
Aus (\ref{FA.4.38.19}) folgt dann, daß 
$$ [0,1] \setminus C \subset \{ t \in [0,1] \, | \, g_1(t) \ne g_2(t) \} \cup R $$
wegen der Abzählbarkeit von $R$, vgl. (\ref{FA.4.38.17}), ebenfalls eine $\mu_1$-Nullmenge ist.
Hieraus ergibt sich nach der Definition von $C$ in (\ref{FA.4.38.18}) die Aussage (ii).

Weiterhin folgt aus (\ref{FA.4.38.11}), (\ref{FA.4.38.29}) und \ref{FA.4.24} (iv)
\begin{equation} \label{FA.4.38.30}
f \in \L_{\R}([0,1],\mu_1) ~~ \mbox{ sowie } ~~ \int_{[0,1]} f \, \d \mu_1 = \int_{[0,1]} g_1 \, \d \mu_1,
\end{equation}
also wegen (\ref{FA.4.38.12}), (\ref{FA.4.38.30}), (\ref{FA.4.38.26})
\begin{equation*} \label{FA.4.38.31}
\int_{[0,1]} f \, \d \mu_1 = \int_0^1 f(t) \, \d t,
\end{equation*}
womit auch (iii) gezeigt ist.

,,(ii) $\Rightarrow$ (i)`` Sei $[0,1] \setminus C$, wobei $C$ wie in (\ref{FA.4.38.18}) sei, eine $\mu_1$-Nullmenge.
Wegen (\ref{FA.4.38.19}) ist dann auch $\{ t \in [0,1] \, | \, g_1(t) \ne g_2(t) \}$ eine $\mu_1$-Nullmenge.
Aus (\ref{FA.4.38.11}) und (\ref{FA.4.38.13}) folgt dann mittels (\ref{FA.4.38.15}) die Riemann-Integrierbaitkeit von $f$ über $[0,1]$. \q

\begin{Bem*}
Leser, die mit der mehrdimensionalen Riemannschen Integrationstheorie, welche z.B.\ in \cite[Kapitel 12]{ElAna} dargestellt wird, vertraut sind, können nun die Ideen des letzten Beweises ausnutzen, um als Übung nachzuweisen, daß Hauptsatz \ref{FA.4.38} auch für beliebiges $n \in \N_+$ richtig ist, wenn man $[a,b]$ durch ein Kompaktum $Q \in \mathfrak{Q}_n$, $\mu_1$ durch $\mu_n$ und das eindimensionale Riemann-Integral von $f$ über $[a,b]$ durch das mehrdimensionale von $f$ über $Q$ ersetzt.
\end{Bem*}

\begin{Bsp*} \index{Integral!Riemann-} \label{FA.4.39}
$\mu_1$-Integrierbarkeit einer reellen Funktion über einem kompakten Intervall impliziert nicht die Riemann-Integrierbarkeit der Funktion über dem Intervall:
Die Funktion $f \: [0,1] \to \R$, die durch
$$ \forall_{t \in [0,1]} \, f(t) := \begin{cases} 0, & \mbox{falls } t \in [0,1] \cap \Q, \\ 1, & \mbox{falls } t \in [0,1] \setminus \Q, \end{cases} $$
gegeben ist, ist $\mu_1$-integriebar über $[0,1]$ -- baechte, daß $[0,1] \cap \Q$ eine $\mu_1$-Nullmenge ist.
$f$ ist jedoch nicht Riemann-in\-te\-grier\-bar über $[0,1]$, da Ober- und Untersumme von $f$ sich für jede Zerlegung von $[0,1]$ unterscheiden.
\end{Bsp*}

\begin{Bem}[Uneigentliche Riemann-Integrale] \index{Integral!Riemann-} \label{FA.4.40}
Es ist im folgenden wichtig, zu betonen, daß wir unter einem \emph{uneigentlichen Riemann-Integral} nur ein solches im Sinne von \cite[6.36]{Henke} und nicht allgemeiner im Sinne von \cite[8.33]{ElAna} verstehen.
D.h.\ genau, daß ein uneigentliches Riemann-Integral für uns per definitionem nur Werte in $\R$ (und nicht etwa in $\widehat{\R} = \R \cup \{ \pm \infty \}$) annimmt.
\end{Bem}

\begin{HS} \index{Integral!Riemann-} \label{FA.4.41} $\,$

\noindent \textbf{Vor.:} Seien $J \in \mathfrak{I}_1$ ein nicht-kompaktes Intervall von $\R$ und $f \: J \to \R$ eine über $J$ uneigentlich Riemann-integrierbare Funktion, d.h.\ insbesondere, daß $f$ über jedem kompakten Teilintervall von $J$ (eigentlich) Riemann-integrierbar ist.

\pagebreak
\noindent \textbf{Beh.:} Die Aussagen (i) und (ii) sind äquivalent:
\begin{itemize}
\item[(i)] $f$ ist $\mu_1$-integrierbar über $J$.
\item[(ii)] $|f|$ ist über $J$ uneigentlich Riemann-integrierbar.
\end{itemize}

Des weiteren folgt aus (i) -- oder äquivalent aus (ii) --
\begin{itemize}
\item[(iii)] Das $\mu_1$-Integral von $f$ über $J$ und das uneigentliche Riemann-Integral von $f$ über $J$ stimmen überein.
\end{itemize}
\end{HS}

Wir bereiten den Beweis des Hauptsatzes durch folgendes Lemma vor.

\begin{Lemma} \label{FA.4.42} $\,$

\noindent \textbf{Vor.:} Seien $J \in \mathfrak{I}_1$ ein nicht-kompaktes Intervall von $\R$, $a := \inf J \in \R \cup \{- \infty\}$ und $b := \sup J \in \R \cup \{+ \infty\}$.
Ferner sei $f \: J \to \R$ eine reellwertige Funktion, deren Beschränkung auf jedes kompakte Teilintervall von $J$ (eigentlich) Riemann-integrierbar ist.

\noindent \textbf{Beh.:} Die folgenden Aussagen sind äquivalent:
\begin{itemize}
\item[(i)] $f$ ist über $J$ uneigentlich Riemann-integrierbar.
\item[(ii)] Die Folge $\D \left( \int_{a_i}^{b_i} f(t) \, \d t \right)_{i \in \N}$ ist für alle Folgen $([a_i, b_i])_{i \in \N}$ kompakter Teilintervalle von $J$ mit $\lim_{i \to \infty} a_i = a $ bzw.\ $\lim_{i \to \infty} b_i = b$, wobei $\forall_{i \in \N} \, a_i = a$, falls $a \in J$, und $\forall_{i \in \N} \, b_i = b$, falls $b \in J$, gelte, in $\R$ konvergent ist -- beachte, daß $\D \int_{a_i}^{b_i} f(t) \, \d t$ nach Voraussetzung wohldefiniert ist --, und für das uneigentliche Riemann-Integral auf der linken Seite der folgenden Gleichung gilt
\begin{equation} \label{FA.4.42.0}
\int_a^b f(t) \, \d t = \lim_{i \to \infty} \int_{a_i}^{b_i} f(t) \, \d t.
\end{equation}
\end{itemize}
\end{Lemma}

\textit{Beweis.} ,,(i) $\Rightarrow$ (ii)`` Sei zunächst $f$ uneigentlich Riemann-integrierbar über $J$ mit uneigentlichem Riemann-Integral $\int_a^b f(t) \, \d t$, d.h.\ unter Beachtung von \ref{FA.4.40} z.B.\ nach \cite[8.33]{ElAna}, daß es eine (nicht eindeutige) Zahl $c \in {]}a,b{[}$ gibt, derart, daß gilt
\begin{equation} \label{FA.4.42.1}
\lim_{i \to \infty} \int_{a_i}^c f(t) \, \d t ~~ \mbox{ und } ~~ \lim_{i \to \infty} \int_c^{b_i} f(t) \, \d t ~~ \mbox{ existieren in $\R$,}
\end{equation}
wobei $([a_i, b_i])_{i \in \N}$ eine beliebige Folge mit der (ii) genannten Eigenschaften sei, und die Limites in (\ref{FA.4.42.1}) hängen nicht von der speziellen Wahl der soeben genannten Folge ab.
Das uneigentliche Riemann-Integral von $f$ über $J$ auf der linken Seite der folgenden Gleichung erfüllt per definitionem
\begin{equation} \label{FA.4.42.2}
\int_a^b f(t) \, \d t = \lim_{i \to \infty} \int_{a_i}^c f(t) \, \d t + \lim_{i \to \infty} \int_c^{b_i} f(t) \, \d t,
\end{equation}
Die Beliebigkeit der genannten Folge und (\ref{FA.4.42.2}) ergeben die Konvergenz von $\left( \int_{a_i}^{b_i} f(t) \, \d t \right)_{i \in \N}$ in $\R$ sowie (\ref{FA.4.42.0}).

,,(ii) $\Rightarrow$ (i)`` Gelte nun für jede Folge $([a_i,b_i])_{i \in \N}$ mit den in (ii) genannten Eigenschaften, daß $\left( \int_{a_i}^{b_i} f(t) \, \d t \right)_{i \in \N}$ in $\R$ konvergiert sowie (\ref{FA.4.42.0}).
Für beliebiges $c \in {]}a,b{[}$ folgt offenbar (\ref{FA.4.42.1}), und die Limites hängen wegen (\ref{FA.4.42.0}) nicht von der Wahl der Folge $([a_i,b_i])_{i \in \N}$ ab.
Aus der Beliebigkeit letztgenannter Folge ergibt sich, daß $f$ uneigentlich Riemann-integrierbar über $J$ ist.
\q
\A
\textit{Beweis des Hauptsatzes.} 
Seien $a := \inf J \in \R \cup \{- \infty\}$, $b := \sup J \in \R \cup \{+ \infty\}$ und $\left( [a_i,b_i] \right)_{i \in \N}$ eine beliebige Folge wie in (ii) des letzten Lemmas.

Wir zeigen zunächst
\begin{equation} \label{FA.4.41.1}
f \ge 0 \Longrightarrow f \in \L_{\R}(J,\mu_1).
\end{equation}

{[} Zu (\ref{FA.4.41.1}): 
Es gelte die zusätliche Eigenschaft
\begin{equation} \label{FA.4.41.2}
\forall_{i \in \N} \, a_{i+1} \le a_i \, \wedge \, b_i \le b_{i+1}.
\end{equation}

Aus $f \ge 0$ und (\ref{FA.4.41.2}) folgt, daß $\left( \chi_{[a_i,b_i]} \, \hat{f} \right)_{i \in \N}$ eine monton wachsende Folge, deren Folgenglieder nach Voraussetzung und \ref{FA.4.38} ,,(i) $\Rightarrow$ (iii)`` Elemente von $\L_{\R}(\R,\mu_1)$ sind, ist, und es gilt -- beachte \ref{FA.4.42} (i) $\Rightarrow$ (ii) --
$$ \forall_{i \in \N} \, \int_{\R} \chi_{[a_i,b_i]} \, \hat{f} \, \d \mu_1 = \int_{[a_i,b_i]} f \, \d \mu_1 = \int_{a_i}^{b_i} f(t) \, \d t \le \int_a^b f(t) \, \d t \in \R. $$
Daher ergibt der Grenzwertsatz von \textsc{Levi}: $\chi_J \, \hat{f} = \lim_{i \to \infty} \chi_{[a_i,b_i]} \, \hat{f} \in \L_{\R}(\R,\mu_1)$, d.h.\ genau $f \in \L_{\R}(J,\mu_1)$. {]}

,,(i) $\Rightarrow$ (ii), (iii)`` Sei $f$ $\mu_1$-in\-te\-grier\-bar über $J$.
Dann ist auch $|f|$ $\mu_1$-in\-te\-grier\-bar über $J$, d.h.\ genau $\chi_J \, \widehat{|f|} \in \L_{\R}(\R,\mu_1)$.
Somit ergibt \ref{FA.4.20} (vii)
\begin{equation*} \label{FA.4.41.3}
\forall_{i \in \N} \, \chi_{[a_i,b_i]} \, \widehat{|f|} \in \L_{\R}(\R,\mu_1), \mbox{ d.h.\ } |f| \in \L_{\R}([a_i,b_i)],\mu_1),
\end{equation*}
und es gilt
\begin{equation*} \label{FA.4.41.4}
\forall_{i \in \N} \, \left| \chi_{[a_i,b_i]} \, \widehat{|f|} \right| \le \chi_J \, \widehat{|f|} \in \L_{\R}(\R,\mu_1).
\end{equation*}
Daher ergibt der Grenzwertsatz von \textsc{Lebesgue} wegen $\lim_{i \to \infty} \chi_{[a_i,b_i]} \, \widehat{|f|} = \chi_J \, \widehat{|f|}$ die Existenz der linken Seite der folgenden Gleichung in $\R$ und
$$ \int_{\R} \chi_J \, \widehat{|f|} \,\, \d \mu_1 = \lim_{i \to \infty} \int_{[a_i,b_i]} |f| \,\, \d \mu_1. $$
Nun ist $|f|$ für jedes $i \in \N$ mit $f$ über $[a_i,b_i]$ -- nach Voraussetzung -- (eigentlich) Rie\-mann-in\-te\-grier\-bar, vgl.\ \cite[8.33 (i) und Satz 8.13]{ElAna}, also folgt aus \ref{FA.4.38}\linebreak ,,(i) $\Rightarrow$ (iii)``
$$ \lim_{i \to \infty} \int_{a_i}^{b_i} |f|(t) \, \d t = \lim_{i \to \infty} \int_{[a_i,b_i]} |f| \,\, \d \mu_1 \in \R. $$
Die Beliebigkeit der Folge $([a_i,b_i])_{i \in \N}$ ergibt zusammen mit \ref{FA.4.42} ,,(ii) $\Rightarrow$ (i)`` die uneigentliche Riemann-Integrierbarkeit von $|f|$ über $J$, d.h.\ es gilt (ii).

Des weiteren folgt aus der $\mu_1$-Integrierbarkeit von $f$ über $J$ mittels \ref{FA.4.20} (vii)
\begin{equation*} \label{FA.4.41.5}
\forall_{i \in \N} \, \chi_{[a_i,b_i]} \, f \in \L_{\R}(\R,\mu_1), \mbox{ d.h.\ } f \in \L_{\R}([a_i,b_i],\mu_1),
\end{equation*}
und es gilt offenbar
\begin{gather*} 
\forall_{i \in \N} \, -\chi_J \, \widehat{|f|} \le \chi_{[a_i,b_i]} \, \hat{f} \le \chi_J \, \widehat{|f|}, \label{FA.4.41.6} \\
\forall_{i \in \N} \, \left| \chi_{[a_i,b_i]} \, \hat{f} \right| \le \chi_ J \, \widehat{|f|} \stackrel{\text{s.o.}}{\in} \L_{\R}(\R,\mu_1). \label{FA.4.41.7}
\end{gather*}
\pagebreak

\noindent Aus dem Grenzwertsatz von \textsc{Lebesgue} folgt wegen $\lim \chi_{[a_i,b_i]} \, \hat{f} = \chi_J \, \hat{f}$ und \ref{FA.4.38} ,,(i) $\Rightarrow$ (iii)`` sowie \ref{FA.4.42} ,,(i) $\Rightarrow$ (ii)``
$$ \int_J f \, \d \mu = \int_{\R^n} \chi_J \, f \, \d \mu_1 = \lim_{i \to \infty} \int_{a_i}^{b_i} f(t) \, \d t  = \int_a^b f(t) \, \d t, $$
also gilt auch (iii).

,,(ii) $\Rightarrow$ (i)`` Sei $|f|$ uneigentlich Riemann-integrierbar über $J$.
Aus (\ref{FA.4.41.1}), angewandt auf $|f|$ anstelle von $f$, ergibt sich $|f| \in \L_{\R}(J,\varphi)$, d.h.\ $\chi_J \, \widehat{|f|} \in \L_{\R}(\R,\varphi)$, also 
\begin{equation*} \label{FA.4.41.8}
\forall_{i \in \N} \, \left| \chi_{[a_i,b_i]} \, \hat{f} \right| \le \chi_J \, \widehat{|f|} \in \L_{\R}(\R,\varphi).
\end{equation*}
Hieraus sowie $\lim_{i \to \infty} \chi_{[a_i,b_i]} \, \hat{f} = \chi_J \, \hat{f}$ und dem Grenzwertsatz von \textsc{Lebesgue} folgt (i). \q

\begin{Bsp*} \label{FA.4.43}
Die Funktion $f \: {[} 1, \infty {[} \to \R$, definiert durch 
$$ \forall_{t \in {[} 1, \infty {[}} \, f(t) := \frac{\sin(t)}{t}, $$
ist uneigentlich Riemann-integrierbar über ${[} 1, \infty {[}$, jedoch nicht $\mu_1$-integrierbar über ${[} 1, \infty {[}$.
\end{Bsp*}

\begin{Def} \label{FA.4.44}
Seien $J \in \mathfrak{I}_1$ mit $J \ne \emptyset$ und $a := \inf J \in \R \cup \{- \infty \}$ sowie $b := \sup J \in \R \cup \{ + \infty \}$.
Ferner sei $f \: J \to \R$ eine über $J$ $\mu_1$-integrierbare Funktion.
Wegen \ref{FA.4.38} ,,(i) $\Rightarrow$ (iii)`` und \ref{FA.4.41} ,,(i) $\Rightarrow$ (iii)`` können wir folgende Notation verwenden
$$ \boxed{\int_a^b f(x) \, \d x} := \int_J f \, \d \mu_1, $$
ohne daß Verwechselungen mit dem Riemann-Integral möglich sind.
\end{Def}

\begin{Satz}[Interpretation der Theorie der gewichteten absolut konvergenten Reihen als Integralrechnung] \index{Reihe!(absolut) konvergente} \label{FA.4.45} $\,$

\noindent \textbf{Vor.:} Seien $M \subset \R^n$ eine nicht-leere Menge ohne Häufungspunkt in $\R^n$, insbesondre ist $M$ -- wie in Beispiel \ref{FA.4.4.B} 4.) erwähnt --, höchstens abzählbar.
Daher existieren $k \in \N_+ \cup \{\infty\}$ und paarweise verschiedene $p_i \in \R^n$, $i \in I := \{ i \in \N \, | \, i \le k \}$, mit $M = \{ p_i \, | \, i \in I \}.$
$m \: M \to \R$ sei ferner eine positivwertige Funktion und $\varphi_m$ das Quadermaß der diskreten Massenverteilung $m$. 
Außerdem seien $D$ eine Teilmenge von $\R^n$ und $f \in \K^D$.

\noindent \textbf{Beh.:} Es gilt
\begin{itemize}
\item[(i)] $\D f \in \mathcal{L}_{\K}(\R^n,\varphi_m) \Longleftrightarrow \left( M \subset D \, \wedge \,\sum_{i=0}^k m(p_i) \, f(p_i) \mbox{ absolut konvergent} \right),$
\item[(ii)] $\D f \in \mathcal{L}_{\K}(\R^n,\varphi_m) \Longrightarrow \left( M \subset D \, \wedge \, \int_{\R^n} f \, \d \varphi_m = \sum_{i=0}^k m(p_i) \, f(p_i) \right).$
\end{itemize}
\end{Satz}

\pagebreak
\textit{Beweis.} 1.) Sei zunächst $f \in \mathcal{L}_{\R}(\R^n,\varphi_m)$ mit $f \ge 0$.
Wir definieren
\begin{gather*}
T_0 := f(p_0) \, \chi_{\{p_0\}}, \label{FA.4.45.1} \\
\forall_{i \in \N} \, \left( 0 < i \le k \Longrightarrow T_i := f(p_i) \, \chi_{\{p_i\}} + T_{i-1} \right), \label{FA.4.45.2} \\
\forall_{i \in \N} \, \left( i > k \Longrightarrow T_i := T_{i-1} \right). \label{FA.4.45.3}
\end{gather*}
Dann ist $(T_i)_{i \in \N}$ eine monoton wachsende Folge in $\mathcal{S}(\R^n)$ mit
\begin{equation*} \label{FA.4.45.4}
\lim_{i \to \infty} T_i =_{\varphi_m} f,
\end{equation*}
beachte, daß $\R^n \setminus M$ eine $\varphi_m$-Nullmenge ist, und
\begin{equation*} \label{FA.4.45.5}
\forall_{i \in \N} \, \int_{\R^n} T_i \, \d \varphi_m \le \int_{\R^n} f \, \d \varphi_m.
\end{equation*}
Der Grenzwertsatz von \textsc{Levi} ergibt nun wegen $\forall_{i \in I} \, \int_{\R^n} \chi_{\{p_i\}} \d \varphi_m = m(p_i)$ die rechte Seite von (ii).

2.) Aus $f \in \L_{\K}(M,\varphi_m)$ folgt $|f| \in \L_{\R}(M,\varphi_m)$.
Daher ergibt 1.) die Konklusion von (i).

3.) Aus der Gültigkeit der rechten Seite von (i) folgt sowohl die absolute Konvergenz des Real- als auch des Imaginärteiles der Reihe $\sum_{i=0}^k m(p_i) \, f(p_i)$.
Daher können wir $\K = \R$ annehmen, und es genügt zu zeigen, daß gilt
$$ f^+, f^- \in \L_{\R}(\R^n,\varphi_m). $$
Wir zeigen dies nur für $f^+$, der Beweis der Aussage für $f^-$ verläuft ebenso.

Sei $I_+ := \{ i \in I \, | \, f(p_i) \ge 0 \}$.
Dann konvergiert auch $\sum_{\substack{i=0 \\ i \in I_+}}^k m(p_i) \, f(p_i)$ absolut, d.h.\ jede endliche Teilsumme besitzt eine obere Schranke.
Eine zu 1.) analoge Argumentation ergibt nun mittels des Grenzwertsatzes von \textsc{Levi}, daß $f^+ \in \L_{\R}(\R^n,\varphi_m)$ gilt.

4.) Sei $f \in \L_{\K}(M,\varphi_m)$.
Wegen 2.) hängt die Summe auf der rechten Seite von (ii) nicht von der Reihenfolge der Summanden ab.
Wir können zunächst wieder ohne Beschränkung der Allgemeinheit $\K = \R$ und wegen $f = f^+ - f^-$ sodann $f \ge 0$ annehmen.
1.) ergibt daher die Aussage (ii).

Mit 1.) - 4.) ist der Satz bewiesen. \q

\begin{Def}[Stieltjes-Integral] \index{Quader!-maß!Stieltjessches} \index{Integral!Stieltjes-} \label{FA.4.46}
Seien $g \: \R \to \R$ eine monoton wachsende Funktion und $\varphi_g$ das gemäß Beispiel \ref{FA.4.4.B} 5.) definierte Stieltjessche Maß.
Weiterhin seien $J$ ein nicht-leeres Intervall von $\R$, $a := \inf J$, $b := \sup J$ und $f \in \L_{\K}(J,\varphi_g)$.
Wir setzen
\begin{equation*} \label{FA.4.46.0}
\boxed{\int_a^b f(x) \, \d \, g(x)} := \int_J f \, \d \, \varphi_g 
\end{equation*}
und nennen dieses Integral das \emph{Stieltjes-Integral von $f$ über $J$ (mit der Gewichtsfunktion $g$)}.

Man beachte, daß im Falle $g = \id_{\R} = x \: \R \to \R$ das $\varphi_g$-Integral mit dem $\mu_1$-Integral übereinstimmt.
\end{Def}

\pagebreak

\begin{Satz} \index{Integral!Stieltjes-} \label{FA.4.47} $\,$

\noindent \textbf{Vor.:} Seien $g \: \R \to \R$ eine monoton wachsende Funktion und $\varphi_g$ das gemäß Beispiel \ref{FA.4.4.B} 5.) definierte Stieltjessche Maß.
Ferner seien $a,b \in \R$ mit $a<b$ und $f \: [a,b] \to \R$ eine über $[a,b]$ $\varphi_g$-integrierbare stetige Funktion.

\noindent \textbf{Beh.:} Ist $g$ stetig differenzierbar, so gilt
\begin{equation*} \label{FA.4.47.0}
\int_a^b f(x) \, \d \, g(x) = \int_a^b f(t) \, g'(t) \, \d t,
\end{equation*}
wobei die rechte Seite der Gleichung als stetige Funktion über dem Kompaktum $[a,b]$ als Riemann-Integral existiert. \q

\begin{Bem*}
Wir werden in \ref{FA.4.57} zeigen, daß aus der Stetigkeit von $f$ auf dem Kompaktum $[a,b]$ bereits die $\varphi_g$-Integrierbarkeit von $f$ über $[a,b]$ und aus der Stetigkeit von $f \, g'$ auf $[a,b]$ die $\mu_1$-Integrierbarkeit von $f \, g'$ über $[a,b]$ folgt.
\end{Bem*}
\end{Satz}

\textit{Beweisskizze.} Ohne Einschränkung können wir $a=0$ und $b=1$ annehmen.
Für alle $k \in \N_+$ und $i \in \{0, \ldots, 2^k\}$ seien dann $a_{k,i}$, $m_{k,i}$, $Q_{k,i}$ und $T_k$ wie in (\ref{FA.4.38.0}), (\ref{FA.4.38.3}), (\ref{FA.4.38.7}) und (\ref{FA.4.38.8}) definiert.
Dann ist $(T_k)_{k \in \N}$ eine monoton wachsende Folge in $\mathcal{S}([0,1])$ mit $\lim_{k \to \infty} T_k = f$ sowie
$$ \forall_{k \in \N} \, \int_0^1 T_k(x) \, \d \, g(x) \le \int_0^1 f(x) \, \d \, g(x) \in \R. $$
Der Grenzwertsatz von \textsc{Levi} ergibt daher
\begin{equation} \label{FA.4.47.1}
\int_0^1 f(x) \, \d \, g(x) = \lim_{k \to \infty} \int_0^1 T_k(x) \, \d \, g(x),
\end{equation}
und es gilt für jedes $k \in \N_+$ wegen der Stetigkeit von $g$
$$ \int_0^1 T_k(x) \, \d \, g(x) = \sum_{i=1}^{2^k} m_{k,i} \left( g(a_{k,i}) - g(a_{k,i-1}) \right), $$
also folgt aus der Differenzierbarkeit von $g$
$$ \int_0^1 T_k(x) \, \d \, g(x) = \sum_{i=1}^{2^k} m_{k,i} \int_{a_{k,i}}^{a_{k,i-1}} g'(t) \, \d t, $$
wobei die rechte Seite für $k \to \infty$ gegen $\D \int_0^1 f(t) \, g'(t) \, \d t$ konvergiert.
Zusammen mit (\ref{FA.4.47.1}) ergibt dies die Behauptung. \q

\subsection*{Lebesgue-integrierbare Mengen} \addcontentsline{toc}{subsection}{Lebesgue-integrierbare Mengen}

\begin{Def} \index{Menge!integrierbare} \label{FA.4.48} 
Es seien $\varphi$ ein Quadermaß auf $\R^n$ und $M$ eine beliebige Teilmenge von $\R^n$.

$M$ heißt \emph{$\varphi$-integrierbar} (oder auch \emph{$\varphi$-summierbar}) genau dann, wenn gilt $\chi_M \in \L_{\R}(\R^n,\varphi)$.
In diesem Falle nennen wir
$$ \boxed{\varphi(M)} := \int_{\R^n} \chi_M \, \d \varphi = \int_M 1_{\R^n} \d \varphi$$
das \emph{$\varphi$-Maß von $M$}.

\begin{Bem*}
Die hier eingeführte Notation ist im Falle $M \in \mathfrak{P}_n$ offenbar mit der vorherigen konsistent.
\end{Bem*}
\end{Def}

Korollar \ref{FA.4.30.K} besagt genau das folgende Resultat.

\begin{Kor} \label{FA.4.30.K2} 
Sind $\varphi$ ein Quadermaß auf $\R^n$ und $M$ eine Teilmenge von $\R^n$, so ist $M$ genau dann eine $\varphi$-Nullmenge, wenn $M$ $\varphi$-integrierbar mit $\varphi(M) = 0$ ist. \q
\end{Kor}

\begin{Satz} \label{FA.4.50} 
Es sei $\varphi$ ein Quadermaß auf $\R^n$.

\begin{itemize}
\item[(i)] Seien $M_1, M_2$ zwei $\varphi$-integrierbare Teilmengen von $\R^n$.

Dann sind $M_1 \cup M_2$, $M_1 \cap M_2$ sowie $M_1 \setminus M_2$ $\varphi$-integrierbar, und es gilt
\begin{gather*}
\varphi(M_1 \cup M_2) = \varphi(M_1) + \varphi(M_2) - \varphi(M_1 \cap M_2), \\
\varphi(M_1 \setminus M_2) = \varphi(M_1) - \varphi(M_1 \cap M_2).
\end{gather*}
\item[(ii)] Für jedes $i \in \N$ sei $M_i$ eine $\varphi$-integrierbare Teilmenge von $\R^n$, und es gelte $\forall_{i \in \N} \, M_i \subset M_{i+1}$.

Dann folgt für $\D M := \bigcup_{i \in \N} M_i$:
\begin{itemize}
\item[(a)] $M$ $\varphi$-integrierbar $\Longrightarrow$ $\D \varphi(M) = \lim_{i \to \infty} \varphi(M_i)$.
\item[(b)] $\D \lim_{i \to \infty} \varphi(M_i) \in \R$ $\Longrightarrow$ $M$ $\varphi$-integrierbar und $\D \varphi(M) = \lim_{i \to \infty} \varphi(M_i)$.
\end{itemize}
\end{itemize}
\end{Satz}

\textit{Beweis.} $\chi_{M_1 \cup M_2} = \sup(\chi_{M_1}, \chi_{M_2}), \, \chi_{M_1 \cap M_2} = \inf(\chi_{M_1}, \chi_{M_2}) \in \L_{\R}(\R^n,\varphi)$, d.h.\ $\chi_{M_1 \setminus M_2} = \chi_{M_1} - \chi_{M_1 \cap M_2} \in \L_{\R}(\R^n,\varphi)$, sowie $\chi_{M_1 \cup M_2} = \chi_{M_1} + \chi_{M_2} - \chi_{M_1 \cap M_2}$ ergeben (i).

(ii) folgt aus \ref{FA.4.61} mit $f := 1_{\R^n}$. \q

\begin{Bsp} \label{FA.4.49} 
Es existiert eine Teilmenge der $\mu_1$-integrierbaren Menge $[0,1]$, die nicht $\mu_1$-integrierbar ist:

Auf $[0,1]$ wird durch
$$ \forall_{t,s \in [0,1]} \, t \sim s : \Longleftrightarrow t-s \in \Q $$
eine Äquivalenzrelation $\sim$ definiert.
$[0,1]_{\sim}$ bezeichne die Menge der Äquivalenzklassen bzgl.\ $\sim$.
Nach dem Auswahlaxiom existiert eine Abbildung
$$ h \: [0,1]_{\sim} \longrightarrow [0,1] $$
mit $\forall_{t \in [0,1]} \, h([t]_{\sim}) \in [t]_{\sim}$.
Setze nun
$$ M := h \left( [0,1]_{\sim} \right), $$
also gilt
\begin{equation} \label{FA.4.49.0}
\R = M \oplus \Q.
\end{equation}

Die Teilmenge $M$ von $[0,1]$ ist nicht $\mu_1$-integrierbar, denn andernfalls wäre für jedes $k \in \N_+$ offenbar auch
$$ M_k := \bigcupdot_{\kappa=1}^k \left( M + \frac{1}{\kappa} \right) \subset [0,2] $$
$\mu_1$-integrierbar mit $\mu_1(M_k) = k \cdot \mu_1(M) \le 2$.\footnote{Hier geht die leicht zu beweisende \emph{$\mu_1$-Translationsinvarianz}
$$ \forall_{M \in \mathfrak{P}(\R)} \forall_{a \in \R} \, M \mbox{ $\mu_1$-integrierbar} \Longrightarrow M + a \mbox{ $\mu_1$-integrierbar und } \mu_1(M + a) = \mu_1(M) $$
ein, die allgemeiner auch für den Beweis von Hauptsatz \ref{FA.4.67} benötigt wird und ein Spezialfall dieses Hauptsatzes ist.}
Wegen $\lim_{k \to \infty} k = \infty$ folgte hieraus, daß $M$ eine $\mu_1$-Nullmenge sei.
Dies widerspäche (\ref{FA.4.49.0}).
\end{Bsp}

\begin{Satz} \label{FA.4.51} 
Seien $\varphi$ ein Quadermaß auf $\R^n$ sowie $U$ eine offene und beschränkte Teilmenge von $\R^n$.

Dann ist $U$ $\varphi$-integrierbar, und es gilt
$$ \varphi(U) = \sup \{ \varphi(P) \, | \, P \in \mathfrak{P}_n \, \wedge \, P \subset U \}. $$
\end{Satz}

\begin{Def} \label{FA.4.52} 
Sind $\varphi$ ein Quadermaß auf $\R^n$ und $M$ eine Teimenge von $\R^n$, so heißt
$$ \boxed{\sup \{ \varphi(P) \, | \, P \in \mathfrak{P}_n \, \wedge \, P \subset M \}} \in [0, \infty ] $$
das \emph{innere Jordan-$\varphi$-Maß von $M$}.
\end{Def}

Wir bereiten den Beweis des Satzes \ref{FA.4.51} durch ein Lemma vor.

\begin{Lemma} \label{FA.4.53} 
Sei $U$ eine offene Teilmenge von $\R^n$.

Dann existiert eine Folge abgeschlossener Quader $(Q_i)_{i \in \N}$ in $\mathfrak{Q}_n$ derart, daß $\D U = \bigcup_{i=0}^{\infty} Q_i$ gilt.
\end{Lemma}

\textit{Beweisskizze.} Für $k \in \N$ definieren wir $\mathcal{W}_k := \left\{ W_{k,a} \, | \, a \in \Z^n \right\}$, wobei
$$ \forall_{a = (a_1, \ldots, a_n) \in \Z^n} \, W_{k,a} := \left[ \frac{a_1}{2^k}, \frac{a_1 + 1}{2^k} \right] \times \ldots \times \left[ \frac{a_n}{2^k}, \frac{a_n + 1}{2^k} \right] $$
sei, und
$$ P_k := \bigcup_{W \in \mathcal{W}_k, \, W \subset U \cap \overline{U_k^{\| \ldots \|_{\infty}}(0)}} W. $$
Dann gilt $U = \bigcup_{k \in \N_+} P_k$. \q
\A
\textit{Beweisskizze des Satzes \ref{FA.4.51}.} Sei $(Q_i)_{i \in \N}$ wie im vorherigen Lemma.
Die Anwendung des Satzes von \textsc{Levi} auf die Folge $(T_k)_{k \in \N}$, gegeben durch
$$ \forall_{k \in \N} \, T_k := \chi_{Q_0 \cup \ldots \cup Q_k}, $$
ergibt die $\varphi$-Integriebarkeit von $U$ und
\begin{equation} \label{FA.4.51.1}
\varphi(U) = \lim_{k \to \infty} \varphi(Q_0 \cup \ldots \cup Q_k),
\end{equation}
beachte, daß die Beschränktheit von $U$ die Existenz eines Quaders $Q \in \mathfrak{Q}_n$ mit $U \subset Q$ liefert, also $\int_{\R^n} T_k \, \d \varphi \le \varphi(Q)$ für alle $k \in \N$.
Aus der $\varphi$-Integrierbarkeit von $U$ folgt
$$ \sup \{ \varphi(P) \, | \, P \in \mathfrak{P}_n \, \wedge \, P \subset U \} \le \varphi(U). $$

Angenommen, in der letzten Ungleichung gilt ,,$<$``.
Dann existiert $\varepsilon \in \R_+$ derart, daß für alle $P \in \mathfrak{P}_n$ mit $P \subset U$ gilt $\varphi(P) \le \varphi(U) - \varepsilon$, also $|\varphi(U) - \varphi(P)| \ge \varepsilon$, im Widerspruch zu (\ref{FA.4.51.1}). \q

\begin{Satz} \label{FA.4.54} 
Es seien $\varphi$ ein Quadermaß auf $\R^n$ sowie $K$ eine kompakte Teilmenge von $\R^n$.

Dann ist $K$ $\varphi$-integrierbar, und es gilt
$$ \varphi(K) = \inf \{ \varphi(P) \, | \, P \in \mathfrak{P}_n \, \wedge \, K \subset P \}. $$
\end{Satz}

\begin{Def} \label{FA.4.55} 
Sind $\varphi$ ein Quadermaß auf $\R^n$ und $M$ eine beschränkte Teilmenge von $\R^n$, so heißt
$$ \boxed{\inf \{ \varphi(P) \, | \, P \in \mathfrak{P}_n \, \wedge \, M \subset P \}} \in {[} 0, \infty {[} $$
das \emph{äußere Jordan-$\varphi$-Maß von $M$}.
\end{Def}

\textit{Beweis des Satzes \ref{FA.4.54}.} Nach der Wahl eines offenen Quaders $Q \in \mathfrak{Q}_n$ mit $K \subset Q$ liefert die Anwendung von \ref{FA.4.51} auf die offenen und beschränkten Mengen $Q \setminus K = Q \cap (\R^n \setminus K)$ sowie $Q$ die $\varphi$-Integrierbarkeit von $Q \setminus K$ sowie $Q$ und
$$ \varphi(Q \setminus K) = \sup \{ \varphi(\widetilde{P}) \, | \, \widetilde{P} \in \mathfrak{P}_n \, \wedge \, \widetilde{P} \subset Q \setminus K \}. $$
Dann ist nach \ref{FA.4.50} (i) auch $K = Q \setminus (Q \setminus K)$ eine $\varphi$-integrierbare Menge mit $\varphi$-Maß
\begin{eqnarray*}
\varphi(K) & = & \varphi(Q) - \varphi(Q \setminus K) \\
& = & \varphi(Q) - \sup \{ \varphi(\widetilde{P}) \, | \, \widetilde{P} \in \mathfrak{P}_n \, \wedge \, \widetilde{P} \subset Q \setminus K \} \\
& = & \inf \{ \varphi(Q) - \varphi(\widetilde{P}) \, | \, \widetilde{P} \in \mathfrak{P}_n \, \wedge \, \widetilde{P} \subset Q \setminus K \} \\
& = & \inf \{ \varphi(\underbrace{\underbrace{Q \setminus \widetilde{P}}_{=: P \in \mathfrak{P}_n}}_{\stackrel{\widetilde{P} \subset Q}{\Rightarrow} \widetilde{P} = Q \setminus P \in \mathfrak{P}_n}) \, | \, \widetilde{P} \in \mathfrak{P}_n \, \wedge \, \underbrace{\widetilde{P} \subset Q \setminus K}_{\Leftrightarrow K \subset P} \} \\
& = & \inf \{ \varphi(P) \, | \, P \in \mathfrak{P}_n \, \wedge \, K \subset P \}.
\end{eqnarray*}
\q

\begin{Bsp} \label{FA.4.56} 
Die \emph{Smith-Volterra-Cantor Menge}\index{Menge!Smith-Volterra-Cantor-} entsteht aus $[0,1]$, indem man im ersten Schritt das in der Mitte liegende offene Intervall ${]} \frac{3}{8} , \frac{5}{8}{[}$ der Länge $\frac{1}{4}$ entfernt, im zweiten Schritt aus jedem verbleibenden Intervall das in der Mitte liegende offene Intervall der Länge $\frac{1}{16} = \left( \frac{1}{4} \right)^2$ entfernt, im $k$-ten Schritt für $k \in \N_+$ aus jedem verbliebenen Intervall das in der Mitte liegende offene Intervall der Länge $\left( \frac{1}{4} \right)^k$ entfernt usw.
Die so entstandene Menge ist eine kompakte (also $\mu_1$-integrierbare) Teilmenge von $[0,1]$ mit leerem offenen Kern, die Null als inneres Jordan-$\mu_1$-Maß und $1 - \lim_{k \to \infty} \sum_{i=1}^k \frac{2^{i-1}}{4^i} = \frac{1}{2}$ als äußeres Jordan-$\mu_1$-Maß besitzt.

Die nicht-negative charakteristische Funktion der überabzählbaren Smith-Volterra-Cantor Menge ist übrigens ein Element von $\L_{\R}(\R,\mu_1) \setminus \mathcal{S}_{\nearrow}(\R,\mu_1)$.
\end{Bsp}

\begin{Def} \label{FA.4.57} 
Seien $M$ eine beschränkte Teilmenge von $\R^n$ und $f \: M \to \R$ eine stetige Funktion. 
\begin{itemize}
\item[(i)] Seien $\delta \in \R_+$ und $T \in \mathcal{S}(\R^n)$.

$T$ heißt \emph{$\delta$-fein für $f$ auf $M$} genau dann, wenn $k \in \N_+$, paarweise disjunkte Quader $Q_1, \ldots, Q_k \in \mathfrak{Q}_n$ derart, daß $M \subset \bigcupdot_{\kappa = 1}^k Q_{\kappa}$ gilt, sowie $\alpha_1, \ldots, \alpha_k \in \R$ existieren mit
\begin{gather*}
\forall_{\kappa \in \{1, \ldots, k\}} \, {\rm diam}^{d_{\infty}} (Q_{\kappa}) = \sup \left\{ \|x-y\|_{\infty} \, | \, x,y \in Q_{\kappa} \right\} < \delta, \label{FA.4.57.1} \\
T  = \sum_{\kappa = 1}^k \alpha_{\kappa} \, \chi_{Q_{\kappa}}, \label{FA.4.57.2} \\
\forall_{\kappa \in \{1, \ldots, k\}} \exists_{x \in \overline{Q_{\kappa}} \cap M} \, \alpha_{\kappa} = f(x). \label{FA.4.57.3}
\end{gather*}
\item[(ii)] Eine Folge $(T_i)_{i \in \N}$ in $\mathcal{S}(\R^n)$ heißt \emph{ausgezeichnete Folge von Treppenfunktionen für $f$ auf $M$}\index{Folgen!ausgezeichnete -- von Treppenfunktionen} genau dann, wenn $(\delta_i)_{i \in \N} \in (\R_+)^{\N}$ mit $\lim_{i \to \infty} \delta_i = 0$ und 
$$ \forall_{i \in \N} \, T_i \mbox{ ist $\delta_i$-fein für $f$ auf $M$} $$
existiert.
\end{itemize}
\end{Def}

\begin{Satz} \label{FA.4.58} $\,$

\noindent \textbf{Vor.:} Es seien $M$ eine beschränkte Teilmenge von $\R^n$ und $f \: M \to \R$ eine stetige Funktion.

\noindent \textbf{Beh.:} Es folgt:
\begin{itemize}
\item[(i)] Es existiert eine ausgezeichnete Folge von Treppenfunktionen für $f$ auf $M$.
\item[(ii)] Ist $(T_i)_{i \in \N}$ eine beliebige ausgezeichnete Folge von Treppenfunktionen für $f$ auf $M$, so gilt
$$ \forall_{x \in M} \, \lim_{i \to \infty} T_i(x) = f(x). $$

\begin{Zusatz}
Im Falle der Kompaktheit von $M$ konvergiert jede ausgezeichnete Folge von Treppenfunktionen für $f$ auf $M$ gleichmäßig auf $M$ gegen $f$.
\end{Zusatz}
\end{itemize}
\end{Satz}

\textit{Beweis.} Zu (i): Es sei $i \in \N$.
Wir definieren Quader durch
$$ \forall_{a = (a_1, \ldots, a_n) \in \Z^n} \, Q_{a,i} := \left[ \frac{a_1}{2^i}, \frac{a_1 + 1}{2^i} \right[ \times \ldots \times \left[ \frac{a_n}{2^i}, \frac{a_n + 1}{2^i} \right[, $$
d.h.\ $\bigcupdot_{a \in \Z^n} Q_{a,i} = \R^n$ sowie
$$ \forall_{a \in \Z^n} \, {\rm diam}^{d_{\infty}} (Q_{a,i}) = \sup \{ \| x-y \|_{\infty} \, | \, x,y \in Q_{a,i} \} = \frac{\sqrt{n}}{2^i} \stackrel{\i \to \infty}{\longrightarrow} 0, $$
und setzen
$$ A_i := \{ a \in \Z^n \, | \, \overline{Q_{a,i}} \cap M \ne \emptyset \}. $$
Wegen der Beschränktheit von $M$ muß $A_i$ endlich sein.
Wir wählen nun für jedes $a \in A_i$ ein $x_{a,i} \in \overline{Q_{a,i}} \cap M$ und definieren
$$ T_i := \sum_{a \in A_i} f(x_{a,i}) \, \chi_{Q_{a,i}} \in \mathcal{S}(\R^n), $$
beachte $Q_{a,i} \cap Q_{\tilde{a},i} = \emptyset$ für $a, \tilde{a} \in \Z^n$ mit $a \ne \tilde{a}$.

Dann erfüllt $(T_i)_{i \in \N}$ offenbar die Behauptung von (i).
\pagebreak

Zu (ii): Seien $(T_i)_{i \in \N}$ eine ausgezeichnete Folge von Treppenfunktionen für $f$ auf $M$, $x_0 \in M$ und $\varepsilon \in \R_+$.
Da $f$ stetig ist, existiert zunächst $\delta \in \R_+$ mit
$$ \forall_{x \in M} \, \left( \| x - x_0 \|_{\infty} < \delta \Longrightarrow |f(x) - f(x_0)| < \varepsilon \right).$$
Sodann können wir $i_0 \in \N$ wählen derart, daß für alle $i \in \N$ mit $i \ge i_0$ gilt
$$ \mbox{$T_i$ ist $\delta$-fein für $f$ auf $M$.} $$
Für jedes solche $i$ gilt dann offenbar
\begin{equation*}
|T_i(x_0) - f(x_0)| = |f(x) - f(x_0)| < \varepsilon,
\end{equation*}
wobei $x \in \overline{Q} \cap M$, $Q \in \mathfrak{Q}_n$ mit ${\rm diam}^{d_{\infty}} (Q) < \delta$, und $T_i(x_0) = f(x)$.

Der Beweis des Zusatzes verläuft analog.
Man nutzt aus, daß wegen der Kompaktheit von $M$ dann $f \: M \to \R$ sogar gleichmäßig stetig ist und $\delta$ sowie $i_0$ daher nicht mehr von $x_0$ abhängen.
\q
 
\begin{HS}[Integrierbarkeit stetiger beschränkter Funktionen über integrierbaren Mengen] \label{FA.4.59} $\,$

\noindent \textbf{Vor.:} Seien $\varphi$ ein Quadermaß auf $\R^n$, $M \subset \R^n$ eine $\varphi$-integrierbare Menge und $f \: M \to \K$ eine stetige beschränkte Funktion.
(Die $\varphi$-Integrierbarkeit von $M$ ist nach \ref{FA.4.54} z.B.\ im Falle eines Kompaktums $M$ erfüllt.
Als stetige Funktion auf einem Kompaktum ist dann die Voraussetzung der Beschränktheit von $f$ evident.)

\noindent \textbf{Beh.:} $f \in \L_{\K}(M,\varphi)$, und im Falle der Beschränktheit von $M$ sowie $\K = \R$ gilt für jede ausgezeichnete Folge von Treppenfunktionen für $f$ auf $M$
$$ \int_M f \, \d \varphi = \lim_{i \to \infty} \int_M  T_i \, \d \varphi. $$
\end{HS}

Wir bereiten den Beweis des Hauptsatzes durch folgendes Lemma vor.

\begin{Lemma} \label{FA.4.60} 
Seien $\varphi$ ein Quadermaß auf $\R^n$, $M \subset \R^n$ eine $\varphi$-integrierbare Menge und $T \in \mathcal{S}(\R^n)$.

Dann gilt $T \in \L_{\R}(M,\varphi)$.
\end{Lemma}

\textit{Beweis.} Für jedes $Q \in \mathfrak{Q}_n$ gilt 
$$ \chi_Q \in \mathcal{S}(\R^n) \subset \L_{\R}(\R^n,\varphi), $$
und es gilt weiterhin nach Voraussetzung
$$ \chi_M \in \L_{\R}(\R^n,\varphi), $$
also folgt
$$ \chi_M \cdot \widehat{\chi_Q} = \chi_M \cdot \chi_Q = \inf ( \chi_M, \chi_Q ) \stackrel{\ref{FA.4.22} (v)}{\in} \L_{\R}(\R^n,\varphi), $$
d.h.\ genau 
$$ \chi_Q \in \L_{\R}(M,\varphi). $$

Hieraus, \ref{FA.4.8} (i) sowie \ref{FA.4.24} (i), (ii) folgt die Behauptung. \q
\A
\textit{Beweis des Hauptsatzes.} Ohne Beschränkung der Allgemeinheit genügt es, den Fall $\K = \R$ und $M \ne \emptyset$ zu betrachten.

1.\ Fall: $M$ beschränkt.
Sei $(T_i)_{i \in \N}$ eine -- nach \ref{FA.4.58} (i) existierende -- ausgezeichnete Folge von Treppenfunktionen für $f$ auf $M$, also gilt nach \ref{FA.4.60}
\begin{equation} \label{FA.4.59.1}
\forall_{i \in \N} \, T_i \in \L_{\R}(M,\varphi)
\end{equation}
und nach Definition \ref{FA.4.57}
\begin{equation} \label{FA.4.59.2}
\forall_{i \in \N} \, T_i(\R^n) \subset f(M).
\end{equation}

$f$ ist beschränkt, also folgt
\begin{equation} \label{FA.4.59.3}
C := \inf f(M) \in \R, ~~ D := \sup f(M) \in \R.
\end{equation}
Des weiteren ist $M$ $\varphi$-integrierbar, d.h.\
$$ \chi_M \cdot \widehat{1_{\R^n}} = \chi_M \in \L_{\R}(\R^n,\varphi), $$
wobei $1_{\R^n} \: \R^n \to \R$ die konstante Funktion vom Wert $1$ bezeichne, also gilt
\begin{equation} \label{FA.4.59.4}
C \cdot 1_{\R^n}, \, D \cdot 1_{\R^n} \in \L_{\R}(M,\varphi).
\end{equation}

Wegen (\ref{FA.4.59.2}), (\ref{FA.4.59.3}) ergibt sich zusätzlich
\begin{equation} \label{FA.4.59.5}
\forall_{i \in \N} \forall_{x \in \R^n} \, C \cdot 1_{\R^n}(x) \le T_i(x) \le D \cdot 1_{\R^n}(x).
\end{equation}

Aus (\ref{FA.4.59.1}), \ref{FA.4.58} (ii), (\ref{FA.4.59.4}), (\ref{FA.4.59.5}) und dem Grenzwertsatz von \textsc{Lebesgue} folgt die Behauptung im 1.\ Fall.

2.\ Fall: $M$ unbeschränkt.
Wir haben zu zeigen, daß $f \in \L_{\R}(M,\varphi)$ gilt, und wir können wegen $f = f^+ - f^-$ ohne Einschränkung $f \ge 0$ annehmen.

Es existiert eine Folge $(Q_i)_{i \in \N}$ in $\mathfrak{Q}_n$ mit
$$ \forall_{i \in \N} \, Q_i \subset Q_{i+1} ~~ \mbox{ und } ~~ \bigcup_{i=0}^{\infty} Q_i = \R^n. $$
Definiere für jedes $i \in \N$ die beschränkte $\varphi$-integrierbare Menge
$$ M_i := M \cap Q_i. $$
Dann gilt
$$ \forall_{i \in \N} \, M_i \subset M_{i+1} \subset M ~~ \mbox{ sowie } ~~ \bigcup_{i=0}^{\infty} M_i = M $$
und nach dem 1.\ Fall
$$ \forall_{i \in \N} \, f \in \L_{\R}(M_i,\varphi). $$
Wegen $f \ge 0$, $\forall_{i \in \N} \, M_i \subset M_{i+1} \subset M$, der Beschränktheit von $f$ und der $\varphi$-In\-te\-grier\-bar\-keit von $M$ ist $\left( \int_{M_i} f \, \d \varphi \right)_{i \in \N}$ eine monoton wachsende reelle Folge, die durch $\sup f(M) \cdot \varphi(M) \in \R$ beschränkt ist und somit in $\R$ konvergiert.
Daher folgt aus \ref{FA.4.61} (ii): $f \in \L_{\R}(M,\varphi)$. \q

\begin{Kor} \label{FA.4.62}
Seien $\varphi$ ein Quadermaß auf $\R^n$ und $K$ eine kompakte Teilmenge von $\R^n$.

Dann gilt $\mathcal{C}(K,\K) \subset \L_{\K}(K,\varphi)$. \q
\end{Kor}

\subsection*{Stetige und differenzierbare Abhängigkeit} \addcontentsline{toc}{subsection}{Stetige und differenzierbare Abhängigkeit}

\begin{Satz}[Stetige Abhängigkeit des Integrales von einem stetigen Parameter im Integranden] \label{FA.4.63} $\,$

\noindent \textbf{Vor.:} Seien $\varphi$ ein Quadermaß auf $\R^n$ und $M$ eine Teilmenge von $\R^n$.
Ferner seien $A$ ein metrischer Raum, $a_* \in A$ und $F \: M \times A \to \K$ eine Funktion mit folgenden Eigenschaften:
\begin{itemize}
\item[(a)] $\forall_{a \in A} \, F_a := F(\ldots,a) \in \L_{\K}(M,\varphi)$.
\item[(b)] Es existiert eine $\varphi$-Nullmenge $N$ derart, daß $F(x,\ldots) \: A \to \K$ für jedes $x \in M \setminus N$ stetig in $a_*$ ist.
\item[(c)] Es existieren $U \in \U(a_*,A)$, $g \in \L_{\R}(M,\varphi)$ und eine $\varphi$-Nullmenge $\widetilde{N}$ mit $\forall_{(x,a) \in (M \setminus \widetilde{N}) \times U} \, |F(x,a)| \le g(x)$.
\end{itemize}

\noindent \textbf{Beh.:} Die Funktion
$$ A \longrightarrow \K, ~~ a \longmapsto \int_M F_a \, \d \varphi $$ 
ist stetig in $a_*$.
\end{Satz}

\begin{Bem*}
Versieht man $M \times A$ mit der Produkttopologie\footnote{\textbf{Definition.} Für topologische Räume $X_1, X_2$ ist die \emph{Produkttopologie für $X_1 \times X_2$}\index{Topologie!Produkt-} durch
$$ \forall_{U \in \mathfrak{P}(X_1 \times X_2)} \, \left( U \mbox{ offen in } X_1 \times X_2 : \Longleftrightarrow \forall_{(x_1,x_2) \in U } \exists_{U_1 \in \U(x_1,X_1)} \exists_{U_2 \in \U(x_2,X_2)} \, U_1 \times U_2 \subset U \right) $$
gegeben.}, so folgen (a), (b) und (c), wenn die Kompaktheit von $M$ und die Stetigkeit sowie Beschränktheit von $F \: M \times A \to \K$ vorausgesetzt wird.
\end{Bem*}

\textit{Beweis des Satzes.} Sei $(a_i)_{i \in \N}$ eine Folge in $A$ mit $\lim_{i \to \infty} a_i = a$.
Ohne Beschränkung der Allgemeinheit gelte $\forall_{i \in \N} \, a_i \in U$.
Wir setzen
$$ \forall_{i \in \N} \, f_i := F_{a_i} \stackrel{(a)}{\in} \L_{\K}(M,\varphi). $$
Dann gilt nach (b)
$$ \forall_{x \in M \setminus N} \, \lim_{i \to \infty} f_i(x) = \lim_{i \to \infty} F(x,a_i) = F(x,a_*) = F_{a_*}(x) $$
und nach (c)
$$ \forall_{i \in \N} \forall_{x \in M \setminus \widetilde{N}} \, |f_i(x)| = |F(x,a_i)| \le g(x). $$
Also folgt aus $g \in \L_{\R}(M,\varphi)$ und dem Grenzwertsatz von \textsc{Lebesgue}
$$ \lim_{i \to \infty} \int_M \underbrace{f_i}_{=F_{a_i}} \, \d \varphi = \int_M F_{a_*} \, \d \varphi, $$
womit die Behauptung gezeigt ist. \q

\begin{Satz}[Differenzierbare Abhängigkeit des Integrales von einem differenzierbaren Parameter im Integranden] \label{FA.4.64} $\,$

\noindent \textbf{Vor.:} Seien $\varphi$ ein Quadermaß auf $\R^n$ und $M$ eine Teilmenge von $\R^n$.
Ferner seien $U$ eine nicht-leere offene Teilmenge eines endlich-dimensionalen normierten $\R$-Vektorraumes $(Y,\|\ldots\|)$ und $F \: M \times U \to \R$ eine Funktion mit folgenden Eigenschaften:
\begin{itemize}
\item[(a)] $\forall_{y \in U} \, F_y := F(\ldots,y) \in \L_{\R}(M,\varphi).$
\item[(b)] Es existiert eine $\varphi$-Nullmenge $N$ derart, daß $F(x,\ldots) \: U \to \R$ für jedes $x \in M \setminus N$ differenzierbar ist.
Das Differential von $F(x,\ldots)$ in $y_* \in U$ sei dann für jedes solche $x$ mit $\D \frac{\partial F}{\partial y}(x,y_*) \in \L_{\R}(Y,\R)$ bezeichnet, d.h.\
$$ \forall_{(x,y_*) \in (M \setminus N) \times U} \forall_{v \in Y} \, \frac{\partial F}{\partial y}(x,y_*) (v) = \lim_{t \to 0} \frac{F(x,y_* + t \, v) - F(x,y_*)}{t}. $$
\item[(c)] Es existieren eine $\varphi$-Nullmenge $\widetilde{N}$ und $g \in \L_{\R}(M,\varphi)$ mit
$$ \forall_{(x,y_*) \in \left( M \setminus (N \cup \widetilde{N}) \right) \times U} \, \left\| \D \frac{\partial F}{\partial y}(x,y_*) \right\| \le g(x). $$
\end{itemize}

\noindent \textbf{Beh.:} $\D \forall_{y_* \in U} \forall_{v \in Y} \, \frac{\partial F}{\partial y}(\ldots,y_*)(v) \in \L_{\R}(M,\varphi)$ und die Funktion
$$ \psi \: U \longrightarrow \R, ~~ y \longmapsto \int_M F_y \, \d \varphi $$ 
ist differenzierbar mit
$$ \forall_{y_* \in U} \forall_{v \in Y} \, \d_y \psi (v) = \int_M \frac{\partial F}{\partial y}(\ldots,y_*) (v) \, \d \varphi.$$
\end{Satz}

\begin{Bem*}
Sind $M$ kompakt, $F$ stetig und $\D \frac{\partial F}{\partial y}$ auf ganz $M \times U$ definiert sowie stetig, so folgen (a), (b) und (c).
\end{Bem*}

\textit{Beweisskizze des Satzes.} 1.) Seien $y_* \in U$, $v \in Y$ und $(t_i)_{i \in \N}$ eine Nullfolge in $\R \setminus \{0\}$.
Wegen der Offenheit von $U$ können wir für jedes $i \in \N$ ohne Einschränkung annehmen
$$ y_* + t_i \, v \in U $$
und $f_i \: M \to \R$ durch
$$ \forall_{x \in M} \, f_i(x) := \frac{F(x,y_* + t_i \, v) - F(x,y_*)}{t_i} $$
definieren, also gilt nach (a): $f_i \in \L_{\R}(M,\varphi)$.
Wegen (b) und dem Mittelwertsatz der Differentialrechnung, siehe z.B.\ \cite[Satz 10.9]{ElAna}, existiert zu $x \in M \setminus (N \cup \widetilde{N})$ eine Zahl $\xi_i \in {]}t_i,0{[} \cup {]}0,t_i{[}$ mit
$$ f_i(x) = \frac{\partial F}{\partial y} (x,y_* + \xi_i \, v)(v), $$
d.h.\ nach (c)
$$ |f_i(x)| \le \left\| \frac{\partial F}{\partial y} (x,y_* + \xi_i \, v) \right \| \, \|v\| \le g(x) \, \|v\| $$
und $\|v\| \, g \in \L_{\R}(M,\varphi)$.
Aus dem Grenzwertsatz von \textsc{Lebesgue} folgt daher
$$ \frac{\partial F}{\partial y}(\ldots,y_*)(v) = \lim_{i \to \infty} f_i \in \L_{\R}(M,\varphi). $$

2.) Sei $y_* \in U$.
Die Funktion
$$ l \: Y \longrightarrow \R, ~~ v \longmapsto \int_M \frac{\partial F}{\partial y}(\ldots,y_*)(v) \, \d \varphi $$
ist nach 1.) offenbar $\R$-linear.
Sei $(y_i)_{i  \in \N}$ eine Folge in $U \setminus \{y_*\}$ mit $\lim_{i \to \infty} y_i = y_*$.
Wir können nun wegen der Offenheit von $U$ die Existenz von $\delta \in \R_*$ derart, daß für alle $i \in \N$ gilt
$$ y_i \in U_{\delta}(y_*,Y) \subset U $$
annehmen, und $g_i \in \L_{\R}(M, \varphi)$ nach (a) und 1.) durch
$$ \forall_{x \in M \setminus N} \, g_i(x) := \frac{F(x,y_i) - F(x,y_*) - \frac{\partial F}{\partial y}(x,y_*)(y_i-y_*)}{\|y_i - y_*\|} $$
definieren.
Dann gilt wegen (b)
\begin{equation} \label{FA.4.64.1}
\forall_{x \in M \setminus N} \, \lim_{i \to \infty} g_i(x) = 0. 
\end{equation}
Aus dem Mittelwertsatz der Differentialrechnung folgt zu jedem $i \in \N$ die Existenz eines $\xi_i \in {]}0,1{[}$ mit
\begin{equation*}
\forall_{x \in M \setminus N} \, F(x,y_i) - F(x,y_*) = \frac{\partial F}{\partial y} \left( x,y_* + \xi_i \, (y_i - y_*) \right) (y_i-y_*),
\end{equation*}
d.h.\ nach (c)
\begin{eqnarray*}
\forall_{x \in M \setminus (N \cup \widetilde{N})} \, |g_i(x)| & = & \left| \frac{\frac{\partial F}{\partial y} \left( x,y_* + \xi_i \, (y_i - y_*) \right) (y_i-y_*) - \frac{\partial F}{\partial y}(x,y_*)(y_i-y_*)}{\|y_i - y_*\|} \right| \\
& \le & \left\| \frac{\partial F}{\partial y} \left( x,y_* + \xi_i \, (y_i - y_*) \right) - \frac{\partial F}{\partial y}(x,y_*) \right\| \le 2 g(x),
\end{eqnarray*}
und mit $g$ ist auch $2 g$ $\varphi$-integrierbar über $M$.
Der Grenzwertsatz von \textsc{Lebesgue} und (\ref{FA.4.64.1}) ergeben nun
$$ \lim_{i \to \infty} \int_M g_i \, d \varphi = 0, $$
also offenbar
$$ \lim_{i \to \infty} \frac{\psi(y_i) - \psi(y_*) - l(y_i-y_*)}{\|y_i - y_*\|} = 0, $$
und der Satz ist nach der Definition von $l$ vollständig bewiesen. \q

\subsection*{Der Satz von \textsc{Fubini} und der Transformationssatz} \addcontentsline{toc}{subsection}{Der Satz von \textsc{Fubini} und der Transformationssatz}

Die Sätze dieses Abschnittes dienen der konkreten Berechnung von Integralen.
Der folgende Hauptsatz ermöglicht u.a.\ die Rückführung der Bestimmung eines $\mu_n$-Integrales auf die Berechnung von $\mu_1$-Integralen.

\begin{HS}[Satz von \textsc{Fubini}]\index{Satz!von \textsc{Fubini}} \label{FA.4.65} $\,$

\noindent \textbf{Vor.:} Seien $n_1, n_2 \in \N_+$ mit $n = n_1 + n_2$, $\varphi_i$ Quadermaße auf $\R^{n_i}$ für $i \in \{1,2\}$ und $f \in \L_{\K}(\R^n,\varphi_1 \times \varphi_2)$.

\noindent \textbf{Beh.:} Es existiert eine $\varphi_1$-Nullmenge $N_1$ derart, daß gilt
$$ \forall_{x_1 \in \R^{n_1} \setminus N_1} \, f_{x_1} := f(x_1,\ldots) \in \L_{\K}(\R^{n_2},\varphi_2), $$
und die Funktion
$$ \R^{n_1} \setminus N_1 \longrightarrow \K, ~~ x_1 \longmapsto \int\limits_{\R^{n_2}} f_{x_1} \, \d \varphi_2, $$
ist $\varphi_1$-integrierbar über $\R^{n_1}$ mit
$$ \int\limits_{\R^n} f \, \d (\varphi_1 \times \varphi_2) = \int\limits_{\R^{n_1}} \left( \; \int\limits_{\R^{n_2}} f_{\ldots} \, \d \varphi_2 \right) \, \d \varphi_1 =: \int\limits_{\R^{n_1}} \left( \; \int\limits_{\R^{n_2}} f(x_1,x_2) \, \d \varphi_2(x_2) \right) \, \d \varphi_1(x_1). $$
\end{HS}

\begin{Bem*}
Völlig analog gilt auch die Version des Satzes von \textsc{Fubini}, die man durch Vertauschung der Reihenfolge der Integrationen erhält.
\end{Bem*} 

\textit{Beweisskizze des Hauptsatzes.}
Ein ausführlicher Beweis im Falle $\varphi_i = \mu_{n_i}$, wobei $i \in \{1,2\}$, findet sich in \cite[Hauptsatz 11.36]{Henke}.
Der Beweis im allgemeinen Falle ist analog zu führen:

1.) Der Satz von \textsc{Fubini} gilt mit $N_1 = \emptyset$ für alle Funktionen $f = \chi_Q$, wobei $Q \in \mathfrak{Q}_n$ sei.

2.) Sind $f, \tilde{f} \in \L_{\R}(\R^n,\varphi)$ und $\lambda, \tilde{\lambda} \in \R$ und gilt der Satz von \textsc{Fubini} für $f$ mit einer $\varphi_1$-Nullmenge $N_1$ und für $\tilde{f}$ mit einer $\varphi_1$-Nullmenge $\widetilde{N_1}$, so gilt er auch für $\lambda \, f + \tilde{\lambda} \, \tilde{f}$ mit der $\varphi_1$-Nullmenge $N_1 \cup \widetilde{N_1}$.

3.) Für jedes Quadermaß $\varphi$ auf $\R^n$ gilt:
Eine Teilmenge $N$ von $\R^n$ ist genau dann eine $\varphi$-Nullmenge, wenn eine Folge $(Q_k)_{k \in \N}$ in $\mathfrak{Q}_n$ mit
\begin{gather}
\sum_{k=0}^{\infty} \varphi(Q_k) \mbox{ ist konvergent in $\R$ und} \label{FA.4.65.1} \\
\forall_{x \in N} \, \# \{ k \in \N \, | \, x \in Q_k \} = \infty \label{FA.4.65.2}
\end{gather}
existiert.

{[} ,,$\Rightarrow$`` Sei $N$ eine $\varphi$-Nullmenge.
Dann existiert zu jedem $i \in \N$ eine Folge $(Q_{ij})_{j \in \N}$ mit 
$$ \forall_{i \in \N} \, \left( N \subset \bigcup_{j=0}^{\infty} Q_{ij} \mbox{ und } \sum_{j=0}^{\infty} \varphi(Q_{ij}) < \frac{1}{2^i} \right). $$
Numeriere die Folge $(Q_{ij})_{(i,j) \in \N \times \N}$ zu einer Folge $(Q_k)_{k \in \N}$ um.
Dann gilt
$$ \sum_{k=0}^{\infty} \varphi(Q_k) = \sum_{i=0}^{\infty} \sum_{j=0}^{\infty} \varphi(Q_{ij}) \le 2 $$
und für jedes $x \in N$
$$ \forall_{i \in \N} \exists_{j \in \N} \, x \in Q_{ij}. $$

,,$\Leftarrow$`` Sei umgekehrt $(Q_k)_{k \in \N}$ eine Folge in $\mathfrak{Q}_n$ mit (\ref{FA.4.65.1}) und (\ref{FA.4.65.2}).
Zu $\varepsilon \in \R_+$ existiert dann $k_0 \in \N$ mit $\sum_{k=k_0}^{\infty} \varphi(Q_k) < \varepsilon$, und es gilt $x \in \bigcup_{k=k_0}^{\infty} Q_k$ für jedes $x \in N$. {]}

4.) Seien $N$ eine $\varphi$-Nullmenge und
$$ \forall_{a_1 \in \R^{n_1}} \, N_{a_1} := \{ x_2 \in \R^{n_2} \, | \, (a_1,x_2) \in N \} \subset \R^{n_2}. $$
Dann existiert eine $\varphi_1$-Nullmenge $N_1$ mit
$$ \forall_{a_1 \in \R^{n_1} \setminus N_1} \, N_{a_1} \mbox{ ist $\varphi_2$-Nullmenge.} $$

{[} Wähle gemäß 3.) ,,$\Rightarrow$`` eine Folge $(Q_k)_{k \in \N}$ in $\mathfrak{Q}_n$ mit (\ref{FA.4.65.1}), wobei $\varphi$ durch $\varphi_1 \times \varphi_2$ ersetzt sei, und (\ref{FA.4.65.2}).
Für jedes $k \in \N$ seien $Q_k' \in \mathfrak{Q}_{n_1}$ und $Q_k'' \in \mathfrak{Q}_{n_2}$ mit $Q = Q_k' \times Q_k''$ und für jedes $i \in \N$
$$ h_i := \sum_{k=0}^i \underbrace{\varphi_2(Q_k'') \cdot \chi_{Q_k'}}_{\ge 0} \in \mathcal{S}(\R^{n_1}) \subset \L_{\R}(\R^{n_1},\varphi_1). $$
Dann gilt $h_i \le h_{i+1}$ und
$$ \int\limits_{\R^{n_1}} h_i \d \varphi_1 = \sum_{i=0}^i (\varphi_1 \times \varphi_2)(Q) \le \sum_{i=0}^{\infty} (\varphi_1 \times \varphi_2)(Q) < \infty. $$
Somit folgt aus dem Grenzwertsatz von \textsc{Levi} die Existenz einer $\varphi_1$-Nullmenge $N_1$ mit
$$ \forall_{a_1 \in \R^{n_1} \setminus N_1} \, \lim_{i \to \infty} h_i(a_1) \in \R. $$
Mittels 3.) ,,$\Leftarrow$`` zeigt man
$$ \forall_{a_1 \in \R^{n_1} \setminus N_1} \, N_{a_1} \mbox{ ist $\varphi_2$-Nullmenge,} $$
womit 4.) bewiesen ist. {]}

5.) Der Hauptsatz gilt, wenn man $\L_{\K}(\R^n,\varphi_1 \times \varphi_2)$ durch $\mathcal{S}_{\nearrow}(\R^n,\varphi_1 \times \varphi_2)$ ersetzt.

{[} Sei $f \in \mathcal{S}_{\nearrow}(\R^n,\varphi_1 \times \varphi_2)$.
Dann existieren eine Folge $(T_i)_{i \in N}$ in $\mathcal{S}(\R^n)$ und eine $(\varphi_1 \times \varphi_2)$-Nullmenge $N$ mit
\begin{gather}
\forall_{i \in \N} \forall_{x \in \R^n} \, T_i(x) \le T_{i+1}, \label{FA.4.65.9} \\
\forall_{x \in \R^n \setminus N} \, \lim_{i \to \infty} T_i(x) = f(x), \label{FA.4.65.10} \\
\int\limits_{\R^n} f \, \d (\varphi_1 \times \varphi_2) = \lim_{i \to \infty} \int\limits_{\R^n} T_i \, \d (\varphi_1 \times \varphi_2). \label{FA.4.65.11}
\end{gather}
Wegen 1.), 2.) gilt der Hauptsatz (mit $N_1 = \emptyset$) für alle Treppenfunktionen auf $\R^n$.
Daher folgt für jedes $i \in \N$ sowie $a_1 \in \R^{n_1}$ wegen (\ref{FA.4.65.9})
$$ \underbrace{\int\limits_{\R^{n_2}} T_i(a_1,\ldots) \, \d \varphi_2}_{=: \int\limits_{\R^{n_2}} T_i(a_1,x_2) \, \d \varphi_2(x_2)} \le \int\limits_{\R^{n_2}} T_{i+1}(a_1,\ldots) \, \d \varphi_2, $$
und $\left( \; \int\limits_{\R^{n_2}} T_i(\ldots, x_2) \, \d \varphi_2(x_2) \right)_{i \in \N}$ ist eine monoton wachsende Folge von Funktionen aus $\L_{\R}(\R^{n_1},\varphi_1)$ derart, daß nach (\ref{FA.4.65.9}), (\ref{FA.4.65.11}) für jedes $i \in \N$ gilt
$$ \int\limits_{\R^{n_1}} \left( \; \int\limits_{\R^{n_2}} T_i(x_1,x_2) \, \d \varphi_2(x_2) \right) \d \varphi_1(x_1) = \int\limits_{\R^n} T_i \, \d (\varphi_1 \times \varphi_2) \le \int\limits_{\R^n} f \, \d (\varphi_1 \times \varphi_2), $$
also folgt aus dem Grenzwertsatz von \textsc{Levi} die Existenz einer $\varphi_1$-Nullmenge $\widetilde{N_1}$ mit
\begin{gather*}
\forall_{a_1 \in \R^{n_1} \setminus \widetilde{N_1}} \, \lim_{i \to \infty} \int\limits_{\R^{n_2}} T_i(a_1,x_2) \, \d \varphi_2(x_2) \in \R, \\
\lim_{i \to \infty} \, \int\limits_{\R^{n_2}} T_i(\ldots,x_2) \, \d \varphi_2(x_2) \in \L_{\R}(\R^{n_1},\varphi_1)
\end{gather*}
\begin{equation} \label{FA.4.65.S}
\mbox{ und }
\end{equation}
\begin{eqnarray*}
\lefteqn{\int\limits_{\R^{n_1}} \left( \lim_{i \to \infty} \, \int\limits_{\R^{n_2}} T_i(\ldots,x_2) \, \d \varphi_2(x_2) \right) = \lim_{i \to \infty} \, \int\limits_{\R^{n_1}} \left( \; \int\limits_{\R^{n_2}} T_i(\ldots,x_2) \, \d \varphi_2(x_2) \right) \, \d \varphi_1} ~~~~~~~~~~~~~~~~~~~~~~~~~~~~~~~~~~~ && \\
&&  = \lim_{i \to \infty} \, \int\limits_{\R^n} T_i \, \d (\varphi_1 \times \varphi_2) \stackrel{(\ref{FA.4.65.11})}{=} \int\limits_{\R^n} f \, \d (\varphi_1 \times \varphi_2).
\end{eqnarray*}

Wir zeigen nun (die Existenz der Integrale in der folgenden Gleichung und) die Gleichung
\begin{equation} \label{FA.4.65.DS}
\int\limits_{\R^{n_2}} f(\ldots,x_2) \, \d \varphi_2(x_2) =_{\varphi_1} \lim_{i \to \infty} \, \int\limits_{\R^{n_2}} T_i(\ldots,x_2) \, \d \varphi_2(x_2).
\end{equation}
Wegen (\ref{FA.4.65.S}) ist 5.) dann gezeigt.

Beweis von (\ref{FA.4.65.DS}): Nach 4.) existiert ein $\varphi_1$-Nullmenge $\widetilde{\widetilde{N_1}}$ derart, daß $N_{a_1}$ für jedes $a_1 \in \R^{n_1} \setminus \widetilde{\widetilde{N_1}}$ eine $\varphi_1$-Nullmenge ist.
Dann ist auch $N_1 := \widetilde{N_1} \cup \widetilde{\widetilde{N_1}}$ eine $\varphi_1$-Nullmenge und für alle $a_1 \in \R^{n_1} \setminus N_1$ gilt -- beachte (\ref{FA.4.65.9}) --, die Gültigkeit des Hauptsatzes für alle Treppenfunktionen (mit $N_1 = \emptyset$) und (\ref{FA.4.65.S}):

$\left( T_i(a_1,\ldots) \right)_{i\in \N}$ ist eine monoton wachsende Folge in $\L_{\R}(\R^{n_2},\varphi_2)$ mit beschränkter (sogar in $\R$ konvergenter) Integralfolge, $N_{a_1}$ ist eine $\varphi_2$-Nullmenge und
$$ \forall_{x_2 \in \R^{n_2} \setminus N_{a_1}} \, (a_1,x_2) \in \R^n \setminus N, $$
d.h.\ nach (\ref{FA.4.65.10})
$$ \lim_{i \to \infty} T_i(a_1,\ldots) =_{\varphi_2} f(a_1,\ldots). $$
Aus dem Satz von \textsc{Levi} folgt hieraus gerade (\ref{FA.4.65.DS}). {]}

6.) Es genügt offenbar, den Fall $\K = \R$ zu betrachten.
Der Satz von \textsc{Fubini} folgt dann aus 5.) und 2.). \q

\begin{Kor}[Prinzip von \textsc{Cavalieri}] \index{Prinzip!von \textsc{Cavalieri}} \label{FA.4.66} $\,$

\noindent \textbf{Vor.:} Es seien $\varphi_i$ ein Quadermaß auf $\R^i$ für $i \in \{1,n\}$, $M$ eine $(\varphi_1 \times \varphi_n)$-integrierbare Menge und für jedes $t \in \R$
$$ M_t := \{ x \in \R^n \, | \, (x,t) \in M \}. $$

\noindent \textbf{Beh.:} Es existiert eine $\varphi_1$-Nullmenge $N$ derart, daß $M_t$ für jedes $t \in \R \setminus N$ eine $\varphi_n$-integrierbare Menge und
$$ \R \setminus N \longrightarrow \R, ~~ t \longmapsto \varphi_n(M_t), $$
eine $\varphi_1$-integrierbare Funktion ist.
Des weiteren gilt
$$ (\varphi_1 \times \varphi_n)(M) = \int_{\R} \varphi_n(M_t) \, \d \varphi_1(t). $$
\q
\end{Kor}

\begin{Bsp*}[Volumina von Rotationskörpern]
Seien $a,b \in \R$, wobei $a < b$ gelte, $\rho \: [a,b] \to \R$ eine stetige Funktion mit $\rho \ge 0$ und
$$ M := \{ (x,y,t) \in \R^3 \, | \, t \in [a,b] \, \wedge \, \sqrt{x^2 + y^2} \le \rho(t) \}. $$

Dann ist $M$ eine $\mu_3$-integrierbare Menge und
$$ \mu_3(M) = \pi \int_a^b \rho(t)^2 \, \d t. $$
\end{Bsp*}

Der folgende Hauptsatz verallgemeinert i.w.\ die Substitutionsregel.
Der Autor möchte betonen, daß es bei der Beweisführung an mehreren Stellen entscheidend ist, das Quadermaß $\mu_n$ (und nicht etwa ein beliebiges Quadermaß auf $\R^n$) zu betrachten.
Z.B.\ wird die $\mu_n$-Nullmengentreue stetig differenzierbarer Abbildungen verwendet:

\begin{Satz} \label{FA.4.NmT}
Seien $U \subset \R^n$ offen und $h \: U \to \R^n$ eine stetig differenzierbare Abbildung.

Dann gilt: $\forall_{N \in \mathfrak{P}(U)} \, \left( \mbox{$N$ ist $\mu_n$-Nullmenge $\Longrightarrow$ $h(N)$ ist $\mu_n$-Nullmenge} \right).$
\end{Satz}

\textit{Beweisskizze.} Überlege zunächst, daß im Falle des Quadermaßes $\mu_n$ in Definition \ref{FA.4.13} äquivalent die Existenz abgeschlossener Quader anstelle beliebiger Quader gefordert werden kann.
Wende dann den Mittelwertabschätzungssatz an. \q

\begin{HS}[Transformationssatz] \index{Satz!Transformations-} \label{FA.4.67} $\,$

\noindent \textbf{Vor.:} Es seien $U$ eine offene Teilmenge von $\R^n$, $h \: U \to \R^n$ eine stetig differenzierbare Abbildung und $M$ eine Teilmenge von $U$ mit folgenden Eigenschaften:
\begin{itemize}
\item[(a)] $M \setminus (M)^{\circ}$ ist eine $\mu_n$-Nullmenge,
\item[(b)] $h|_{M^{\circ}}$ ist injektiv, und
\item[(c)] $\forall_{x \in M^{\circ}} ~ \d_x h \in \L_{\R}(\R^n,\R^n)$ ist bijektiv.
\end{itemize}
Ferner seien $D \subset \R^n$ und $f \in \R^D$.

\pagebreak
\noindent \textbf{Beh.:}
\begin{itemize}
\item[(i)] $f \in \L_{\R}(h(M),\mu_n) \Longleftrightarrow (f \circ h) \, | \det \d h | \in \L_{\R}(M,\mu_n)$.
\item[(ii)] Im Falle der Gültigkeit einer der äquivalenten Aussagen von (i) folgt
\begin{equation} \label{FA.4.67.S}
\int_{h(M)} f \, \d \mu_n = \int_M (f \circ h) \, | \det \d h | \, d \mu_n.
\end{equation}
\end{itemize}
\end{HS}

\textit{Beweisskizze.} 1.) Mittels des Satzes von \textsc{Levi} zeigt man, daß (i) ,,$\Rightarrow$`` und (\ref{FA.4.67.S}) für alle Translationen $t_a \: \R^n \to \R^n, \, x \mapsto x + a,$ anstelle von $h$ gilt.

2.) Seien $a_1 = \left( \begin{array}{c} a_{11} \\ \vdots \\ a_{1n} \end{array} \right), \ldots, a_n = \left( \begin{array}{c} a_{n1} \\ \vdots \\ a_{nn} \end{array} \right) \in \R^n$ und $P$ \emph{das von $a_1, \ldots, a_n$ aufgespannte Parallelotop}, d.h.\ genau $P = \sum_{i=1}^n [0,1] \, a_i$.

Dann folgt aus 1.): $\mu_n(P) = | \det (a_{ij})_{i,j \in \{1, \ldots, n\}} |$.

3.) Aus 2.) ergibt sich $\forall_{Q \in \mathfrak{Q}_n} \forall_{h \in \L_{\R}(\R^n,\R^n)} \, \mu_n(h(Q)) = | \det h | \cdot \mu_n(Q)$.

4.) Mittels \ref{FA.4.NmT},  1.) und 3.) sieht man ein, daß für alle offenen Teilmengen $U, \widetilde{U}$ von $\R^n$ und jeden $\mathcal{C}^1$-Diffeomorphismus\footnote{Ein \emph{$\mathcal{C}^1$-Diffeomorphismus zwischen offenen Teilmengen von $\R^n$} ist per definitionem eine bijektive stetig differenzierbare Abbildung zwischen offenen Teilmengen von $\R^n$ mit stetig differenzierbarer Umkehrabbildung.} $h \: U \to \widetilde{U}$ gilt:
$$ \forall_{Q \in \mathfrak{Q}_n} \, \left( \chi_{h(Q)} \in \mathcal{S}_{\nearrow}(\R^n,\varphi) \, \wedge \, \mu_n(h(Q)) \le \int_Q | \det \d h | \, \d \mu_n \right) . $$

5.) Seien $U$ eine offene Teilmenge von $\R^n$ und $f \in \mathcal{C}(U,\R)$ mit $f \ge 0$.

Dann existiert eine Folge $(T_i)_{i \in \N}$ in $\mathcal{S}(\R^n)$ derart, daß für alle $i \in \N$ gilt $0 \le T_i \le T_{i+1}$, ${\rm Tr}(T_i) \subset U$ und $\lim_{i \to \infty} T_i = f$ auf $U$.

6.) $\left( f \in \mathcal{S}_{\nearrow}(\R^n,\varphi) \, \wedge \, \widetilde{U} \in {\rm Top}(\R^n) \right) \Longrightarrow \chi_{\widetilde{U}} \, f \in \mathcal{S}_{\nearrow}(\R^n,\varphi)$.

7.) Mittels obiger Vorbereitungen zeigt man den Hauptsatz. \q

\begin{Bem*}
Obwohl in \cite[Kapitel 12]{ElAna} die $n$-dimensionale Riemannsche Integrationstheorie studiert wird, können Anwendungen des Transformationssatzes analog zu den dort in 12.55 bis 12.58 dargestellten eingesehen werden.
Außerdem ist der Hauptsatz durch getrennte Anwendung auf Real- und Imaginärteil (mit ggf.\ unterschiedlichem $h$) auch bei der Integration komplexwertiger Funktionen hilfreich.
\end{Bem*}

\subsection*{Lebesgue-Meßbarkeit und der Satz von \textsc{Tonelli}} \addcontentsline{toc}{subsection}{Lebesgue-Meßbarkeit und der Satz von \textsc{Tonelli}}

Der Begriff der Meßbarkeit wird in der Literatur unterschiedlich verwendet.
Wir verallgemeinern hier die in \cite{Henke} und \cite{Hirz} genannten Definitionen.

\begin{Def}[Lebesgue-meßbare Funktionen und Mengen] \index{Funktion!meßbare} \index{Menge!meßbare} \label{FA.4.68} $\,$
Es sei $\varphi$ ein Quadermaß auf $\R^n$.
\begin{itemize}
\item[(i)] 
\begin{itemize}
\item[(a)] Eine auf einer Teilmenge von $\R^n$ definierte reellwertige Funktion $f$ heißt genau dann \emph{(reellwertig) $\varphi$-meßbar über $\R^n$}, wenn eine Folge $(T_i)_{i \in \N}$ in $\mathcal{S}(\R^n)$ und eine $\varphi$-Nullmenge $N$ existieren derart, daß gilt 
$$ \forall_{x \in \R^n \setminus N} \, \lim_{i \to \infty} T_i(x) = f(x). $$
Die Menge solcher Funktionen bezeichnen wir mit $\boxed{\mathcal{M}_{\R}(\R^n,\varphi)}$.
\pagebreak
\item[(b)] Eine auf einer Teilmenge von $\R^n$ definierte komplexwertige Funktion $f$ heißt genau dann \emph{$\varphi$-meßbar über $\R^n$}, wenn sowohl ihr Real- als auch ihr Imaginärteil Elemente von $\mathcal{M}_{\R}(\R^n,\varphi)$ sind.
Die Menge jener Funktionen bezeichnen wir mit $\boxed{\mathcal{M}_{\C}(\R^n,\varphi)}$.
\item[(c)] Ist $M$ eine Teilmenge von $\R^n$, so heißt eine auf einer Teilmenge von $\R^n$ definierte $\K$-wertige Funktion $f$ genau dann \emph{(reellwertig, falls\linebreak $\K=\R$,)} \emph{$\varphi$-meßbar über $M$}, wenn gilt $\chi_M \, \hat{f} \in \mathcal{M}_{\K}(\R^n,\varphi)$.
Letzteres ist offenbar gleichbedeutend mit der Existenz zweier Folgen $(T_i)_{i \in \N}, (S_i)_{i \in \N}$ in $\mathcal{S}(\R^n)$ und einer $\varphi$-Null\-menge $N$ derart, daß gilt
$$ \forall_{x \in M \setminus N} \, \lim_{i \to \infty} \left( T_i(x) + \i \, S_i(x) \right) = f(x), $$
insbes.\ gehören alle Punkte von $M$ mit Ausnahme einer $\varphi$-Nullmenge zum Definitionsbereich von $f$.
Die Menge solcher Funktionen $f$ sei mit $\boxed{\mathcal{M}_{\K}(M,\varphi)}$ bezeichnet.
\end{itemize}
\item[(ii)] Ist $M$ eine Teilmenge von $\R^n$, so heißt $M$ genau dann \emph{$\varphi$-meßbar}, wenn gilt $\chi_M \in \mathcal{M}_{\R}(\R^n,\varphi)$.
\end{itemize}
\end{Def}

\begin{Satz} \label{FA.4.69} 
Sind $\varphi$ ein Quadermaß auf $\R^n$ und $M$ eine Teilmenge von $\R^n$, so folgt:
\begin{itemize}
\item[(i)] $\mathcal{M}_{\K}(M,\varphi)$ ist eine assoziative $\K$-Algebra, insbes.\ also ein $\K$-Vektorraum.
\item[(ii)] $\L_{\K}(M,\varphi) \subset \mathcal{M}_{\K}(M,\varphi)$.
\item[(iii)] $f, \tilde{f} \in \mathcal{M}_{\R}(M,\varphi)$ $\Longrightarrow$ $\sup(f,\tilde{f}), \, \inf(f,\tilde{f}), \, f^+, \, f^-, \, |f| \, \in \mathcal{M}_{\R}(M,\varphi).$
\item[(iv)] $f \in \mathcal{M}_{\C}(M,\varphi) \Longrightarrow |f| \in \mathcal{M}_{\R}(M,\varphi)$.
\item[(v)] $(f \in \mathcal{M}_{\K}(M,\varphi) \, \wedge \, f \ne_{\varphi} 0 \mbox{ auf $M$})$ $\Longrightarrow$ $\D \frac{1}{f} \in \mathcal{M}_{\K}(M,\varphi)$.
\item[(vi)] Sind $f \in \mathcal{M}_{\K}(M,\varphi)$ und $f^*$ eine auf einer Teilmenge von $\R^n$ definierte $\K$-wertige Funktion mit $f^* =_{\varphi} f$ auf $M$, so folgt $f^* \in \mathcal{M}_{\K}(M,\varphi)$.
\item[(vii)] $f \in \mathcal{M}_{\K}(M,\varphi) \Longleftrightarrow \forall_{\widetilde{M} \in \mathfrak{P}(M)} \, f \in \mathcal{M}_{\K}(\widetilde{M},\varphi)$.
\end{itemize}
\end{Satz}

\textit{Beweis.} (i) sowie (iii) werden aus Stetigkeitsgründen durch die entsprechenden Eigenschaften der Treppenfunktionen in \ref{FA.4.8} (iv) sowie (v) impliziert, und (vi) sowie (vii) sind trivial.

Zu (ii): Es genügt, den Fall $\K = \R$ zu betrachten.
Sei $f \in \L_{\R}(M,\varphi)$.
Dann existieren $g,h \in \mathcal{S}_{\nearrow}(\R^n,\varphi)$ mit $\chi_M \, \hat{f} =_{\varphi} g - h$ auf $\R^n$, also auch zwei Folgen $(T_i)_{i \in \N}, (S_i)_{i \in \N}$ in $\mathcal{S}(\R^n)$ derart, daß auf $\R^n$ gilt $\chi_M \, \hat{f} =_{\varphi} \lim_{i \to \infty} (T_i - S_i)$, d.h.\ $f \in \mathcal{M}_{\R}(M,\varphi)$.

Zu (iv): Es existierten also Folgen $(T_i)_{i \in \N}$ bzw.\ $(S_i)_{i \in \N}$ in $\mathcal{S}(\R^n)$, die $\varphi$-fast überall auf $M$ gegen ${\rm Re} \, f$ bzw.\ ${\rm Im} \, f$ konvergieren.
Aus \ref{FA.4.8} (iii) folgt dann offenbar, daß $(\sqrt{{T_i}^2 + {S_i}^2})_{i \in \N}$ eine Folge in $\mathcal{S}(\R^n)$ ist, die aus Stetigkeitsgründen $\varphi$-fast überall auf $M$ gegen $|f|$ konvergiert.

Zu (v): Es existierten wieder Folgen $(T_i)_{i \in \N}$ bzw.\ $(S_i)_{i \in \N}$ in $\mathcal{S}(\R^n)$, die $\varphi$-fast überall auf $M$ gegen ${\rm Re} \, f$ bzw.\ ${\rm Im} \, f$ konvergieren.
Beachte im folgenden, daß die Folge $(S_i)_{i \in \N}$ im Falle $\K = \R$ als konstant vom Wert Null gewählt werden kann.
Nach \ref{FA.4.8} (iii) werden durch
\begin{gather*}
\forall_{i \in \N} \forall_{x \in \R^n} \, \widetilde{T}_i(x) := \left\{ \begin{array}{cc} \frac{T_i(x)}{T_i(x)^2 + S_i(x)^2}, & \mbox{falls } (T_i(x),S_i(x)) \ne (0,0), \\ 0, & \mbox{falls } (T_i(x),S_i(x)) = (0,0), \end{array} \right. \\
\forall_{i \in \N} \forall_{x \in \R^n} \, \widetilde{S}_i(x) := \left\{ \begin{array}{cc} \frac{S_i(x)}{T_i(x)^2 + S_i(x)^2}, & \mbox{falls } (T_i(x),S_i(x)) \ne (0,0), \\ 0, & \mbox{falls } (T_i(x),S_i(x)) = (0,0), \end{array} \right.
\end{gather*}
Folgen $(\widetilde{T}_i)_{i \in \N}, (\widetilde{S}_i)_{i \in \N}$ in $\mathcal{S}(\R^n)$ definiert.
Dann gilt aus Stetigkeitsgründen wegen $f \ne_{\varphi} 0$ auf $M$
$$ \lim_{i \to \infty} \widetilde{T}_i - \i \, \widetilde{S}_i =_{\varphi} \frac{1}{f} \mbox{ auf $M$}. $$

Damit ist der Satz vollständig bewiesen. \q

\begin{HS} \label{FA.4.70} 
Seien $\varphi$ ein Quadermaß auf $\R^n$, $M$ eine Teilmenge von $\R^n$ und $f \in \mathcal{M}_{\K}(M,\varphi)$.

Dann gilt:
\begin{itemize}
\item[(i)] $(\exists_{g \in \L_{\R}(M,\varphi)} \, |f| \le_{\varphi} g \mbox{ auf } M) \Longrightarrow f \in \L_{\K}(M,\varphi)$.
\item[(ii)] $|f| \in \L_{\R}(M,\varphi) \Longrightarrow f \in \L_{\K}(M,\varphi)$.
\item[(iii)] Sind $C \in \R$ mit $|f| \le_{\varphi} C$ auf $M$ und $g \in \L_{\K}(M,\varphi)$, so folgt $f \, g \in \L_{\K}(M,\varphi)$.
\item[(iv)] $\widetilde{M}$ $\varphi$-meßbare Teilmenge von $M$ und $f \in \L_{\K}(M,\varphi)$ \\
$\Longrightarrow$ $f \in \L_{\K}(\widetilde{M},\varphi)$.

Insbesondere ist jede $\varphi$-meßbare Teilmenge einer $\varphi$-integrierbaren Menge $\varphi$-meßbar.
\end{itemize}
\end{HS}

\textit{Beweis.} (ii) folgt trivial aus (i), und (iv) ergibt sich aus (iii).

Zu (i): Wegen $|{\rm Re} f| \le |f|$ und $|{\rm Im} f| \le |f|$ genügt es offenbar, den Fall $\K = \R$ zu betrachten.
Außerdem können wir durch Übergang von $f$ und $g$ zu $\chi_M \, \hat{f}$ und $\chi_M \, \hat{g}$ ohne Einschränkung $M = \R^n$ annehmen.
Dann existieren eine $\varphi$-Nullmenge $N \subset \R^n$ und eine Folge $(T_i)_{i \in \N}$ in $\mathcal{S}(\R^n)$ mit
\begin{gather}
\forall_{x \in \R^n \setminus N} \, \lim_{i \to \infty} T_i(x) = f(x), \label{FA.4.70.1} \\
\forall_{x \in \R^n \setminus N} \, f(x) \le |f(x)| \le g(x). \label{FA.4.70.2}
\end{gather}

Für jedes $i \in \N$ definieren wir eine Funktion $f_i \: \R^n \setminus N \to \R$ durch 
$$ \forall_{x \in \R^n \setminus N} \, f_i(x) := \sup \{ -g(x), \inf \{ g(x), T_i(x) \} \}. $$
Dann folgt für alle $i \in \N$
\begin{gather} 
f_i \in \L_{\R}(\R^n,\varphi), \label{FA.4.70.3} \\
|f_i| \le_{\varphi} g \mbox{ auf } \R^n. \label{FA.4.70.4}
\end{gather}
Es gilt
\begin{equation} \label{FA.4.70.5}
\forall_{x \in \R^n \setminus N} \, \lim_{i \to \infty} f_i(x) = f(x).
\end{equation}

{[} Zu (\ref{FA.4.70.5}): Für alle $x \in \R^n \setminus N$ gilt wegen (\ref{FA.4.70.1}), der Stetigkeit der Funktion $\inf \{ g(x), \ldots \} \: \R \to \R$ in $f(x)$ und (\ref{FA.4.70.2}) sowie der Stetigkeit der Funktion $\sup \{ -g(x), \ldots \} \: \R \to \R$ in $f(x)$
$$ \lim_{i \to \infty} f_i(x) = \sup \{ -g(x), f(x) \}, $$
also nach (\ref{FA.4.70.2}): $\lim_{i \to \infty} f_i(x) = f(x)$. {]}

(\ref{FA.4.70.3}) - (\ref{FA.4.70.5}) und der Grenzwertsatz von \textsc{Lebesgue} ergeben $f \in \L_{\R}(\R^n,\varphi)$.

Zu (iii): Aus $f \in \mathcal{M}_{\K}(M,\varphi)$ und $g \in \L_{\K}(M,\varphi) \subset \mathcal{M}_{\K}(M,\varphi)$ folgt zunächst mittels \ref{FA.4.69} (i)
$$ f \, g \in \mathcal{M}_{\K}(M,\varphi) $$
und sodann aus $|f \, g| \le_{\varphi} |C| \, |g| \in \L_{\R}(M,\varphi)$ und (i) die Konklusion der ersten Aussage von (iii). 
Setzt man $f := 1_{\R^n} \in \L_{\R}(M,\varphi)$, so folgt die zweite Aussage von (iii) aus der ersten. \q

\begin{Bem*} 
Die Voraussetzung der $\varphi$-Beschränktheit in (iii) ist notwendig:
Aus \ref{FA.4.41} folgt nämlich $\D f := \frac{1}{\sqrt{x}} \in \L_{\R}({]}0,1{]}, \mu_1)$ und $f^2 \notin \L_{\R}({]}0,1{]}, \mu_1)$.
\end{Bem*}

\begin{Bsp} \label{FA.4.49.M} 
Die Teilmenge $M$ von $[0,1]$ mit $\R = M \oplus \Q$ wie in \ref{FA.4.49} ist nach \ref{FA.4.70} (iv) nicht $\mu_1$-meßbar.
\end{Bsp}

Der folgende Hauptsatz zeigt, daß eine weitere Grenzwertbildung fast überall konvergenter Folgen Lebesgue-meßbarer Funktionen keine größere Menge Le\-bes\-gue-meßbarer Funktionen liefert.
Insbesondere ist die Definition der Lebesgue-Meßbarkeit über $\R^n$ à la \textsc{Dombrowski} in \cite{Dom} gleichbedeutend mit der in \ref{FA.4.68} gegebenen.\footnote{Sei $\varphi$ ein Quadermaß auf $\R^n$.
Eine reellwertige Funktion $f$ heißt in \cite{Dom} genau dann \emph{$\varphi$-meßbar über $\R^n$}, wenn eine Folge in $\mathcal{L}_{\R}(\R^n,\varphi)$ existiert, die $\varphi$-fast überall auf $\R^n$ gegen $f$ konvergiert.
Für eine beliebige Teilmenge $M$ von $\R^n$ ist die Definition der $\varphi$-Meßbarkeit von $f$ über $M$ in loc.\ cit.\ allerdings spezieller als unsere; es wird nämlich zusätzlich gefordert, daß $M$ eine $\varphi$-meßbare Menge ist.}

\begin{HS} \label{FA.4.69.viii} 
Seien $\varphi$ ein Quadermaß auf $\R^n$, $M$ eine Teilmenge von $\R^n$ und $(f_i)_{i \in \N}$ eine $\varphi$-konvergente Folge auf $M$ in $\mathcal{M}_{\K}(M,\varphi)$.

Dann gilt $\D \lim_{i \to \infty} f_i \in \mathcal{M}_{\K}(M,\varphi)$.
\end{HS}

Wir bereiten den Beweis des Hauptsatzes durch folgendes Lemma vor.

\begin{Lemma} \label{FA.4.69.L} 
Ist $\varphi$ ein Quadermaß auf $\R^n$, so existiert $h \: \R^n \to \R_+$ derart, daß $h \in \mathcal{L}_{\R}(\R^n,\varphi)$ gilt.
\end{Lemma}

\textit{Beweisskizze des Lemmas.} Wähle eine Folge $(Q_i)_{i \in \N}$ paarweise disjunkter Quader mit 
$$ \forall_{i \in \N} \, \left( Q_i = \left[ a_{i_1}, a_{i_1} + 1 \right[ \times \ldots \times \left[ a_{i_n}, a_{i_n} + 1 \right[ ~ \wedge ~ \forall_{j \in \{1, \ldots, n\}} \, a_{i_j} \in \Z \right) $$
und $\bigcupdot_{i \in \N} Q_i = \R^n$.
Setze für $j \in \N$
$$ \alpha_j := \left\{ \begin{array}{cc} 0, & \mbox{falls } \varphi(Q_j) = 0, \\ \frac{1}{2^j \, \varphi(Q_j)}, & \mbox{falls } \varphi(Q_j) \ne 0. \end{array} \right. $$
Dann erfüllt $h := \lim_{i \to \infty} T_i$, wobei $\forall_{i \in \N} \, T_i := \sum_{j=0}^i \alpha_j \chi_{Q_j}$, nach dem Grenzwertsatz von \textsc{Levi} die Behauptung des Lemmas. \q
\A
\textit{Beweis des Hauptsatzes.} Es genügt, den Fall $\K = \R$ zu betrachten.
Wir wählen $h$ wie im letzten Lemma und definieren auf $M$ außerhalb einer $\varphi$-Nullmenge
$$ \forall_{i \in \N} \,  g_i := \frac{h \, f_i}{h + |f_i|} ~~ \mbox{ und } ~~ g := \frac{h \, f}{h + |f|}. $$
Der Reihe nach folgt $\lim_{i \to \infty} g_i =_{\varphi} g$, $\forall_{i \in \N} \, |g_i| <_{\varphi} h$, $\forall_{i \in \N} \, g_i \in \L_{\R}(M,\varphi)$ mittels \ref{FA.4.70} (i), also $g \in \L_{\R}(M,\varphi)$ nach dem Grenzwersatz von \textsc{Lebesgue}, $|g| <_{\varphi} h$ und schließlich $f =_{\varphi} \frac{h \, g}{h - |g|} \in \mathcal{M}_{\R}(M,\varphi)$. \q

\begin{Satz} \label{FA.4.69.F} 
Es sei $\varphi$ ein Quadermaß auf $\R^n$.

Dann gilt:
\begin{itemize}
\item[(i)] Jede $\varphi$-integrierbare Teilmenge von $\R^n$ ist $\varphi$-meßbar.

\item[(ii)] Jede offene Teilmenge von $\R^n$ ist $\varphi$-meßbar.
\item[(iii)] Für je zwei $\varphi$-meßbare Mengen $\widetilde{M}, M$ ist auch $\widetilde{M} \setminus M$ $\varphi$-meßbar; insbesondere ergibt (i), daß jede abgeschlossene Teilmenge von $\R^n$ eine $\varphi$-meßbare Menge ist.
\item[(iv)] Ist $M_i$ für jedes $i \in \N$ eine $\varphi$-meßbare Menge, so sind $\D \bigcup_{i \in \N} M_i$ und $\D \bigcap_{i \in \N} M_i$ $\varphi$-meßbare Mengen.
\end{itemize}
\end{Satz}

\textit{Beweis.} (i) folgt sofort aus \ref{FA.4.69} (ii).

Zu (ii): Seien $M \subset \R^n$ offen und für jedes $i \in \N_+$
$$ Q_i := U_i^{\| \ldots \|_{\infty}}(0) = {{]} -i, i {[}}\,^n, $$
also ist $M \cap Q_i$ offen sowie beschränkt und somit nach \ref{FA.4.51} $\varphi$-integrierbar, d.h.\ es gilt $\chi_{{M \cap Q_i}} \in \L_{\R}(\R^n,\varphi) \subset \mathcal{M}_{\R}(\R^n,\varphi)$.
Aus \ref{FA.4.69.viii} und $\lim_{i \to \infty} \chi_{{M \cap Q_i}} = \chi_M$ ergibt sich (ii).

Zu (iii): Die Behauptung folgt aus $\chi_{{{\widetilde{M}} \setminus M}} = (\chi_{\widetilde{M}} - \chi_M)^+$ und der Stetigkeit der Funktion $\ldots^+ \: \R \to \R$.
Beachte, daß wegen der $\varphi$-Meßbarkeit von $\widetilde{M}, M$ eine Folge $(T_i)_{i \in \N}$ in $\mathcal{S}(\R^n)$ mit $\chi_{\widetilde{M}} - \chi_M =_{\varphi} \lim_{i \to \infty} T_i$ auf $\R^n$ existiert.

Zu (iv): Wir definieren Folgen $(f_i)_{i \in \N}$ und $(g_i)_{i \in \N}$ durch
\begin{gather*}
\forall_{i \in \N} \, f_i := \chi_{\bigcup_{\iota=1}^i M_{\iota}} = \sup( \chi_{M_0}, \ldots, \chi_{M_i} ) \stackrel{\ref{FA.4.69} (iii)}{\in} \mathcal{M}_{\R}(M,\varphi), \\
\forall_{i \in \N} \, g_i := \chi_{\bigcap_{\iota=1}^i M_{\iota}} = \inf( \chi_{M_0}, \ldots, \chi_{M_i} ) \stackrel{\ref{FA.4.69} (iii)}{\in} \mathcal{M}_{\R}(M,\varphi),
\end{gather*}
für die $\lim_{i \to \infty} f_i = \chi_{\bigcup_{i=1}^{\infty} M_i}$ sowie $\lim_{i \to \infty} g_i = \chi_{\bigcap_{i=1}^{\infty} M_i}$ gilt.
\ref{FA.4.69.viii} ergibt die Konklusion von (iv). \q

\begin{Satz}[Lebesgue-Meßbarkeit stetiger Funktionen über Lebesgue-meßbaren Mengen] \label{FA.4.70.F} $\,$

\noindent \textbf{Vor.:} Seien $\varphi$ ein Quadermaß auf $\R^n$, $M$ eine $\varphi$-meßbare Teilmenge von $\R^n$ und $f \: M \to \K$ eine stetige Funktion.

\noindent \textbf{Beh.:} $f \in \mathcal{M}_{\K}(M,\varphi)$.
\end{Satz}

\textit{Beweis.} Es genügt erneut, den Fall $\K = \R$ zu betrachten.
Wegen der $\varphi$-Meßbarkeit von $M$ ist
$$ \mbox{$\chi_M \in \mathcal{M}_{\R}(\R^n,\varphi)$ beschränkt.} $$
Des weiteren ist für jedes $i \in \N$
$$ Q_i := U^{\| \ldots \|_{\infty}}_i(0) = {]} -i, i {[}^n \in \mathfrak{Q}_n $$
eine $\varphi$-integrierbare Menge, d.h.\
$$ \chi_{Q_i} \in \L_{\R}(\R^n,\varphi). $$
Somit folgt aus \ref{FA.4.70} (iii)
$$ \chi_{{M \cap Q_i}} = \chi_M \, \chi_{Q_i} \in \L_{\R}(\R^n,\varphi), $$
also ist $M \cap Q_i$ eine beschränkte $\varphi$-integrierbare Menge.
Aus Stetigkeitsgründen ergibt \ref{FA.4.59} nun zusammen mit \ref{FA.4.22} (v)
$$ g_i := \sup( \inf( f|_{M \cap Q_i}, i), -i ) \in \L_{\R}(M \cap Q_i, \varphi) \mbox{ beschränkt}, $$
d.h.\
$$ \underbrace{\chi_{{M \cap Q_i}} \, \hat{g}_i}_{\stackrel{i \to \infty}{\longrightarrow} \chi_M \, \hat{f}} \in \L_{\R}(\R^n, \varphi) \subset \mathcal{M}_{\R}(\R^n,\varphi), $$
womit die Behauptung aus \ref{FA.4.69.viii} folgt. \q
\A
In gewisser Weise stellt der folgende Hauptsatz die Umkehrung des Satzes von \textsc{Fubini} dar.

\begin{HS}[Satz von \textsc{Tonelli}] \index{Satz!von \textsc{Tonelli}} \label{FA.4.71} $\,$

\noindent \textbf{Vor.:} Es seien $n_1, n_2 \in \N_+$ mit $n = n_1 + n_2$, $\varphi_i$ Quadermaße auf $\R^{n_i}$ für $i \in \{1,2\}$ und $f \in \mathcal{M}_{\K}(\R^n,\varphi_1 \times \varphi_2)$.
Ferner existiere
\begin{equation} \label{FA.4.71.S}
\int\limits_{\R^{n_1}} \left( \; \int\limits_{\R^{n_2}} |f(x_1,x_2)| \, \d \varphi_2(x_2) \right) \, \d \varphi_1(x_1) \in \R.
\end{equation}

\noindent \textbf{Beh.:} $f \in \L_{\K}(\R^n,\varphi_1 \times \varphi_2)$.
\end{HS}

\textit{Beweis.} 1.\ Fall: $\K = \R$.
Wir wählen gemäß \ref{FA.4.69.L} eine Abbildung 
\begin{equation} \label{FA.4.71.1}
g \in \L_{\R}(\R^n,\varphi_1 \times \varphi_2) \cap \R^{\R^n} \mbox{ mit } g > 0
\end{equation}
und definieren für jedes $i \in \N$
\begin{equation} \label{FA.4.71.2}
f_i := \inf ( |f|, i \, g ).
\end{equation}

Es gilt
\begin{gather}
\lim_{i \to \infty} f_i =_{\varphi_1 \times \varphi_2} |f|, \label{FA.4.71.3} \\
\forall_{i \in \N} \, f_i \in \L_{\R}(\R^n,\varphi_1 \times \varphi_2), \label{FA.4.71.4} \\
\forall_{i \in \N} \, f_i \le_{\varphi_1 \times \varphi_2} f_{i+1} \le_{\varphi_1 \times \varphi_2} |f|, \label{FA.4.71.5} \\
\exists_{C \in \R} \forall_{i \in \N} \, \int\limits_{\R^n} f_i \, \d (\varphi_1 \times \varphi_2) \le C. \label{FA.4.71.6}
\end{gather}

{[} (\ref{FA.4.71.3}) sowie (\ref{FA.4.71.5}) sind wegen (\ref{FA.4.71.2}) und (\ref{FA.4.71.1}) klar.

Zu (\ref{FA.4.71.4}): Aus (\ref{FA.4.71.1}), (\ref{FA.4.71.2}) folgt zunächst mittels \ref{FA.4.69} (ii), (iii) für jedes $i \in \N$ $f_i \in \mathcal{M}_{\R}(M,\varphi_1 \times \varphi_2)$ und sodann mittels \ref{FA.4.70} (i) die Aussage (\ref{FA.4.71.4}).

Zu (\ref{FA.4.71.6}): Aus (\ref{FA.4.71.4}) und dem Satz von \textsc{Fubini} ergibt sich wegen der Existenz des Integrales (\ref{FA.4.71.S}) und (\ref{FA.4.71.5}) die Aussage (\ref{FA.4.71.6}), wenn $C$ als die reelle Zahl in (\ref{FA.4.71.S}) definiert wird. {]}

(\ref{FA.4.71.3}) - (\ref{FA.4.71.6}) und der Satz von \textsc{Levi} ergeben $|f| \in \L_{\R}(\R^n,\varphi_1 \times \varphi_2)$.
Hieraus folgt mittels \ref{FA.4.70} (ii) die Behauptung des Hauptsatzes.

2.\ Fall: $\K = \C$.
Die Behauptung folgt durch Anwendung des 1.\ Falles auf ${\rm Re} \, f$ und ${\rm Im} \, f$, die betraglich kleiner oder gleich $|f|$ sind. \q

\subsection*{Anhang: Meßbarkeit im Sinne der Maßtheorie} \addcontentsline{toc}{subsection}{Anhang: Meßbarkeit im Sinne der Maßtheorie}

Zum Abschluß dieses Kapitels gehen wir auf die Definition der Meßbarkeit, die in der \emph{Maßtheorie} verwendet wird, ein.

\begin{Def}[(Borel-)$\sigma$-Algebren, Borel-Mengen] \label{FA.4.M1} 
Sei $X$ eine nicht-leere Menge.
\begin{itemize}
\item[(i)] Eine Teilmenge $\mathcal{A}$ von $\mathfrak{P}(X)$ heißt \emph{$\sigma$-Algebra über $X$}\index{Algebra!$\sigma$-} genau dann, wenn gilt
\begin{itemize}
\item[(a)] $X \in \mathcal{A}$,
\item[(b)] $\forall_{A \in \mathcal{A}} \, X \setminus A \in \mathcal{A}$ und
\item[(c)] $\D \forall_{(A_i)_{i \in \N} \in \mathcal{A}^{\N}} \, \bigcup_{i \in \N} A_i \in \mathcal{A}$.
\end{itemize}
Dann gilt offenbar auch
\begin{itemize}
\item[(d)] $\emptyset = X \setminus X \in \mathcal{A}$,
\item[(e)] $\D \forall_{(A_i)_{i \in \N} \in \mathcal{A}^{\N}} \, \bigcap_{i \in \N} A_i = X \setminus \left( \bigcup_{i \in \N} (X \setminus A_i) \right) \in \mathcal{A}$ sowie
\item[(f)] $\forall_{A,B \in \mathcal{A}} \, A \setminus B = A \cap (X \setminus B) \in \mathcal{A}$.
\end{itemize}
\begin{Bsp*} $\,$
\begin{itemize}
\item[1.)] $\{ \emptyset, X \}$ und $\mathfrak{P}(X)$ sind $\sigma$-Algebren über $X$.
\item[2.)] Der Schnitt beliebig vieler $\sigma$-Algebren über $X$ ist ebenfalls eine $\sigma$-Algebra über $X$.
\item[3.)] Sind $\mathcal{A}$ eine $\sigma$-Algebra über $X$ und $M \in \mathcal{A} \setminus \{ \emptyset \}$, so ist $\mathcal{A} \cap \mathfrak{P}(M)$ eine $\sigma$-Algebra über $M$.
\end{itemize}
\end{Bsp*}
Für eine $\sigma$-Algebra $\mathcal{A}$ über $X$ nennt man das Paar $(X,\mathcal{A})$ einen \emph{Meßraum}\index{Meßraum}\index{Raum!Meß-}.
\item[(ii)] Eine nicht-leere Teilmenge $\mathfrak{S}$ von $\mathfrak{P}(X)$ erzeugt offenbar kanonisch eine minimale $\sigma$-Algebra $\boxed{\sigma(\mathfrak{S})}$ über $X$, die $\mathfrak{S}$ enthält, nämlich
$$ \sigma(\mathfrak{S}) := \bigcap_{\substack{\mathcal{A} \text{ $\sigma$-Algebra über $X$} \\ \mathfrak{S} \subset \mathcal{A}}} \mathcal{A}. $$
\item[(iii)] Ist $X$ sogar ein toplogischer Raum, so heißt $\boxed{\B(X)} := \sigma({\rm Top}(X))$ die \emph{Borel-$\sigma$-Algebra von $X$}\index{Algebra!$\sigma$-!Borel-}, deren Elemente \emph{Borel-Mengen von $X$}\index{Menge!Borel-} -- kurz \emph{Borelsch} -- oder \emph{meßbar}\index{Menge!meßbare} heißen.
\begin{Bsp*}
Sind $X$ ein topologischer Raum und $M \in \B(X) \setminus \{ \emptyset \}$, so gilt $\B(M) \subset \B(X) \cap \mathfrak{P}(M)$.
\end{Bsp*}
\item[(iv)] Ist $X = \R^n$, so setzen wir $\boxed{\B^n} := \B(\R^n)$.
(Im Falle $n=1$ schreibt man natürlich $\boxed{\B}$ anstelle von $\B^1$.)
\end{itemize}
\end{Def}

\begin{Def}[Meßbare Abbildungen] \index{Abbildung!meßbare} \label{FA.4.M2} 
Es sei $(X,\mathcal{A})$ ein Meßraum.
\begin{itemize}
\item[(i)] Ist $(Y,\mathcal{B})$ ein weiterer Meßraum, so heißt eine Abbildung $f \: X \to Y$ genau dann \emph{$(\mathcal{A}, \mathcal{B})$-meßbar}, wenn gilt $\forall_{B \in \mathcal{B}} \, \overline{f}^1(B) \in \mathcal{A}$.
\begin{Bsp*}
Seien $(X,\mathcal{A}), (Y,\mathcal{B}), (Z,\mathcal{C})$ Meßräume.
Dann ist die Komposition $f \circ g \: X \to Z$ $(\mathcal{A}, \mathcal{B})$-meßbarer bzw.\ $(\mathcal{B}, \mathcal{C})$-meßbarer Abbildungen $g \: X \to Y$ bzw.\ $f \: Y \to Z$ eine $(\mathcal{A}, \mathcal{C})$-meßbare Abbildung.
\end{Bsp*}
\item[(ii)] Ist $Y \ne \emptyset$ ein topologischer Raum, so heißt eine Abbildung $f \: X \to Y$ genau dann \emph{$\mathcal{A}$-meßbar}, wenn $f \: X \to Y$ eine $(\mathcal{A}, \B(Y))$-meßbare Abbildung ist.

Falls auch $X$ ein topologischer Raum und $\mathcal{A} = \B(X)$ ist, so nennen wir $f$ \emph{meßbar}.
\end{itemize}
\end{Def}

\begin{Satz} \label{FA.4.M2S}
Es seien $(X,\mathcal{A})$ ein Meßraum, $Y$ eine nicht-leere Menge und $\mathfrak{S}$ eine nicht-leere Teilmenge von $\mathfrak{P}(Y)$.
Ferner sei $f \: X \to Y$ eine Abbildung mit $\forall_{S \in \mathfrak{S}} \, \overline{f}^1(S) \in \mathcal{A}$.

Dann ist $f$ $(\mathcal{A},\sigma(\mathfrak{S}))$-meßbar.
\end{Satz}

\textit{Beweis.} $\mathcal{B} := \{ B \in \mathfrak{P}(Y) \, | \, \overline{f}^1(B) \in \mathcal{A} \}$ ist offenbar eine $\sigma$-Algebra über $Y$, die nach Voraussetzung $\mathfrak{S}$ enthält.
Daher gilt $\sigma(\mathfrak{S}) \subset \mathcal{B}$.
Hieraus folgt die Behauptung. \q

\begin{Kor} \label{FA.4.M2K} 
Es seien $X, Y$ nicht-leere topologische Räume und $f \: X \to Y$ eine stetige Abbildung.

Dann ist $f$ meßbar. \q
\end{Kor}

\begin{Satz} \label{FA.4.M6.L} 
Es seien $(X,\mathcal{A})$ ein Meßraum und $f \in \R^X$.

Dann sind die folgenden Aussagen paarweise äquivalent:
\begin{itemize}
\item[(i)] $f$ ist $\mathcal{A}$-meßbar.
\item[(ii)] $\forall_{b \in \R} \, \overline{f}^1( {]} -\infty, b {[} ) \in \mathcal{A}$.
\item[(iii)] $\forall_{b \in \R} \, \overline{f}^1( {]} -\infty, b {]} ) \in \mathcal{A}$.
\end{itemize}
\end{Satz}

\textit{Beweis.} Seien $\mathfrak{S}_1 := \{ \, {]} -\infty, b {[} \, | \, b \in \R \}$ und $\mathfrak{S}_2 := \{ \, {]} -\infty, b {]} \, | \, b \in \R \}$ sowie $\mathcal{B}_i := \sigma(\mathfrak{S}_i)$ für $i \in \{1,2\}$, d.h.\ insbesondere $\mathcal{B}_1 \subset \B$.
Da für jedes $b \in \R$ gilt
\begin{gather*}
{]} - \infty, b {]} = \bigcap_{k \in \N_+} {]} - \infty, b + \frac{1}{k}{[}  \in \mathcal{B}_1, \\
{]} - \infty, b {[} = \bigcup_{k \in \N_+} {]} - \infty, b - \frac{1}{k}{]}  \in \mathcal{B}_2, 
\end{gather*}
folgt $\mathcal{B}_1 = \mathcal{B}_2$.
Ferner gilt für alle $a,b \in \R$ mit $a<b$
$$ {]} a, b {[} = {]} - \infty, b {[} \setminus {]} - \infty, a {]} \in \mathcal{B}_1 = \mathcal{B}_2, $$
also folgt zunächst ${\rm Top}(\R) \subset \mathcal{B}_1 \subset \B$ und sodann $\B = \mathcal{B}_1 = \mathcal{B}_2$.
Der Satz ergibt sich nun aus \ref{FA.4.M2S}. \q
\A
Wir wollen nun natürlich ein Quadermaß auf $\R^n$ ins Spiel bringen und orientieren uns bei der Entwickelung der Theorie an \cite[Abschnitte 5.5 und 5.6]{Reck}.

\begin{Def}[Lebesgue-meßbare Mengen i.S.d.\ Maßtheorie] \index{Menge!meßbare} \label{FA.4.M3} 
Sei $\varphi$ ein Quadermaß auf $\R^n$.
$\boxed{\mathfrak{I}(\R^n,\varphi)}$ bezeichne die Menge der $\varphi$-integrierbaren Teilmengen von $\R^n$.\index{Menge!integrierbare}

Wir setzen $\boxed{\B(\R^n,\varphi)} := \sigma(\mathfrak{I}(\R^n,\varphi))$ und nennen die Elemente dieser Menge \emph{$\varphi$-meßbar i.S.d.\ Maßtheorie}.

\begin{Bem*} 
Um Verwechselungen zu vermeiden, werden wir eine i.S.d.\ Maßtheorie $\varphi$-meßbare Menge zunächst nicht abkürzend ,,$\varphi$-meßbar`` nennen.
In \ref{FA.4.M6.K2} werden wir sehen, daß die $\varphi$-meßbaren Teilmengen von $\R^n$ genau die i.S.d.\ Maßtheorie $\varphi$-meßbaren Teilmengen von $\R^n$ sind.
\end{Bem*}
\end{Def}

\begin{Satz} \label{FA.4.72} 
Es sei $\varphi$ ein Quadermaß auf $\R^n$.

Dann ist die Menge $\boxed{\mathfrak{M}(\R^n,\varphi)}$ der $\varphi$-meßbaren Teilmengen von $\R^n$ eine $\sigma$-Algebra über $\R^n$, und es gilt $\B^n \subset \B(\R^n,\varphi) \subset \mathfrak{M}(\R^n,\varphi)$.\index{Menge!meßbare}
\end{Satz}

\textit{Beweis.} Daß $\mathfrak{M}(\R^n,\varphi)$ eine $\sigma$-Algebra ist folgt aus \ref{FA.4.69.F} (ii) - (iv).
\ref{FA.4.69} (i) und die Definition von $\B(\R^n,\varphi)$ ergeben dann $\B(\R^n,\varphi) \subset \mathfrak{M}(\R^n,\varphi)$.
Ferner ist nach \ref{FA.4.53} offenbar jede offene Teilmenge von $\R^n$ eine i.S.d.\ Maßtheorie $\varphi$-meßbare Teilmenge von $\R^n$, also ergibt die Definition von $\B^n$ auch $\B^n \subset \B(\R^n, \varphi)$. \q

\begin{Bem*}  
Aus \ref{FA.4.49.M} folgt $\mathfrak{M}(\R,\mu_1) \subsetneqq \mathfrak{P}(\R)$.
Es gilt jedoch
\begin{equation} \label{FA.4.72.S}
\# \B(\R,\mu_1) = \# \mathfrak{P}(\R).
\end{equation}

{[} Beweis hiervon: $C$ bezeichne das Cantorsche Diskontinuum wie in \ref{FA.4.13.B} 4.), also ist $C$ eine überabzählbare $\mu_1$-Nullmenge.
Dann ist auch jede Teilmenge von $C$ eine $\mu_1$-Nullmenge und somit nach \ref{FA.4.30.K2} eine $\mu_1$-integrierbare Menge.
Daher folgt $\# \mathfrak{P}(\R) = \# \mathfrak{P}(C) \le \# \mathfrak{I}(\R,\mu_1) \le \# \B(\R,\mu_1) \le \# \mathfrak{P}(\R)$. {]}

Mittels \emph{transfiniter Induktion} kann man 
$$ \# \B = \# \R $$
zeigen, siehe z.B.\ \cite[Korollar II.8.6]{Els}.\footnote{Die in \cite{Els} mit $\mathfrak{B}^p$ bzw.\ $\mathfrak{L}^p$ bezeichneten Mengen stimmen in unserem Sinne mit $\B^p$ bzw.\ $\B(\R^p,\mu_p) = \mathfrak{M}(\R^p,\mu_p)$ überein, und die in loc.\ cit.\ ,,Lebesgue-meßbar`` genannten Teilmengen des $\R^p$ heißen bei uns ,,$\mu_p$-meßbar (i.S.d.\ Maßtheorie)``, vgl. \cite[Sätze III.2.2 und II.7.4]{Els} sowie \ref{FA.4.M4F}, \ref{FA.4.M6.K2} und \ref{FA.4.M4S} unten.}
Der Beweis von (\ref{FA.4.72.S}) impliziert dann, daß es eine Teilmenge des Cantorschen Diskontinuums gibt, die $\mu_1$-integrierbar und nicht Borelsch ist.
\end{Bem*}

\begin{Satz}[Charakterisierung i.S.d.\ Maßtheorie Lebesgue-meßbarer Mengen] \label{FA.4.M4} $\,$

\noindent \textbf{Vor.:} Es sei $\varphi$ ein Quadermaß auf $\R^n$.

\noindent \textbf{Beh.:} Eine Teilmenge $M$ von $\R^n$ ist genau dann $\varphi$-meßbar i.S.d.\ Maßtheorie, wenn $M \cap U$ für jede offene und beschränkte Teilmenge $U$ von $\R^n$ $\varphi$-integrierbar ist.
\end{Satz}

\textit{Beweis.} $\mathcal{B}$ bezeichne die Menge aller Teilmengen $M$ von $\R^n$, die die rechte Seite der Behauptung erfüllen.
Wir zeigen
\begin{gather}
{\rm Top}(\R^n) \subset \mathcal{B}, \label{FA.4.M4.1} \\
\mathcal{B} \mbox{ ist eine $\sigma$-Algebra}, \label{FA.4.M4.2} \\
\B(\R^n,\varphi) \subset \mathcal{B}, \label{FA.4.M4.3} \\
\mathcal{B} \subset \B(\R^n,\varphi), \label{FA.4.M4.4}
\end{gather}
womit der Satz bewiesen ist.

{[} Zu (\ref{FA.4.M4.1}): Seien $M \in {\rm Top}(\R^n)$ sowie $U$ eine offene und beschränkte Teilmenge von $\R^n$.
Dann ist $M \cap U$ eine offene und beschränkte -- also nach \ref{FA.4.51} auch $\varphi$-integrierbare -- Teilmenge von $\R^n$.

Zu (\ref{FA.4.M4.2}): $\R^n \in \mathcal{B}$ folgt aus (\ref{FA.4.M4.1}).
Des weiteren gilt
$$ M \in \mathcal{B} \Longrightarrow \R^n \setminus M \in \mathcal{B}, $$
da für alle $M \in \mathcal{B}$ sowie jede offene und beschränkte Teilmenge $U$ von $\R^n$ aus \ref{FA.4.51}, \ref{FA.4.50} (i) folgt: $(\R^n \setminus M) \cap U = U \setminus (M \cap U) \in \mathfrak{I}(\R^n,\varphi)$.
Zum Nachweis von (\ref{FA.4.M4.2}) bleibt zu zeigen, daß mit $M_i \in \mathcal{B}$ für jedes $i \in \N$ auch $\bigcup_{i \in \N} M_i \in \mathcal{B}$ gilt:
Wir setzen für alle $\iota \in \N$
$$ \widetilde{M}_{\iota} := \bigcup_{i=0}^{\iota} M_i. $$
Für jede offene und beschränkte Teilmenge $U \subset \R^n$ ist $\widetilde{M}_{\iota} \cap U = \bigcup_{i=1}^{\iota} (M_i \cap U)$ dann wegen \ref{FA.4.51}, \ref{FA.4.50} (i) $\varphi$-in\-tegrier\-bar mit $\varphi$-Maß $\le \varphi(U)$, also ist nach \ref{FA.4.50} (ii) (b) offenbar auch $\bigcup_{\iota \in \N} (\widetilde{M}_{\iota}  \cap U ) = \left( \bigcup_{\iota \in \N} \widetilde{M}_{\iota} \right) \cap U = \left( \bigcup_{i \in \N} M_i \right) \cap U$ $\varphi$-integrierbar.

Zu (\ref{FA.4.M4.3}): Der Schnitt zweier $\varphi$-integrierbarer Mengen ist $\varphi$-integrierbar, siehe \ref{FA.4.50} (i), also folgt aus \ref{FA.4.51}, daß $\mathfrak{I}(\R^n,\varphi) \subset \mathcal{B}$ und somit (\ref{FA.4.M4.3}) gilt.

Zu (\ref{FA.4.M4.4}): Sei $M \in \mathcal{B}$.
Dann ist $M \cap {]}-i,i{[}^n$ für jedes $i \in \N$ eine $\varphi$-integrierbare Menge und folglich $M = \bigcup_{i \in \N} (M \cap {]}-i,i{[}^n) \in \sigma(\mathfrak{I}(\R^n,\varphi)) = \B(\R^n,\varphi)$. {]} \q

\begin{Def} \label{FA.4.M4D} 
Es seien $\varphi$ ein Quadermaß auf $\R^n$ und $M$ eine i.S.d.\ Maßtheorie $\varphi$-meßbare Teilmenge von $\R^n$.
Wir setzen dann
$$ \boxed{\varphi(M)} := \lim_{i \to \infty} \varphi(M \cap {]}-i,i{[}^n) \in {[} 0, \infty {]}. $$

\begin{Bem*} 
Im Falle $M \in \mathfrak{I}(\R^n,\varphi)$ stimmt diese Definition von $\varphi(M)$ wegen des letzten Satzes und \ref{FA.4.50} (i) (a) mit der von $\varphi(M) < \infty$ in \ref{FA.4.48} überein.
Nach \ref{FA.4.50} (ii) (b) gilt daher
\begin{equation} \label{FA.4.M4D.S}
\forall_{M \in \B(\R^n,\varphi)} \, \varphi(M) < \infty \Longleftrightarrow M \in \mathfrak{I}(\R^n,\varphi).
\end{equation}
\end{Bem*}
\end{Def}

Wir vereinbaren im folgenden, daß mit $\infty$ naheliegend gerechnet wird.
Für jedes $t \in \R$ setzen wir $t + \infty := \infty + t := \infty - t := \infty$ sowie $\infty + \infty := \infty$.

\begin{Satz} \label{FA.4.M4F} 
Seien $\varphi$ ein Quadermaß auf $\R^n$ und $M_i \in \B(\R^n,\varphi)$ für jedes $i \in \N$.

Dann gilt:
\begin{itemize}
\item[(i)] $M_1 \cup M_2, M_1 \cap M_2, M_1 \setminus M_2 \in \B(\R^n,\varphi)$ und
\begin{gather*}
\varphi(M_1 \cup M_2) + \varphi(M_1 \cap M_2) = \varphi(M_1) + \varphi(M_2), \\
\varphi(M_1) = \varphi(M_1 \setminus M_2) + \varphi(M_1 \cap M_2).
\end{gather*}
\item[(ii)] $M_1 \subset M_2$ $\Longrightarrow$ $\varphi(M_1) \le \varphi(M_2)$.
\item[(iii)] $\D M:= \bigcup_{i \in \N} M_i \in \B(\R^n,\varphi)$ und
\begin{gather*}
\forall_{i \in \N} \, M_i \subset M_{i+1} \Longrightarrow  \varphi(M) = \lim_{i \to \infty} \varphi(M_i), \nonumber \\
\varphi(M) \le \sum_{i=0}^{\infty} \varphi(M_i) \mbox{ mit Gleichheit, falls } \forall_{i,j \in \N, \, i \ne j} \, M_i \cap M_j = \emptyset \mbox{ gilt.} \label{FA.4.M4F.1}
\end{gather*}
\end{itemize}
\end{Satz}

\textit{Beweis als Übung.} \q

\begin{Def}[Lebesgue-meßbare Funktionen i.S.d.\ Maßtheorie] \index{Funktion!meßbare} \label{FA.4.M5} 
Es seien $\varphi$ ein Quadermaß auf $\R^n$ und $M \ne \emptyset$ eine i.S.d.\ Maßtheorie $\varphi$-meßbare Teilmenge von $\R^n$.
Wir definieren eine $\sigma$-Algebra über $M$ durch
$$ \boxed{\B(M, \varphi)} := \B(\R^n, \varphi) \cap \mathfrak{P}(M) $$
und nennen eine Funktion $f \: M \to \K$ \emph{(reellwertig, falls $\K = \R$,) $\varphi$-meßbar über $M$ i.S.d.\ Maßtheorie} genau dann, wenn sie $\B(M, \varphi)$-meßbar ist.
Die Menge solcher Funktionen bezeichnen wir mit $\boxed{\widetilde{\mathcal{M}}_{\K}(M,\varphi)}$.
\pagebreak

\begin{Bem*} $\,$
\begin{itemize}
\item[1.)] Um diese Definition geben zu können, ist es notwendig, $M \ne \emptyset$ als $\varphi$-meßbar i.S.d.\ Maßtheorie vorauszusetzen, da $\B(\R^n, \varphi) \cap \mathfrak{P}(M)$ andernfalls i.a.\ keine $\sigma$-Algebra über $M$ ist.
\item[2.)] Um Verwechselungen zu vermeiden, werden wir eine i.S.d.\ Maßtheorie $\varphi$-meßbare Funktion \underline{nie} abkürzend ,,$\varphi$-meßbar`` nennen.
\end{itemize}
\end{Bem*}
\end{Def}

\begin{Satz} \label{FA.4.M5S} 
Es seien $\varphi$ ein Quadermaß auf $\R^n$ und $M \ne \emptyset$ eine i.S.d.\ Maßtheorie $\varphi$-meßbare Teilmenge von $\R^n$.

Dann gilt:
\begin{itemize}
\item[(i)] $\mathcal{C}(M,\K) \subset \widetilde{\mathcal{M}}_{\K}(M,\varphi)$.
\item[(ii)] $\chi_M \in \widetilde{\mathcal{M}}_{\R}(\R^n,\varphi)$.
\item[(iii)] Sind $f \in \widetilde{\mathcal{M}}_{\K}(M,\varphi)$ und $g \in \K^M$ mit $f =_{\varphi} g$, so gilt $g \in \widetilde{\mathcal{M}}_{\K}(M,\varphi)$.
\end{itemize}
\end{Satz}

\textit{Beweis.} Zu (i): Seien $f \in \mathcal{C}(M,\K)$ und $U \subset \K$ offen.
Dann ist $\overline{f}^1(U)$ eine offene Teilmenge von $M$.
Wegen
$$ {\rm Top}_{\R^n}(M) = {\rm Top}(\R^n) \cap \{M\} \subset \B^n \cap \mathfrak{P}(M) \subset \B(\R^n,\varphi) \cap \mathfrak{P}(M) = \B(M,\varphi) $$
folgt (i) aus \ref{FA.4.M2S}.

Zu (ii): Da für jedes $B \in \B$ gilt
$$ \overline{\chi_M}^1(B) = \left\{\begin{array}{cl} \R^n, & \mbox{falls } \{0,1\} \subset B, \\ \emptyset, & \mbox{falls } \{0,1\} \cap B = \emptyset, \\ M, & \mbox{falls } 0 \notin B \wedge 1 \in B, \\ \R^n \setminus M, & \mbox{falls } 0 \in B \wedge 1 \notin B, \end{array} \right. $$
ergibt sich (ii).

Zu (iii): Sei $N$ eine $\varphi$-Nullmenge mit
$$ \forall_{x \in M \setminus N} \, f(x) = g(x). $$
Dann gilt für jedes $B \in \B(\K)$
\begin{eqnarray*}
\overline{g}^1(B) & = & \left( \overline{g}^1(B) \cap \left( M \setminus N \right) \right) \cup \left( \overline{g}^1(B) \cap N \right) \\
& = & \left( \overline{f}^1(B) \cap \left( M \setminus N \right) \right) \cup \left( \overline{g}^1(B) \cap N \right)
\end{eqnarray*}
Aus $f \in \widetilde{\mathcal{M}}_{\K}(M,\varphi)$, $g \in \K^M$, \ref{FA.4.30.K2} und $\B^n \subset \B(\R^n,\varphi)$ folgt offenbar, daß die rechte Seite der letzten Gleichung ein Element von $\B(M,\varphi)$ ist. \q

\begin{Satz} \label{FA.4.A3} 
Es seien $\varphi$ ein Quadermaß auf $\R^n$, $M \ne \emptyset$ eine i.S.d.\ Maßtheorie $\varphi$-meßbare Teilmenge von $\R^n$ und $m \in \N_+$ sowie $f = (f_1, \ldots, f_m) \: M \to \K^m$ eine Abbildung.

Dann ist $f$ genau dann $\B(M,\varphi)$-meßbar, wenn für jedes $i \in \{1, \ldots, m\}$ die Funktion $f_i \: M \to \K$ i.S.d.\ Maßtheorie $\varphi$-meßbar über $M$ ist.
\end{Satz}

\textit{Beweis.} ,,$\Rightarrow$`` folgt sofort aus \ref{FA.4.M5S} (i) und der Stetigkeit der Projektionen $\pi_i \: \K^m \to \K$, $i \in \{1, \ldots, m \}$.

,,$\Leftarrow$`` Wir setzen
$$ \mathfrak{S} := \left\{ U_{\varepsilon}^{\|\ldots\|_{\infty}}(x) \, | \, x \in \K^m \, \wedge \, \varepsilon \in \R_+ \right\}. $$
Dann gilt $\sigma(\mathfrak{S}) = \B(\K^m)$ und für alle $x = (x_1, \ldots, x_m) \in \K^m$ sowie $\varepsilon \in \R_+$
$$ U_{\varepsilon}^{\|\ldots\|_{\infty}}(x) = \bigtimes_{i = 1}^m U_{\varepsilon}^{|\ldots|} (x_i), $$
also ist
$$ \overline{f}^1 \left( U_{\varepsilon}^{\|\ldots\|_{\infty}} (x) \right) = \bigcap_{i=1}^m \overline{f_i}^1\left( (U_{\varepsilon}^{|\ldots|} (x_i) \right) $$
nach Voraussetzung der rechten Seite eine i.S.d.\ Maßtheorie $\varphi$-meßbare Menge.
Aus der Beliebigkeit von $x \in \K^m$ und $\varepsilon \in \R_+$ sowie \ref{FA.4.M2S} folgt die $\B(M,\varphi)$-Meßbarkeit von $f$. \q

\begin{Kor} \label{FA.4.A3K} 
Seien $\varphi$ ein Quadermaß auf $\R^n$ und $M \ne \emptyset$ eine i.S.d.\ Maßtheorie $\varphi$-meßbare Teilmenge von $\R^n$.

Dann gilt $\D \widetilde{\mathcal{M}}_{\C}(M,\varphi) = \{ f + \i \, g \, | \, f,g \in \widetilde{\mathcal{M}}_{\R}(M,\varphi) \}$. \q
\end{Kor}

\begin{Lemma} \label{FA.4.A4} $\,$

\noindent \textbf{Vor.:} Es seien $\varphi$ ein Quadermaß auf $\R^n$ und $M \ne \emptyset$ eine i.S.d.\ Maßtheorie $\varphi$-meßbare Teilmenge von $\R^n$.
Ferner seien $m \in \N_+$, $U$ eine offene Teilmenge von $\K^m$ sowie $f = (f_1, \ldots, f_m) \: M \to \K^m$ eine $\B(M,\varphi)$-meßbare Abbildung mit $f(M) \subset U$ und $g \: U \to \K$ eine stetige Funktion.

\noindent \textbf{Beh.:} $g \circ f \: M \to \K$ ist eine i.S.d.\ Maßtheorie $\varphi$-meßbare Funktion über $M$.
\end{Lemma}

\textit{Beweis.} Sei $V$ eine offene Teilmenge von $\K$.
Die Stetigkeit von $g$ ergibt, daß $\overline{g}^1(V)$ eine offene Teilmenge der offenen Teilmenge $U$ von $\K^n$ -- also selbst offen in $\K^n$ -- ist.
Aus der $\B(M,\varphi)$-Meßbarkeit der Abbildung $f$ folgt dann $\overline{g \circ f}^1(V) = \overline{f}^1(\overline{g}^1(V)) \in \B(M,\varphi)$.
Somit gilt die Behauptung nach \ref{FA.4.M2S}. \q

\begin{Satz} \label{FA.4.A4F} $\,$

\noindent \textbf{Vor.:} Seien $\varphi$ ein Quadermaß auf $\R^n$ und $M$ eine nicht-leere i.S.d.\ Maßtheorie $\varphi$-meßbare Teilmenge von $\R^n$.

\noindent \textbf{Beh.:}
\begin{itemize}
\item[(i)] $\D \widetilde{\mathcal{M}}_{\K}(M,\varphi)$ ist eine assoziative $\K$-Unteralgebra von $\K^M$.
\item[(ii)] Für je zwei Funktionen $f, \tilde{f} \in \widetilde{\mathcal{M}}_{\R}(M,\varphi)$ gilt
$$ \sup(f,\tilde{f}), \, \inf(f,\tilde{f}), \, f^+, \, f^-, \, |f| \, \in \widetilde{\mathcal{M}}_{\R}(M,\varphi). $$
\item[(iii)] $\D f \in \widetilde{\mathcal{M}}_{\C}(M,\varphi) \Longrightarrow |f| \in \mathcal{M}_{\R}(M,\varphi)$.
\item[(iv)] $\D f \in \widetilde{\mathcal{M}}_{\K}(M,\varphi)$ \mbox{nicht konstant vom Wert Null auf $M$} \\
$\Longrightarrow$ $\D \overline{f}^1(\K \setminus \{0\}) \in \B(M,\varphi)$ \mbox{und} $\D \frac{1}{f} \in \widetilde{\mathcal{M}}_{\K}(\overline{f}^1(\K \setminus \{0\}),\varphi)$.
\end{itemize}
\end{Satz}

\textit{Beweis.} Der Satz folgt sofort aus dem vorherigen Lemma. \q
\pagebreak

\begin{HS} \label{FA.4.M6} $\,$

\noindent \textbf{Vor.:} Seien $\varphi$ ein Quadermaß auf $\R^n$ und $M$ eine nicht-leere i.S.d.\ Maßtheorie $\varphi$-meßbare Teilmenge von $\R^n$.

\noindent \textbf{Beh.:}
\begin{itemize}
\item[(i)] Es seien $\D (f_i)_{i \in \N}$ eine Folge in $\D \widetilde{\mathcal{M}}_{\K}(M,\varphi)$ und $f \in \K^M$ derart, daß gilt $\D f =_{\varphi} \lim_{i \to \infty} f_i$ auf $M$.\\
Dann folgt $f \in \widetilde{\mathcal{M}}_{\K}(M,\varphi)$.
\item[(ii)] Sei $f \in \widetilde{\mathcal{M}}_{\K}(M,\varphi)$. \\
Dann existiert eine Folge $\D (f_i)_{i \in \N}$ in $\mathcal{L}_{\K}(M,\varphi) \cap \K^M$, die gegen $f$ konvergiert.
Im Falle $\K = \R$ und $f \ge 0$ kann $\D (f_i)_{i \in \N}$ sogar als monoton wachsend gewählt werden.
\end{itemize}
\end{HS}

\textit{Beweis.} Wegen \ref{FA.4.A3K} können wir ohne Beschränkung der Allgemeinheit $\K = \R$ annehmen.

Zu (i): Nach eventuellem Übergang von $f$, $(f_i)_{i \in \N}$ zu $\chi_M \, \hat{f}$, $(\chi_M \, \widehat{f_i})_{i \in \N}$ gilt ohne Einschränkung $M = \R^n$ und
$$ f = \lim_{i \to \infty} f_i = \limsup_{i \to \infty} f_i = \lim_{i \to \infty} \underbrace{\sup \{ f_k \, | \, \k \in \N \, \wedge \, k \ge i \}}_{=: g_i} = \inf \{ g_i(x) \, | \, i \in \N \}, $$
beachte $\forall_{i \in \N} \, g_{i+1} \le g_i$.
Es gilt für jedes $i \in \N$ und jedes $b \in \R$
$$ \forall_{x \in M} \, \left( g_i(x) \le b \Longleftrightarrow \forall_{k \in \N, \, k \ge i} \, f_k(x) \le b \right), $$
also nach Voraussetzung
$$ \overline{g_i}^1( {]} - \infty, b {]} ) = \bigcap_{\substack{k \in \N \\ k \ge i}} \overline{f_k}^1( {]} - \infty, b {]} ) \in \B(M,\varphi), $$
und $g_i$ ist nach \ref{FA.4.M6.L} somit $\B(M,\varphi)$-meßbar.
Daher folgt aus
$$ \forall_{b \in \R} \forall_{x \in M} \, \left( f(x) < b \Longleftrightarrow \exists_{i \in \N} \, g_i(x) < b \right), $$
daß für jedes $b \in \R$ gilt
$$ \overline{f}^1( {]} - \infty, b {[} )= \bigcup_{i \in \N} \overline{g_i}^1( {]} - \infty, b {[} ) \in \B(M,\varphi). $$
Erneut nach \ref{FA.4.M6.L} ist $f$ damit eine $\B(M,\varphi)$-meßbare Funktion.

Zu (ii): Wegen $f = f^+ - f^-$ können wir $f \ge 0$ annehmen.
Für alle $i,k \in \N$ sei
$$ M_{i,k} := \overline{f}^1 \left( \left[ \frac{k}{2^i}, \frac{k+1}{2^i}  \right[ \right) \cap {]} -i, i {[}^n \in \B(M,\varphi). $$
Wegen $\varphi(M_{i,k}) \le \varphi({]} -i, i {[}) < \infty$ gilt $M_{i,k} \in \mathfrak{I}(\R^n,\varphi)$, also $\chi_{M_{i,k}} \in \L(\R^n,\varphi)$.
Daher wird durch
$$ \forall_{i \in \N} \, f_i := \sum_{k=0}^{i \,2^i} \frac{k}{2^i} \, \chi_{M_{i,k}} $$
eine Folge $(f_i)_{i \in \N}$ in $\L(\R^n,\varphi)$ definiert, die die in der Konklusion von (ii) genannten Eigenschaften hat. \q

\begin{Kor} \label{FA.4.M6.K1} 
Seien $\varphi$ ein Quadermaß auf $\R^n$ und $M$ eine nicht-leere i.S.d.\ Maßtheorie $\varphi$-meßbare Teilmenge von $\R^n$.

Dann gilt:
\begin{itemize}
\item[(i)] $\D \mathcal{L}_{\K}(M,\varphi) \cap \K^M \subset \widetilde{\mathcal{M}}_{\K}(M,\varphi)$,
\item[(ii)] $\D \widetilde{\mathcal{M}}_{\K}(M,\varphi) = \mathcal{M}_{\K}(M,\varphi) \cap \K^M$.
\end{itemize}
\end{Kor} 

\textit{Beweis.} Nach \ref{FA.4.A3K} können wir $\K = \R$ annehmen.

Wegen 
$$ \forall_{T \in \mathcal{S}(\R^n)} \left( \left( \forall_{\alpha \in \R \setminus \{0\}} \, \overline{T}^1(\{\alpha\}) \in \mathfrak{P}_n \right) \wedge \left( \exists_{P \in \mathfrak{P}_n} \, \overline{T}^1(\{0\}) = \R^n \setminus P \in \B(\R^n,\varphi) \right) \right) $$ gilt $\mathcal{S}(\R^n) \subset \widetilde{\mathcal{M}}_{\R}(\R^n,\varphi)$, also folgen (i) bzw.\ (ii) ,,$\supset$``  aus \ref{FA.4.M6} (i), da zu jedem $f \in \mathcal{L}_{\R}(M,\varphi) \cap \R^M$ bzw.\ $f \in \mathcal{M}_{\R}(M,\varphi) \cap \R^M$ offenbar jeweils eine Folge $(T_i)_{i \in \N}$ in $\mathcal{S}(\R^n)$ mit $\hat{f} = \chi_M \, \hat{f} =_{\varphi} \lim_{i \to \infty} T_i$ auf $\R^n$, d.h.\ $f =_{\varphi} \lim_{i \to \infty} T_i$ auf $M$, existiert.

(ii) ,,$\subset$`` ergibt sich aus \ref{FA.4.M6} (ii), $\L_{\R}(M,\varphi) \subset \mathcal{M}_{\R}(M,\varphi)$ und \ref{FA.4.69.viii}. \q

\begin{Kor} \label{FA.4.M6.K2} 
Sei $\varphi$ ein Quadermaß auf $\R^n$.

Dann gilt $\B(\R^n,\varphi) = \mathfrak{M}(\R^n,\varphi)$.
\end{Kor}

\textit{Beweis.} ,,$\subset$`` ist bereits nach \ref{FA.4.72} klar.
Sei daher $M \in \mathfrak{M}(\R^n,\varphi)$, d.h.\ genau $\chi_M \in \mathcal{M}_{\R}(\R^n,\varphi) \cap \R^{\R^n} = \widetilde{\mathcal{M}}_{\R}(\R^n,\varphi)$.
Dann gilt $M = \overline{\chi_M}^1(\{1\}) \in \B(\R^n,\varphi)$. \q

\begin{Satz}[Charakterisierung Lebesgue-meßbarer Mengen] \label{FA.4.M4S} $\,$
Es seien $\varphi$ ein Quadermaß auf $\R^n$ und $M$ eine Teilmenge von $\R^n$.

Dann ist $M$ genau dann $\varphi$-meßbar (i.S.d.\ Maßtheorie), wenn zu jedem $\varepsilon \in \R_+$ eine abgeschlossene Teilmenge $A$ des $\R^n$ sowie eine offene Teilmenge $U$ des $\R^n$ mit $A \subset M \subset U$ und $\varphi(U \setminus A) < \varepsilon$ existieren.
\end{Satz}

Zum Nachweis des Satzes beweisen wir zunächst ein weiteres Lemma.

\begin{Lemma} \label{FA.4.M4L} 
Sind $\varphi$ ein Quadermaß auf $\R^n$ und $M$ eine (i.S.d.\ Maßtheorie) $\varphi$-meßbare Teilmenge von $\R^n$, so gilt
\begin{eqnarray*}
\varphi(M) & = & \inf \{ \varphi(U) \, | \, M \subset U \mbox{ und $U$ offen in $\R^n$} \} \\
& = & \sup \{ \varphi(A) \, | \, A \subset M \mbox{ und $A$ abgeschlossen in $\R^n$} \}.
\end{eqnarray*}
\end{Lemma}

\textit{Beweis.} (I) 1.\ Fall: $M$ $\varphi$-integrierbar.
Zu $\varepsilon \in \R_+$ existieren gemäß \ref{FA.4.26} dann $g,h \in \mathcal{S}_{\nearrow}(\R^n,\varphi)$ mit $h \ge_{\varphi} 0$ auf $\R^n$ und
\begin{gather}
\chi_M =_{\varphi} g-h \mbox{ auf $\R^n$}, \label{FA.4.M4L.1} \\
\int_{\R^n} h \, \d \varphi < \frac{\varepsilon}{3}, \label{FA.4.M4L.2}
\end{gather}
insbes.\ gilt $0 \le_{\varphi} h \le_{\varphi} g$ auf $\R^n$.
Der Leser überlege sich als einfache Übung, daß wir ohne Einschränkung annehmen können, daß $g, h$ auf $\R^n$ definiert sind und $0 \le h \le g$ auf $\R^n$ gilt.
Es gibt nun eine $\varphi$-Nullmenge $N$ derart, daß für jedes $x \in \R^n \setminus N$ gilt
\begin{equation*}
x \in M \Longleftrightarrow \chi_M(x) = 1 \Longleftrightarrow \underbrace{h(x)}_{\ge 0} = g(x) - 1,
\end{equation*}
und wir definieren eine Obermenge $\widetilde{M}$ von $M$ durch
\begin{equation} \label{FA.4.M4L.3}
\widetilde{M} := \overline{g}^1 \left( \left] \frac{1}{2}, \infty \right[ \right) \cup N \in \B(\R^n,\varphi),
\end{equation}
beachte, daß aus \ref{FA.4.M6.K1} (i) folgt: $g \in \mathcal{S}_{\nearrow}(\R^n,\varphi) \cap \R^{\R^n} \subset \widetilde{\mathcal{M}}_{\R}(\R^n,\varphi)$.
Es gilt
\begin{equation} \label{FA.4.M4L.4}
\varphi \left( \widetilde{M} \setminus M \right) \le \frac{2 \varepsilon}{3}.
\end{equation}

{[} Zu (\ref{FA.4.M4L.4}): 1.1.\ Fall: $\widetilde{M}$ $\varphi$-integrierbar.
Dann folgt
\begin{eqnarray*}
\varphi \left( \widetilde{M} \setminus M \right) & = & \int_{\R^n} \chi_{\widetilde{M} \setminus M} \d \varphi \stackrel{2 \, g|_{\widetilde{M}} \, >_{\varphi} 1}{\le} 2 \int_{\R^n} \chi_{\widetilde{M} \setminus M} \,g \, \d \varphi \\
& \stackrel{M \subset \widetilde{M}}{=} & 2 \left( \int_{\R^n} \chi_{\widetilde{M}} \, g \, \d \varphi \, - \int_{\R^n} \chi_M \, g \, \d \varphi \right) \\
& \stackrel{g|_M \ge_{\varphi} 1}{\le} & 2 \left( \int_{\R^n} \chi_{\widetilde{M}} \, g \, \d \varphi \, - \int_{\R^n} \chi_{M} \, \d \varphi \right) \\
& \stackrel{(\ref{FA.4.M4L.1})}{=} & 2 \left( \int_{\R^n} \chi_{\widetilde{M}} \, g \, \d \varphi \, - \int_{\R^n} g \, \d \varphi \, + \int_{\R^n} h \, \d \varphi \right) \\
& = & 2 \left( \int_{\R^n} \left( \chi_{\widetilde{M}} - 1_{\R^n} \right) \, g \, \d \varphi \, + \int_{\R^n} h \, \d \varphi \right) \stackrel{g \ge 0}{\le} 2 \int_{\R^n} h \, \d \varphi \\
& \stackrel{(\ref{FA.4.M4L.2})}{<} & \frac{2 \varepsilon}{3}
\end{eqnarray*}

1.2.\ Fall: $\widetilde{M}$ nicht $\varphi$-integrierbar.
Dann gilt $\widetilde{M} \setminus M \in \B(\R^n,\varphi)$ sowie 
$$ \widetilde{M}_i := \widetilde{M} \cap {]}-i, i{[}^n, \, M_i := M \cap {]}-i, i{[}^n \in \mathfrak{I}(\R^n,\varphi) $$
für jedes $i \in \N$, und durch Anwendung des Beweises des ersten Falles auf $\widetilde{M}_i, M_i$ anstelle von $\widetilde{M}, M$ erhält man
$$ \varphi \left( (\widetilde{M} \setminus M) \cap {]} -i, i {[}^n \right) = \varphi \left( \widetilde{M}_i \setminus M_i \right) < \frac{2 \varepsilon}{3}, $$
also $\varphi \left( \widetilde{M} \setminus M \right) = \lim_{i \to \infty} \varphi \left( (\widetilde{M} \setminus M) \cap {]} -i, i {[}^n \right) \le \frac{2 \varepsilon}{3}$. {]}

Wegen $g \in \mathcal{S}_{\nearrow}(\R^n,\varphi)$ existiert eine monoton wachsende Folge $(T_i)_{i \in \N}$ in $\mathcal{S}(\R^n)$ mit
$$ g =_{\varphi} \lim_{i \to \infty} T_i. $$
Dann folgt für jedes $i \in \N$
$$ P_i := \overline{T_i}^1 \left( \left] \frac{1}{2}, \infty \right[ \right) \in \mathfrak{P_n}, $$
und es existiert eine $\varphi$-Nullmenge $\widetilde{N}$, die $N$ enthält und deren Offenheit in $\R^n$ wir wegen \ref{FA.4.13.S} (iii) ohne Einschränkung annehmen können, mit
\begin{equation} \label{FA.4.M4L.5}
\widetilde{M} \cup \widetilde{N} \stackrel{(\ref{FA.4.M4L.3})}{=} \overline{g}^1 \left( \left] \frac{1}{2}, \infty \right[ \right) \cup \widetilde{N} = \left( \bigcup_{i \in \N} P_i \right) \cup \widetilde{N}.
\end{equation}

Wegen der Regularität von $\varphi$ gibt es zu jedem $i \in \N$ eine offene parkettierbare Obermenge $P'_i$ von $P_i$ derart, daß gilt
\begin{equation} \label{FA.4.M4L.6}
\forall_{i \in \N} \, \varphi \left( P'_i \setminus P_i \right) < \frac{\varepsilon}{3 \cdot 2^{i+3}} ~~ \mbox{ und somit } ~~ \sum_{i=0}^{\infty} \varphi \left( P'_i \setminus P_i \right) \le \frac{\varepsilon}{6}. 
\end{equation}
Für die nach (\ref{FA.4.M4L.5}) offene Obermenge $U := \left( \bigcup_{i \in \N} P'_i \right) \cup \widetilde{N}$ von $\widetilde{M}$ ergibt sich nun mittels \ref{FA.4.M4F} (ii), (iii)
\begin{eqnarray*}
\varphi(U) - \varphi \left( \widetilde{M} \right) & \stackrel{(\ref{FA.4.M4L.5})}{=} & \varphi \left( \left( \bigcup_{i \in \N} P'_i \right) \setminus \left( \bigcup_{i \in \N} P_i \right) \right) = \varphi \left( \bigcap_{j \in \N} \left( \left( \bigcup_{i \in \N}  P'_i  \right) \setminus P_j \right) \right) \\
& \le & \varphi \left( \bigcup_{i \in \N} \left( P'_i  \setminus P_i \right) \right) \le \sum_{i=0}^{\infty} \varphi \left( P'_i  \setminus P_i \right) \stackrel{(\ref{FA.4.M4L.6})}{\le} \frac{\varepsilon}{6} < \frac{\varepsilon}{3},
\end{eqnarray*}
also nach (\ref{FA.4.M4L.4}) auch
$$ \varphi(U) - \varphi(M) = \underbrace{\varphi(U) - \varphi \left( \widetilde{M} \right)}_{< \frac{\varepsilon}{3}} + \underbrace{\varphi \left( \widetilde{M} \right) - \varphi(M)}_{\le \frac{2 \varepsilon}{3}} < \varepsilon. $$

2.\ Fall: $M$ $\varphi$-meßbar i.S.d.\ Maßtheorie.
Sei $\varepsilon \in \R_+$.
Dann folgt zu jedem $i \in \N$ aus \ref{FA.4.M4} und dem 1.\ Fall die Existenz einer offenen Obermenge $U_i$ von $M \cap {]} -i, i {[}^n$ mit $\varphi ( U_i \setminus ( M \cap {]} -i, i {[}^n ) ) < \frac{\varepsilon}{2^{i+2}}$.
Somit ist $U := \bigcup_{i \in \N} U_i$ eine offene Obermenge von $M$ derart, daß nach \ref{FA.4.M4F} (iii), (ii) gilt
$$ \varphi(U \setminus M) = \varphi \left( \bigcup_{i \in \N} ( U_i \setminus M ) \right) \le \sum_{i=0}^{\infty} \varphi ( U_i \setminus M ) \le \sum_{i=0}^{\infty} \varphi ( U_i \setminus ( M \cap {]} -i, i {[}^n ) ) < \varepsilon. $$

(II) Nach (I) existiert zu $\varepsilon \in \R_+$ eine offene Obermenge $U$ von $\R^n \setminus M$ mit $\varphi(U \setminus (\R^n \setminus M)) < \varepsilon$.
Dann ist $A := \R^n \setminus U$ eine abgeschlossene Teilmenge von $\R^n$\linebreak mit $A \subset M$ und $\varphi(M \setminus A) = \varphi(M \setminus (\R^n \setminus U)) = \varphi(M \cap U) = \varphi(U \setminus (\R^n \setminus M)) < \varepsilon$.

Aus (I) und (II) ergibt sich das Lemma.
\q
\A
\textit{Beweis des Satzes.} ,,$\Rightarrow$`` folgt sofort aus obigem Lemma.

,,$\Leftarrow$`` Wähle zu jedem $i \in \N$ eine offene Obermenge $U_i$ von $M$ und eine abgeschlossene Teilmenge $A_i$ von $M$ derart, daß gilt $\varphi(U_i \setminus A_i) < \frac{1}{i+1}$.
Dann folgt
\begin{equation} \label{FA.4.M4S.1}
B := \bigcap_{i \in \N} U_i, \, C := \bigcup_{i \in \N} A_i \in \B^n \subset \B(\R^n,\varphi) ~~ \mbox{ und } ~~ C \subset M \subset B
\end{equation}
sowie
\begin{eqnarray*}
\varphi(B \setminus C) & = & \varphi \left( \bigcap_{i \in \N} (B \setminus A_i) \right) = \varphi \left( \bigcap_{i \in \N} \left( \left( \bigcap_{j \in \N} U_j \right) \setminus A_i \right) \right) \\
& \le & \varphi \left( \bigcap_{i \in \N} \left( U_i \setminus A_i \right) \right) \le \inf\{ \varphi(U_i \setminus A_i) \, | \, i \in \N \} = 0,
\end{eqnarray*}
also ist $M \setminus C$ als Teilmenge von $B \setminus C$ eine $\varphi$-Nullmenge und somit nach \ref{FA.4.30.K2} $\varphi$-integrierbar.
Daher ergibt sich aus (\ref{FA.4.M4S.1}): $M = C \cup (M \setminus C) \in \B(\R^n,\varphi)$. \q

\subsection*{Übungsaufgaben}

\begin{UA}
Zeige, daß die in \ref{FA.4.4.B} 5.) definierten Stieltjeschen Maße tatsächlich Quadermaße sind.
\end{UA}

\begin{UA}
Sei $N$ eine Teilmenge von $\R^n$.

Beweise, daß $N$ genau dann eine $\mu_n$-Nullmenge ist, wenn es zu jedem $\varepsilon \in \R_+$\linebreak eine Folge $(Q_i)_{i \in \N}$ abgeschlossener Quader des $\R^n$ mit $N \subset \bigcup_{i=0}^{\infty} Q_i$ und\linebreak $\sum_{i=0}^{\infty} \mu_n(Q_i) < \varepsilon$ gibt.
\end{UA}

\begin{UA}
Sei $U$ eine offene Teilmenge von $\R^n$.

Zeige, daß eine Folge $(Q_i)_{i \in \N}$ abgeschlossener Quader des $\R^n$ mit $U = \bigcup_{i=1}^{\infty} Q_i$ und $\forall_{i,j \in \N, \, i \ne j} \, \mu_n (Q_i \cap Q_j) = 0$ existiert.
\end{UA}

\begin{UA}
$M$ bezeichne die Smith-Volterra-Cantor Menge, vgl.\ \ref{FA.4.56}.

Beweise $\chi_M \in \L_{\R}(M,\mu_1) \setminus \mathcal{S}_{\nearrow}(\R^n,\mu_1)$.
\end{UA}

\begin{UA}
Beweise das Beispiel zu \ref{FA.4.66}.
\end{UA}

\begin{UA}
Führe den Beweis des Satzes \ref{FA.4.NmT} in allen Einzelheiten aus.
\end{UA}

\begin{UA} \label{FA.4.exp}
Zeige $\e^{-(x_1^2 + \ldots + x_n^2)} \in \L_{\R}(\R^n,\mu_n)$ und 
$$ \int_{\R^n} \e^{-(x_1^2 + \ldots + x_n^2)} \, \d \mu_n = \pi^{\frac{n}{2}}. $$
\end{UA}

Tip: Betrachte zunächst den Fall $n=2$.

\begin{UA}
Zeige, daß $\frac{1}{x + y} \: \R^2 \setminus \{ (\alpha , \beta) \in \R^2 \, | \, \alpha + \beta = 0 \} \to \R$ $\mu_2$-in\-te\-grier\-bar über ${[}0,1{]}^2$ ist, und berechne das $\mu_2$-Integral.
(Zur Kontrolle sei erwähnt, daß letzteres gleich $2 \, \ln(2)$ ist.)
\end{UA}

\begin{UA}
Sei $\rho \in \R_+$.
Es bezeichnen $Z_1 \subset \R^3$ den Vollzylinder vom Radius $\rho$ um die $z$-Achse und $Z_2 \subset \R^3$ den Vollzylinder vom Radius $\rho$ um die $x$-Achse.

Zeige, daß $K := Z_1 \cap Z_2$ kompakt ist, und berechne das $\mu_3$-Volumen $\mu_3(K)$.
(Zur Kontrolle erwähnen wir bereits $\mu_3(K) = \frac{16}{3} \, \rho^3$.)

Der Rand von $\left({[}-\rho, \rho{]}^2 \times {[}0, \rho{]}\right) \cap K$ ist übrigens ein \emph{Kreuzgewölbe}.
\end{UA}

\begin{UA}
Berechne das $\mu_3$-Volumen der kompakten Vollteilmenge $K$ des $\R^3$, die von der $(x,y)$-Ebene, dem Zylinder $\{ (x,y,z) \in \R^3 \, | x^2 + y^2 = 2 \, x \}$ und dem Kegel $\{ (x,y,z) \in \R^3 \, | \, z = \sqrt{x^2 + y^2}\}$ berandet wird.
(Zur Kontrolle: Die Lösung lautet $\frac{32}{9}$.)
Skizziere $K$!
\end{UA}

\begin{UA}
Es sei $\rho \in \R_+$.
$B_{\rho}(0)$ bezeichne die Vollkugel vom Radius $\rho$ in $(\R^n, \| \ldots \|_2)$.

Zeige
$$ \mu_n(B_{\rho}(0)) = 
\left\{ 
\begin{array}{cl} \pi \cdot \frac{\pi}{2} \cdot \frac{\pi}{3} \cdot \cdots \cdot \frac{\pi}{\frac{n}{2}} \cdot \rho^n, & \mbox{ falls } n \equiv 0 (2), \\ 
                  2 \cdot \frac{2 \, \pi}{3} \cdot \frac{2 \, \pi}{5} \cdot \cdots \cdot \frac{2 \, \pi}{n} \cdot \rho^n, & \mbox{ falls } n \equiv 1 (2). 
\end{array}
\right. $$
\end{UA}

\begin{UA}[Die allgemeine Version der Sätze von \textsc{Fubini} und \textsc{Tonelli}]\index{Satz!von \textsc{Fubini}}\index{Satz!von \textsc{Tonelli}}
Es seien $n_1, n_2 \in \N_+$ mit $n = n_1 + n_2$, $M$ eine Teilmenge von $\R^n$ und $f \in \K^M$ sowie $\varphi_i$ Quadermaße auf $\R^{n_i}$ für $i \in \{1,2\}$.
Wir definieren
\begin{gather*}
\pi_1(M) := \{ x_1 \in \R^{n_1} \, | \, \exists_{x_2 \in \R^{n_2}} \, (x_1,x_2) \in M \} \subset \R^{n_1}
\end{gather*}
und für jedes $x_1 \in \pi_1(M)$
\begin{gather*}
M_{x_1} := \{ x_2 \in \R^{n_2} \, | \, (x_1,x_2) \in M \} \subset \R^{n_2}, \\
f_{x_1} := f(x_1, \ldots) \in \K^{M_{x_1}}.
\end{gather*}
Des weiteren sei $N_1 := \{ x_1 \in \pi_1(M) \, | \, f_{x_1} \notin \L(M_{x_1},\varphi_2)\}$.

Zeige:
\begin{itemize}
\item[(i)] Ist $f \in \L_{\K}(M,\varphi_1 \times \varphi_2)$, so ist $N_1$ eine $\varphi_1$-Nullmenge, und die Funktion
$$ \pi_1(M) \setminus N_1 \longrightarrow \K, ~~ x_1 \longmapsto \int_{M_{x_1}} f_{x_1} \, \d \varphi_2, $$
ist $\varphi_1$-integrierbar über $\pi_1(M)$.
Des weiteren gilt
\begin{eqnarray*}
\int_M f \, \d (\varphi_1 \times \varphi_2) & = &\int_{\pi_1(M)} \left( \int_{M_{\ldots}} f_{\ldots} \, \d \varphi_2 \right) \d \varphi_1 \\
& =: & \int_{\pi_1(M)} \left( \int_{M_{x_1}} f_{x_1}(x_2) \, \d \varphi_2(x_2) \right) \, \d \varphi_1(x_1).
\end{eqnarray*}
\item[(ii)] Sind $f \in \mathcal{M}_{\K}(M,\varphi_1 \times \varphi_2)$, $N_1$ eine $\varphi_1$-Nullmenge und
$$ \pi_1(M) \setminus N_1 \longrightarrow \R, ~~ x_1 \longmapsto \int_{M_{x_1}} |f_{x_1}| \, \d \varphi_2, $$
$\varphi_1$-integrierbar über $\pi_1(M)$, so gilt $f \in \L_{\K}(M,\varphi_1 \times \varphi_2)$.
\end{itemize}
\end{UA}

\begin{UA}
Beweise Satz \ref{FA.4.M4F}.
\end{UA}

\cleardoublepage
\section{Lebesguesche Räume} \label{FAna5}

Sei im folgenden für jedes $p \in \R_+$ stets 
$$ \boxed{0^p} := 0, $$
beachte, daß die differenzierbare Funktion
$$ \R_+ \longrightarrow \R, ~~ t \longmapsto t^p, $$
in Null durch Null stetig fortgesetzt wird.
Wie oben erwähnt, übernehmen wir des weiteren in diesem Kapitel noch die Voraussetzungen des vorherigen.

\subsection*{Einführung} \addcontentsline{toc}{subsection}{Einführung}

\begin{Def} \label{FA.5.1}
Seien $\varphi$ ein Quadermaß auf $\R^n$ und $M$ eine Teilmenge von $\R^n$.

Für $p \in \R_+$ setzen wir
$$ \boxed{\L_{\K}^p(M,\varphi)} := \{ f \in \mathcal{M}_{\K}(M,\varphi) \, | \, |f|^p \in \L_{\R}(M,\varphi) \} $$
und
$$ \forall_{f \in \L_{\K}^p(M,\varphi)} \, \boxed{\|f\|_p} := \left( \int_M |f|^p \, \d \varphi \right)^{\frac{1}{p}}. $$

Wir definieren des weiteren
$$ \boxed{\L_{\K}^{\infty}(M,\varphi)} := \{ f \in \mathcal{M}_{\K}(M,\varphi) \, | \, \exists_{C \in \R_+} \, |f| \le_{\varphi} C \} $$
sowie
$$ \forall_{f \in \L_{\K}^{\infty}(M,\varphi)} \, \boxed{\|f\|_{\infty}} := \inf\{ C \in \R_+ \, | \, |f| \le_{\varphi} C \}. $$

Ist $p \in {]}0, \infty{]}$, so heißen die Elemente von $\L_{\K}^p(M,\varphi)$ \emph{(reellwertig, falls $\K=\R$) $p$-integrierbar über $M$ bzgl.\ $\varphi$} und im Falle $p=2$ auch \emph{(reellwertig, falls $\K=\R$) $\varphi$-quadrat\-in\-te\-grier\-bar über $M$ bzgl.\ $\varphi$}.\index{Funktion!integrierbare!$p$-}
Es gilt offenbar für jedes $f \in \mathcal{M}_{\K}(M,\varphi)$
\begin{gather*}
f \in \L_{\K}^p(M,\varphi) \Longleftrightarrow |f| \in \L_{\R}^p(M,\varphi), \\
f \in \L_{\K}^p(M,\varphi) \Longleftrightarrow {\rm Re} \, f, \, {\rm Im} \, f \in \L_{\R}^p(M,\varphi) \Longleftrightarrow \overline{f} \in \L_{\K}^p(M,\varphi).
\end{gather*}

\begin{Bem*} $\,$
\begin{itemize}
\item[1.)] Aus \ref{FA.4.24} (iii) und \ref{FA.4.70} (ii) folgt $\L_{\K}^1(M,\varphi) = \L_{\K}(M,\varphi)$.
\item[2.)] Für alle $f \in \L_{\K}^{\infty}(M,\varphi)$ gilt $|f| \le_{\varphi} \|f\|_{\infty}$, denn 
$$ \{ x \in M \, | \, |f(x)| > \|f\|_{\infty} \} = \bigcup_{k \in \N} \left\{ x \in M \left| \, |f(x)| > \|f\|_{\infty} + \frac{1}{k+1} \right. \right\} $$
ist eine $\varphi$-Nullmenge.
\end{itemize}
\end{Bem*}
\end{Def}

\begin{Satz} \label{FA.5.2}
Seien $\varphi$ ein Quadermaß auf $\R^n$ und $M$ eine Teilmenge von $\R^n$.

Dann gilt:
\begin{itemize}
\item[(i)] Für jedes $p \in {]}0, \infty{]}$ ist $\L_{\K}^p(M,\varphi)$ ein $\K$-Untervektorraum von $\mathcal{M}_{\K}(M,\varphi)$.
\item[(ii)] Für jedes $p \in {[}1, \infty{]}$ ist $\| \ldots \|_p$ eine Halbnorm auf $\L_{\K}^p(M,\varphi)$.
\end{itemize}
\end{Satz}
\pagebreak

\textit{Beweis.} Zu (i): Daß $\mathcal{M}_{\K}(M,\varphi)$ ein $\K$-Vektorraum ist, haben wir in \ref{FA.4.69} (i) bereits eingesehen.
Daher bleibt zu zeigen
\begin{equation} \label{FA.5.2.1}
\forall_{\lambda \in \K} \, \forall_{f,g \in \mathcal{L}_{\K}^p(M,\varphi)} \, \lambda \, f , f+g \in \L_{\K}^p(M,\varphi).
\end{equation}

{[} Zu (\ref{FA.5.2.1}): 1.\ Fall: $p \in \R_+$.
Seien $\lambda \in \K$ und $f,g \in \mathcal{M}_{\K}(M,\varphi)$ zwei Funktionen mit $|f|^p, |g|^p \in \L_{\R}(M,\varphi)$.
Dann folgt aus \ref{FA.4.22} (ii) sofort $| \lambda \, f |^p \in \L_{\R}(M,\varphi)$.
Des weiteren gilt
\begin{equation} \label{FA.5.2.2}
| f+g |^p \le \left( |f| + |g| \right)^p \le \left( 2 \, \sup( |f|, |g| ) \right)^p = 2^p \, \sup \left( |f|^p, |g|^p \right).
\end{equation}
Nun folgt aus $f,g \in \mathcal{M}_{\K}(M,\varphi)$ auch $f+g \in \mathcal{M}_{\K}(M,\varphi)$, also wegen der Stetigkeit von $| \ldots |^p \: \R \to \R$ mit $|0|^p = 0$ sowie \ref{FA.4.8} (iii) offenbar die $\varphi$-Meßbarkeit von $|f+g|^p = | \ldots |^p \circ (f+g)$ über $M$.
Außerdem gilt mit $|f|^p, |g|^p \in \L_{\R}(M,\varphi)$ auch $2^p \, \sup ( |f|^p, |g|^p ) \in \L_{\R}(M,\varphi)$, vgl.\ \ref{FA.4.22} (ii), (v).
Somit ergeben (\ref{FA.5.2.2}) und \ref{FA.4.70} (i):\linebreak $|f+g|^p \in \L_{\R}(M,\varphi)$.

2.\ Fall: $p = \infty$.
Seien $\lambda \in \K$ und $f,g \in \mathcal{M}_{\K}(M,\varphi)$ \emph{$\varphi$-beschränkt}, d.h.\ per definitionem $\varphi$-fast überall beschränkt.
Trivialerweise sind dann auch $\lambda \, f$ sowie $f+g$ $\varphi$-beschränkt. {]}

Zu (ii): Die Homogenität von $\| \ldots \|_p \: \R \to {[}0, \infty{[}$ ist trivial, und die Minkowskische Ungleichung, die wir unten in \ref{FA.5.11} beweisen werden, besagt genau die Dreiecksungleichung. \q

\begin{Bem*} \label{FA.5.3}
Es seien $\varphi$ ein Quadermaß auf $\R^n$, $M$ eine Teilmenge von $\R^n$ und $p \in {]}0,1{[}$.
Wir werden in \ref{FA.5.13} und \ref{FA.5.14} (ii) sehen, daß dann zwar 
\begin{equation}
\forall_{f,g \in \L_{\R}^p(M,\varphi)} \, \left( f,g \ge_{\varphi} 0 \, \wedge \, f+g >_{\varphi} 0 \right) \Longrightarrow \| f + g \|_p \ge \|f\|_p + \|g\|_p \label{FA.5.3.1}
\end{equation}
aber immerhin
\begin{equation}
\forall_{f,g \in \L_{\K}^p(M,\varphi)} \, \| f + g \|_p \le 2^{\frac{1}{p}-1} \left( \|f\|_p + \|g\|_p \right) \label{FA.5.3.2}
\end{equation}
gilt.
(\ref{FA.5.3.1}) läßt vermuten, daß $\| \ldots \|_p$ i.a.\ keine Halbnorm auf $\L_{\K}^p(M,\varphi)$ ist, und dies stimmt:
Sind z.B.\ $\varphi:= \mu_1$ das eindimensionale Volumen, $M := [0,1]$, $p := \frac{1}{2}$ und $f := \chi_{{[}0, \frac{1}{2}{]}}, g := \chi_{{]}\frac{1}{2}, 1{]}} \in \L_{\R}^{\frac{1}{2}}([0,1],\mu_1)$, so ergibt sich
$$ \| f + g \|_{\frac{1}{2}} = 1 > \frac{1}{2} = 2 \left( \frac{1}{2} \right)^2 = \|f\|_{\frac{1}{2}} + \|g\|_{\frac{1}{2}}. $$
\end{Bem*}

\begin{Def} \label{FA.5.4}
Seien $\varphi$ ein Quadermaß auf $\R^n$, $M$ eine Teilmenge von $\R^n$ und $p \in {]}0, \infty{]}$.

Dann ist
$$ \boxed{\mathcal{N}_{\K}(M,\varphi)} := \{ f \in \mathcal{M}_{\K}(M,\varphi) \, | \, f =_{\varphi} 0 \} $$
offenbar ein $\K$-Untervektorraum von $\L_{\K}^p(M,\varphi)$, und wir definieren den $\K$-Vek\-tor\-raum
$$ \boxed{L_{\K}^p(M,\varphi)} := \L_{\K}^p(M,\varphi) / \mathcal{N}_{\K}(M,\varphi), $$
dessen Elemente ebenso wie die von $\L_{\K}^p(M,\varphi)$ die \emph{(reellwertigen, falls $\K=\R$) $p$-integrierbaren Funktionen über $M$ bzgl.\ $\varphi$}\index{Funktion!integrierbare!$p$-} heißen.
\end{Def}

\begin{Lemma} \label{FA.5.5}
Seien $\varphi$ ein Quadermaß auf $\R^n$, $M$ eine Teilmenge von $\R^n$ und $p \in {]}0, \infty{]}$.

Dann gilt $\mathcal{N}_{\K}(M,\varphi) = \{ f \in \L_{\K}^p(M,\varphi) \, | \, \|f\|_p = 0 \}$.
\end{Lemma}

\textit{Beweis.} ,,$\subset$`` und im Falle $p = \infty$ ,,$\supset$`` sind trivial.

Zu ,,$\supset$`` im Falle $p \in \R_+$: Sei $f \in \L_{\K}^p(M,\varphi)$ mit $\|f\|_p = 0$.
Dann folgt zum einen $|f|^p \in \mathcal{M}_{\R}(M,\varphi)$, also wie oben wegen der Stetigkeit von $| \ldots |^{\frac{1}{p}} \: \R \to \R$ mit $|0|^{\frac{1}{p}} = 0$ sowie \ref{FA.4.8} (iii) auch $|f| \in \mathcal{M}_{\R}(M,\varphi)$, d.h.\ $f \in \mathcal{M}_{\K}(M,\varphi)$ nach \ref{FA.4.69} (iii),\linebreak und zum anderen $f =_{\varphi} 0$. \q

\begin{Def}[Lebesguesche Räume] \label{FA.5.6} \index{Raum!Lebesguescher}
Seien $\varphi$ ein Quadermaß auf $\R^n$, $M$ eine Teilmenge von $\R^n$ und $p \in {]}0, \infty{]}$.
$\pi \: \L_{\K}^p(M,\varphi) \to L_{\K}^p(M,\varphi)$ bezeichne den kanonischen Epimorphismus.

\begin{itemize}
\item[(i)] Durch
$$ \forall_{f \in \L_{\K}^p(M,\varphi)} \, \boxed{\| \pi(f) \|_p} := \| f \|_p $$
wird eine Abbildung 
$$ \boxed{\| \ldots \|_p \: L_{\K}^p(M,\varphi) \longrightarrow {[}0, \infty {[}} $$
definiert, die wir also ebenfalls mit $\| \ldots \|_p$ bezeichnen, so daß nach dem letzten Lemma und \ref{FA.5.2} (ii) gilt  \index{Norm!$p$-}
\begin{equation} \label{FA.5.6.S}
p \in {[} 1, \infty {]} \Longrightarrow  \left( L_{\K}^p(M,\varphi), \| \ldots \|_p \right) \mbox{ ist ein normierter $\K$-Vektorraum.}
\end{equation}
\item[(ii)] Für $p \in {[} 1, \infty {]}$ heißen die normierten $\K$-Vektorräume $L_{\K}^p(M,\varphi)$ wie in (\ref{FA.5.6.S}) die \emph{Lebesgueschen Räume auf $M$ bzgl.\ $\varphi$}.
\item[(iii)] Für $p \in {]}0,1{[}$ ist $\| \ldots \|_p$ zwar i.a.\ keine Norm auf $L_{\K}^p(M,\varphi)$ (nur die Dreiecksungleichung ist ggf.\ nicht erfüllt), aber wegen (\ref{FA.5.3.2}) kann man analog zu \ref{FA.2.2} (i), \ref{FA.1.2} eine Topologie für $L_{\K}^p(M,\varphi)$ definieren, welche offenbar durch die Metrik $\boxed{{d_p}^p}$ auf $L_{\K}^p(M,\varphi)$, gegeben durch
$$ \forall_{\mathbf{f}, \mathbf{g} \in L_{\K}^p(M,\varphi)} \, {d_p}^p(\mathbf{f}, \mathbf{g}) := {\| \mathbf{f} - \mathbf{g} \|_p}^p, $$
beachte \ref{FA.5.14} (i) unten, metrisiert wird, die $L_{\K}^p(M,\varphi)$ zu einem topologischen $\K$-Vektorraum macht -- letzteres sieht man analog zu \ref{FA.2.2} (iv) ein.
Wir versehen die $\K$-Vek\-tor\-räume $L_{\K}^p(M,\varphi)$ mit der Metrik ${d_p}^p$ und nennen diese Räume ebenfalls \emph{Lebesguesche Räume auf $M$ bzgl.\ $\varphi$}.
\end{itemize}

\begin{Bem*}
Sei $p \in {]}0,1{[}$.
Dann zeigen die Lebesgueschen Räume auf $M$ bzgl.\ $\varphi$ einige Pathologien.
Z.B.\ ist für $\mu_n$-meßbares $M$ nach \cite[Theorem 1]{Day} jede stetige lineare Abbildung $L_{\R}^p(M,\mu_n) \to \R$ konstant vom Wert Null.
Vergleiche dies mit \ref{FA.5.17} (ii) unten!
Der Trennungssatz von \textsc{Hahn-Banach} \ref{FA.2.24} (ii) zeigt nun übrigens, daß $L_{\R}^p(M,\mu_n)$, aufgefaßt als topologischer Raum, für solche $M$ im Falle $M \ne \emptyset$ nicht normierbar ist.
\end{Bem*}
\end{Def}

\begin{Bsp} \label{FA.5.6.Bsp} $\,$
\begin{itemize}
\item[1.)] Seien $m := 1_{\{1, \ldots, n\}} \: \{1, \ldots, n\} \to \R_+$ und $\varphi_m$ das zugehörige Quadermaß der diskreten Massenverteilung $m$, welches ein Quadermaß auf $\R$ ist, vgl.\ Beispiel \ref{FA.4.4.B} 4.).
Offenbar können wir dann $\mathcal{M}_{\K}(\{1, \ldots, n\},\varphi_m)$ mit $\K^n$ identifizieren, und es gilt des weiteren für dessen $\K$-Un\-ter\-vek\-tor\-raum $\mathcal{N}(\{1, \ldots, n\},\varphi_m) = \{0\}$.
U.a.\ aus \ref{FA.4.45} folgt daher, daß die in Beispiel 1.) zu \ref{FA.2.1} mit $\| \ldots \|_p$ für $p \in \{1, 2, \infty\}$ bezeichneten Normen auf $\K^n$ mit den ebenso bezeichneten auf $L^p_{\K}(\{1, \ldots, n\},\varphi_m) \equiv \K^n$ übereinstimmen.
\item[2.)] Sind $m := 1_{\N} \: \N \to \R_+$ und $\varphi_m$ das zugehörige Quadermaß der diskreten Massenverteilung $m$, so können wir $\mathcal{M}_{\K}(\N,\varphi_m)$ analog zu 1.) mit $\K^{\N}$ identifizieren, und es gilt auch jetzt $\mathcal{N}(\N,\varphi_m) = \{0\}$.
Somit sind offenbar folgende Identifikationen möglich
$$ \ell^1_{\K} \equiv L_{\K}^1(\N,\varphi_m) ~~ \mbox{ und } ~~ \ell^{\infty}_{\K} \equiv L_{\K}^{\infty}(\N,\varphi_m), $$
vgl.\ Beispiel \ref{FA.3.B}.
\end{itemize}
\end{Bsp}

\begin{Def}[Lebesguesche Folgenräume] \label{FA.5.7} \index{Folgen!-räume!Lebesguesche}
Es seien $m := 1_{\N} \: \N \to \R_+$, $\varphi_m$ das zugehörige Quadermaß der diskreten Massenverteilung $m$ und $p \in {]}0, \infty{]}$.

Wir definieren dann
$$ \boxed{\ell^p_{\K}} := L_{\K}^p(\N,\varphi_m), $$
d.h.\ $\ell^p_{\K}$ ,,ist`` im Falle $p < \infty$ der $\K$-Vektorraum aller Folgen $(x_i)_{i \in \N}$ in $\K$, für welche die Reihe $\sum_{i=0}^{\infty} {x_i}^p$ absolut konvergiert, und im Falle $p = \infty$ der $\K$-Vektorraum aller beschränkten Folgen in $\K$ zusammen mit der Metrik
$$ \forall_{(x_i)_{i \in \N}, (y_i)_{i \in \N} \in \K^{\N}} \, {d_p}^p \left( (x_i)_{i \in \N}, (y_i)_{i \in \N} \right) = \sum_{i=0}^{\infty} \left| x_i - y_i \right|^p , ~~ \mbox{falls $p \in {]}0,1{[}$}, $$
bzw.\ der Norm
$$ \forall_{(x_i)_{i \in \N} \in \K^{\N}} \, \| (x_i)_{i \in \N} \|_p = \left\{ \begin{array}{cl} \D \left( \sum_{i=0}^{\infty} |x_i|^p \right)^{\frac{1}{p}}, & \mbox{falls $p \in {[}1, \infty{[}$}, \\ \D \sup \left\{ |x_i| \, | \, i \in \N \right\}, & \mbox{falls $p = \infty$}. \end{array} \right. $$
Diese Räume heißen die \emph{Lebesgueschen Folgenräume}.

\begin{Bem*} \label{FA.5.7.B} $\,$
\begin{itemize}
\item[1.)] $\| \ldots \|_p$ ist im Falle $p \in {]}0,1{[}$ keine Norm auf $\ell^p_{\K}$, denn es gilt
$$ \| (1,0,0,\ldots) + (0,1,0,\ldots) \|_p = 2^{\frac{1}{p}} > 2 = \| (1,0,0,\ldots) \|_p + \| (0,1,0,\ldots) \|_p . $$
\item[2.)] Aus der Jensenschen Ungleichung, vgl.\ \ref{FA.5.10} unten, folgt offenbar für alle $p, \tilde{p} \in \R_+$ mit $p \le \tilde{p}$: $\ell_{\K}^p \subset \ell_{\K}^{\tilde{p}} \subset \ell_{\K}^{\infty}$ und $\| \ldots \|_{\infty} \le \| \ldots \|_{\tilde{p}} \le \| \ldots \|_p$ auf $\ell_{\K}^p$.
\item[3.)] In der aktuellen mathematischen Forschung spielen die Lebesgueschen Folgenräume für $p \in {]}0,1{]}$ in der Theorie der \emph{Verdichteten Mathematik} von \textsc{Clausen} und \textsc{Scholze} eine Rolle, vgl.\ hierzu \cite[Seite 17]{Scholze}.
Das Ziel dieser Theorie besteht u.a.\ darin, die Funktionalanalysis zu algebraisieren, d.h.\ zu einem Teil der \emph{Kommutativen Algebra} zu machen.
\end{itemize}
\end{Bem*}
\end{Def}

\subsection*{Ungleichungen} \addcontentsline{toc}{subsection}{Ungleichungen}

Ohne Beweis haben wir im letzten Abschnitt mehrere Ungleichungen verwendet, u.a.\ deren Beweise sollen nun erbracht werden.
Wir weisen zunächst auf das folgende Lemma hin.

\begin{Lemma} \label{FA.5.pq}
Für alle $p,q \in \R \setminus \{0\}$ gilt
\begin{eqnarray*}
\lefteqn{\frac{1}{p} + \frac{1}{q} = 1} \\
& & \Longleftrightarrow p + q = p \, q \Longleftrightarrow (p-1) \, (q-1) = 1 \Longleftrightarrow p \, (q-1) = q \Longleftrightarrow  q \, (p-1) = p.
\end{eqnarray*}
\q
\end{Lemma}

\begin{Satz}[Höldersche Ungleichung] \index{Ungleichung!Höldersche} \label{FA.5.8}
Es seien $\varphi$ ein Quadermaß auf $\R^n$, $M$ eine Teilmenge von $\R^n$, $p,q \in [1,\infty]$ mit $\frac{1}{p} + \frac{1}{q} = 1$, wobei $\frac{1}{\infty}$ als $0$ zu lesen ist, und $f \in \L^p_{\K}(M,\varphi)$ sowie $g \in \L^q_{\K}(M,\varphi)$.

Dann gilt $f \, g \in \L^1_{\K}(M,\varphi) = \L_{\K}(M,\varphi)$ und
$$ \| f \, g \|_1 \le \|f\|_p \, \|g\|_q. $$
\end{Satz}

\textit{Beweis.} Wir können ohne Beschränkung der Allgemeinheit annehmen, daß gilt 
$$ \|f\|_p \, \|g\|_q \ne 0, $$
denn sonst ist die Behauptung trivial.

Seien zunächst $p,q \in {]} 1, \infty {[}$.
Zeige als Übung die \emph{Youngsche Ungleichung}\index{Ungleichung!Youngsche}
\begin{equation} \label{FA.5.8.1}
\forall_{\alpha, \beta \in \R} \, \left( \alpha, \beta \ge 0 \Longrightarrow \alpha \, \beta \le \frac{\alpha^p}{p} + \frac{\beta^q}{q} \right).\footnote{Tip: Die Funktion
$$ {[} 0, \infty {[} \longrightarrow \R, ~~ t \longmapsto \frac{\alpha^p}{p} + \frac{t^q}{q} - \alpha \, t, $$
besitzt für festes $\alpha \in {[} 0, \infty {[}$ in $\alpha^{p-1}$ ein absolutes Minimum mit Funktionswert Null.}
\end{equation}
Aus (\ref{FA.5.8.1}) folgt
\begin{equation} \label{FA.5.8.2}
\frac{|f|}{\|f\|_p} \, \frac{|g|}{\|g\|_q} \le \frac{1}{p} \, \frac{|f|^p}{{\|f\|_p}^p} + \frac{1}{q} \, \frac{|g|^q}{{\|g\|_q}^q},
\end{equation}
wobei nach Voraussetzung auf der rechten Seite eine über $M$ $\varphi$-integrierbare Funktion und auf der linken Seite eine über $M$ $\varphi$-meßbare steht.
Nach \ref{FA.4.70} (i) ist daher auch die linke Seite von (\ref{FA.5.8.2}) $\varphi$-integrierbar über $M$, also wegen \ref{FA.4.70} (ii):\linebreak $f \, g \in \L_{\K}(M,\varphi)$.
Integration der Ungleichung (\ref{FA.5.8.2}) ergibt nun
\begin{eqnarray*}
\frac{1}{\|f\|_p \, \|g\|_q} \, \int_M |f \, g| \, \d \varphi & \le & \frac{1}{p} \, \frac{\int_M |f|^p \, \d \varphi}{{\|f\|_p}^p} + \frac{1}{q} \, \frac{\int_M |g|^q \, \d \varphi}{{\|g\|_q}^q} \\
& = & \frac{1}{p} + \frac{1}{q} = 1, 
\end{eqnarray*}
und dies ist gleichbedeutend mit der Behauptung. 

Für den Rest des Beweises können wir nun ohne Einschränkung annehmen, daß gilt $p = 1$ und $q = \infty$.
Dann ist $f$ $\varphi$-integrierbar über $M$ und $g$ $\varphi$-meßbar über $M$ sowie $\varphi$-beschränkt auf $M$, also folgt aus \ref{FA.4.70} (iii) die $\varphi$-In\-te\-grier\-bar\-keit von $f \, g$ über $M$.
$\| f \, g \|_1 \le \|f \|_1 \, \| g \|_{\infty}$ ist klar. \q

\begin{Bsp} \label{FA.5.9} 
Seien $p,q \in {]}1,\infty{[}$ mit $\frac{1}{p} + \frac{1}{q} = 1$. \\
Dann gilt $(x_i \, y_i)_{i \in \N} \in \ell_{\K}^1$ für alle $(x_i)_{i \in \N} \in \ell_{\K}^p$ und alle $(y_i)_{i \in \N} \in \ell_{\K}^q$ sowie
$$ \sum_{i=0}^{\infty} |x_i| \, |y_i| \le \left( \sum_{i=0}^{\infty} |x_i|^p \right)^{\frac{1}{p}} \left( \sum_{i=0}^{\infty} |y_i|^q \right)^{\frac{1}{q}}. $$
\end{Bsp}

\begin{Satz}[Jensensche Ungleichung] \index{Ungleichung!Jensensche} \label{FA.5.10}
Es seien $p, \tilde{p} \in \R_+$ mit $p \le \tilde{p}$ und $(x_i)_{i \in \N}$ eine Folge in $\K$ derart, daß die Reihe $\sum_{i=0}^{\infty} {x_i}^p$ absolut konvergiert.

Dann konvergiert auch die Reihe $\sum_{i=0}^{\infty} {x_i}^{\tilde{p}}$ absolut, und es folgt
$$ \sup \{ |x_i| \, | \, i \in \N \} \le \left( \sum_{i=0}^{\infty} |x_i|^{\tilde{p}} \right)^{\frac{1}{\tilde{p}}} \le \left( \sum_{i=0}^{\infty} |x_i|^p \right)^{\frac{1}{p}}. $$
\end{Satz}

\textit{Beweis.} Im Falle $\sum_{i=0}^{\infty} |x_i|^p = 0$ ist die Behauptung klar.
Daher können wir ohne Einschränkung annehmen, daß $\sum_{i=0}^{\infty} |x_i|^p = 1$ gilt.
Dann gilt $\forall_{i \in \N} \, |x_i| \le 1$.
Hieraus ergibt sich zum einen $\sup \{ |x_i| \, | \, i \in \N \} \le 1 = ( \sum_{i=0}^{\infty} |x_i|^p )^{\frac{1}{p}}$, und dies gilt auch im Falle $p = \tilde{p}$, also gilt die erste Ungleichung.
Zum anderen folgt wegen $\tilde{p} \ge p$: $\forall_{i \in \N} \, |x_i|^{\tilde{p}} \le |x_i|^p$.
Die absolute Konvergenz der Reihe $\sum_{i=0}^{\infty} {x_i}^p$ ergibt somit die der Reihe $\sum_{i=0}^{\infty} {x_i}^{\tilde{p}}$ und $\sum_{i=0}^{\infty} |x_i|^{\tilde{p}} \le \sum_{i=0}^{\infty} |x_i|^p = 1$.
Letzteres impliziert die zweite Ungleichung: $( \sum_{i=0}^{\infty} |x_i|^{\tilde{p}} )^{\frac{1}{\tilde{p}}} \le 1 = ( \sum_{i=0}^{\infty} |x_i|^p )^{\frac{1}{p}}$. \q

\begin{Satz}[Minkowskische Ungleichung] \index{Ungleichung!Minkowskische}\index{Minkowski!-Ungleichung} \label{FA.5.11}
Es seien $\varphi$ ein Quadermaß auf $\R^n$, $M$ eine Teilmenge von $\R^n$, $p \in {[}1,\infty{]}$ und $f,g \in \L^p_{\K}(M,\varphi)$.

Dann folgt 
$$ \| f + g \|_p \le \| f \|_p + \| g \|_p. $$
\end{Satz}

\textit{Beweis.} Für $p=1$ ist wegen
$$ \| f + g \|_1 = \int_M | f + g | \, \d \varphi \le \int_M (|f| + |g|) \, \d \varphi = \|f\|_1 + \|g\|_1 $$
nichts zu zeigen, und für $p = \infty$ ist die Behauptung wegen 
$$ | f + g | \le |f| + |g| $$
klar.

Sei also $p \in {]}1, \infty{[}$.
Ferner sei $q \in {]}1, \infty{[}$ mit $\frac{1}{p} + \frac{1}{q} = 1$.
Dann gilt
\begin{gather}
| f + g |^p = | f + g | \, | f + g |^{p-1} \le |f| \, | f + g |^{p-1} + |g| \, | f + g |^{p-1}, \label{FA.5.11.1} \\
| f + g |^p \in \L_{\R}(M,\varphi), \label{FA.5.11.2} \\
|f|, |g| \in \L_{\R}^p(M,\varphi), \label{FA.5.11.3} \\
| f + g|^{p-1} \in \L^q_{\R}(M,\varphi). \label{FA.5.11.4}
\end{gather}

{[} (\ref{FA.5.11.1}) - (\ref{FA.5.11.3}) sind klar.

Zu (\ref{FA.5.11.4}): Die $\varphi$-Meßbarkeit von $| f + g |^{p-1}$ über $M$ folgt aus (\ref{FA.5.11.2}) und der Stetigkeit der Funktion $| \ldots |^{\frac{p-1}{p}} \: \R \to \R$ mit $0^{\frac{p-1}{p}}$ sowie \ref{FA.4.8} (iii).
Des weiteren ist wiederum wegen (\ref{FA.5.11.2}) auch $| f + g |^{(p-1)q} \stackrel{\ref{FA.5.pq}}{=} |f+g|^p$ $\varphi$-integrierbar über $M$. {]}

(\ref{FA.5.11.1}) - (\ref{FA.5.11.4}) und die Höldersche Ungleichung \ref{FA.5.8} ergeben offenbar
$$ \int_M | f + g |^p \, \d \varphi \le \|f\|_p \, \|f+g\|_q + \|g\|_p \, \|f+g\|_q. $$
Ist $|f+g| =_{\varphi} 0$, so ist die Behauptung trivial.
Andernfalls ergibt sie sich aus der letzten Ungleichung durch Division durch $\|f+g\|_q = ( \int_M | f + g |^p \, \d \varphi )^{\frac{1}{q}}$ wegen $1 - \frac{1}{q} = \frac{1}{p}$.
\q
\A
Die nächsten drei Sätze behandeln i.w.\ Varianten der Ungleichungen von \textsc{Hölder} und \textsc{Minkowski} für den Fall $p \in {]}0,1{[}$.

\begin{Satz} \label{FA.5.12}
Seien $\varphi$ ein Quadermaß auf $\R^n$, $M$ eine Teilmenge von $\R^n$ und 
$$ \mbox{$p \in {]}0,1{[}$ sowie $q \in \R_-$ mit $\frac{1}{p} + \frac{1}{q} = 1$.} $$
Ferner seien $f \in \L_{\K}^p(M,\varphi)$ und $g \in \mathcal{M}_{\K}(M,\varphi)$ mit $g \ne_{\varphi} 0$ sowie $|g|^q \in \L_{\R}(M,\varphi)$ derart, daß gilt $f \, g \in \L_{\K}(M,\varphi)$.\footnote{$f \, g \in \L_{\K}(M,\varphi)$ ist nach \ref{FA.4.70} (iii) z.B.\ erfüllt, wenn $g \in \mathcal{M}_{\K}(M,\varphi)$ $\varphi$-beschränkt ist, d.h.\ $g \in \L_{\K}^{\infty}(M,\varphi)$.}

Dann gilt
$$ \underbrace{\int_M | f \, g | \, \d \varphi}_{= \| f \, g \|_1} \ge \underbrace{\left( \int_M |f|^p \, \d \varphi \right)^{\frac{1}{p}}}_{= \|f\|_p} \left( \int_M |g|^q \, \d \varphi \right)^{\frac{1}{q}}. $$
\end{Satz}

\textit{Beweis.} Seien $\tilde{p} := \frac{1}{p}$ und $\tilde{q} := \frac{1}{1-p}$.
Dann folgt $\tilde{p}, \tilde{q} \in {]}1, \infty{[}$ sowie $\frac{1}{\tilde{p}} + \frac{1}{\tilde{q}} = 1$ und
\begin{equation} \label{FA.5.12.1}
-p \, \tilde{q} = \frac{p}{p-1} \stackrel{\ref{FA.5.pq}}{=} q.
\end{equation}
Wir behaupten
\begin{gather} 
\tilde{f} := |f \, g|^p \in \L_{\R}^{\tilde{p}}(M,\varphi), \label{FA.5.12.2} \\
\tilde{g} := |g|^{-p} \in \L_{\R}^{\tilde{q}}(M,\varphi). \label{FA.5.12.3}
\end{gather}

{[} (\ref{FA.5.12.2}) bedeutet genau $f \, g \in \L_{\K}(M,\varphi)$, und (\ref{FA.5.12.3}) ergibt sich sofort aus der Voraussetzung an $g$ und (\ref{FA.5.12.1}). {]}

Die Höldersche Ungleichung \ref{FA.5.8} -- angewandt auf $\tilde{p}, \tilde{q}, \tilde{f}, \tilde{g}$ -- ergibt nun
$$ \int_M |f|^p \, \d \varphi \le \left( \int_M |f \, g| \, \d \varphi \right)^{\frac{1}{\tilde{p}}} \left( \int_M |g|^{-p \tilde{q}} \, \d \varphi \right)^{\frac{1}{\tilde{q}}}, $$
also wegen $g \ne_{\varphi} 0$ und der Definition von $\tilde{q}, \tilde{p}$
\begin{gather*}
\left( \int_M |f|^p \, \d \varphi \right) \left( \int_M |g|^{-p \tilde{q}} \, \d \varphi \right)^{p-1} \le \left( \int_M |f \, g| \, \d \varphi \right)^p, \\
\int_M |f \, g| \, \d \varphi \ge \left( \int_M |f|^p \, \d \varphi \right)^{\frac{1}{p}} \left( \int_M |g|^{-p \tilde{q}} \, \d \varphi \right)^{\frac{p-1}{p}},
\end{gather*}
womit der Satz wegen (\ref{FA.5.12.1}) bewiesen ist. \q

\begin{Satz} \label{FA.5.13}
Es seien $\varphi$ ein Quadermaß auf $\R^n$, $M$ eine Teilmenge von $\R^n$, $p \in {]}0,1{[}$ und $f,g \in \L^p_{\R}(M,\varphi)$ mit $f,g \ge_{\varphi} 0$ sowie $f+g >_{\varphi} 0$.

Dann folgt 
$$ \| f + g \|_p \ge \| f \|_p + \| g \|_p. $$
\end{Satz}

\textit{Beweis.} Wir wählen $q \in \R_-$ mit $\frac{1}{p} + \frac{1}{q} = 1$.
Analog zum Beweis der Minkowskischen Ungleichung \ref{FA.5.11} (dort im Falle $p \in {]} 1, \infty {[}$) gilt
\begin{gather*}
( f + g )^p = ( f + g ) \, ( f + g )^{p-1} =_{\varphi} f \, ( f + g )^{p-1} + g \, ( f + g )^{p-1}, \label{FA.5.13.1} \\
( f + g )^p \in \L_{\R}(M,\varphi), \label{FA.5.13.2} \\
f, g \in \L_{\R}^p(M,\varphi), \label{FA.5.13.3} \\
( f + g )^{(p-1)} \in \mathcal{M}_{\R}(M,\varphi), \label{FA.5.13.4} \\
| f + g |^{(p-1)q} \in \L_{\R}(M,\varphi), \label{FA.5.13.5}
\end{gather*}
also wegen $f+g >_{\varphi} 0$ nach \ref{FA.5.12}
\begin{eqnarray*}
\lefteqn{\int_M | f + g |^p \, \d \varphi} \\
& & \ge \left( \int_M |f|^p \, \d \varphi \right)^{\frac{1}{p}} \left( \int_M | f + g |^p \, \d \varphi \right)^{\frac{1}{q}} + \left( \int_M |g|^p \, \d \varphi \right)^{\frac{1}{p}} \left( \int_M | f + g |^p \, \d \varphi \right)^{\frac{1}{q}}
\end{eqnarray*}
und Division durch $( \int_M | f + g |^p \, \d \varphi )^{\frac{1}{q}}$ beweist den Satz.
\q

\begin{Bem*} 
Die Prämisse $f,g \ge_{\varphi} 0$ und $f+g >_{\varphi} 0$ in \ref{FA.5.13} ist für die Konklusion von \ref{FA.5.13} notwendig:
Sind nämlich $f := \chi_{[0,1]}, \, g_1:= - \frac{f}{2}, \, g_2 := -f \in \L^{\frac{1}{2}}_{\R}([0,1],\mu_1)$, so gilt $f > 0$, $g_1 < 0$, $f + g_1 > 0$ sowie $\underbrace{f + g_2 = 0}_{\Rightarrow \, g_2 = -f}$ und
\begin{gather*}
\| f + g_1 \|_{\frac{1}{2}} = \frac{1}{2} < 1 + \frac{1}{2} = \|f\|_{\frac{1}{2}} + \|g_1\|_{\frac{1}{2}}, \\
\| f + g_2 \|_{\frac{1}{2}} = 0 < 2 = \|f\|_{\frac{1}{2}} + \|g_2\|_{\frac{1}{2}}.
\end{gather*} 
Dies zeigt auch, daß ein zu \ref{FA.5.13} analoges Resultat für $f,g \in \L^p_{\C}(M,\varphi)$ mit $f,g \ne_{\varphi} 0$ und $f+g \ne_{\varphi} 0$ i.a.\ nicht gelten kann.
\end{Bem*}

\begin{Satz} \label{FA.5.14}
Es seien $\varphi$ ein Quadermaß auf $\R^n$, $M$ eine Teilmenge von $\R^n$, $p \in {]}0,1{[}$ und $f,g \in \L^p_{\K}(M,\varphi)$.

Dann gilt:
\begin{itemize}
\item[(i)] ${\|f+g\|_p}^p \le {\|f\|_p}^p + {\|g\|_p}^p$,
\item[(ii)] $\| f + g \|_p \le 2^{\frac{1}{p}-1} \left( \|f\|_p + \|g\|_p \right)$.
\end{itemize}
\end{Satz}

\textit{Beweis.} Zu (i): Zunächst gilt
\begin{equation} \label{FA.5.14.1}
\forall_{\alpha, \beta \in {[}0, \infty{[}} \, (\alpha + \beta)^p \le \alpha^p + \beta^p.
\end{equation} 

{[} Zu (\ref{FA.5.14.1}): Die Funktion
$$ {[}0, \infty{[} \longrightarrow \R, ~~ t \longmapsto 1 + t^p - (1 + t)^p, $$
ist wegen $p-1 < 0$ streng monoton fallend und bildet Null auf Null ab.
Im Falle $\beta \ne 0$ erhält man die Behauptung mit $t := \frac{\alpha}{\beta}$, beachte $1 + \frac{\alpha}{\beta} = \frac{\alpha + \beta}{\beta}$, und im Falle $\beta = 0$ ist sie trivial. {]}

Aus (\ref{FA.5.14.1}) folgt $|f+g|^p \stackrel{\text{klar}}{\le} (|f|+|g|)^p \le |f|^p + |g|^p$, also durch Integration auch die Ungleichung (i).

Zu (ii): Es gilt
\begin{equation} \label{FA.5.14.2}
\forall_{\alpha, \beta \in {[}0, \infty{[}} \forall_{\tilde{p} \in {]} 1, \infty {[}} \, (\alpha + \beta)^{\tilde{p}} \le 2^{\tilde{p}-1} \left( \alpha^{\tilde{p}} + \beta^{\tilde{p}} \right).
\end{equation} 

{[} Zu (\ref{FA.5.14.2}): Zeige als Übung, daß die Funktion
$$ {[}0, \infty{[} \longrightarrow \R, ~~ t \longmapsto \frac{1 + t^{\tilde{p}}}{(1 + t)^{\tilde{p}}}, $$ 
für $\tilde{p} \in {]} 1, \infty {[}$ in Eins ein absolutes Minimum mit Funktionswert $2^{\tilde{p}-1}$ besitzt.
Nun argumentiert man wie im Beweis von (\ref{FA.5.14.1}). {]} 

Mit $\alpha := \int_M |f|^p \, \d \varphi$, $\beta := \int_M |g|^p \, \d \varphi$ und $\tilde{p} := \frac{1}{p}$ gilt nach (i) 
$$ \|f+g\|_p = \left( \int_M |f+g|^p \, \d \varphi \right)^{\tilde{p}} \le \left( \int_M |f|^p \, \d \varphi + \int_M |g|^p \, \d \varphi \right)^{\tilde{p}} = (\alpha + \beta)^{\tilde{p}}, $$
also folgt die Ungleichung (ii) aus (\ref{FA.5.14.2}). \q
\A
Die im folgenden Satz genannten Ungleichungen, welche in \cite{Hirz} und \cite{Reck} \emph{Par\-al\-le\-lo\-gramm\-un\-gleich\-ungen}\index{Ungleichung!Parallelogramm-}\index{Parallelogramm!-ungleichung} genannt werden, benötigen wir im siebenten Kapitel, um den Satz von \textsc{Clarkson} zu beweisen.

\begin{Satz}[Clarksonsche Ungleichungen] \index{Ungleichung!Clarksonsche} \label{FA.5.15} $\,$

\noindent \textbf{Vor.:} Seien $\varphi$ ein Quadermaß auf $\R^n$, $M$ eine Teilmenge von $\R^n$, $p \in {]} 1, \infty {[}$ und $f,g \in \L^p_{\K}(M,\varphi)$.

\noindent \textbf{Beh.:}
\begin{itemize}
\item[(i)] $\D p \in {[}2, \infty{[} \Longrightarrow {\|f+g\|_p}^p + {\|f-g\|_p}^p \le 2^{p-1} \left( {\|f\|_p}^p + {\|g\|_p}^p \right)$.
\item[(ii)] Ist $p \in {]} 1, 2 {[}$ und $q := \frac{1}{1 - \frac{1}{p}} \in {]} 2, \infty {[}$, d.h.\ $\frac{1}{p} + \frac{1}{q} = 1$, so gilt
$$ {\|f+g\|_p}^q + {\|f-g\|_p}^q \le 2 \left( {\|f\|_p}^p + {\|g\|_p}^p \right)^{q-1}. $$
\end{itemize}
\end{Satz}

Um diesen Satz zu beweisen, benötigen wir zwei Lemmata.
Wir erinnern zunächst daran, daß für jedes $\zeta \in \C$ die Folge $\left( {\zeta \choose k} \right)_{k \in \N}$ der \emph{Binominalkoeffizienten}\index{Binominal!-koeffizienten} rekursiv durch
\begin{equation*} \label{FA.5.16.S}
\boxed{{\zeta \choose 0}} := 1 ~~ \mbox{ sowie  } ~~ \forall_{k \in \N_+} \, \boxed{{\zeta \choose k}} := {\zeta \choose k - 1} \, \frac{\zeta - (k - 1)}{k},
\end{equation*}
d.h.\
\begin{equation} \label{FA.5.16.DS}
\forall_{k \in \N} \, {\zeta \choose k} = \frac{ \prod_{\kappa = 1}^k (\zeta - (\kappa - 1))}{k!},
\end{equation}
definiert ist.
Aus (\ref{FA.5.16.DS}) folgt leicht
\begin{gather}
\forall_{k \in \N_+} \forall_{\zeta \in \C \setminus \{2k\}} \, {\zeta \choose 2k} = \frac{\zeta \, (\zeta - 1)}{2 k \, (\zeta - 2 k)} \cdot \frac{\prod_{\kappa=2}^{2k} (\zeta - \kappa)}{(2 k - 1)!} = \frac{\zeta \, (\zeta - 1)}{2 k \, (2k - \zeta)} \cdot \frac{\prod_{\kappa=2}^{2k} (\kappa - \zeta)}{(2 k - 1)!}, \label{FA.5.16.DS1} \\
\forall_{k \in \N_+} \forall_{\zeta \in \C} \, {\zeta - 1 \choose 2k} = \frac{\zeta - 1}{2 k} \cdot \frac{\prod_{\kappa=2}^{2k} (\zeta - \kappa)}{(2 k - 1)!} = - \frac{\zeta - 1}{2 k} \cdot \frac{\prod_{\kappa=2}^{2k} (\kappa - \zeta)}{(2 k - 1)!}, \label{FA.5.16.DS2} \\
\forall_{k \in \N_+} \forall_{\zeta \in \C \setminus \{2k\}} \, {\zeta - 1 \choose 2k - 1} = \frac{\zeta - 1}{\zeta - 2 k} \cdot \frac{\prod_{\kappa=2}^{2k} (\zeta - \kappa)}{(2 k - 1)!} = \frac{\zeta - 1}{2k - \zeta} \cdot \frac{\prod_{\kappa=2}^{2k} (\kappa - \zeta)}{(2 k - 1)!}. \label{FA.5.16.DS3}
\end{gather}

\begin{Lemma}[Binominalentwicklung] \index{Binominal!-entwicklung} \label{FA.5.16.BE}
Für alle $z \in \C$ mit $|z| < 1$ und $p \in \R_+$ gilt
$$ (1+z)^p = \sum_{k=0}^{\infty} {p \choose k} \, z^k, $$
wobei die Potenzreihe auf der rechten Seite der Gleichung absolut konvergent ist.
\end{Lemma}

\textit{Beweis.} Das Lemma folgt natürlich sofort aus dem (mittels des Cauchyschen Integralsatzes zu beweisenden) Taylorschen Satz der Funktionentheorie durch Taylorentwicklung der holomorphen Funktion 
$$ (1+ \ldots)^p = \exp( p \, \log(1+ \ldots) ) \: {\C \setminus {]}- \infty, -1{]}} \longrightarrow \C $$
in Null, wobei $\log \: \C \setminus{]}- \infty, 0 {]} \to \C$ den \emph{Hauptzweig des Logarithmus} bezeichne, d.h.\ per definitionem
$$ \forall_{r \in \R_+} \forall_{\phi \in {]} - \pi, \pi {[}} \, \log( r \, e^{\i \phi}) = \ln(r) + \i \, \phi, $$
aber für $z \in {]}-1,1{[}$ auch aus dem Taylorschen Satz der Analysis I, vgl.\ ggf.\ \cite[Hauptsatz 6.24]{ElAna}.
Dieses einzusehen, überlassen wir dem Leser als Übung. \q

\begin{Lemma} \label{FA.5.16}
Seien $p \in {]} 1, 2 {[}$ und $q \in {]} 2, \infty {[}$ mit $\frac{1}{p} + \frac{1}{q} = 1$.

Dann gilt $\D \forall_{\tau \in {[}0,1{]}} \, (1 + \tau^q)^{p-1} \le \frac{1}{2} \left( \left( 1 + \tau \right)^p +  \left( 1 - \tau \right)^p \right)$.
\end{Lemma}

\textit{Beweis.} Im Falle $\tau \in \{0, 1\}$ ist die Behauptung trivial.
Sei im folgenden daher $\tau \in {]}0,1{[}$.
Aus \ref{FA.5.16.BE} folgt dann
\begin{gather*}
(1 + \tau)^p = \sum_{k=0}^{\infty} {p \choose k} \, \tau^k, \\
(1 - \tau)^p = \sum_{k=0}^{\infty} (-1)^k \, {p \choose k} \, \tau^k, \\
(1 + \tau^q)^{p-1} = \sum_{k=0}^{\infty} {p-1 \choose k} \, \tau^{q k}
\end{gather*}
\pagebreak
und
\begin{eqnarray*}
(1 + \tau)^p + (1 - \tau)^p & = & 2 \sum_{k=0}^{\infty} \underbrace{{p \choose 2 k} \, \tau^{2k}}_{= 1 \text{ für } k=0}, \\
(1 + \tau^q)^{p-1} & = & \sum_{k=0}^{\infty} \underbrace{{p-1 \choose 2 k} \, \tau^{2 k q}}_{= 1 \text{ für } k=0} + \sum_{k=1}^{\infty} {p-1 \choose 2 k - 1} \, \tau^{(2k-1) q},
\end{eqnarray*}
also wegen $p \notin 2 \N_+$ sowie (\ref{FA.5.16.DS1}) - (\ref{FA.5.16.DS3}) auch
\begin{eqnarray*}
\lefteqn{\frac{1}{2} \left( \left( 1 + \tau \right)^p +  \left( 1 - \tau \right)^p \right) - (1 + \tau^q)^{p-1}} \\
&& = \sum_{k=1}^{\infty} \left( {p \choose 2 k} - {p - 1 \choose 2 k} \, \tau^{2 k (q-1)} -  {p - 1 \choose 2 k - 1} \, \tau^{2 k (q-1) - q} \right) \tau^{2k} \\
&& = \sum_{k=1}^{\infty} \underbrace{\frac{\prod_{\kappa=2}^{2k} (\kappa - p)}{(2 k - 1)!}}_{> 0, \text{ da } p < 2} \left( \frac{p \, (p - 1)}{2 k \, (2k - p)} + \frac{p-1}{2 k} \, \tau^{2 k (q-1)} - \frac{p-1}{2k - p} \, \tau^{2 k (q-1) - q} \right) \underbrace{\tau^{2k}}_{> 0}.
\end{eqnarray*}
Zum Nachweis des Lemmas bleibt daher zu zeigen
\begin{equation} \label{FA.5.16.1}
\forall_{k \in \N_+} \, \frac{p \, (p - 1)}{2 k \, (2k - p)} + \frac{p-1}{2 k} \, \tau^{2 k (q-1)} - \frac{p-1}{2k - p} \, \tau^{2 k (q-1) - q} \ge 0. 
\end{equation}

Zu (\ref{FA.5.16.1}):
Sei $k \in \N_+$.
Zunächst gilt
\begin{equation} \label{FA.5.16.2}
\frac{p \, (p - 1)}{2 k \, (2k - p)} + \frac{p-1}{2 k} \, \tau^{2 k (q-1)} - \frac{p-1}{2k - p} \, \tau^{2 k (q-1) - q}
= \frac{1 - \tau^{\frac{2k - p}{p-1}}}{\frac{2k - p}{p-1}} - \frac{1 - \tau^{\frac{2k}{p-1}}}{\frac{2k}{p-1}}.
\end{equation}

{[} Die rechte Seite in (\ref{FA.5.16.2}) ist nach \ref{FA.5.pq} gleich
\begin{eqnarray*}
\lefteqn{\frac{p-1}{2k \, (2k-p)} \left( 2k \left( 1 - \tau^{2k (q-1) -q} \right) - (2k - p) \left( 1 - \tau^{2k (q-1)} \right) \right)} \\
&& = \frac{p-1}{2k \, (2k-p)} \left( p + (2k - p) \, \tau^{2k (q-1)} - 2k \, \tau^{2k (q-1) -q} \right)
\end{eqnarray*}
und somit gleich der linken Seite in (\ref{FA.5.16.2}). {]}

Erneut sei es dem Leser überlassen, zu zeigen, daß die Funktion
$$ \R_+ \longrightarrow \R, ~~ t \longmapsto \frac{1 - \tau^{\frac{t}{p-1}}}{\frac{t}{p-1}}, $$
wegen $\tau \in {]}0,1{[}$ und $p>1$ streng monoton fallend ist.
Hieraus und aus (\ref{FA.5.16.2}) folgt dann (\ref{FA.5.16.1}) für $k$. \q
\A
\textit{Beweis der Clarksonschen Ungleichungen \ref{FA.5.15}.} 
Zu (i): Es sei $p \in {[} 2, \infty {[}$.
Wir wollen
$$ \forall_{f,g \in \L_{\K}^p(M, \varphi)} \, {\|f+g\|_p}^p + {\|f-g\|_p}^p \le 2^{p-1} \left( {\|f\|_p}^p + {\|g\|_p}^p \right) $$
zeigen und behaupten zunächst
\begin{equation} \label{FA.5.15.1}
\forall_{\zeta, \xi \in \C} \, | \zeta + \xi |^p + | \zeta - \xi |^p \le 2^{p-1} \left( |\zeta|^p + |\xi|^p \right).
\end{equation}

{[} Nachweis hiervon: Seien $\zeta, \xi \in \C$.
Dann existiert $t \in {[} 0 , 2 \,\pi {]}$ mit
\begin{eqnarray*}
| \zeta \pm \xi |^2 & = & ( \zeta \pm \xi ) \, ( \overline{\zeta \pm \xi} ) = \zeta \, \overline{\zeta} \pm \zeta \, \overline{\xi} \pm \overline{ \overline{\zeta}} \, \xi + \xi \, \overline{\xi} = |\zeta|^2 + |\xi|^2 \pm 2 \, {\rm Re}(\zeta \, \overline{\xi}) \\
& = & |\zeta|^2 + |\xi|^2 \pm 2 \, |\zeta| \, |\xi| \, \cos(t),
\end{eqnarray*}
also gilt
$$ | \zeta + \xi |^p + | \zeta - \xi |^p = ( |\zeta|^2 + |\xi|^2 + 2 \, |\zeta| \, |\xi| \, \cos(t))^{\frac{p}{2}} - ( |\zeta|^2 + |\xi|^2 - 2 \, |\zeta| \, |\xi| \, \cos(t))^{\frac{p}{2}}. $$
Wir überlassen es wieder dem Leser, nachzuweisen, daß die Funktion 
\begin{eqnarray*}
\R & \longrightarrow & \R \\
t & \longmapsto & ( |\zeta|^2 + |\xi|^2 + 2 \, |\zeta| \, |\xi| \, \cos(t))^{\frac{p}{2}} - ( |\zeta|^2 + |\xi|^2 - 2 \, |\zeta| \, |\xi| \, \cos(t))^{\frac{p}{2}}
\end{eqnarray*}
wegen $p \ge 2$ für jedes $k \in \Z$ auf ${[} (k-1) \, \pi , k \, \pi - \frac{\pi}{2} {]}$ monoton fallend und auf ${[} k \, \pi - \frac{\pi}{2} , k \, \pi {]}$ monoton wachsend ist.
Ihre Beschränkung auf ${[} 0 , 2 \, \pi {]}$ nimmt daher ihr globales Maximum in $0$, $\pi$ oder $2 \, \pi$ mit einem Funktionswert gleich
\begin{eqnarray*}
\lefteqn{( |\zeta|^2 + |\xi|^2 \pm 2 \, |\zeta| \, |\xi|)^{\frac{p}{2}} - ( |\zeta|^2 + |\xi|^2 \mp 2 \, |\zeta| \, |\xi|)^{\frac{p}{2}}} && \\
&& = ( (|\zeta| \pm |\xi|)^2 )^{\frac{p}{2}} - ( (|\zeta| \mp |\xi|)^2 )^{\frac{p}{2}} =  | \, |\zeta| \pm |\xi| \, |^p - | \, |\zeta| \mp |\xi| \, |^p \\
&& = \pm ( | \, |\zeta| + |\xi| \, |^p - | \, |\zeta| - |\xi| \, |^p ),
\end{eqnarray*}
an.
Ohne Einschränkung gelte $|\zeta| \ge |\xi|$.
Dann folgt aus den bisherigen Überlegungen
$$ | \zeta + \xi |^p + | \zeta - \xi |^p \le | \, |\zeta| + |\xi| \, |^p - | \, |\zeta| - |\xi| \, |^p, $$
und zum Beweis von (\ref{FA.5.15.1}) bleibt zu zeigen, daß für $\alpha := |\zeta|, \beta := |\xi| \in \R$ gilt
\begin{equation} \label{FA.5.15.1s}
| \alpha + \beta |^p + | \alpha - \beta |^p \le 2^{p-1} \left( |\alpha|^p + |\beta|^p \right).
\end{equation}

Zu (\ref{FA.5.15.1s}): Wegen $p \ge 2$ folgt aus der Jensenschen Ungleichung \ref{FA.5.10} -- angewandt auf die Folge $(\alpha + \beta, \alpha - \beta, 0, \ldots )$ --
$$ \left( | \alpha + \beta |^p + | \alpha - \beta |^p  \right)^{\frac{1}{p}} \le \left( | \alpha + \beta |^2 + | \alpha - \beta |^2  \right)^{\frac{1}{2}} = \sqrt{2} \left( \alpha^2 + \beta^2 \right)^{\frac{1}{2}}. $$
Hieraus und aus der noch zu zeigenden Ungleichung
\begin{equation} \label{FA.5.15.2}
\alpha^2 + \beta^2 \le 2^{\frac{p-2}{p}} \left( |\alpha|^p + |\beta|^p \right)^{\frac{2}{p}}
\end{equation}
ergibt sich
\begin{equation*}
\left( | \alpha + \beta |^p + | \alpha - \beta |^p  \right)^{\frac{1}{p}} \le \sqrt{2} \, 2^{\frac{p-2}{2 p}} \left( |\alpha|^p + |\beta|^p \right)^{\frac{1}{p}} = 2^{\frac{p-1}{p}} \left( |\alpha|^p + |\beta|^p \right)^{\frac{1}{p}},
\end{equation*}
also gilt (\ref{FA.5.15.1s}).

Zu (\ref{FA.5.15.2}): Im Falle $p=2$ ist die Behauptung trivial.
Sei daher $p > 2$.
Wir setzen  nun $\tilde{q} := \frac{p}{p-2} \in {]} 1, \infty {[}$, d.h.\ $\frac{1}{\frac{p}{2}} + \frac{1}{\tilde{q}} = 1$.
Beispiel \ref{FA.5.9} -- angewandt auf die Folgen $(\alpha^2, \beta^2, 0, \ldots)$ und $(1,1,0, \ldots)$ mit $\tilde{p} := \frac{p}{2}, \, \tilde{q}$ anstelle von $p,q$ -- ergibt
$$ \alpha^2 + \beta^2 \le \left( |\alpha|^p + |\beta|^p \right)^{\frac{2}{p}} \, (1+1)^{\frac{p-2}{p}}, $$
und auch (\ref{FA.5.15.2}) ist gezeigt. {]}
\pagebreak

Aus (\ref{FA.5.15.1}) folgt für alle $f,g \in \L_{\K}^p(M,\varphi)$
$$ | f + g |^p + | f - g |^p \le 2^{p-1} \left( |f|^p + |g|^p \right), $$
und (i) ergibt sich durch Integration.

Zu (ii): Seien $p \in {]}1,2{[}$ und $q \in {]} 2, \infty {[}$ mit $\frac{1}{p} + \frac{1}{q} = 1$.
Nun wollen wir
$$ \forall_{f,g \in \L_{\K}^p(M, \varphi)} \, {\|f+g\|_p}^q + {\|f-g\|_p}^q \le 2 \left( {\|f\|_p}^p + {\|g\|_p}^p \right)^{q-1} $$
zeigen.
Es gilt
\begin{equation} \label{FA.5.15.3}
\forall_{\zeta, \xi \in \C} \, | \zeta + \xi |^q + | \zeta - \xi |^q \le 2 \left( |\zeta|^p + |\xi|^p \right)^{q-1}.
\end{equation}

{[} Beweis hiervon: Analog zum Beweis von (\ref{FA.5.15.1}) begründet man, daß es wegen $q > 2$ genügt, für $\alpha := |\zeta|, \beta := |\xi| \in \R$
\begin{equation} \label{FA.5.15.3s}
| \alpha + \beta |^q + | \alpha - \beta |^q \le 2 \left( |\alpha|^p + |\beta|^p \right)^{q-1}
\end{equation}
zu zeigen.
Wir definieren $\tilde{\alpha}, \tilde{\beta} \in \R$ durch
$$ \alpha + \beta = 2 \tilde{\alpha} ~~ \mbox{ und } ~~ \alpha - \beta = 2 \tilde{\beta}, $$
d.h.\ $\alpha = \tilde{\alpha} + \tilde{\beta}$ sowie $\beta = \tilde{\alpha} - \tilde{\beta}$.
Dann bedeutet (\ref{FA.5.15.3s}) genau
\begin{equation} \label{FA.5.15.4}
2^q \left( |\tilde{\alpha}|^q  + |\tilde{\beta}|^q \right) \le 2 \left( |\tilde{\alpha} + \tilde{\beta}|^p + |\tilde{\alpha} - \tilde{\beta}|^p \right)^{q-1},
\end{equation}
und ohne Beschränkung der Allgemeinheit können wir zum Nachweis dieser Ungleichung $\tilde{\alpha} \ge \tilde{\beta} > 0$ annehmen.
Sei weiterhin
$$ \tau := \frac{\tilde{\beta}}{\tilde{\alpha}} \in {]} 0, 1 {]}. $$
Indem wir (\ref{FA.5.15.4}) mit $(\frac{1}{\tilde{\alpha}})^q$ multiplizieren, sehen wir unter Berücksichtigung von $\frac{q}{(q-1) \, p} \stackrel{\ref{FA.5.pq}}{=} 1$, daß (\ref{FA.5.15.4}) nun mit
$$ (1 + \tau^q) \le  \frac{1}{2^{q-1}}  \left( (1 + \tau)^p - (1 - \tau)^p \right)^{q-1}, $$
also wegen $(p-1) \, (q-1) \stackrel{\ref{FA.5.pq}}{=} 1$ auch mit
$$ (1 + \tau^q)^{p-1} \le \frac{1}{2} \left( (1 + \tau)^p - (1 - \tau)^p \right) $$
gleichbedeutend ist, und letzteres ist nach \ref{FA.5.16} klar. {]} 

Aus (\ref{FA.5.15.3}) folgt für alle $f,g \in \L_{\K}^p(M,\varphi)$
\begin{equation} \label{FA.5.15.5}
\begin{array}{c}
\underbrace{| f + g |^q}_{\stackrel{\ref{FA.5.pq}}{\in} \L_{\R}^{p-1}(M,\varphi)} + \underbrace{| f - g |^q}_{\stackrel{\ref{FA.5.pq}}{\in} \L_{\R}^{p-1}(M,\varphi)} \le ~ \underbrace{2 \left( |f|^p + |g|^p \right)^{q-1}}_{\stackrel{\ref{FA.5.pq}}{\in} \L_{\R}^{p-1}(M,\varphi)},
\end{array}
\end{equation}
also mittels \ref{FA.5.pq} auch
\begin{equation} \label{FA.5.15.6}
\begin{array}{rcl}
\D \left\| |f + g |^q + | f - g |^q \right\|_{p-1} & = & \D \left( \int_M \left( |f + g |^q + | f - g |^q \right)^{p-1} \d \varphi \right)^{\frac{1}{p-1}} \\
& \le & \D 2 \left( \int_M \left( |f|^p + |g|^p \right) \d \varphi \right)^{q-1} = 2 \left( {\|f\|_p}^p + {\|g\|_p}^p \right)^{q-1}.
\end{array}
\end{equation}
Des weiteren gilt
\begin{equation} \label{FA.5.15.7}
\| |f + g|^q \|_{p-1} + \| |f - g|^q \|_{p-1} \le \| |f + g|^q + |f - g|^q \|_{p-1},
\end{equation}
denn wegen $p-1 \in {]}0,1{[}$ und $| f \pm g |^q \stackrel{(\ref{FA.5.15.5})}{\in} \L_{\R}^{p-1}(M,\varphi)$ ergibt \ref{FA.5.13} -- angewandt auf $p-1$, $|f+g|^q, |f-g|^q$ anstelle von $p, f ,g$ -- (\ref{FA.5.15.7}) im Falle $|f+g|^q + |f-g|^q >_{\varphi} 0$, und im Falle $|f+g|^q + |f-g|^q =_{\varphi} 0$, d.h.\ $|f \pm g|^q =_{\varphi} 0$, ist (\ref{FA.5.15.7}) trivial.
$$ {\| f \pm g \|_p}^q = \left( \int_M |f \pm g|^p \, \d \varphi \right)^{\frac{q}{p}} \stackrel{\ref{FA.5.pq}}{=} \left( \int_M |f \pm g|^{q \, (p-1)} \, \d \varphi \right)^{\frac{1}{p-1}} = \| |f \pm g|^q \|_{p-1}, $$
(\ref{FA.5.15.7}) und (\ref{FA.5.15.6}) beweisen schließlich (ii). \q

\subsection*{Dichte Teilmengen und Vollständigkeit} \addcontentsline{toc}{subsection}{Dichte Teilmengen und Vollständigkeit}

Wir haben in Bemerkung 2.) zu \ref{FA.5.7} bereits gesehen, daß für alle $p, \tilde{p} \in \R_+$ mit $p < \tilde{p}$ gilt: $\ell_{\K}^p \subset \ell_{\K}^{\tilde{p}} \subset \ell_{\K}^{\infty}$ und $\| \ldots \|_{\infty} \le \| \ldots \|_{\tilde{p}} \le \| \ldots \|_p$ auf $\ell_{\K}^p$.
Zusätzlich gilt der folgende Satz.

\begin{Satz} \label{FA.5.18}
Sei $p \in \R_+$.
\begin{itemize}
\item[(i)] Für alle $\tilde{p} \in \R_+$ mit $\tilde{p} > p$ ist $\left| \ell_{\K}^p \right|$ dicht in $\ell_{\K}^{\tilde{p}}$.
\item[(ii)] $\left| \ell_{\K}^p \right|$ ist nicht dicht in $\ell_{\K}^{\infty}$.
\end{itemize}
\end{Satz}

\textit{Beweis.} Zu (i): Seien $\tilde{p} \in \R_+$ mit $\tilde{p} > p$ und $(x_i)_{i \in \N} \in \K^n$ derart, daß $\sum_{i=0}^{\infty} |x_i|^{\tilde{p}}$ konvergiert, d.h. genau $(x_i)_{i \in \N} \in \ell_{\K}^{\tilde{p}}$.

1.\ Fall: $\tilde{p} \ge 1$. 
Zu $\varepsilon \in \R_+$ existiert dann $i_0 \in \N$ mit $\sum_{i=i_0}^{\infty} |x_i|^{\tilde{p}} < \varepsilon^{\tilde{p}}$.
Dann gilt $(x_0, \ldots, x_{i_0}, 0, \ldots) \in \ell_{\K}^p$ und
$$ \| (x_i)_{i \in \N} - (x_0, \ldots, x_{i_0}, 0, \ldots) \|_{\tilde{p}} = \left( \sum_{i=i_0}^{\infty} |x_i|^{\tilde{p}} \right)^{\frac{1}{{\tilde{p}}}} < \varepsilon. $$

2.\ Fall: $\tilde{p} \in {]}0,1{[}$.
Dann existiert zu $\varepsilon \in \R_+$ ein $i_0 \in \N$ mit $\sum_{i=i_0}^{\infty} |x_i|^{\tilde{p}} < \varepsilon$, und es gilt
$$ {d_{\tilde{p}}}^{\tilde{p}} \left( (x_i)_{i \in \N}, (x_0, \ldots, x_{i_0}, 0, \ldots) \right) = \sum_{i=i_0}^{\infty} |x_i|^{\tilde{p}} < \varepsilon. $$

Zu (ii): Sei $(x_i)_{i \in \N} \in \ell_{\K}^p$.
Dann folgt $\lim_{i \to \infty} x_i = 0$, also 
$$ \| \underbrace{(1,1,\ldots)}_{\in \ell_{\K}^{\infty}} - (x_i)_{i \in \N} \|_{\infty} = \sup \{ |1 - x_i| \, | \, i \in \N \} \ge 1. $$
Damit ist auch (ii) bewiesen. \q
\A
Für jedes $p \in {]}0, \infty{[}$ gilt $\ell_{\K}^p = L^p_{\K}(\N,\varphi_m)$, wobei $\varphi_m$ das Quadermaß zur diskreten Massenverteilung $m := 1_{\N} \: \N \to \R_+$ sei.
$\N$ ist eine $\varphi_m$-meßbare Teilmenge von $\R$ mit $\varphi_m(\N) = \infty$.
Für ein beliebiges Quadermaß und eine meßbare Teilmenge von $\R^n$, deren Maß endlich ist, erhält man ein zum letzten Satz ,,umgekehrtes`` Resultat.

\begin{Satz} \label{FA.5.19} $\,$

\noindent \textbf{Vor.:} Es seien $\varphi$ ein Quadermaß auf $\R^n$ und $M$ eine $\varphi$-meßbare Teilmenge von $\R^n$ derart, daß $\varphi(M) < \infty$ gilt.
Ferner seien $p, \tilde{p} \in {]}0, \infty{]}$ mit $p < \tilde{p}$.

\noindent \textbf{Beh.:}
\begin{itemize}
\item[(i)] $\L_{\K}^{\tilde{p}}(M,\varphi) \subset \L_{\K}^p(M,\varphi)$ und
\begin{equation} \label{FA.5.19.S}
\forall_{f \in \L_{\K}^{\tilde{p}}(M,\varphi)} \, \|f\|_p \le \varphi(M)^{\frac{1}{p} - \frac{1}{\tilde{p}}} \, \|f\|_{\tilde{p}},
\end{equation}
wobei $\frac{1}{\infty}$ wieder als $0$ zu lesen ist.
\item[(ii)] $\big| L_{\K}^{\tilde{p}}(M,\varphi) \big|$ ist dicht in $L_{\K}^p(M,\varphi)$.
\end{itemize}
\end{Satz}

\textit{Beweis.} Zu (i): Sei $f \in \L_{\K}^{\tilde{p}}(M,\varphi)$.

1.\ Fall: $\tilde{p} < \infty$.
Wir setzen $r := \frac{\tilde{p}}{p} \in {]}1, \infty{[}$ und wählen $s \in {]}1, \infty {[}$ mit $\frac{1}{r} + \frac{1}{s} = 1$. 
Dann gilt $\tilde{f} := |f|^p \in \L_{\R}^r(M,\varphi)$ sowie $\tilde{g} := 1_M \in \L_{\R}^s(M,\varphi)$ -- beachte, daß $M$ als $\varphi$-meßbare Menge mit endlichem $\varphi$-Maß eine $\varphi$-integrierbare Menge ist.
Somit folgt aus der Hölderschen Ungleichung \ref{FA.5.8} -- angewandt auf $r, s, \tilde{f}, \tilde{g}$ anstelle von $p, q, f, g$ --, daß gilt $|f|^p \in \L_{\R}(M,\varphi)$, also auch $f \in \L_{\K}^p(M,\varphi)$, und
\begin{gather*}
\| |f|^p \|_1 \le \| |f|^p \|_r \, \| 1_M \|_s, \\
\int_M |f|^p \, \d \varphi \le \left( \int_M |f|^p \, \d \varphi \right)^{\frac{p}{\tilde{p}}} \varphi(M)^{\frac{1}{s}},
\end{gather*}
also gilt (\ref{FA.5.19.S}) wegen $\frac{1}{p \, s} = \frac{r-1}{p \, r} = \frac{1}{p} - \frac{1}{\tilde{p}}$.

2.\ Fall: $\tilde{p} = \infty$.
Dann gilt analog zum 1.\ Fall $\tilde{f} := |f|^p \in \L_{\R}^{\infty}(M,\varphi)$ sowie $\tilde{g} := 1_M \in \L_{\R}^1(M,\varphi)$, und aus der Hölderschen Ungleichung folgt $|f|^p \in \L_{\R}(M,\varphi)$, $f \in \L_{\K}^p(M,\varphi)$ und 
$$ \| |f|^p \|_1 \le \| |f|^p \|_{\infty} \, \| 1_M \|_1 = \varphi(M) \, {\|f\|_{\infty}}^p. $$
Diese Ungleichung ist äquivalent zu (\ref{FA.5.19.S}).

Zu (ii): Sei $f \in \L_{\K}^p(M,\varphi) \subset \mathcal{M}_{\K}(M,\varphi)$.
Wir zeigen, daß eine Folge $(f_i)_{i \in \N}$ in $\L_{\K}^{\widetilde{p}}(M,\varphi)$ existiert derart, daß gilt
\begin{equation*} \label{FA.5.19.0}
\lim_{i \to \infty} \underbrace{\int_M |f - f_i|^p \, \d \varphi}_{= {\| f - f_i \|_p}^p} = 0,
\end{equation*}
d.h.\ auch
$$ \lim_{i \to \infty} \| f - f_i \|_p = 0. $$
Dann folgt die Behauptung im Falle $p \in {]}0,1{[}$ aus der vorletzten und im Falle $p \in {[}1, \infty{[}$ aus der letzten Gleichung.

1.\ Fall: $\K = \R$.
Für jedes $i \in \N$ definieren wir die ,,horizontal bei $i$ und $-i$ gestutzte`` (beschränkte) Funktion
$$ f_i := \sup ( \inf (f, i \, \chi_M ) , - i \, \chi_M ) \in \mathcal{M}_{\R}(M,\varphi), $$
beachte, daß $M$ als $\varphi$-meßbar vorausgesetzt ist, und behaupten
\begin{equation} \label{FA.5.19.1}
f_i \in \L_{\R}^{\tilde{p}}(M,\varphi) \stackrel{(i)}{\subset} \L_{\R}^p(M,\varphi).
\end{equation}
\pagebreak

{[} Zu (\ref{FA.5.19.1}): Im Falle $\tilde{p} = \infty$ ist dies klar, sei also $\tilde{p} < \infty$.
Dann folgt
\begin{gather}
|f|^{\tilde{p}} \in \mathcal{M}_{\R}(M,\varphi), \label{FA.5.19.2} \\
i^{\tilde{p}} \, \chi_M  \in \L_{\R}(M,\varphi), \label{FA.5.19.3} \\
|f_i|^{\tilde{p}} = \inf ( |f|^{\tilde{p}}, i^{\tilde{p}} \, \chi_M ) \in \L_{\R}(M,\varphi), \label{FA.5.19.4}
\end{gather}
also ergibt sich (\ref{FA.5.19.1}) aus (\ref{FA.5.19.2}) und (\ref{FA.5.19.4}).

Zu (\ref{FA.5.19.2}): Daß mit $f$ auch $|f|^{\tilde{p}}$ wegen der Stetigkeit von $| \ldots |^{\tilde{p}} \: \R \to \R$ mit $|0|^{\tilde{p}} = 0$ und \ref{FA.4.8} (iii) $\varphi$-meßbar über $M$ ist, haben wir schon mehrfach erwähnt.

Zu (\ref{FA.5.19.3}): Wie bereits im Beweis von (i) verwendet, folgt aus der $\varphi$-Meßbar\-keit von $M$ mit $\varphi(M) < \infty$ die $\varphi$-Integrierbarkeit von $M$.
Damit  ist (\ref{FA.5.19.3}) klar.

Zu (\ref{FA.5.19.4}): Aus (\ref{FA.5.19.2}), (\ref{FA.5.19.3}) folgt $|f_i|^{\tilde{p}} \in \mathcal{M}_{\R}(M,\varphi)$.
Ferner ist $|f_i|^{\tilde{p}}$ durch die nach (\ref{FA.5.19.2}) über $M$ $\varphi$-integrierbare Funktion $i^{\tilde{p}} \, \chi_M$ beschränkt.
\ref{FA.4.70} (i) ergibt daher (\ref{FA.5.19.4}). {]}

Nach Definition der $f_i$ gilt $\forall_{i \in \N} \, |f - f_i| \le |f|$, also auch
$$ \forall_{i \in \N} \, |f - f_i|^p \le |f|^p \in \L_{\R}(M,\varphi), $$
und $\lim_{i \to \infty} f_i = f$.
Hieraus, $f \in \L_{\R}^p(M,\varphi)$ sowie (\ref{FA.5.19.1}) folgt des weiteren
$$ (|f_i - f|^p)_{i \in \N} \mbox{ ist eine auf $M$ gegen Null $\varphi$-konvergente Folge in } \L_{\R}(M,\varphi). $$
Der Grenzwertsatz von \textsc{Lebesgue} ergibt nun
$$ \lim_{i \to \infty} \int_M |f - f_i|^p \, \d \varphi = 0, $$
und der Beweis ist wegen (\ref{FA.5.19.1}) erbracht.

2.\ Fall: $\K = \C$.
Die Anwendung des ersten Falles jeweils auf Real- und Imaginärteil von $f$ liefert eine Folge $(f_i)_{i \in \N}$ in $\L_{\C}^{\widetilde{p}}(M,\varphi)$ mit 
$$ \lim_{i \to \infty} \int_M |f - f_i|^p \, \d \varphi \le \lim_{i \to \infty} \int_M ({\rm Re} \, |f - f_i|)^p \, \d \varphi + \lim_{i \to \infty} \int_M ({\rm Im} \,|f - f_i|)^p \, \d \varphi = 0. $$
\q

\begin{Bem*} $\,$
\begin{itemize}
\item[1.)] Die Voraussetzung $\varphi(M) < \infty$ in \ref{FA.5.19} ist notwendig:
Z.B. ist die Funktion $x^{-1} \: \R \setminus \{0\} \to \R$ ein Element von $\L^2_{\R}({[}1,\infty{[},\mu_1) \setminus \L^1_{\R}({[}1,\infty{[},\mu_1)$.
\item[2.)] Wegen $\mu_1([0,1]) = \infty = \varphi_m(\N)$, wobei $\varphi_m$ das Quadermaß zu diskreten Massenverteilung $m:= 1_{\N}$ bezeichne, könnte man nun aufgrund von Satz \ref{FA.5.18} vermuten, daß $\L^1_{\R}({[}0,1{]},\mu_1)$ eine Teilmenge von $\L^2_{\R}({[}0,1{]},\mu_1)$ ist.
Letzteres stimmt jedoch nicht, denn die Funktion $x^{-\frac{1}{2}} \: \R_+ \to \R$ ist ein Element von $\L^1_{\R}({[}0,1{]},\mu_1) \setminus \L^2_{\R}({[}0,1{]},\mu_1)$.
\end{itemize}
\end{Bem*}

\begin{HS}[Satz von \textsc{Riesz-Fischer}] \index{Satz!von \textsc{Riesz-Fischer}} \label{FA.5.20} $\,$

\noindent \textbf{Vor.:} Seien $\varphi$ ein Quadermaß auf $\R^n$, $M$ eine Teilmenge von $\R^n$ und $p \in {]}0, \infty{]}$.

\noindent \textbf{Beh.:} $L^p_{\K}(M,\varphi)$ ist vollständig.
\end{HS}

\textit{Beweis.} 1.\ Fall: $p \in {[}1, \infty{]}$.
Es genügt offenbar zu zeigen, daß jede ,,Cauchyfolge`` in $(\L^p_{\K}(M,\varphi), \| \ldots \|_p)$ ,,konvergiert``.
Sei also $(f_k)_{k \in \N}$ ein Folge in $\L^p_{\K}(M,\varphi)$ mit
\begin{equation} \label{FA.5.20.0}
\forall_{\varepsilon \in \R_+} \exists_{m \in \N} \forall_{k,l \in \N, \, k,l \ge m} \, \| f_k - f_l \|_p < \varepsilon.
\end{equation}

1.1.\ Fall: $p = \infty$.
Wegen Bemerkung 2.) zu \ref{FA.5.1} existiert zu $k,l \in \N$ eine $\varphi$-Nullmenge $N_{k,l}$ mit
$$ \forall_{x \in M \setminus N_{k,l}} \, |f_k(x) - f_l(x)| \le \| f_k - f_l \|_{\infty}. $$
Dann ist auch die abzählbare Vereinigung
$$ N := \bigcup_{k,l \in \N} N_{k,l} $$ 
eine $\varphi$-Nullmenge und $(f_k)_{k \in \N}$ bzgl.\ der Supremumsnorm eine Cauchyfolge in $\mathcal{B}(M \setminus N, \K)$, die nach \ref{FA.3.2} (iii) im Falle $M \setminus N \ne \emptyset$ bzgl.\ der Supremumsnorm gegen eine Funktion $\tilde{f} \in \mathcal{B}(M \setminus N, \K)$ konvergiert.
Da $(f_k)_{k \in \N}$ auch eine Folge in $\mathcal{M}_{\K}(M,\varphi)$ ist, impliziert \ref{FA.4.69.viii} daher, daß die Folge $(f_k)_{k \in \N}$ in $\L_{\K}^{\infty}(M,\varphi)$ bzgl.\ $\| \ldots \|_{\infty}$ gegen $f \in \L_{\K}^{\infty}(M,\varphi)$, definiert durch
$$ \forall_{x \in M} \, f(x) := \left\{ \begin{array}{cl} \tilde{f}(x), & \mbox{falls } x \in M \setminus N, \\ 0, & \mbox{falls } x \in N, \end{array} \right. $$
,,konvergiert``.

1.2.\ Fall: $p \in {[}1, \infty{[}$.
Wegen (\ref{FA.5.20.0}) existiert eine Teilfolge $(f_{i_k})_{k \in \N}$ von $(f_k)_{k \in \N}$ mit
\begin{equation} \label{FA.5.20.1}
\forall_{k,l \in \N \, l \ge i_k} \, \| f_{i_k} - f_l \|_p < \frac{1}{2^k}.
\end{equation}
Wir definieren eine Folge $(g_k)_{k \in \N}$ in $\L^p_{\K}(M,\varphi)$ durch
\begin{equation} \label{FA.5.20.2}
\forall_{k \in \N} \, g_k := f_{i_{k+1}} - f_{i_k} \in \L^p_{\K}(M,\varphi).
\end{equation}
Dann ist $((\sum_{j=0}^k |g_j|)^p)_{k \in \N}$ eine nicht-negative Folge in $\L_{\R}(M,\varphi)$ derart, daß für jedes $k \in \N$ gilt
\begin{gather*}
\left( \sum_{j=0}^k |g_j| \right)^p \le_{\varphi} \left( \sum_{j=0}^{k+1} |g_j| \right)^p, \\
\int_M \left( \sum_{j=0}^k |g_j| \right)^p \, \d \varphi = \left( \left\| \sum_{j=0}^k |g_j| \right\|_p \right)^p \le \left( \sum_{j=0}^k \|g_j\|_p \right)^p \stackrel{(\ref{FA.5.20.2}), (\ref{FA.5.20.1})}{<} 2^p.
\end{gather*}
Aus dem Grenzwertsatz von \textsc{Levi} \ref{FA.4.25} folgt daher, daß $((\sum_{j=0}^k |g_j|)^p)_{k \in \N}$ -- und damit auch $(\sum_{j=0}^k g_j)_{k \in \N}$ -- auf $M$ $\varphi$-konvergent ist.
Somit impliziert
\begin{equation*} 
\forall_{k \in \N} \, f_{i_k} \stackrel{(\ref{FA.5.20.2})}{=_{\varphi}} f_{i_0} + \sum_{j=0}^{k-1} g_j
\end{equation*}
die $\varphi$-Konvergenz von $\left( f_{i_k} \right)_{k \in \N}$ auf $M$, d.h.\ es existiert eine auf einer Teilmenge von $\R^n$ definierte Funktion $f$ mit
\begin{equation} \label{FA.5.20.3}
f =_{\varphi} \lim_{k \to \infty} f_{i_k} \stackrel{\ref{FA.4.69.viii}}{\in} \mathcal{M}_{\K}(M,\varphi).
\end{equation}
Die Folge $(|f_{i_k}|^p)_{k \in \N}$ in $\L_{\R}(M,\varphi)$ ist folglich auf $M$ $\varphi$-konvergent gegen $|f|^p$, und die Integralfolge $(\int_M |f_{i_k}|^p \, \d \varphi)_{k \in \N}$ ist nach oben beschränkt, denn nach (\ref{FA.5.20.0}) existiert $m \in \N$ derart, daß gilt
$$ \forall_{k \in \N, k \ge m} \, \int_M |f_{i_k}|^p \, \d \varphi = {\|f_{i_m}\|_p}^p \le {\|f_{i_m}\|_p}^p + {\| f_{i_k} - f_{i_m}\|_p}^p < {\|f_{i_m}\|_p}^p + 1. $$
Also ergibt das Lemma von \textsc{Fatou} \ref{FA.4.33} die $\varphi$-Integrierbarkeit von $|f|^p$ über $M$, d.h.\ wegen (\ref{FA.5.20.3})
$$ f \in \L^p_{\K}(M,\varphi). $$ 

Zu zeigen bleibt 
$$ \lim_{k \to \infty} \|f_k - f\|_p = 0. $$
Hierzu seien $k,l \in \N$ mit $l \ge k$, also auch mit $i_l \ge i_k \ge k$.
Dann gilt
$$ \| f_k - f_{i_l} \|_p \stackrel{(\ref{FA.5.20.1})}{<} \frac{1}{2^k}, $$ 
also
\begin{equation} \label{FA.5.20.4}
{\| f_k - f_{i_l} \|_p}^p < \frac{1}{2^{k p}}
\end{equation}
und somit
$$ \sup \left\{ \int\limits_M |f_k - f_{i_l}|^p \, \d \mu \, | \, l \in \N \right\} < \infty. $$
Hieraus, der $\varphi$-Konvergenz der Folge $\left( |f_l - f_{i_l}|^p \right)_{l \in \N}$ auf $M$ gegen $|f_k - f|^p$ und dem Lemma von \textsc{Fatou} \ref{FA.4.33} folgt die $\varphi$-Integrierbarkeit von $|f_k - f|^p$ über $M$ sowie
$$ \int_M |f_k - f|^p \, \d \varphi \, \le \, \liminf_{l \to \infty} \int_M |f_k - f_{i_l}|^p \, \d \varphi = \liminf_{l \to \infty} {\|f_k - f_{i_l}\|}^p \stackrel{(\ref{FA.5.20.4})}{\le} \frac{1}{2^{k p}}, $$
also auch $\lim_{k \to \infty} \|f_k - f\|_p = 0$.

2.\ Fall: $p \in {]}0,1{[}$.
Der Leser adaptiere den Beweis des 1.2.\ Falles, um als Übung zu zeigen, daß zu jeder Folge $(f_k)_{k \in \N}$ in $\L^p_{\K}(M,\varphi)$ mit
\begin{equation*} \label{FA.5.20.5}
\forall_{\varepsilon \in \R_+} \exists_{m \in \N} \forall_{k,l \in \N, \, k,l \ge m} \, {\| f_k - f_l \|_p}^p < \varepsilon
\end{equation*}
eine Funktion $f \in \L^p_{\K}(M,\varphi)$ mit $\lim_{k \to \infty} {\|f_k - f\|_p}^p = 0$ existiert, womit der Satz dann offenbar bewiesen ist.
\q

\begin{Bem} \label{FA.5.20.B}
Der soeben geführte Beweis zeigt, daß für alle Quadermaße $\varphi$ auf $\R^n$ und alle Teilmengen $M$ von $\R^n$ sowie alle $p \in \R_+$ gilt:

\textit{Jede Cauchyfolge in $L^p(M,\varphi)$ besitzt eine Teilfolge, die $\varphi$-fast überall punktweise gegen ein Element von $L^p(M,\varphi)$ konvergiert.}

Wie das folgende Beispiel zeigt, konvergiert die gesamte Folge u.U.\ in keinem Punkt.
Betrachte das Volumenmaß $\mu_1$ auf $\R$ und numeriere die Intervalle
$$ I_{i,j} := \left[ \frac{i}{j+1}, \frac{i+1}{j+1} \right], ~~ \mbox{ wobei } i,j \in \N \mbox{ mit } i \le j, $$
so zu einer Folge $(I_k)_{k \in \N}$, daß gilt
$$ (I_k)_{k \in \N} = ( I_{0,0}, I_{0,1}, I_{1,1}, I_{0,2}, I_{1,2}, I_{2,2}, \ldots ). $$
Dann folgt
$$ \forall_{k \in \N} \, \chi_{I_k} \in \L_{\R}^p([0,1],\mu_1), ~~ \lim_{k \to \infty} \| \chi_{I_k} \|_p = 0, $$
und $( \chi_{I_k} )_{k \in \N}$ konvergiert in keinem Punkt von $[0,1]$, da die Folge dort sowohl $0$ als auch $1$ unendlich oft als Wert annimmt.
\end{Bem}

Wir verallgemeinern als nächstes \ref{FA.4.37} auf den Fall $p \in \R_+$ anstelle von $p=1$.

\begin{Satz} \label{FA.5.21} $\,$

\noindent \textbf{Vor.:} Seien $\varphi$ ein Quadermaß auf $\R^n$, $J \in \mathfrak{I}_n$ ein nicht-leeres Intervall von $\R^n$, $p \in \R_+$ und $f \in \mathcal{L}_{\K}^p(J,\varphi)$.

\noindent \textbf{Beh.:} Zu jedem $\varepsilon \in \R_+$ existieren $T,S \in \mathcal{S}(J)$ mit $\D \int_J |f - (T + \i \, S)|^p ~ \d \varphi < \varepsilon$.
\end{Satz}

\begin{Zusatz}
Im Falle $\K = \R$ kann $S = 0$ gewählt werden, und im Falle $\K = \C$ existieren $k \in \N_+$, $\zeta_1, \ldots, \zeta_k \in \C$ und paarweise disjunkte Quader $Q_1, \ldots, Q_k$ mit $T + \i \, S = \sum_{i=1}^k \zeta_i \, \chi_{Q_i}^J$.
\end{Zusatz}

\textit{Beweis.} Wegen \ref{FA.5.2} (ii), falls $p \ge 1$, bzw.\ \ref{FA.5.14} (ii), falls $p<1$, genügt es offenbar, den Fall $\K = \R$ zu betrachten.
Der Zusatz ist dann klar, beachte auch \ref{FA.4.8} (ii).

1.\ Fall: $f \in \mathcal{L}_{\R}^p(J,\varphi)$ ist beschränkt und besitzt einen kompakten Träger.
Dann existiert $C \in \R_+$ mit
\begin{equation} \label{FA.5.21.1.1.S}
|f| \le C ~~ \mbox{ und } ~~ {\rm Tr}(f) \subset [-C,C]^n.
\end{equation}
Da $f$ insbesondere $\varphi$-meßbar über $J$ ist, existiert offenbar eine Folge $(T_i)_{i \in \N}$ in $\mathcal{S}(J)$ mit
\begin{equation} \label{FA.5.21.0}
\lim_{i \to \infty} T_i =_{\varphi} f \mbox{ auf } J,
\end{equation}
und wir können ohne Einschränkung
\begin{gather} 
\forall_{i \in \N} \, \left( | T_i | \le C \, \wedge \, {\rm Tr}(T_i) \subset [-C,C]^n \right) \label{FA.5.21.1} 
\end{gather}
annehmen.
Es gilt
\begin{gather*} 
\forall_{i \in \N} \, f_i := |T_i - f|^p \in \L_{\R}(J,\varphi), \label{FA.5.21.2.5} \\
\lim_{i \to \infty} \, f_i \stackrel{(\ref{FA.5.21.0})}=_{\varphi} 0 \mbox{ auf } J, \label{FA.5.21.3} \\
\forall_{i \in \N} \, |f_i| \le \left( |T_i| + |f| \right)^p \stackrel{(\ref{FA.5.21.1.1.S}),(\ref{FA.5.21.1})}{\le} (2 \, C)^p \, \chi_{[-C,C]} \in \L(J,\varphi). \label{FA.5.21.4}
\end{gather*}
Aus dem Grenzwertsatz von \textsc{Lebesgue} folgt daher
$$ \lim_{i \to \infty} \int_J f_i \, \d \varphi = 0, $$
und dies ergibt die Behauptung im 1.\ Falle.

2.\ Fall: $f \in \mathcal{L}_{\R}^p(J,\varphi)$ beliebig.
Wegen des 1. Falles genügt es, zu zeigen, daß es eine Folge $(f_i)_{i \in \N}$ beschränkter Funktionen mit kompaktem Träger in $\L_{\R}^p(J,\varphi)$ mit
\begin{equation} \label{FA.5.21.5}
\lim_{i \to \infty} \int_J |f - f_i|^p \, \d \varphi = 0
\end{equation}
gibt.

Beweis hiervon: Für jedes $i \in \N$ definieren wir die ,,horizontal und vertikal bei $i$ und $-i$ gestutzte`` (beschränkte) Funktion mit kompaktem Träger
$$ f_i := \sup( \inf(f, i), -i ) \, \chi_{[-i,i]^n} \in \mathcal{M}_{\R}(J,\varphi), $$
und es gilt sogar
$$ f_i \in \L^p_{\R}(J,\varphi), $$
beachte $|f_i|^p \le_{\varphi} i^p \, \chi_{[-i,i]^n} \in \L_{\R}(J,\varphi)$ sowie \ref{FA.4.70} (i).
Es folgt
\begin{gather*}
| f - f_i |^p \in \L_{\R}(J,\varphi), \\
| f - f_i |^p \le |f|^p \in \L_{\R}(J,\varphi),
\end{gather*}
und der Grenzwertsatz von \textsc{Lebesgue} zusammen mit der Definition von $f_i$ ergibt (\ref{FA.5.21.5}).
\q

\begin{Bem*}
Es gilt $1_{\R^n} \in \L^{\infty}_{\R}(\R^n,\mu_n)$ und für jedes $T \in \mathcal{S}(\R^n)$
$$ \|1_{\R^n} - T \|_{\infty} \ge 1. $$
\end{Bem*}

\begin{Kor} \label{FA.5.22}
Seien $\varphi$ ein Quadermaß auf $\R^n$, $J \in \mathfrak{I}_n$ ein nicht-leeres Intervall von $\R^n$ und $p \in \R_+$.

Dann ist $L^p_{\K}(J,\varphi)$ separabel, d.h.\ per definitionem, daß eine höchstens abzählbare dichte Teilmenge von $L^p_{\K}(J,\varphi)$ existiert. 
\end{Kor}

\textit{Beweis.} Die Restklassen der Funktionen $\sum_{i=1}^k \zeta_i \, \chi_{Q_i}^J$ wie im Zusatz mit rationalen Koeffizienten im Falle $\K=\R$ bzw.\ $(\Q + \i \, \Q)$-wertigen Koeffizienten im Falle $\K=\C$ und $\Q^n$-wertigen ,,Eckpunkten`` der Quader leisten offenbar das Gewünschte.
\q

\begin{Satz} \label{FA.5.23} $\,$

\noindent \textbf{Vor.:} Seien $p \in \R_+$ und $f \in \mathcal{L}_{\K}^p(\R^n,\mu_n)$.

\noindent \textbf{Beh.:} Zu jedem $\varepsilon \in \R_+$ existiert $g \in \mathcal{C}_c(\R^n,\K)$ mit $\D \int_{\R^n} |f - g|^p \, \d \mu_n < \varepsilon$.
\end{Satz}

\textit{Beweis.} Nach \ref{FA.5.21} und \ref{FA.5.2} (ii), falls $p \ge 1$, bzw.\ \ref{FA.5.14} (ii), falls $p<1$, genügt es zu zeigen, daß im Falle $\K=\R$ zu jeder Treppenfunktion $T \in \mathcal{S}(\R)$ und jeder Zahl $\varepsilon \in \R_+$ eine stetige Funktion $g \: \R^n \to \R$ mit kompaktem Träger mit der Eigenschaft 
$$ \int_{\R^n} |T - g|^p \, \d \mu_n < \varepsilon $$
existiert.
Wiederum wegen \ref{FA.5.2} (ii), falls $p \ge 1$, bzw.\ \ref{FA.5.14} (ii), falls $p<1$, können wir ohne Beschränkung der Allgemeinheit $T = \chi_{Q}$ mit $Q \in \mathfrak{Q}_n$ annehmen, und der Fall $Q^{\circ} = \emptyset$ ist trivial, da dann $\chi_Q =_{\mu_n} 0 \in \mathcal{C}_c(\R^n,\R)$ gilt.

Seien also $\varepsilon \in \R_+$ und $Q \in \mathfrak{Q}_n$ mit $Q^{\circ} \ne \emptyset$.
Wegen $\overline{Q} \in \mathfrak{Q}_n$, der Regularität von $\mu_n$ und  $\mu_n(\overline{Q}) = \mu_n(Q)$ existiert ein in $\R^n$ offener Quader $Q' \in \mathfrak{Q}_n$ mit $\overline{Q} \subset Q'$ und
\begin{equation} \label{FA.5.23.1}
\mu_n(Q') \le \mu_n(Q) + \frac{\varepsilon}{2}.\footnote{Hier geht entscheidend ein, das Quadermaß $\mu_n$ (und nicht etwa ein beliebiges Quadermaß auf $\R^n$) zu betrachten.}
\end{equation}
Wir beweisen unten, daß es eine stetige Funktion $g \: \R^n \to \R$ mit
\begin{equation} \label{FA.5.23.2}
g|_Q = 1 \, \wedge \, g|_{\R^n \setminus Q'} = 0 \, \wedge \, 0 \le g - \chi_Q \le 1,
\end{equation}
gibt.
Offensichtlich gilt $g \in \mathcal{C}_c(\R^n,\R)$.
Es folgt
$$ 0 \le g - \chi_Q \le \chi_{Q'} - \chi_Q, $$
also wegen $(\chi_{Q'} - \chi_Q)(\R) = \{0,1\}$
\begin{gather*}
|g - \chi_Q|^p \le \chi_{Q'} - \chi_Q, \\
\int_{\R^n} |g - \chi_Q|^p \, \d \mu_n \le \mu_n(Q') - \mu_n(Q) \stackrel{(\ref{FA.5.23.1})}{<} \varepsilon.
\end{gather*}

Zu zeigen bleibt die Existenz einer stetigen Funktion $g \: \R^n \to \R$ mit (\ref{FA.5.23.2}):
Wegen $Q^{\circ} \ne \emptyset$ gibt es eine $\R$-Vektorraum-Affinität $h \: \R^n \to \R^n$ (d.h.\ per definitionem $h = a + \tilde{h}$, wobei $a \in \R^n$ und $\tilde{h} \in {\rm Hom}_{\R}(\R^n,\R^n) = \L_{\R}(\R^n,\R^n)$ ein $\R$-Vektorraum-Isomorphismus sind) derart, daß 
$$ h(\overline{Q}) = {[} -1,1 {]}^n =  \{ x \in \R^n \, | \, \|x\|_{\infty} \le 1 \} $$
gilt.
Da $h$ insbesondere ein Homöomorphismus ist, ist $h(Q')$ eine in $\R^n$ offene Obermenge von ${[}-1,1{]}^n$, und es existiert folglich eine Zahl $\delta \in {]}1, \infty{[}$ mit
$$ \{ x \in \R^n \, | \, \|x\|_{\infty} < \delta \} = {]} - \delta, \delta{[}^n \subset h(Q'). $$
Wir definieren nun eine stetige Funktion $\psi \: {[}0, \infty{[} \to \R$ durch
$$ \forall_{t \in {[}0, \infty{[}} \, \psi(t) := \left\{ \begin{array}{cl} 1, & \mbox{falls } 0 \le t \le 1, \\ 1 - \frac{t-1}{\delta - 1}, & \mbox{falls } 1 < t < \delta, \\ 0, & \mbox{falls } t \ge \delta. \end{array} \right. $$
Dann leistet $g \: \R^n \to \R$, gegeben durch
$$ \forall_{x \in \R^n} \, g(x) := \psi(\|h(x)\|_{\infty}), $$
offenbar das Gewünschte.
\q

\begin{Bem*} $\,$
\begin{itemize}
\item[1.)] Betrachte das Volumenmaß $\mu_n$ auf $\R^n$, und sei $p \in \R_+$.
Dann definiert ${d_p}^p$ eine Metrik und $\| \ldots \|_p$ für $p \in {[}1, \infty{[}$ eine Norm auf $\mathcal{C}_c(\R^n,\K)$ bzgl.\ derer $\mathcal{C}_c(\R^n,\K)$ aber nicht vollständig ist.
Aufgrund des letzten Satzes und des Satzes von \textsc{Riesz-Fischer} kann $L^p_{\K}(\R^n,\mu_n)$ im Falle $p \in {]} 0, 1 {[}$ als Vervollständigung von $(\mathcal{C}_c(\R^n,\K),{d_p}^p)$ und im Falle $p \in {[}1, \infty{[}$ von $(\mathcal{C}_c(\R^n,\K),\| \ldots \|_p)$ aufgefaßt werden.
\item[2.)] Aus $1_{\R^n} \in \L_{\K}^{\infty}(\R^n,\mu_n) \setminus \mathcal{C}_c(\R^n,\K)$ folgt leicht, daß Satz \ref{FA.5.23} i.a.\ falsch ist, wenn man $p = \infty$ zuläßt.
\end{itemize}
\end{Bem*}

\begin{Kor} \label{FA.5.23.K} $\,$

\noindent \textbf{Vor.:} Seien $p \in \R_+$, $K$ eine kompakte Teilmenge von $\R^n$ und $f \in \mathcal{L}_{\K}^p(K,\mu_n)$.

\noindent \textbf{Beh.:} Zu jedem $\varepsilon \in \R_+$ existiert $g \in \mathcal{C}(K,\K)$ mit $\D \int_{K} |f - g|^p \, \d \mu_n < \varepsilon$. \q
\end{Kor}

\begin{Bem*}
Das Korollar ist für $p = \infty$ i.a.\ falsch.
\end{Bem*}

Zum Abschluß dieses Kapitels geben wir im Falle $p \in {[}1, \infty{]}$ eine Teilmenge des Dualraumes eines Lebesgueschen Raumes an.

\begin{Satz} \label{FA.5.17}
Es seien $\varphi$ ein Quadermaß auf $\R^n$, $M$ eine Teilmenge von $\R^n$ und $p,q \in {[} 1, \infty {]}$ mit $\frac{1}{p} + \frac{1}{q} = 1$, wobei $\frac{1}{\infty}$ wieder als $0$ zu lesen ist.
\begin{itemize}
\item[(i)] Sei $\mathbf{f} \in L^q_{\K}(M,\varphi)$.
Dann definiert
\begin{equation} \label{FA.5.17.0}
\forall_{\mathbf{g} \in L^p_{\K}(M,\varphi)} \, \Lambda_{\mathbf{f}} \, \mathbf{g} := \int_M f \, g \, \d \varphi, \mbox{ wobei $f \in \mathbf{f}$ und $g \in \mathbf{g}$ beliebig seien},
\end{equation}
ein stetiges Funktional $\Lambda_{\mathbf{f}}$ von $L^p_{\K}(M,\varphi)$, d.h.\ $\Lambda_{\mathbf{f}} \in L^p_{\K}(M,\varphi)'$, und es gilt
\begin{equation} \label{FA.5.17.4}
\| \Lambda_{\mathbf{f}} \| = \| \mathbf{f} \|_q.
\end{equation}
\item[(ii)] Durch
\begin{equation*} \label{FA.5.17.2}
\Lambda \: L^q_{\K}(M,\varphi) \longrightarrow L^p_{\K}(M,\varphi)', ~~ \mathbf{f} \longmapsto \Lambda_{\mathbf{f}},
\end{equation*}
ist ein isometrischer $\K$-Vektorraum-Homomorphismus gegeben.
$\Lambda$ ist also insbesondere injektiv und stetig.
\end{itemize}
\end{Satz}

\textit{Beweis.} Zu (i): Aus der Hölderschen Ungleichung \ref{FA.5.8} folgt offenbar, daß (\ref{FA.5.17.0}) wohldefiniert ist und
\begin{equation} \label{FA.5.17.3}
\forall_{\mathbf{g} \in L^p_{\K}(M,\varphi)} \, | \Lambda_{\mathbf{f}} \, \mathbf{g} | \le \| \Lambda_{\mathbf{f}} \, \mathbf{g} \|_1 \le \|\mathbf{f}\|_q \, \|\mathbf{g}\|_p,
\end{equation}
also ist die $\K$-lineare Funktion $\Lambda_{\mathbf{f}} \: L^p_{\K}(M,\varphi) \to \K$ nach \ref{FA.2.7} ,,(iv) $\Rightarrow$ (ii)`` stetig.

Zu (\ref{FA.5.17.4}): Sei $f \in \mathbf{f}$.
Im Falle $f =_{\varphi} 0$ ist die Behauptung trivial.
Daher können wir ohne Beschränkung der Allgemeinheit 
$$ \forall_{x \in M} \, f(x) \ne 0 $$
annehmen -- insbes.\ ist $f$ auf $M$ definiert --, also ist nach \ref{FA.4.69} (iv), (v) die \emph{,,Vorzeichenfunktion von $f$``}
\begin{equation} \label{FA.5.17.5}
\sigma := \frac{f}{|f|} \: M \longrightarrow \K \mbox{ $\varphi$-meßbar über $M$ mit $\| \sigma \|_{\infty} = 1$.}
\end{equation}

1.\ Fall: $q = 1$. 
Dann gilt $p = \infty$, $g := \sigma \in \L^{\infty}_{\K}(M,\varphi)$, $\| g \|_{\infty} = 1$ und
$$ \left| \int_M f \, g \, \d \varphi \right| = \int_M |f| \, \d \varphi = \|f\|_1 = \|f\|_1 \, \|g\|_{\infty}. $$
Wegen (\ref{FA.5.17.3}) und der letzten Aussage in \ref{FA.2.8} (ii) ist daher (\ref{FA.5.17.4}) erfüllt.

2.\ Fall: $q \in {]}1, \infty{[}$.
Dann gilt auch $p \in {]}1, \infty{[}$ sowie $p \, (q-1) \stackrel{\ref{FA.5.pq}}{=} q$, also nach Voraussetzung $|f|^{q-1} \in \L^p_{\K}(M,\varphi)$, und daher wegen (\ref{FA.5.17.5}) und \ref{FA.4.70} (iii)
$$ g := \sigma \, |f|^{q-1} \in \L^p_{\K}(M,\varphi). $$
Außerdem folgt
\begin{gather} 
f \, g = |f|^q, \label{FA.5.17.6} \\
|f|^q = |\sigma|^p \, |f|^{(q-1) \, p} = |g|^p \label{FA.5.17.7}
\end{gather}
und somit
$$ \left| \int_M f \, g \, \d \varphi \right| \stackrel{(\ref{FA.5.17.6})}{=} \int_M |f|^q \, \d \varphi = \left( \int_M |f|^q \, \d \varphi \right)^{\frac{1}{q} + \frac{1}{p}} \stackrel{(\ref{FA.5.17.7})}{=} \|f\|_q \, \|g\|_p, $$
d.h.\ erneut wegen (\ref{FA.5.17.3}) und der letzten Aussage in \ref{FA.2.8} (ii) auch (\ref{FA.5.17.4}).

3.\ Fall: $q = \infty$.
Dann gilt $p = 1$.
Wegen $f \ne_{\varphi} 0$ folgt aus der Definition von $\|f\|_{\infty}$ für jedes $\varepsilon \in {]} 0, \|f\|_{\infty} {[} \ne \emptyset$, daß
\begin{equation*} \label{FA.5.17.8}
M_{\varepsilon} := \{ x \in M \, | \, |f(x)| > \| f \|_{\infty} - \varepsilon \} = \overline{|f|}^1 \left( {]} \| f \|_{\infty} - \varepsilon, \infty {[} \right) \in \B(M,\varphi)
\end{equation*}
keine $\varphi$-Nullmenge ist, und nach \ref{FA.4.M4}, (\ref{FA.4.M4D.S}) existiert $i \in \N$ derart, daß für die Teilmenge
\begin{equation*} \label{FA.5.17.9}
\widetilde{M}_{\varepsilon} := M_{\varepsilon} \cap {]}-i,i{[} \in \B(M,\varphi)
\end{equation*}
von $M$ gilt 
$$ 0 < \varphi(\widetilde{M}_{\varepsilon}) < \infty, $$
d.h.\ genau: $\widetilde{M}_{\varepsilon}$ ist $\varphi$-integrierbare Nicht-Nullmenge.
Wir setzen nun
$$ g_{\varepsilon} := \chi_{\widetilde{M}_{\varepsilon}} \, \sigma \in \L^1_{\K}(M,\varphi), $$
beachte, daß $\chi_M \, \widehat{g_{\varepsilon}} = \chi_{\widetilde{M}_{\varepsilon}} \, \hat{\sigma}$ u.a.\ nach (\ref{FA.5.17.5}) und \ref{FA.4.70} (iii) eine über $\R^n$ $\varphi$-in\-te\-grier\-bare Funktion ist.
(\ref{FA.5.17.5}) ergibt weiterhin
$$ \| g_{\varepsilon} \|_1 = \int_{\widetilde{M}_{\varepsilon}} |\sigma| \, \d \varphi = \varphi \left( \widetilde{M}_{\varepsilon} \right) \ne 0 $$
und 
$$ \left| \int_M f \, g_{\varepsilon} \, \d \varphi \right| = \int_M |f| \, |g_{\varepsilon}| \, \d \varphi = \int_{\widetilde{M}_{\varepsilon}} |f| \, |\sigma| \, \d \varphi \ge \left( \|f\|_{\infty} - \varepsilon \right) \, \| g_{\varepsilon} \|_1 , $$
also folgt aus \ref{FA.2.8} (i)
$$ \| \Lambda_{\mathbf{f}} \| \ge \|f\|_{\infty} - \varepsilon. $$
Die Beliebigkeit von $\varepsilon \in {]} 0, \|f\|_{\infty} {[} \ne \emptyset$ und (\ref{FA.5.17.3}) ergeben nun (\ref{FA.5.17.4}).

Zu (ii): Da die $\K$-Linearität der Abbildung 
$$ \Lambda \: L^q_{\K}(M,\varphi) \longrightarrow L^p_{\K}(M,\varphi)' $$
klar ist, folgt (ii) aus (i).
\q

\begin{Bem*}
Wir werden im siebenten Kapitel den \emph{Rieszschen Darstellungssatz für Lebesguesche Räume} beweisen, der besagt, daß $\Lambda \: L^q_{\K}(M,\varphi) \to L^p_{\K}(M,\varphi)'$ im Falle einer $\varphi$-meßbaren Teilmenge $M$ von $\R^n$ und $p \in {[} 1, \infty {[}$ sogar eine $\K$-Vek\-tor\-raum-Isometrie ist.
\end{Bem*}

\subsection*{Übungsaufgaben}

\begin{UA} $\,$
\begin{itemize}
\item[(i)] Beweise die Youngsche Ungleichung (\ref{FA.5.8.1}) mittels des in loc.\ cit.\ angegebenen Tips, welcher ebenfalls zu zeigen ist.
\item[(ii)] Sei $p \in {]}0,1{[}$.
Zeige, daß die beim Nachweis von (\ref{FA.5.14.1}) verwendete Funktion
$$ {[}0, \infty{[} \longrightarrow \R, ~~ t \longmapsto 1 + t^p - (1 + t)^p, $$
streng monoton fallend ist, und daß sie Null auf Null abbildet.
\item[(iii)] Zeige, daß die beim Nachweis von (\ref{FA.5.14.2}) verwendete Funktion
$$ {[}0, \infty{[} \longrightarrow \R, ~~ t \longmapsto \frac{1 + t^{\tilde{p}}}{(1 + t)^{\tilde{p}}}, $$ 
für $\tilde{p} \in {]} 1, \infty {[}$ in Eins ein absolutes Minimum mit Funktionswert $2^{\tilde{p}-1}$ besitzt.
\item[(iv)] Zeige, daß die im Beweis von \ref{FA.5.16} verwendete Funktion
$$ \R_+ \longrightarrow \R, ~~ t \longmapsto \frac{1 - \tau^{\frac{t}{p-1}}}{\frac{t}{p-1}}, $$
für $\tau \in {]}0,1{[}$ und $p \in {]}1, \infty{[}$ streng monoton fallend ist.
\item[(v)] Seien $\zeta, \xi \in \C$ und $p \in {[}2,\infty{[}$.
Zeige, daß die im Beweis der Clarksonschen Ungleichung \ref{FA.5.15} verwendete Funktion 
\begin{eqnarray*}
\R & \longrightarrow & \R \\
t & \longmapsto & ( |\zeta|^2 + |\xi|^2 + 2 \, |\zeta| \, |\xi| \, \cos(t))^{\frac{p}{2}} - ( |\zeta|^2 + |\xi|^2 - 2 \, |\zeta| \, |\xi| \, \cos(t))^{\frac{p}{2}}
\end{eqnarray*}
für jedes $k \in \Z$ auf ${[} (k-1) \, \pi , k \, \pi - \frac{\pi}{2} {]}$ monoton fallend und auf ${[} k \, \pi - \frac{\pi}{2} , k \, \pi {]}$ monoton wachsend ist.
\end{itemize}
\end{UA}

\begin{UA} \label{FA.5.UA1}
Seien $\varphi$ ein Quadermaß auf $\R^n$, $M$ eine Teilmenge von $\R^n$, $p \in {]}0, \infty{]}$, $f \in \L_{\K}^p(M,\varphi)$ und $r \in \R_+$.

Zeige $|f|^r \in \L_{\K}^{\frac{p}{r}}(M,\varphi)$ und $\| |f|^r \|_{\frac{p}{r}} = {\| f \|_p}^r$, wobei $\frac{\infty}{r}$ als $\infty$ zu lesen ist.
\end{UA}

\begin{UA}[Verallgemeinerte Höldersche Ungleichung] \index{Ungleichung!Höldersche}
Es seien $\varphi$ ein Quadermaß auf $\R^n$, $M$ eine Teilmenge von $\R^n$, $k \in \N_+$, $p, q_1, \ldots, q_k, r \in {]} 0, \infty {]}$ mit $\frac{1}{p} + \sum_{\kappa = 1}^k \frac{1}{q_{\kappa}} = \frac{1}{r}$, wobei $\frac{1}{\infty}$ wie oben als $0$ zu lesen ist, und $f \in \L_{\K}^p(M,\varphi)$ sowie $g_{\kappa} \in \L_{\K}^{q_{\kappa}}(M,\varphi)$ für $\kappa \in \{1, \ldots, k\}$.

Beweise durch vollständige Induktion über $k \in \N_+$: $f \cdot \prod_{\kappa=1}^k g_{\kappa} \in \L^r(M,\varphi)$ und
$$ \left\| f \cdot \prod_{\kappa=1}^k g_{\kappa} \right\|_r \le \|f\|_p \cdot \prod_{\kappa=1}^k \| g_{\kappa} \|_{q_{\kappa}}. $$
\end{UA}

Tip: Verwende beim Induktionsanfang im Falle $r \in \R_+$ zunächst die vorherige Aufgabe \ref{FA.5.UA1} und sodann die Höldersche Ungleichung \ref{FA.5.8}.

\begin{UA} \label{FA.5.UA2}
Seien $\varphi$ ein Quadermaß auf $\R^n$, $M$ eine Teilmenge von $\R^n$, $p_1, p_2 \in {]}0, \infty{]}$ mit $p_1 < p_2$ und $f \in \L_{\K}^{p_1}(M,\varphi) \cap \L_{\K}^{p_2}(M,\varphi)$.

Zeige, daß dann für jedes $p \in {]}p_1,p_2{[}$ auch $f \in \L_{\K}^p(M,\varphi)$ gilt sowie 
$$ \|f\|_p \le {\|f\|_{p_1}}^{1-\tau} \, {\|f\|_{p_1}}^{\tau}, $$
wobei $\tau \in {]}0,1{[}$ die eindeutig bestimmte Zahl mit
$$ \frac{1}{p} = \frac{1 - \tau}{p_1} + \frac{\tau}{p_2} $$
sei und $\frac{\tau}{\infty}$ als $0$ zu lesen ist.
Begründe auch die Existenz und Eindeutigkeit eines solchen $\tau$.
\end{UA}

\begin{Def*}
Seien $X$ ein $\K$-Vektorraum, $C$ eine konvexe Teilmenge von $X$ und $\psi \: C \to \R$ eine Funktion.
\begin{itemize}
\item[(i)] $\psi$ heißt \emph{konvex} $: \Longleftrightarrow$ $\forall_{x,y \in C} \forall_{t \in {[}0,1{]}} \, \psi((1-t) \, x + t \, y) \le (1-t) \, \psi(x) + t\, \psi(y)$.
\item[(ii)] Gilt zusätzlich $\psi(C) \subset \R_+$, so heißt $\psi$ genau dann \emph{logarithmisch konvex} wenn $\ln \circ \psi \: C \to \R$ konvex ist.
\end{itemize}
\end{Def*} 

\begin{UA}
Seien $\varphi$ ein Quadermaß auf $\R^n$ und $M$ eine Teilmenge von $\R^n$.
Nach Aufgabe \ref{FA.5.UA2} wird für jedes $f \in \mathcal{M}_{\K}(M,\varphi)$ durch
$$ I(f) := \{ p \in {]}0,\infty{]} \, | \, f \in \L^p_{\K}(M,\varphi) \} $$
ein ggf.\ leeres Intervall von $\widehat{\R} = \R \cup \{ - \infty, \infty \}$ definiert.
\begin{itemize}
\item[(i)] Es sei weiterhin $f \in \mathcal{M}_{\K}(M,\varphi) \setminus \mathcal{N}_{\K}(M,\varphi)$ mit $I(f) \ne \emptyset$.
Dann ist $I(f)^{-1} = \{ t \in {[}0, \infty{[} \, | \, \frac{1}{t} \in I(f) \}$, wobei wir $\frac{1}{0}$ als $\infty$ lesen, ein nicht-leeres Intervall von $\R$.

Zeige, daß
$$ I(f)^{-1} \longrightarrow \R_+, ~~ t \longmapsto \|f\|_{\frac{1}{t}} $$
logarithmisch konvex und daher (!) konvex ist.
Folgere hieraus die Stetigkeit der Funktion
$$ (I(f) \setminus \{ \infty \})^{\circ} \longrightarrow \R_+, ~~ p \longmapsto \|f\|_p. $$
\item[(ii)] Sei $f \in \mathcal{M}_{\K}(M,\varphi)$ derart, daß $\{\infty\}$ eine echte Teilmenge von $I(f)$ ist.

Zeige $\D \lim_{p \to \infty} \|f\|_p = \|f\|_{\infty}$.
\end{itemize}
\end{UA}

Tip zu (ii): Für geeignete $p,p_0$ gilt $\|f\|^p \le |f|^{p_0} \, {\|f\|_{\infty}}^{p-p_0}$ und $\alpha^p \cdot \chi_N \le |f|^p$ für jedes $\alpha < \|f\|_{\infty}$, wobei $N := \{ x \in M \, | \, |f(x)| > \alpha \}$ sei.
Leite hieraus die Ungleichung $\limsup_{p \to \infty} \|f\|_p \le \|f\|_{\infty} \le \liminf_{p \to \infty} \|f_p\|$ her. 

\begin{UA}
Zeige die Binominalentwicklung \ref{FA.5.16.BE} für $z \in {]}-1,1{[}$ mittels des Taylorschen Satzes der Analysis I.
\end{UA}

\begin{UA}
Führe den Beweis des Satzes von \textsc{Riesz-Fischer} \ref{FA.5.20} im Falle $p \in {]}0,1{[}$ in allen Einzelheiten aus.
\end{UA}

\begin{UA} \label{FA.5.CcRnpnichtvollst}
Sei $p \in \R_+$.

Beweise, daß $C_c(\R^n,\K)$ nicht abgeschlossen (und somit auch nicht vollständig) in $(\L_{\K}^p(\R^n,\mu_n), \| \ldots \|_p)$ ist.
\end{UA}

Tip: Wegen Aufgabe \ref{FA.4.exp} gilt offenbar $\sqrt[p]{\e^{-(x_1^2 + \ldots + x_n^2)}} \in \L_{\K}^p(\R^n,\mu_n) \setminus C_c(\R^n,\K)$.

\begin{Bem*}
Da wir bereits wissen, daß $C_c(\R^n,\K)$ nach Satz \ref{FA.5.23} für jede positive natürliche Zahl $n \in \N_+$ sowie alle $p \in \R_+$ als dichte Teilmenge von $L_{\K}^p(\R^n,\mu_n)$ aufgefaßt werden kann und $C_c(\R^n,\K)$ trivialerweise eine echte Teilmenge von $C_0(\R^n,\K)$ ist, können wir uns aufgrund der letzten Aufgabe berechtigterweise fragen, ob $C_0(\R^n,\K)$ nicht bereits gleich $\L_{\K}^p(\R^n,\mu_n) \cap \mathcal{C}(\R^n,\K)$ ist.
Tatsächlich ist dies zumindest für $p=2$ nicht der Fall.
\textsc{Richtmyer} möchte in \cite[S.\ 85]{Richt} zwar eine unbeschränkte stetige quadratintegrierbare Funktion nennen, verschreibt sich aber, indem er in loc.\ cit.\ $x^2 \, \exp(-x^8 \, \sin^2(x))$ anstelle von $|x| \, \exp(-x^8 \, \sin^2(x))$ angibt.
Der Leser sei hierzu auf die Adresse \verb+https://math.stackexchange.com/questions/355447+ hingewiesen.
Summa summarum gibt es jedenfalls ein Gegenbeispiel.
\end{Bem*}

\begin{UA} \label{FA.5.nichtdicht}
Zeige die folgenden beiden Aussagen:
\begin{itemize}
\item[(i)] Der $\K$-Banachraum $(\mathcal{C}_0(\R^n,\K), \| \ldots \|_{\infty})$ ist keine dichte Teilmenge von $(\L_{\K}^{\infty}(\R^n,\mu_n), \| \ldots \|_{\infty})$.
\item[(ii)] Der $\R$-Banachraum $(\mathcal{C}({[}0,1{]},\R), \| \ldots \|_{\infty})$ ist keine dichte Teilmenge von $(\L^{\infty}_{\R}({[}0,1{]},\mu_1), \| \ldots \|_{\infty})$.
\end{itemize}
\end{UA}

\cleardoublepage
\section{Reflexive Räume} \label{FAna6}
\subsection*{Schwache Topologien} \addcontentsline{toc}{subsection}{Schwache Topologien}

\begin{Def} \label{FA.6.1}
Sei $X$ eine Menge. 
\begin{itemize}
\item[(i)] Sind $\mathcal{T}_1$ sowie $\mathcal{T}_2$ Topologien für $X$ mit $\mathcal{T}_1 \subset \mathcal{T}_2$, so heißt \emph{$\mathcal{T}_1$ gröber als $\mathcal{T}_2$}\index{Topologie!gröbere} und \emph{$\mathcal{T}_2$ feiner als $\mathcal{T}_1$}\index{Topologie!feinere}.

\begin{Bsp*}
Die triviale Topologie für $X$ ist die gröbste und die diskrete die feinste Topologie für $X$.
\end{Bsp*}
\item[(ii)] Sei $\mathfrak{S}$ eine Teilmenge von $\mathfrak{P}(X)$.\\
Dann erhält man aus $\mathfrak{S}$ durch Hinzunahme aller endlichen Schnitte von Elementen von $\mathfrak{S}$ und beliebig vieler Vereinigungen solcher Mengen zu $\mathfrak{S}$ offenbar eine gröbste Topologie $\boxed{\mathcal{T}(\mathfrak{S})}$ für $X$, die $\mathfrak{S}$ enthält.\footnote{Man beachte, daß der leere Schnitt von Teilmengen als der ganze Raum und die leere Vereinigung als die leere Menge definiert sind.}
$\mathcal{T}(\mathfrak{S})$ heißt \emph{die von $\mathfrak{S}$ erzeugte Topologie für $X$}\index{Topologie!erzeugte}.

\begin{Bem*}
Analog zu \ref{FA.4.M1} (ii) könnten wir auch folgende äquivalente Definition geben
$$ \mathcal{T}(\mathfrak{S}) := \bigcap_{\substack{\mathcal{T} \text{ Topologie für $X$} \\ \mathfrak{S} \subset \mathcal{T}}} \mathcal{T}. $$
Die vorher genannte Definition hat allerdings den Vorteil, daß sie eine explizite Konstruktion für $\mathcal{T}(\mathfrak{S})$ beinhaltet.
Eine solche anzugeben, ist (natürlich im Falle $\mathfrak{S} \ne \emptyset$) für die minimale $\sigma$-Algebra $\sigma(\mathfrak{S})$, die $\mathfrak{S}$ enthält, nicht möglich.
\end{Bem*}
\end{itemize}
\end{Def}

\begin{Def}[Die schwache Topologie bzgl.\ einer Menge von Abbildungen] \index{Topologie!schwache} \label{FA.6.2}
Es seien $X, I$ Mengen und für jedes $i \in I$ des weiteren $Y_i$ ein topologischer Raum sowie $f_i \: X \to Y_i$ eine Abbildung.

Dann heißt die von
$$ \left\{ \overline{f_i}^1(V_i) \, | \, i \in I \, \wedge \, V_i \in  {\rm Top}(Y_i) \right\} \subset \mathfrak{P}(X) $$
erzeugte Topologie für $X$ die \emph{schwache Topologie für $X$ bzgl.\ $\{ f_i \, | \, i \in I \}$}.
Dies ist offenbar die gröbste Topologie für $X$ bzgl.\ derer alle $f_i$, $i \in I$, stetig sind.
\end{Def}

\begin{Satz} \label{FA.6.2.S} $\,$

\noindent \textbf{Vor.:} Seien $X, I$ Mengen und für jedes $i \in I$ des weiteren $Y_i$ ein topologischer Raum sowie $f_i \: X \to Y_i$ eine Abbildung.
Ferner bezeichne $\mathcal{T}$ die schwache Topologie für $X$ bzgl.\ $\{ f_i \, | \, i \in I \}$.

\noindent \textbf{Beh.:}
\begin{itemize}
\item[(i)] Sind $\{ f_i \, | \, i \in I \}$ \emph{punktetrennend}, d.h.\ $\forall_{x, \tilde{x} \in X, \, x \ne \tilde{x}} \exists_{i_0 \in I} \, f_{i_0}(x) \ne f_{i_0}(\tilde{x})$, und $Y_i$ für jedes $i \in I$ hausdorffsch, so ist auch $(X,\mathcal{T})$ hausdorffsch.
\item[(ii)] Seien $(x_j)_{j \in \N}$ eine Folge in $X$ und $x \in X$.\\
Dann konvergiert $(x_j)_{j \in \N}$ genau dann in $(X,\mathcal{T})$ gegen $x$, wenn für jedes $i \in I$ die Folge $(f_i(x_j))_{j \in \N}$ in $Y_i$ gegen $f_i(x)$ konvergiert.
\item[(iii)] Seien $Z$ ein weiterer topologischer Raum und $g \in X^Z$.\\
Dann ist $g \: Z \to (X,\mathcal{T})$ genau dann stetig, wenn für jedes $i \in I$ die Abbildung $f_i \circ g \: Z \to Y_i$ stetig ist.
\end{itemize}
\end{Satz}

\textit{Beweis.} Zu (i): Seien $x, \tilde{x} \in X$ derart, daß $x \ne \tilde{x}$ gilt.
Dann existieren wegen der Voraussetzung in (i) ein $i_0 \in I$ mit $f_{i_0}(x) \ne f_{i_0}(\tilde{x})$ und $V \in \U(f_{i_0}(x),Y_{i_0})$ sowie $\widetilde{V} \in \U(f_{i_0}(\tilde{x}),Y_{i_0})$ mit $V \cap \widetilde{V} = \emptyset$.
Hieraus folgt $x \in U := \overline{f_{i_0}}^1(V) \in \mathcal{T}$, $\tilde{x} \in \widetilde{U} := \overline{f_{i_0}}^1(\widetilde{V}) \in \mathcal{T}$ und $U \cap \widetilde{U} = \emptyset$.

Zu (ii): ,,$\Rightarrow$`` ist wegen der Stetigkeit der $f_i$, $i \in \N$, klar, vgl.\ (\ref{FA.1.St.S}).

,,$\Leftarrow$`` Sei $U \in \U(x,(X,\mathcal{T}))$.
Zu zeigen ist, daß für alle hinreichend großen $j \in \N$ gilt: $x_j \in U$.

Die Definition von $\mathcal{T}$ impliziert offenbar die Existenz einer Zahl $k \in \N_+$, $i_1, \ldots, i_k \in I$ sowie $V_{\kappa} \in \U(f_{i_{\kappa}}(x),Y_{i_{\kappa}})$ für $\kappa \in \{1, \ldots, k\}$ mit
$$ \bigcap_{\kappa=1}^k \overline{f_{i_{\kappa}}}^1(V_{\kappa}) \subset U, $$
und nach Voraussetzung der rechten Seite gilt $x_j \in \overline{f_{i_{\kappa}}}^1(V_{\kappa})$ für alle $\kappa \in \{1, \ldots, k\}$ sowie hinreichend große $j \in \N$.

Zu (iii): ,,$\Rightarrow$`` ist wieder wegen der Stetigkeit der $f_i$, $i \in \N$, klar. 

,,$\Leftarrow$`` Sei $U \in \mathcal{T}$.
Aus der Definition von $\mathcal{T}$ folgt, daß sich $U$ als beliebige Vereinigung endlicher Schnitte von Mengen der Form
$$ \overline{f_{i_1,i_2}}^1(V_{i_1,i_2}) \mbox{ mit } i_1,i_2 \in I \mbox{ und } V_{i_1,i_2} \in {\rm Top}(Y_{i_1,i_2}) $$
schreiben läßt. Wegen 
$$ \overline{g}^1 \left( \bigcup_{i_1} \bigcap_{i_2} \overline{f_{i_1,i_2}}^1(V_{i_1,i_2}) \right) = \bigcup_{i_1} \overline{g}^1 \left( \bigcap_{i_2} \overline{f_{i_1,i_2}}^1(V_{i_1,i_2}) \right) = \bigcup_{i_1} \bigcap_{i_2} \overline{g}^1 \left( \overline{f_{i_1,i_2}}^1(V_{i_1,i_2}) \right) $$
genügt es zum Nachweis der Offenheit von $\overline{g}^1(U)$ in $Z$ zu zeigen, daß die
$$ \overline{g}^1 \left( \overline{f_{i_1,i_2}}^1(V_{i_1,i_2}) \right) = \overline{f_{i_1,i_2} \circ g}^1(V_{i_1,i_2}) $$
in $Z$ offen sind, und dies ist nach Voraussetzung der rechten Seite von (ii) klar. \q
\A
Für die Funktionalanalysis von großer Bedeutung ist der folgende Spezialfall der letzten Definition.

\begin{Def}[Die schwache Topologie eines normierten Vektorraumes] \index{Topologie!schwache} \label{FA.6.3}
Sei $X$ ein normierter $\K$-Vektorraum.

Dann heißt die schwache Topologie für $X$ bzgl.\ seines topologischen Dualraumes $X'$ die \emph{schwache Topologie für $X$}.
Für den topologischen Raum, der aus der Menge $X$ zusammen mit dieser Topologie besteht, schreiben wir $\boxed{X_s}$.\footnote{$X_s$ in unserem Sinne wird in der Literatur häufig mit $X_w$ bezeichnet.}
\end{Def}

\begin{Bem} \label{FA.6.4}
Ist $X$ ein normierter $\K$-Vektorraum, so ist die schwache Topologie für $X$ gröber als seine Normtopologie, denn bzgl.\ der Normtopologie sind alle Elemente von $X'$ stetig.
Daher wird die Normtopologie für $X$ in der Literatur gelegentlich auch die \emph{starke Topologie für $X$} genannt.
\end{Bem}

\begin{Satz} \label{FA.6.4.S} $\,$

\noindent \textbf{Vor.:} Sei $X$ ein endlich-dimensionaler normierter $\K$-Vektorraum.

\noindent \textbf{Beh.:} Die Normtopologie und die schwache Topologie für $X$ stimmen überein.

\end{Satz}

\textit{Beweis.} Wegen \ref{FA.6.4} müssen wir zeigen, daß die Normtopologie in der schwachen Topologie für $X$ enthalten ist.
Seien also $U \in {\rm Top}(X,\| \ldots\|)$ und $x_0 \in U$ beliebig.
Dann existiert eine Zahl $\varepsilon \in \R_+$ mit $U_{\varepsilon}(x_0) \subset U$.

Ohne Beschränkung der Allgemeinheit gelte $X \ne \{0\}$.
Nach Voraussetzung existieren $n \in \N_+$ und eine Basis $\{b_1, \ldots, b_n\}$ von $X$ mit $\forall_{i \in \{1, \ldots, n\}} \, \|b_i\| = 1$.
$\{b_1^*, \ldots, b_n^*\}$ bezeichne die dazu duale Basis von $X^* = X'$ -- beachte $\dim_{\K} X < \infty$.
Nun gilt
$$ V := \bigcap_{i=1}^n \overline{b_i^*}^1\left(U_{\frac{\varepsilon}{n}}(b_i^*(x_0)) \right) \in \U(x_0,X_s) $$
und wegen
$$ \forall_{x \in X} \, \| x - x_0 \| = \left\| \sum_{i=1}^n b_i^*(x-x_0) \, b_i \right\| \le \sum_{i=1}^n |b_i^*(x-x_0)| \, \|b_i\| = \sum_{i=1}^n |b_i^*(x-x_0)| $$
offenbar auch $V \subset U_{\varepsilon}(x_0) \subset U$, d.h.\ $x_0$ ist innerer Punkt von $U$ bzgl.\ der schwachen Topologie für $X$.
Aus der Beliebigkeit von $x_0 \in U$ folgt, daß $U$ in $X_s$ offen ist. \q
\A
Das folgende Lemma, das sofort aus \ref{FA.2.20} (iii) und \ref{FA.6.2.S} (i) folgt, ermöglicht es, den Grenzwert einer in $X_s$ konvergenten Folge, wobei $X$ ein normierter $\K$-Vek\-tor\-raum sei, zu definieren, vgl.\ \ref{FA.1.5}. 
Eine solche Folge nennt man dann \emph{schwach konvergent in $X$ gegen ihren (eindeutig bestimmten) Grenzwert}.\index{Konvergenz!schwache} 

\begin{Lemma} \label{FA.6.5}
Ist $X$ ein normierter $\K$-Vektorraum, so ist $X_s$ hausdorffsch. \q
\end{Lemma}

\begin{Satz} \label{FA.6.6} $\,$

\noindent  \textbf{Vor.:} Es seien $X$ ein normierter $\K$-Vektorraum, $(x_j)_{j \in \N}$ eine Folge in $X$ und $x \in X$.

\noindent  \textbf{Beh.:} $\D \lim_{j \to \infty} x_j = x \mbox{ in } X_s \Longleftrightarrow  \forall_{F \in X'} \lim_{j \to \infty} F(x_j) = F(x) \mbox{ in } \K$.
\end{Satz}

\textit{Beweis klar nach \ref{FA.6.2.S} (ii).} \q

\begin{Kor} \label{FA.6.6.K}
Es seien $X$ ein normierter $\K$-Vektorraum und $(x_j)_{j \in \N}$ eine Folge in $X$.

Ist $(x_j)_{j \in \N}$ schwach konvergent, so ist $(\|x_j\|)_{j \in \N}$ beschränkt.
\end{Kor}

\textit{Beweis.} Das Korollar folgt sofort aus dem letzten Satz und \ref{FA.2.28}. \q

\begin{Bem} \label{FA.6.6.B}
Sei $(x_j)_{j \in \N}$ eine Folge in einem normierten $\K$-Vek\-tor\-raum $X$.
Dann gilt wegen \ref{FA.6.4}:

\emph{Konvergiert $(x_j)_{j \in \N}$ bzgl.\ der Normtopologie gegen ein $x \in X$, so konvergiert $(x_j)_{j \in \N}$ schwach gegen $x$.}

Die Umkehrung dieser Aussage ist i.a.\ falsch, siehe \ref{FA.6.6.B.B} unten.
\end{Bem}

\begin{Satz} \label{FA.6.7}
Seien $X$ ein normierter $\K$-Vektorraum und $C$ eine konvexe Teilmenge von $X$.

Dann stimmt die abgeschlossene Hülle $\overline{C}$ von $C$ in $X$ mit der abgeschlossenen Hülle $\overline{C}^s$ von $C$ in $X_s$ überein.
\end{Satz}

\textit{Beweis.} Nach \ref{FA.6.4} ist $\overline{C} \subset \overline{C}^s$ (auch ohne die Voraussetzung der Konvexität) trivial -- beachte, daß die abgeschlossene Hülle der Schnitt aller abgeschlossenen Obermengen ist.

Zu ,,$\overline{C}^s \subset \overline{C}$``:
Ohne Beschränkung der Allgemeinheit können wir $C \ne \emptyset$ annehmen.
Sei $x_0 \in X \setminus \overline{C}$.
Dann gilt $d(\{x_0\},\overline{C}) > 0$ -- beachte die Fußnote auf Seite \pageref{FA.2.20.0} --, und der Trennungssatz von \textsc{Hahn-Banach} \ref{FA.2.24} (ii) ergibt die Existenz eines $F \in X' \setminus \{0\}$ sowie einer Zahl $\gamma \in \R$ mit
$$ \forall_{x \in \overline{C}} \, ({\rm Re} \, F)(x_0) < \gamma < ({\rm Re} \, F)(x). $$
Es folgt 
\begin{gather*}
U := \{ x \in X \, | \, {\rm Re} \, F(x) < \gamma \} = \left\{ \begin{array}{c} \overline{F}^1( {]} - \infty, \gamma {[} ) \mbox{ für } \K = \R \\ \overline{F}^1( {]} - \infty, \gamma {[} \times \R ) \mbox{ für } \K = \C \end{array} \right\} \in \U(x_0,X_s), \\
U \cap C = \emptyset,
\end{gather*} 
also gilt $x_0 \notin \overline{C}^s$. \q

\begin{Kor} \label{FA.6.7.K}
Seien $X$ ein normierter $\K$-Vektorraum und $C$ eine konvexe Teilmenge von $X$.
Dann gilt:
\begin{itemize}
\item[(i)] $C$ ist genau dann eine abgeschlossene Teilmenge von $X$, wenn $C$ eine abgeschlossene Teilmenge von $X_s$ ist.
\item[(ii)] $C$ ist genau dann dicht in $X$, wenn $C$ dicht in $X_s$ ist. \q
\end{itemize}
\end{Kor}

Wir zeigen in \ref{FA.6.7.B} 1.), daß die Einheitssphäre eines unendlich-dimensionalen normierten $\K$-Vektorraumes nicht schwach abgeschlossen ist.
Vorher führen wir die folgenden Begriffe ein, obwohl wir nur den ersten an dieser Stelle benötigen.

\begin{Def}[Umgebungsbasis, Basis einer Topologie, erstes und zweites Abzählbarkeitsaxiom] \label{FA.6.UB} \index{Basis!Umgebungs-} \index{Umgebungsbasis} \index{Basis!einer Topologie} \index{Abzählbarkeitsaxiome}
Es sei $X$ ein topologischer Raum. 
\begin{itemize}
\item[(i)] Sei $x_0 \in X$.
Dann heißt eine Teilmenge $\mathfrak{U}$ von $\U(x_0,X)$ eine \emph{Umgebungsbasis von $x_0$ (in $X$)} genau dann, wenn gilt $\forall_{U \in \U(x_0,X)} \exists_{V \in \mathfrak{U}} \, V \subset U$.
\item[(ii)] Eine Teilmenge $\mathfrak{U}$ von ${\rm Top}(X)$ heißt \emph{\mbox{Basis von ${\rm Top}(X)$}} genau dann, wenn gilt $\forall_{U \in {\rm Top}(X)} \, U = \bigcup_{V \in \mathfrak{U}, \, V \subset U} V$.
\item[(iii)] $X$ erfüllt per definitionem das \emph{erste Abzählbarkeitsaxiom}, wenn jedes Element $x_0 \in X$ eine höchstens abzählbare Umgebungsbasis besitzt.
\item[(iv)] $X$ erfüllt per definitionem das \emph{zweite Abzählbarkeitsaxiom}, wenn ${\rm Top}(X)$ eine höchstens abzählbare Basis besitzt.\\
(Aus dem zweiten Abzhählbarkeitsaxiom folgt das erste und die Separabilität\index{Raum!topologischer!separabler} des Raumes.)
\end{itemize}
\end{Def}

\begin{Lemma} \label{FA.6.7.L}
Seien $X$ ein normierter $\K$-Vektorraum und $x_0 \in X$.

Dann ist 
$$ \left\{ \bigcap_{\kappa = 1}^k \overline{F_{\kappa}}^1(U_{\varepsilon}(F_{\kappa}(x_0))) \, | \, k \in \N_+ \, \wedge \, \forall_{\kappa \in \{1, \ldots, k\}} \, F_{\kappa} \in X' \, \wedge \, \varepsilon \in \R_+ \right\} $$ 
eine Umgebungsbasis von $x_0$ in $X_s$.
\end{Lemma}

\textit{Beweis als Übung.} \q

\begin{Bsp} \label{FA.6.7.B}
Sei $X$ ein unendlich-dimensionaler normierter $\K$-Vektorraum.
\begin{itemize}
\item[1.)] Die Einheitssphäre $S := \{ x \in X \, | \, \|x\| = 1 \}$ ist eine nicht-abgeschlossene Teilmenge von $X_s$, genauer ist die abgeschlossene Hülle $\overline{S}^s$ von $S$ in $X_s$ gleich der Einheitsvollkugel $B_1(0)$.
\item[2.)] $U_1(0)$ ist eine nicht-offene Teilmenge von $X_s$, genauer ist der offene Kern $U_1(0)^{\circ s}$ von $U_1(0)$ in $X_s$ leer.
\end{itemize}

{[} Zu 1.): Wir beweisen zunächst
\begin{equation} \label{FA.6.7.B.1}
\forall_{x_0 \in U_1(0)} \forall_{U \in \U(x_0,X_s)} \, U \cap S \ne \emptyset.
\end{equation}

Zu (\ref{FA.6.7.B.1}): Seien $x_0 \in U_1(0)$ und $U \in \U(x_0,X_s)$.
Aus \ref{FA.6.7.L} folgt, daß es genügt, den Fall
\begin{equation} \label{FA.6.7.B.2}
U = \bigcap_{\kappa = 1}^k \overline{F_{\kappa}}^1(U_{\varepsilon}(F_{\kappa}(x_0))),
\end{equation}
wobei $k \in \N_+$, $\forall_{\kappa \in \{1, \ldots, k\}} \, F_{\kappa} \in X'$ und $ \varepsilon \in \R_+$ seien, zu betrachten.
Nun existiert $y_0 \in X \setminus \{0\}$ mit
\begin{equation} \label{FA.6.7.B.3}
\forall_{\kappa \in \{1, \ldots, k\}} \, F_{\kappa}(y_0) = 0,
\end{equation}
denn andernfalls wäre die $\K$-lineare Abbildung
$$ X \longrightarrow \K^k, ~~ x \longmapsto (F_1(x), \ldots, F_k(x)), $$
injektiv, im Widerspruch zu $\dim_{\K} X = \infty$. 

Wir definieren eine stetige Funktion $\varphi \: {[}0, \infty{[} \to {[}0, \infty{[}$ durch
$$ \forall_{t \in {[}0, \infty{[}} \, \varphi(t) := \|x_0 + t \, y_0\|. $$
Dann gilt $\varphi(0) = \|x_0\| < 1$ und $\lim_{t \to \infty} \varphi(t) = \infty$ -- beachte $x_0 \in U_1(0)$ und $y_0 \ne 0$ sowie $(C \in \R_+ \, \wedge \, |\varphi| \le C) \Rightarrow \forall_{t \in \R_+} \, t \, \|y_0\| \le C + \|x_0\|$ --, also folgt aus dem Zwischenwertsatz der Analysis I die Existenz einer Zahl $t_0 \in \R_+$ mit $\varphi(t_0) = 1$, d.h.\
$$ x_0 + t_0 \, y_0 \in S. $$
Ferner gilt $F_{\kappa}(x_0 + t_0 \, y_0 - x_0) \stackrel{(\ref{FA.6.7.B.3})}{=} 0$ für jedes $\kappa \in \{1, \ldots, k\}$, also wegen (\ref{FA.6.7.B.2}) auch
$$ x_0 + t_0 \, y_0 \in U. $$
Damit ist (\ref{FA.6.7.B.1}) gezeigt.

Aus (\ref{FA.6.7.B.1}) ergibt sich $S \subset B_1(0) \subset \overline{S}^s$, also mittels \ref{FA.6.7}
$$ \overline{S}^s \subset B_1(0) \subset \overline{S}^s. $$

Zu 2.): Angenommen, es existierte $x_0 \in U_1(0)^{\circ s} \subset U_1(0)$.
Dann gäbe es auch $U \in \U(x_0,X_s)$ mit $U \subset U_1(0)$.
Analog zu Beweis von (\ref{FA.6.7.B.1}) folgte die Existenz von $y_0 \in X \setminus \{0\}$ und $t_0 \in \R_+$ mit
$$ \|x_0 + t_0 \, y_0\| = 1 ~~ \mbox{ und } ~~ x_0 + t_0 \, y_0  \in U, $$
im Widerspruch zu $U \subset U_1(0)$. {]}
\end{Bsp}

\begin{Satz} \label{FA.6.St}
Seien $X,Y$ normierte $\K$-Vektorräume und $T \: X \to Y$ eine $\K$-lineare Abbildung.
Dann gilt:
$$ T \in \L_{\K}(X,Y) \Longleftrightarrow T \: X_s \to Y_s \mbox{ ist stetig.} $$
\end{Satz}

\textit{Beweis.} ,,$\Rightarrow$`` Nach \ref{FA.6.2.S} (iii) ist zu zeigen, daß für jedes $G \in Y'$ die $\K$-lineare Funktion
\begin{equation*} \label{FA.6.St.S}
G \circ T \: X_s \longrightarrow \K
\end{equation*}
stetig ist.
Seien daher $G \in Y'$, $\zeta \in \K$ und $\varepsilon \in \R_+$.
Aus der Voraussetzung der linken Seite folgt $G \circ T \in X'$, also ist $U := \overline{(G \circ T)}^1(U_{\varepsilon}(\zeta))$ eine in $X_s$ offene Menge mit $(G \circ T)(U) \subset U_{\varepsilon}(\zeta)$.
Hieraus folgt die Stetigkeit von $G \circ T \: X_s \to \K$. 

,,$\Leftarrow$`` \ref{FA.6.2.S} (iii) und die Voraussetzung der rechten Seite ergeben die Stetigkeit von $G \circ T \: X_s \to \K$ -- und damit offenbar auch die von $G \circ T \: X \to \K$ -- für jedes $G \in Y'$.

Angenommen, $T \: X \to Y$ ist nicht stetig, d.h.\ nach \ref{FA.2.7} ,,(ii) $\Leftrightarrow$ (iv)``, daß $T \: X \to Y$ kein beschränkter Operator ist.
Dann kann $T(B_1(0))$ keine beschränkte Teilmenge von $Y$ sein, also existiert nach \ref{FA.2.28} ein $G \in Y'$ derart, daß $G(T(B_1(0)))$ eine unbeschränkte Teilmenge von $Y$ ist, und $G \circ T \: X \to \K$ ist erneut nach \ref{FA.2.7} ,,(ii) $\Leftrightarrow$ (iv)`` nicht stetig, Widerspruch! \q
\A
Wir erinnern daran, daß wir zu einen normierten $\K$-Vektorraum $X$ in \ref{FA.2.TB.6} durch
$$ \forall_{x \in X} \forall_{F \in X'} \, \left( i_X(x) \right) (F) := F(x) $$
einen kanonischen isometrischen $\K$-Vektorraum-Homomorphismus 
$$ \boxed{i_X \: X \longrightarrow X''} $$
erhalten haben, der folglich injektiv und stetig ist.

\begin{Def}[Die schwache-$*$-Topologie des topologischen Dualraumes] \label{FA.6.8} \index{Topologie!schwache-$*$-}
Ist $X$ ein normierter $\K$-Vektorraum, so heißt die schwache Topologie für $X'$ bzgl.\ $i_X(X) \subset X''$ \emph{die schwache-$*$-Topologie für $X'$}.
Den topologischen Raum, der aus der Menge $X'$ zusammen mit dieser Topologie besteht, bezeichnen wir mit $\boxed{X'_{s*}}$.\footnote{\textbf{Warnung.} In \cite{Reck} wird $X'_{s*}$ in unserem Sinne mit $X'_s$ bezeichnet.}
\end{Def}

\begin{Bem*} \label{FA.6.9} $\,$
Sei $X$ ein normierter $\K$-Vektorraum.
\begin{itemize}
\item[1.)] ${\rm Top}(X'_{s*})$ ist gröber als ${\rm Top}(X'_s)$, und letztere ist gröber als die durch die Operatornorm induzierte Topologie für $X'$.
\item[2.)] Die Bezeichnung \emph{schwache-$*$-Topologie} rührt daher, daß viele Autoren den topologischen Dualraum mit $X^*$ bezeichnen.
Konsequenterweise müßten wir hier eigentlich von der ,,schwachen-$'$-Topologie'' reden.
Dies wäre allerdings vollkommen unüblich.
\end{itemize}
\end{Bem*}

Analog zu oben ermöglicht das folgende Lemma, vom Grenzwert einer in $X'_{s*}$ konvergenten Folge zu sprechen, wobei $X$ ein normierter $\K$-Vektorraum sei.
Eine derartige Folge heißt dann \emph{schwach-$*$-konvergent in $X'$ gegen ihren (eindeutig bestimmten) Grenzwert}.\index{Konvergenz!schwache-$*$-}

\begin{Lemma} \label{FA.6.10}
Ist $X$ ein normierter $\K$-Vektorraum, so ist $X'_{s*}$ hausdorffsch.
\end{Lemma}

\textit{Beweis.} Seien $F_1, F_2 \in X'$ mit $F_1 \ne F_2$.
Dann existiert $x \in X$ derart, daß gilt
$$ \left( i_X(x) \right) (F_1) = F_1(x) \ne F_2(x) = \left( i_X(x) \right) (F_2). $$
Daher folgt das Lemma aus \ref{FA.6.2.S} (i). \q

\begin{Satz} \label{FA.6.11} $\,$

\noindent \textbf{Vor.:} Es seien $X$ ein normierter $\K$-Vektorraum, $(F_j)_{j \in \N}$ eine Folge in $X'$ und $F \in X'$.

\noindent \textbf{Beh.:} $\D \lim_{j \to \infty} F_j = F \mbox{ in } X'_{s*} \Longleftrightarrow \forall_{x \in X} \, \lim_{j \to \infty} F_j(x) = F(x) \mbox{ in } \K$.
\end{Satz}

\textit{Beweis klar nach \ref{FA.6.2.S} (ii).} \q

\begin{Bem*}
Sei $X$ ein normierter $\K$-Vektorraum.
Der letzte Satz rechtfertigt, daß die Topologie für $X'_{s*}$ in \cite{Reck} die \emph{Topologie der punktweisen Konvergenz für $X'$}\index{Topologie!der punktweisen Konvergenz} genannt wird.
\end{Bem*}

\begin{Kor} \label{FA.6.11.K}
Es seien $X$ ein $\K$-Banachraum und $(F_j)_{j \in \N}$ eine Folge in $X'$.

Ist $(F_j)_{j \in \N}$ schwach-$*$-konvergent, so ist $(\|F_j\|)_{j \in \N}$ beschränkt.
\end{Kor}

\textit{Beweis.} Seien $(F_j)_{j \in \N}$ schwach-$*$-konvergent und $\mathcal{H} := \{ F_j \, | \, j \in \N \}$.
Nach \ref{FA.6.11} ist $\mathcal{H}(x)$ für jedes $x \in X$ beschränkt, also folgt aus \ref{FA.2.25}, daß $\mathcal{H}$ beschränkt in $X'$ ist.
\q

\begin{HS}[Satz von \textsc{Banach-Alaoglu}] \label{FA.6.BA} \index{Satz!von \textsc{Banach-Alaoglu}} $\,$

\noindent \textbf{Vor.:} Sei $X$ ein normierter $\K$-Vektorraum.

\noindent \textbf{Beh.:} Die Einheitsvollkugel $B':= \{ F \in X' \, | \, \|F\| \le 1 \}$ ist kompakt in $X'_{s*}$.
\end{HS}

\textit{Beweis.} Wir setzen für jedes $x \in X$
$$ X_x := \K $$
und versehen $\bigtimes_{\mathfrak{x} \in X} X_{\mathfrak{x}}$ mit der Produkttopologie, vgl.\ \ref{FA.C.P}.
Die Topologie für $\bigtimes_{\mathfrak{x} \in X} X_{\mathfrak{x}} = \K^X$ ist also die gröbste Topologie bzgl.\ derer für jedes $x \in X$ die Abbildung
$$ \pi_x \: \bigtimes_{\mathfrak{x} \in X} X_{\mathfrak{x}} \longrightarrow \K, ~~ f \longmapsto f(x), $$
stetig ist.
Da ${\rm Top}(X'_{s*})$ die gröbste Topologie für $X' \subset \K^X = \bigtimes_{\mathfrak{x} \in X} X_{\mathfrak{x}}$ ist bzgl.\ derer für jedes $x \in X$ die Abbildung
$$ i_X(x) \: X' \longrightarrow \K, ~~ F \longmapsto F(x), $$
stetig ist, stimmt der topologische Teilraum mit zugrundeliegender Menge $X'$ von $\bigtimes_{\mathfrak{x} \in X} X_{\mathfrak{x}}$ mit $X'_{s*}$ überein.

Definieren wir nun für jedes $x \in X$ die kompakte Menge
$$ B_x := B_{\|x\|}(0) \subset \K, $$
so ist $\bigtimes_{\mathfrak{x} \in X} B_{\mathfrak{x}}$ nach dem Satz von \textsc{Tychonoff} \ref{FA.C.T} eine kompakte Teilmenge von $\bigtimes_{\mathfrak{x} \in X} X_{\mathfrak{x}}$, und zum Nachweis des Hauptsatzes bleibt wegen \ref{FA.1.24} (i) zu zeigen, daß $B'$ eine abgeschlossene Teilmenge von $\bigtimes_{\mathfrak{x} \in X} B_{\mathfrak{x}}$ ist.

Wir behaupten
\begin{equation} \label{FA.6.BA.S}
B' = \bigcap_{x,y \in X, \, \lambda \in \K} \left\{ f \in \bigtimes_{\mathfrak{x} \in X} B_{\mathfrak{x}} \, |\, f(x+y) = f(x) + f(y) \, \wedge f(\lambda \, x) = \lambda \, f(x) \right\}.
\end{equation}

{[} Zu (\ref{FA.6.BA.S}): ,,$\subset$`` ist wegen wegen
$$ \forall_{F \in B'} \forall_{x \in X} \, |F(x)| \le \|F\| \, \|x\| \le \|x\| $$
klar.

,,$\supset$`` Jedes Element $f$ der rechten Seite von (\ref{FA.6.BA.S}) ist $\K$-linear, und es gilt wegen $f \in \bigtimes_{\mathfrak{x} \in X} B_{\mathfrak{x}}$
$$ \forall_{x \in X} \, \, |f(x)| \le \|x\|, $$
d.h., daß $f$ stetig ist sowie $\|f\| \le 1$, also $f \in B'$. {]}

Wegen der Stetigkeit der Abbbildung
$$ \pi_{\tilde{x}}|_{\bigtimes_{\mathfrak{x} \in X} B_{\mathfrak{x}}} \: \bigtimes_{\mathfrak{x} \in X} B_{\mathfrak{x}} \longrightarrow \K $$
für jedes $\tilde{x} \in X$ ist jede der im Schnitt der rechten Seite in (\ref{FA.6.BA.S}) auftretenden Menge als Schnitt zweier abgeschlossener Mengen selbst abgeschlossen.
Somit ist $B'$ nach (\ref{FA.6.BA.S}) eine abgeschlossene Teilmenge von $\bigtimes_{\mathfrak{x} \in X} B_{\mathfrak{x}}$.
\q

\begin{Satz} \label{FA.6.12}
Sei $X$ ein separabler normierter $\K$-Vektorraum, d.h.\ genau, daß $X$ eine höchstens abzählbare (bzgl.\ der Normtopologie) dichte Teilmenge besitzt.

Dann erfüllt $X'_{s*}$ das erste Abzählbarkeitsaxiom.
\end{Satz}

\textit{Beweisskizze.} Sei $M$ eine höchstens abzählbare dichte Teilmenge von $X$.
Zu vorgegebenem $F \in X'$ bilden die Mengen
\begin{eqnarray*}
U \left( F; x_1, \ldots,x_k; \frac{1}{m} \right) & := & \left\{ \widetilde{F} \in X' \, | \, \forall_{\kappa \in \{1, \ldots, k\}} \, |\widetilde{F}(x_{\kappa}) - F(x_{\kappa})| < \frac{1}{m} \right\} \\
& = & \bigcap_{\kappa = 1}^k \overline{i_X(x_{\kappa})}^1 (U_{\frac{1}{m}}(F(x_{\kappa})) ~ \in ~ \U(F,X'_{s*}),
\end{eqnarray*}
wobei $k,m \in \N_+$ und $x_1, \ldots, x_k \in M$ seien, eine abzählbare Umgebungsbasis von $F$ in $X'_{s*}$. \q

\subsection*{Charakterisierung reflexiver Räume} \addcontentsline{toc}{subsection}{Charakterisierung reflexiver Räume}

Der Vollständigkeit halber sei hier noch einmal erwähnt, daß ein normierter $\K$-Vek\-tor\-raum $X$ genau dann \emph{reflexiv} heißt, wenn $i_X \: X \to X''$ sogar eine $\K$-Vek\-tor\-raum-Iso\-metrie ist.

\begin{Bem*}
Wegen des folgenden Satzes \ref{FA.6.21} (iv) ,,$\Rightarrow$`` ist jeder reflexive $\K$-Vek\-tor\-raum ein $\K$-Ba\-nach\-raum.
\textsc{James} \cite{James} hat gezeigt, daß ein $\K$-Banachraum, der isometrisch zu seinem topologischen Bidualraum ist, nicht reflexiv zu sein braucht.
\end{Bem*}

\begin{Satz} \label{FA.6.21}
Sei $X$ ein normierter $\K$-Vek\-tor\-raum.
Dann gilt:
\begin{itemize}
\item[(i)] $X$ reflexiv $\Longrightarrow$ $X'_s = X'_{s*}$.
\item[(ii)] Ist $X$ reflexiv und $W$ ein abgeschlossener $\K$-Unter\-vek\-tor\-raum von $X$, so ist auch $W$ reflexiv.
\item[(iii)] Sind $Y$ ein weiterer normierter $\K$-Vek\-tor\-raum und $T \: X \to Y$ ein $\K$-Vek\-tor\-raum-Iso\-morphismus, so ist $X$ genau dann reflexiv, wenn $Y$ reflexiv ist.
\item[(iv)] $X$ reflexiv $\Longleftrightarrow$ $X$ $\K$-Banachraum und $X'$ reflexiv.
\end{itemize}
\end{Satz}

\textit{Beweis.} (i) ist wegen der Definitionen von $X'_s$ und $X'_{s*}$ klar.

Zu (ii): Sei $X$ reflexiv, und bezeichne $j \: W \hookrightarrow X$ die Inklusionsabbildung.
Nach \ref{FA.2.TB.9} kommutiert das folgende Diagramm:
%
%\begin{diagram}
%W             & \rEmbed^{j} & X                 \\
%\dEmbed^{i_W} &             & \dEmbedonto_{i_X} \\
%W''           & \rTo^{j''}  & X''               \\
%\end{diagram}
%
\begin{center}
\begin{tikzpicture}[node distance=2cm, every arrow/.style={thick,->}]
    % Knoten definieren
    \node (W) {$W$};
    \node (X) [right of=W] {$X$};
    \node (W'') [below of=W] {$W''$};
    \node (X'') [right of=W''] {$X''$};
    
    % Pfeile zeichnen
    \draw[>->] (W) -- node[above] {$j$} (X);
    \draw[>->] (W) -- node[left] {$i_W$} (W'');
    \draw[>-{>>}] (X) -- node[right] {$i_X$} (X'');
    \draw[->] (W'') -- node[above] {$j''$} (X'');
    \end{tikzpicture}
\end{center}
Wir haben zu zeigen, daß $i_W \: W \to W''$ surjektiv ist.
Sei daher $\Psi \in W''$.
Dann existiert (wegen der Reflexivität von $X$) ein Element $x \in X$ mit $j''(\Psi) = i_X(x)$, also
\begin{equation} \label{FA.6.21.1}
\forall_{F \in X'} \, \Psi(F|_W) = \Psi(F \circ j) \stackrel{\ref{FA.2.TB.4} (ii)}{=} (j''(\Psi))(F) = (i_X(x))(F) = F(x).
\end{equation}
Es folgt $x \in W$, denn andernfalls gäbe es wegen der Abgeschlossenheit von $W$ nach \ref{FA.2.20} (i) ein Funktional $F \in X'$ derart, daß $F(x) > 0$ und $F|_W = 0$, d.h.\ $\Psi(F|_W) = 0$, gilt, im Widerspruch zu (\ref{FA.6.21.1}).
Nunmehr zeigen wir
$$ j_W(x) = \Psi. $$

Zum Nachweis hiervon seien $f \in W'$ beliebig und $F \in X'$ mit $F|_W = f$ gemäß der ,,funktionalanalytischen Version`` des Fortsetzungssatzes von \textsc{Hahn-Banach} \ref{FA.2.19} gewählt.
Dann gilt
$$ \Psi(f) \stackrel{(\ref{FA.6.21.1})}{=} F(x) \stackrel{x \in W}{=} f(x) = (i_W(x))(f). $$

\pagebreak %nur mit tikz notwendig
Zu (iii): Sei $X$ reflexiv. Daß das Diagramm
%
%\begin{diagram}
%X                 & \rEmbedonto^{T} & Y             \\
%\dEmbedonto^{i_X} &                 & \dEmbed_{i_Y} \\
%X''               & \rTo^{T''}      & Y''           \\
%\end{diagram}
%
\begin{center}
\begin{tikzpicture}[node distance=2cm, every arrow/.style={thick,->}]
    % Knoten definieren
    \node (X) {$X$};
    \node (Y) [right of=X] {$Y$};
    \node (X'') [below of=X] {$X''$};
    \node (Y'') [right of=X''] {$Y''$};
    
    % Pfeile zeichnen
    \draw[>-{>>}] (X) -- node[above] {$T$} (Y);
    \draw[>-{>>}] (X) -- node[left] {$i_X$} (X'');
    \draw[>->] (Y) -- node[right] {$i_Y$} (Y'');
    \draw[->] (X'') -- node[above] {$T''$} (Y'');
    \end{tikzpicture}
\end{center}
dann kommutiert, folgt erneut aus \ref{FA.2.TB.9}, und somit ist $i_Y$ surjektiv, also $Y$ reflexiv.
Dies genügt offenbar zum Nachweis von (iii).

Zu (iv): ,,$\Rightarrow$`` Sei $i_X \: X \to X''$ eine $\K$-Vektorraum-Isometrie.
Wegen \ref{FA.2.11} ist somit $X$ ein $\K$-Banachraum.
Zu zeigen bleibt die Surjektivität des Operators $i_{X'} \: X' \to X''' := (X'')'$.

Sei also $\mathfrak{f} \in X'''$.
Dann gilt offenbar
\begin{equation} \label{FA.6.21.2}
F := \mathfrak{f} \circ i_X \in X',
\end{equation}
und zu jedem $\Phi \in X''$ existiert (genau) ein $x \in X$ mit
\begin{equation} \label{FA.6.21.3}
\Phi = i_X(x),
\end{equation}
also folgt
$$ \mathfrak{f}(\Phi) \stackrel{(\ref{FA.6.21.3})}{=} \mathfrak{f}(i_X(x)) \stackrel{(\ref{FA.6.21.2})}{=} F(x) = (i_X(x))(F) \stackrel{(\ref{FA.6.21.3})}{=} \Phi(F) =  (i_{X'}(F))(\Phi), $$
d.h.\ $i_{X'}(F) = \mathfrak{f}$.

,,$\Leftarrow$`` Sei $X'$ reflexiv.
Dann folgt aus der bereits bewiesenen Richtung ,,$\Rightarrow$``, daß $X''$ reflexiv ist.
Des weiteren ist $i_X(X)$ eine abgeschlossene Teilmenge von $X''$ -- beachte, daß $i_X \: X \to X''$ ein isometrischer $\K$-Vektorraum-Homomorphismus und $X$ nach Voraussetzung der rechten Seite vollständig ist, sowie \ref{FA.1.11} (ii). 
Nun ergibt (ii) die Reflexivität von $i_X(X)$ und somit (iii) die von $X$. \q 
\A
Aufgrund des folgenden Satzes interessieren wir uns für das Bild eines normierten $\K$-Vektorraumes $X$ unter $i_X \: X \to X''$ in $\boxed{X''_{s*}} := (X')'_{s*}$.

\begin{Satz} \label{FA.6.22.A}
Sei $X$ ein normierter $\K$-Vektorraum, und $(i_X(X))_{s*}$ bezeichne den topologischen Teilraum mit zugrundeliegender Menge $i_X(X)$ von $X''_{s*}$.

Dann ist $i_X \: X_s \to (i_X(X))_{s*}$ ein Homöomorphismus.
\end{Satz}

\textit{Beweis.} 1.) Trivialerweise ist $i_X \: X \to i_X(X)$ bijektiv.

2.) Jede in $(i_X(X))_{s*}$ offene Menge läßt sich als beliebige Vereinigung endlicher Schnitte von Mengen der Form
\begin{eqnarray*}
\overline{i_{X'}(F)}^1(V) \cap i_X(X) & = & \{ \Phi \in X'' \, | \, \overbrace{(i_{X'}(F))(\Phi)}^{= \Phi(F)} \in V \} \cap i_X(X) \\
& = & \{ i_X(x) \, | \, x \in X \, \wedge \, \underbrace{(i_X(x))(F)}_{= F(x)} \in V \},
\end{eqnarray*}
wobei $F \in X'$ und $V \in {\rm Top}(\K)$ sind, schreiben, welche
$$ \{ x \in X \, | \, F(x) \in V \} = \overline{F}^1(V) \in {\rm Top}(X_s) $$
als Urbild unter $i_X$ besitzen.
Hieraus folgt offenbar die Stetigkeit der Abbildung $i_X \: X_s \to (i_X(X))_{s*}$.

3.) Jede in $X_s$ offene Menge läßt sich als beliebige Vereinigung endlicher Schnitte von Mengen der Form
$$ \overline{F}^1(V) = \{ x \in X \, | \, F(x) \in V \} $$ 
wobei $F \in X'$ und $V \in {\rm Top}(\K)$ sind, schreiben, welche
$$ \{ i_X(x) \, | \, x \in X \, \wedge \, F(x) \in V \} \stackrel{\text{s.o.}}{=} \overline{i_{X'}(F)}^1(V) \cap i_X(X) \in {\rm Top}((i_X(X))_{s*}) $$
als Bild unter $i_X$ besitzen.
Wegen der Injektivität von $i_X$ folgt hieraus offenbar, daß $i_X \: X_s \to (i_X(X))_{s*}$ eine offene Abbildung ist -- beachte u.a., daß für injektive Abbildungen der Schnitt der Bilder gleich dem Bilde der Schnitte ist. 

Mit 1.) - 3.) ist der Satz bewiesen. \q

\begin{Satz} \label{FA.6.22} $\,$

\noindent \textbf{Vor.:} Seien $X$ ein normierter $\K$-Vektorraum und $B'' := \{ \Phi \in X'' \, | \, \|\Phi\| \le 1 \}$.

\noindent \textbf{Beh.:}
\begin{itemize}
\item[(i)] $i_X(B_1(0))$ ist in dem topologischen Teilraum $B''_{s_*}$ mit zugrundeliegender Menge $B''$ von $X''_{s*}$ dicht.
\item[(ii)] $i_X(X)$ ist dicht in $X''_{s*}$.
\end{itemize}
\end{Satz}

Hätten wir die wir den Trennungssatz von \textsc{Hahn-Banach} \ref{FA.2.24} (ii) auf dem Niveau hausdorffscher lokal-konvexer $\K$-Vektorräume bewiesen, so folgte hieraus Teil (i) der Behauptung, vgl.\ \cite[Proposition V.4.1]{Con}.
Wir benötigen zum Nachweis allerdings das folgende Lemma, das wir auch beim Beweis des u.g.\ Hauptsatzes \ref{FA.6.Char2} ausnutzen werden.

\begin{Lemma} \label{FA.6.22.L}
Seien $X$ ein normierter $\K$-Vektorraum, $\Phi \in X''$ mit $\| \Phi \| \le 1$, $k \in \N_+$ und $F_1, \ldots, F_k \in X'$.
Ferner sei $h \: X \to \R$ definiert durch
$$ \forall_{x \in X} \, h(x) := \sum_{\kappa = 1}^k |\Phi(F_{\kappa}) - F_{\kappa}(x)|^2. $$

Dann gilt $\inf \{ h(x) \, | \, x \in B_1(0) \} = 0$.
\end{Lemma}

\textit{Beweis.} Seien
\begin{equation} \label{FA.6.22.L.1}
I := \inf \{ h(x) \, | \, x \in B_1(0) \} = \inf \left\{ \sum_{\kappa = 1}^k |\Phi(F_{\kappa}) - F_{\kappa}(x)|^2 \, | \, x \in B_1(0) \right\} \ge 0
\end{equation}
und $(x_l)_{l \in \N}$ eine Folge in $B_1(0)$ mit
\begin{equation} \label{FA.6.22.L.2}
I = \lim_{l \to \infty} \sum_{\kappa = 1}^k |\Phi(F_{\kappa}) - F_{\kappa}(x_l)|^2.
\end{equation}
Es gilt für jedes $\kappa \in \{1, \ldots k\}$
$$ \forall_{l \in \N} \, \|F_{\kappa}(x_l)\| \le \|F_{\kappa}\| \, \|x_l\| \le \|F_{\kappa}\|, $$
also existieren nach dem Satz von \textsc{Bolzano-Weierstraß} ein gewisses $\xi_{\kappa} \in \K$ und eine Teilfolge von $(F_{\kappa}(x_l))_{l \in \N}$, die gegen $\xi_{\kappa}$ konvergiert.
Wir können offenbar ohne Einschränkung annehmen
\begin{equation} \label{FA.6.22.L.3}
\forall_{\kappa \in \{1, \ldots, k\}} \, \lim_{l \to \infty} F_{\kappa}(x_l) = \xi_{\kappa}
\end{equation}
und setzen
\begin{equation} \label{FA.6.22.L.4}
\forall_{\kappa \in \{1, \ldots, k\}} \, \zeta_{\kappa} := \Phi(F_{\kappa}) - \xi_{\kappa} \in \K,
\end{equation}
d.h.\
\begin{equation} \label{FA.6.22.L.5}
I \stackrel{(\ref{FA.6.22.L.2}) - (\ref{FA.6.22.L.4})}{=} \sum_{\kappa = 1}^k |\zeta_{\kappa}|^2
\end{equation}
Es folgt
\begin{equation} \label{FA.6.22.L.6}
\forall_{x \in B_1(0)} \, \sum_{\kappa = 1}^k {\rm Re} ( \overline{\zeta_{\kappa}} \, (F_{\kappa}(x) - \xi_{\kappa}) )  \le 0.
\end{equation}

{[} Zu (\ref{FA.6.22.L.6}): Seien $x \in B_1(0)$, $t \in {]}0,1{]}$ und $l \in \N$.
Dann gilt $t \, x + (1-t) \, x_l \in B_1(0)$ sowie
\begin{eqnarray*}
I & \stackrel{(\ref{FA.6.22.L.1})}{\le} & \sum_{\kappa = 1}^k |\Phi(F_{\kappa}) - F_{\kappa}(t \, x + (1-t) \, x_l)|^2 \\
& = & \sum_{\kappa = 1}^k {|} {(} \underbrace{\Phi(F_{\kappa}) - F_{\kappa}(x_l)}_{=: v_{\kappa}} {)} - \underbrace{t \,( F_{\kappa}(x) - F_{\kappa}(x_l)}_{=: w_{\kappa}} ) {|}^2 \\
& = & \sum_{\kappa = 1}^k (v_{\kappa} - w_{\kappa}) \, (\overline{v_{\kappa} - w_{\kappa}})
  =   \sum_{\kappa = 1}^k  v_{\kappa} \, \overline{v_{\kappa}} - v_{\kappa} \, \overline{w_{\kappa}} - \overline{v_{\kappa} \, \overline{w_{\kappa}}} + w_{\kappa} \, \overline{w_{\kappa}} \\
& = & \sum_{\kappa = 1}^k  |v_{\kappa}|^2 - 2 \, {\rm Re} ( v_{\kappa} \, \overline{w_{\kappa}} ) + |w_{\kappa}|^2 \\
& = & \sum_{\kappa = 1}^k |\Phi(F_{\kappa}) - F_{\kappa}(x_l)|^2 - 2 t \, \sum_{\kappa = 1}^k {\rm Re} ((F_{\kappa}(x) - F_{\kappa}(x_l)) \, (\overline{\Phi(F_{\kappa}) - F_{\kappa}(x_l)})) \\ && + \, t^2 \, \sum_{\kappa = 1}^k |F_{\kappa}(x) - F_{\kappa}(x_l)|^2.
\end{eqnarray*}
Die letzte Ungleichung sowie (\ref{FA.6.22.L.3}) - (\ref{FA.6.22.L.5}) ergeben durch Grenzwertbildung für $l \to \infty$
$$ I \le I - 2 t \, \sum_{\kappa = 1}^k {\rm Re} (F_{\kappa}(x) - \xi_{\kappa}) \, \overline{\zeta_{\kappa}}) + t^2 \, \sum_{\kappa = 1}^k |F_{\kappa}(x) - \xi_{\kappa}|^2, $$
d.h.\
$$ \sum_{\kappa = 1}^k {\rm Re} ((F_{\kappa}(x) - \xi_{\kappa}) \, \overline{\zeta_{\kappa}}) \le \frac{t}{2} \, \sum_{\kappa = 1}^k |F_{\kappa}(x) - \xi_{\kappa}|^2, $$
und wegen der Beliebigkeit von $t \in {]}0,1{]}$ muß (\ref{FA.6.22.L.6}) gelten. {]}

Sei nun
\begin{equation} \label{FA.6.22.L.F}
F := \sum_{\kappa = 1}^k \overline{\zeta_{\kappa}} \, F_{\kappa} \in X'.
\end{equation}
Dann folgt aus (\ref{FA.6.22.L.6})
\begin{equation} \label{FA.6.22.L.7} 
\forall_{x \in B_1(0)} \, {\rm Re} F(x) \le {\rm Re} \left( \sum_{\kappa = 1}^k \overline{\zeta_{\kappa}} \, \xi_{\kappa} \right),
\end{equation}
und wir behaupten 
\begin{equation} \label{FA.6.22.L.8}
\| F \| \le {\rm Re} \left( \sum_{\kappa = 1}^k \overline{\zeta_{\kappa}} \, \xi_{\kappa} \right).
\end{equation}

{[} Zu (\ref{FA.6.22.L.8}): Seien $\tilde{x} \in B_1(0)$ und $\phi \in {[}0, 2 \pi{[}$ mit
$$ F(\tilde{x}) = \e^{\i \phi} \, |F(\tilde{x})|, $$
d.h.\
\begin{gather*}
\underbrace{|F(\tilde{x})|}_{\in {[}0, \infty{[}} = F(\e^{- \i \phi} \, \tilde{x}) = {\rm Re} \, F(\e^{- \i \phi} \, \tilde{x}) \stackrel{(\ref{FA.6.22.L.7})}{\le} {\rm Re} \left( \sum_{\kappa = 1}^k \overline{\zeta_{\kappa}} \, \xi_{\kappa} \right).
\end{gather*}
Hieraus, der Beliebigkeit von $\tilde{x} \in B_1(0)$ sowie (\ref{FA.2.8.Stern}) ergibt sich (\ref{FA.6.22.L.8}). {]}

Außerdem konvergiert $(|F(x_l)|)_{l \in \N}$ wegen (\ref{FA.6.22.L.F}), (\ref{FA.6.22.L.3}) gegen $| \sum_{\kappa = 1}^k \overline{\zeta_{\kappa}} \, \xi_{\kappa} |$.
Somit folgt aus (\ref{FA.2.8.Stern}) und $\forall_{l \in \N} \, x_l \in B_1(0)$
$$ \| F \| \ge \left| \sum_{\kappa = 1}^k \overline{\zeta_{\kappa}} \, \xi_{\kappa} \right|. $$
Letzteres sowie (\ref{FA.6.22.L.8}) implizieren
$$ \sqrt{ \left( {\rm Re} \left( \sum_{\kappa = 1}^k \overline{\zeta_{\kappa}} \, \xi_{\kappa} \right) \right)^2 + \left( {\rm Im} \left( \sum_{\kappa = 1}^k \overline{\zeta_{\kappa}} \, \xi_{\kappa} \right) \right)^2 } \le \| F \| \le {\rm Re} \left( \sum_{\kappa = 1}^k \overline{\zeta_{\kappa}} \, \xi_{\kappa} \right), $$
also gilt
\begin{equation} \label{FA.6.22.L.9}
\| F \| = \sum_{\kappa = 1}^k \overline{\zeta_{\kappa}} \, \xi_{\kappa}.
\end{equation}

Schließlich ergibt sich
\begin{eqnarray*}
I & \stackrel{(\ref{FA.6.22.L.5})}{=} & \sum_{\kappa = 1}^k \overline{\zeta_{\kappa}} \, \zeta_{\kappa} \stackrel{(\ref{FA.6.22.L.4})}{=} \sum_{\kappa = 1}^k \overline{\zeta_{\kappa}} \, (\Phi(F_{\kappa}) - \xi_{\kappa}) \stackrel{(\ref{FA.6.22.L.F}), (\ref{FA.6.22.L.9})}{=} \Phi(F) - \|F\| \le |\Phi(F)| - \|F\| \\
& \le & \|F\| \, (\|\Phi\| - 1) \stackrel{\|\Phi\| \le 1}{\le} 0
\end{eqnarray*}
und somit wegen (\ref{FA.6.22.L.1}): $\inf \{ h(x) \, | \, x \in B_1(0) \} = 0$.
\q
\A
\textit{Beweis des Satzes.} Zu (i): Zunächst bemerken wir $i_X(B_1(0)) \subset B''$.
Es seien $\Phi \in B''$ und $U \in \U(\Phi,X''_{s*}).$
Wegen \ref{FA.1.8} (ii) und \ref{FA.1.TT} ist $(U \cap B'') \cap i_X(B_1(0)) \ne \emptyset$ zu zeigen, d.h.\
$$ U \cap i_X(B_1(0)) \ne \emptyset. $$

Als Umgebung von $\Phi$ in $X''_{s*}$ enthält $U$ eine Menge der Form
$$ \bigcap_{\kappa = 1}^k \overline{i_{X'}(F_{\kappa})}^1(V_{\kappa}) = \{ \Psi \in X'' \, | \, \forall_{\kappa \in\{1, \ldots, k\}} \, \Psi(F_{\kappa}) \in V_{\kappa} \}, $$
wobei $k \in \N_+$, $F_1, \ldots, F_k \in X'$ sowie $V_1 \in \U(\Phi(F_1),\K), \ldots, V_k \in \U(\Phi(F_k),\K)$ sind, also existiert eine Zahl $\varepsilon \in \R_+$ mit
$$ U \cap i_X(B_1(0)) \supset \bigg\{ i_X(x) \, | \, x \in B_1(0) \, \wedge \, \underbrace{\forall_{\kappa \in\{1, \ldots, k\}} \, | \Phi(F_{\kappa}) - F_{\kappa}(x) | < \sqrt{\frac{\varepsilon}{k}}}_{\Rightarrow \, \sum_{\kappa = 1}^k | \Phi(F_{\kappa}) - F_{\kappa}(x) |^2 < \varepsilon} \bigg\}, $$
und die rechte Seite ist wegen des vorherigen Lemmas \ref{FA.6.22.L} nicht leer.

Zu (ii): Für jedes $t \in \R_+$ ist der zu
$$ S_t \: X \longrightarrow X, ~~ x \longmapsto t \, x, $$
bitransponierte Operator, vgl.\ \ref{FA.2.TB.4} (ii), auch als Abbildung
$$ S''_t \: X''_{s*} \longrightarrow X''_{s*}, ~~ \Phi \longmapsto t \, \Phi, $$
stetig -- beachte, daß ${\rm Top}(X''_{s*})$ als schwache Topologie für $X''$ bzgl.\ $i_{X'}(X')$ definiert ist -- mit stetiger Umkehrabbildung $S''_{t^{-1}} \: X''_{s*} \to X''_{s*}$, also ein Homöomorphismus, und es gilt nach \ref{FA.2.TB.9}
$$ S''_t \circ i_X = i_X \circ S_t. $$
Daher folgt aus (i)
$$ S''_t(B'') = S''_t(\overline{i_X(B_1(0)}^{s*}) = \overline{S''_t(i_X(B_1(0))}^{s*} = \overline{i_X(S_t(B_1(0))}^{s*} \subset \overline{i_X(X)}^{s*}, $$
wobei $\overline{A}^{s*}$ die abgeschlossene Hülle einer Teilmenge $A$ von $X''$ in $X''_{s*}$ bezeichne.
Somit ergibt sich
$$ X'' = \bigcup_{t \in \R_+} S''_t(B'') \subset \overline{i_X(X)}^{s*}. $$
\q

\begin{HS} \label{FA.6.Char2} $\,$

\noindent \textbf{Vor.:} Sei $X$ ein normierter $\K$-Vektorraum.

\noindent \textbf{Beh.:} $X$ ist genau dann reflexiv, wenn die Einheitsvollkugel $B_1(0)$ in dem topologischen Raum $X_s$ kompakt ist.
\end{HS}

\textit{Beweis.} $B''$ bezeichne erneut die Einheitsvollkugel in $X''$.

,,$\Rightarrow$`` Sei $X$ reflexiv.
Dann folgt zunächst $B'' = i_X(B_1(0))$ sowie aus \ref{FA.6.22.A}, daß $i_X \: X_s \to X''_{s*}$ ein Homöomorphismus ist, und sodann wegen der Kompaktheit von $B''$ in $X''_{s*}$ (nach dem Satz von \textsc{Banach-Alaoglu} \ref{FA.6.BA} -- angewandt auf $X'$ anstelle von $X$ --), daß $B_1(0)$ kompakt in $X_s$ ist. 

,,$\Leftarrow$`` Sei $B_1(0)$ in $X_s$ kompakt.
Zu zeigen ist die Surjektivität des isometrischen $\K$-Vektorraum-Homomorphismus $i_X \: X \to X''$.
Hierfür genügt es, nachzuweisen, daß $i_X|_{B_1(0)} \: B_1(0) \to B''$ surjektiv ist.

Sei also $\Phi \in B''$.
Für jedes $F \in X'$ setzen wir
$$ A_F := \overline{F}^1(\{\Phi(F)\}) \cap B_1(0) = \{ x \in B_1(0) \, | \, \underbrace{F(x)}_{= (i_X(x))(F)} = \Phi(F) \} $$
und behaupten:
\begin{equation} \label{FA.6.Char2.0}
\mbox{$(A_F)_{F \in X'}$ ist ein zentriertes System in $B_1(0) \subset X_s$ abgeschlossener Mengen.}
\end{equation}
Dann folgt aus der Kompaktheit von $B_1(0)$ in $X_s$ sowie \ref{FA.1.22.Durchschnitt} ,,(i) $\Rightarrow$ (ii)``
$$ \bigcap_{F \in X'} A_F \ne \emptyset, $$
d.h.\ genau $\exists_{x \in B_1(0)} \, i_X(x) = \Psi$, und dies wollen wir gerade zeigen.

Zu (\ref{FA.6.Char2.0}): 1.) Daß $A_F$ für jedes $F \in X'$ eine in dem topologischen Teilraum $B_1(0)$ von $X_s$ abgeschlossene Teilmenge ist, folgt sofort daraus, daß ${\rm Top}(X_s)$ die gröbste Topologie für $X$ ist, so daß alle $F \in X'$ stetig sind.

2.) Seien $k \in \N_+$ und $F_1, \ldots, F_k \in X'$.
Mit $F_1, \ldots, F_k$ ist dann auch $h \: X_s \to \R$, definiert durch
$$ \forall_{x \in X} \, h(x) := \sum_{\kappa = 1}^k |\Phi(F_{\kappa}) - F_{\kappa}(x)|^2, $$
eine stetige Funktion, die folglich auf der kompakten Teilmenge $B_1(0)$ von $X_s$ ihr Minimum annimmt.
Somit existiert wegen \ref{FA.6.22.L} ein $x_0 \in B_1(0)$ mit $h(x_0) = 0$, d.h.\ $F_{\kappa}(x_0) = \Phi(x_0)$ für alle $\kappa \in \{1, \ldots, k\}$, also gilt
$$ \bigcap_{\kappa = 1}^k A_{F_{\kappa}} \ne \emptyset. $$
Daher ist $(A_F)_{F \in X'}$ ein zentriertes System.
\q

\begin{Kor} \label{FA.6.Char2.K} $\,$

\noindent \textbf{Vor.:} Sei $X$ ein normierter $\K$-Vektorraum.

\noindent \textbf{Beh.:} $X$ reflexiv $\Longleftrightarrow$ $X$ $\K$-Banachraum und $X'_s = X'_{s*}$.
\end{Kor}

\textit{Beweis.} ,,$\Rightarrow$`` ist bereits nach \ref{FA.6.21} (iv) ,,$\Rightarrow$`` und (i) klar.

,,$\Leftarrow$`` Aus $X'_s = X'_{s*}$ und dem Satz von \textsc{Banach-Alaoglu} \ref{FA.6.BA} folgt, daß die Einheitsvollkugel $B' := \{ F \in X' \, | \, \|F\| \le 1 \}$ eine kompakte Teilmenge von $X'_s$ ist.
Nach dem letzten Hauptsatz \ref{FA.6.Char2} ist daher $X'$ reflexiv.
Hieraus, der Voraussetzung, daß $X$ ein $\K$-Banachraum ist, sowie \ref{FA.6.21} (iv) ,,$\Leftarrow$`` ergibt sich die Reflexivität von $X$.
\q
\A
\ref{FA.6.Char2} und der u.g.\ Satz von \textsc{Eberlein-\v{S}mulian} \ref{FA.6.ES} ergeben eine weitere Charakterisierung der Reflexivität -- beachte, daß die Einheitsvollkugel eines normierten $\K$-Vektorraumes als konvexe abgeschlossene Teilmenge nach \ref{FA.6.7.K} (i) auch bzgl.\ der schwachen Topologie abgeschlossen ist.

\begin{HS} \label{FA.6.Char3} $\,$

\noindent \textbf{Vor.:} Sei $X$ ein normierter $\K$-Vektorraum.

\noindent \textbf{Beh.:} $X$ ist genau dann reflexiv, wenn die Einheitsvollkugel $B_1(0)$ in dem topologischen Raum $X_s$ folgenkompakt ist. \q
\end{HS}

\begin{Bem*} $\,$
\begin{itemize}
\item[1.)]
Daß aus der Reflexivität eines $\K$-Vektorraumes die schwache Folgenkompaktheit der Einheitsvollkugel folgt, ergibt sich auch ohne Verwendung des Satzes von \textsc{Eberlein-\v{S}mulian} aus \ref{FA.6.Char2}, vgl.\ z.B.\ \cite[Satz 14.9]{Hirz}.
In loc.\ cit.\ wird begründet, daß es genügt, separable $\K$-Vektorräume zu betrachten, und dann \ref{FA.6.22.A}, \ref{FA.6.21} (i) sowie \ref{FA.6.12} verwendet.
Kompakta, die das erste Abzählbarkeitsaxiom erfüllen, sind nach \cite[Anhang I.4]{Hirz} folgenkompakt.
\item[2.)] Die schwache Folgenkompaktheit der Einheitsvollkugel besagt übrigens genau, daß jede beschränkte Folge eine schwach konvergente Teilfolge besitzt.
Man könnte \ref{FA.6.Char3} ,,$\Rightarrow$`` daher als \textit{schwachen Satz von \textsc{Bolzano-Weierstrass}} bezeichnen.
\end{itemize}
\end{Bem*}

\subsection*{Die Approximationsaufgabe} \addcontentsline{toc}{subsection}{Die Approximationsaufgabe}

Seien $X$ ein normierter $\K$-Vekrorraum, $a \in X$ und $C$ eine nicht-leere Teilmenge von $X$.
Die \emph{Approximationsaufgabe\index{Approximationsaufgabe} $\boxed{A(a,C)}$} lautet, (mindestens) ein $x \in C$ mit
$$ \| a - x \| = \inf \{ \| a - y \| \, | \, y \in C \} \stackrel{\text{Def.}}{=} d( \{a\} , C) $$
zu finden.
Wir betrachten im folgenden nur solche Teilmengen $C$ von $X$, die sowohl abgeschlossen als auch konvex sind.

\begin{Bem} \label{FA.6.A.B} $\,$
\begin{itemize}
\item[1.)] Soll $A(a,C)$ für jedes $a \in X$ (mindestens) eine Lösung besitzen, so muß $C$ abgeschlossen in $X$ sein.
Ist $C$ nämlich nicht abgeschlossen, so existiert $a_* \in \overline{C} \setminus C$, also gilt $d(\{a_*\},C) = 0$, und $A(a_*,C)$ besitzt keine Lösung.
\item[2.)] Es ist anschaulich klar, daß es im Falle $n \in \N_+$ und $X = \R^n$ nicht-konvexe Teilmengen $C$ von $X$ und $a \in X$ gibt so, daß $A(a,C)$ mehrere Lösungen besitzt.
Dies gilt z.B.\ für $C = \R \setminus {]}-1,1{[} \subset \R$ und $a = 0$.
\end{itemize}
\end{Bem}

\begin{Satz} \label{FA.6.Approx} \index{Approximationsaufgabe} $\,$

\noindent \textbf{Vor.:} Seien $X$ ein reflexiver $\K$-Vekrorraum, $a \in X$ und $C$ eine nicht-leere abgeschlossene konvexe Teilmenge von $X$.

\noindent \textbf{Beh.:} $A(a,C)$ besitzt (mindestens) eine Lösung.
\end{Satz}

Wir bereiten den Beweis des Satzes durch das folgende Lemma vor.

\begin{Lemma} \label{FA.6.Approx.L}
Seien $X$ ein normierter $\K$-Vektorraum, $(x_j)_{j \in \N}$ eine Folge in $X$ und $x \in X$.

Dann gilt:
$$ \lim_{j \to \infty} x_j = x \mbox{ in } X_s \Longrightarrow \|x\| \le \liminf_{j \to \infty} \|x_j\|. $$
\end{Lemma}

\textit{Beweis.} Sei also $(x_j)_{j \in \N}$ schwach konvergent gegen $x$.
Ohne Einschränkung gelte $x \ne 0$.
Dann existiert nach \ref{FA.2.20} (ii) ein Funktional $F \in X'$ mit $F(x) = \|x\|$ sowie $\|F\| = 1$, und für jedes $j \in \N$ gilt $|F(x_j)| \le \|F\| \, \|x_j\| = \|x_j\|$.
Daher folgt $\|x\| = |F(x)| \stackrel{\ref{FA.6.6}}{=} \lim_{j \to \infty} |F(x_j)| = \liminf_{j \to \infty} |F(x_j)| \le \liminf_{j \to \infty} \|x_j\|$. \q
\A
\textit{Beweis des Satzes.} Sei $(y_j)_{j \in \N}$ eine Folge in $C$ mit
\begin{equation} \label{FA.6.Approx.1}
\lim_{j \to \infty} \|a - y_j\| = d(\{a\},C).
\end{equation}
Dann ist $(\|y_j\|)_{j \in \N}$ beschränkt, also können wir ohne Einschränkung annehmen, daß $\forall_{j \in \N} \, y_j \in B_1(0)$ gilt.
Nach \ref{FA.6.Char3} existiert daher wegen der Reflexivität von $X$ eine Teilfolge von $(y_j)_{j \in \N}$, die gegen ein gewisses $x \in B_1(0)$ in $X_s$ konvergiert.
Ohne Beschränkung der Allgemeinheit gelte 
\begin{equation}
\lim_{j \to \infty} y_j = x \mbox{ in } X_s, \label{FA.6.Approx.2}
\end{equation}
d.h.\ (z.B.\ nach \ref{FA.6.6})
\begin{gather} 
\lim_{j \to \infty} a - y_j = a - x \mbox{ in } X_s. \label{FA.6.Approx.3}
\end{gather}
Da $(y_j)_{j \in \N}$ eine Folge in $C$ ist, folgt aus (\ref{FA.6.Approx.2})
\begin{equation} \label{FA.6.Approx.4} 
x \in \overline{C}^s \stackrel{\ref{FA.6.7}}{=} \overline{C} = C,
\end{equation}
beachte, daß $C$ eine abgeschlossene konvexe Teilmenge von $X$ ist.
Außerdem ergeben (\ref{FA.6.Approx.3}) und \ref{FA.6.Approx.L}
\begin{equation} \label{FA.6.Approx.5}
\| a - x \| \le \liminf_{j \to \infty} \| a - y_j \|.
\end{equation}
Schließlich gilt
$$ d(\{a\},C) \stackrel{(\ref{FA.6.Approx.4})}{\le} \| a - x \| \stackrel{(\ref{FA.6.Approx.5})}{\le} \liminf_{j \to \infty} \| a - y_j \| \stackrel{(\ref{FA.6.Approx.1})}{=} d(\{a\},C), $$
und $x \in C$ ist somit eine Lösung von $A(a,C)$.
\q

\subsection*{Anhang: Der Satz von \textsc{Eberlein-\v{S}mulian}} \addcontentsline{toc}{subsection}{Anhang: Der Satz von \textsc{Eberlein-\v{S}mulian}}

I.a.\ ist die schwache Topologie eines normierten $\K$-Vektorraumes nicht metrisierbar.
Dennoch gilt das folgende Resultat.

\begin{HS}[Satz von \textsc{Eberlein-\v{S}mulian}] \index{Satz!von \textsc{Eberlein-\v{S}mulian}} \label{FA.6.ES} $\,$

\noindent \textbf{Vor.:} Sei $X$ ein normierter $\K$-Vektorraum.

\noindent \textbf{Beh.:} Eine Teilmenge von $X$ ist genau dann kompakt in $X_s$, wenn sie sowohl abgeschlossen in $X_s$ als auch folgenkompakt in $X_s$ ist.
\end{HS}

Für $\K$-Banachräume haben \textsc{\v{S}mulian} 1940 in \cite{Smu} für \emph{abzählbar kompakte} Teilmengen die Hinrichtung und \textsc{Eberlein} 1947 in \cite{Eb} die Rückrichtung der Behauptung bewiesen, daher rührt der Name des Satzes.
Es gibt (wie bei fast jedem mathematischen Satz) Verallgemeinerungen, die ebenso genannt werden.
Wir skizzieren den Beweis des nächstgenannten Hauptsatzes -- jener folgt der Beweisführung in \cite[Chapter III]{Die}, die ursprünglich 1967 von \textsc{Whitley} \cite{Whit} präsentiert wurde, -- und zeigen, daß \ref{FA.6.ES} leicht daraus folgt.
Wir weisen an dieser Stelle bereits ausdrücklich auf die erste Bemerkung auf Seite \pageref{FA.6.ES.B} hin!

\begin{Def}[Relative Folgenkompaktheit] \index{Kompaktheit!Folgen-!relative} \label{FA.6.31}
Es seien $X$ ein topologischer Raum und $K$ eine Teilmenge von $X$.

$K$ heißt \emph{relativ folgenkompakt} genau dann, wenn jede Folge in $K$ eine in $X$ konvergente Teilfolge besitzt.
\end{Def}

\begin{HS} \label{FA.6.32} $\,$

\noindent \textbf{Vor.:} Seien $X$ ein normierter $\K$-Vektorraum und $K$ eine Teilmenge von $X$.

\noindent \textbf{Beh.:} $K \subset \subset X_s$ $\Longleftrightarrow$ $K$ relativ folgenkompakt in $X_s$.
\end{HS}

\textit{Beweisskkizze.} 1.) Sei $W$ ein normierter $\K$-Vektorraum.
Eine Teilmenge $\mathcal{H}$ von $W'$ heißt \emph{total} genau dann, wenn $(\forall_{F \in \mathcal{H}} \, F(w) = 0) \Rightarrow w = 0$ für alle $w \in W$ gilt. 

2.) Ist $W$ ein nicht-nulldimensionaler separabler normierter $\K$-Vektorraum, so besitzt $W'$ eine abzählbare totale Teilmenge, die in $\{ F \in W' \, | \, \|F\| = 1 \}$ enthalten ist.

{[} Ist nämlich $\{ w_j \, | \, j \in \N \}$ eine abzählbare dichte Teilmenge von $W$ -- beachte, daß eine dichte Teilmenge von $W$ nicht endlich sein kann --, so kann man offenbar $\forall_{j \in \N} \, w_j \ne 0$ annehmen, und es existiert nach \ref{FA.2.20} (ii)
zu jedem $j \in \N$ ein $F_j \in W'$ mit $F_j(w_j) = \|w_j\|$ sowie $\|F_j\| = 1$.
Dann ist $\{ F_j \, | \, j \in \N \}$ total, denn wegen $\overline{\{ w_j \, | \, j \in \N \}} = W$ existiert zu jedem $w \in W$ mit $\forall_{j \in \N} \, F_j(w) = 0$ sowie jedem $\varepsilon \in \R_+$ ein $j_0 \in \N$ mit
$$ \| w - w_{j_0} \| < \frac{\varepsilon}{2}, $$
also folgt
$$ \| w_{j_0} \| = | F_{j_0}(w_{j_0}) | = | F_j (w - w_{j_0}) | \le \| F_{j_0} \| \, \| w - w_{j_0} \| < \frac{\varepsilon}{2} $$
und $\| w \| \le \| w - w_{j_0} \| + \| w_{j_0} \| < \varepsilon$. {]}

3.) Sei $W$ ein normierter $\K$-Vektorraum derart, daß $W'$ eine abzählbare totale Teilmenge, die in $\{ F \in W' \, | \, \|F\| = 1 \}$ enthalten ist, besitzt.
Dann ist jeder kompakte topologische Teilraum von $W_s$ metrisierbar.

{[} Ist $\{ F_j \, | \, j \in \N \}$ eine abzählbare totale Teilmenge von $\{ F \in W' \, | \, \|F\| = 1 \}$, die dicht in $W'$ liegt, so wird durch
$$ \forall_{w,\tilde{w} \in W} \, d(w,\tilde{w}) := \sum_{j=0}^{\infty} \frac{1}{2^{j+1}} \, | F_j(w - \tilde{w}) | ~ ( \, \le \| w - \tilde{w} \| < \infty ) $$
eine Metrik $d$ auf $W$ definiert. (Hier geht die Totalität von $\{ F_j \, | \, j \in \N \}$ ein!)

Sei nun $\widetilde{K}$ eine Teilmenge von $W$ derart, daß der topologischer Teilraum $\widetilde{K}_s$ mit zugrundeliegender Menge $\widetilde{K}$ von $W_s$ kompakt ist.
Zu zeigen bleibt, daß
$$ \id_{\widetilde{K}} \: \widetilde{K}_s \longrightarrow (\widetilde{K},d_{\widetilde{K} \times \widetilde{K}}) $$
ein Homöomorphismus ist.
$\widetilde{K}_s$ ist kompakt und $(\widetilde{K},d_{\widetilde{K} \times \widetilde{K}})$ als metrischer Raum hausdorffsch, also genügt es, die Stetigkeit der Abbildung $\id_{\widetilde{K}} \: \widetilde{K}_s \to (\widetilde{K},d_{\widetilde{K} \times \widetilde{K}})$ zu beweisen, beachte \ref{FA.1.26}.

Zunächst ist $\widetilde{K}$ bzgl.\ $\| \ldots \|$ beschränkt.
(Dies folgt aus \ref{FA.2.28}, denn für alle $F \in W'$ ist $F \: W_s \to \K$ nach der Definition der schwachen Topologie stetig, und die Kompaktheit von $\widetilde{K}$ in $W_s$ ergibt mittels \ref{FA.1.25} die Kompaktheit -- und damit die Beschränktheit -- von $F(K) \subset \K$.)
Ohne Beschränkung der Allgemeinheit gelte daher
$$ \widetilde{K} \subset \{ w \in W \, | \, \| w \| \le 1 \}. $$
Seien $w_0 \in \widetilde{K}$ und $\varepsilon \in \R_+$.
Dann existiert also eine Zahl $j_0 \in \N$ mit
$$ \forall_{w \in \widetilde{K}} \, \sum_{j=j_0+1}^{\infty} \frac{1}{2^{j+1}} \underbrace{\|F_j\|}_{= \, 1} \underbrace{\|w-w_0\|}_{\le \, 2} \le \sum_{j=j_0+1}^{\infty} \frac{1}{2^j} = \frac{1}{2^{j_0}} \sum_{j=1}^{\infty} \frac{1}{2^j} = \frac{1}{2^{j_0}} < \frac{\varepsilon}{2}. $$
Des weiteren gilt
\begin{eqnarray*}
V & := & \left\{ w \in \widetilde{K} \, | \, \forall_{j \in \{0, \ldots,j_0\}} \, |F_j (w-w_0)| < \frac{\varepsilon}{2 \cdot 2^{j_0+1}} \right\} \\
& = & \left( \bigcap_{j=0}^{j_0} \overline{F_j}^1(U^{\K}_{\frac{\varepsilon}{2 \cdot 2^{j_0+1}}}(0)) \right) \cap \widetilde{K} ~ \in ~ \U(w_0,\widetilde{K}_s)
\end{eqnarray*}
und
$$ \forall_{w \in V} \, \sum_{j=0}^{j_0} \frac{1}{2^{j+1}} |F_j (w - w_0)| \le 2^{j_0+1} |F_j(w-w_0)| < \frac{\varepsilon}{2}. $$
Insgesamt ergibt sich für jedes $w \in V$
\begin{eqnarray*}
d(w,w_0) & = & \sum_{j=0}^{\infty} \frac{1}{2^{j+1}} |F_j (w - w_0)| \\
& \le & \sum_{j=0}^{j_0} \frac{1}{2^{j+1}} |F_j (w - w_0)| + \sum_{j=j_0+1}^{\infty} \frac{1}{2^{j+1}} \|F_j\| \|w-w_0\| < \varepsilon,
\end{eqnarray*}
d.h.\ $\id_{\widetilde{K}}(V) \subset \{ w \in K \, | \, d(w,w_0) < \varepsilon \}$.
Folglich ist $\id_{\widetilde{K}} \: \widetilde{K}_s \to (\widetilde{K},d_{\widetilde{K} \times \widetilde{K}})$ stetig. {]}

4.) Beweis von ,,$\Rightarrow$``: Seien $K \subset \subset X_s$ und $(x_j)_{j \in \N}$ eine Folge in $K$.
Ohne Einschränkung sei $(x_j)_{j \in \N}$ nicht konstant.
Setze
$$ W := \overline{{\rm Span}_{\K} \{ x_j \, | \, j \in \N \}}^{\| \ldots \|} \, ( \ne \{0\} ), $$
d.i.\ die abgeschlossene Hülle der Menge alle Linearkombinationen $\sum_{j=0}^{\infty} \lambda_j \, x_j$, wobei $\lambda_j \in \K$ nur für endlich viele $j \in \N$ ungleich Null ist, bzgl.\ der Normtopologie.
Dann ist $W$ nach \ref{FA.2.4} (ii) ein abgeschlossener $\K$-Untervektorraum von $X$ und wegen seiner Konvexität -- beachte Beispiel 5.) zu \ref{FA.N.2} (iv) -- sowie \ref{FA.6.7.K} (i) auch eine abgeschlossene Teilmenge von $X_s$.
Aus $K \subset \subset X_s$ folgt mittels \ref{FA.1.24} (i) des weiteren $K \cap W \subset \subset X_s$, also insgesamt
$$ \{ x_j \, | \, j \in \N \} \subset K \cap W \subset \subset W_s, $$
denn der topologische Teilraum mit zugrundeliegender Menge $W$ von $X_s$ stimmt mit $W_s$ überein.
Letzteres liegt daran, daß ein Element von ${\rm Top}_{X_s}(W)$ bzw.\ ${\rm Top}(W_s)$ sich als beliebige Vereinigung endlicher Schnitte von Mengen der Form
$$ \overline{F}^1(V) \cap W = \overline{F|_W}^1(V) \mbox{ mit } F \in X' \mbox{ (also insbes.\  $F|_W \in W'$)  und } V \in {\rm Top}(\K) $$
bzw.\
$$ \overline{f}^1(V) \mbox{ mit } f \in W' \mbox{ und } V \in {\rm Top}(\K) $$
schreiben läßt, und zu jedem $f \in W'$ nach \ref{FA.2.19} ein $F \in X'$ mit $F|_W = f$ existiert.
\pagebreak
Da $W \ne \{0\}$ als normierter $\K$-Vektorraum separabel ist (weil
$$ \left\{ \sum_{j=0}^k \lambda_j \, x_j \, | \, k \in \N \mbox{ und $\lambda_j \in \left\{ \begin{array}{c} \Q \mbox{ im Falle } \K = \R \\ \Q + \i \, \Q \mbox{ im Falle } \K = \C \end{array} \right\}$ für $j \in \{0, \ldots, k\}$} \right\} $$
dicht in $W$ liegt), ergeben 2.) und 3.) nun, daß der kompakte Teilraum $\overline{K \cap W}^s$ von $W_s$ metrisierbar ist, also besitzt $(x_j)_{j \in \N}$ nach \ref{FA.1.30} ,,(i) $\Rightarrow$ (ii)`` eine in $W_s \subset X_s$ konvergente Teilfolge.
Damit ist die relative Folgenkompaktheit von $K$ in $X_s$ gezeigt.

5.) Sei $E$ ein endlich-dimensionaler $\K$-Untervektorraum von $X''$.
Dann existieren $k \in \N_+$ und $F_1, \ldots, F_k \in S' := \{ F \in X' \, | \, \|F\| = 1 \}$ mit
$$ \forall_{\Psi \in E} \, \frac{\|\Psi\|}{2} \le \max \{ |\Psi(F_{\kappa})| \, | \, \kappa \in \{1, \ldots, k\} \}. $$

{[} $E$ ist in kanonischer Weise ein endlich-dimensionaler normierter $\K$-Vek\-tor\-raum, also ergibt \ref{FA.2.E.11} ,,(i) $\Rightarrow$ (iii)`` die Kompaktheit der Einheitsvollkugel von $E$.
Offenbar ist dann auch die Einheitssphäre $S_E := \{ \Psi \in E \, | \, \|\Psi\| = 1 \}$ kompakt, und es existieren $k \in \N_+$ sowie $\Psi_1, \ldots, \Psi_k \in S_E$ derart, daß gilt $S_E \subset \bigcup_{\kappa=1}^k U_{\kappa}$, wobei $U_{\kappa} := \{ \Psi \in E \, | \, \| \Psi - \Psi_{\kappa} \| < \frac{1}{4} \}$ für $\kappa \in \{1, \ldots, k\}$ sei.
Zu jedem solchen $\kappa$ gibt es des weiteren ein $F_{\kappa} \in S'$ mit $|\Psi_{\kappa}(F_{\kappa})| > \frac{3}{4}$.
Nun existiert zu jedem $\Psi \in S_E$ eine Zahl $\kappa_0 \in \{1, \ldots, k\}$ mit $\| \Psi - \Psi_{\kappa_0} \| < \frac{1}{4}$, d.h.\
\begin{gather*}
|\Psi(F_{\kappa_0})| \le |(\Psi - \Psi_{\kappa_0})(F_{\kappa_0})| + |\Psi_{\kappa_0}(F_{\kappa_0})| \le \| \Psi - \Psi_{\kappa_0} \| \, \|F_{\kappa_0} \| + | \Psi_{\kappa_0}(F_{\kappa_0}) |, \\
|\Psi_{\kappa_0}(F_{\kappa_0})| \ge |\Psi(F_{\kappa_0})| - \| \Psi - \Psi_{\kappa_0} \| \, \|F_{\kappa_0} \| > \frac{3}{4} - \frac{1}{4} = \frac{1}{2}, \\
\max \{ |\Psi(F_{\kappa})| \, | \, \kappa \in \{1, \ldots, k\} \} > \frac{1}{2}.
\end{gather*}
Für beliebiges $\Psi \in E \setminus \{0\}$ folgt hieraus das Gewünschte, indem man von $\Psi$ zu $\frac{\Psi}{\|\Psi\|} \in S_E$ übergeht, und für $\Psi = 0$ ist nichts zu zeigen. {]}

6.) Beweis von ,,$\Leftarrow$``: Sei $K$ relativ folgenkompakt in $X_s$.
Aus \ref{FA.6.6} folgt dann, daß $F(K)$ für jedes $F \in X'$ relativ folgenkompakt in $\K$ und somit offenbar beschränkt ist.
\ref{FA.2.28} ergibt dann die Beschränktheit von $K \subset X$, also ist auch $i_X(K) \subset X''$ beschränkt.
Somit folgt aus dem Satz von \textsc{Banach-Alaoglu} \ref{FA.6.BA} und \ref{FA.1.24} (i), daß die abgeschlossene Hülle $\overline{i_X(K)}^{s*}$ eine kompakte Teilmenge von $X''_{s*} \stackrel{\text{Def.}}{=} (X')'_{s*}$ ist.
Zu zeigen bleibt $\overline{i_X(K)}^{s*} \subset i_X(X)$, denn dann liefert \ref{FA.6.22.A} die Kompaktheit von $\overline{K}^s$ in $X_s$.

Sei also $\Phi \in \overline{i_X(K)}^{s*}$.
Ferner seien $k_0 := 0$ und aus technischen Gründen $k_{-1} := 0$.
Für jedes $j \in \N$ werden nun rekursiv endlich-dimensionale $\K$-Untervektorräume von $X''$
\begin{equation} \label{FA.6.ES.0} 
E_j := {\rm Span}_{\K} \{ \Phi, i_X(x_0), \ldots, i_X(x_{k_{j-1}}) \}
\end{equation}
definiert und $k_j \in \N_+$ mit $k_j > k_{j-1}$, $F_{k_{j-1}+1}, \ldots, F_{k_j} \in S' = \{ F \in X' \, | \, \|F\| = 1 \}$ derart, daß
\begin{equation} \label{FA.6.ES.1}
\forall_{\Psi \in E_j} \, \frac{\|\Psi\|}{2} \le \max \{ |\Psi(F_{\kappa})| \, | \, \kappa \in \{0, \ldots, k_j\} \}
\end{equation}
gilt, sowie $x_j \in K$ mit
\begin{equation} \label{FA.6.ES.2}
\forall_{\kappa \in \{0, \ldots, k_j\}} \, |(i_X(x_j))(F_{\kappa}) - \Phi(F_{\kappa})| < \frac{1}{j+1}
\end{equation}
gewählt.
\pagebreak

{[} Wähle für $j=0$ ein $F_0 \in S'$ mit $\frac{\|\Phi\|}{2} \le |\Phi(F_0)|$.
Gemäß der Definition (\ref{FA.6.ES.0}) ist (\ref{FA.6.ES.1}) dann trivialerweise erfüllt.
Es gilt weiterhin
$$ \{ \Psi \in X'' \, | \, |\underbrace{(i_{X'}(F_0))(\Psi)}_{= \Psi(F_0)} - \Phi(F_0)| < 1 \} = \overline{i_{X'}(F_0)}^1(U_1^{\K}(\Phi(F_0))) \in \U(\Phi,X''_{s*}) $$
und (weil $\Phi$ ein Berührungspunkt von $i_X(K)$ in $X''_{s*}$ ist)
\begin{eqnarray*}
\emptyset & \ne & \{ \Psi \in X'' \, | \, |\Psi(F_0) - \Phi(F_0)| < 1 \} \cap i_X(K) \\
& & \, = \{ i_X(x) \, | \, x \in K \, \wedge \, |(i_X(x))(F_0) - \Phi(F_0)| < 1 \},
\end{eqnarray*}
also existiert $x_0 \in K$ mit
$$ |(i_X(x_0))(F_0) - \Phi(F_0)| < 1, $$
d.h.\ auch (\ref{FA.6.ES.2}) ist erfüllt.

Sei nun $j \in \N_+$ und die Konstruktion für $j-1$ bereits erfolgt.
Der $\K$-Unter\-vek\-tor\-raum $E_j \stackrel{(\ref{FA.6.ES.0})}{=} {\rm Span}_{\K} \{ \Phi, i_X(x_0), \ldots, i_X(x_{k_j - 1}) \}$ von $X''$ ist endlich-dimensional.
Nach 5.) gibt es zusätzlich zu den bereits vorhandenen $F_0, \ldots, F_{k_{j-1}}$ ein $k_j \in \N_+$, $k_j > k_{j-1}$, und $F_{k_{j-1}+1}, \ldots, F_{k_j} \in S'$ mit (\ref{FA.6.ES.1}).
Analog zu oben gilt weiterhin
\begin{eqnarray*}
\lefteqn{\left\{ \Psi \in X'' \, | \, \forall_{\kappa \in \{0, \ldots, k_j\}} \, |\Psi(F_{\kappa}) - \Phi(F_{\kappa})| < \frac{1}{j+1} \right\}} && \\
&& ~~~~~~~~~~~~~~~~ = \bigcap_{\kappa=0}^{k_j} \overline{i_{X'}(F_{\kappa})}^1(U_{\frac{1}{j+1}}^{\K}(\Phi(F_{\kappa}))) \in \, \U(\Phi,X''_{s*}), 
\end{eqnarray*}
und (wieder weil $\Phi$ ein Berührungspunkt von $i_X(K)$ in $X''_{s*}$ ist)
\begin{eqnarray*}
\emptyset & \ne & \left\{ \Psi \in X'' \, | \, \forall_{\kappa \in \{0, \ldots, k_j\}} \, |\Psi(F_{\kappa}) - \Phi(F_{\kappa})| < \frac{1}{j+1} \right\} \cap i_X(K) \\
& & \, = \left\{ i_X(x) \, | \, x \in K \, \wedge \, \forall_{\kappa \in \{0, \ldots, k_j\}} \, |(i_X(x))(F_{\kappa}) - \Phi(F_{\kappa})| < \frac{1}{j+1} \right\},
\end{eqnarray*}
also existiert $x_j \in K$ mit (\ref{FA.6.ES.2}). {]}

Nun folgt aus der relativen Folgenkompaktheit von $K$ in $X_s$ die Existenz einer Teilfolge von $(x_j)_{j \in \N}$, die in $X_s$ gegen ein gewisses $x \in X$ konvergiert.
Ohne Beschränkung der Allgemeinheit gelte
\begin{equation} \label{FA.6.ES.3}
\lim_{j \to \infty} x_j = x \mbox{ in } X_s.
\end{equation}
Wie zu Beginn von 4.) begründet man, daß $\overline{{\rm Span}_{\K} \{ x_j \, | \, j \in \N \}}^{\|\ldots\|}$ in $X_s$ abgeschlossen ist, also folgt
\begin{gather*}
x \in \overline{{\rm Span}_{\K} \{ x_j \, | \, j \in \N \}}^{\|\ldots\|}, \\
i_X(x) - \Phi \in \overline{{\rm Span}_{\K} \{ \Phi, i_X(x_j) \, | \, j \in \N \}}^{\|\ldots\|} \subset X'',
\end{gather*}
und (\ref{FA.6.ES.0}), (\ref{FA.6.ES.1}) implizieren
\begin{equation} \label{FA.6.ES.4}
\frac{\| i_X(x) - \Phi \|}{2} \le \sup \{  | (i_X(x) - \Phi)(F_{\kappa}) | \, | \, \kappa \in \N \}.
\end{equation}
Für alle $\kappa, j \in \N$ gilt allerdings
$$ | (i_X(x) - \Phi)(F_{\kappa}) | \le | (i_X(x_j) - \Phi)(F_{\kappa}) | + | (i_X(x_j) - i_X(x))(F_{\kappa}) |, $$ 
also folgt aus (\ref{FA.6.ES.2}), (\ref{FA.6.ES.3})  und (\ref{FA.6.ES.4}): $i_X(x) = \Phi$.
\q

\begin{Bem*} \label{FA.6.ES.B}
Der Autor dieser Seiten hat mit einer gewissen Überraschung festgestellt, daß die Beweisführung, der er folgt -- nämlich der in \cite{Die}, die wie o.g.\ auf \cite{Whit} zurückgeht --, an keiner Stelle erfordert, daß der betrachtete Raum ein $\K$-Banachraum ist; für einen normierten $\K$-Vektorraum ,,funktioniert`` der Beweis seines Erachtens ebenso.
Es ist aber selbstverständlich, daß der Autor nicht frei von Fehlern ist.
Für einen Hinweis auf solche ist er dankbar.
\end{Bem*}

\textit{Beweis des Satzes von \textsc{Eberlein-\v{S}mulian} \ref{FA.6.ES}.} Es sei $K$ eine Teilmenge von $X$.

,,$\Rightarrow$'' Ist $K$ kompakt in $X_s$, so ist $K$ nach \ref{FA.6.5} und \ref{FA.1.24} (ii) auch abgeschlossen in $X_s$, d.h.\ $K \subset \subset X_s$.
\ref{FA.6.32} ,,$\Rightarrow$`` ergibt folglich die relative Folgenkompaktheit von $K$ in $X_s$, und dieselbe Argumentation wie im Beweis von \ref{FA.1.9} ,,$\Leftarrow$`` zeigt wegen der Abgeschlossenheit von $K$ in $X_s$, daß $K$ folgenkompakt in $X_s$ ist.

,,$\Leftarrow$`` Sei $K$ abgeschlossen in $X_s$ und folgenkompakt in $X_s$.
Dann ist $K$ natürlich auch relativ folgenkompakt in $X_s$, also gilt nach \ref{FA.6.32} ,,$\Leftarrow$``: $K \subset \subset X_s$.
Hieraus und der Abgeschlossenheit von $K$ in $X_s$ folgt, daß $K$ eine kompakte Teilmenge von $X_s$ ist.
\q

\begin{Bem*}
Obwohl \textsc{W.\ F.\ Eberlein} dies nicht bewiesen hat, wird der Satz von \textsc{Eberlein-\v{S}mulian} in der Literatur gelegentlich folgendermaßen formuliert:

\textit{Eine Teilmenge eines $\K$-Banachraumes ist genau dann schwach kompakt, wenn sie schwach folgenkompakt ist.}

Der Autor vermutet, daß die letztgenannte Aussage falsch ist.
Um die obige Beweisführung aufrecht zu erhalten, müßte eine schwach folgenkompakte Teilmenge eines $\K$-Banachraumes stets schwach abgeschlossen sein.
Für die Nennung eines Gegenbeispieles wäre der Autor dankbar.
\end{Bem*}

\subsection*{Übungsaufgaben}

\begin{UA}
Beweise Lemma \ref{FA.6.7.L}.
\end{UA}

\begin{UA}
Seien $X$ ein normierter $\K$-Vektorraum, $(x_j)_{j \in \N}$ eine Folge in $X$ und $x \in X$.

Zeige
$$ \lim_{j \to \infty} x_j = x \mbox{ in } X_s \Longleftrightarrow \lim_{j \to \infty} i_X(x_j) = i_X(x) \mbox{ in } X''_{s*} \stackrel{\text{Def.}}{=} (X')'_{s*}. $$
\end{UA}

\begin{UA}
Seien $X$ ein normierter $\K$-Vektorraum, $(F_j)_{j \in \N}$ eine Folge in $X'$ und $F \in X$.

Zeige
$$ \lim_{j \to \infty} F_j = F \mbox{ in } X'_{s*} \Longrightarrow \|F\| \le \liminf_{j \to \infty} \|F_j\|. $$
\end{UA}

\begin{UA}
Seien $X$ ein normierter $\K$-Vektorraum, $(x_j)_{j \in \N}$ eine Folge in $X$ und $x \in X$ sowie $(F_j)_{j \in \N}$ eine Folge in $X'$ und $F \in X'$.
Ferner gelte entweder $\D \lim_{j \to \infty} x_j = x$ in $X$ und $\D \lim_{j \to \infty} F_j = F$ in $X'_{s*}$ oder $\D \lim_{j \to \infty} x_j = x$ in $X_s$ und $\D \lim_{j \to \infty} F_j = F$ in $X'$.

Zeige, daß dann $\D \lim_{j \to \infty} F_j(x_j) = F(x)$ gilt.
\end{UA}

\begin{UA}
Schwach-$*$-konvergente Folgen sind i.a.\ nicht beschränkt.
\end{UA}

\begin{UA}
Seien $X$ ein normierter $\K$-Vektorraum und $Y$ ein endlich-di\-men\-sio\-naler $\K$-Untervektorraum von $X$ sowie $a \in X$.

Zeige, daß die Approximationsaufgabe $A(a,Y)$ (mindestens) eine Lösung besitzt.
\end{UA}

\cleardoublepage
\section{Gleichmäßig konvexe Räume} \label{FAna7}

Wir weisen vor Beginn der Diskussion auf folgendes Lemma hin.

\begin{Lemma} \label{FA.7.1}
Sei $X$ ein normierter $\K$-Vektorraum.

Dann sind die folgenden Aussagen paarweise äquivalent:
\begin{itemize}
\item[(i)] $\forall_{\varepsilon \in \R_+} \exists_{\delta \in \R_+} \forall_{x,y \in B_1(0)} \, \left( \| x - y \| \ge \varepsilon \Longrightarrow \frac{\|x + y \|}{2} \le 1 - \delta \right)$.
\item[(ii)] $\forall_{\varepsilon \in \R_+} \exists_{\delta \in \R_+} \forall_{x,y \in X, \, \|x\| = \|y\| = 1} \, \left( \| x - y \| \ge \varepsilon \Longrightarrow \frac{\|x + y \|}{2} \le 1 - \delta \right)$.
\item[(iii)] $\forall_{\varepsilon \in \R_+} \exists_{\delta \in \R_+} \forall_{x,y \in X, \, \|x\| = \|y\| = 1} \, \left( \frac{\|x + y \|}{2} > 1 - \delta \Longrightarrow \|x - y \| < \varepsilon \right)$.
\item[(iv)] Für je zwei Folgen $(x_j)_{j \in \N}$, $(y_j)_{j \in \N}$ in $S_1(0) := \{ x \in X \, | \, \|x\| = 1\}$ gilt
$$ \lim_{j \to \infty} \frac{\|x_j + y_j \|}{2} = 1 \Longrightarrow \lim_{j \to \infty} (x_j - y_j) = 0. $$
\item[(v)] Für je zwei Folgen $(x_j)_{j \in \N}$, $(y_j)_{j \in \N}$ in $X$ mit
\begin{equation} \label{FA.7.1.0}
\limsup_{j \to \infty} \|x_j\| \le 1, ~~ \limsup_{j \to \infty} \|y_j\| \le 1
\end{equation}
gilt
$$ \lim_{j \to \infty} \frac{\|x_j + y_j \|}{2} = 1 \Longrightarrow \lim_{j \to \infty} (x_j - y_j) = 0. $$
\end{itemize}
\end{Lemma}

\textit{Beweis.} ,,(i) $\Rightarrow$ (ii)``, ,,(ii) $\Rightarrow$ (iii)`` und ,,(iii) $\Rightarrow$ (iv)`` sind trivial.

Zu ,,(iv) $\Rightarrow$ (v)``: Seien also $(x_j)_{j \in \N}$, $(y_j)_{j \in \N}$ Folgen in $X$ mit (\ref{FA.7.1.0}) und
\begin{equation} \label{FA.7.1.1}
\lim_{j \to \infty} \frac{\|x_j + y_j \|}{2} = 1.
\end{equation}
Zunächst gilt
\begin{equation} \label{FA.7.1.2}
\lim_{j \to \infty} \|x_j\| = \lim_{j \to \infty} \|y_j\| = 1.
\end{equation}

{[} Zu (\ref{FA.7.1.2}): Aus Symmetriegründen genügt es, $\lim_{j \to \infty} \|x_j\| = 1$ zu zeigen.
Aus der Annahme, dies sei falsch, und (\ref{FA.7.1.0}) folgte die Existenz einer Teilfolge $(x_{i_j})_{j \in \N}$ von $(x_j)_{j \in \N}$ mit $\lim_{j \to \infty} \|x_{i_j}\| < 1$, also gälte
$$ 2 \stackrel{(\ref{FA.7.1.1})}{=} \lim_{j \to \infty} \|x_{i_j} + y_{i_j} \| = \limsup_{j \to \infty} \|x_{i_j} + y_{i_j} \| \le \underbrace{\limsup_{j \to \infty} \|x_{i_j}\|}_{< 1} + \limsup_{j \to \infty} \|y_{i_j}\| \stackrel{(\ref{FA.7.1.0})}{<} 2, $$
Widerspruch! {]}

Daher existiert $j_0 \in \N$ derart, daß für alle $j \in \N$ mit $j \ge j_0$ gilt
\begin{gather*}
x_j \ne 0 \, \wedge \, y_j \ne 0, \\
\left| \, \left\| \frac{x_j}{\|x_j\|} + \frac{y_j}{\|y_j\|} \right\| - \|x_j + y_j \| \, \right| \stackrel{\ref{FA.2.2} (iii), (\rm N3)}{\le} \left\| \frac{x_j}{\|x_j\|} - x_j \right\| + \left\| \frac{y_j}{\|y_j\|} - y_j \right\| \stackrel{(\ref{FA.7.1.2})}{\stackrel{j \to \infty}{\longrightarrow} 0}
\end{gather*}
und somit
$$ \lim_{j_0 \le j \to \infty} \left\| \frac{x_j}{\|x_j\|} + \frac{y_j}{\|y_j\|} \right\| = \lim_{j_0 \le j \to \infty} \|x_j + y_j\| \stackrel{(\ref{FA.7.1.1})}{=} 2. $$
Hieraus und aus (iv) folgt
$$ \lim_{j_0 \le j \to \infty} \frac{x_j}{\|x_j\|} - \frac{y_j}{\|y_j\|} = 0, $$
also wegen
$$ \forall_{j \in \N, \, j \ge j_0} \, x_j - y_j = \underbrace{\left( \frac{x_j}{\|x_j\|} - \frac{y_j}{\|y_j\|} \right)}_{\stackrel{j \to \infty}{\longrightarrow} 0} - \underbrace{\left( \left( \frac{x_j}{\|x_j\|} - x_j \right) - \left( \frac{y_j}{\|y_j\|} - y_j \right) \right)}_{\stackrel{(\ref{FA.7.1.2})}{\stackrel{j \to \infty}{\longrightarrow} 0}} $$
die Konklusion von (v).

Zu ,,(v) $\Rightarrow$ (i)``: Angenommen, (i) gilt nicht.
Dann existieren $\varepsilon \in \R_+$ und Folgen $(x_j)_{j \in \N}$, $(y_j)_{j \in \N}$ in $B_1(0)$ mit
$$ \forall_{j \in \N} \, \| x_j - y_j \| \ge \varepsilon \, \wedge \, \frac{\|x_j + y_j \|}{2} > 1 - \frac{1}{j+1}, $$
im Widerspruch zu (v).
\q

\begin{Def}[Gleichmäßig und strikt konvexe Vektorräume] \label{FA.7.2} 
Sei $X$ ein normierter $\K$-Vektorraum.
\begin{itemize}
\item[(i)] $X$ heißt genau dann \emph{gleichmäßig konvex}\index{Raum!gleichmäßig konvexer Vektor-}, wenn eine der äquivalenten Aussagen (i) - (v) des letzten Lemmas gilt.
\item[(ii)] $X$ heißt genau dann \emph{strikt konvex}\index{Raum!strikt konvexer Vektor-}, wenn gilt
$$ \forall_{x,y \in X, \, \|x\| = \|y\| = 1} \, \left( x \ne y \Longrightarrow \frac{\|x+y\|}{2} < 1 \right). $$
\end{itemize}
\end{Def}

\begin{Satz} \label{FA.7.2.S1}
Ein gleichmäßig konvexer $\K$-Vektorraum ist strikt konvex.
\end{Satz}

\textit{Beweis.} Klar, denn \ref{FA.7.1} (ii) impliziert die Bedingung in \ref{FA.7.2} (ii). \q

\begin{Satz}[Charakterisierung der strikten Konvexität] \label{FA.7.2.S2}
Sei $X$ ein normierter $\K$-Vektorraum.

Dann sind die folgenden Aussagen paarweise äquivalent.
\begin{itemize}
\item[(i)] $X$ ist strikt konvex.
\item[(ii)] $X$ ist \emph{strikt normiert}\index{Raum!normierter Vektor-!strikt}, d.h.\ per definitionem
$$ \forall_{x,y \in X \setminus \{0\}} \left( \| x + y \| = \|x\| + \|y\| \Longrightarrow \exists_{\lambda \in \R_+} \, x = \lambda \, y \right). $$
\item[(iii)] $\forall_{x,y \in X} \, \left( (\|x\| = \|y\| = 1 \, \wedge \, x \ne y) \Longrightarrow \forall_{t \in {]}0,1{[}} \, \| (1-t) \, x + t \, y \| < 1 \right)$.
\end{itemize}
\end{Satz}

\textit{Beweis.} Zu ,,(i) $\Rightarrow$ (ii)``: Seien $x,y \in X \setminus \{0\}$ mit $\|x + y\| = \|x\| + \|y\|$.
Ohne Beschränkung der Allgemeinheit gelte $\|x\|=1$.
Setze $\tilde{x} := \frac{x+y}{1+\|y\|} = \frac{x+y}{\|x+y\|}$.
Dann gilt $\|x+y-\tilde{x}\| = \|y\|$ und somit $\|\frac{x-\tilde{x}}{2} + y\| \le \|y\|$, also folgt nacheinander
\begin{gather*}
1 + \|y\| = \|x + y\| = \left\| \frac{x+\tilde{x}}{2} + \frac{x-\tilde{x}}{2} + y \right\| \le \frac{\| x + \tilde{x} \|}{2} + \left\| \frac{x-\tilde{x}}{2} + y \right\| \le 1 + \|y\|, \\
\frac{\| x + \tilde{x} \|}{2} = 1.
\end{gather*}
Die letzte Gleichung und die strikte Konvexität von $X$ ergeben $x = \tilde{x}$, d.h.\ $x = \frac{1}{\|y\|} \, y$.

Zu ,,(ii) $\Rightarrow$ (iii)``: Seien $x,y \in X$ mit $\|x\|=\|y\|=1$ und $x\ne y$ sowie $t \in {]}0,1{[}$.
Gälte $\| (1-t) \, x + t \, y \| = 1 \stackrel{\text{klar}}{=} \|(1-t) \, x\| + \| t \, y \|$, so folgte aus (ii): $x=y$, Widerspruch!

Zu ,,(iii) $\Rightarrow$ (i)``: Setzt man $t:= \frac{1}{2}$ in (iii), so gilt (i).
\q

\begin{Bem*}
Die Einheitsvollkugel $B_1(0)$ eines normierten $\K$-Vektorraumes ist konvex, und sie (genauer gesagt ihr offener Kern $U_1(0)$) bestimmt die Norm nach \ref{FA.N.E} eindeutig.
Daher nennt man $B_1(0)$ auch den \emph{Eichkörper}\index{Eichkörper} des normierten $\K$-Vektorraumes.
Geometrisch präziser wäre es, in \ref{FA.7.2} (ii) von der ,,strikten Konvexität des Eichkörpers $B_1(0)$`` zu reden.
Sie besagt nach \ref{FA.7.2.S2} ,,(i) $\Rightarrow$ (iii)``, daß der Rand des Eichkörpers keine echten Strecken enthält, d.h., daß er überall ,,gekrümmt`` ist.
Falls der $\K$-Vektorraum sogar gleichmäßig konvex ist, so bezieht sich der Begriff der Gleichmäßigkeit auf das ,,Krümmungsverhalten`` des Randes des Eichkörpers:
Es soll überall eine gewisse ,,Mindestkrümmung`` vorliegen.
\end{Bem*}

\begin{Bsp} \label{FA.7.3} $\,$
\begin{itemize}
\item[1.)] Mit einem normierten $\K$-Vektorraum ist natürlich auch jeder seiner $\K$-Unter\-vek\-tor\-räume strikt konvex bzw.\ gleichmäßig konvex.
\item[2.)] Ein endlich-dimensionaler strikt konvexer $\K$-Vektorraum ist gleichmäßig konvex.

{[} Sei $X$ ein endlich-dimensionaler strikt konvexer $\K$-Vektorraum.
Versehe $X^2$ mit der Maximumsnorm.
Sei $\varepsilon \in \R_+$.
Dann ist
$$ K := \{ (x,y) \in X^2 \, | \, \|x\| = \|y\| = 1 \, \wedge \, \|x-y\| \ge \varepsilon \} $$
eine beschränkte und abgeschlossene (also kompakte) Teilmenge von $X^2$.
Die stetige Funktion
$$ X^2 \longrightarrow \R, ~~ (x,y) \longmapsto \frac{\|x+y\|}{2}, $$
nimmt daher auf $K$ ein Maximum an, welches wegen der strikten Konvexität von $X$ kleiner als Eins sein muß.
Hieraus folgt die gleichmäßige Konvexität von $X$. {]}
\item[3.)] Sei $n \in \N_+$.
Dann ist $(\R^n, \| \ldots \|_p)$ für jedes $p \in {]}1, \infty{[}$ gleichmäßig konvex.\footnote{Dies folgt auch aus dem u.g.\ Satz von \textsc{Clarkson} \ref{FA.7.C}.}
$(\R^n, \| \ldots \|_1)$ und $(\R^n, \| \ldots \|_{\infty})$ sind nicht strikt konvex.
\item[4.)] Auf $\mathcal{C}({[}0,1{]},\R)$ wird durch
$$ \forall_{f \in \mathcal{C}({[}0,1{]},\R)} \, \|f\| := \|f\|_{\infty} + \|f\|_2 $$
eine Norm definiert.
Es gilt:
\begin{itemize}
\item[(a)] $\| \ldots\|$ und $\| \ldots \|_{\infty}$ sind äquivalent und erzeugen somit dieselbe Topologie.
\item[(b)] $(\mathcal{C}({[}0,1{]},\R), \| \ldots \|_{\infty})$ ist nicht strikt konvex.
\item[(c)] $(\mathcal{C}({[}0,1{]},\R), \| \ldots \|)$ ist strikt konvex aber nicht gleichmäßig konvex. 
\end{itemize}

{[} (a) folgt i.w.\ aus \ref{FA.5.19} (i) mit $\varphi = \mu_1$, und (c) zeigen wir auf Seite \pageref{FA.7.Beweis} unter Verwendung des Satzes von \textsc{Clarkson} sowie des Satzes von \textsc{Milman}.

Zu (b): Betrachte die Funktionen $f_1, f_2 \: {[}0,1{]} \to \R$, definiert durch
$$ \forall_{t \in {[}0,1{]}} \, f_1(t) := 1 \, \wedge \, f_2(t) := t. $$
Dann folgt $\|f_1\|_{\infty} = \|f_2\|_{\infty} = 1$ und $\|f_1 + f_2 \|_{\infty} = 2$. {]}
\end{itemize}
\end{Bsp}

Wir kommen nun erneut auf die auf Seite \pageref{FA.6.Approx} gestellte \emph{Approximationsaufgabe} zu sprechen.

\begin{Satz} \label{FA.7.A} \index{Approximationsaufgabe} $\,$

\noindent \textbf{Vor.:} Seien $X$ ein normierter $\K$-Vekrorraum, $a \in X$ und $C$ eine nicht-leere abgeschlossene konvexe Teilmenge von $X$.

\noindent \textbf{Beh.:}
\begin{itemize}
\item[(i)] Ist $X$ strikt konvex, so besitzt $A(a,C)$ höchstens eine Lösung.
\item[(ii)] Ist $X$ gleichmäßig konvex und $C$ sogar vollständig, so besitzt $A(a,C)$ genau eine Lösung.
\end{itemize}
\end{Satz}

\textit{Beweis.} Im Falle $d := d(\{a\},C) = 0$ gilt $a \in C$ -- beachte, daß $C$ abgeschlossen ist --, und $A(a,C)$ besitzt die eindeutig bestimmte Lösung $a \in C$.
Ohne Beschränkung der Allgemeinheit gelte daher $d > 0$.

Zu (i): Angenommen, $x, y \in C$ sind zwei verschiedene Lösungen von $A(a,C)$, d.h.\ insbesondere $\| a - x \| = d = \| a - y \|$.
Dann sind $\tilde{x} := \frac{a - x}{d}$ sowie $\tilde{y} := \frac{a - y}{d}$ zwei verschiedene Einheitsvektoren, und die strikte Konvexität von $X$ ergibt $\| \frac{x+y}{2} - a \| = d \, \| \frac{\tilde{x} + \tilde{y}}{2} \| < d$, im Widerspruch zur Definition von $d$, da $C$ als konvexe Menge $\frac{x+y}{2}$ enthält.

Zu (ii): Sei $(y_j)_{j \in \N}$ eine Folge in $C$ mit $\lim_{j \to \infty} \|a - y_j\| = d$.
Wir behaupten:
\begin{equation} \label{FA.7.A.1}
(y_j)_{j \in \N} \mbox{ ist eine Cauchyfolge.}
\end{equation}

{[} Zu (\ref{FA.7.A.1}): Wegen $d>0$ gilt $a \notin C$, also auch $d_j := \|a - y_j \| > 0$.
Da $C$ konvex ist, gilt für alle $j,k \in \N$: $\frac{y_j + y_k}{2} \in C$, also
$$ d \le \left\| a - \frac{y_j + y_k}{2} \right\| = \left\| \frac{(a - y_j) + (a - y_k)}{2} \right\|. $$

Für jedes $j \in \N$ sei $x_j := \frac{a - y_j}{d_j} \in S_1(0) := \{  x \in X \, | \, \|x\| = 1 \}$.
Aus der letzten Ungleichung folgt dann für alle $j,k \in \N$
\begin{eqnarray}
1 & \le & \left\| \frac{(a - y_j) + (a - y_k)}{2 \, d} \right\| = \frac{1}{2} \, \left\| \frac{d_j}{d} \, x_j + \frac{d_k}{d} \, x_k \right\| \nonumber \\
& = & \frac{1}{2} \, \left\| (x_j + x_k) - \left( \left( 1 - \frac{d_j}{d} \right) \, x_j + \left( 1 - \frac{d_k}{d} \right) \, x_k \right) \right\| \label{FA.7.A.2} \\
& \le & \left\| \frac{x_j + x_k}{2} \right\| +  \frac{1}{2 \, d} \, ( |d - d_j| + |d - d_k| ). \nonumber
\end{eqnarray}

Hieraus ergibt sich zunächst, daß $(x_j)_{j \in \N}$ eine Cauchyfolge ist, denn zu $\varepsilon \in \R_+$ existiert wegen der gleichmäßigen Konvexität von $X$ eine Zahl $\delta \in \R_+$ mit
\begin{equation} \label{FA.7.A.3}
\forall_{x,y \in S_1(0)} \, \left( \left\| \frac{x+y}{2} \right\| > 1 - \delta \Longrightarrow \| x - y \| < \varepsilon \right),
\end{equation}
und wegen $\lim_{j \to \infty} d_j = 0$ existiert $j_0 \in \N$ derart, daß für alle $j,k \in \N$ mit $j,k \ge j_0$ gilt
$$ \frac{1}{2 \, d} \, ( |d - d_j| + |d - d_k| ) < \delta, $$
also nacheinander nach (\ref{FA.7.A.2})
$$ \left\| \frac{x_j + x_k}{2} \right\| > 1 - \delta $$
sowie (\ref{FA.7.A.3}): $\|x_j - x_k \|< \varepsilon$.

Da für alle $j,k \in \N$ gilt
\begin{eqnarray*}
\|y_j - y_k \| & = & \| (a - y_k) - (a - y_j) \| = \| d_k \, x_k - d_j \, x_j \| \\
& = & \| d \, (x_k + x_j) + ((d_k - d) \, x_k + (d_j - d) \, x_j) \| \\
& \le & d \, \|x_k + x_j\| + ( |d_k - d| + |d_j - d| ),
\end{eqnarray*}
ist $(y_j)_{j \in \N}$ mit $(x_j)_{j \in \N}$ wegen $\lim_{j \to \infty} d_j = d$ eine Cauchyfolge. {]}

Nun folgt aus (\ref{FA.7.A.1}) und der Vollständigkeit von $C$ die Existenz eines $x \in C$ mit $\lim_{j \to \infty} y_j = x$, d.h.\ auch $d = \lim_{j \to \infty} \|a - y_j\| = \|a - x\|$, also ist $x \in C$ eine Lösung von $A(a,C)$.
\ref{FA.7.2.S1} und (i) ergeben die Eindeutigkeit dieser Lösung.
\q
\A
Für den Rest dieses Kapitels liegt der Schwerpunkt erneut auf den Lebesgueschen Räumen.
Die wesentlichen Zutaten zum Nachweis der kommenden Hauptsätze haben wir bereits in den vorherigen Kapiteln vorbereitet.

\begin{HS}[Satz von \textsc{Clarkson}] \label{FA.7.C} \index{Satz!von \textsc{Clarkson}} $\,$

\noindent \textbf{Vor.:} Seien $n \in \N_+$, $\varphi$ ein Quadermaß auf $\R^n$, $M$ eine Teilmenge von $\R^n$ und $p \in {]}1, \infty{[}$.

\noindent \textbf{Beh.:} $L^p_{\K}(M,\varphi)$ ist gleichmäßig konvex.
\end{HS}

\textit{Beweis.} Es seien $(f_j)_{j \in \N}$ und $(g_j)_{j \in \N}$ Folgen in $\L^p_{\K}(M,\varphi)$ mit 
$$ \forall_{j \in \N} \, \|f_j\|_p = \|g_j\|_p = 1. $$
Aus den Ungleichungen von \textsc{Clarkson} \ref{FA.5.15} folgt für jedes $j \in \N$ im Falle $p \ge 2$
$$ {\|f_j + g_j\|_p}^p + {\|f_j -g_j\|_p}^p \le 2^p $$
und im Falle $p < 2$
$$ {\|f_j + g_j\|_p}^q + {\|f_j - g_j\|_p}^q \le 2^q $$
mit $q := \frac{1}{1 - \frac{1}{p}} \in {]} 2, \infty {[}$.

Gelte zusätzlich $\lim_{j \to \infty} \frac{\|f_j + g_j\|_p}{2} = 1$.
Dann implizieren die letzten beiden Ungleichungen $\lim_{j \to \infty} \|f_j - g_j\|_p = 0$, d.h., daß \ref{FA.7.1} (iv) in $\L_{\K}^p(M,\varphi)$ bzgl.\ der Halbnorm $\| \ldots \|_p$ und somit auch in $L_{\K}^p(M,\varphi)$ gilt.
Daher ist $L_{\K}^p(M,\varphi)$ gleichmäßig konvex. 
\q

\begin{Bem*}
$L^{\infty}_{\R}([0,1],\mu_1)$ ist nach \ref{FA.7.3} 1.) und 4.) (b) nicht strikt konvex und somit nicht gleichmäßig konvex.\footnote{An dieser Stelle muß der Leser sich natürlich überlegen, daß die Norm für $L^{\infty}_{\R}([0,1],\mu_1)$ auf der Menge der stetigen Funktionen $[0,1] \to \R$ mit der Supremumsnorm übereinstimmt.}
\end{Bem*}

Die Reflexivität eines normierten $\K$-Vektorraumes hängt nach \ref{FA.6.21} (iii) nur von der $\K$-Vektorraum-Isomorphieklasse ab, insbesondere handelt es sich um eine topologische Eigenschaft des normierten $\K$-Vektorraumes.
Der nun folgende Satz von \textsc{Milman} besagt, daß diese für einen $\K$-Banachraum aus der geometrischen Eigenschaft der gleichmäßigen Konvexität folgt.
Man beachte, daß die Topologie nur von der Äquivalenzklasse der Norm abhängt und diese Äquivalenzklasse gleichzeitig ,,gleichmäßig konvexe Normen`` und ,,nicht gleichmäßig konvexe Normen`` enthalten kann, vgl.\ \ref{FA.7.3} 3.). 

\begin{HS}[Satz von \textsc{Milman}] \label{FA.7.M} \index{Satz!von \textsc{Milman}} $\,$

\noindent \textbf{Vor.:} Sei $X$ ein gleichmäßig konvexer $\K$-Banachraum.

\noindent \textbf{Beh.:} $X$ ist reflexiv.
\end{HS}

Wir bereiten den Beweis des Hauptsatzes durch folgendes Lemma vor.

\begin{Lemma} \label{FA.7.4}
Es seien $X$ ein gleichmäßig konvexer $\K$-Vektorraum und $(x_k)_{k \in \N}$ ein Folge in $B_1(0)$ mit $\lim_{k,l \to \infty} \|x_k + x_l\| = 2$.\footnote{$\lim_{k,l \to \infty} \|x_k + x_l\| = 2 : \Longleftrightarrow \forall_{\varepsilon \in \R_+} \exists_{k_0 \in \N} \forall_{k,l \in \N, \, k,l \ge k_0} \, | \|x_k + x_l\| - 2 | < \varepsilon$.}

Dann gilt $\lim_{k \to \infty} \|x_k\| = 1$, und $(x_k)_{k \in \N}$ ist eine Cauchyfolge in $X$.
\end{Lemma}

\textit{Beweis.} Die erste Behauptung ist trivial.
Angenommen, es es gibt ein $\varepsilon \in \R_+$ derart, daß zu allen $k_0 \in \N$ natürliche Zahlen $k,l \in \N$ mit $k,l \ge k_0$ und $\|x_k - x_l \| \ge \varepsilon$ existieren.
Dann existieren Teilfolgen $(x_{i_k})_{k \in \N}$, $(x_{j_k})_{k \in \N}$ von $(x_k)_{k \in \N}$ mit
$$ \limsup_{k \to \infty} \| x_{i_k} \| \le 1, ~~ \limsup_{k \to \infty} \| x_{j_k} \| \le 1, ~~ \lim_{k \to \infty} \frac{\| x_{i_k} + x_{j_k} \|}{2} = 1  $$
und
$$ \liminf_{k \to \infty} \| x_{i_k} + x_{j_k} \| \ge \varepsilon, $$
im Widerspruch zur Bedingung \ref{FA.7.1} (iv) der gleichmäßigen Konvexität.
Daher ist $(x_k)_{k \in \N}$ eine Cauchyfolge in $X$. \q
\A
\textit{Beweis des Hauptsatzes.} Sei $\Phi \in X''$.
Zu zeigen ist, daß es ein $x \in X$ mit $i_X(x) = \Phi$, d.h.\ $\forall_{F_0 \in X'} \, F_0(x) = \Phi(F_0)$, gibt. 
Ohne Beschränkung der Allgemeinheit können wir zunächst $\Phi \ne 0$ (denn sonst leistet $x = 0$ das Gewünschte) und sodann $\|\Phi\| = 1$ annehmen.
Daher existiert eine Folge $(F_k)_{k \in \N_+}$ in $B' := \{ F \in X' \| \|F\| \le 1 \}$ mit $1 - \frac{1}{k+1} < |\Phi(F_k)|$ für jedes $k \in \N_+$.
Wegen $\forall_{k \in \N_+} \forall_{t \in {[}0, 2 \pi {[}} \, |\Phi(\e^{\i t} \, F_k)| = |\Phi(F_k)|$ kann die Folge so gewählt werden, daß sogar gilt
\begin{equation} \label{FA.7.M.1}
\forall_{k \in \N_+} \, 1 - \frac{1}{k+1} < \Phi(F_k) \in \R.
\end{equation}
Nach \ref{FA.6.22.L} existiert dann offenbar weiterhin eine Folge $(x_k)_{k \in \N_+}$ in $B_1(0)$ mit
\begin{equation} \label{FA.7.M.2}
\forall_{k \in \N_+} \forall_{\kappa \in \{1, \ldots, k\}} \, |\Phi(F_{\kappa}) - F_{\kappa}(x_k)| < \frac{1}{2 \, (k+1)}.
\end{equation}
Für alle $k \in \N_+$ und $\kappa \in \{1, \ldots, k\}$ folgt aus (\ref{FA.7.M.2}), (\ref{FA.7.M.1})
\begin{equation*} \label{FA.7.M.3}
\begin{aligned}
{\rm Re}(F_{\kappa}(x_k)) & = {\rm Re}(F_{\kappa}(x_k) - \Phi(F_{\kappa})) + {\rm Re}(\Phi(F_{\kappa})) \\
& > - \frac{1}{2 \, (k+1)} + 1 - \frac{1}{k+1} = 1 - \frac{3}{2 \, (k+1)},
\end{aligned}
\end{equation*}
also für jedes $l \in \N$ mit $l \ge k$
\begin{equation*} \label{FA.7.M.4}
\begin{aligned}
2 - \frac{3}{k+1} & \le 1 - \frac{3}{2 \, (k+1)} + 1 - \frac{3}{2 \, (l+1)} \\
& < {\rm Re}(F_{\kappa}(x_k) + F_{\kappa}(x_l))  \le |F_{\kappa}(x_k) + F_{\kappa}(x_l)| \le \|F_{\kappa}\| \, \|x_k + x_l\| \\
& \le \|x_k + x_l\| \le 2.
\end{aligned}
\end{equation*}
Das vorherige Lemma ergibt nun, daß $(x_k)_{k \in \N_+}$ eine Cauchyfolge in dem $\K$-Banach\-raum $X$ ist und $\lim_{k \to \infty} \|x_k\| = 1$ gilt.
$x \in S_1(0) := \{ \tilde{x} \in X \, | \, \| \tilde{x} = 1\}$ bezeichne den Grenzwert der Folge $(x_k)_{k \in \N_+}$.
Aus der Stetigkeit von $F_{\kappa}$ für $\kappa \in \N_+$ und (\ref{FA.7.M.2}) folgt daher
\begin{equation} \label{FA.7.M.S}
\forall_{\kappa \in \N_+} \, F_{\kappa}(x) = \Phi(F_{\kappa}).
\end{equation}

Wir behaupten, daß $x$ nach Wahl der Folge $(F_{\kappa})_{\kappa \in \N_+}$ eindeutig durch (\ref{FA.7.M.S}) und $\|x\| = 1$ bestimmt ist:
Erfüllt nämlich $x' \in S_1(0)$ ebenfalls (\ref{FA.7.M.S}) (mit $x'$ anstelle von $x$), so genügt die Folge $(\tilde{x}_k)_{k \in \N}$ in $S_1(0)$, definiert durch
$$ \tilde{x}_1 := x \, \wedge \, \forall_{k \in \N_+} \, (\tilde{x}_{2k} := x' \, \wedge \, \tilde{x}_{2k+1} := x), $$
der Ungleichung (\ref{FA.7.M.2}) und ist analog zu oben eine Cauchyfolge in dem $\K$-Banachraum $X$, also folgt $x = x'$.

Sei nun $F_0 \in X'$.
Dann existiert wie oben eine Folge $(x_k)_{k \in \N}$, deren Teilfolge $(x_k)_{k \in \N_+}$ sich ggf.\ von der ursprünglichen Folge $(x_k)_{k \in \N_+}$ unterscheidet, für die wieder
(\ref{FA.7.M.2}) und zusätzlich
\begin{equation} \label{FA.7.M.5}
\forall_{k \in \N} \, |\Phi(F_0) - F_0(x_k)| < \frac{1}{2 \, (k+1)}
\end{equation}
gilt.
Die bisherige Diskussion zeigt, daß auch die ,,neue`` Folge gegen (dasselbe) $x \in S_1(0)$ konvergiert, also folgt aus (\ref{FA.7.M.5}): $\Phi(F_0) = F_0(x)$.
\q

\begin{Bem*}
\textsc{Day} \cite{Day2} hat gezeigt, daß es separable reflexive strikt konvexe $\K$-Ba\-nach\-räume gibt, die zu keinem gleichmäßig konvexen $\K$-Vektorraum isomorph sind.
\end{Bem*}

Das nächste Korollar ergibt sich sofort aus dem Satz von \textsc{Clarkson} \ref{FA.7.C}, dem Satz von \textsc{Riesz-Fischer} \ref{FA.5.20} und dem Satz von \textsc{Milman} \ref{FA.7.M}.

\begin{Kor} \label{FA.7.5}
Es seien $n \in \N_+$, $\varphi$ ein Quadermaß auf $\R^n$, $M$ eine Teilmenge von $\R^n$ und $p \in {]}1, \infty{[}$.

Dann ist $L^p_{\K}(M,\varphi)$ ist reflexiv. \q
\end{Kor}

Wir können nun das o.g.\ Beispiel \ref{FA.7.3} 4.) (c) beweisen.
Dafür zeigen wir zunächst die folgenden beiden Lemmata.

\begin{Lemma} \label{FA.7.6}
Seien $X$ ein $\K$-Vektorraum und $\| \ldots \|_*$, $\| \ldots \|_s$ Normen auf $X$ derart, daß
$(X, \| \ldots \|_s)$ strikt konvex ist.

Dann ist $\| \ldots \| := \| \ldots \|_* + \| \ldots \|_s$ eine Norm auf $X$, und $(X, \| \ldots \|)$ ist ebenfalls strikt konvex.
\end{Lemma}

\textit{Beweis.} Seien $x,y \in X \setminus \{0\}$ $\K$-linear unabhängig.
Dann gilt nach Voraussetzung und \ref{FA.7.2.S2} ,,(i) $\Rightarrow$ (ii)``
$$ \| x + y \|_s < \| x + y \|_s, $$
also auch
$$ \| x + y \| = \| x + y \|_* + \| x + y \|_s < \|x\|_* + \|y\|_* + \|x\|_s + \|y\|_s = \|x\| + \|y\|. $$
Damit ist $(X, \| \ldots \|)$ nach \ref{FA.7.2.S2} (ii) $\Rightarrow$ (i) strikt konvex. \q

\begin{Lemma} \label{FA.7.7}
Der $\R$-Banachraum $(\mathcal{C}([0,1],\R),\|\ldots\|_{\infty})$ ist nicht reflexiv.
\end{Lemma}

\textit{Beweis.} Nach \ref{FA.3.B} ist der $\R$-Banachraum $\mathbf{c_0}_{\R}$ nicht reflexiv.
Wegen \ref{FA.6.21} (ii) genügt es, einen isometrischen $\K$-Vektorraum-Homomorphismus 
$$ T \: \mathbf{c_0}_{\R} \longrightarrow (\mathcal{C}([0,1],\R),\|\ldots\|_{\infty}), ~~ x \longmapsto T_x, $$
anzugeben.
Sei $x = (x_j)_{j \in \N} \in \mathbf{c_0}_{\R}$.
Setze dann 
$$ \forall_{j \in \N} \, T_x \left( \frac{1}{j+1} \right) := x_j ~~ \mbox{ sowie } ~~ T_x(0) := 0, $$ 
und erweitere dies zu einer stückweise linearen Funktion ${[}0,1{]} \to \R$, d.h.\ für alle $j \in \N$ sowie $t \in {]} \frac{1}{j+2}, \frac{1}{j+1} {[}$
$$ T_x(t) := T_x \left( \frac{1}{j+2} \right) + \left( t - \frac{1}{j+2} \right) (j+1) (j+2) \left( T_x \left( \frac{1}{j+1} \right) - T_x \left( \frac{1}{j+2} \right) \right). $$
Das so definierte $T$ leistet das Gewünschte. 
\q
\A
\textit{Beweis des Beispiels \ref{FA.7.2} 4.) (c).}\label{FA.7.Beweis} Aus dem Satz von \textsc{Clarkson} \ref{FA.7.C} und \ref{FA.7.3} 1.)\linebreak folgt, daß mit $L^2([0,1],\mu_1)$ auch $(\mathcal{C}({[}0,1{]},\R), \| \ldots \|_2)$ gleichmäßig konvex und somit strikt konvex ist.
Damit ist auch $(\mathcal{C}({[}0,1{]},\R), \| \ldots \|)$ nach \ref{FA.7.6} strikt konvex, beachte $\| \ldots \| = \| \ldots \|_{\infty} + \| \ldots \|_2$.
Des weiteren kann der letztgenannte Raum nicht gleichmäßig konvex sein, da $(\mathcal{C}({[}0,1{]},\R), \| \ldots \|_{\infty})$ sonst wegen \ref{FA.7.3} 4.) (a) und dem Satz von \textsc{Milman} \ref{FA.7.M} reflexiv wäre, im Widerspruch zu \ref{FA.7.7}. \q

\begin{HS}[Rieszscher Darstellungssatz für Lebesguesche Räume] \label{FA.7.RD} \index{Satz!Darstellungs- von \textsc{Riesz}} $\,$

\noindent \textbf{Vor.:} Seien $n \in \N_+$, $\varphi$ ein Quadermaß auf $\R^n$, $M$ eine $\varphi$-meßbare Teilmenge von $\R^n$ und $p \in {[} 1, \infty {[}$ sowie $q \in {]} 1, \infty {]}$ mit $\frac{1}{p} + \frac{1}{q} = 1$, wobei $\frac{1}{\infty}$ als $0$ zu lesen ist.

\noindent \textbf{Beh.:} Die Abbildung
$$ \Lambda \: L^q_{\K}(M,\varphi) \longrightarrow L^p_{\K}(M,\varphi)', ~~ \mathbf{f} \longmapsto \Lambda_{\mathbf{f}}, $$
wobei $\Lambda_{\mathbf{f}} \: L^p_{\K}(M,\varphi) \to \K$ für jedes $\mathbf{f} \in L^q_{\K}(M,\varphi)$ durch
$$ \forall_{\mathbf{g} \in L^p_{\K}(M,\varphi)} \, \Lambda_{\mathbf{f}} \, \mathbf{g} := \int_M f \, g \, \d \varphi \mbox{ mit beliebigen $f \in \mathbf{f}$ und $g \in \mathbf{g}$} $$
geben sei, ist eine $\K$-Vektor\-raum-Iso\-metrie.
\end{HS}

\textit{Beweis.} Wegen \ref{FA.5.17} bleibt nur zu zeigen, daß $\Lambda \: L^q_{\K}(M,\varphi) \to L^p_{\K}(M,\varphi)'$ surjektiv ist.

1.\ Fall: $p \in {]}1, \infty{[}$.
Wir können ohne Einschränkung $\varphi(M) > 0$ annehmen, denn andernfalls gilt $L^q_{\K}(M,\varphi) = \{0\}$ sowie $L^p_{\K}(M,\varphi)' = \{0\}$ und $\Lambda$ ist trivialerweise surjektiv.
\ref{FA.4.M6.K2} und \ref{FA.4.M4} ergeben dann $\forall_{k \in \N} \, M_k := M \cap {]}-k,k{[}^n \in \mathfrak{I}(\R^n,\varphi)$, also $\varphi(M_k) \in \R_+$ für hinreichend großes $k \in \N$.
Dies impliziert offenbar
\begin{equation} \label{FA.7.RD.0}
L^q_{\K}(M,\varphi) \ne \{0\}.
\end{equation}

Nach dem Satz von \textsc{Riesz-Fischer} \ref{FA.5.20} ist $L_{\K}^q(M,\varphi)$ und wegen \ref{FA.5.17} (ii) somit auch $\Lambda(L_{\K}^q(M,\varphi))$ vollständig.
Daher folgt aus \ref{FA.1.11} (ii) die Abgeschlossenheit von $\Lambda(L_{\K}^q(M,\varphi))$ in $L_{\K}^p(M,\varphi)'$.

Existierte $G \in L_{\K}^p(M,\varphi)' \setminus \Lambda(L_{\K}^q(M,\varphi))$, so gäbe es folglich wegen \ref{FA.2.20} (i)\linebreak ein $\Phi \in L_{\K}^p(M,\varphi)'' \setminus \{0\}$ mit
\begin{equation} \label{FA.7.RD.1}
\Phi|_{\Lambda(L_{\K}^q(M,\varphi))} = 0.
\end{equation}
Da $L_{\K}^p(M,\varphi)$ nach \ref{FA.7.5} reflexiv ist, existierte weiterhin ein $\mathbf{g} \in L_{\K}^p(M,\varphi) \setminus \{0\}$ derart, daß gälte $i_{L_{\K}^p(M,\varphi)}(\mathbf{g}) = \Phi$, d.h.\
\begin{equation} \label{FA.7.RD.2}
\forall_{\widetilde{G} \in L_{\K}^p(M,\varphi)'} \, \widetilde{G}(\mathbf{g}) = \Phi(\widetilde{G}),
\end{equation}
also auch für jedes $\mathbf{f} \in L_{\K}^q(M,\varphi)$
$$ 0 \stackrel{(\ref{FA.7.RD.1})}{=} \Phi(\Lambda_{\mathbf{f}}) \stackrel{(\ref{FA.7.RD.2})}{=} \Lambda_{\mathbf{f}} \, \mathbf{g} = \int_M f \, g \, \d \varphi \mbox{ mit $f \in \mathbf{f}$ und $g \in \mathbf{g}$}. $$
Wegen (\ref{FA.7.RD.0}) bedeutete dies $\mathbf{g} = 0$, Widerspruch!

2.\ Fall: $p = 1$, d.h.\ $q = \infty$.
Sei $G \: L^1_{\K}(M,\varphi) \to \K$ eine stetiges Funktional.
Für jedes $k \in \N$ ist $M_k := M \cap {]}-k,k{[}^n$ eine $\varphi$-integrierbare Teilmenge von $M$, also folgt $\varphi(M_k) < \infty$.
\ref{FA.5.19} (i) ergibt daher $|L_{\K}^2(M_k,\varphi)| \subset L_{\K}^1(M_k,\varphi)$ und die Stetigkeit der natürlichen Abbildung
\begin{equation*} \label{FA.7.RD.3}
\iota_k \: L^2_{\K}(M_k,\varphi) \longrightarrow L_{\K}^1(M_k,\varphi).
\end{equation*}
Des weiteren gibt es einen kanonischen isometrischen $\K$-Vek\-tor\-raum-Ho\-mo\-mor\-phis\-mus
\begin{equation} \label{FA.7.RD.4}
\mu_k \: L_{\K}^1(M_k,\varphi) \longrightarrow L_{\K}^1(M,\varphi) ~~ \mbox{ (mit $\|\mu_k\|= 1$)}
\end{equation}
denn
\begin{equation*} \label{FA.7.RD.5}
(\L_{\K}^1(M_k,\varphi), \| \ldots \|_1) \longrightarrow (\L_{\K}^1(M,\varphi), \| \ldots \|_1), ~~ g \longmapsto \chi_{M_k} \, \hat{g} = \chi_M \, ( \chi_{M_k} \, \hat{g} ),
\end{equation*}
ist eine ,,isometrische`` $\K$-lineare Abbildung -- beachte, daß $M$ eine $\varphi$-meß\-bare Teilmenge von $\R^n$ ist, und \ref{FA.4.70} (iii).
Daher ist auch
\begin{equation*} \label{FA.7.RD.6}
(G \circ \mu_k \circ \iota_k) \: L_{\K}^2(M_k,\varphi) \longrightarrow \K
\end{equation*}
ein stetiges Funktional, und aus dem bereits bewiesenen 1.\ Fall (angewandt auf $M_k$ anstelle von $M$ sowie $p=q=2$) folgt die Existenz eines $f_k \in \L_{\K}^2(M_k,\varphi)$ derart, daß gilt
\begin{equation} \label{FA.7.RD.7}
\forall_{\mathbf{g} \in L^2_{\K}(M_k,\varphi)} \, (G \circ \mu_k \circ \iota_k)(\mathbf{g}) = \int_{M_k} f_k \, g \, \d \varphi \mbox{ mit beliebigem $g \in \mathbf{g}$}. 
\end{equation}
Wir behaupten
\begin{gather}
\L^2_{\K}(M_{k+1},\varphi) \subset \L^2_{\K}(M_k,\varphi), \label{FA.7.RD.8} \\
\forall_{h \in \L^2_{\K}(M_{k+1},\varphi)} \, \hat{h}|_{M_k} \in \L^2_{\K}(M_k,\varphi), \label{FA.7.RD.9} \\
\widehat{f_{k+1}}|_{M_k} =_{\varphi} \widehat{f_k}|_{M_k}, \label{FA.7.RD.10} \\
f_k \in \L_{\K}^{\infty}(M_k, \varphi) \, \wedge \, \|f_k\|_{\infty} \le \|G\|. \label{FA.7.RD.13}
\end{gather}

{[} Zu (\ref{FA.7.RD.8}): Sei $h \in \L^2_{\K}(M_{k+1},\varphi)$.
Da $M_k$ insbesondere eine $\varphi$-meßbare Teilmenge von $\R^n$ ist, ergibt \ref{FA.4.70} (iii)
$$ \chi_{M_k} \, \widehat{|h|^2} = \chi_{M_k} \left( \chi_{M_{k+1}} \, \widehat{|h|^2} \right) \in \L_{\K}(\R^n,\varphi), $$
also $h \in \L^2_{\K}(M_k,\varphi)$.

Zu (\ref{FA.7.RD.9}): Sei $h \in \L^2_{\K}(M_{k+1},\varphi)$.
Dann gilt $h \in \L^2_{\K}(M_k,\varphi)$ nach (\ref{FA.7.RD.8}), d.h.\
$$ \chi_{M_k} \, \left| \widehat{\stackrel{\wedge}{h}|_{M_k}} \right|^2 = \chi_{M_k} \, \widehat{|h|^2} \in \L_{\K}(\R^n,\varphi), $$
und (\ref{FA.7.RD.9}) folgt.

Zu (\ref{FA.7.RD.10}): Sei $g \in \L^2_{\K}(M_{k+1},\varphi)$, also $\hat{g}|_{M_k} \in \L^2_{\K}(M_{k+1},\varphi) \subset \L^2_{\K}(M_k,\varphi)$ nach (\ref{FA.7.RD.9}) sowie (\ref{FA.7.RD.8}).
Bezeichnen $\pi_{\kappa}^{\tilde{p}} \: \L_{\K}^{\tilde{p}}(M_{\kappa},\varphi) \to L_{\K}^{\tilde{p}}(M_{\kappa},\varphi)$ für $\kappa \in \{k,k+1\}$, $\tilde{p} \in \{1,2\}$ und $\pi \: \L_{\K}^1(M, \varphi) \to L_{\K}^1(M, \varphi)$ die kanonischen Projektionen, so gilt
\begin{gather*}
\mu_{k+1}(\iota_{k+1}(\pi_{k+1}^2(\hat{g}|_{M_k}))) = \mu_{k+1}(\pi_{k+1}^1(\hat{g}|_{M_k})) =  \pi\Big(\overbrace{\chi_{M_{k+1}} \, \widehat{\stackrel{\wedge}{g}|_{M_k}}}^{= \, \chi_{M_k} \, \widehat{\stackrel{\wedge}{g}|_{M_k}}}\Big), \label{FA.7.RD.11.1} \\
\mu_{k}(\iota_{k}(\pi_{k}^2(\hat{g}|_{M_k}))) = \mu_{k}(\pi_{k}^1(\hat{g}|_{M_k})) = \pi\Big(\chi_{M_{k}} \, \widehat{\stackrel{\wedge}{g}|_{M_k}}\Big). \label{FA.7.RD.11.2}
\end{gather*}
Daher ergibt (\ref{FA.7.RD.7}) offenbar
\begin{equation*} \label{FA.7.RD.12}
\int\limits_{M_{k}} (f_{k+1} - f_{k}) \, g \, \d \varphi = 0,
\end{equation*}
und diese Gleichung gilt insbesondere für $g := \overline{f_{k+1} - f_{k}} \in \L^2_{\K}(M_{k+1},\varphi)$, also folgt aus \ref{FA.4.30}: $f_{k+1} =_{\varphi} f_k$ auf $M_k$.
Letzteres bedeutet genau (\ref{FA.7.RD.10}).  
\pagebreak

Zu (\ref{FA.7.RD.13}): Im Falle $f_k =_{\varphi} 0$ ist nichts zu zeigen.
Gelte daher ohne Einschränkung
\begin{equation*} \label{FA.7.RD.14}
\forall_{x \in M_k} \, f_k(x) \ne 0,
\end{equation*}
insbesondere ist $f$ auf $M$ definiert.
Dann ist wegen \ref{FA.4.69} (iv), (v) die \emph{,,Vorzeichenfunktion von $f_k$``}
\begin{equation} \label{FA.7.RD.15}
\sigma := \frac{f_k}{|f_k|} \: M_k \longrightarrow \K \mbox{ $\varphi$-meßbar über $M_k$ mit $\| \sigma \|_{\infty} = 1$.}
\end{equation}

Angenommen, für die nach \ref{FA.4.M6.K1} (ii), \ref{FA.4.M5} offenbar $\varphi$-meßbare Teilmenge
\begin{equation} \label{FA.7.RD.16}
\widetilde{M} := \{ x \in M_k \, | \, |f_k(x)| > \|G\| \}
\end{equation}
von $M_k$, die wegen $\varphi(M_k) < \infty$ sogar $\varphi$-integrierbar ist, gälte $\varphi(\widetilde{M}) > 0$.
Hieraus, (\ref{FA.7.RD.16}) und $f_k \in \L_{\K}^2(\widetilde{M},\varphi)$ -- beachte $\chi_{\widetilde{M}} \, \widehat{|f_k|}^2 = \chi_{\widetilde{M}} \, (\chi_{M_k} \, \widehat{|f_k|}^2)$ und \ref{FA.4.70} (iii) -- folgte die echte Ungleichung in der folgenden Abschätzung
\begin{eqnarray*}
\varphi(\widetilde{M}) \, \|G\| & = & \int_{\widetilde{M}} \|G\| \, \d \varphi \\
& < & \int_{\widetilde{M}} |f_k| \, \d \varphi \stackrel{(\ref{FA.7.RD.15})}{=} \int_{M_k} \chi_{\widetilde{M}} \, \frac{f_k}{\sigma} \, \d \varphi \\
& \stackrel{(\ref{FA.7.RD.7})}{\le} & \|G\| \, \|\mu_k\| \, \left\| \frac{\chi_{\widetilde{M}}}{\sigma} \right\|_1 \stackrel{(\ref{FA.7.RD.4}), (\ref{FA.7.RD.15})}{=} \|G\| \, \varphi(\widetilde{M}), 
\end{eqnarray*}
also ein Widerspruch. {]}

Wegen (\ref{FA.7.RD.10}) und (\ref{FA.7.RD.13}) existiert nun eine Funktion $f \in \L^{\infty}_{\K}(M,\varphi)$ derart, daß gilt
$$ \forall_{k \in \N} \, \hat{f}|_{M_k} =_{\varphi} \widehat{f_k}|_{M_k}. $$
Sei schließlich $g \in \L^1_{\K}(M,\varphi)$.
Dann folgt
\begin{gather*}
\forall_{k \in \N} \ g_k := \chi_{M_k} \, g = \chi_{M_k} \, (\chi_M \, g) \in \L^1_{\K}(M,\varphi), \\
\lim_{k \to \infty} g_k = g \mbox{ in } (\L^1_{\K}(M,\varphi),\|\ldots\|_1),
\end{gather*}
und $(f_k \, g_k)_{k \in \N}$ ist eine auf $M$ $\varphi$-konvergente Folge in $\L_{\K}^1(M,\varphi)$ -- beachte \ref{FA.5.8} --\linebreak mit $\forall_{k \in \N} \, |f_k \, g_k| \stackrel{(\ref{FA.7.RD.13})}{\le_{\varphi}} \|G\| \, |g| \in \L_{\K}^1(M,\varphi)$ sowie $\lim_{k \to \infty} f_k \, g_k = f \, g$ auf $M$.
Der Grenzwertsatz von \textsc{Lebesgue} ergibt nun
$$ G(\pi(g)) = \lim_{k \to \infty} G(\pi(g_k)) = \lim_{k \to \infty} \int_M f_k \, g_k \, \d \varphi = \int_M f \, g \, \d \varphi, $$
wobei $\pi \: \L^1_{\K}(M,\varphi) \to L^1_{\K}(M,\varphi)$ wieder den kanonischen Epimorphismus bezeichne.
Damit ist der Rieszsche Darstellungssatz für Lebesguesche Räume vollständig bewiesen. 
\q

\begin{Kor} \label{FA.7.8}
Die $\K$-Banachräume $\ell^1_{\K}$ und $\ell^{\infty}_{\K}$ sind nicht reflexiv.
\end{Kor}

\textit{Beweis.} Wir wissen bereits aus \ref{FA.3.B}, daß der $\K$-Banachraum $\mathbf{c_0}_{\K}$ nicht reflexiv und sein Dualraum isometrisch zu $\ell_{\K}^1$ ist.
Daher kann $\ell_{\K}^1$ nach \ref{FA.6.21} (iv), (iii) nicht reflexiv sein.
Da der $\K$-Banachraum $\ell_{\K}^{\infty}$ nach dem Rieszschen Darstellungssatz für Lebesguesche Räume isometrisch zu $(\ell_{\K}^1)'$ ist, sieht man analog, daß jener ebenfalls nicht reflexiv ist. \q
\A
Wir geben schließlich das in \ref{FA.6.6.B} angekündigte Beispiel an.

\begin{Bsp} \label{FA.6.6.B.B}
Für jedes $k \in \N$ bezeichne $e_k := (e_{k,j})_{j \in \N} \in \ell^2_{\R}$ die Folge mit $\forall_{j \in \N} \, e_{k,j} = \delta_{kj}$.
Dann ist $(e_k)_{k \in \N}$ wegen $\forall_{k,l \in \N, \, k \ne l} \, \|e_k - e_l\|_2 = \sqrt{2}$ keine Cauchyfolge in $\ell^2_{\R}$ und somit nicht konvergent.
Allerdings existiert nach dem Rieszschen Darstellungssatz zu jedem $F \in (\ell^2_{\R})'$ eine Folge $(x_j)_{j \in \N} \in \ell^2_{\R}$ derart, daß gilt
$$ \forall_{(y_j)_{j \in \N} \in \ell^2_{\R}} \, F((y_j)_{j \in \N}) = \sum_{j=0}^{\infty} x_j \, y_j, $$
also $\forall_{k \in \N} \, F(e_k) = x_k$.
Da die Elemente von $\ell^2_{\R}$ insbesondere Nullfolgen sind, erhalten wir $\lim_{k \to \infty} F(e_k) = 0$, und aus \ref{FA.6.6} sowie der Beliebigkeit von $F$ folgt, daß $(e_k)_{k \in \N}$ in $(\ell^2_{\R})_s$ gegen Null konvergiert.
\end{Bsp}

\subsection*{Übungsaufgabe}

\begin{UA}
Sei $n \in \N_+$.

Zeige, daß $\Lambda \: L_{\K}^1(\R^n,\mu_n) \to L_{\K}^{\infty}(\R^n,\mu_n)'$ nicht surjektiv ist.
\end{UA}

Tip: Gemäß dem topologischen Fortsetzungssatz \ref{FA.2.19} von \textsc{Hahn-Banach} existiert eine Erweiterung des stetigen (!) Funktionals
$$ \delta \: \mathcal{C}_b(\R^n,\K) \longrightarrow \K, ~~ f \longmapsto f(0), $$
zu $\Delta \in L_{\K}^{\infty}(\R^n,\mu_n)'$ mit $\Delta \ne 0$ (!).
Angenommen, es existierte $f \in \L^1(\R^n,\mu_n)$ mit
$$ \forall_{\mathbf{g} \in L_{\K}^{\infty}(\R^n,\mu_n)'} \forall_{g \in \mathbf{g}} \, \Delta(\mathbf{g}) = \int_{\R^n} f \, g \, \d \mu_n, $$
so folgte insbesondere
$$ \forall_{g \in \mathcal{C}_b(\R^n,\K)} \, g(0) = \delta(g) = \int_{\R^n} f \, g \, \d \mu_n, $$
also auch $f =_{\mu_n} 0$ (!), im Widerspruch zu $\Delta \ne 0$.

\cleardoublepage
\section{Hilberträume} \label{FAna8}
\subsection*{Semiskalarprodukte und induzierte Halbnormen} \addcontentsline{toc}{subsection}{Semiskalarprodukte und induzierte Halbnormen}

\begin{Def}[Sesquilinearformen, hermitesche Sesquilinearformen, Semi\-ska\-lar\-pro\-dukte und Skalarprodukte] \label{FA.8.1}
Sei $X$ ein $\K$-Vektorraum.
\begin{itemize}
\item[(i)] Wir nennen eine Abbildung $s \: X \times X \to \K$ eine \emph{Sesquilinearform auf $X$}\index{Sesquilinearform} genau dann, wenn für alle $\lambda \in \K$ und $x,\tilde{x},y,\tilde{y} \in X$ gilt
$$
\begin{array}{lcl}
s( \lambda \, x , y ) = \lambda \, s( x , y ), && s( x + \tilde{x} , y ) = s( x , y ) + s( \tilde{x} , y ), \\
s( x , \lambda \, y ) = \overline{\lambda} \, s( x , y ) & \mbox{und} & s( x , y + \tilde{y} ) = s( x , y ) + s( x , \tilde{y} ),
\end{array}
$$
d.h.\ insbesondere $s(0,0) = 0$.

\begin{Bem*} $\,$
\begin{itemize} 
\item[1.)] Im Falle $\K = \R$ sind die Sesquilinearformen auf $X$ wegen $\overline{\lambda} = \lambda$ für $\lambda \in \R$ genau die $\R$-Bilinearformen auf $X$ und im Falle $\K = \C$ spezielle $\R$-Bilinearformen auf $X$.
\item[2.)] In der physikalischen Literatur, also z.B.\ in der Quantenmechanik, ist es üblich, eine Sequilinearform als linear im zweiten und semilinear im ersten Argument zu definieren.
\end{itemize}
\end{Bem*}
Die Menge aller Sesquilinearformen auf $X$ bildet in kanonischer Weise einen $\K$-Vektorraum, den wir mit $\boxed{{\rm SQF}(X)}$ bezeichnen.
\item[(ii)] Eine \emph{hermitesche Sesquilinearform auf $X$}\index{Sesquilinearform!hermitesche} ist per definitionem eine Sesquilinearform $s \: X \times X \to \K$ auf $X$ derart, daß gilt
$$ \forall_{x,y \in X} \, s( x , y ) = \overline{s( y , x )}, $$
d.h.\ insbesondere $\forall_{x \in X} \, s( x , x ) \in \R$.

\begin{Bem*}
Im Falle $\K = \R$ ist eine hermitesche Sesquilinearform auf $X$ also eine symmetrische $\R$-Bilinearform auf $X$ und umgekehrt.
\end{Bem*}
Die Menge aller hermiteschen Sesquilinearformen auf $X$ stellt einen $\K$-Unter\-vek\-tor\-raum von ${\rm SQF}(X)$ dar, den wir mit $\boxed{{\rm HF}(X)}$ bezeichnen.
\item[(iii)] Eine hermitesche Sesquilinearform $\langle \ldots, \ldots \rangle \: X \times X \to \K$ auf $X$ heißt \emph{positiv semi-definit} bzw.\ \emph{positiv definit} genau dann, wenn gilt
$$ \forall_{x \in X} \, \langle x , x \rangle \ge 0 ~~ \mbox{ bzw. } ~~ \forall_{x \in X \setminus \{0\}} \, \langle x , x \rangle > 0. $$
Wir nennen $\langle \ldots, \ldots \rangle$ in diesem Falle auch ein \emph{Semiskalarprodukt auf $X$}\index{Skalarprodukt!Semi-} bzw.\ ein \emph{Skalarprodukt}\index{Skalarprodukt} oder \emph{inneres Produkt auf $X$}.

\begin{Bem*}
Der Begriff ,,Semiskalarprodukt auf $X$`` wird in der Literatur auch allgemeiner für \emph{semi-innere Produkte auf $X$} verwendet.
Sie wurden 1961 von \textsc{G.\ Lumer} \cite{Lumer} eingeführt und sind per definitionem im ersten Argument $\K$-lineare Abbildungen $\langle \ldots, \ldots \rangle \: X \times X \to \K$, die positiv definit sind und die Cauchy-Schwarzsche Ungleichung, vgl.\ (\ref{FA.8.CS.S}) unten, erfüllen.
Zur Abgrenzung sei auch noch erwähnt, daß man im Falle $\K = \R$ unter einem \emph{Pseudoskalarprodukt auf $X$} eine symmetrische $\R$-Bilinearform $\langle \ldots, \ldots \rangle \: X \times X \to \R$, die nicht-entartet ist (d.h.\ genau $\forall_{x \in X} \, (\langle x, \ldots \rangle = 0 \Rightarrow x = 0)$, versteht.
\end{Bem*}

\begin{Bsp*} Sei $n \in \N_+$.
\begin{itemize}
\item[1.)] 
Seien $\varphi$ ein Quadermaß auf $\R^n$ und $M$ eine Teilmenge von $\R^n$.
Dann definiert
$$ \forall_{f,g \in \L^2_{\K}(M,\varphi)} \, \boxed{\langle f , g \rangle_2} := \int_M f \, \overline{g} \, \d \varphi $$
ein Semiskalarprodukt auf $\L^2_{\K}(M,\varphi)$, das in kanonischer Weise ein Skalarprodukt $\boxed{\langle \ldots, \ldots \rangle_2}$ auf $L^2_{\K}(M,\varphi)$ induziert.
Wir versehen diese beiden Räume im folgenden stets mit jenen Semiskalarprodukten.

\item[2.)] Ist in 1.) speziell $M := \{1, \ldots, n\} \subset \R$ (das dortige $n$ ist also gleich Eins) und $\varphi := \varphi_m$, wobei $\varphi_m$ das Quadermaß der diskreten Massenverteilung $m := 1_{\{1, \ldots,n\}} \: \{1,\ldots,n\} \to \R_+$ bezeichne, vgl.\ Beispiel \ref{FA.4.4.B} 4.), so können wir $L^2_{\K}(\{1, \ldots, n\}, \varphi_m)$ mit $\K^n$ identifizieren, siehe ggf.\ Beispiel \ref{FA.5.6.Bsp} 1.), und erhalten das kanonische Skalarprodukt auf $\K^n$, das durch
$$ \forall_{x=(x_1,\ldots,x_n), y=(y_1,\ldots,y_n) \in \K^n} \, \langle x, y \rangle_2 = \sum_{i=1}^n x_i \, \overline{y_i} $$
gegeben ist.
\item[3.)] Analog zu 2.) erhält man aus 1.) mit $M := \N$ und $m := 1_{\N} \: \N \to \R_+$ sowie $\varphi := \varphi_m$ das folgende Skalarprodukt auf $\ell^2_{\K}$:
$$ \forall_{x = (x_k)_{k \in \N}, y = (y_k)_{k \in \N} \in \ell_{\K}^2} \, \langle x, y \rangle_2 = \sum_{i=1}^{\infty} x_i \, \overline{y_i}. $$
\item[4.)] Informativ erwähnen wir, daß auf $\R^4$ durch
$$ \forall_{x=(x_1,\ldots,x_4), y=(y_1,\ldots,y_4) \in \R^4} \, \langle x, y \rangle := - x_1 \, y_1 + x_2 \, y_2 + x_3 \, y_3 + x_4 \, y_4 $$
ein Pseudoskalarprodukt definiert wird, das in der \emph{Speziellen Relativitätstheorie} von Bedeutung ist.
\end{itemize}
\end{Bsp*}
\end{itemize}
\end{Def}

\begin{Satz} \label{FA.8.2}
Sind $X$ ein $\K$-Vektorraum und $\langle \ldots, \ldots \rangle$ ein Semiskalarprodukt auf $X$ bzw.\ ein Skalarprodukt auf $X$, so definiert
\begin{equation} \label{FA.8.2.S}
\boxed{\forall_{x \in X} \, \|x\| := \sqrt{\langle x , x \rangle}}
\end{equation}
eine Halbnorm auf $X$ bzw.\ eine Norm auf $X$, die sog.\ \emph{durch $\langle \ldots, \ldots \rangle$ induzierte Halbnorm} bzw.\ \emph{Norm auf $X$}.
\end{Satz}

Zum Nachweis des Satzes benötigen wir die folgende (wichtige) Ungleichung.

\begin{Lemma}[Cauchy-Schwarzsche Ungleichung] \index{Ungleichung!Cauchy-Schwarzsche} \label{FA.8.CS}
Seien $X$ ein $\K$-Vek\-tor\-raum, und $\langle \ldots, \ldots \rangle$ ein Semiskalarprodukt auf $X$.
\begin{itemize}
\item[(i)] Für alle $x,y \in X$ gilt
\begin{equation} \label{FA.8.CS.S}
| \langle x, y \rangle |^2 \le \langle x, x \rangle \, \langle y, y \rangle
\end{equation}
mit Gleichheit, falls $x$ und $y$ $\K$-linear abhängig sind.
\item[(ii)] Ist $\langle \ldots, \ldots \rangle$ sogar ein Skalarprodukt, so gilt für $x,y \in X$ in (\ref{FA.8.CS.S}) genau dann Gleichheit, wenn $x$ und $y$ $\K$-linear abhängig sind.
\end{itemize}
\end{Lemma}

\textit{Beweis.} Zu (i): Seien $x,y \in X$ und $\lambda \in \K$.
Da $\langle \ldots, \ldots \rangle$ ein Semiskalarprodukt ist, folgt
\begin{equation} \label{FA.8.CS.0}
\begin{aligned}
0 \le \langle x + \lambda \, y, x + \lambda \, y \rangle & = \langle x, x \rangle + \overline{\lambda} \, \langle x, y \rangle + \lambda \, \langle y, x \rangle + \lambda \, \overline{\lambda} \, \langle y, y \rangle \\
& = \langle x, x \rangle + \overline{\lambda} \, \langle x, y \rangle + \lambda \, \overline{\langle x, y \rangle} + |\lambda|^2 \, \langle y, y \rangle.
\end{aligned}
\end{equation} 
Hieraus erhält man im Falle $\langle y , y \rangle \ne 0$ mit $\lambda := - \frac{\langle x , y \rangle}{\langle y , y \rangle}$
\begin{gather*}
0 \le \langle x, x \rangle - \frac{|\langle x , y \rangle|^2}{\langle y , y \rangle} - \frac{|\langle x , y \rangle|^2}{\langle y , y \rangle} + \frac{|\langle x , y \rangle|^2}{\langle y , y \rangle} = \langle x , x \rangle - \frac{|\langle x , y \rangle|^2}{\langle y , y \rangle},
\end{gather*}
und dies bedeutet genau (\ref{FA.8.CS.S}).

Durch Vertauschung der Rollen von $x$ und $y$ sieht man (\ref{FA.8.CS.S}) auch im Falle $\langle x , x \rangle \ne 0$ ein. 

Gilt schließlich $\langle x , x \rangle = \langle y , y \rangle = 0$, so ergibt (\ref{FA.8.CS.0}) mit $\lambda := - \langle x , y \rangle$, daß $\langle x, y \rangle = 0$ gilt, also ebenfalls die zu zeigende Ungleichung (\ref{FA.8.CS.S}).

Daß für beliebige $\K$-linear abhängige $x,y \in \K$ in (\ref{FA.8.CS.S}) Gleichheit gilt, ist klar.

Zu (ii): Seien $x,y \in \K$ derart, daß
\begin{equation} \label{FA.8.CS.1}
|\langle x, y \rangle|^2 = \langle x, x \rangle \, \langle y, y \rangle.
\end{equation}
Im Falle $\langle x, y \rangle = 0$ folgt $x = 0$ oder $y = 0$, weil $\langle \ldots, \ldots \rangle$ ein Skalarprodukt ist, und $x$ und $y$ sind $\K$-linear abhängig.
Gelte daher $\langle x, y \rangle \ne 0$, also insbesondere $y \ne 0$, d.h.\ wiederum, da $\langle \ldots, \ldots \rangle$ ein Skalarprodukt ist: $\langle y, y \rangle \ne 0$.
Mit $\lambda := - \frac{\langle x, y \rangle}{\langle y, y \rangle} \ne 0$ ist dann die rechte Seite von (\ref{FA.8.CS.0}) wegen (\ref{FA.8.CS.1}) gleich Null, also ergibt (\ref{FA.8.CS.0})
$$ \langle x + \lambda \, y , x + \lambda \, y \rangle = 0. $$
Da $\langle \ldots, \ldots \rangle $ ein Skalarprodukt ist, folgt hieraus $x + \lambda \, y = 0$, d.h.\ $x$ und $y$ sind $\K$-linear abhängig.
Wegen (i) genügt dies zu Nachweis von (ii).
\q
\A
\textit{Beweis des Satzes.} Sei $\langle \ldots, \ldots \rangle$ ein Semiskalarprodukt auf $X$.
Dann ist $\| \ldots \| \: X \to {[}0, \infty{[}$ durch (\ref{FA.8.2.S}) wohldefiniert, und (N2) ist trivial.

Zu (N3): Es seien $x,y \in X$.
Dann gilt
\begin{equation*}
\begin{aligned}
\| x + y \|^2 & = \langle x + y, x + y \rangle = \langle x, x \rangle + \langle x, y \rangle + \langle y, x \rangle + \langle y, y \rangle \\
& = \langle x, x \rangle + \langle x, y \rangle + \overline{\langle x, y \rangle} + \langle y, y \rangle = \|x\|^2 + 2 \, {\rm Re}(\langle x, y \rangle) + \|y\|^2 \\
& \le \|x\|^2 + 2 \, |\langle x, y \rangle| + \|y\|^2  \stackrel{\ref{FA.8.CS}}{\le} \|x\|^2 + 2 \, \|x\| \, \|y\| + \|y\|^2 = ( \|x\| + \|y \| )^2,
\end{aligned}
\end{equation*}
und die Dreiecksungleichung folgt.

Ist $\langle \ldots, \ldots \rangle$ sogar ein Skalarprodukt, so ist auch (N1) klar. \q
\pagebreak

\begin{Satz} \label{FA.8.5}
Es seien $X$ ein $\K$-Vektorraum und $\langle \ldots, \ldots \rangle$ ein Semiskalarprodukt auf $X$.

Dann gilt für die durch $\langle \ldots, \ldots \rangle$ induzierte Halbnorm $\| \ldots \|$ auf $X$
\begin{itemize}
\item[(i)] (Satz des \textsc{Pythagoras})\index{Satz!des \textsc{Pythagoras}}
$$ \forall_{x,y \in X} \, \|x+y\|^2 = \|x\|^2 + 2 \, {\rm Re}(\langle x, y \rangle) + \|y\|^2, $$
\item[(ii)] (Parallelogrammgleichung)\index{Parallelogramm!-gleichung}\index{Gleichung!Parallelogramm-}
$$ \forall_{x,y \in X} \, \|x+y\|^2 + \|x-y\|^2 = 2 \, \|x\|^2 + 2 \, \|y\|^2. $$
\end{itemize}
\end{Satz}

\textit{Beweis.} (i) haben wir beim Beweis des letzten Satzes mitbewiesen, und (ii) zeige der Leser als einfache Übung. \q

\begin{Bem*} $\,$
\begin{itemize}
\item[1.)] Teil (i) des letzten Satzes den \emph{Satz des \textsc{Pythagoras}} zu nennen, ist dadurch gerechtfertigt, daß zwei Elemente $x,y$ eines $\K$-Vektorraumes $X$, auf dem ein Semiskalarprodukt $\langle \ldots, \ldots \rangle$ gegeben ist, \emph{zueinander orthogonal}\index{Orthogonal} heißen (i.Z.\ $x \perp y$), wenn $\langle x, y \rangle = 0$ gilt, denn dann ergibt sich
$$ \forall_{x,y \in X} \, \left( x \perp y \Longrightarrow \|x+y\|^2 = \|x\|^2 + \|y\|^2 \right). $$
\item[2.)] Anschaulich besagt die Parallelogrammgleichung, daß die Summe der Quadrate der Seiten eines Parallelogrammes gleich der Summe der Quadrate der Diagonalen ist.
\end{itemize}
\end{Bem*}

Wir haben jedem Semiskalarprodukt, das auf einem vorgegebenen $\K$-Vek\-tor\-raum $X$ definiert ist, eine Funktion, nämlich die induzierte Halbnorm, zugeordnet.
In ähnlicher Weise werden wir nun jeder Sesquilinearform $s$ auf $X$ eine Funktion $Q_s$, bei der es sich im Falle $\K = \R$ um eine quadratische und im Falle $\K = \C$ um eine hermitesche Form handelt, zuordnen, welche insbesondere im Falle einer hermiteschen Sesquilinearform $s$ bereits die gesamte Information von $s$ enthält.

\begin{Satz} \label{FA.8.3}
Sei $X$ ein $\K$-Vektorraum.

Dann ist
$$ Q \: {\rm SQF}(X) \longrightarrow \K^X, ~~ s \longmapsto (Q_s \: X  \to \K, \, x \mapsto s(x,x)), $$
eine $\K$-lineare Abbildung, die im Falle $\K = \C$ injektiv ist, und deren Einschränkung auf ${\rm HF}(X)$ im Falle $\K = \R$ injektiv ist.
Für jedes $s \in {\rm SQL}(X)$ gilt nämlich im Falle $\K = \C$ die \emph{Polarisierungsformel}\index{Polarisierungsformel}
\begin{align*}
\forall_{x,y \in X} \, s(x,y) & = \frac{1}{4} \sum_{k=0}^3 \i^k \, Q_s(x + \i^k y) \\ & = \frac{1}{4} \left( Q_s(x+y) -  Q_s(x-y) + \i (Q_s (x + \i \, y) - Q_s(x - \i \, y)) \right)
\end{align*}
und für jedes $s \in {\rm HF}(X)$ im Falle $\K = \R$ die \emph{Polarisierungsformel}\index{Polarisierungsformel}
$$ \forall_{x,y \in X} \, s(x,y) = \frac{1}{4} \left( Q_s(x+y) -  Q_s(x-y) \right). $$
\end{Satz}

\textit{Beweis als Übung.} \q
\pagebreak

\begin{Kor} \label{FA.8.4}
Seien $X$ ein $\K$-Vektorraum und $\| \ldots \|$ eine Halbnorm auf $X$.

Dann existiert höchstens ein Semiskalarprodukt $\langle \ldots, \ldots \rangle$ auf $X$, das die Halbnorm $\| \ldots \|$ induziert, und im Falle der Existenz ist $\langle \ldots, \ldots \rangle$ durch $\| \ldots \|$ mittels der Polarisierungsformel bestimmbar -- beachte $Q_{\langle \ldots, \ldots \rangle} = \| \ldots \|^2$. \q
\end{Kor}

\subsection*{Prähilberträume} \addcontentsline{toc}{subsection}{Prähilberträume}

\begin{Def}[Prähilberträume] \index{Raum!Hilbert-!Prä-} \index{Hilbert!-raum!Prä-} \label{FA.8.6}
Ein Paar $(X, \langle \ldots, \ldots \rangle)$, für das wir auch kurz $X$ schreiben, bestehend aus einem $\K$-Vektorraum $X$ und einem Skalarprodukt $\langle \ldots , \ldots \rangle$ auf $X$, heißt ein \emph{$\K$-Prähilbertraum}.

\begin{Bem*}
Wir betrachten einen $\K$-Prähilbertraum im folgenden stets auch als einen normierten $\K$-Vektorraum mit der durch das Skalarprodukt induzierten Norm.
\end{Bem*}
\end{Def}

\begin{Lemma} \label{FA.8.7}
Ist $X$ ein $\K$-Prähilbertraum, so ist
$$ \langle \ldots, \ldots \rangle \: X \times X \longrightarrow \K $$
eine stetige Funktion, wobei wir $X \times X$ mit der Maximumsnorm versehen.
\end{Lemma}

\textit{Beweis.} Seien $x,y \in X$ und $(x_k)_{k \in \N}$ bzw.\ $(y_k)_{k \in \N}$ Folgen in $X$, die gegen $x$ bzw.\ $y$ konvergieren -- beachte, daß Konvergenz bzgl.\ der Maximumsnorm gleichbedeutend mit der Konvergenz der Komponenten ist.
Dann gilt für $k \in \N$ nach der Cauchy-Schwarzschen Ungleichung \ref{FA.8.CS}
\begin{align*}
\left| \langle x, y \rangle - \langle x_k, y_k \rangle \right| & = \left| \langle x - x_k, y \rangle - \langle x_k, y_k - y \rangle \right| \le \left| \langle x - x_k, y \rangle \right| + \left| \langle x_k, y_k - y \rangle \right|\\
& \le \| x - x_k \| \, \|y\| + \| x_k \| \, \| y_k - y \| \stackrel{k \to \infty}{\longrightarrow} 0,
\end{align*}
und die Behauptung ist gezeigt. \q

\begin{Kor} \label{FA.8.7.K}
Seien $X$ ein $\K$-Prähilbertraum und $y \in X$.

Dann ist $\langle \ldots, y \rangle \: X \to \K$ eine stetige $\K$-Linearform. \q
\end{Kor}

\begin{Satz} \label{FA.8.8}
Sei $X$ ein normierter $\K$-Vektorraum.

Dann ,,ist`` $X$ genau dann ein $\K$-Prähilbertraum, d.h.\ per definitionem, daß ein Skalarprodukt auf $X$ existiert, das die Norm von $X$ induziert, wenn die in \ref{FA.8.5} (ii) genannte \emph{Parallelogrammgleichung}\index{Parallelogramm!-gleichung}\index{Gleichung!Parallelogramm-} gilt:
$$ \forall_{x,y \in X} \, \|x+y\|^2 + \|x-y\|^2 = 2 \, \|x\|^2 + 2 \, \|y\|^2. $$
\end{Satz}

\textit{Beweisskizze des Satzes.} ,,$\Rightarrow$`` ist bereits nach \ref{FA.8.5} (ii) klar.

,,$\Leftarrow$`` Wir skizzieren den in \cite[S.\ 85]{Hirz} angegebenen Beweis.
Zu zeigen ist, daß durch die Polarisationsformel in \ref{FA.8.3} ein Skalarprodukt $s = \langle \ldots, \ldots \rangle$ auf $X$ definiert wird, wobei $Q_s = \| \ldots \|^2$ sei.
Zunächst bemerken wir, daß $\langle \ldots, \ldots \rangle$ dann bzgl.\ der Maximumsnorm stetig ist, da die Norm stetig ist.
Die Fälle $\K = \R$ und $\K = \C$ sind nun getrennt voneinander zu betrachten.
Man rechnet jeweils unter Verwendung der Parallelogrammgleichung nach, daß $\langle \ldots, \ldots \rangle$ in der ersten Komponente additiv ist.
Hieraus folgt nun, daß für alle $\lambda \in \Q$ und $x,y \in X$ gilt $\langle \lambda \, x, y \rangle = \lambda \, \langle x, y \rangle$.
Aus Stetigkeitsgründen gilt dies dann auch für alle $\lambda \in \R$.
Um die Gültigkeit der obigen Gleichung im Falle $\K = \C$ für alle $\lambda \in \C$ zu beweisen, genügt es nun den Fall $\lambda := \i$ zu betrachten, der sich durch simples Nachrechnen verifizieren läßt.
Abschließend überzeugt man sich durch eine erneute Rechnung davon, daß $\forall_{x,y \in X} \, \langle x, y \rangle = \overline{\langle y, x \rangle}$ gilt. \q

\begin{Kor} \label{FA.8.8.K}
Es seien $n \in \N_+$, $\varphi$ ein Quadermaß auf $\R^n$, $M$ eine Teilmenge von $\R^n$, die (mindestens) zwei disjunkte $\varphi$-integrierbare Teilmengen mit positivem $\varphi$-Maß besitzt, und $p \in {[}1, \infty{]}$.

Dann ist $L^p_{\K}(M,\varphi)$ genau dann ein $\K$-Prähilbertraum, wenn $p=2$ gilt. 
\end{Kor}

\textit{Beweis als Übung.} \q

\begin{Bem*} \label{FA.8.6.B}
Seien $n \in \N_+$ und $p \in {[}1, \infty{]} \setminus \{2\}$.
Das letzte Korollar zeigt, daß die zu der induzierten Norm $\|\ldots\|_2$ des $\K$-Prähilbertraumes $(\K^n, \langle \ldots, \ldots \rangle_2)$ äquivalente Norm $\|\ldots\|_p$ auf $\K^n$ nicht durch ein Skalarprodukt induziert wird.
\end{Bem*}

\begin{Satz} \label{FA.8.9}
Jeder $\K$-Prähilbertraum ist gleichmäßig konvex.
\end{Satz}

\textit{Beweis.} Es seien $X$ ein $\K$-Prähilbertraum, $\varepsilon \in {]}0, \frac{1}{2}{]}$ sowie $x,y \in X$ derart, daß $\|x\| = \|y\| = 1$ und $\|x-y\| \ge \varepsilon$ gilt.
Dann folgt zunächst aus der nach dem letzten Satz gültigen Parallelogrammgleichung $\|x+y\|^2 + \|x-y\|^2 = 4$ und sodann $\frac{\|x+y\|}{2} \le 1 - \delta$ mit $\delta := 1 - \sqrt{1 - \frac{\varepsilon^2}{4}} \in \R_+$. \q
\A
Die folgende Definition kann man natürlich auch in $\K$-Vek\-tor\-räu\-men, die lediglich mit einem Semiskalarprodukt -- und somit einer Halbnorm -- anstelle eines Skalarproduktes -- und somit einer Norm -- ausgestattet sind, geben.
Für jene Räume haben wir allerdings keine Topologie eingeführt.
Man kann eine solche analog zu der eines normierten $\K$-Vektorraumes definieren und einen lokal-konvexen topologischen $\K$-Vektorraum, der im Falle einer Halbnorm, die keine Norm ist, nicht hausdorffsch ist, erhalten.
Wenn es Mehraufwand bedeutet, möchte der Autor in einem Kapitel über Hilberträume aber keine Nicht-Hausdorff-Räume miteinbeziehen. 

\begin{Def}[Orthogonalität]  \label{FA.8.12}
Sei $X$ ein $\K$-Prähilbertraum.
\begin{itemize}
\item[(i)] Zwei Elemente $x,y \in X$ heißen \emph{zueinander orthogonal}\index{Orthogonal} (i.Z.\ $\boxed{x \perp y}$) genau dann, wenn $\langle x, y \rangle = 0$ gilt.
\item[(ii)] Zwei Teilmengen $A, B$ von $X$ heißen \emph{zueinander orthogonal}\index{Orthogonal} (i.Z.\ $\boxed{A \perp B}$) genau dann, wenn $\forall_{x \in A} \forall_{y \in B} \, \langle x, y \rangle = 0$ gilt.
\item[(iii)] Sind $I$ eine Menge und $X_i$ für jedes $i \in I$ ein $\K$-Untervektorraum von $X$, so nennt man $X$ genau dann die \emph{orthogonale Summe von $X_i$, $i \in I$,} (i.Z.\ $\boxed{X = {\bigobot}_{i \in I} X_i}$), wenn $X = \bigoplus_{i \in I} X_i$ und $\forall_{i,j \in I} \, (i \ne j \Rightarrow X_i \perp X_j)$ gilt.
\item[(iv)] Für eine Teilmenge $A$ von $X$ heißt
\begin{equation} \label{FA.8.12.S}
A^{\perp} := \{ x \in X \, | \, \forall_{y \in A} \, \langle x, y \rangle = 0 \} = \bigcap_{y \in A} \overline{\langle \ldots, y \rangle}^{1}(\{0\})
\end{equation}
das \emph{Orthogonalkomplement\index{Orthogonal!-komplement} von $A$}.
\end{itemize}
\end{Def}

\begin{Satz} \label{FA.8.13}
Seien $X$ ein $\K$-Prähilbertraum und $A$ eine Teilmenge von $X$.
\begin{itemize}
\item[(i)] $A^{\perp}$ ist ein abgeschlossener $\K$-Untervektorraum von $X$, für den $A \perp A^{\perp}$ gilt.
\item[(ii)] $\overline{{\rm Span}_{\K} A} \subset A^{\perp \perp} := (A^{\perp})^{\perp}$.
\item[(iii)] $\emptyset^{\perp} = \{0\}^{\perp} = X$ und $X^{\perp} = \{0\}$.
\item[(iv)] Ist $B$ eine weitere Teilmenge von $X$, so gilt $A \subset B \Rightarrow B^{\perp} \subset A^{\perp}$.
\item[(v)] $A^{\perp \perp \perp} = A^{\perp} = (\overline{{\rm Span}_{\K} A})^{\perp}.$
\end{itemize}
\end{Satz}

\textit{Beweis.} (i) folgt sofort aus (\ref{FA.8.12.S}) sowie \ref{FA.8.7.K}, und (iv) ist trivial.

Zu (ii): Sei $a \in A$.
Für jedes $x \in A^{\perp}$ gilt dann $\langle x, a \rangle = 0$, also $\langle x, a \rangle = 0$.
Daher folgt zunächst $A \subset A^{\perp \perp}$ und sodann aus (i) die Behauptung von (ii).

Zu (iii): Die erste Aussage ist trivial, und die zweite folgt daraus, daß ein Skalarprodukt positiv definit ist.

Zu (v): $(\overline{{\rm Span}_{\K} A})^{\perp} \stackrel{(iv)}{\subset} A \stackrel{(ii)}{\subset} A^{\perp \perp \perp} \stackrel{(ii), (iv)}{\subset} (\overline{{\rm Span}_{\K} A})^{\perp}$. \q
\A
Da ein $\K$-Prähilbertraum in kanonischer Weise ein normierter $\K$-Vektorraum ist, kann man ,,seine`` Vervollständigung betrachten -- vgl.\ Satz \ref{FA.2.13} --, in der der Ausgangsraum somit dicht liegt.
Aus Stetigkeitsgründen handelt es sich bei jener nach Satz \ref{FA.8.8} wieder um einen $\K$-Prä\-hil\-bert\-raum, dessen Skalarprodukt wegen Korollar \ref{FA.8.4} die eindeutig bestimmte Fortsetzung des Ausgangsskalarproduktes ist.

\begin{Def}[Hilberträume] \index{Raum!Hilbert-} \index{Hilbert!-raum} \label{FA.8.10}
Ein $\K$-Prähilbertraum, der als normierter $\K$-Vek\-tor\-raum ein $\K$-Banachraum ist, heißt ein \emph{$\K$-Hilbertraum}\index{Raum!Hilbert-}\index{Hilbert!-raum}.
\end{Def}

\begin{Bsp} \label{FA.8.10.B} $\,$
\begin{itemize}
\item[1.)] Ein endlich-dimensionaler $\K$-Prähilbertraum ist nach Satz \ref{FA.2.E.7} ein $\K$-Hil\-bert\-raum.
\item[2.)] Sei $n \in \N_+$.
Die in den Beispielen 1.) bis 3.) zu \ref{FA.8.1} (iii) genannten Räume $L_{\K}^2(M,\varphi)$, $\K^n$ und $\ell_{\K}^2$ sind zusammen mit den dort genannten jeweiligen Skalarprodukten $\langle \ldots, \ldots \rangle_2$ nach dem Satz von \textsc{Riesz-Fischer} \ref{FA.5.20} sämtlich $\K$-Hil\-bert\-räume.
Der nach Aufgabe \ref{FA.5.CcRnpnichtvollst} nicht-vollständige $\K$-Untervektorraum $\mathcal{C}_c(\R^n,\K)$ von $\L_{\K}^2(\R^n,\mu_n)$ besitzt zwar eine Struktur als $\K$-Prä\-hil\-bert\-raum, kann aber kein $\K$-Hilbertraum sein.
Wir haben des weiteren in Korollar \ref{FA.8.8.K} bereits eingesehen, daß die $\K$-Banachräume $L^p_{\K}(\R^n,\mu_n)$ für alle $p \in [1, \infty] \setminus \{2\}$ keine $\K$-Prähilberträume ,,sind``.
\end{itemize}
\end{Bsp}

\begin{Bem*}
Insbesondere in den Anfängen der Entwicklung der \emph{Funktionalanalysis} war es üblich, $\ell^2_{\R}$ den \emph{Hilbertschen Raum}\index{Raum!Hilbertscher} zu nennen.
Diese Bezeichnung findet sich daher häufig in der Literatur der ersten Hälfte des letzten Jahrhunderts, siehe z.B. \cite[S.\ 287]{Hausdorff}.
Wir werden den Begriff eines \emph{Hilbertschen Raumes} unten etwas allgemeiner einführen.
\end{Bem*}

Das folgende Korollar ergibt sich unmittelbar aus Satz \ref{FA.8.9} und dem Satz von \textsc{Milman} \ref{FA.7.M}.

\begin{Kor} \label{FA.8.11}
Jeder $\K$-Hilbertraum ist reflexiv. \q
\end{Kor}

\subsection*{Der Rieszsche Darstellungssatz für Hilbert\-räume} \addcontentsline{toc}{subsection}{Der Rieszsche Darstellungssatz für Hilbert\-räume}

\begin{Def}[Projektionen] \label{FA.8.14.DP}
Es seien $X$ ein $\K$-Vektorraum und $P \: X \to X$ eine $\K$-lineare Abbildung.

$P$ heißt  genau dann \emph{idempotent} oder eine \emph{Projektion}\index{Projektion}, wenn $P \circ P = P$ gilt.
\end{Def}

\begin{Lemma} \label{FA.8.14.LP}
Es seien $X$ ein $\K$-Vektorraum und $P \: X \to X$ eine $\K$-lineare Abbildung.

Dann gilt:
\begin{itemize}
\item[(i)] $P$ Projektion $\Longleftrightarrow$ $\id_X - P$ Projektion.
\item[(ii)] $P$ Projektion \newline $\Longrightarrow$ $P(X) = \{ x \in X \, | \, P(x) = x \} = {\rm Kern} \, (\id_X - P) \, \wedge \, X = P(X) \oplus {\rm Kern} \, P$.
\item[(iii)] Sind $Y,Z$ $\K$-Untervektorräume von $X$ mit $X = Y \oplus Z$, so existiert genau eine Projektion $P \in {\rm Hom}_{\K}(X,X)$ mit $P(X) = Y$ und ${\rm Kern} \, P = Z$.
\end{itemize}
\end{Lemma}

\begin{Zusatz} $\,$
\begin{itemize}
\item[1.)] Wenn $X$ ein normierter $\K$-Vektorraum und $P$ zusätzlich stetig ist, dann sind die Teilmengen $P(X)$ und ${\rm Kern} \, P$ von $X$ in der Konklusion von (ii) abgeschlossene $\K$-Un\-ter\-vek\-tor\-räume.
\item[2.)] Sind $X$ ein $\K$-Banachraum und $Y,Z$ abgeschlossene $\K$-Untervektorräume von $X$ mit $X = Y \oplus Z$, so ist die Projektion $P$ in (iii) stetig.
\end{itemize}
\end{Zusatz}

\textit{Beweis als Übung.} \q

\begin{Def}[Orthogonalprojektionen] \label{FA.8.14.DOP}
Sei $X$ ein $\K$-Prähilbertraum.

Eine Projektion $P \: X \to X$ heißt \emph{Orthogonalprojektion}\index{Projektion!Orthogonal-}\index{Orthogonal!-projektion} genau dann, wenn sie \emph{selbstadjungiert} ist, d.h.\ per definitiomen, daß $\forall_{x, \tilde{x} \in X} \, \langle P \, x, \tilde{x} \rangle = \langle x, P \, \tilde{x} \rangle$ gilt.
\end{Def}

\begin{Kor} \label{FA.8.14.KOP}
Sind $X$ ein $\K$-Prähilbertraum und $P \: X \to X$ eine Orthogonalprojektion, so folgt $X = P(X) \obot {\rm Kern} \, P$. \q  
\end{Kor}

\begin{Def}[Orthogonalprojektionen von Hilberträumen auf abgeschlossene Untervektorräume] \index{Orthogonal!-projektion} \index{Projektion!Orthogonal-} \label{FA.8.14}
Seien $X$ ein $\K$-Hilbertraum und $Y$ ein abgeschlossener $\K$-Unter\-vek\-tor\-raum von $X$.

Wir definieren die \emph{Orthogonalprojektion $\boxed{\pi_Y \: X \to X}$ von $X$ auf $Y$}\footnote{Die Namenswahl ist dadurch gerechtfertigt, daß der nächste Satz u.a.\ $\pi_Y \in {\rm Hom}_{\K}(X,X)$, $\pi_Y \circ \pi_Y = \pi_Y$ und $\forall_{x,\tilde{x} \in X} \, \langle \pi_Y(x), \tilde{x} \rangle = \langle x, \pi_Y(\tilde{x}) \rangle$ besagt.} durch
\begin{equation*} \label{FA.8.14.DOPx}
\forall_{x \in X} \, \| x - \pi_Y(x) \| = \inf \{ \| x - y \| \, | \, y \in Y \} \stackrel{\text{Def.}}{=} d( \{x\} , Y),
\end{equation*}
beachte, daß $X$ nach \ref{FA.8.9} ein gleichmäßig konvexer $\K$-Banachraum ist und die Approximationsaufgabe $A(x,Y)$ somit nach \ref{FA.7.A} (ii) offenbar für jedes $x \in X$ eine eindeutig bestimmte Lösung $\pi_Y(x) \in Y$ besitzt.
\end{Def}

\begin{Bem*}
Daß die Approximationsaufgabe $A(x,C)$ in $\K$-Hilberträumen $X$ für jede nicht-leere abgeschlossene konvexe Teilmenge $C$ von $X$ und jedes $x \in X$ genau eine Lösung besitzt, kann man auch aus der Parallelogrammgleichung herleiten, ohne die gleichmäßige Konvexität von $X$ auszunutzen.
\end{Bem*}

\begin{Satz} \label{FA.8.15}
Seien $X$ ein $\K$-Hilbertraum und $Y$ ein abgeschlossener $\K$-Un\-ter\-vek\-tor\-raum von $X$ (also ist $Y^{\perp}$ nach \ref{FA.8.13} (i) ebenfalls ein abgeschlossener $\K$-Unter\-vek\-tor\-raum von $X$).

Dann folgt
\begin{itemize}
\item[(i)] $\pi_Y, \pi_{Y^{\perp}} = \id_X - \pi_Y \in \L_{\K}(X,X)$, $\|\pi_Y\| \le 1$ und $\|\pi_{y^{\perp}}\| \le 1$ mit Gleicheit in der vorletzten bzw.\ letzten Ungleichung genau dann, wenn $Y \ne \{0\}$ bzw.\ $Y^{\perp} \ne \{0\}$ gilt,
\item[(ii)] $X = Y \obot Y^{\perp}$ (\emph{orhogonale Spaltung}) und $Y = Y^{\perp \perp}$,
\item[(iii)] $\pi_Y(X) = Y = {\rm Kern} \, \pi_{Y^{\perp}}$ und ${\rm Kern} \, \pi_Y = Y^{\perp} = \pi_{Y^{\perp}}(X)$,
\item[(iv)] $\pi_Y \circ \pi_Y = \pi_Y$, $\pi_{Y^{\perp}} \circ \pi_{Y^{\perp}} = \pi_{Y^{\perp}}$ und $\pi_Y \circ \pi_{Y^{\perp}} = \pi_{Y^{\perp}} \circ \pi_Y = 0$ sowie
\item[(v)] $\forall_{x,\tilde{x} \in X} \, \langle \pi_Y(x), \tilde{x} \rangle = \langle x, \pi_Y(\tilde{x}) \rangle \, \wedge \, \langle \pi_{Y^{\perp}}(x), \tilde{x} \rangle = \langle x, \pi_{Y^{\perp}}(\tilde{x}) \rangle$.
\end{itemize}
\end{Satz}

\textit{Beweisskizze.} Die Definition der Orthogonalprojektion $\pi_Y$ von $X$ auf $Y$ impliziert $\pi_Y(Y) \subset Y$ und $\forall_{y \in Y} \, \pi_Y(y) = y$, also folgt
\begin{gather}
\pi_Y(Y) = Y, \label{FA.8.15.1} \\
\pi_Y \circ \pi_Y = \pi_Y, \label{FA.8.15.2}
\end{gather}
Wir definieren $\tilde{\pi} := \id_X - \pi_Y \: X \to X$, d.h.\ genau
\begin{equation} \label{FA.8.15.3}
\id_X = \pi_Y + \tilde{\pi},
\end{equation}
und behaupten
\begin{equation} \label{FA.8.15.4}
\tilde{\pi}(X) \subset Y^{\perp}.
\end{equation}

{[} Zum Nachweis von (\ref{FA.8.15.4}) ist zu zeigen
$$ \forall_{x \in X} \forall_{y \in Y} \, \langle x - \pi_Y(x), y  \rangle = 0. $$
Seien also $x \in X$ und $y \in Y$.
Da $Y$ ein $\K$-Untervektorraum von $X$ ist, gilt dann $\forall_{t \in \R} \, \pi_Y(x) + t \, y \stackrel{(\ref{FA.8.15.1})}{\in} Y$.
Somit besagt die Definition von $\pi_Y$, daß die differenzierbare Funktion $q \: \R \to \R$, die durch
\begin{eqnarray*}
\forall_{t \in \R} \, q(t) & := & \|x - (\pi_Y(x) + t \, y) \|^2 = \|(x - \pi_Y(x)) + t \, y) \|^2 \\
& \stackrel{\ref{FA.8.5} (i)}{=} & \|x - \pi_Y(x)\|^2 + 2 t \, {\rm Re}(\langle x - \pi_Y(x), y \rangle) +t^2 \, \|y\|^2
\end{eqnarray*}
gegeben ist, in Null ein Minimum annimmt.
Hieraus folgt $q'(0) = 0$, und dies bedeutet genau
\begin{equation} \label{FA.8.15.5}
{\rm Re}(\langle x - \pi_Y(x), y \rangle) = 0. 
\end{equation}
Im Falle $\K = \R$ ist der Beweis von (\ref{FA.8.15.4}) damit erbracht.
Weil $Y$ ein $\K$-Un\-ter\-vek\-tor\-raum von $Y$ ist, gilt (\ref{FA.8.15.5}) im Falle $\K = \C$ wegen der Beliebigkeit von $y \in Y$ zusätzlich auch für $\i \, y$ anstelle von $y$, d.h.
$$ 0 = {\rm Re}(\langle x - \pi_Y(x), \i \, y \rangle) =  {\rm Re}(- \i \langle x - \pi_Y(x), y \rangle) = {\rm Im}(\langle x - \pi_Y(x), y \rangle), $$
und (\ref{FA.8.15.4}) ist ebenfalls gezeigt. {]}
\pagebreak

Da $\langle \ldots, \ldots \rangle$ positiv definit ist, gilt $Y \cap Y^{\perp} = \{0\}$.
Hieraus, (\ref{FA.8.15.1}), (\ref{FA.8.15.4}) und (\ref{FA.8.15.3}) folgt zunächst die erste Aussage in (ii) und sodann (iii) mit $\tilde{\pi}$ anstelle von $\pi_{Y^{\perp}}$, (iv) mit $\tilde{\pi}$ anstelle von $\pi_{Y^{\perp}}$ sowie die $\K$-Linearität von $\pi_Y$ und $\tilde{\pi}$.
Zusammen mit (\ref{FA.8.15.2}) ergibt Lemma \ref{FA.8.14} (i), (iii) nun offenbar $\tilde{\pi} = \pi_{Y^{\perp}}$.
Daher gelten (iii) und (iv).

Des weiteren gilt für jedes $x \in X$ nach dem Satz des \textsc{Pythagoras} wegen $\pi_Y(x) \perp \pi_{Y^{\perp}}(x)$
$$ \|\pi_Y(x)\|^2 \le \|\pi_Y(x)\|^2 + \|\pi_{Y^{\perp}}(x)\|^2 = \|\pi_Y(x) + \pi_{Y^{\perp}}(x)\|^2 = \|x\|^2, $$
also ist $\pi_Y$ stetig mit $\|\pi\| \le 1$.
Im Falle $Y \ne \{0\}$ existiert $y \in Y \setminus \{0\}$ mit $\pi_Y(y) = y$, d.h.\ $\|\pi_Y(y)\| = \|y\| \ne 0$, und es folgt $\|\pi_Y\| = 1$.
Analog schließt man für $\pi_{Y^{\perp}}$, und (i) ist wegen $\pi_{Y^{\perp}} = \id_X - \pi_Y$ gezeigt.

Zum Beweis des Satzes bleiben $Y = Y^{\perp \perp}$ und (v) zu zeigen.
$Y \subset Y^{\perp \perp}$ ist nach \ref{FA.8.13} (ii) klar.
Zum Nachweis von $Y^{\perp \perp} \subset Y$ sei $y \in Y$.
Dann gilt für alle $\tilde{y} \in Y^{\perp} = \pi_{Y^{\perp}}(X)$
$$ 0 = \langle y, \tilde{y} \rangle = \langle \pi_Y(y) + \pi_{Y^{\perp}}(y), \tilde{y} \rangle = \langle \pi_{Y^{\perp}}(y), \tilde{y} \rangle, $$
also auch für $\tilde{y} = \pi_{Y^{\perp}}(y)$: $\langle \pi_{Y^{\perp}}(y), \pi_{Y^{\perp}}(y) \rangle = 0$.
Dies aber ist äquivalent zu $\pi_{Y^{\perp}}(y) = 0$ bzw.\ $y \in {\rm Kern} \, \pi_{Y^{\perp}} = Y$.

Zu (v): Seien $x, \tilde{x} \in X$. Dann gilt
\begin{align*}
\langle \pi_Y(x), \tilde{x} \rangle & = \langle \pi_Y(x), \pi_Y(\tilde{x}) + \pi_{Y^{\perp}}(\tilde{x}) \rangle = \langle \pi_Y(x), \pi_Y(\tilde{x}) + \pi_{Y^{\perp}}(\tilde{x}) \rangle \\
& = \langle \pi_Y(x), \pi_Y(\tilde{x}) \rangle = \langle \pi_Y(x) + \pi_{Y^{\perp}}(x), \pi_Y(\tilde{x}) \rangle = \langle x, \pi_Y(\tilde{x}) \rangle.
\end{align*}
Analog argumentiert man für $\pi_{Y^{\perp}}$. \q

\begin{HS}[Rieszscher Darstellungssatz für Hilberträume] \label{FA.8.RD} \index{Satz!Darstellungs- von \textsc{Riesz}} $\,$

\noindent \textbf{Vor.:} Sei $X$ ein $\K$-Hilbertraum.

\noindent \textbf{Beh.:} Die Abbildung
$$ \Lambda \: X \longrightarrow X', ~~ x \longmapsto \Lambda_y := \langle \ldots, y \rangle, $$
ist eine semilineare\footnote{\textbf{Definition.} Eine Abbildung $T \: X \to Y$ zwischen $\K$-Vektorräumen heißt \emph{semilinear}\index{Abbildung!semilineare} genau dann, wenn $\forall_{x,y \in X} \, T(x + y) = T(x) + T(y)$ und $\forall{\lambda \in X} \forall_{x \in X} \, T(\lambda \, x) = \overline{\lambda} \, x$ gilt.} isometrische\footnote{\textbf{Definition.} Eine semilineare Abbildung $T \: X \to Y$ zwischen normierten $\K$-Vektorräumen heißt \emph{isometrisch}\index{Abbildung!isometrische} genau dann, wenn $\forall_{x \in X} \, \|T(x)\| = \|x\|$ gilt.} Bijektion.
\end{HS}

\textit{Beweis.} Korollar \ref{FA.8.7.K} impliziert, daß $\Lambda$ wohldefiniert ist, und die Semilinearität ist trivial.
Weiterhin gilt für alle $y \in X$ und jedes $x \in X$ nach der Cauchy-Schwarzschen Ungleichung $\|\Lambda_y \, x\| = \|\langle x, y \rangle| \le \|y\| \|x\|$, also $\|\Lambda_y\| \le \|y\|$.
Wegen $\|\Lambda_y \, y\| = \|y\| \, \|y\|$ gilt folglich $\|\Lambda_y\| = \|y\|$, und $\Lambda$ ist eine isometrische Abbildung, die somit auch injektiv ist.

Zu zeigen bleibt, daß $\Lambda$ surjektiv ist.
Sei $F \in X'$.
Ohne Beschränkung der Allgemeinheit gelte $F \ne 0$.
Dann ist $Y := {\rm Kern} \, F$ ein abgeschlossener $\K$-Un\-ter\-vek\-tor\-raum, also folgt $X = Y \oplus Y^{\perp}$ aus \ref{FA.8.15} (ii), und es gilt $\dim_{\K} Y^{\perp} = 1$,
denn wegen $F \ne 0$ besagt der Homomorphisatz \ref{FA.2.H} (ii), daß $X/Y$ eindimensional ist.
Daher existiert ein Einheitsvektor $\tilde{y} \in X$ mit $Y^{\perp} = \K \, \tilde{y}$.
Sei dann $y := \overline{F(\tilde{y})} \, \tilde{y}$.
Es folgt $\Lambda_y|_Y = 0$ und $\Lambda_y(\tilde{y}) = F(a)$, also $\Lambda_y = F$. \q

\begin{Bem*}
Die soeben bewiesene Version des Rieszschen Darstellungssatzes enthält mit der Surjektivität von $\Lambda$ genau die Aussage \cite[Satz B]{Riesz34}, die \textsc{Friedrich Riesz}\footnote{Es gab auch einen bedeutenden Mathematiker \textsc{Marcel Riesz}. Er war ein jüngerer Bruder von \textsc{Friedrich Riesz}.} 1934 ursprünglich bewiesen hat.
Jene beinhaltet i.w.\ die Version \ref{FA.7.RD} des Rieszschen Darstellungssatzes für Lebesguesche Räume im Falle $p=2$, denn $L_{\K}^2(M,\varphi)$, wobei $n \in \R_+$, $\varphi$ ein Quadermaß auf $\R^n$ und $M$ eine $\varphi$-meßbare Teilmenge von $\R^n$ seien, ist nach Beispiel 2.) in \ref{FA.8.10.B} ein $\K$-Hilbertraum.\footnote{Man beachte, daß in \ref{FA.7.RD} die $\K$-Vek\-tor\-raum-Iso\-me\-trie
$$ L_{\K}^2(M,\varphi) \longrightarrow L_{\K}^2(M,\varphi)', \, \mathbf{f} \longmapsto \left( L_{\K}^2(M,\varphi) \to \K, \, \mathbf{g} \mapsto \int_M f \, g \, \d \varphi, \, f \in \mathbf{f}, g \in \mathbf{g}, \right), $$
und in \ref{FA.8.RD} die antilineare isometrische Bijektion
$$ L_{\K}^2(M,\varphi) \longrightarrow L_{\K}^2(M,\varphi)', \, \mathbf{f} \longmapsto \left( L_{\K}^2(M,\varphi) \to \K, \, \mathbf{g} \mapsto \int_M f \, \overline{g} \, \d \varphi, \, f \in \mathbf{f}, g \in \mathbf{g}, \right), $$
betrachtet wurde.}
Für $p \in {[}1,\infty{]} \setminus \{2\}$ anstelle von $p=2$ ist letzteres wegen \ref{FA.8.8.K} i.a.\ nicht der Fall. 
Der Rieszsche Darstellungssatz für die Räume der quadrat\-integrierbaren Funktionen kann übrigens umgekehrt auch nicht als Spezialfall desjenigen für Hilberträume aufgefaßt werden, weil einerseits
$$ \# L^2_{\K}(M,\varphi) \le \# \left( \K^M \right) = (\# \K)^{\# M} \le (\# \K)^{\# \R^n} = (2^{\aleph_0})^{(2^{\aleph_0})} = 2^{(\aleph_0 \cdot 2^{\aleph_0})} = 2^{(2^{\aleph_0})} $$
gilt und andererseits zu jeder vorgegebenen Kardinalzahl $\kappa > 1$ ein $\K$-Hil\-bert\-raum, der als Menge die Mächtigkeit $\kappa^{\aleph_0} \ge \kappa$ hat, existiert, vgl.\ \ref{FA.8.H.K} unten.
\end{Bem*}

\subsection*{Hilbertbasen} \addcontentsline{toc}{subsection}{Hilbertbasen}

Jeder $\K$-Vektorraum ist bis auf $\K$-Vektorraum-Isomorphie durch die Mächtigkeit einer Basis eindeutig bestimmt.
Wir wollen in diesem Abschnitt ein Analogon für $\K$-Hilberträume beweisen.
Hierfür müssen wir den Begriff der \emph{Basis} geeignet anpassen und benötigen einige Vorbereitungen.
Wir beginnen, indem wir über \emph{Summen} über eine ggf.\ überabzählbare Indexmenge sprechen.  

\begin{Nr}[Summierbare Familien] \label{FA.8.16}
Seien $X$ ein normierter $\K$-Vektorraum, $I$ eine nicht-leere Menge und $(x_i)_{i \in I} \in X^I$ eine Familie von Elementen von $X$.
\\[2mm]
\noindent \textbf{Definition 1.} Wir nennen $(x_i)_{i \in I}$ \emph{(in $X$) summierbar}\index{Summierbarkeit} genau dann, wenn es ein $x \in X$ gibt derart, daß zu jedem $\varepsilon \in \R_+$ eine endliche Teilmenge $J_{\varepsilon}$ von $I$ existiert, so daß für alle endlichen Mengen $J$ mit $J_{\varepsilon} \subset J \subset I$
$$ \left\| \left( \sum_{i \in J} x_i \right) - x \right\| < \varepsilon $$
gilt, d.h.\ per definitionem \emph{$(x_i)_{i \in I}$ ist summierbar zu $x$}.

\begin{Bsp*}
$\# I < \infty$ $\Longrightarrow$ $(x_i)_{i \in I}$ summierbar.
\end{Bsp*}

\begin{Lemma*}
Seien $x,y \in X$.
Ist $(x_i)_{i \in I}$ sowohl zu $x$ als auch zu $y$ summierbar, so gilt $x=y$.
\end{Lemma*}

\textit{Beweis.} Angenommen, es gilt $x \ne y$.
Zu $\varepsilon := \|x-y\| \in \R_+$ existieren dann endliche Teilmengen $J_{\varepsilon},\widetilde{J}_{\varepsilon}$ von $I$ derart, daß für alle endlichen Teilmengen $J,\widetilde{J}$ von $I$ mit $J_{\varepsilon} \subset J$ und $\widetilde{J}_{\varepsilon} \subset \widetilde{J}$
$$ \left\| \left( \sum_{i \in J} x_i \right) - x \right\| < \frac{\varepsilon}{2} ~~ \wedge ~ \left\| \left( \sum_{i \in J'} x_i \right) - y \right\| < \frac{\varepsilon}{2} $$
gilt.
Dann ist $J:= J_{\varepsilon} \cap \widetilde{J}_{\varepsilon}$ eine endliche Teilmenge von $I$, die sowohl $J_{\varepsilon}$ als auch $\widetilde{J}_{\varepsilon}$ enthält, und es folgt
$$ \|x-y\| \le \left\| \left( \sum_{i \in J} x_i \right) - x \right\| + \left\| \left( \sum_{i \in J} x_i \right) - y \right\| < \varepsilon, $$
im Widerspruch zur Definition von $\varepsilon$. \q
\\[2mm]
\noindent \textbf{Definition 2.} Ist $(x_i)_{i \in I}$ summierbar, so heißt das nach dem Lemma eindeutig bestimmte Element von $X$, zu dem $(x_i)_{i \in I}$ summierbar ist, die \emph{Summe der $(x_i)_{i \in I}$} und wird mit $\boxed{{\sum}_{i \in I} x_i}$ bezeichnet.

\begin{Bem*}
Wir werden in Beispiel \ref{FA.8.22.B} sehen, daß es im Falle $I = \N$ wichtig ist, zwischen den Notationen $\sum_{i \in \N} x_i$ und $\sum_{i=0}^{\infty} x_i$ zu unterscheiden!
\end{Bem*}
\end{Nr}

\begin{Satz} \label{FA.8.17} $\,$

\noindent \textbf{Vor.:} Seien $X$ ein normierter $\K$-Vektorraum und $I$ eine nicht-leere Menge.

\noindent \textbf{Beh.:}
\begin{itemize}
\item[(i)] $(x_i)_{i \in I}$ summierbare Familie von Elementen von $X$ 
\newline $\Longrightarrow$ $\sum_{i \in I} x_i \in \overline{{\rm Span}_{\K} \{ x_i \, | \, i \in I \}}$.
\item[(ii)] $(x_i)_{i \in I}, (y_i)_{i \in I}$ summierbare Familien von Elementen von $X$ sowie $\lambda \in \R$ 
\newline $\Longrightarrow$ $(x_i + y_i)_{i \in I}, (\lambda \, x_i)_{i \in I}$ summierbar und $\sum_{i \in I} (x_i + y_i) = \sum_{i \in I} x_i + \sum_{i \in I} y_i$ $~~~~~~~\, \mbox{sowie $\sum_{i \in I} \lambda \, x_i = \lambda \, \sum_{i \in I} x_i$}.$
\item[(iii)] Ist $Y$ ein weiterer normierter $\K$-Vektorraum, $T \: X \to Y$ ein beschränkter Operator und $(x_i)_{i \in I}$ eine summierbare Familie von Elementen von $X$, so ist $(T \, x_i)_{i \in I}$ summierbar, und es gilt $\sum_{i \in I} T \, x_i = T (\sum_{i \in I} x_i)$.
\item[(iv)] Ist $X$ sogar ein $\K$-Prähilbertraum und $(x_i)_{i \in I}$ eine summierbare Familie von Elementen von $X$ sowie $y \in X$, so ist $(\langle x_i, y \rangle)_{i \in I}$ summierbar, und es gilt $\sum_{i \in I} \langle x_i, y \rangle = \langle \sum_{i \in I} x_i, y \rangle$.
\end{itemize}
\end{Satz}

\textit{Beweis als Übung.} \q

\begin{Satz} \label{FA.8.22.L}
Seien $X$ ein normierter $\K$-Vektorraum, $I$ eine nicht-leere Menge und $(x_i)_{i \in I}$ eine summierbare Familie von Elementen von $X$.

Dann ist $\{ i \in I \, | \, x_i \ne 0 \}$ höchstens abzählbar.
\end{Satz}

\textit{Beweis.} Sei $x := \sum_{i \in I} x_i$.
Zu jedem $k \in \N$ existiert dann eine endliche Teilmenge $J_k$ von $I$ derart, daß für endlichen Teilmengen $J$ von $I$ mit $J_k \subset J$ gilt $\| (\sum_{i \in J} x_i) - x \| < \frac{1}{2(k+1)}$.
Somit ist $\widetilde{J} := \bigcup_{k \in \N} J_k$ eine höchstens abzählbare Teilmenge von $J$, und für jedes $i_0 \in I \setminus \widetilde{J}$ folgt
$$ \forall_{k \in \N} \, \|x_{i_0}\| = \left\| \left( \sum_{i \in J_k \cup \{i_0\}} x_i \right) - \left( \sum_{i \in J_k} x_i \right) \right\| \le \left\| \sum_{i \in J_k \cup \{i_0\}} x_i \right\| + \left\| \sum_{i \in J_k} x_i \right\| < \frac{1}{(k+1)}, $$
also $\|x_{i_0}\| = 0$ und schließlich $x_{i_0} \in I \setminus \{ i \in I \, | \, x_i \ne 0 \}$.
Damit ist gezeigt, daß $\{ i \in I \, | \, x_i \ne 0 \}$ eine Teilmenge der höchstens abzählbaren Menge $\widetilde{J}$ ist. \q

\begin{Satz} \label{FA.8.18}
Es seien $I$ eine nicht-leere Menge und $(a_i)_{i \in I}$ eine Familie nicht-negativer reeller Zahlen.
\begin{itemize}
\item[(i)] $(a_i)_{i \in I}$ ist genau dann summierbar, wenn gilt
$$ s := \sup \left\{ \sum_{i \in J} a_i \, | \, J \subset I \, \wedge \, \# J < \infty \right\} < \infty. $$
\item[(ii)]
Im Falle der Gültigkeit einer der beiden äquivalenten Aussagen in (i) gilt $\sum_{i \in I} a_i = s$.
\end{itemize} 
\end{Satz}

\textit{Beweis.} Zu (i): ,,$\Rightarrow$ `` Ist $(a_i)_{i \in I}$ summierbar, so existieren $a \in \R$ und zu $\varepsilon = 1$ eine endliche Teilmenge $J_1$ von $I$ derart, daß für alle endlichen Mengen $\widetilde{J}$ mit $J_1 \subset \widetilde{J} \subset I$ gilt
$$ \left| \left( \sum_{i \in \widetilde{J}} a_i \right) - a \right| < 1. $$
Hieraus folgt für jede endliche Teilmenge $J$ von $I$
$$ \sum_{i \in J} a_i = \left(\sum_{i \in J_1 \cup J} a_i \right) - a - \underbrace{\left(\sum_{i \in J_1 \setminus J} a_i \right)}_{\ge 0} + a \le \left| \left(\sum_{i \in J_1 \cup J} a_i \right) - a \right| + |a| < 1 + |a| < \infty. $$

,,$\Leftarrow$`` Zu $\varepsilon \in \R_+$ existiert eine endliche Teilmenge $J_{\varepsilon}$ von $I$ mit $s - \varepsilon < \sum_{i \in J_{\varepsilon}} a_i$.
Daher gilt wegen $\forall_{i \in I} \, a_i \ge 0$ für jede endliche Obermenge $J$ von $J_{\varepsilon}$, die in $I$ enthalten ist,
$$ s - \varepsilon < \sum_{i \in J} a_i \le s, $$
und es folgt $\sum_{i \in J} a_i = s$.

Zu (ii): Der Beweis der Rückrichtung von (i) zeigt auch (ii). \q

\begin{Bem} \label{FA.8.19}
Es sei $I$ eine nicht-leere Menge.
Aufgrund des letzten Satzes ist es sinnvoll, für eine nicht-summierbare Familie $(a_i)_{i \in I}$ nicht-negativer reeller Zahlen $\sum_{i \in I} a_i := \infty$ zu setzen.
\end{Bem}

Der letzte Satz verwendet durch die Bildung des Supremums die Vollständigkeit der Menge der reellen Zahlen.
Es ist offensichtlich, daß wir bei der Rückrichtung des folgenden Summierbarkeitskriteriums ebenfalls die Vollständigkeit des betrachteten normierten $\K$-Vek\-tor\-raumes benötigen.

\begin{Satz}[Cauchysches Summierbarkeitskriterium] \label{FA.8.20} $\,$

\noindent \textbf{Vor.:} Seien $X$ ein $\K$-Banachraum, $I$ eine nicht-leere Menge und $(x_i)_{i \in I}$ eine Familie von Elementen von $X$.

\noindent \textbf{Beh.:} $(x_i)_{i \in I}$ ist genau dann summierbar, wenn zu jedem $\varepsilon \in \R_+$ eine endliche Teilmenge $J_{\varepsilon}$ von $I$ existiert derart, daß für alle endlichen Teilmengen $J$ von $I$ mit $J \cap J_{\varepsilon} = \emptyset$ gilt $\| \sum_{i \in J} x_i \| < \varepsilon$.
\end{Satz}

\textit{Beweis.} ,,$\Rightarrow$`` Existiere also $x \in X$ derart, daß $(x_i)_{i \in I}$ summierbar zu $x$ ist, d.h.\ zu jedem $\varepsilon \in \R_+$ gibt es eine endliche Teilmenge $J_{\varepsilon}$ von $I$, so daß für alle Obermengen $J$ von $I$, die $J_{\varepsilon}$ enthalten, gilt $\| (\sum_{i \in J} x_i) - x \| < \frac{\varepsilon}{2}$.
Für ein solches $J$, das zusätzlich $J \cap J_{\varepsilon} = \emptyset$ erfüllt, ergibt sich dann 
$$ \left\| \sum_{i \in J} x_i \right\| = \left\| \left( \sum_{i \in J \cup J_{\varepsilon}} x_i \right) - \left( \sum_{i \in J_{\varepsilon}} x_i \right) \right\| \le \left\| \sum_{i \in J \cup J_{\varepsilon}} x_i \right\| + \left|\ \sum_{i \in J_{\varepsilon}} x_i \right\| < \varepsilon. $$

,,$\Leftarrow$`` Für alle $k \in \N$ wählen wir zu $\varepsilon_k := \frac{1}{k+1}$ eine endliche Teilmenge $J_{\varepsilon_k}$ von $I$ gemäß der rechten Seite der Behauptung.
Ohne Beschränkung der Allgemeinheit gelte $\forall_{k \in \N} \, J_{\varepsilon_k} \subset J_{\varepsilon_{k+1}}$.
Dann ist $(y_k)_{k \in \N}$, definiert durch $\forall_{k \in \N} \, y_k := \sum_{i \in J_{\varepsilon_k}} x_i$, eine Cauchyfolge in $X$ (denn für $k,l \in \N$ mit $k \ge l$ gilt $(J_{\varepsilon_k} \setminus J_{\varepsilon_l}) \cap J_{\varepsilon_l} = \emptyset$ und somit $\| y_k - y_l \| = \| \sum_{i \in J_{\varepsilon_k} \setminus J_{\varepsilon_l}} x_i \| < \varepsilon_l = \frac{1}{l+1}$), die wegen der Vollständigkeit von $X$ gegen ein gewisses $y \in X$ konvergiert.
Wir behaupten $\sum_{i \in I} x_i = y$.

Zum Nachweis hiervon sei $\varepsilon \in \R_+$.
Dann existiert eine natürliche Zahl $k \in \N$ mit $\|y_k - y\| < \varepsilon_k = \frac{1}{k+1} < \frac{\varepsilon}{2}$.
Daher gilt für jede endliche Teilmenge $J$ von $I$, die $J_{\varepsilon_k}$ enthält, $(J \setminus J_{\varepsilon_k}) \cap J_{\varepsilon_k} = \emptyset$, also auch
$$ \left\| \left(\sum_{i \in J} x_i \right) - y \right\| = \left\| y_k + \left( \sum_{i \in J \setminus J_{\varepsilon_k}} x_i \right) - y \right\| \le \| y_k - y \| + \left\| \sum_{i \in J \setminus J_{\varepsilon_k}} x_i \right\| < \varepsilon, $$
womit die Behauptung bewiesen ist. \q

\begin{Kor}[Majorantenkriterium] \label{FA.8.21}
Seien $X$ ein $\K$-Banachraum, $I$ eine nicht-leere Menge, $(x_i)_{i \in I}$ eine Familie von Elementen von $X$ und $(a_i)_{i \in I}$ eine summierbare Familie reeller Zahlen derart, daß $\forall_{i \in I} \, \|x_i\| \le a_i$ gilt.

Dann ist $(x_i)_{i \in I}$ summierbar und es gilt $\| \sum_{i \in I} x_i\| \le \sum_{i \in I} a_i$. \q
\end{Kor}

\begin{Kor} \label{FA.8.21.K}
Seien $X$ ein $\K$-Banachraum, $I$ eine nicht-leere Menge und $(x_i)_{i \in I}$ eine Familie von Elementen von $X$.

Ist $(x_i)_{i \in I}$ \emph{absolut summierbar}\index{Summierbarkeit!absolute}, d.h.\ per definitionem, daß $(\|x_i\|)_{i \in I}$ in den reellen Zahlen summierbar ist, so ist $(x_i)_{i \in I}$ in $X$ summierbar (und es gilt $\| \sum_{i \in I} x_i \| \le \sum_{i \in I} \|x_i\|$). \q
\end{Kor}

\begin{Bem*}
Ist $X$ ein endlich-dimensionaler normierter $\K$-Vektorraum, so überlassen wir es dem Leser als Übung, zu zeigen, daß eine summierbare Familie von Elementen von $X$ auch absolut summierbar ist.
\end{Bem*}

\begin{Satz} \label{FA.8.22} $\,$

\noindent \textbf{Vor.:} Es seien $X$ ein $\K$-Banachraum, $I$ eine nicht-leere Menge und $(x_i)_{i \in I}$ eine Familie von Elementen von $X$.
Ferner sei $I_* := \{ i \in I \, | \, x_i \ne 0 \}$ keine endliche Menge.

\pagebreak
\noindent \textbf{Beh.:} Die folgenden beiden Aussagen sind äquivalent:
\begin{itemize}
\item[(i)] $(x_i)_{i \in I}$ ist summierbar.
\item[(ii)] $I_*$ ist abzählbar, und es existiert ein $x \in X$ derart, daß die Reihe $\sum_{k=0}^{\infty} x_{j(k)}$ für jede Bijektion $j \: \N \to I_*$ gegen $x$ konvergiert.\footnote{Es sei an dieser Stelle darauf hingewiesen, daß man so etwas wie ,,unbedingte Summierbarkeit einer Familie von Elementen eines normierten $\K$-Vektorraumes`` nicht definieren kann!}
\end{itemize}
\begin{Zusatz}
Aus der Gültigkeit einer der beiden äquivalenten Aussagen (i) oder (ii) folgt $\sum_{i \in I} x_i = \sum_{k=0}^{\infty} x_k$.
\end{Zusatz}
\end{Satz}

\textit{Beweis.} ,,(i) $\Rightarrow$ (ii), Zusatz`` Daß $I_*$ abzählbar ist, ergibt sich sofort aus Satz \ref{FA.8.22.L}.
Seien also $x := \sum_{i \in I} x_i$, $j \: \N \to I_*$ eine Bijektion und $\varepsilon \in \R_+$.
Dann existiert eine endliche Teilmenge $J_{\varepsilon}$ von $I$ derart, daß für alle endlichen Teilmengen $J$ von $I$, die $J_{\varepsilon}$ enthalten, gilt $\| ( \sum_{i \in I} x_i ) - x \| < \varepsilon$.
Ohne Beschränkung der Allgemeinheit können wir $J_{\varepsilon} \subset I_*$ annehmen, und wir setzen $k_0 := \max \{ k \in \N \, | \, j(k) \in J_{\varepsilon} \}$.
Für alle $l, m \in \N$ mit $m > l \ge k_0$ gilt dann $\| (\sum_{k=l+1}^m x_{j(k)} ) - x \| < \varepsilon$, und wir haben $\sum_{k=0}^{\infty} x_{j(k)} = \sum_{i \in I} x_i$ gezeigt.

,,(ii) $\Rightarrow$ (i)`` Angenommen, $(x_i)_{i \in I}$ ist nicht summierbar, d.h.\ insbesondere, daß auch $(x_i)_{i \in I_*}$ nicht summierbar ist.
Aus dem Cauchyschen Summierbarkeitskriterium folgt dann offenbar die Existenz einer Zahl $\varepsilon \in \R_+$ und einer Folge $(J_{\kappa})_{\kappa \in \N}$ paarweise disjunkter nicht-leerer endlicher Teilmengen von $I_*$, so daß 
\begin{equation} \label{FA.8.22.1}
\forall_{k \in \N} \, \left\| \sum_{i \in J_k} x_i \right\| \ge \varepsilon
\end{equation}
gilt.
Wegen der Abzählbarkeit von $I_*$ ist $J_* := I_* \setminus \left( \bigcupdot_{\kappa \in \N} J_{\kappa} \right)$ höchstens abzählbar.

1.\ Fall: $\# J_* < \infty$.
Definiere dann eine bijektive Abbildung $j \: \N \to I_*$, indem man zunächst $J_* \setminus \left( \bigcupdot_{k \in \N} J_k \right)$ und sodann $J_0, J_1, J_2 \ldots$ in dieser Reihenfolge durchzählt.
Wegen der Voraussetzung von (ii) ist dann die Reihe $\sum_{k=0}^{\infty} x_{j(k)}$ konvergent in $X$ und somit eine Cauchyfolge in $X$, im Widerspruch zu (\ref{FA.8.22.1}).

2.\ Fall: $\# J_* = \infty$.
Definiere dann eine bijektive Abbildung $j \: \N \to I_*$, indem man zunächst ein Element $x_{j(0)} \in J_* \setminus \left( \bigcupdot_{k \in \N} J_k \right)$ auswählt, dann $J_0$ abzählt, nun ein Element aus $J_* \setminus \left( \bigcupdot_{k \in \N} J_k \right)$, das ungleich $x_{j(0)}$ ist, auswählt, sodann $J_1$ abzählt usw.
Wie im ersten Falle ergibt sich ein Widerspruch zu (\ref{FA.8.22.1}). \q

\begin{Bem*} $\,$
\begin{itemize}
\item[1.)] Der Beweis der Richtung ,,(i) $\Rightarrow$ (ii)`` des Satzes \ref{FA.8.22} verwendet die Vollständigkeit von $X$ nicht.
\item[2.)] Ein Analogon zu Satz \ref{FA.8.22} ist im Falle $\# I_* < \infty$ trivialerweise für jeden lediglich normierten $\K$-Vektorraum richtig.
\end{itemize}
\end{Bem*}

\begin{Kor} \label{FA.8.22.K}
In einem $\K$-Banachraum $X$ ist ein Familie $(x_i)_{i \in \N}$ von Elementen von $X$ genau dann summierbar, wenn $\sum_{i=0}^{\infty} x_i$ unbedingt konvergiert.
In diesem Falle gilt $\sum_{i \in \N} x_i = \sum_{i=0}^{\infty} x_i$. \q
\end{Kor}

\begin{Bsp} \label{FA.8.22.B}
Die alternierende harmonische Reihe konvergiert gegen $\ln(2)$, aber $\left( \frac{(-1)^i}{i+1} \right)_{i \in \N}$ ist nicht summierbar.
Es gilt also $\sum_{i=0}^{\infty} \frac{(-1)^i}{i+1} = \ln (2)$, und das Symbol $\sum_{i\in \N} \frac{(-1)^i}{i+1}$ ist gar nicht definiert.
\end{Bsp}

\pagebreak
Wir werden unten sehen, daß jeder $\K$-Hilbertraum bis auf Isomorphie von $\K$-Hilberträumen\footnote{\textbf{Definition.} Ein \emph{$\K$-(Prä-)Hilbertraum-Isomorphismus} zwischen $\K$-(Prä-)Hilberträumen ist ein skalarprodukttreuer $\K$-Vek\-tor\-raum-Isomorphismus (also insbesondere eine Isometrie) zwischen ihnen.} von der Gestalt eines \emph{Hilbertschen Raumes}\index{Raum!Hilbertscher} im Sinne der folgenden Definition ist.

\begin{Def}[Hilbertsche Räume] \label{FA.8.H}
Es seien $I$ eine nicht-leere Menge und der $\K$-Vektorraum $\boxed{\ell^2_{\K}(I)}$ der sog.\ \emph{Hilbertsche Raum}\index{Raum!Hilbertscher} aller Familien $(x_i)_{i \in I}$ von Elementen von $\K$ derart, daß $\left( |x_i|^2 \right)_{i \in I}$ summierbar ist, versehen mit dem durch
$$ \forall_{(x_i)_{i \in I}, (y_i)_{i \in I} \in \ell^2_{\K}(I)} \, \langle (x_i)_{i \in I}, (y_i)_{i \in I} \rangle_2 := \sum_{i \in I} x_i \, \overline{y_i}, $$
definierten Skalarprodukt $\boxed{\langle \ldots, \ldots \rangle_2}$.

{[} Seien $(x_i)_{i \in I}, (y_i)_{i \in I} \in \boxed{\ell^2_{\K}(I)}$ sowie $\lambda \in \K$.
Dann gilt für jede endliche Teilmenge $J$ von $I$
\begin{gather*}
\sum_{i \in J} |\lambda \, x_i|^2 = |\lambda|^2 \, \sum_{i \in J} |x_i|^2 \le |\lambda|^2 \, \sum_{i \in I} |x_i|^2, \\
\sum_{i \in J} |x_i + y_i|^2  \le  \sum_{i \in J} |x_i|^2 + \sum_{i \in J} |y_j|^2 \le  \sum_{i \in I} |x_i|^2 + \sum_{i \in I} |y_j|^2
\end{gather*}
und wegen der Cauchy-Schwarzschen Ungleichung in $(\R^n, \langle \ldots, \ldots \rangle_2^{\R^n})$, wobei $n := \# J \in \N$ sei, 
\begin{align*}
\sum_{i \in J} |x_i \, \overline{y_i}| & \, = \, \sum_{i \in J} |x_i| \, |y_i| \, \le \, \sqrt{\sum_{i \in J} |x_i|^2} \, \sqrt{\sum_{i \in J} |y_i|^2} \, \le \,\sqrt{\sum_{i \in I} |x_i|^2} \, \sqrt{\sum_{i \in I} |y_i|^2} \, < \infty,
\end{align*}
also folgt aus Satz \ref{FA.8.18} (i) sowohl $(\lambda \, x_i)_{i \in I}, (x_i + y_i)_{i \in I} \in \ell^2_{\K}(I)$ als auch die absolute Summierbarkeit der Familie $(x_i \, \overline{y_i})_{i \in I}$ von Elementen von $\K$.
Zusammen mit Korollar \ref{FA.8.21.K} zeigt dies, daß $\ell^2_{\K}(I)$ ein $\K$-Prähilbertraum ist, denn $\langle \ldots, \ldots, \rangle_2$ erfüllt trivialerweise die Eigenschaften eines Skalarproduktes. {]}
\end{Def}

\begin{Satz} \label{FA.8.H.S}
Für jede nicht-leere Menge $I$ ist der Hilbertsche Raum $\ell_{\K}^2(I)$ ein $\K$-Hil\-bert\-raum.
\end{Satz}

\textit{Beweis.} Zum Nachweis des Satzes bleibt nur noch zu zeigen, daß $\ell_{\K}^2(I)$ vollständig ist:
Sei $\left( (x_{k,i})_{i \in I} \right)_{k \in \N}$ eine Cauchyfolge in $\ell_{\K}^2(I)$.
Für jedes $i \in I$ ist dann $(x_{k,i})_{k \in \N}$ offenbar eine Cauchyfolge in $\K$, die somit in $\K$ gegen ein gewisses $x_i \in \K$ konvergiert, und man sieht sofort, daß $\left( (x_{k,i})_{i \in I} \right)_{k \in \N}$ dann bzgl.\ der durch $\langle \ldots, \ldots \rangle_2$ induzierten Norm gegen $(x_i)_{i \in I}$ konvergiert, d.h.\ genau $\lim_{k \to \infty} \sqrt{ \sum_{i \in I} |x_{k,i} - x_i|^2} = 0$.
Hieraus folgt $\sup \{ \sum_{i \in I} |x_i - x_{k,i}|^2 \, | \, k \in \N \} < \infty$, also gilt für jede endliche Teilmenge $J$ von $I$ und jedes $k \in \N$
$$ \sum_{i \in J} |x_i|^2 \le \sum_{i \in J} |x_i - x_{k,i}|^2 + \sum_{i \in J} |x_{k,i}|^2 \le \sum_{i \in I} |x_i - x_{k,i}|^2 + \sum_{i \in I} |x_{k,i}|^2 < \infty $$
und somit $(x_i)_{i \in I} \in \ell_{\K}^2(I)$ nach Satz \ref{FA.8.18} (i).
Damit haben wir bewiesen, daß $\left( (x_{k,i})_{i \in I} \right)_{k \in \N}$ in $\ell_{\K}^2(I)$ konvergiert. \q

\begin{Bem*}
Als $\K$-Banachräume sind $\ell_{\K}^2$ und $\ell_{\K}^2(\N)$ einander gleich.
\end{Bem*}
\pagebreak
Das folgende Korollar, dessen Beweis Kenntnisse über Kardinalzahlen verwendet, die der Leser z.B.\ in \cite[Chapter 5 ff.]{HJ} findet, beinhaltet eine Aussage, die in der Bemerkung zum Rieszschen Darstellungssatz für Hilberträume \ref{FA.8.RD} nur erwähnt wurde.

\begin{Kor} \label{FA.8.H.K}
Sei $\kappa$ eine Kardinalzahl mit $\kappa > 1$.

Dann ist $\ell_{\K}^2(\kappa)$ ein $\K$-Hilbertraum, der als Menge die Mächtigkeit $\kappa^{\aleph_0}$ besitzt.
\end{Kor}

\textit{Beweisskizze.} Im Falle $1 < \kappa < \aleph_0$ ist die Behauptung klar.
Im folgenden gelte daher $\kappa \ge \aleph_0$.

Weil $\kappa$ wohlgeordnet\footnote{Eine Menge heißt \emph{wohlgeordnet}\index{Wohlordnung} genau dann, wenn sie totalgeordnet ist und jede nicht-leere Teilmenge ein kleinstes Element besitzt.} ist, können wir $\ell_{\K}^2(\kappa)$ in kanonischer Weise injektiv in $(\K \times \kappa)^{\N}$ abbilden, also gilt offenbar $\mathfrak{c} \le \# \ell_{\K}^2(\kappa) \le \max\{ \mathfrak{c}, {\kappa}^{\aleph_0} \}$, wobei $\mathfrak{c} := \# \R$ wie üblich die Mächtigkeit des Kontinuums bezeichne.

Außerdem definiert jede Abbildung $f \: \aleph_0 \to \kappa$ via
$$ \forall_{\alpha \in \kappa} \, x_{f,\alpha} := \left\{ \begin{array}{cl} \sum_{i \in \overline{f}^1(\{\alpha\})} 2^{-i}, & \mbox{falls } \alpha \in f(\aleph_0), \\ 0, & \mbox{falls } \alpha \notin f(\aleph_0), \end{array} \right\} $$ 
ein Element $(x_{f,\alpha})_{\alpha \in \kappa}$ von $\ell_{\K}^2(\kappa)$, beachte
\begin{align*}
\sum_{\alpha \in \kappa} \left| x_{f,\alpha} \right|^2 & = \sum_{\alpha \in f(\aleph_0)} \bigg| \sum_{i \in \overline{f}^1(\{\alpha\})} 2^{-i} \bigg|^2 = \sum_{\alpha \in f(\aleph_0)} \sum_{i_1, i_2 \in \overline{f}^1(\{\alpha\})} 2^{-(i_1+i_2)} \\
& = \sum_{\alpha \in f(\aleph_0)} \bigg( \sum_{i \in \overline{f}^1(\{\alpha\})} 2^{-(2 i)} + ~ 2 \sum_{\substack{i_1,i_2 \in \overline{f}^1(\{\alpha\}) \\ i_1 < i_2}} 2^{-(i_1+i_2)} \bigg) \\
& \le \sum_{\alpha \in f(\aleph_0)} \bigg( \sum_{i \in \overline{f}^1(\{\alpha\})} 2^{-(2 i)} + ~ 2 \sum_{i_1 \in \overline{f}^1(\{\alpha\})} 2^{-i_1} \bigg) \mbox{ mit } \bigcupdot_{\alpha \in f(\aleph_0)} \overline{f}^1(\{\alpha\}) = \aleph_0  \\
& = \sum_{i=0}^{\infty} 4^{-i} + 2 \sum_{i=0}^{\infty} 2^{-i} = \frac{16}{3},
\end{align*}
und die Abbildung
\begin{equation*} \label{FA.8.H.K.1}
M := \left\{ f \in \kappa^{\aleph_0} \, | \, \forall_{\alpha \in \kappa} \forall_{i \in \overline{f}^1(\{\alpha\})} \exists_{j \in \aleph_0, j \ge i} \, j \notin \overline{f}^1(\{\alpha\}) \right\} \longrightarrow \ell^2_{\K}(\kappa), ~~ f \longmapsto (x_{f,\alpha})_{\alpha \in \kappa},
\end{equation*}
ist wegen der Eindeutigkeit der dyadischen Entwicklung einer jeden reellen Zahl $x = \sum_{i=0}^{\infty} a_i \, 2^{-i} \in {[}0,2{[}$, wobei $(a_i)_{i \in \N}$ eine Folge in $\{0,1\}$ mit
$$ \forall_{i \in \N} \exists_{j \in \N, j \ge i} \, a_j \ne 1 $$
ist, injektiv.
Ferner ist die Komplementärmenge von $M$ in $\kappa^{\aleph_0}$
$$ \bigcup_{\alpha \in \kappa} \bigcup_{i \in \aleph_0} \left\{ f \in \kappa^{\aleph_0} \, | \, \forall_{j \in \aleph_0, j \ge i} \, f(j) = \alpha \right\} $$
von der Mächtigkeit kleiner oder gleich 
$$ \kappa \cdot (\sum_{i \in \aleph_0} \kappa^i) = \kappa \cdot \aleph_0 \cdot \sup \{ \underbrace{\kappa^i}_{=\kappa} \, | \, i \in \aleph_0 \} = \kappa \cdot \aleph_0 \cdot \kappa = \kappa, $$
d.h.\ $\kappa^{\aleph_0} \le \# M + \kappa = \max\{ \# M, \kappa \}$ bzw.\ $\# M = \kappa^{\aleph_0}$, also auch ${\kappa}^{\aleph_0} \le \# \ell_{\K}^2(\kappa)$. 
\q

\pagebreak
Mit dem oben beschrieben Rüstzeug über summierbare Familien, das inhaltlich eigentlich Kapitel \ref{FAna2} zuzuordnen ist, können wir uns nun dem Begriff einer \emph{Hilbertbasis} zuwenden.
Solch eine wird per definitionem immer ein \emph{Orthonormalsystem} im folgenden Sinne sein.

\begin{Def}[Orthogonal- und Orthonormalsysteme] \label{FA.8.23}
Es seien $X$ ein $\K$-Prä\-hil\-bert\-raum und $I$ eine nicht-leere Menge sowie $(x_i)_{i \in I}$ eine Familie von Elementen von $X$.

$(x_i)_{i \in I}$ heißt genau dann ein \emph{Orthogonal-}\index{Orthogonal!-system} bzw.\ \emph{Orthonormalsystem von $X$}\index{Orthonormal!-system}, wenn $\forall_{i,j \in I} \, (i \ne j \Rightarrow x_i \perp x_j) \mbox{ bzw. } \langle x_i, x_j \rangle = \delta_{ij}$ gilt.
\end{Def}

\begin{Bsp*} $\,$
\begin{itemize}
\item[1.)] Sei $I$ eine nicht-leere Menge.
Für jedes $i \in I$ bezeichne $e_i \in \K^I$ die Funktion, die durch $\forall_{j \in I} \, e_i(j) := \delta_{ij}$ gegeben ist und die wir in kanonischer Weise mit $(e_i(j))_{j \in I} \in \ell^2_{\K}(I)$ identifizieren. 

Dann ist $(e_i)_{i \in I}$ ein Orthonormalsystem von $\ell^2_{\K}(I)$, das im Falle einer unendlichen Menge $I$ wegen \ref{FA.8.22} ,,(i) $\Rightarrow$ (ii)`` und $1_I \equiv (1_I(j))_{j \in I} \notin \ell^2_{\K}$ nicht in $\ell^2_{\K}(I)$ summierbar ist.
\item[2.)] $\big( \frac{1}{\sqrt{2 \pi}} \, \left. \e^{\i k x} \right|_{{[}- \pi, \pi{]}} \big)_{k \in \Z}$ ist ein Orthonormalsystem von $\L^2_{\C}({[}- \pi, \pi{]}, \mu_1)$.
\end{itemize}
\end{Bsp*}

\begin{Lemma} \label{FA.8.24}
Sind $X$ ein $\K$-Prähilbertraum, $I$ eine nicht-leere Menge und $(x_i)_{i \in I}$ ein Orthogonalsystem von $X$, so ist $(x_i)_{i \in I}$ $\K$-linear unabhängig, d.h. per definitionem, daß die Elemente jeder endlichen Teilmenge von $\{ x_i \, | \, i \in I \}$ $\K$-linear unabhängig sind.
\end{Lemma}

\textit{Beweis als Übung.} \q
\A
In $\K$-Hilberträumen können wir den Satz des \textsc{Pythagoras} auf Summen über beliebige anstelle von endlichen Indexmengen erweitern.

\begin{Satz} \label{FA.8.25} $\,$

\noindent \textbf{Vor.:} Seien $X$ ein $\K$-Hilbertraum, $I$ eine nicht-leere Menge und $(x_i)_{i \in I}$ ein Orthogonalsystem von $X$.

\noindent \textbf{Beh.:} 
\begin{itemize}
\item[(i)] $(x_i)_{i \in I}$ in $X$ summierbar $\Longleftrightarrow$ $\underbrace{(\|x_i\|^2)_{i \in I} \mbox{ in $\R$ summierbar}}_{\text{d.h.\ genau } (\|x_i\|)_{i \in I} \in \ell^2_{\K}(I)}$.
\item[(ii)](Satz des \textsc{Pythagoras})\index{Satz!des \textsc{Pythagoras}}
\newline
Im Falle der Gültigkeit einer der äquivalenten Aussagen in (i) folgt
$$ \left\| \sum_{i \in I} x_i \right\|^2 = \sum_{i \in I} \|x_i\|^2. $$
\end{itemize}
\end{Satz}

\textit{Beweis.} Zu (i): Für jede endliche Teilmenge $J$ von $I$ gilt aufgrund der Orthogonalität der $x_i$, $i \in I$, nach dem Satz des \textsc{Pythagoras} \ref{FA.8.5} (ii) für $\K$-Vek\-tor\-räume, die mit einem Semiskalarprodukt ausgestattet sind, 
$$ \left\| \sum_{i \in J} x_i \right\|^2 = \sum_{i \in J} \| x_i \|^2, $$ 
also folgt (i) offenbar aus dem Cauchyschen Summierbarkeitskriterium \ref{FA.8.20}.

\pagebreak
Zu (ii): Sei $x := \sum_{i \in I} x_i$.
Zweimalige Anwendung von \ref{FA.8.17} (iv) und die Orthogonalität der $x_i$, $i \in I$, ergeben dann 
$$ \langle x , x \rangle = \sum_{i \in I} \langle x_i , x \rangle = \sum_{i \in I} \overline{\langle x , x_i \rangle} = \sum_{i \in I} \overline{\sum_{j \in i}  \langle x_j, x_i \rangle} = \sum_{i \in I} \langle x_i, x_i \rangle, $$
und dies bedeutet genau $\| \sum_{i \in I} x_i \|^2 = \sum_{i \in I} \|x_i\|^2$. \q

\begin{Satz} $\,$ \label{FA.8.26}

\noindent \textbf{Vor.:} Es seien $X$ ein $\K$-Hilbertraum, $I$ eine nicht-leere Menge, $(e_i)_{i \in I}$ ein Orthonormalsystem von $X$ und $Y := \overline{{\rm Span}_{\K}\{e_i \, | \, i \in I\}}$.
($Y$ ist also in kanonischer Weise ebenfalls ein $\K$-Hilbertraum.)

\noindent \textbf{Beh.:}
\begin{itemize}
\item[(i)] Für jede Familie $(a_i)_{i \in I}$ von Elementen von $\K$ ist $(a_i \, e_i)_{i \in I}$ genau dann in $Y$ summierbar, wenn $(|a_i|^2)_{i \in I}$ in $\R$ summierbar ist (wobei letzteres wiederum genau $(a_i)_{i \in I} \in \ell^2_{\K}(I)$ bedeutet), und im Falle der Summierbarkeit gilt mit $y := \sum_{i \in I} a_i \, e_i \in Y$
\begin{equation*}
\forall_{i \in I} \, a_i = \langle y, e_i \rangle \label{FA.8.26.1}
\end{equation*}
-- die Koeffizienten der \emph{Fourierentwicklung von $y$}\index{Fourier!-entwicklung} $= \sum_{i \in I} a_i \, e_i$ sind also durch $y$ eindeutig bestimmt -- und
\begin{equation}
\|y\|^2 = \sum_{i \in I} |a_i|^2. \label{FA.8.26.2}
\end{equation}
\item[(ii)](Besselsche Ungleichung)\index{Ungleichung!Besselsche}
\newline
Ist $x \in X$, so ist $(|\langle x, e_i \rangle|^2)_{i \in I}$ in $\R$ summierbar, und es gilt
$$ \sum_{i \in I} |\langle x, e_i \rangle|^2 \le \|x\|^2. $$
\item[(iii)] Bezeichne $\pi_Y \: X \to X$ die Orthogonalprojektion von $X$ auf $Y$, vgl.\ \ref{FA.8.14}.
Für jedes $x \in X$ ist dann $(\langle x, e_i \rangle \, e_i)_{i \in I}$ in $Y$ summierbar, und es gilt
$$ \pi_Y(x) = \sum_{i \in I} \langle x, e_i \rangle \, e_i. $$
\item[(iv)](Gaußapproximation) \index{Gaußapproximation}
\newline
Für jedes $x \in X$ sind die \emph{Fourierkoeffizienten}\index{Fourier!-koeffizienten} $\langle x, e_i \rangle \in \K$ diejenigen Zahlen $a_i \in \K$, $i \in I$, für die $\|(\sum_{i \in I} a_i \, e_i) - x\|$ minimal wird.
\item[(v)] Für jedes $x \in X$ gilt genau dann $x \in Y$ (und d.i.\ wegen (iii) äquivalent zu $x = \sum_{i \in I} \langle x, e_i \rangle \, e_i$), wenn die \emph{Parsevalsche Gleichung}\index{Gleichung!Parsevalsche}
$$ \sum_{i \in I} |\langle x, e_i \rangle|^2 = \|x\|^2 $$
erfüllt ist.
\item[(vi)] Die Abbildung
$$ Y \longrightarrow \ell^2_{\K}(I), ~~ y \longmapsto (\langle y, e_i \rangle)_{i \in I}, $$
ist ein $\K$-Hilbertraum-Isomorphismus, d.h.\ genau ein $\K$-Vektorraum-Iso\-mor\-phis\-mus mit
$$ \forall_{y, \tilde{y} \in Y} \, \langle y, \tilde{y} \rangle = \sum_{i \in I} \langle y, e_i \rangle \, \overline{\langle \tilde{y}, e_i \rangle}. $$
\end{itemize}
\end{Satz}

\pagebreak
\textit{Beweis.} Zu (i): Die erste Aussage sowie (\ref{FA.8.26.2}) sind nach \ref{FA.8.25} (i) sowie (ii) klar, beachte \ref{FA.8.17} (i), und für jedes $i \in I$ gilt $\langle y, e_i \rangle = \sum_{j \in I} a_j \, \langle e_j, e_i \rangle = a_i$, wenn $y = \sum_{j \in I} a_j \, e_j$ existiert. 

Zu (ii): Seien $x \in X$ und $J$ eine endliche Teilmenge von $I$.
Dann gilt
\begin{align*}
0 & \le \| x - \sum_{i \in J} \langle x, e_i \rangle \, e_i \|^2 = \langle x - \sum_{i \in J} \langle x, e_i \rangle \, e_i, x - \sum_{j \in J} \langle x, e_j \rangle \, e_j \rangle \\
& = \langle x, x \rangle - \sum_{j \in J} \overline{\langle x, e_j \rangle} \, \langle x, e_j \rangle - \sum_{i \in J} \langle x, e_i \rangle \, \langle e_i, x \rangle + \sum_{i,j \in J} \langle x, e_i \rangle \, \overline{\langle x, e_j \rangle} \, \langle e_i, e_j \rangle \\
& = \|x\|^2 - \sum_{i \in J} | \langle x, e_i \rangle |^2,
\end{align*}
also $\sum_{i \in J} | \langle x, e_i \rangle |^2 \le \|x\|^2$.
Hieraus und aus \ref{FA.8.18} folgt (ii).

Zu (iii): Für jedes $x \in X$ ist $( \langle x, e_i \rangle \, e_i )_{i \in I}$ nach (ii), (i) in $Y$ summierbar, und für alle $j \in I$ gilt
\begin{align*}
\left\langle \left( \sum_{i \in I} \langle x, e_i \rangle \, e_i \right) - x, e_j \right\rangle & = \left\langle \sum_{i \in I} \langle x, e_i \rangle \, e_i, e_j \right\rangle - \langle x, e_j \rangle = \left( \sum_{i \in I} \langle x, e_i \rangle \, \langle e_i, e_j \rangle \right) - \langle x, e_j \rangle \\
& = \langle x, e_j \rangle - \langle x, e_j \rangle = 0,
\end{align*}
also folgt $\D \left(\sum_{i \in I} \langle x, e_i \rangle \, e_i \right) - x \in \{e_j \, | \, j \in I \}^{\perp} \stackrel{\ref{FA.8.13} (v)}{=} \left( \overline{{\rm Span}_{\K} \{e_j \, | \, j \in I \}} \right)^{\perp} \stackrel{\text{Def.}}{=} Y^{\perp}$ bzw.\
\begin{equation} \label{FA.8.26.3}
\sum_{i \in I} \langle x, e_i \rangle \, e_i \in Y \cap \left( x + Y^{\perp} \right).
\end{equation}

Des weiteren gilt nach Definition von $\pi_Y$ in \ref{FA.8.14} auch $\pi_Y(x) \in Y$ und für jedes $j \in J$ wegen $\pi_Y(x) + \langle x - \pi_Y(x), e_j \rangle \, e_j \in Y$
\begin{eqnarray*}
\| x - \pi_Y(x) \|^2 & \le & \| x - \pi_Y(x) - \langle x - \pi_Y(x), e_j \rangle \, e_j \|^2 \\
& = & \|x - \pi_Y(x) \|^2 - \overline{\langle x - \pi_Y(x), e_j \rangle} \, \langle x - \pi_Y(x), e_j \rangle \\
&   & - \, \langle x - \pi_Y(x), e_j \rangle \, \langle e_j, x - \pi_Y(x) \rangle \\
&   & + \, \langle x - \pi_Y(x), e_j \rangle \, \overline{\langle x - \pi_Y(x), e_j \rangle} \\
& = & \|x - \pi_Y(x) \|^2 - |\langle x - \pi_Y(x), e_j \rangle|^2,
\end{eqnarray*}
also $|\langle x - \pi_Y(x), e_j \rangle| = 0$ bzw.\ $\langle \pi_Y(x) - x, e_j \rangle = 0$, und d.h.\ analog zu oben
\begin{equation} \label{FA.8.26.4}
\pi_Y(x) \in Y \cap \left( x + Y^{\perp} \right).
\end{equation}

Wegen (\ref{FA.8.26.3}) und (\ref{FA.8.26.4}) genügt es zum Nachweis von (iii) zu zeigen, daß $Y \cap \left( x + Y^{\perp} \right)$ einelementig ist.
Seien also $y, \tilde{y} \in Y \cap \left( x + Y^{\perp} \right)$.
Dann ist direkt klar, daß $y - \tilde{y} \in Y$ gilt, und es existieren $z, \tilde{z} \in Y^{\perp}$ mit $y = x + z$ sowie $\tilde{y} = x + \tilde{z}$, also ergibt sich
$y - \tilde{y} = z - \tilde{z} \in Y \cap Y^{\perp} = \{0\}$ aus \ref{FA.8.15} (ii), d.h.\ $y = \tilde{y}$.

Zu (iv) - (vi): (iv) folgt sofort aus (iii) und der Definition von $\pi_Y$ in \ref{FA.8.14}.
(v) ist wegen $\| x - \sum_{i \in I} \langle x, e_i \rangle \, e_i \|^2 = \| x \|^2 - \sum_{i \in I} |\langle x, e_i \rangle|^2$ klar.
Die Aussage der $\K$-Vektorraum-Isomorphie in (vi) ergibt sich aus (i) und (v), und eine einfache Rechnung bestätigt die Skalarprodukttreue.
\q

\begin{Satz}[Gram-Schmidtsches Orthonormalisierungsverfahren]\index{Gram-Schmidtsches Or\-tho\-nor\-ma\-li\-sie- \newline rungs\-ver\-fahren}\index{Orthonormalisierungsverfahren nach\newline \textsc{Gram-Schmidt}} \label{FA.8.27} $\,$

\noindent \textbf{Vor.:} Es seien $X$ ein $\K$-Hilbertraum, $I$ entweder gleich $\{1, \ldots, n\}$ für ein $n \in \N_+$ oder gleich $\N_+$ und $(x_i)_{i \in I}$ eine Familie $\K$-linear unabhängiger Elemente von $X$, d.h.\ per definitionem, daß jeweils endlich viele Elemente von $(x_i)_{i \in I}$ $\K$-linear unabhängig sind.
Nach Korollar \ref{FA.2.E.8} ist übrigens $Y_i := {\rm Span}_{\K}\{x_1, \ldots, x_i\}$ für jedes $i \in I$ ein abgeschlossener $\K$-Untervektorraum von $X$.

\noindent \textbf{Beh.:} Es gilt
\begin{equation} \label{FA.8.27.1}
\forall_{i \in I, i>1} \, \pi_{{(Y_{i-1})}^{\perp}}(x_i) = x_i - \sum_{j = 1}^{i-1} \langle x_i, e_j \rangle \, e_j,
\end{equation}
und durch
\begin{gather*}
e_1 := \frac{x_1}{\|x_1\|} ~~ \mbox{ und } ~~ \forall_{i \in I, i>1} \, e_i := \frac{\pi_{{(Y_{i-1})}^{\perp}}(x_i)}{\|\pi_{{(Y_{i-1})}^{\perp}}(x_i)\|}
\end{gather*}
wird rekursiv ein Orthonormalsystem $(e_i)_{i \in I}$ von $X$ mit 
\begin{gather}
\forall_{i \in I} \, {\rm Span}_{\K} \{e_1, \ldots, e_i \} = Y_i \mbox{ sowie} \label{FA.8.27.2} \\
{\rm Span}_{\K} \{e_i \, | \, i \in I \} = {\rm Span}_{\K} \{x_i \, | \, i \in I \} \label{FA.8.27.3}
\end{gather}
definiert.
\end{Satz}

\textit{Beweis.} (\ref{FA.8.27.1}) ist wegen \ref{FA.8.15} (i) und \ref{FA.8.26} (iii) klar.

Wir zeigen induktiv für jedes $i \in I$, daß $(e_j)_{j \in I, j \le i}$ ein Orthonormalsystem von $X$ ist und daß (\ref{FA.8.27.2}) für $i$ gilt:
Im Falle $i=1$ ist dies trivial.
Sei nun $i \in I$ mit $i > 1$ und gelte die Behauptung für $i-1$.
Dann folgt aus
$$ \forall_{j \in I, j \le i-1} \, \langle e_j, e_i \rangle = \left\langle \underbrace{e_j}_{\in Y_{i-1}}, \frac{\pi_{{(Y_{i-1})}^{\perp}}(x_i)}{\|\pi_{{(Y_{i-1})}^{\perp}}(x_i)\|} \right\rangle = 0, $$
daß $(e_j)_{j \in I, j \le i}$ ein Orthonormalsystem von $X$ ist, also sind $e_1, \ldots, e_i$ nach \ref{FA.8.24} ebenso wie $x_1, \ldots, x_i$ $\K$-linear unabhängig.
Letztere bilden eine Basis von $Y_i$ und sind nach Induktionsvoraussetzung sowie (\ref{FA.8.27.1}) Elemente von ${\rm Span}_{\K} \{e_1, \ldots, e_i \}$.
Aus Dimensionsgründen muß daher (\ref{FA.8.27.2}) für $i$ gelten.

Ist $I$ endlich, so ist der Satz vollständig bewiesen.
Andernfalls bleibt (\ref{FA.8.27.3}) zu zeigen.
Dies ist aber klar, da wegen Aufgabe \ref{FA.8.27} (v) und (\ref{FA.8.27.2}) dann
\begin{align*}
{\rm Span}_{\K} \{e_i \, | \, i \in \N_+ \} & = \bigcup_{i \in \N_+} {\rm Span}_{\K} \{e_1, \ldots, e_i \} = \bigcup_{i \in \N_+} {\rm Span}_{\K} \{x_1, \ldots, x_i \} \\
& = {\rm Span}_{\K} \{x_i \, | \, i \in \N_+ \}
\end{align*}
gilt. \q
\A
Das Gram-Schmidtsche Orthonormalisierungsverfahren läßt sich modifiziert auch für eine positive Ordinalzahl anstelle von $I$ mittels \emph{transfiniter Induktion} beweisen.
Die Modifikation ist i.w.\ erforderlich, weil es sogar im Falle einer abzählbaren Familie $(x_i)_{i \in I}$ $\K$-linear unabhängiger Elemente eines $\K$-Hil\-bert\-raum\-es passieren kann, daß sich nicht jedes Element von $\overline{{\rm Span}_{\K} \{ x_i \, | \, i \in I \}}$ als $\sum_{i \in I} a_i \, x_i$ mit $a_i \in \K$ für alle $i \in I$ schreiben läßt (und dies bedeutet offenbar genau $\{ \sum_{i \in I} a_i \, x_i \, | \, \forall_{i \in I} \, a_i \in \K \mbox{ und } (a_i \, x_i)_{i \in I} \mbox{ summierbar} \} \subsetneqq \overline{{\rm Span}_{\K} \{ x_i \, | \, i \in I \}}$).
Ein Beispiel hierzu findet der Leser in \cite[Section 17.5 Example 3]{NS}.

\pagebreak
Zum Verständnis des nächsten Satzes sind Grundkenntnisse über Ordinalzahlen und transfinite Induktion erforderlich, die sich in \cite[Chapter 6]{HJ} finden.

\begin{Satz}[Allgemeines Gram-Schmidtsches Orthonormalisierungsverfahren]\index{Gram-Schmidtsches Or\-tho\-nor\-ma\-li\-sie- \newline rungs\-ver\-fahren}\index{Orthonormalisierungsverfahren nach\newline \textsc{Gram-Schmidt}} \label{FA.8.27.M}

\noindent \textbf{Vor.:} Seien $X$ ein nicht-nulldimensionaler $\K$-Hilbertraum, $\kappa$ eine positive Kardinalzahl sowie $\{ x_{\gamma} \, | \, \gamma \in \kappa \}$ eine (nicht notwendigerweise $\K$-linear unabhängige) Familie von Elementen von $X \setminus \{0\}$ und $Y := \overline{{\rm Span}_{\K} \{ x_{\gamma} \, | \, \gamma \in \kappa \}}$.

\noindent \textbf{Beh.:}
\begin{itemize}
\item[(i)] Für jede Ordinalzahl $\alpha$ wird wie folgt ein Orthonormalsystem $\mathfrak{E}_{\alpha}$ von $Y$ definiert derart, daß
\begin{equation} \label{FA.8.27.M.1}
\overline{{\rm Span}_{\K} \mathfrak{E}_{\alpha}} \subset Y ~~ \mbox{ und } ~~ \forall_{\tilde{\beta} \in \alpha} \, \mathfrak{E}_{\tilde{\beta}} \subset \mathfrak{E}_{\alpha}
\end{equation}
gilt:
Seien $\alpha$ eine Ordinalzahl und für alle $\beta \in \alpha$ bereits Orthonormalsysteme $\mathfrak{E}_{\beta}$ von $Y$ konstruiert, so daß (\ref{FA.8.27.M.1}) mit $\beta$ anstelle von $\alpha$ gilt.

1.\ Fall: $\overline{{\rm Span}_{\K} \bigcup_{\beta \in \alpha} \mathfrak{E}_{\beta}} = Y$.
Dann ist $\mathfrak{E}_{\alpha} := \bigcup_{\beta \in \alpha} \mathfrak{E}_{\beta}$ ein Orthonormalsystem von $Y$ mit (\ref{FA.8.27.M.1}).

2.\ Fall: $\overline{{\rm Span}_{\K} \bigcup_{\beta \in \alpha} \mathfrak{E}_{\beta}} \subsetneqq Y = \overline{{\rm Span}_{\K} \{ x_{\gamma} \, | \, \gamma \in \kappa \}}$.
Weil $\kappa$ wohlgeordnet ist, existiert eine kleinste Kardinalzahl $\gamma \in \kappa$ mit $x_{\gamma} \notin \overline{{\rm Span}_{\K} \bigcup_{\beta \in \alpha} \mathfrak{E}_{\beta}}$.
Dann gilt
$$ e := \frac{\pi_{\, \overline{{\rm Span}_{\K} \bigcup_{\beta \in \alpha} \mathfrak{E}_{\beta}}^{\perp} \cap Y}(x_{\gamma})}{\big\| \pi_{\, \overline{{\rm Span}_{\K} \bigcup_{\beta \in \alpha} \mathfrak{E}_{\beta}}^{\perp} \cap Y}(x_{\gamma}) \big\|} \in {\overline{{\rm Span}_{\K} \mathfrak{E}_*}}^{\perp} \cap Y, $$
und $\mathfrak{E}_{\alpha} := (\bigcup_{\beta \in \alpha} \mathfrak{E}_{\beta}) \cup \{e\}$ ist ein Orthonormalsystem von $Y$ mit (\ref{FA.8.27.M.1}).
\item[(ii)] Es existiert eine Ordinalzahl $\alpha_* \le \kappa$ mit $\overline{{\rm Span}_{\K} \mathfrak{E}_{\alpha_*}} = Y$.
\end{itemize}
\end{Satz}

\textit{Beweis.} (i) folgt sofort mittels transfiniter Induktion, und (ii) ist trivial, da bei der Rekursion in (i) spätestens für $\alpha = \kappa$ der 1.\ Fall vorliegt. \q 

\begin{Bem*}
Setzt man sogar voraus, daß $\kappa$ eine positive Ordinalzahl anstelle einer positiven Kardinalzahl ist, so hat man für $x_{\gamma}$ im 2. Falle in (i) eine Wahl zu treffen. 
\end{Bem*}

\begin{Def}[Hilbertbasen] \label{FA.8.28}
Seien $X$ ein $\K$-Hilbertraum und $I$ eine nicht-leere Menge.

Eine Orthonormalsystem $(e_i)_{i \in I}$ von $X$ mit
$$ \overline{{\rm Span}_{\K} \{ e_i \, | \, i \in I \}} = X$$
\index{Orthonormal!-system!vollständiges --} heißt eine \emph{Hilbertbasis\index{Basis!Hilbert-}\index{Hilbert!-Basis} von $X$}.\footnote{In der Literatur werden Hilbertbasen eines Hilbertraumes in unserem Sinne auch als \emph{vollständige Orthonormalsysteme} oder einfach als \emph{Orthonormalbasen}\index{Orthonormal!-basis}\index{Basis!Orthonormal-} bezeichnet.}
\end{Def}

\begin{Bsp} \label{FA.8.28.B}
Wir überlassen es dem Leser als Übung, zu zeigen, daß der $\K$-Hilbertraum $L^2_{\K}({[}- \pi, \pi{]}, \mu_1)$ eine abzählbare Hilbertbasis besitzt.
Beispielsweise induziert die Familie $\big( \frac{1}{\sqrt{2 \pi}} \, \left. \e^{\i k x} \right|_{{[}- \pi, \pi{]}} \big)_{k \in \Z}$ von Elementen des $\C$-Prä\-hilbert\-raumes $\L^2_{\C}({[}- \pi, \pi{]}, \mu_1)$ eine Hil\-bert\-basis des $\C$-Hil\-bert\-raumes $L^2_{\C}({[}- \pi, \pi{]}, \mu_1)$.
Abzählbare Hilbertbasen von $L^2_{\K}({[}- \pi, \pi{]}, \mu_1)$ können keine Hamelbasen von $L^2_{\K}({[}- \pi, \pi{]}, \mu_1)$ sein, denn wir haben bereits in Satz \ref{FA.2.E.9} eingesehen, daß eine unendliche Hamelbasis eines $\K$-Ba\-nach\-raumes überabzählbar ist.
\end{Bsp}

Daß jeder $\K$-Hilbertraum tatsächlich eine Hilbertbasis besitzt, ergibt sich sofort aus dem letzten Satz, indem man eine dichte Teilmenge wählt und ausnutzt, daß diese gleichmächtig zu einer Kardinalzahl ist.
Man beachte, daß damit der zum Auswahlaxiom äquivalente \emph{Wohlordnungssatz} verwendet wird, der besagt, daß sich jede Menge wohlordnen läßt.
Im Rahmen der \emph{Axiomatischen Mengenlehre nach \textsc{Zermelo} und \textsc{Fraenkel}} läßt sich nämlich ohne Verwendung des Auswahlaxiomes nur zeigen, daß jede wohlgeordnete Menge gleichmächtig zu einer Anfangsordinalzahl (und d.i.\ per definitionem eine Kardinalzahl) ist.
Wir führen im folgenden einen eleganteren Beweis, der das (ebenfalls zum Auswahlaxiom äquivalente) Zornsche Lemma und nicht das allgemeine Gram-Schmidtsche Orthonormalisierungsverfahren verwendet.

\begin{Satz} \label{FA.8.30}
Jeder nicht-nulldimensionale $\K$-Hilbertraum besitzt eine Hilbert-\linebreak basis.
\end{Satz}

\textit{Beweis.} Seien $X$ ein nicht-nulldimensionaler $\K$-Hilbertraum, $\mathcal{E} \subset \mathfrak{P}(X)$ die (nicht-leere) Menge aller Orthonormalsysteme von $X$ und $M$ eine bzgl.\ der Teilmengenrelation totalgeordnete Teilmenge von $\mathcal{E}$.
Dann ist $\bigcup_{\mathfrak{E} \in M} \mathfrak{E} \in \mathcal{E}$ eine obere Schranke von $M$.
Nach dem Zornschen Lemma \ref{D.3.11} existiert daher ein maximales Element $\mathfrak{E}_*$ von $M$. 

Angenommen, es gälte $\overline{{\rm Span}_{\K} \mathfrak{E}_*} \subsetneqq X$.
Dann folgte\footnote{An dieser Stelle würde der Beweis eines analogen Resultates für i.a.\ nicht-vollständige $\K$-Prähilberträume zusammenbrechen.} ${\overline{{\rm Span}_{\K} \mathfrak{E}_*}}^{\perp} \ne \{0\}$ aus der orthogonalen Spaltung in Satz \ref{FA.8.15} (ii), und für jedes $x \in {\overline{{\rm Span}_{\K} \mathfrak{E}_*}}^{\perp}$ mit $\|x\| = 1$ wäre $\mathfrak{E_*} \cup \{x\}$ ein maximales Element von $M$, im Widerspruch zur Maximalität von $\mathfrak{E}_*$. \q

\begin{HS}[Charakterisierung einer Hilbertbasis] $\,$ \label{FA.8.29}

\noindent \textbf{Vor.:} Seien $X$ ein $\K$-Hilbertraum und $(e_i)_{i \in I}$ ein Orthonormalsystem von $X$.

\noindent \textbf{Beh.:} Die folgenden Aussagen sind paarweise äquivalent:
\begin{itemize}
\item[(i)] $(e_i)_{i \in I}$ ist eine Hilbertbasis von $X$.
\item[(ii)](Fourierentwicklung)\index{Fourier!-entwicklung}
\newline
Für jedes $x \in X$ ist $(\langle x, e_i \rangle \, e_i)_{i \in I}$ in $X$ summierbar, und es gilt
$$ x = \sum_{i \in I} \langle x, e_i \rangle \, e_i. $$
\item[(iii)] Für alle $x, \tilde{x} \in X$ ist $(\langle x, e_i \rangle \, \overline{\langle \tilde{x}, e_i \rangle})_{i \in I}$ in $\K$ summierbar, und es gilt
$$ \langle x, y \rangle = \sum_{i \in I} \langle x, e_i \rangle \, \overline{\langle \tilde{x}, e_i \rangle}. $$
\item[(iv)] Für jedes $x \in X$ ist $(|\langle x, e_i \rangle|^2)_{i \in I}$ in $\R$ summierbar, und es gilt die \emph{Parsevalsche Gleichung}\index{Gleichung!Parsevalsche}
$$ \|x\|^2 = \sum_{i \in I} |\langle x, e_i \rangle|^2. $$
\item[(v)] Die Abbildung
$$ X \longrightarrow \ell^2_{\K}(I), ~~ x \longmapsto (\langle x, e_i \rangle)_{i \in I}, $$
ist ein $\K$-Hilbertraum-Isomorphismus, d.h.\ genau ein $\K$-Vektorraum-Iso\-mor\-phis\-mus mit
$$ \forall_{x, \tilde{x} \in X} \, \langle x, \tilde{x} \rangle = \sum_{i \in I} \langle x, e_i \rangle \, \overline{\langle \tilde{x}, e_i \rangle}. $$
\item[(vi)] $\forall_{x \in X} \, \left( \forall_{i \in I} \, \langle x, e_i \rangle = 0 \Longrightarrow x = 0 \right)$.
\end{itemize}
\end{HS}

\textit{Beweis.} Sei $Y := {\rm Span}_{\K} \{ e_i \, | \, i \in I \}$.
Dann ist (i) gleichbedeutend mit $X = Y$, und letzteres ist nach Satz \ref{FA.8.26} (v) sowohl äquivalent zu (ii) als auch zu (iv).
(iii) folgt aus (ii), und (iv) folgt aus (iii).
Ferner besagt Satz \ref{FA.8.26} (vi), daß (i) die Aussage (v) impliziert, welche wiederum trivialerweise (vi) zur Folge hat.
Schließlich sieht man sofort ein, daß ,,$\neg$ (iv) $\Rightarrow$ $\neg$ (vi)`` gilt, womit der Hauptsatz dann vollständig bewiesen ist. \q

\begin{Satz} \label{FA.8.31}
Es seien $X$ ein $\K$-Hilbertraum und $I, J$ nicht-leere Mengen sowie $(e_i)_{i \in I}$, $(\tilde{e}_j)_{j \in J}$ Hilbertbasen von $X$.

Dann gilt $\# I = \#J$.
\end{Satz}

\textit{Beweis.} Ist $X$ ein endlich-dimensionaler $\K$-Vektorraum, so ist die Behauptung klar, weil $(e_i)_{i \in I}$ und $(\tilde{e}_j)_{j \in J}$ dann Hamelbasen von $X$ sind.
Seien $I$ und $J$ daher im folgenden keine endliche Mengen.
Wir setzen für jedes $i \in I$
$$ J_i := \{ j \in J \, | \, \forall_{i \in I} \, \langle e_i, \tilde{e}_j \rangle \ne 0 $$
und behaupten
\begin{gather} 
\forall_{i \in I} \, \# J_i \le \aleph_0, \label{FA.8.31.1} \\
J \subset \bigcup_{i \in I} J_i . \label{FA.8.31.2}
\end{gather}

{[} (\ref{FA.8.31.1}) ist trivial, da die Fourierentwicklung von $e_i = \sum_{j \in J} \langle e_i, \tilde{e}_j \rangle \, \tilde{e}_j$ für jedes $i \in I$ nur höchstens abzählbar viele nicht-verschwindende Summanden enthält.

Zu (\ref{FA.8.31.2}): Sei $j \in J$. Dann gilt
$$ 0 \ne \tilde{e}_j = \sum_{i \in I} \langle \tilde{e}_j, e_i \rangle \, e_i, $$
also existiert $i_0 \in I$ mit
$$ 0 \ne \langle \tilde{e}_j, e_{i_0} \rangle = \overline{\langle e_{i_0}, \tilde{e}_j \rangle}, $$
und dies bedeutet $j \in J_{i_0} \subset \bigcup_{i \in I} J_i$. {]}

Aus (\ref{FA.8.31.2}) und (\ref{FA.8.31.1}) folgt wegen der Unendlichkeit der Menge $I$
$$ \# J \le \# I \cdot \aleph_0 = \# I, $$
und durch Vertauschung der Rollen von $I$ und $J$ sieht man auch $\# I \le \# J$ ein.
\q

\begin{Def}[Die Hilbertdimension eines $\K$-Hilbertraumes] \label{FA.8.32}
Sei $X$ ein $\K$-Hilbertraum.
Wir definieren dann $\boxed{{\rm H}\,\text{-}\dim_{\K} X}$ als die Mächtigkeit einer beliebigen Hilbertbasis von $X$ und nennen diese Kardinalzahl die \emph{Hilbertdimension von $X$}\index{Hilbert!-dimension}.
Besitzt $X$ keine Hilbertbasis, d.h.\ genau $X = \{0\}$, so setzen wir die \emph{Hilbertdimension von $X$} ferner als $\boxed{{\rm H}\,\text{-}\dim_{\K} X} := 0$. 
\end{Def}

\begin{Bem*}
Im Falle eines endlich-dimensionalen $\K$-Hilbertraumes, ist die Mächtigkeit einer Hamelbasis stets gleich der Mächtigkeit einer Hilbertbasis.
Für unendlich-dimensionale $\K$-Hilberträume gilt dies i.a.\ nicht mehr, wie wir in Beispiel \ref{FA.8.28.B} gesehen haben.
Weil sich jede $\K$-linear unabhängige Teilmenge eines $\K$-Vektorraumes zu einer Hamelbasis vervollständigen läßt\footnote{Dies beweist man analog zum Beweis von Satz \ref{FA.A.2}, indem man die Menge aller $\K$-linear unabhängigen Teilmengen, die die vorgegebene $\K$-linear unabhängige Teilmenge enthalten, betrachtet.}, ist jede Hilbertbasis eines $\K$-Hil\-bert\-raumes $X$ Teilmenge einer Hamelbasis von $X$, und es folgt 
$$ {\rm H}\,\text{-}\dim_{\K} X \le \dim_{\K} X, $$
wobei $\dim_{\K} X$ hier die Mächtigkeit einer (beliebigen) Hamelbasis von $X$ bezeichne, also im unendlich-dimensionalen Falle nicht etwa das Symbol ,,$\infty$``.
\end{Bem*}

\begin{HS} \label{FA.8.33}
Zwischen je zwei $\K$-Hil\-bert\-räumen existiert genau dann ein $\K$-Hil\-bert\-raum-Isomorphismus, wenn sie dieselbe Hilbertdimension haben.
\end{HS}

\textit{Beweis.} Der nulldimensionale Fall ist trivial.

,,$\Rightarrow$`` Seien $T \: X \to Y$ ein $\K$-Hil\-bert\-raum-Isomorphismus zwischen zwei $\K$-Hil\-bert\-räumen $X$ und $Y$ sowie $(e_i)_{i \in I}$ eine Hilbertbasis von $X$.
Dann ist $(T e_i)_{i \in I}$ eine Hilbertbasis von $Y$, und $X$ und $Y$ haben somit die gleiche Hilbertdimension.

,,$\Leftarrow$`` folgt offenbar aus Hauptsatz \ref{FA.8.29} ,,(i) $\Rightarrow$ (v)``.
\q

\begin{Kor} \label{FA.8.34}
Wenn $\kappa$ alle positiven Kardinalzahlen durchläuft, dann durchläuft $\ell^2_{\K}(\kappa)$ bis auf $\K$-Hilbertraum-Isomorphie alle nicht-nulldimensionalen $\K$-Hil\-bert\-räume (der Hilbertdimension $\kappa$). \q
\end{Kor}

\begin{Bsp} \label{FA.8.34.B}
Aus Korollar \ref{FA.8.34} und Beispiel \ref{FA.8.28.B} folgt, daß $\ell^2_{\K}(\aleph_0)$ als $\K$-Hil\-bert\-raum isomorph zu dem nach Korollar \ref{FA.5.22} separablen $\K$-Hilbertraum $L^2_{\K}({[}- \pi, \pi{]}, \mu_1)$ ist.
\end{Bsp}

\begin{Kor} \label{FA.8.35}
Ein nicht-nulldimensionaler $\K$-Hilbertraum ist genau dann separabel, wenn er eine höchstens abzählbare Hilbertbasis besitzt.
M.a.W.\ ist jeder separabele $\K$-Hil\-bert\-raum bis auf $\K$-Hilbertraum-Isomorphie entweder gleich $(\K^n, \langle \ldots, \ldots \rangle_2)$ für ein $n \in \N$ oder gleich $\ell_{\K}^2(\N)$.
\end{Kor}

\textit{Beweis.} Da die letztgenannten $\K$-Hilberträume separabel sind -- beachte Beispiel \ref{FA.8.34.B} --, ist zu zeigen, daß $\ell^2_{\K}(\kappa)$ für jede Kardinalzahl $\kappa > \aleph_0$ nicht separabel ist.
Angenommen, es existierte eine höchstens abzählbare dichte Teilmenge von $\ell^2_{\K}(\kappa)$ für eine Kardinalzahl $\kappa > \aleph_0$, d.h.\ ${\rm H}\,\text{-}\dim_{\K} \ell^2_{\K}(\kappa) > \aleph_0$.
Aus der Konstruktion im allgemeinen Gram-Schmidtschen Orthonormalisierungsverfahren \ref{FA.8.27.M} folgte dann ${\rm H}\,\text{-}\dim_{\K} \ell^2_{\K}(\kappa) \le \aleph_0$, Widerspruch! \q
\A
Das folgende Korollar ergibt sich offenbar aus Korollar \ref{FA.5.22} und dem vorgenannten Korollar.

\begin{Kor} \label{FA.8.36}
Seien $n \in \N_+$, $\varphi$ ein Quadermaß auf $\R^n$ und $J \in \mathfrak{I}_n$ ein nicht-entartetes Intervall von $\R^n$.

Dann ist $L^2_{\K}(J, \varphi)$ als $\K$-Hilbertraum isomorph zu $\ell^2(\N)$, wenn $L^2_{\K}(J, \varphi)$ unendlich-dimensional ist.
(Letzteres ist z.B.\ im Falle $\varphi = \mu_n$ immer erfüllt.) \q
\end{Kor}

\subsection*{Übungsaufgaben}

\begin{UA}[Charkterisierung der stetigen $\R$-bilinearen Abbildungen] \label{FA.8.A.St}
Es seien $X,Y,Z$ normierte $\K$-Vektorräume und $b \: X \times Y \to Z$ eine $\R$-bilineare Abbildung.
Ferner sei $X \times Y$ mit der Maximumsnorm versehen.
Den $\K$-Vektorraum aller stetigen $\R$-bilinearen Abbildung bezeichnen wir mit $\boxed{\L^2_{\R}(X,Y;Z)}$ und im Falle $X=Y$ mit $\boxed{\L^2_{\R}(X;Z)}$.

Zeige, daß die folgenden Aussagen paarweise äquivalent sind:
\item[(i)] $b \: X \times Y \to Z$ ist stetig, d.h.\ $b \in \L_{\R}^2(X,Y;Z)$.
\item[(ii)] $b \: X \times Y \to Z$ ist stetig in $(0,0)$.
\item[(iii)] $\exists_{C \in \R} \forall_{x \in X} \forall_{y \in Y} \, \|b(x,y)| \le C \, \|x\| \, \|y\|$.
\end{UA}

\begin{UA} \label{FA.8.A.St.N}
Es seien $X,Y,Z$ normierte $\K$-Vektorräume und des weiteren für jedes $b \in \L^2_{\R}(X,Y;Z)$
$$ \|b\| := \left\{ \frac{\|b(x,y)\|}{\|x\| \, \|y\|} \, | \, x \in X \setminus \{0\} \, \wedge \, y \in Y \setminus \{0\} \right\}  (:= 0, \mbox{ falls } X = \{0\} \mbox{ oder } Y = \{0\}). $$

Beweise, daß hierdurch eine Norm für $\L^2_{\R}(X,Y;Z)$ definiert wird und 
\begin{align*}
\|b\| & = \sup \{ \|b(x,y)\| \, | \, x \in X \, \wedge \, y \in Y \, \wedge \, \|x\| \le 1 \, \wedge \|y\| \le 1 \} \\
& = \underbrace{\sup \{ \|b(x,y)\| \, | \, x \in X \, \wedge \, y \in Y \, \wedge \, \|x\| = 1 \, \wedge \|y\| = 1 \}}_{:= 0, \mbox{ falls } X = \{0\} \mbox{ oder } Y = \{0\}}
\end{align*}
gilt.
Zeige ferner $\forall_{x \in X} \forall_{y \in Y} \, \|b(x,y)| \le \|b\| \, \|x\| \, \|y\|$ und daß $\|b\|$ die kleinste reelle Zahl mit dieser Eigenschaft ist.
\end{UA}

\begin{UA}[Fortsetzungssatz für stetige Sesquilinearformen]
Es seien $X$ ein normierter $\K$-Vektorraum sowie $Y$ ein dichter $\K$-Untervektorraum von $X$ und $s \: Y \times Y \to \K$ eine stetige Sesquilinearform.

Zeige, daß genau eine stetige Sesquilinearform $\overline{s} \: X \times X \to \K$ mit $\overline{s}|_{Y \times Y} = s$ existiert.
\end{UA}

Tip: Unter Ausnutzung der Aufgabe \ref{FA.8.A.St} und des Fortsetzungssatzes für stetige lineare Abbildungen setze man $s$ zunächst auf $X \times Y$ und anschließend auf $X \times X$ fort.

\begin{UA}
Beweise Satz \ref{FA.8.5} (ii).
\end{UA}

\begin{UA}
Beweise Satz \ref{FA.8.3}.
\end{UA}

\begin{UA}
Führe den Beweis des Satzes \ref{FA.8.8} in allen Einzelheiten aus.
\end{UA}

\begin{UA}
Beweise Korollar \ref{FA.8.8.K}.
\end{UA}

\begin{UA}
Seien $X$ ein $\K$-Prähilbertraum, $C$ eine nicht-leere vollständige konvexe Teilmenge von $X$ und $a \in X$.

Folgere aus der Parallelogrammgleichung, daß die Approximationsaufgabe $A(a,C)$ genau eine Lösung besitzt (ohne die gleichmäßige Konvexität von $X$ zu verwenden).
\end{UA}

\begin{UA}
Seien $X$ ein $\K$-Prähilbertraum und $Y$ ein vollständiger $\K$-Un\-ter\-vek\-tor\-raum von $X$.

Zeige $X = Y \obot Y^{\perp}$.
\end{UA}

\begin{UA}
Seien $X$ $\K$-Hilbertraum und $Y$ ein $\K$-Untervektorraum von $X$.

Beweise, daß $Y$ genau dann abgeschlossen ist, wenn $Y^{\perp \perp} = Y$ gilt.
\end{UA}

\begin{UA} \label{FA.8.13.U}
Seien $X$ ein $\K$-Hilbertraum und $A$ eine Teilmenge von $X$.

Zeige, daß sich Satz \ref{FA.8.13} (i) und (ii) dann wie folgt verschärfen lassen.
\begin{itemize}
\item[(i)] $A^{\perp}$ ist ein vollständiger $\K$-Untervektorraum von $X$.
\item[(ii)] $\overline{{\rm Span}_{\K} A} = A^{\perp \perp}$.
\end{itemize}
\end{UA}

\begin{UA}
Beweise Lemma \ref{FA.8.14.LP} inkl.\ Zusatz.
\end{UA}

Tip zu Zusatz 2.): $Y \times Z \to X, \, (y,z) \mapsto y+z,$ ist stetig sowie bijektiv und hat daher (!) eine stetige Umkehrabbildung. Hieraus folgt die Stetigkeit von $P$.

\begin{UA}[Satz von \textsc{Hellinger-Toeplitz}] \index{Satz!von \textsc{Hellinger-Toeplitz}}
Seien $X$ ein $\K$-Hilbertraum und $T \: X \to X$ ein \emph{selbstadjungierter} Operator, d.h.\ per definitionem, daß für den Operator $T \: X \to X$ gilt:
$\forall_{x, \tilde{x} \in X} \, \langle T \, x, \tilde{x} \rangle = \langle x, T \, \tilde{x} \rangle$.

Der Leser verwende den Satz vom abgeschlossenen Graphen, um zu zeigen, daß $T \: X \to X$ stetig ist.
\end{UA}

\begin{UA}
Seien $X$ ein $\K$-Hilbertraum und $P \: X \to X$ eine Projektion.

Zeige, daß die folgenden Aussagen paarweise äquivalent sind:
\begin{itemize}
\item[(i)] Es gibt einen abgeschlossenen $\K$-Untervektorraum $Y$ von $X$ mit $P = \pi_Y$.
\item[(ii)] $P \in {\rm Hom}_{\K}(X,X)$ ist selbstadjungiert, vgl.\ die vorherige Übungsaufgabe.
\item[(iii)] $P \in \L_{\K}(X,X)$ und $\|P\| \le 1$.
\item[(iv)] $P(X) \perp {\rm Kern} \, P$.
\end{itemize}
\end{UA}

\begin{UA}
Beweise Satz \ref{FA.8.17}.
\end{UA}

\begin{UA}[Summierbarkeit von Teilfamilien] 
Seien $X$ ein $\K$-Ba\-nach\-raum, $I$ eine nicht-leere Menge und $(x_i)_{i \in I}$ eine summierbare Familie von Elementen von $X$.
Beweise die folgenden Aussagen:
\begin{itemize}
\item[(i)] Ist $I_0$ eine nicht-leere Teilmenge von $I$, so ist auch $(x_i)_{i \in I_0}$ eine summierbare Familie von Elementen von $X$.
\item[(ii)] Sind $J$ eine Menge und $I_j$ für jedes $j \in J$ eine nicht-leere Teilmenge von $I$ mit $I = \bigcupdot_{j \in J} I_j$ sowie $\forall_{j \in J} \, y_j := \sum_{i \in I_j} x_i$ -- beachte (i) --, so ist $(y_j)_{k \in K}$ summierbar und das \emph{Assoziativgesetz}
$$ \sum_{j \in J} y_j = \sum_{i \in I} x_i $$
gilt.
\end{itemize}
\end{UA}

\begin{UA}
Seien $X$ ein $\K$-Ba\-nach\-raum, $I$ eine nicht-leere Menge und $(x_i)_{i \in I}$ eine Familie von Elementen von $X$.
Ferner seien $J$ eine Menge und $I_j$ für jedes $j \in J$ eine nicht-leere Teilmenge von $I$ mit $I = \bigcupdot_{j \in J} I_j$ derart, daß $(y_i)_{i \in J_j}$ für alle $j \in J$ summierbar ist.

Zeige, daß $(x_i)_{i \in I}$ dann nicht summierbar sein muß.
\end{UA}

Tip: $X := \R$, $I := J := \N$, $x_i := \left\{ \begin{array}{rl} -2, & \mbox{falls } i \equiv 0 (5), \\ 1, & \mbox{sonst,} \end{array} \right\}$ für alle $i \in I = \N$ und $I_{2j} := \{ 5j + 1, 5j + 2, 5(2j + 1)\}$ sowie $I_{2j+1} := \{5j + 3, 5j + 4, 10j\}$ für alle $j \in J = \N$.

\begin{UA}
Es seien $n \in \N_+$, $X$ ein $n$-dimensionaler normierter $\K$-Vek\-tor\-raum, $I$ eine nicht-leere Menge und $(x_i)_{i \in I}$ eine summierbare Familie von Elementen von $X$.

Zeige, daß $(x_i)_{i \in I}$ absolut summierbar ist.
\end{UA}

Tip: Begründe zunächst, daß man $\K = \R$ annehmen kann.
Zerlege $I$ dann im Falle $X = \R$ in $I_+ := \{ i \in I \, | \, x_i \ge 0\}$ sowie $I_- := \{ i \in I \, | \, x_i < 0\}$ und betrachte die Summe der Familien $(x_i)_{i \in I_1}$ und $(x_i)_{i \in I_2}$.
Überlege anschließend zum Beweis des allgemeinen Falles, daß es genügt, die euklidische Norm $\| \ldots \|_2$ auf $X = \R^n$ zu betrachten.

\begin{UA}
Beweise Lemma \ref{FA.8.24}.
\end{UA}

\begin{UA}
Seien $X$ ein $\K$-Prähilbertraum, $I$ eine nicht-leere Menge und $(e_i)_{i \in I}$ ein Orthonormalsystem.
Zeige die
\begin{itemize}
\item[(i)] (Besselsche Ungleichung) \index{Ungleichung!Besselsche}
\newline $\forall_{x \in X} \, \sum_{i \in I} |\langle x, e_i \rangle|^2 \le \|x\|^2$
\item[(ii)] und für $x \in X$ genau dann die \emph{Parsevalsche Gleichung}\index{Gleichung!Parsevalsche} $\sum_{i \in I} |\langle x, e_i \rangle|^2 = \|x\|^2$, wenn $x = \sum_{i \in I} \langle x, e_i \rangle \, e_i$ erfüllt ist.
\end{itemize}
\end{UA}

Tip: Für jedes $x \in X$ und jede endliche Teilmenge $J$ von $I$ gilt $$ \big\langle \sum_{j \in J} \langle x, e_j \rangle \, e_j, x - \sum_{j \in J} \langle x, e_j \rangle \, e_j \big\rangle = 0, $$ also folgt aus dem Satz des \textsc{Pythagoras} ...

\begin{UA}\label{FA.8.17.B} 
Seien $X$ ein $\K$-Vektorraum, $A, B$ Teilmengen von $X$, $I$ eine beliebige Menge und $(A_i)_{i \in I}$ eine Familie von Teilmengen von $X$.

Beweise die folgenden Aussagen:
\begin{itemize}
\item[(i)] ${\rm Span}_{\K} \, \emptyset = \{0\}$,
\item[(ii)] $A \subset {\rm Span}_{\K} A$, und es gilt Gleichheit, wenn $A$ ein $\K$-Untervektorraum von $X$ ist,
\item[(iii)] $A \subset B \Longrightarrow {\rm Span}_{\K} A \subset {\rm Span}_{\K} B$,
\item[(iv)] $A \subset {\rm Span}_{\K} B \Longrightarrow {\rm Span}_{\K} A \subset {\rm Span}_{\K} B$,
\item[(v)] $\left( \forall_{i,j \in I} \exists_{k \in I} \, A_i \cup A_j \subset A_k \right) \Longrightarrow {\rm Span}_{\K} \left( \bigcup_{i \in I} A_i \right) = \bigcup_{i \in I} {\rm Span}_{\K} A_i$,
\item[(vi)] ${\rm Span}_{\K} \left( \bigcup_{i \in I} A_i \right) = {\rm Span} \left( \bigcup_{i \in I} {\rm Span}_{\K} A_i \right),$
\item[(vii)] ${\rm Span}_{\K} A = \bigcap_{\substack{Y \text{ $\K$-Untervektor-} \\ \text{raum von } X \text{ mit } A \subset Y}} Y$.
\end{itemize}
\end{UA}

\begin{UA}[Gram-Schmidtsches Orthonormalaisierungsverfahren]\index{Gram-Schmidtsches Or\-tho\-nor\-ma\-li\-sie- \newline rungs\-ver\-fahren}\index{Orthonormalisierungsverfahren nach\newline \textsc{Gram-Schmidt}}
Es seien $X$ ein $\K$-Prähilbertraum, $I$ entweder gleich $\{1, \ldots, n\}$ für ein $n \in \N_+$ oder gleich $\N_+$ und $(x_i)_{i \in I}$ eine $\K$-linear unabhängige Familie von Elementen von $X$.

Zeige, daß durch
$$ e_1 := \frac{x_1}{\|x_1\|} ~~ \mbox{ und } ~~ \forall_{i \in I, i > 1} \, e_i := \frac{x_i - \sum_{j=1}^{i-1} \langle x_i, e_j \rangle \, e_j}{\| x_i - \sum_{j=1}^{i-1} \langle x_i, e_j \rangle \, e_j \|} $$
ein Orthonormalsystem $(e_i)_{i \in I}$ von $X$ mit ${\rm Span}_{\K} \{e_i \, | \, i \in I \} = {\rm Span}_{\K} \{x_i \, | \, i \in I \}$ definiert wird.
\end{UA}

\begin{UA} \label{FA.8.VONS}
Beweise die folgenden beiden Aussagen:
\begin{itemize}
\item[(i)] Es seien $f_0 := \frac{1}{\sqrt{2 \pi}} \, 1_{{[}-\pi, \pi{]}}$ und für jedes $k \in \N_+$ 
$$ f_{-k} := \frac{1}{\sqrt{\pi}} \, \cos(k \, x)|_{{[}-\pi,\pi{]}} ~~ \mbox{ sowie } ~~ f_k := \frac{1}{\sqrt{\pi}} \, \sin(k \, x)|_{{[}-\pi,\pi{]}}. $$
\newline Dann ist $(\mathbf{f}_k)_{k \in \Z}$ eine Hilbertbasis von $L^2_{\R}({[}-\pi,\pi{]}, \mu_1)$, wobei $\mathbf{f}_k$ für jedes $k \in \Z$ die Äquivalenzklasse von $f_k \in \L^2_{\R}({[}-\pi,\pi{]}, \mu_1)$ bezeichne.

\emph{Tip: Aufgabe \ref{FA.3.W2}, Satz \ref{FA.5.19} (ii) und Korollar \ref{FA.5.23.K}.}
\item[(ii)] Das Orthonormalsystem $\big( \frac{1}{\sqrt{2 \pi}} \, \left. \e^{\i k x} \right|_{{[}- \pi, \pi{]}} \big)_{k \in \Z}$ von $\L^2_{\C}({[}-\pi,\pi{]}, \mu_1)$ induziert eine Hilbertbasis von $L^2_{\C}({[}-\pi,\pi{]}, \mu_1)$.
\end{itemize}
\end{UA}

\begin{UA}
Seien $n \in \N_+$ und $J \in \mathfrak{I}_n$ ein nicht-entartetes Intervall.

Zeige, daß der separable $\K$-Hilbertraum $L^2_{\K}(J,\mu_n)$ unendlich-dimensional ist.
\end{UA}

\cleardoublepage
\begin{appendix}

\section{Das Lemma von \textsc{Zorn}} \label{FAnaA}

\begin{Def} \label{D.3.10} 
Sei $\mathcal{X}$ eine Menge.
\begin{itemize}
\item[(i)] Eine \emph{Ordnung von $\mathcal{X}$}\index{Ordnung} ist eine zweistellige Relation $\le$ (d.h.\ eine Teilmenge von $\mathcal{X} \times \mathcal{X}$, wobei wir $x \le y$ für $(x,y) \in \,\, \le$ schreiben) derart, daß für alle $x,y,z \in \mathcal{X}$ gilt:

$\begin{array}{ll}
x \le x, & ~~~~~~~ \mbox{(Reflexivität)} \\
x \le y \, \wedge \, y \le x \Longrightarrow x = y, & ~~~~~~~ \mbox{(Antisymmetrie)} \\
x \le y \, \wedge \, y \le z \Longrightarrow x \le z. & ~~~~~~~ \mbox{(Transitivität)}
\end{array}$

$x,y \in \mathcal{X}$ heißen \emph{miteinander vergleichbar}, wenn gilt $x \le y$ oder $y \le x$.

Eine Ordnung heißt \emph{Totalordnung}\index{Totalordnung}, wenn je zwei Elemente miteinander vergleichbar sind.
\item[(ii)] Ist $\le$ eine Ordnung von $\mathcal{X}$, so definieren wir für alle $x,y \in \mathcal{X}$
\begin{eqnarray*}
x \ge y & : \Longleftrightarrow & y \le x, \\
x < y & : \Longleftrightarrow & ( x \le y \, \wedge \, x \neq y ), \\
x > y & : \Longleftrightarrow & y < x.
\end{eqnarray*}
\item[(iii)] Seien $\le$ eine Ordnung von $\mathcal{X}$ und $a \in \mathcal{X}$.

$a$ heißt \emph{größtes Element von $\mathcal{X}$}\index{Element!größtes} $: \Longleftrightarrow$ $\forall_{x \in \mathcal{X}} \, a \ge x$.
\newline \noindent
Insbesondere ist ein größtes Element mit allen $x \in \mathcal{X}$ vergleichbar.

$a$ heißt \emph{maximales Element von $\mathcal{X}$}\index{Element!maximales} $: \Longleftrightarrow$ $\forall_{x \in \mathcal{X}} \, ( x \ge a \Rightarrow x = a)$,
\newline \noindent d.h.\ ein maximales Element ist $\ge$ allen Elementen, mit denen es vergleichbar ist.

Es gilt:
\begin{equation*}
\begin{minipage}{0.935\textwidth}
$a$ größtes Element von $\mathcal{X}$ $\Longrightarrow$ $a$ maximales Element von $\mathcal{X}$.

Die Umkehrung ist i.a.\ falsch.
\end{minipage}
\end{equation*}
\begin{equation*}
\begin{minipage}{0.935\textwidth}
Es existiert höchstens ein größtes Element von $\mathcal{X}$. Dagegen können viele maximale Elemente von $\mathcal{X}$ existieren, s.u.\ Beispiele.
\end{minipage}
\end{equation*}

Analog definiert man die Begriffe \emph{kleinstes Element von $\mathcal{X}$}\index{Element!kleinstes} und \emph{minimales Element von $\mathcal{X}$}\index{Element!minimales}.
\item[(iv)] Seien $\le$ eine Ordnung von $\mathcal{X}$, $M$ eine nicht-leere Teilmenge von $\mathcal{X}$ und $a \in \mathcal{X}$.

$a$ heißt \emph{obere Schranke von $M$ (in $\mathcal{X}$)}\index{Schranke!obere}, wenn gilt: $\forall_{x \in M} \, a \ge x$.

$a$ heißt \emph{kleinste obere Schranke von $M$ (in $\mathcal{X}$)}\index{Schranke!kleinste obere}, wenn $a$ obere Schranke von $M$ ist und für jede obere Schranke $b \in \mathcal{X}$ von $M$ gilt $a \le b$, d.h.\ $a$ ist kleinstes -- nicht nur minimales !! -- Element der Menge aller oberen Schranken von $M$.

Es gilt:
\begin{equation*}
\begin{minipage}{0.935\textwidth}
Es existiert höchstens eine kleinste obere Schranke von $M$, die wir im Falle ihrer Existenz mit $\boxed{\sup M}$ bezeichnen.
\end{minipage}
\end{equation*}

Analog definiert man die Begriffe \emph{untere Schranke}\index{Schranke!untere}, \emph{größte untere Schranke}\index{Schranke!größte untere} und das Symbol $\boxed{\inf M}$.
\end{itemize}

\begin{Bsp*} $\,$
\begin{itemize}
\item[1.)] Auf $\mathcal{X} := \N_+$ wird durch $$\forall_{m,n \in \N_+} \, m \le_* n : \Longleftrightarrow m | n $$ eine Ordnung definiert, die keine Totalordnung ist, da z.B.\ weder $5 | 7$ noch $7 | 5$ gilt. 
$1$ ist kleinstes Element von $\N_+$ bzgl $\le_*$.
\item[2.)] $\N_+ \setminus \{1\}$ ist durch $\le_*$ ebenfalls geordnet, und $\N_+ \setminus \{1\}$ besitzt kein kleinstes Element.
Die Menge der minimalen Elemente ist die Menge aller Primzahlen.
\item[3.)] Die übliche Ordnung $\le$ von $\R$ ist eine Totalordnung.
\end{itemize}
\end{Bsp*}
\end{Def}

\begin{HS}[Lemma von \textsc{Zorn}]\index{Lemma!von \textsc{Zorn}} \label{D.3.11}
Sei $\mathcal{X} \neq \emptyset$ eine bzgl.\ $\le$ geordnete Menge mit der Eigenschaft, daß jede bzgl.\ $\le$ totalgeordnete nicht-leere Teilmenge von $\mathcal{X}$ eine obere Schranke in $\mathcal{X}$ besitzt.

Dann existiert (mindestens) ein maximales Element von $\mathcal{X}$.
\end{HS}

Die folgende Formulierung ist zum Lemma von \textsc{Zorn} äquivalent:

\emph{Ist $\mathcal{X} \neq \emptyset$ eine geordnete Menge ohne maximales Element, so existiert eine totalgeordnete Teilmenge $M \neq \emptyset$ von $\mathcal{X}$, die keine obere Schranke besitzt.}

Wir bereiten den Beweis des Hauptsatzes durch das folgende Lemma vor:

\begin{Lemma} \label{FA.A.C.3}
Sei $\mathcal{X} \neq \emptyset$ eine bzgl.\ $\le$ geordnete Menge mit der Eigenschaft, daß jede bzgl.\ $\le$ totalgeordnete nicht-leere Teilmenge von $\mathcal{X}$ eine kleinste obere Schranke in $\mathcal{X}$ besitzt.
Ferner sei $f\: \mathcal{X} \to \mathcal{X}$ eine Abbildung mit $\forall_{x \in \mathcal{X}} \, x \le f(x)$.

Dann existiert $x_* \in \mathcal{X}$ mit $f(x_*) = x_*$.
\end{Lemma}

\textit{Beweis des Lemmas.} Nach Voraussetzung existiert $a \in \mathcal{X}$.
Wir setzen 
$$ \mathcal{Y} := \{ y \in \mathcal{X} \, | \, y \ge a \}. $$
Dann gilt $a \in \mathcal{Y}$ ($a$ ist sogar kleinstes Element von $\mathcal{Y}$) und $\forall_{y \in \mathcal{Y}} \, a \le y \le f(y)$, also $f(\mathcal{Y}) \subset \mathcal{Y}$.
$\mathcal{Y}$ ist daher wegen der Voraussetzung des Lemmas im folgenden Sinne zulässig.

Eine Teilmenge $A$ von $\mathcal{Y}$ heißt \emph{zulässig}, wenn gilt
\begin{itemize}
\item[a)] $a \in A$,
\item[b)] $f(A) \subset A$,
\item[c)] $\emptyset \neq B \subset A$ und $B$ totalgeordnet bzgl.\ $\le$ $\Longrightarrow$ $\sup B \in A$,
\newline beachte, daß $\sup B$ nach Voraussetzung in $\mathcal{X}$ existiert.
\end{itemize}

Wir setzen 
$$ M := \bigcap_{A \subset \mathcal{Y}, A \text{ zulässig}} A $$
und werden zeigen:
\begin{gather}
\mbox{$M$ ist zulässig,} \label{D.3.11.L1} \\
\mbox{$M$ ist totalgeordnet.} \label{D.3.11.L2}
\end{gather}
Hieraus und der Voraussetzung des Lemmas folgt die Existenz einer kleinsten oberen Schranke $x_* := \sup M$ von $M$ in $\mathcal{X}$ mit $x_* \in M$.
Dann gilt für alle $y \in M$
$$ y \le \underbrace{f(y)}_{\stackrel{(\ref{D.3.11.L1})}{\in} M} \le x_* $$
und somit wegen $x_* \in M$: $f(x_*) = x_*$.

Zu zeigen bleiben (\ref{D.3.11.L1}) und (\ref{D.3.11.L2}).

Zu (\ref{D.3.11.L1}): a) $\mathcal{Y}$ ist zulässig, und jede zulässige Menge enthält $a$, also gilt $a \in M$.

b) Sei $y \in M$, also folgt aus der Definition von $M$, daß $\forall_{A \subset \mathcal{Y}, A \text{ zulässig}} \, y \in A$ gilt.
Wegen der Eigenschaft b) in der Definition der Zulässigkeit gilt dann auch $\forall_{A \subset \mathcal{Y}, A \text{ zulässig}} \, f(y) \in A$, d.h.\ $f(y) \in M$.
Damit ist gezeigt $f(M) \subset M$.

c) Sei $B \neq \emptyset$ eine totalgeordnete Teilmenge von $M = \bigcap_{A \subset \mathcal{Y}, A \text{ zulässig}} A$.
Dann gilt $\forall_{A \subset \mathcal{Y}, A \text{ zulässig}} \, B \subset A$, also nach Eigenschaft c) in der Definition der Zulässigkeit $\forall_{A \subset \mathcal{Y}, A \text{ zulässig}} \, \sup B \in A$, und somit folgt aus der Definition von $M$: $\sup B \in A$.

Damit ist (\ref{D.3.11.L1}) gezeigt.

Zu (\ref{D.3.11.L2}): Wir nennen ein Element $e \in M$ \emph{extremales Element}, wenn
$$ \forall_{y \in M} \, ( y < e \Longrightarrow f(y) \le e ) $$
gilt, und setzen für jedes extremale Element $e \in M$
$$ M_e := \{ y \in M \, | \, y \le e \, \vee \, f(e) \le y \}. $$
Dann folgen nacheinander 
\begin{gather}
e \in M \mbox{ extremal } \Longrightarrow M_e = M, \label{D.3.11.L3} \\
\forall_{y \in M} \, y \mbox{ extremal.} \label{D.3.11.L4}
\end{gather}

{[} Zu (\ref{D.3.11.L3}): Sei $e \in M$ extremal.
Wir zeigen, daß $M_e$ zulässig ist:

a) $a$ ist kleinstes Element von $\mathcal{Y}$, und nach (\ref{D.3.11.L1}) gilt $a \in M \subset \mathcal{Y}$.
Hieraus folgt $a \le e$, d.h.\ $a \in M_e$.

b) Sei $y \in M_e$, also $y \le e$ oder $f(e) \le y$.

Im Falle $y < e$ folgt aus der Extremalität von $e$, daß gilt $f(y) \le e$, also $f(y) \in M_e$.

Im Falle $y = e$ gilt trivialerweise $f(y) = f(e)$, also $f(e) \le f(y)$ und somit $f(y) \in M_e$.

Im Falle $f(e) \le y$ gilt wegen der Voraussetzung an $f$
$$ f(e) \le y \le f(y), $$
also wieder $f(y) \in M_e$.

Damit ist gezeigt $f(M_e) \subset M_e$.

c) Sei $B \neq \emptyset$ eine totalgeordnete Teilmenge von $M_e$.

Im Falle $\forall_{b \in B} \, b \le e$ folgt $\sup B \le e$ (denn in diesem Fall ist $e$ eine obere Schranke von $B$ und $\sup B$ ist die kleinste solche) und somit $\sup B \in M_e$.

Existiert andernfalls $b_0 \in B \subset M_e$ mit $b_0 > e$, so folgt aus der Defintion von $M_e$, daß gilt $f(e) \le b_0 \le \sup B$, also $\sup B \in M_e$.

Insgesamt haben wir die Zulässigkeit der Teilmenge $M_e$ von $M$ bewiesen.

Nun ist $M$ der Schnitt aller zulässigen Mengen, also $M \subset M_e$, daher folgt $M = M_e$.

Zu (\ref{D.3.11.L4}): Wir setzen
$$ E := \{ e \in M \, | \, e \mbox{ extremal} \} $$
und zeigen, daß $E$ zulässig ist:

a) $a$ ist kleinstes Element von $\mathcal{Y} \supset M$, d.h.\ die Prämisse der Aussage
$$ y < a \Longrightarrow f(y) \le a $$
ist für jedes $y \in M$ falsch, d.h.\ die Aussage ist wahr.
Hieraus folgt $a \in E$.

b) Es sei $e \in E$. 
Wir wollen nachweisen, daß dann auch $f(e) \in E$ gilt, d.h.\ $\forall_{y \in M} \, y < f(e) \Rightarrow f(y) \le f(e)$.
Sei also $y \in M$ mit $y < f(e)$.
Nach (\ref{D.3.11.L3}) gilt $y \in M_e$, d.h.\ entweder $y < e$ oder $y = e$ oder $f(e) \le y$, und letzteres kann nicht eintreten, da $y < f(e)$.

Ist $y < e$, so folgt aus der Extremaleigenschaft von $e$ und der Voraussetzung an $f$, daß gilt $f(y) \le e \le f(e)$.

Ist $y = e$, so gilt trivialerweise $f(y) = f(e)$.

c) Sei $B \neq \emptyset$ eine totalgeordnete Teilmenge von $E$.
Zu zeigen ist $\sup B \in E$, d.h.\
$$ \forall_{y \in M} \, y < \sup B \Longrightarrow f(y) \le \sup B. $$
Sei also $y \in M$ mit $y < \sup B$.

Im Falle $\forall_{b \in B} \, f(b) \le y$ gälte auch $\forall_{b \in B} \, b \le f(b) \le y$, d.h.\ $y$ wäre eine obere Schranke von $B$.
Daher folgte (weil $\sup B$ kleinste obere Schranke von $B$ ist) $\sup B \le y$, im Widerspruch zu $y < \sup B$.

Daher existiert $b_0 \in B \subset E$ mit $f(b_0) > y$.
Nach (\ref{D.3.11.L3}) gilt $M = M_{b_0}$, d.h.\ $y \in M_{b_0}$, also $y \le b_0 \, \vee \, f(b_0) \le y$.
Daher muß $y \le b_0$ gelten.

Entweder gilt nun $y < b_0$ und (da $b_0$ extremal) somit $f(y) \le b_0 \le \sup B$, oder es gilt $y = b_0$, also $f(y) = f(b_0)$.
Aus $B \subset E \subset M$ und (\ref{D.3.11.L1}) folgt dann $$\sup B \in M \stackrel{(\ref{D.3.11.L3})}{=} M_{b_0},$$ und somit gilt nach Definition von $M_{b_0}$ wegen $b_0 = y < \sup B$
$$f(y) = f(b_0) \le \sup B.$$

\pagebreak
Analog zum Ende des Beweises von (\ref{D.3.11.L3}) folgt aus der Zulässigkeit von $E$ und der Definition von $M$: $M = E$.
Die Definition von $E$ ergibt (\ref{D.3.11.L4}). {]}

Aus (\ref{D.3.11.L3}), (\ref{D.3.11.L4}) folgt jetzt, daß $M$ totalgeordnet ist:

Seien nämlich $y_1, y_2 \in M$. 
Nach (\ref{D.3.11.L4}) ist $y_1$ extremal, also folgt aus (\ref{D.3.11.L3}) $y_2 \in M_{y_1}$, d.h.\
$$ y_2 \le y_1 ~~~ \mbox{ oder } ~~~ \underbrace{f(y_1) \le y_2}_{\Rightarrow \, y_1 \, \le \, f(y_1) \, \le \, y_2}. $$

Somit ist auch (\ref{D.3.11.L2}) bewiesen, und hierauf hatten wir das Lemma zurückgeführt. \q
\A
\textit{Beweis des Hauptsatzes.} 1.\ Fall: Jede nicht-leere totalgeordnete Teilmenge von $\mathcal{X}$ hat eine kleinste obere Schranke in $\mathcal{X}$.
Wir nehmen dann an, daß $\mathcal{X}$ kein maximales Element besitzt.
Dann gilt
$$ \forall_{x \in \mathcal{X}} \, \mathcal{X}_{>x} := \{ y \in \mathcal{X} \, | \, y > x \} \neq \emptyset. $$
Aus dem \underline{Auswahlaxiom} folgt daher $\bigtimes_{x \in \mathcal{X}} \mathcal{X}_{>x} \neq \emptyset$, d.h.\ es existiert eine Abbildung $f \: \mathcal{X} \to \mathcal{X}$ mit $\forall_{x \in \mathcal{X}} \, f(x) \in \mathcal{X}_{>x}$.
Es gilt also
\begin{equation} \label{D.3.11.S}
\forall_{x \in \mathcal{X}} \, f(x) > x.
\end{equation}
Aus (\ref{D.3.11.S}) und \ref{FA.A.C.3} folgt die Existenz von $x_* \in \mathcal{X}$ mit $f(x_*) = x_*$, im Widerspruch zu (\ref{D.3.11.S}).

2.\ Fall: Jede nicht-leere totalgeordnete Teilmenge von $\mathcal{X}$ hat eine obere Schranke in $\mathcal{X}$. (D.i.\ ist die Voraussetzung des Hauptsatzes.)
Sei dann 
$$ \widetilde{\mathcal{X}} :=  \left\{ A \in \mathfrak{P}(\mathcal{X}) \setminus \{ \emptyset \} \, | \, \mbox{$A$ totalgeordnet bzgl.\ $\le$} \right\}. $$

Durch 
$$ \forall_{A,B \in \widetilde{\mathcal{X}}} \, A \preceq B : \Longleftrightarrow A \subset B $$ 
wird eine Ordnung auf $\widetilde{\mathcal{X}}$ definiert, und jede nicht-leere bzgl.\ $\preceq$ totalgeordnete Teilmenge von $\widetilde{\mathcal{X}}$ besitzt eine kleinste obere Schranke in $\widetilde{\mathcal{X}}$.

{[} Denn ist $(A_i)_{i \in I}$ eine nicht-leere bzgl.\ $\preceq$ totalgeordnete Teilmenge von $\widetilde{\mathcal{X}}$, so ist $\bigcup_{i \in I} A_i$ bzgl.\ $\le$ totalgeordnet%
\footnote{Seien $x,y \in \bigcup_{i \in I} A_i$. Dann existieren $i_0,i_1 \in I$ mit $x \in A_{i_0}$ und $y \in A_{i_1}$.
Da $(A_i)_{i \in I}$ totalgeordnet bzgl.\ $\preceq$ ist, gilt ohne Einschränkung $A_{i_0} \subset A_{i_1}$, also $x,y \in A_{i_1}$.
Nun ist aber $A_{i_1}$ als Element von $\widetilde{\mathcal{X}}$ totalgeordnet bzgl.\ $\le$, d.h.\ $x \le y$ oder $y \le x$.}
und offenbar (bzgl.\ $\preceq$) die kleinste obere Schranke von $(A_i)_{i \in I}$ in $\widetilde{\mathcal{X}}$. {]}

Aus dem bereits bewiesenen 1. Fall folgt nun die Existenz eines (bzgl.\ $\preceq$) maximalen Elementes $\widetilde{\mathcal{X}}_0$ von $\widetilde{\mathcal{X}}$.
Insbesondere gilt $\widetilde{\mathcal{X}}_0 \neq \emptyset$, und $\widetilde{\mathcal{X}}_0$ ist eine bzgl.\ $\le$ totalgeordnete Teilmenge von $\mathcal{X}$.
Nach Voraussetzung des Hauptsatzes besitzt $\widetilde{\mathcal{X}}_0$ eine obere Schranke $a \in \mathcal{X}$.
Dann ist $a$ maximales Element von $\mathcal{X}$ bzgl $\le$:

Sei nämlich $x \in \mathcal{X}$ mit $x \ge a$.
Dann ist $\widetilde{\mathcal{X}}_0 \cup \{ x \}$ totalgeordnet bzgl.\ $\le$ mit $\widetilde{\mathcal{X}}_0 \preceq \widetilde{\mathcal{X}}_0 \cup \{ x \}$.
Wegen der Maximalität von $\widetilde{\mathcal{X}}_0$ bzgl.\ $\preceq$ gilt dann $\widetilde{\mathcal{X}}_0 = \widetilde{\mathcal{X}}_0 \cup \{ x \}$, d.h.\ $x \in \widetilde{\mathcal{X}}_0$.
Nun ist $a$ obere Schranke von $\widetilde{\mathcal{X}}_0$, also folgt $a \ge x$. 
Nach Wahl von $x$ gilt aber auch $x \ge a$, also folgt $a = x$.

Damit ist das Lemma von \textsc{Zorn} vollständig bewiesen. \q
\A Wir haben das Zornsche Lemma unter Verwendung des Auswahlaxiomes bewiesen.
Tatsächlich ist es zum Auswahlaxiom äquivalent.

Es gilt die folgende Verschärfung des Lemmas von \textsc{Zorn}:

\begin{Satz} \label{FA.A.0}
Sei $\mathcal{X} \neq \emptyset$ eine bzgl.\ $\le$ geordnete Menge mit der Eigenschaft, daß jede bzgl.\ $\le$ totalgeordnete nicht-leere Teilmenge von $\mathcal{X}$ eine obere Schranke in $\mathcal{X}$ besitzt.
Ferner sei $M$ eine solche Teilmenge von $\mathcal{X}$.

Dann existiert (mindestens) ein maximales Element von $\mathcal{X}$, welches eine obere Schranke von $M$ ist.
\end{Satz}

\textit{Beweis.} Wir setzen
$$ \widetilde{\mathcal{X}} := \{ x \in \mathcal{X} \, | \, \forall_{y \in M} \, y \le x \}. $$ 
Da $M$ nach Voraussetzung eine obere Schranke besitzt, gilt $\widetilde{\mathcal{X}} \ne \emptyset$.
Als Teilmenge von $\mathcal{X}$ besitzt jede bzgl.\ $\le$ totalgeordnete nicht-leere Teilmenge $\widetilde{M}$ von $\widetilde{\mathcal{X}}$ eine obere Schranke $a \in \mathcal{X}$, d.h.\ $\forall_{x \in \widetilde{M}} \, x \le a$, also gilt wegen $\widetilde{M} \subset \widetilde{\mathcal{X}}$
$$ \forall_{x \in \widetilde{M}} \forall_{y \in M} \, y \le x \le a $$
und somit $a \in \widetilde{\mathcal{X}}$.

Damit erfüllt $\widetilde{\mathcal{X}}$ die Voraussetzung des Lemmas von \textsc{Zorn} \ref{D.3.11}, besitzt also ein maximales Element $\tilde{x} \in \widetilde{\mathcal{X}}$.
Daß $\tilde{x}$ eine obere Schranke von $M$ ist, folgt sofort aus der Definition von $\widetilde{\mathcal{X}}$.
\q
\A Typische Anwendungen des Zornschen Lemmas sind die unten genannten Sätze.

\begin{Def}[(Hamel-)Basis eines Vektorraumes] \label{FA.A.1}
Es seien $\k$ ein Körper und $X$ ein $\k$-Vektorraum.

Eine Teilmenge $\mathfrak{B}$ von $X$ heißt eine \emph{Basis von $X$}\index{Basis} genau dann, wenn jeweils endlich viele Elemente von $\mathfrak{B}$ $\k$-linear unabhängig sind und sich jedes Element von $X \setminus \{0\}$ bis auf die Reihenfolge der Summanden eindeutig als \emph{$\k$-Linearkombination von Elementen aus $\mathfrak{B}$} schreiben läßt.
Letzteres bedeutet per definitionem, daß zu jedem $x \in X \setminus \{0\}$ sowohl $k \in \N_+$ als auch $\lambda_1, \ldots, \lambda_k \in \k \setminus \{0\}$ sowie $b_1, \ldots, b_k \in \mathfrak{B}$ mit $x = \sum_{i=1}^k \lambda_i \, b_i$ existieren.
(Linearkombinationenen sind also \underline{endliche} Summen! Es sei darauf hingewiesen, daß sogar der Begriff der ,,Konvergenz einer Reihe`` i.a.\ auch gar keinen Sinn machen würde.)

\begin{Bem*} $\,$
\begin{itemize}
\item[1.)] Um Verwechselungen zu vermeiden, nennt man eine Basis im gerade definierten Sinne gelegentlich eine \emph{Hamelbasis von $X$}.\index{Basis!Hamel-}
Wir verwenden die Begriffe ,,Basis`` und ,,Hamelbasis`` synonym.
\item[2.)] Ist $X$ ein normierter $\K$-Vektorraum, so heißt eine Folge $(b_i)_{i \in \N}$ in $X$ eine \emph{Schauderbasis}\index{Basis!Schauder-} genau dann, wenn zu jedem $x \in X$ genau eine Folge $(\lambda_i)_{i \in \N}$ in $\K$ mit $x = \sum_{i=0}^{\infty} \lambda_i \, b_i$ existiert.
\end{itemize}
\end{Bem*}
\end{Def}

\begin{Satz} \label{FA.A.2}
Seien $\k$ ein Körper und $X$ ein $\k$-Vektorraum.

Dann besitzt $X$ eine Basis. 
\end{Satz}

\textit{Beweisskizze.} Sei $\mathcal{B} \subset \mathfrak{P}(X)$ die Menge aller $\k$-linear unabhängigen Teilmengen von $X$, d.h.\ per definitionem, daß für alle $B \in \mathcal{B}$ jeweils endlich viele Elemente von $B$ $\k$-linear unabhängig sind.
Dann ist $\mathcal{B} \ne \emptyset$ mit der Teilmengenrelation $\subset$ eine geordnete Menge, und jede bzgl.\ $\subset$ totalgeordnete Menge $M \subset \mathcal{B}$ besitzt $\bigcup_{B \in M} B$ als obere Schranke.
Nach dem Zornschen Lemma existiert ein bzgl.\ $\subset$ maximales Element $\mathfrak{B}$ von $\mathcal{B}$, welches offenbar eine Basis von $X$ ist. \q

\begin{Bem*}
Wir haben soeben bewiesen, daß der $\K$-Vektorraum $\K^{\N}$ aller Folgen in $\K$ eine Hamelbasis besitzt.
$ \{ (1,0,0, \ldots), (0,1,0, \ldots), \ldots \}$ ist keine solche!
\end{Bem*}

Aus dem letzten Satz folgt übrigens das überraschende Resultat, daß es additive Funktionen $f \: \R \to \R$ gibt\footnote{Eine Funktion $f \: \R \to \R$ heißt \emph{additiv}\index{Funktion!additive}, wenn für alle $s,t \in \R$ gilt $f(s+t) = f(s) + f(t)$.}, die nicht von der Gestalt $f = a \, x$ mit $a \in \R$ sind.

\begin{Satz} \label{FA.A.3}
Es seien $R$ ein Ring mit Einselement, der nicht der Nullring ist, und $\mathfrak{a}$ ein echtes beidseitiges Ideal von $R$.

Dann existiert ein maximales beidseitiges Ideal von $R$, welches $\mathfrak{a}$ enthält.
\end{Satz}

\textit{Beweisskizze.} Bzgl.\ der Teilmengenrelation ist die Menge aller echten beidseitigen Ideale von $R$ eine geordnete Menge.
Sind $I$ eine nicht-leere Menge und $\{ \mathfrak{a}_i \, | \, i \in I \}$ eine totalgeordnete Teilmenge der Menge aller echten beidseitigen Ideale, die $\mathfrak{a}$ enthalten, so ist $\bigcup_{i \in I} \mathfrak{a}_i$ ein echtes beidseitiges Ideal von $R$, das $\mathfrak{a}$ enthält.
Daher folgt die Behauptung aus dem Zornschen Lemma. \q

\subsection*{Übungsaufgaben}

\begin{UA}
Zeige, daß aus dem Lemma von \textsc{Zorn} das Auswahlaxiom folgt.
\end{UA}

\begin{UA}
Zeige, daß eine additive Funktion $f \: \R \to \R$, die nicht von der Gestalt $f = a \, x$ mit $a \in \R$ ist, existiert.
\end{UA}

\cleardoublepage

\section{Der Satz von \textsc{Tychonoff}} \label{FAnaC}

\begin{Def}[Produktopologie] \index{Topologie!Produkt-} \label{FA.C.P}
Seien $I$ eine Menge und $X_i$ für jedes $i \in I$ ein topologischer Raum.

Die \emph{Produktopolologie für 
$$ \bigtimes_{\iota \in I} X_{\iota} = \left\{ f \in \left( \bigcup_{\iota \in I} X_{\iota} \right)^I \, | \, \forall_{i \in I} \, f(i) \in X_i \right\} $$}
ist per definitionem die schwache Topologie für $\bigtimes_{\iota \in I} X_{\iota}$ bzgl.\ $\{ \pi_i \, | \, i \in I \}$, wobei
\begin{equation} \label{FA.C.P.S}
\pi_i \: \bigtimes_{\iota \in I} X_{\iota} \longrightarrow X_i, ~~ f \longmapsto f(i),
\end{equation}
für jedes $i \in I$ die kanonische Projektion bezeichne, d.h.\ genau, daß diese Topologie die gröbste für $\bigtimes_{\iota \in I} X_{\iota}$ ist, derart, daß alle $\pi_i$, $i \in I$, stetig sind, vgl.\ \ref{FA.6.2}.

\begin{Bem*}
Schreibt man $X$ anstelle von $I$ und sind alle $X_{\iota}$, $\iota \in I$, gleich einem (festen) topologischen Raum $Y$, so gilt
$$ \bigtimes_{\iota \in I} X_{\iota} = Y^X. $$
Aus \ref{FA.6.2.S} (ii) folgt dann, daß eine Folge $(f_k)_{k \in \N}$ in $Y^X$ genau dann bzgl.\ der Produkttopologie für $\bigtimes_{\iota \in I} X_{\iota}$ gegen $f \in Y^X$ konvergiert, wenn für jedes $x \in X$ die Folge $(f_k(x))_{k \in \N}$ in $Y$ gegen $f(x)$ konvergiert.
Daher nennt man die Produkttopologie für $\bigtimes_{\iota \in I} X_{\iota}$ in diesem Falle die \emph{Topologie der punktweisen Konvergenz für $Y^X$}.\index{Topologie!der punktweisen Konvergenz}
\end{Bem*}
\end{Def}

\begin{Satz} \label{FA.C.2}
Seien $I$ eine Menge sowie $X_{\iota}$ für alle $\iota \in I$ ein topologischer Raum.
Ferner seien $\pi_i \: \bigtimes_{\iota \in I} X_{\iota} \to X_i$ für jedes $i \in I$ die kanonische Projektion und $\bigtimes_{\iota \in I} X_{\iota}$ mit der Prodokttopologie versehen.
Dann gilt:
\begin{itemize}
\item[(i)] Eine Teilmenge $U$ von $\bigtimes_{\iota \in I} X_{\iota}$ ist genau dann offen, wenn zu jedem $f_0 \in U$ eine endliche Teilmenge $J$ von $I$ und Umgebungen $U_j \in \U(f_0(j),X_j)$ für alle $j \in J$ existieren derart, daß der ,,Zylinder``
$$ \bigcap_{j \in J} \overline{\pi_j}^1(U_j) = \left\{ f \in \bigtimes_{\iota \in I} X_{\iota} \, | \, \forall_{j \in J} \, f(j) \in U_j \right\} \in \U \left( f_0, \bigtimes_{\iota \in I} X_{\iota} \right) $$
eine Teilmenge von $U$ ist.
\item[(ii)] $\pi_i \: \bigtimes_{\iota \in I} X_{\iota} \to X_i$ ist für jedes $i \in I$ eine offene Abbildung.
\end{itemize}
\end{Satz}

\textit{Beweis.} Sei $U$ eine offene Teilmenge von $\bigtimes_{\iota \in I} X_{\iota}$, d.h.\ per definitionem, daß $U$ sich als beliebige Vereinigung von Mengen der Form
$$ \bigcap_{j \in J} \overline{\pi_j}^1(U_j) \mbox{ mit } J \subset I \mbox{ endlich sowie } \forall_{j \in J} \, U_j \in {\rm Top}(X_j) $$
schreiben läßt.
Hieraus folgt offenbar (i).

Sei nun $i \in I$.
Da für alle Teilmengen $A_{\lambda}$ von $\bigtimes_{\iota \in I} X_{\iota}$, $\lambda \in \Lambda$, wobei $\Lambda$ eine beliebige Menge sei, gilt
$$ \pi_i \left( \bigcup_{\lambda \in \Lambda} A_{\lambda} \right) = \bigcup_{\lambda \in \Lambda} \pi_i \left( A_{\lambda} \right), $$
genügt es zum Nachweis von (ii), zu zeigen, daß die Menge
$$ \pi_i \left( \bigcap_{j \in J} \overline{\pi_j}^1(U_j) \right) = \pi_i \left( \left\{ f \in \bigtimes_{\iota \in I} X_{\iota} \, | \, \forall_{j \in J} \, f(j) \in U_j \right\} \right) = \left\{ \begin{array}{cc} U_i, & \mbox{ falls } i \in J, \\ X_i, & \mbox{ falls } i \notin J, \end{array} \right. $$
für alle $J$ und $U_j$, $j \in J$, wie oben offen in $X_i$ ist, und dies ist trivial. \q

\begin{Bem*}
Seien $n \in \N_+$ und speziell $I = \{1, \ldots, n\} \equiv n$.
Des weiteren seien $X_1, \ldots, X_n$ topologische Räume und $\bigtimes_{\iota = 1}^n X_{\iota}$ mit der Produkttopologie versehen.
\begin{itemize}
\item[1.)] Eine Teilmenge $U$ von $\bigtimes_{\iota = 1}^n X_{\iota}$ ist genau dann offen, wenn zu jedem Element $(x_1, \ldots, x_n) \in U$ für jedes $i \in \{1, \ldots,n\}$ ein $U_i \in \U(x_i,X_i)$ existiert derart, daß $\bigtimes_{\iota = 1}^n U_{\iota} \subset U$ gilt.

Insbesondere ist die Definition der Produkttopologie in \ref{FA.C.P} mit der in der Fußnote auf Seite \pageref{FA.4.63} konform.
\item[2.)] Sind $X_1, \ldots, X_n$ sogar normierte $\K$-Vektorräume, so wird die Produkttopologie für $\bigtimes_{\iota = 1}^n X_{\iota}$ durch die Maximumsnorm auf $\bigtimes_{\iota = 1}^n X_{\iota}$ induziert.
\end{itemize}
\end{Bem*}

Für endliches $I$ kann der u.g.\ Hauptsatz auch ohne Verwendung des Auswahlaxiomes gezeigt werden, siehe z.B.\ \cite[Theorem 26.7]{Munkres}.
Alle Beweise benötigen jedoch in irgendeiner Form das Auswahlaxiom, wenn $I$ eine unendliche Menge ist, und \textsc{Kelley} \cite{Kelley} hat bewiesen, daß das Auswahlaxiom aus dem Hauptsatz folgt, also gilt Äquivalenz.

\begin{HS}[von \textsc{Tychonoff}] \label{FA.C.T} \index{Satz!von \textsc{Tychonoff}} $\,$

\noindent \textbf{Vor.:} Es seien $I$ eine Menge und $X_i$ für jedes $i \in I$ ein nicht-leerer kompakter topologischer Raum.

\noindent \textbf{Beh.:} $\bigtimes_{\iota \in I} X_{\iota}$ versehen mit der Produkttopologie ist kompakt.
\end{HS}

Zum Nachweis des Hauptsatzes benötigen wir einige Vorbereitungen.

\begin{Def}[(Maximales) Zentriertes System] \index{zentriertes System} \label{FA.C.ZS}
Sei $X$ eine Menge.
\begin{itemize}
\item[(i)] Eine Teilmenge $\mathcal{A}$ von $\mathfrak{P}(X)$ heißt ein \emph{zentriertes System} genau dann, wenn für alle $n \in \N$ und $A_1, \ldots, A_n \in \mathcal{A}$ gilt $A_1 \cap \ldots \cap A_n \ne \emptyset$.

\begin{Bsp*} $\,$
\begin{itemize}
\item[1.)] Aus der Existenz eines zentrierten Systemes folgt $X \ne \emptyset$, denn der leere Schnitt von Teilmengen ist per definitionem der ganze Raum.
\item[2.)] $\emptyset \subset \mathfrak{P}(X)$ ist ein zentriertes System.
\item[3.)] Ein zentriertes System kann $\emptyset$ nicht als Element besitzen.
\end{itemize}
\end{Bsp*}
\item[(ii)] Ein zentriertes Systen $\mathcal{A} \subset \mathfrak{P}(X)$ heißt \emph{maximal} genau dann, wenn jede echte Obermenge $\mathcal{B}$ von $\mathcal{A}$ mit $\mathcal{B} \subset \mathfrak{P}(X)$ kein zentriertes System ist.
\end{itemize}
\end{Def}

\begin{Lemma} \label{FA.C.3}
Seien $X$ eine nicht-leere Menge und $\mathcal{A} \subset \mathfrak{P}(X)$ ein zentriertes System.

Dann existiert ein maximales zentriertes System $\mathcal{A}_{\max} \subset \mathfrak{P}(X)$ derart, daß $\mathcal{A} \subset \mathcal{A}_{\max}$ gilt.
\end{Lemma}

\textit{Beweis.} Wir definieren die nicht-leere Menge
$$ \mathfrak{A} := \left\{ \widetilde{\mathcal{A}} \in \mathfrak{P}(\mathfrak{P}(X)) \, | \, \mbox{$\mathcal{A} \subset \widetilde{\mathcal{A}}$ und $\widetilde{\mathcal{A}}$ zentriertes System} \right\} $$
und versehen $\mathfrak{A}$ mit der Teilmengenrelation $\subset$ als Ordnung von $\mathfrak{A}$.
Wegen des Lemmas von \textsc{Zorn} \ref{D.3.11} bleibt zu zeigen, daß jede bzgl.\ $\subset$ totalgeordnete nicht-leere Teilmenge $\mathfrak{B}$ von $\mathfrak{A}$ eine obere Schranke in $\mathfrak{A}$ besitzt.
Hierfür wiederum genügt der Nachweis davon, daß für alle solchen $\mathfrak{B}$ gilt
$$ \bigcup_{\mathcal{B} \in \mathfrak{B}} \mathcal{B} \in \mathfrak{A}, $$ 
d.h.\ genau
\begin{gather}
\mathcal{A} \subset \bigcup_{\mathcal{B} \in \mathfrak{B}} \mathcal{B} \subset \mathfrak{P}(X), \label{FA.C.3.1} \\
\bigcup_{\mathcal{B} \in \mathfrak{B}} \mathcal{B} \mbox{ ist ein zentriertes System.} \label{FA.C.3.2}
\end{gather}

Sei also $\mathfrak{B}$ wie oben.
(\ref{FA.C.3.1}) folgt dann sofort daraus, daß $\mathfrak{B}$ eine nicht-leere Teilmenge von $\mathfrak{A}$ ist.

Zu (\ref{FA.C.3.2}): Seien $n \in \N$ und $B_1, \ldots, B_n \in \bigcup_{\mathcal{B} \in \mathfrak{B}} \mathcal{B}$.
Ohne Einschränkung gelte $n \in \N_+$, beachte daß der leere Schnitt von Teilmengen von $X$ als $X \ne \emptyset$ definiert ist.
Wähle $\mathcal{B}_1, \ldots, \mathcal{B}_n \in \mathfrak{B}$ mit 
$$ \forall_{\nu \in \{1, \ldots,n\}} \, B_{\nu} \in \mathcal{B}_{\nu}. $$
Als Teilmenge von $\mathfrak{B}$ ist $\{\mathcal{B}_1, \ldots, \mathcal{B}_n\}$ eine bzgl.\ $\subset$ totalgeordnete Menge, die zusätzlich endlich ist.
Daher existiert offenbar $n_0 \in \{1, \ldots, n\}$ mit
$$ \forall_{\nu \in \{1, \ldots, n\}} \, \mathcal{B}_{\nu} \subset \mathcal{B}_{n_0}, $$
d.h.\
$$ \forall_{\nu \in \{1, \ldots, n\}} \, B_{\nu} \in \mathcal{B}_{n_0}. $$
Wegen $\mathcal{B}_{n_0} \in \mathfrak{B} \subset \mathfrak{A}$ ist $\mathcal{B}_{n_0}$ ein zentriertes System, also folgt $B_1 \cap \ldots \cap B_n \ne \emptyset$.
\q

\begin{Lemma} \label{FA.C.4} $\,$

\noindent \textbf{Vor.:} Seien $X$ eine nicht-leere Menge und $\mathcal{A} \subset \mathfrak{P}(X)$ ein maximales zentriertes System.

\noindent \textbf{Beh.:} 
\begin{itemize}
\item[(i)] Sind $n \in \N_+$ und $A_1, \ldots, A_n \in \mathcal{A}$, so gilt $A_1 \cap \ldots \cap A_n \in \mathcal{A}$.
\item[(ii)] Ist $B \in \mathfrak{P}(X) \setminus \{\emptyset\}$ mit
\begin{equation} \label{FA.C.4.S}
\forall_{A \in \mathcal{A}} \, B \cap A \ne \emptyset,
\end{equation}
so gilt $B \in \mathcal{A}$.
\end{itemize}
\end{Lemma}

\textit{Beweis.} Zu (i): Seien $n \in \N_+$, $A_1, \ldots, A_n \in \mathcal{A}$ sowie $B := A_1 \cap \ldots \cap A_n$ und
$$ \mathcal{B} := \mathcal{A} \cup \{B\}. $$
Wir weisen nach, daß $\mathcal{B}$ ein zentriertes System ist.
Dann folgt aus der Maximalität von $\mathcal{A}$, daß gilt $\mathcal{B} = \mathcal{A}$, also $B \in \mathcal{A}$.

Ohne Einschränkung seien daher $m \in \N_+$ und $B_1, \ldots, B_m \in \mathcal{B}$.

Im Falle $\forall_{\mu \in \{1, \ldots, m\}} \, B_{\mu} \ne B$ gilt $B_1, \ldots, B_m \in \mathcal{A}$ und somit $B_1 \cap \ldots \cap B_m \ne \emptyset$, weil $\mathcal{A}$ ein zentriertes System ist.
 
Im Falle $\exists_{m_0 \in \{1, \ldots, m\}} \, B_{m_0} = B$ können wir nach eventueller Umnumerierung der $B_1, \ldots, B_m$ ohne Beschränkung der Allgemeinheit annehmen, daß gilt
$$ B_1 = \ldots = B_{m_0} = B ~~ \mbox{ und } ~~ B_{m_0 + 1}, \ldots, B_m \in \mathcal{A}, $$
also folgt
$$ B_1 \cap \ldots \cap B_m = B \cap B_{m_0 + 1} \cap \ldots \cap B_m = A_1 \cap \ldots \cap A_n \cap B_{m_0 + 1} \cap \ldots \cap B_m \ne \emptyset, $$
da es sich hierbei um einen Schnitt von Elementen des zentrierten Systemes $\mathcal{A}$ handelt.

Zu (ii): Sei $B \in \mathfrak{P}(X) \setminus \{\emptyset\}$ mit (\ref{FA.C.4.S}).
Analog zum Beweis von (i) genügt es zu zeigen, daß
$$ \mathcal{B} := \mathcal{A} \cup \{B\} $$
ein zentriertes System ist.

Ohne Einschränkung seien also $m \in \N_+$ und $B_1, \ldots, B_m \in \mathcal{B}$.

Im Falle $\forall_{\mu \in \{1, \ldots, m\}} \, B_{\mu} \ne B$ gilt wieder $B_1, \ldots, B_m \in \mathcal{A}$ und $B_1 \cap \ldots \cap B_m \ne \emptyset$, weil $\mathcal{A}$ ein zentriertes System ist.

Im Falle $\exists_{m_0 \in \{1, \ldots, m\}} \, B_{m_0} = B$ können wir erneut ohne Beschränkung der Allgemeinheit annehmen, daß gilt
$$ B_1 = \ldots = B_{m_0} = B ~~ \mbox{ und } ~~ B_{m_0 + 1}, \ldots, B_m \in \mathcal{A}. $$
Hieraus folgt zunächst
$$ B_{m_0 + 1} \cap \ldots \cap B_m = \left\{ \begin{array}{cc} \stackrel{(i)}{\in} \mathcal{A}, & \text{ falls } m_0 < m, \\ = X, & \text{ falls } m_0 = m, \end{array} \right. $$
und sodann wegen (\ref{FA.C.4.S}) bzw.\ $B \ne \emptyset$
$$ B_1 \cap \ldots \cap B_m = B \cap B_{m_0 + 1} \cap \ldots \cap B_m\ne \emptyset. $$
Damit ist das Lemma bewiesen.
\q

\pagebreak
\textit{Beweis des Hauptsatzes.} Es genügt offenbar, den Fall $I \ne \emptyset$ zu betrachten.
Seien dann $X_i$ für jedes $i \in I$ ein nicht-leerer kompakter topologischer Raum und
$$ X := \bigtimes_{\iota \in I} X_{\iota} \ne \emptyset. $$
Ferner bezeichne
$$ \pi_i \: X \longrightarrow X_i $$
für jedes $i \in I$ die kanonische Projektion wie in (\ref{FA.C.P.S}).
Wir versehen $X$ mit der Produkttoplogie und behaupten, daß für jedes zentrierte System $\mathcal{A} \subset \mathfrak{P}(X)$ gilt
$$ \bigcap_{A \in \mathcal{A}} \overline{A} \ne \emptyset. $$
Hieraus folgt dann offenbar mittels \ref{FA.1.22.Durchschnitt} ,,(ii) $\Rightarrow$ (i)`` die Kompaktheit von $X$.

Seien $\mathcal{A}$ wie oben und $\mathcal{A}_{\max}$ ein gemäß \ref{FA.C.3} existierendes maximales zentriertes System mit $\mathcal{A} \subset \mathcal{A}_{\max}$.
Wir zeigen sogar
$$ \bigcap_{A \in \mathcal{A}_{\max}} \overline{A} \ne \emptyset $$
und können $\mathcal{A}_{\max} \ne \emptyset$ annehmen, denn sonst ist die letzte Aussage trivial.

Sei $i \in I$.
Dann ist
$$ \{ \pi_i(A) \, | \, A \in \mathcal{A}_{\max} \} $$
ein zentriertes System, denn für alle $n \in \N$ und $A_1, \ldots, A_n \in \mathcal{A}_{\max}$ gilt
$$ \pi_i(A_1) \cap \ldots \cap \pi_i(A_n) \supset \pi_i \left( A_1 \cap \ldots ... \cap A_n \right) \ne \emptyset, $$
weil $\mathcal{A}_{\max}$ ein zentriertes System ist.
Da $X_i$ des weiteren ein nicht-leerer kompakter topologischer Raum ist, folgt hieraus und aus \ref{FA.1.22.Durchschnitt} ,,(i) $\Rightarrow$ (ii)`` die Existenz eines $x_i \in X_i$ mit
\begin{equation} \label{FA.C.T.0}
x_i \in \bigcap_{A \in \mathcal{A}_{\max}} \pi_i(A) ~ \subset \bigcap_{A \in \mathcal{A}_{\max}} \overline{\pi_i(A)}.
\end{equation}

Wir setzen nun
$$ x := (x_i)_{i \in I} \in X $$
und zeigen schließlich
\begin{equation} \label{FA.C.T.1}
\forall_{A \in \mathcal{A}_{\max}} \, x \in \overline{A},
\end{equation}
womit der Satz von \textsc{Tychonoff} bewiesen ist.

Zum Nachweis von (\ref{FA.C.T.1}) weisen wir zunächst nach, daß gilt
\begin{equation}
\forall_{i \in I} \forall_{V_i \in {\rm Top}(X_i)} \, \left( x \in \overline{\pi_i}^1(V_i) \Longrightarrow \forall_{A \in \mathcal{A}_{\max}} \, \overline{\pi_i}^1(V_i) \cap A \ne \emptyset \right). \label{FA.C.T.2} 
\end{equation}

{[} Zu (\ref{FA.C.T.2}): Es seien $i \in I$, $V_i \in {\rm Top}(X_i)$ mit $x \in \overline{\pi_i}^1(V_i)$ und $A \in \mathcal{A}_{\max}$.
Dann gilt $V_i \in \U(x_i,X_i)$ und wegen (\ref{FA.C.T.0}) außerdem $x_i \in \overline{\pi_i(A)}$.
\ref{FA.1.8} (ii) impliziert $V_i \cap \pi_i(A) \ne \emptyset$, also $\overline{\pi_i}^1(V_i) \cap A \ne \emptyset$. {]}

Aus (\ref{FA.C.T.2}) und \ref{FA.C.4} (ii) folgt
$$ \forall_{i \in I} \forall_{V_i \in {\rm Top}(X_i)} \, \left( x \in \overline{\pi_i}^1(V_i) \Longrightarrow \overline{\pi_i}^1(V_i) \in \mathcal{A}_{\max} \right), $$
also ergibt \ref{FA.C.4} (i), daß jeder positiv endliche Schnitt von Mengen der Form
$$ \overline{\pi_j}^1(V_j) \in \U(x,X), \mbox { wobei } V_j \in {\rm Top}(X_j) \mbox{ sei,} $$
ein Element von $\mathcal{A}_{\max}$ ist.
Nun ist $\mathcal{A}_{\max}$ ein zentriertes System, und somit folgt aus der Definition der Produkttopologie offenbar
$$ \forall_{A \in \mathcal{A}_{\max}} \forall_{U \in \U(x,X)} \, U \cap A \ne \emptyset, $$
d.h.\ nach \ref{FA.1.8} (ii): $\forall_{A \in \mathcal{A}_{\max}} \, x \in \overline{A}$. \q

\cleardoublepage

\section{Das Lemma von \textsc{Urysohn}} \label{FAnaB}

Üblicherweise wird das Lemma von \textsc{Urysohn} für normale topologische Räume formuliert, worauf wir unten eingehen.
In Kapitel \ref{FAna3} benötigen wir eine Version für lokal-kompakte Hausdorff-Räume, die wir nun beweisen.

\begin{HS}[Lemma von \textsc{Urysohn} in der Version für lokal-kompakte Haus\-dorff-Räu\-me] \index{Lemma!von \textsc{Urysohn}} \label{FA.B.1}$\,$ 

\noindent \textbf{Vor.:} Seien $X$ ein lokal-kompakter Hausdorff-Raum sowie $K$ eine kompakte und $U$ eine offene Teilmenge von $X$ derart, daß gilt $K \subset U$.

\noindent \textbf{Beh.:} Es existiert $f \in \mathcal{C}_c(X,[0,1])$ mit $f|_K = 1$ und ${\rm Tr}(f) \subset U$.
\end{HS}

Entscheidend für den Beweis des Hauptsatzes ist die folgende Umformulierung von \ref{FA.1.47}.

\begin{Satz} \label{FA.B.2}
Es seien $X$ ein lokal-kompakter Hausdorff-Raum, $K$ eine kompakte und $U$ eine offene Teilmenge von $X$ derart, daß gilt $K \subset U$.

Dann existiert eine offene Teilmenge $V$ von $X$ mit $K \subset V \subset \overline{V} \subset U$ und $V \subset \subset X$. \q
\end{Satz}

\textit{Beweis des Hauptsatzes.} Wir konstruieren eine stetige Funktion $g \: X \to [0,1]$ so, daß $f := 1-g$ die Behauptung erfüllt.

Wenn ein solches $g$ existiert, so gilt für die offenen Mengen $U_{\rho} := \overline{g}^1({]} -\infty, \rho{[})$, wobei $\rho \in \R_+$: $\forall_{\rho, \sigma \in \R_+, \, \rho < \sigma} \, K \subset U_{\rho} \subset \overline{g}^1({]} -\infty, \rho{]}) \subset U_s$.
Hieraus folgt wegen der Abgeschlossenheit von $\overline{g}^1({]} -\infty, \rho{]})$ zum einen  
$$ \forall_{\rho, \sigma \in \R_+, \, \rho < \sigma} \, K \subset U_{\rho} \subset \overline{U_{\rho}} \subset \overline{g}^1({]} -\infty, \rho{]}) \subset U_{\sigma} ~ \wedge ~ \left( \rho > 1 \Rightarrow U_{\rho} = X \right) $$
und zum anderen $\forall_{x \in X} \, g(x) = \inf\{ \sigma \in \R_+ \, | \, x \in U_{\sigma} \}$.
Wir wollen nun umgekehrt für $\rho$ aus einer dichten Teilmenge $I$ von ${[}0, \infty{[}$ offene Teilmengen $U_{\rho}$ von $X$ mit
\begin{equation} \label{FA.B.1.0}
\forall_{\rho, \sigma \in I, \, \rho < \sigma} \, K \subset U_{\rho} \subset \overline{U_{\rho}} \subset U_{\sigma} ~ \wedge ~ \left( \rho > 1 \Rightarrow U_{\rho} = X \right)
\end{equation}
finden und $g$ durch $\forall_{x \in X} \, g(x) := \inf\{ \rho \in I \, | \, x \in U_{\rho} \}$ definieren.

Wegen \ref{FA.B.2} existiert eine relativ kompakte offene Teilmenge $U_1$ von $X$ mit $K \subset U_1 \subset \overline{U_1} \subset U$.
Wir wählen für jedes $\rho \in \{ \frac{k}{2^n} \, | \, n \in \N_+ \, \wedge \, k \in \{1, \ldots, 2^n-1\} \}$, d.i.\ die Menge der dyadisch endlich dargestellten (rationalen) Zahlen aus ${]}0,1{[}$, relativ kompakte offene Teilmengen $U_{\rho}$ von $X$ mit
\begin{equation} \label{FA.B.1.1}
\forall_{n \in \N_+} \, K \subset U_{\frac{1}{2^n}} \subset \overline{U_{\frac{1}{2^n}}} \subset U_{\frac{2}{2^n}} \subset \ldots \subset \overline{U_{\frac{2^n-2}{2^n}}} \subset U_{\frac{2^n-1}{2^n}} \subset \overline{U_{\frac{2^n-1}{2^n}}} \subset U_1
\end{equation}
wie folgt:

Im Falle $n=1$ ergibt \ref{FA.B.2} die Existenz einer relativ kompakten offenen Teilmenge $U_{\frac{1}{2}}$ mit (\ref{FA.B.1.1}).

Seien $n \in \N_+$ und relativ kompakte offene Teilmengen $U_{\frac{1}{2^n}}, \ldots, U_{\frac{2^n-1}{2^n}}$ mit (\ref{FA.B.1.1}) bereits gewählt, d.h.\
$$ K \subset U_{\frac{2}{2^{n+1}}} \subset \overline{U_{\frac{2}{2^{n+1}}}} \subset U_{\frac{4}{2^{n+1}}} \subset \ldots \subset \overline{U_{\frac{2^{n+1}-4}{2^{n+1}}}} \subset U_{\frac{2^{n+1}-2}{2^{n+1}}} \subset \overline{U_{\frac{2^{n+1}-2}{2^{n+1}}}} \subset U_1. $$
Dann liefert $2^n$-malige Anwendung von \ref{FA.B.2} die Existenz relativ kompakter offener Teilmengen $U_{\frac{1}{2^{n+1}}}, U_{\frac{3}{2^{n+1}}}, \ldots, U_{\frac{2^{n+1}-3}{2^{n+1}}}, U_{\frac{2^{n+1}-1}{2^{n+1}}}$ derart, daß (\ref{FA.B.1.1}) für $n+1$ anstelle von $n$ gilt.

Ferner seien $I := \{ \frac{k}{2^n} \, | \, n \in \N_+ \, \wedge \, k \in \{1, \ldots, 2^n-1\} \} \cup {[}1, \infty{[}$ und $U_{\rho} := X$ für $\rho \in {]}1,\infty{[}$.
Damit sind für alle $\rho \in I$ offene Teilmengen $U_{\rho}$ mit (\ref{FA.B.1.0}) definiert.

Nun können wir $g \: X \to \R$ durch $\forall_{x \in X} \, g(x) := \inf \{ \rho \in I \, | \, x \in U_{\rho} \}$ definieren.
Es gilt offenbar $g(X) \subset [0,1]$, 
$$ \forall_{x \in K} \forall_{\rho \in I} \, x \in U_{\rho}, $$
also $g|K = 0$, und
$$ \forall_{x \in X} \, \left( x \notin U_1 \Longrightarrow \left( \forall_{\rho \in I, \, \rho < 1} \, x \notin U_{\rho} \right) \Longrightarrow g(x) \ge 1 \Longrightarrow g(x) = 1 \right). $$
Somit ist $\overline{\{ x \in X \, | \, g(x) \ne 1 \}}$ eine abgeschlossene Teilmenge des Kompaktums $U_1$, also selbst kompakt.
Zu zeigen bleibt, daß $g$ stetig ist, denn $f := 1-g \: X \to [0,1]$ leistet dann das Gewünschte.

Wir zeigen zunächst für jedes $x \in X$ und $\rho \in I$
\begin{gather}
x \in \overline{U_{\rho}} \Longrightarrow g(x) \le \rho, \label{FA.B.1.2} \\
x \notin U_{\rho} \Longrightarrow g(x) \ge \rho. \label{FA.B.1.3}
\end{gather}

{[} Beweis hiervon: Aus $x \in \overline{U_{\rho}}$ folgt $\forall_{\sigma \in I, \, \sigma \ge \rho} \, x \in U_{\sigma}$, also $g(x) \le \rho$, und aus $x \notin U_{\rho}$ folgt $\forall_{\sigma \in I, \, \sigma \le \rho} \, x \notin U_{\sigma}$, also $g(x) \ge \rho$. {]}

Zum Nachweis der Stetigkeit von $g$ sei nun $x_0 \in X$.

1. Fall: $g(x_0) \ne 0$. 
Seien $a,b \in \R_+$ mit $a < g(x_0) < b$.
Gesucht ist $V \in \U(x_0,X)$ mit $g(V) \subset {]}a,b{[}$.
Da $I$ in ${[}0, \infty{[}$ dicht liegt, existieren $\rho, \sigma \in I$ derart, daß $a < \rho < g(x_0) < \sigma < b$ gilt, und $V := U_{\sigma} \setminus \overline{U_{\rho}} = U_{\sigma} \cap (X \setminus \overline{U_{\rho}})$ ist eine offene Teilmenge von $X$.
Wegen $g(x_0) < \sigma$, (\ref{FA.B.1.3}) und $g(x_0) > \rho$, (\ref{FA.B.1.2}) gilt $x_0 \in V$.
Ferner gilt $g(V) \subset {]}a,b{[}$, denn für jedes $x \in V$ folgt aus $x \in U_{\sigma} \subset \overline{U_{\sigma}}$ und (\ref{FA.B.1.2}), daß gilt $g(x) \le \sigma < b$, und aus $x \notin \overline{U_{\sigma}} \supset U_{\sigma}$ und (\ref{FA.B.1.3}), daß gilt $g(x) \ge \rho < a$.
Damit ist offenbar die Stetigkeit von $g \: X \to {[}0, \infty{[}$ in $x_0$ gezeigt.

2. Fall: $g(x_0) = 0$. 
Sei $b \in \R_+$.
Dann existiert $\sigma \in I$ mit $g(x_0) < \sigma < b$, und es gilt $U_{\sigma} \in \U(x_0,X)$ sowie $g(U_{\sigma}) \subset {[}0,b{[}$, also folgt die Stetigkeit von $g \: X \to {[}0, \infty{[}$ in $x_0$. \q

\begin{Kor} \label{FA.B.3}
Seien $X$ ein lokal-kompakter Hausdorff-Raum, $K$ eine kompakte und $B$ eine abgeschlossene Teilmenge von $X$ derart, daß $K \cap B = \emptyset$ gilt. 

Dann existiert $f \in \mathcal{C}_c(X,[0,1])$ mit $f|_K = 1$ und $f|_B = 0$. \q
\end{Kor}

\textit{Beweis.} Setze $U := X \setminus B$ in \ref{FA.B.1}. \q
\A
Wir kommen nun auf die Version des Urysohnschen Lemmas für normale topologische Räume zu sprechen.

\begin{Def}[Normale topologische Räume] \label{FA.B.4}
Ein topologischer Raum $X$ heißt \emph{normal}\index{Raum!topologischer!normaler} genau dann, wenn $X$ hausdorffsch ist und für alle disjunkten abgeschlossenen Teilmengen $A,B$ von $X$ disjunkte offene Teilmengen $U,V$ von $X$ mit $A \subset U$ und $B \subset V$ existieren.
\end{Def}

\begin{Bem*}
Das \emph{vierte Trennungsaxiom}\index{Trennungs!-axiome} für einen topologischen Raum besagt, daß zu je zwei disjunkten abgeschlossenen Teilmengen disjunkte offene Obermengen existieren.
\end{Bem*}

\begin{Lemma} \label{FA.B.5}
Sei $X$ ein Hausdorff-Raum.

Dann ist $X$ genau dann normal, wenn zu jeder abgeschlossenen Teilmenge $A$ von $X$ sowie jeder in $X$ offenen Obermenge $U$ von $A$ eine in $X$ offene Menge $V$ mit $A \subset V \subset \overline{V} \subset U$ existiert.
\end{Lemma}

\textit{Beweis.} ,,$\Rightarrow$`` Seien $A$ eine abgeschlossene und $U$ eine offene Teilmenge von $X$ mit $A \subset U$.
Dann ist $B := X \setminus U$ abgeschlossen, und es gilt $A \cap B = \emptyset$.
Nach Voraussetzung existieren disjunkte offene Mengen $V,W$ mit $A \subset V$ und $B \subset W$.
Es folgt $\overline{V} \cap B = \emptyset$, denn für jedes $x \in B$ gilt $W \in \U(x,X)$ und $V \cap W = \emptyset$.
Damit ist gezeigt $\overline{V} \subset X \setminus B = U$.

,,$\Leftarrow$`` Seien $A,B$ disjunkte abgeschlossene Teilmengen von $X$.
Dann ist $X \setminus B$ eine offene Obermenge von $A$, und nach Voraussetzung existiert eine offene Menge $V$ mit $A \subset V \subset \overline{V} \subset X \setminus B$.
$W := X \setminus \overline{V}$ ist dann offen und leistet $V \cap W = \emptyset$ sowie $B \subset W$. \q

\begin{HS}[Lemma von \textsc{Urysohn}] \index{Lemma!von \textsc{Urysohn}} \label{FA.B.6} $\,$ 

\noindent \textbf{Vor.:} Seien $X$ ein normaler topologischer Raum und $A,B$ disjunkte abgeschlossene Teilmengen von $X$.

\noindent \textbf{Beh.:} Es existiert $f \in \mathcal{C}(X,[0,1])$ mit $f|_A = 1$ und $f|_B = 0$.
\end{HS}

\textit{Beweis.} $U:= X \setminus B$ ist offene Obermenge von $A$.
Es genügt $f \in \mathcal{C}(X,[0,1])$ mit $f|_A = 1$ und ${\rm Tr}(f) \subset U$ zu finden.
Dies läuft analog zum Beweis von \ref{FA.B.1}, indem man $K$ durch $A$ ersetzt und anstelle von \ref{FA.B.2} auf \ref{FA.B.5} ,,$\Rightarrow$`` verweist. \q

\begin{Bem*} $\,$
\begin{itemize}
\item[1.)] Ein toplologischer Raum $X$, in dem zu je zwei disjunkten abgeschlossenen Teilmengen $A,B$ von $X$ eine stetige Funktion wie in der Behauptung des letzten Hauptsatzes existiert, heißt \emph{vollständig regulär}.
Jeder normale topologische Raum ist also vollständig regulär.
Dies wird beim Beweis des Metrisierbarkeitssatzes von \textsc{Urysohn}, der besagt, daß ein normaler topologischer Raum mit abzählbarer Basis der Topologie metrisierbar ist, benötigt.
Der Leser kann hierzu mehr in \cite[Chapter 4]{Munkres} finden.
\item[2.)] Jeder metrische Raum $(X,d)$ ist vollständig regulär.
Dies sieht man auch direkt ein, ohne das Lemma von \textsc{Urysohn} zu verwenden.
Sind nämlich $A,B$ disjunkte abgeschlossene Teilmengen von $X$, so wird durch
$$ \forall_{x \in X} \, f(x) := \frac{d(\{x\},B)}{d(\{x\},A) + d(\{x\},B)} $$
eine stetige Funktion $f \in \mathcal{C}(X,[0,1])$ mit $f|_A = 1$ und $f|_B = 0$ definiert.
\end{itemize}
\end{Bem*}

\cleardoublepage
\end{appendix}
\cleardoublepage
\addcontentsline{toc}{section}{Index} 
\printindex
\end{document}